\documentclass[12pt]{article}
\usepackage{amssymb,amsmath,amsthm,tikz,multirow,wasysym}
\usetikzlibrary{arrows,calc}

\title{Tilings of the Sphere by Congruent Pentagons IV: Edge Combination $a^4b$}
\author{Ho Man Cheung, Hoi Ping Luk, Min Yan\thanks{Research was supported by Hong Kong RGC General Research Fund 16303515.} \\
Hong Kong University of Science and Technology}

\usepackage[hidelinks, hyperindex]{hyperref}

\hyperref[sec:function]{}

\newcommand{\mc}{\mathcal}
\newcommand{\bb}{\mathbb}
\newcommand{\ssum}{\textstyle \sum}

\newcommand{\thin}{\hspace{0.1em}\rule{0.7pt}{0.8em}\hspace{0.1em}}
\newcommand{\thick}{\hspace{0.1em}\rule{1.5pt}{0.8em}\hspace{0.1em}}

\newcommand{\arcThroughThreePoints}[4][]{
\coordinate (middle1) at ($(#2)!.5!(#3)$);
\coordinate (middle2) at ($(#3)!.5!(#4)$);
\coordinate (aux1) at ($(middle1)!1!90:(#3)$);
\coordinate (aux2) at ($(middle2)!1!90:(#4)$);
\coordinate (center) at ($(intersection of middle1--aux1 and middle2--aux2)$);
\draw[#1] 
 let \p1=($(#2)-(center)$),
      \p2=($(#4)-(center)$),
      \n0={veclen(\p1)},       
      \n1={atan2(\y1,\x1)}, 
      \n2={atan2(\y2,\x2)},
      \n3={\n2>\n1?\n2:\n2+360}
    in (#2) arc(\n1:\n3:\n0);
}

\newtheorem{theorem}{Theorem}
\newtheorem{lemma}[theorem]{Lemma}

\newtheorem{proposition}{Proposition}
\newtheorem*{theorem*}{Theorem}

\theoremstyle{definition}
\newtheorem*{definition*}{Definition}

\newtheorem*{case*}{Case}

\newtheorem*{subcase*}{Subcase}

\newtheorem*{subsubcase*}{Subsubcase}

\theoremstyle{remark}

\numberwithin{equation}{section}

\begin{document}

\maketitle

\begin{abstract}
We classify edge-to-edge tilings of the sphere by congruent almost equilateral pentagons, in which four edges have the same length. Together with our earlier classifications of edge-to-edge tilings of the sphere by congruent equilateral pentagons of other types, and our classification of edge-to-edge tilings of the sphere by congruent quadrilaterals or triangles, we complete the classification of edge-to-edge tilings of the sphere by congruent polygons.

\end{abstract}

\tableofcontents

\section{Introduction}

\begin{figure}[htp]
\centering
\begin{tikzpicture}

\foreach \c in {0,1}
{
\begin{scope}[xshift=2.5*\c cm]

\foreach \a in {0,...,4}
\draw[rotate=72*\a]
	(18:0.8) -- (90:0.8);

\draw[line width=1.2]
	(-54:0.8) -- (234:0.8);
	
\foreach \a in {0,...,3}
\node at (-18+72*\a:0.8) {\small $a$};

\node at (-90:0.85) {\small $b$};

\end{scope}
}

\node at (90:0.5) {\small $\alpha$};
\node at (162:0.5) {\small $\beta$};
\node at (18:0.5) {\small $\gamma$};
\node at (234:0.5) {\small $\delta$};
\node at (-54:0.5) {\small $\epsilon$};

\begin{scope}[xshift=2.5cm]

\node at (90:0.5) {\small $\alpha$};
\node at (162:0.5) {\small $\gamma$};
\node at (18:0.5) {\small $\beta$};
\node at (234:0.5) {\small $\epsilon$};
\node at (-54:0.5) {\small $\delta$};

\node at (0,0) {\small $-$};

\end{scope}

\end{tikzpicture}
\caption{Almost equilateral pentagon.}
\label{pentagon}
\end{figure}
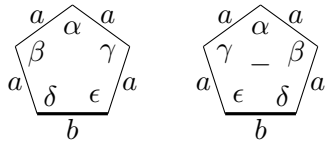

In an edge-to-edge tiling of the sphere by congruent pentagons, the pentagon has five possible edge combinations \cite[Lemma 9]{wy1}: $a^2b^2c,a^3bc,a^3b^2,a^4b,a^5$. We classified tilings for $a^2b^2c,a^3bc,a^3b^2,a^5$ in \cite{ay,awy,gsy,wy1,wy2}. In this paper, we classify tilings for the edge combination $a^4b$. We call such a pentagon {\em almost equilateral}. See Figure \ref{pentagon}, in which the two pentagons are congruent but with different orientations. The pentagon is positively oriented if $\alpha\to\beta$ is counterclockwise, and is negatively oriented (indicated by the $-$ sign) if $\alpha\to\beta$ is clockwise.

In an edge-to-edge tiling of the sphere by congruent polygons, the polygon must be triangle, quadrilateral, or pentagon. This paper and \cite{ay,awy,gsy,wy1,wy2} complete the classification pentagonal tilings. Together with the companion paper \cite{cly}, therefore, we completely classify all the edge-to-edge tilings of the sphere by congruent polygons.

The following is the main result of this paper. We denote the number of tiles by $f$.

\begin{theorem*}
Edge-to-edge tilings of the sphere by congruent pentagons with the edge combination $a^4b$ ($a,b$ distinct) are the following:
\begin{enumerate}
\item Three one parameter families {\rm $P_{\pentagon}P_4,P_{\pentagon}P_8,P_{\pentagon}P_{20}$} of the pentagonal subdivisions of the Platonic solids.
\item An earth map tiling {\rm $E_{\pentagon}1$} by symmetric pentagons for each $f\ge 12$ satisfying $f=0$ mod $4$, and two flip modifications {\rm $F_1E_{\pentagon}1,F_2E_{\pentagon}1$} for $f=4$ mod $8$.
\item A one parameter family of earth map tilings {\rm $E_{\pentagon}2$} for each $f\ge 12$ satisfying $f=0$ mod $4$, and a rotation modification {\rm $RE_{\pentagon}2$} for $f=4$ mod $8$, and the following modifications in two special cases:
\begin{itemize}
\item For $\alpha=\beta$ and $f=4$ mod $8$, two flip modifications {\rm $F_1E_{\pentagon}2$, $F_2E_{\pentagon}2$}; and {\rm $F'_2E_{\pentagon}2,F''_2E_{\pentagon}2$} obtained by further shift and flip of certain four-tile patch ${\mc A}$ in {\rm $F_2E_{\pentagon}2$}.
\item For $\alpha=\gamma$ and $f=16$, the flip modification {\rm $S_{16\pentagon}$} of the southern half tiling.
\end{itemize}
\end{enumerate}
\end{theorem*} 

The tilings are explained in Section \ref{construction}. The pentagonal subdivision tilings are are given in Figure \ref{subdivision_tiling}, and are the reductions of the pentagonal subdivision tilings introduced in \cite{wy1}. The two earth map tilings $E_{\pentagon}1,E_{\pentagon}2$ are given in Figure \ref{emt}, and are the repetitions of the timezones consisting of four shaded tiles in Figure \ref{emt}. Therefore the number $f$ of tiles in an earth map tiling is a multiple of $4$. 

The flip modifications $F_1E_{\pentagon}1,F_2E_{\pentagon}1$ are given in Figure \ref{emt1flip}. The rotation modification $RE_{\pentagon}2$ is given in Figure \ref{emt2mod2}. The flip modifications $F_1E_{\pentagon}2,F_2E_{\pentagon}2$ are given in Figure \ref{emt2mod3}. The tiling $S_{16\pentagon}$ is given in Figure \ref{emt2mod7}. 

The patch ${\mc A}$ is given in Figure \ref{emt2mod8}. The two flip modifications $F_1E_{\pentagon}2,F_2E_{\pentagon}2$ are distinguished by that $F_1E_{\pentagon}2$ does not contain ${\mc A}$, and $F_2E_{\pentagon}2$ contains a single ${\mc A}$. We may shift and flip ${\mc A}$ in $F_2E_{\pentagon}2$ to get $F'_2E_{\pentagon}2$ that still contains a single ${\mc A}$. On the other hand, if we first flip the single ${\mc A}$ in $F_2E_{\pentagon}2$, then we get a tiling with two disjoint ${\mc A}$. We may then independently shift and flip the two ${\mc A}$ to get $F''_2E_{\pentagon}2$ that still contains two disjoint ${\mc A}$.

The main theorem and \cite{ay,awy,gsy,wy1,wy2} give complete classification of all the edge-to-edge tilings of the sphere by congruent pentagons. The tilings can be divided into the Platonic type and the earth map type. In  the summary below, the pentagon is unique if we do not mention free parameters.
\begin{enumerate}
\item Platonic: 
\begin{itemize}
\item Pentagonal subdivisions $P_{\pentagon}P_4,P_{\pentagon}P_8,P_{\pentagon}P_{20}$ of the Platonic solids, which allow two free parameters, and has $f=12, 24, 60$; and flip modifications $FP_{\pentagon}P_8,F_1P_{\pentagon}P_{20},F_2P_{\pentagon}P_{20}$ for special edge combinations and angle values.  
\item Double pentagonal subdivisions $D_{\pentagon}P_8,D_{\pentagon}P_{20}$ of the Platonic solids, with $f=48, 120$.
\end{itemize}
\item Earth map:
\begin{itemize}
\item Earth map tiling $E_{\pentagon}1$ by symmetric pentagons, for each $f=0$ mod $4$; and flip modifications $F_1E_{\pentagon}1,F_2E_{\pentagon}1$, for $f=4$ mod 8.
\item Earth map tiling $E_{\pentagon}2$, which allows one free parameter, for each $f=0$ mod $4$; and a rotation modification $RE_{\pentagon}2$ for $f=4$ mod 8; and flip modifications $F_1E_{\pentagon}2,F_2E_{\pentagon}2$, for $f=4$ mod 8 and special angle values; and further shift and flip modifications $F'_2E_{\pentagon}2,F''_2E_{\pentagon}2$.
\item A sporadic tiling $S_{16\pentagon}$ with $16$ tiles, which can be regarded as a flip modification $S_{16\pentagon}$ of $E_{\pentagon}2$.
\end{itemize}
\end{enumerate}

Now we outline the content of the paper. 

Section \ref{construction} describes all the tilings in the the main theorem. At the end of the section, we also describe several tilings obtained in the paper, that are combinatorially possible but not geometrically possible. These tilings are not included in the main theorem.

Section \ref{basic} develops the tools for handling tilings by almost equilateral pentagons. The classifications for the edge combinations $a^2b^2c,a^3bc,a^3b^2$ in \cite{wy1,wy2} mainly use the strong combinatorial information provided by the variations of edge lengths. The classification for the edge combination $a^5$ (equilateral pentagon) in \cite{awy} cannot use any edge information, and is basically calculations in spherical trigonometry. The difficulty for the edge combination $a^4b$ is that, on the one hand, it does not have enough edge variation to provide sufficient combinatorial information. On the other hand, it is a bit more flexible than the equilateral case, that the pentagon cannot be easily determined by calculations in spherical trigonometry. Therefore more sophisticated technical tools in combinatorics and geometry are needed. The sophistication is illustrated by nineteen technical lemmas in Section \ref{basic}.

Sections \ref{symmetric_tiling}, \ref{distinct}, \ref{ndistinct} are the classification of tilings. Any pentagonal tiling of the sphere has at least two degree $3$ vertices with distinct angle combinations. We divide the classification into individual cases according to the angle combinations at degree $3$ vertices. Here we need to carefully select the most relevant pairs of degree $3$ vertices, and carry out the classifications in particular order. Moreover, we need to be concerned with the ambiguity caused by  distinct angles with the same angle values. As explained in Section \ref{count}, we deal with the concern by considering three broad cases: symmetric pentagon (Section \ref{symmetric_tiling}), distinct $\alpha,\beta,\gamma$ values (Section \ref{distinct}), and two of $\alpha,\beta,\gamma$ have the same value (Section \ref{ndistinct}). The whole classification is covered by thirty propositions.

\section{Construction}
\label{construction}

\subsubsection*{Pentagonal Subdivision}

Let $P_4,P_6,P_8,P_{12},P_{20}$ be the tetrahedron, cube, octahedron, dodecahedron, and icosahedron. The general pentagonal subdivision tilings $P_{\pentagon}P_n$ and double pentagonal subdivision tilings $D_{\pentagon}P_n$ are obtained in \cite{wy1}. The general combinatorial structures of the pentagonal and double pentagonal subdivisions are discussed in \cite{yan2}. 

The pentagonal and double pentagonal subdivisions of a Platonic solid and its dual are the same: $P_{\pentagon}P_6=P_{\pentagon}P_8$, $P_{\pentagon}P_{12}=P_{\pentagon}P_{20}$, $D_{\pentagon}P_6=D_{\pentagon}P_8$, $D_{\pentagon}P_{12}=D_{\pentagon}P_{20}$. Moreover, we remark that $P_{\pentagon}P_4=P_{12}$ is the deformed dodecahedron in \cite{ay,gsy}, and $D_{\pentagon}P_4$ is a reduction of $P_{\pentagon}P_8$. Therefore we only get three pentagonal subdivision tilings with $f=12,24,60$, and two double pentagonal subdivision tilings with $f=48,120$. 

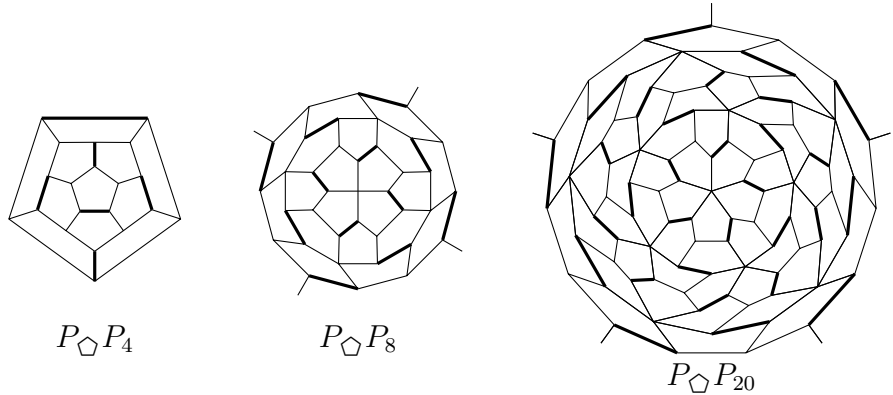
\begin{figure}[htp]
\centering
\begin{tikzpicture}[>=latex,scale=1]


\draw
	(234:0.33) -- (162:0.33) -- (90:0.33)
	(162:0.33) -- (162:0.65)
	(18:0.65) -- (54:0.8) -- (90:0.65)
	(54:0.8) -- (54:1.2)
	(-18:0.8) -- (-54:0.65) -- (-90:0.8)
	(-54:0.65) -- (-54:0.33)
	(126:1.2) -- (198:1.2) -- (270:1.2)
	(198:1.2) -- (198:0.8);

\draw
	(-54:0.33) -- (18:0.33) -- (90:0.33)
	(18:0.33) -- (18:0.65)
	(198:0.8) -- (234:0.65) -- (270:0.8)
	(234:0.33) -- (234:0.65)
	(-90:1.2) -- (-18:1.2) -- (54:1.2)
	(-18:1.2) -- (-18:0.8)
	(90:0.65) -- (126:0.8) -- (162:0.65)
	(126:0.8) -- (126:1.2);
	
\draw[line width=1.2]
	(-54:0.33) -- (234:0.33)
	(90:0.33) -- (90:0.65) 
	(18:0.65) -- (-18:0.8)
	(162:0.65) -- (198:0.8)
	(-90:1.2) -- (-90:0.8)
	(54:1.2) -- (126:1.2);

\node at (0,-2) {$P_{\pentagon}P_4$};


\foreach \a in {0,1,2,3}
{

\begin{scope}[xshift=3.5cm,rotate=90*\a]

\draw
	(0.4,0) -- (0.6,0.25) -- (0.25,0.6)
	(0.6,0.25) -- (15:1)
	(0:1.3) -- (30:1.3) -- (60:1.3)
	(45:1) -- (30:1.3);

\draw
	(0,0) -- (0.4,0)
	(-45:1) -- (-15:1) -- (15:1)
	(0.6,-0.25) -- (-15:1) -- (0:1.3)
	(60:1.3) -- (60:1.6);

\draw[line width=1.2]
	(0.25,0.6) -- (0,0.4)
	(15:1) -- (45:1)
	(60:1.3) -- (90:1.3);

\end{scope}
}

\node at (3.5,-2) {$P_{\pentagon}P_8$};


\foreach \a in {0,...,4}
{
\begin{scope}[xshift=8.2cm, rotate=72*\a]

\draw
	(0,0) -- (18:0.45) -- (36:0.7) -- (72:0.7)
	(0:0.7) -- (6:1.1)
	(30:1.1) -- (36:0.7)
	(-18:1.1) -- (6:1.1) -- (30:1.1)
	(-9:1.4) -- (6:1.1) -- (21:1.4) -- (39:1.4) -- (54:1.1)
	(-9:1.4) -- (0:1.6) -- (12:1.6)
	(39:1.4) -- (48:1.6) 
	(12:1.6) -- (30:1.9) -- (48:1.6)
	(54:1.9) -- (30:1.9) -- (6:1.9)
	(-18:1.9) -- (0:1.6)
	(6:1.9) -- (-6:2.2)
	(30:1.9) -- (42:2.2)
	(-30:2.2) -- (-6:2.2) -- (18:2.2) -- (18:2.5);

\draw[line width=1.2]
	(18:0.45) -- (0:0.7)
	(30:1.1) -- (54:1.1)
	(-9:1.4) -- (-24:1.6)
	(12:1.6) -- (21:1.4)
	(6:1.9) -- (-18:1.9)
	(-54:2.2) -- (-30:2.2);

\draw
	(0,0) -- (18:0.45)
	(0:0.7) -- (6:1.1)
	(-18:1.1) -- (6:1.1) -- (30:1.1)
	(-9:1.4) -- (6:1.1) -- (21:1.4)
	(30:1.9) -- (42:2.2)
	(12:1.6) -- (30:1.9) -- (48:1.6)
	(54:1.9) -- (30:1.9) -- (6:1.9)
	(18:2.2) -- (18:2.5);

\end{scope}
}

\node at (8.2,-2.5) {$P_{\pentagon}P_{20}$};
	
\end{tikzpicture}
\caption{Pentagonal subdivision tilings for $a^4b$.}
\label{subdivision_tiling}
\end{figure}

In general, the pentagons in pentagonal subdivision tilings have edge length combination $a^2b^2c$. When $a=b$, the pentagon becomes almost equilateral, and we get the pentagonal subdivision tilings in Figure \ref{subdivision_tiling} (all tiles have the same orientation), which are the ones in the main theorem.

The two free parameters for pentagonal subdivision tilings are described in \cite{wy2}. The reduction $a=b$ cuts two free parameters down to one free parameter.

Next we explain the notations $FP_{\pentagon}P_8,F_1P_{\pentagon}P_{20},F_2P_{\pentagon}P_{20}$ in the whole classification in the introduction. The flip modifications apply to the edge combination $a^3b^2$ and assume special angle values. In fact, the flip is applied to a cluster of tiles \cite[Figure 5]{wy2}: In $FP_{\pentagon}P_8$, the flip is applied to the cluster of four faces of the octahedron $P_8$ around a vertex. In $F_1P_{\pentagon}P_{20},F_2P_{\pentagon}P_{20}$, under two different angle assumptions, the flip is applied to the cluster of five faces of the icosahedron $P_{20}$ around a vertex. If there are several disjoint clusters in the tiling, then we may apply the flip to the clusters simultaneously. 

The pentagonal subdivision $P_{\pentagon}P_8$ is the union of two disjoint clusters. If we apply the flip to one cluster, then we get the tiling $T(4\beta\gamma\epsilon^2,2\epsilon^4)$ in the notation of \cite{wy2}. If we apply the flip to both clusters, then we get $T(4\beta^2\gamma^2,2\epsilon^4)$. Both are illustrated in \cite[Figure 6]{wy2}.

In the pentagonal subdivision $P_{\pentagon}P_{20}$, we may pick only one cluster, or two disjoint clusters (there are two possibilities), or three disjoint clusters \cite[Figures 7 and 8]{wy2}. Therefore each of $F_1P_{\pentagon}P_{20},F_2P_{\pentagon}P_{20}$ consists of four tilings:
\begin{itemize}
\item $F_1P_{\pentagon}P_{20}$: $T(5\beta\gamma\epsilon^3,7\epsilon^5)$, $T(10\beta\gamma\epsilon^3,2\epsilon^5)$, $T(2\beta^2\gamma^2\epsilon, 6\beta\gamma\epsilon^3,4\epsilon^5)$,
$T(6\beta^2\gamma^2\epsilon, 3\beta\gamma\epsilon^3,3\epsilon^5)$.
\item $F_2P_{\pentagon}P_{20}$: $T(5\beta\gamma\epsilon^2,5\delta\epsilon^3,7\epsilon^5)$,
$T(10\beta\gamma\epsilon^2,10\delta\epsilon^3,2\epsilon^5)$, $T(10\beta\gamma\epsilon^2,6\delta\epsilon^3,4\epsilon^5)$, 
$T(15\beta\gamma\epsilon^2,3\delta\epsilon^3,3\epsilon^5)$. 
\end{itemize}

\subsubsection*{Earth Map Tilings and Flip Modifications}

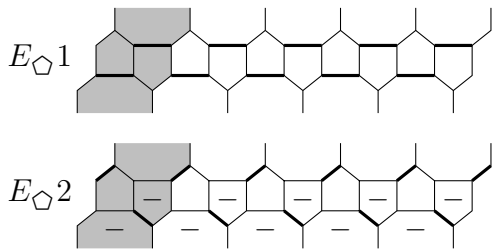
\begin{figure}[htp]
\centering
\begin{tikzpicture}[>=latex]


\foreach \b in {0,1}
{
\begin{scope}[yshift=1.8*\b cm]

\fill[gray!50]
	(-0.5,-0.7) -- (-0.5,-0.4) -- (-0.25,-0.2) -- (-0.25,0.2) -- (0,0.4) -- (0,0.7) -- (1,0.7) -- (1,0.4) -- (0.75,0.2) -- (0.75,-0.2) -- (0.5,-0.4) -- (0.5,-0.7);
	
\foreach \a in {0,...,5}
\draw[xshift=\a cm]
	(-0.5,-0.7) -- (-0.5,-0.4) -- (-0.25,-0.2) -- (-0.25,0.2) -- (0,0.4) -- (0,0.7);
	
\foreach \a in {0,...,4}
\draw[xshift=\a cm]
	(-0.25,-0.2) -- (0.25,-0.2) -- (0.5,-0.4) 
	(0,0.4) -- (0.25,0.2) -- (0.75,0.2)
	(0.25,0.2) -- (0.25,-0.2);	
	
\end{scope}
}


\begin{scope}[yshift=1.8 cm]

\foreach \a in {0,...,4}
\draw[line width=1.2, xshift=\a cm]
	(-0.25,-0.2) -- (0.25,-0.2)
	(0.25,0.2) -- (0.75,0.2);

\node at (-1,0) {$E_{\pentagon}1$};

\end{scope}


\foreach \a in {0,...,5}
\draw[line width=1.2, xshift=\a cm]	
	(-0.25,0.2) -- ++(0.25,0.2);
	
\foreach \a in {0,...,4}
{
\draw[line width=1.2, xshift=\a cm]
	(0.25,-0.2) -- ++(0.25,-0.2);

\node at (0.5+\a, -0.05) {\small $-$};
\node at (\a, -0.45) {\small $-$};
}

\node at (-1,0) {$E_{\pentagon}2$};
	
\end{tikzpicture}
\caption{Earth map tilings for $a^4b$.}
\label{emt}
\end{figure}

There are two earth map tilings, given by Figure \ref{emt}. Combinatorially, they are the earth map tilings of distance 5 in \cite{yan1}. The shaded part that consists of four tiles is one timezone, and the earth map tiling is the repetition of the timezone. The number of timezones can be any integer $\ge 3$. Since each timezone consists of four tiles, the number of tiles $f\ge 12$ and satisfies $f=0$ mod $4$. All the edges pointing upwards converge to a single vertex (north pole), and all the edges pointing downwards converge to another vertex (south pole). The left and right of the picture are glued together to form the spherical tiling.

In $E_{\pentagon}1$, the pentagon is symmetric, with angle values
\[
\alpha=\tfrac{8}{f}\pi,\;
\beta=\gamma=(1-\tfrac{4}{f})\pi,\;
\delta=\epsilon=(\tfrac{1}{2}+\tfrac{2}{f})\pi.
\]
For each $f$, this determines a unique pentagon. Figure \ref{symmetricA} shows part of $E_{\pentagon}1$ with all angles indicated. 

In $E_{\pentagon}2$, the angles satisfy
\[
\alpha+\delta+\epsilon=2\pi,\;
\beta=(1-\tfrac{4}{f})\pi,\;
\gamma=\tfrac{8}{f}\pi.
\]
For each $f$, this determines a one parameter family of pentagons. Moreover, the pentagon is not symmetric for $f>12$. In the unlabelled tiles, the direction of $\alpha\to\beta$ is in counterclockwise. In the tiles labelled by $-$, the direction is clockwise. Figure \ref{4adeD1} shows part of $E_{\pentagon}2$ with all angles indicated. 

If $f=8q+4$, then both earth map tilings have flip modifications. In Figure \ref{emt1flip}, we divide $E_{\pentagon}1$ into two identical {\em half earth map tilings} ${\mc H}$. Each ${\mc H}$ has $\alpha^{q+1}$ at one end and $\alpha^q$ at the other end. The first of the second row is another view of $E_{\pentagon}1$, decomposed into inner ${\mc H}$ and outer ${\mc H}$, with the angles indicated around the boundary between the two halves. In the second of the second row, we indicate the angle values, with $\bar{\beta}=2\pi-\beta$. Then we see many possible flips and rotations of the inner ${\mc H}$ that still give tilings. Many of these modifications are the same, or equivalent to exchanging the inner and outer ${\mc H}$. At the end, we only get two distinct flip modifications $F_1E_{\pentagon}1,F_2E_{\pentagon}1$ with respect to the two gray lines.

\begin{figure}[htp]
\centering
\begin{tikzpicture}[>=latex]

\foreach \b in {1,-1}
{
\begin{scope}[scale=\b]

\draw[<->]
	(0,-0.9) -- node[fill=white] {\small $\alpha^{q+1}$} ++(4,0);

\draw[<->]
	(0.5,0.9) -- node[fill=white] {\small $\alpha^q$} ++(3,0);
	
\foreach \a in {1,...,4}
\draw[xshift=-0.5cm+\a cm]
	(-0.5,-0.7) -- (-0.5,-0.4) -- (-0.25,-0.2) -- (0.25,-0.2) -- (0.5,-0.4) -- (0.5,-0.7)
	(-0.25,-0.2) -- (-0.25,0.2) -- (0,0.4) -- (0.25,0.2) -- (0.25,-0.2)
	(0,0.4) -- (0,0.7);
	
\foreach \a in {1,...,4}
	\draw[line width=1.2, xshift=-0.5cm+\a cm]
(-0.25,-0.2) -- (0.25,-0.2);

\foreach \a in {1,...,3}
\draw[line width=1.2, xshift=-0.5cm+\a cm]
	(0.25,0.2) -- (0.75,0.2);

\node at (4.3,0) {\small ${\mc H}$};

\end{scope}
}

\begin{scope}[yshift=-3.3cm]

\foreach \a in {0,...,9}
\foreach \x in {1,-1}
\draw[xshift=2.5*\x cm, rotate=36*\a]
	(-18:1.6) -- (18:1.6);

	
\begin{scope}[xshift=-2.5cm]

\foreach \x in {1,-1}
{
\begin{scope}[xscale=\x]

\draw[line width=1.2]
	(18:1.6) -- (18:2)
	(-18:1.6) -- (-18:1.2);

\draw
	(54:1.6) -- (54:1.2)
	(-54:1.6) -- (-54:2);

\node at (-62:1.8) {\small $\beta$};
\node at (-48:1.8) {\small $\alpha$};
\node at (-18:1.8) {\small $\beta$};
\node at (12:1.8) {\small $\delta$};
\node at (24:1.8) {\small $\delta$};
\node at (54:1.8) {\small $\beta$};

\node at (-54:1.35) {\small $\beta$};
\node at (-27:1.35) {\small $\delta$};
\node at (-9:1.35) {\small $\delta$};
\node at (16:1.35) {\small $\beta$};
\node at (63:1.35) {\small $\beta$};
\node at (46:1.35) {\small $\alpha$};

\end{scope}
}

\node at (90:1.8) {\small $\alpha^{q+1}$};
\node at (-90:1.8) {\small $\alpha^q$};

\node at (-90:1.25) {\small $\alpha^{q+1}$};
\node at (90:1.35) {\small $\alpha^q$};

\node at (0,0) {\small ${\mc H}$};
\node at (-35:2.2) {\small ${\mc H}$};

\end{scope}

\begin{scope}[xshift=2.5cm]

\draw[gray]
	(54:2.2) -- (54:-2.2)
	(18:2.2) -- (18:-2.2);

\node at (54:2.4) {\small $1$};
\node at (18:2.4) {\small $2$};

\foreach \a in {0,...,4}
{
\node at (-18+72*\a:1.8) {\small $\beta$};
\node at (18+72*\a:1.85) {\small $\bar{\beta}$};
\node at (18+72*\a:1.35) {\small $\beta$};
\node at (-18+72*\a:1.3) {\small $\bar{\beta}$};
}

\end{scope}

\end{scope}

\end{tikzpicture}
\caption{Flips modifications $F_1E_{\pentagon}1,F_2E_{\pentagon}1$ of $E_{\pentagon}1$.}
\label{emt1flip}
\end{figure}

We remark that a symmetric pentagon can be divided into two congruent quadrilaterals. If we divide all the symmetric tiles in this way, then we get the {\em simple quadrilateral subdivision} $E_{\square}2$ of $E_{\pentagon}1$ in \cite{cly}.

Similarly, we may divide $E_{\pentagon}2$ into two identical half earth map tilings in Figure \ref{emt2mod1}. In the second of the second row, we indicate the angle values between the two halves, with $\bar{\theta}=2\pi-\theta$. Then we may flip the inner half with respect to the gray line to get a new tiling.

\begin{figure}[htp]
\centering
\begin{tikzpicture}[>=latex]

\foreach \b in {1,-1}
{
\begin{scope}[scale=\b]

\draw[<->]
	(0,-0.9) -- node[fill=white] {\small $\gamma^{q+1}$} ++(4,0);

\draw[<->]
	(0.5,0.9) -- node[fill=white] {\small $\gamma^q$} ++(3,0);
	
\foreach \a in {1,...,4}
\draw[xshift=-0.5cm+\a cm]
	(-0.5,-0.7) -- (-0.5,-0.4) -- (-0.25,-0.2) -- (0.25,-0.2) -- (0.5,-0.4) -- (0.5,-0.7)
	(-0.25,-0.2) -- (-0.25,0.2) -- (0,0.4) -- (0.25,0.2) -- (0.25,-0.2)
	(0,0.4) -- (0,0.7);

\foreach \a in {1,...,3}
\draw[xshift=-0.5cm+\a cm]
	(0.25,0.2) -- (0.75,0.2);
	
\node at (4.3,0) {\small ${\mc H}$};

\end{scope}
}

\foreach \a in {1,...,4}
{
\draw[line width=1.2, xshift=-0.5cm+\a cm]
	(0.25,-0.2) -- (0.5,-0.4)
	(-0.25,0.2) -- (0,0.4);

\draw[line width=1.2, xshift=-5cm+\a cm]
	(0.25,-0.2) -- (0.5,-0.4)
	(0.75,0.2) -- (1,0.4);

\node at (-4.5+\a, -0.05) {\small $-$};
\node at (-0.5+\a, -0.45) {\small $-$};

}

\foreach \a in {1,...,3}
{

\node at (0+\a, -0.05) {\small $-$};
\node at (-4+\a, -0.45) {\small $-$};

}


\begin{scope}[yshift=-3.3cm]

\foreach \a in {0,...,9}
\foreach \x in {1,-1}
\draw[xshift=2.5*\x cm, rotate=36*\a]
	(-18:1.6) -- (18:1.6);

	
\begin{scope}[xshift=-2.5cm]

\draw
	(18:1.6) -- (18:2)
	(-54:1.6) -- (-54:2)
	(162:1.6) -- (162:2)
	(126:1.6) -- (126:1.2)
	(-18:1.6) -- (-18:1.2)
	(198:1.6) -- (198:1.2);
	
\draw[line width=1.2]
	(234:1.6) -- (234:2)
	(54:1.6) -- (54:1.2)
	(-18:1.6) -- (-54:1.6)
	(162:1.6) -- (126:1.6);

\node at (-60:1.8) {\small $\alpha$};
\node at (-48:1.8) {\small $\delta$};
\node at (-18:1.8) {\small $\epsilon$};
\node at (12:1.8) {\small $\gamma$};
\node at (24:1.8) {\small $\beta$};
\node at (54:1.8) {\small $\alpha$};

\node at (240:1.8) {\small $\epsilon$};
\node at (228:1.8) {\small $\delta$};
\node at (198:1.8) {\small $\beta$};
\node at (168:1.8) {\small $\alpha$};
\node at (156:1.8) {\small $\delta$};
\node at (126:1.8) {\small $\epsilon$};

\node at (-54:1.35) {\small $\epsilon$};
\node at (-27:1.35) {\small $\delta$};
\node at (-9:1.35) {\small $\alpha$};
\node at (16:1.35) {\small $\beta$};
\node at (63:1.35) {\small $\epsilon$};
\node at (46:1.35) {\small $\delta$};

\node at (234:1.35) {\small $\alpha$};
\node at (207:1.35) {\small $\beta$};
\node at (189:1.35) {\small $\gamma$};
\node at (164:1.35) {\small $\epsilon$};
\node at (117:1.35) {\small $\alpha$};
\node at (134:1.35) {\small $\delta$};

\node at (90:1.8) {\small $\gamma^{q+1}$};
\node at (-90:1.8) {\small $\gamma^q$};

\node at (-90:1.25) {\small $\gamma^{q+1}$};
\node at (90:1.35) {\small $\gamma^q$};

\node at (0,0) {\small ${\mc H}$};
\node at (-35:2.2) {\small ${\mc H}$};

\end{scope}

\foreach \b in {1,-1}
{
\begin{scope}[xshift=2.5cm, rotate=-36, xscale=\b]

\draw[gray]
	(90:2.2) -- (90:-2.2);

\draw[line width=1.2]
	(-18:1.6) -- (18:1.6);

\node at (54:1.35) {\small $\beta$};
\node at (54:1.85) {\small $\bar{\beta}$};
\node at (-54:1.3) {\small $\bar{\beta}$};
\node at (-54:1.8) {\small $\beta$};

\node at (18:1.35) {\small $\bar{\epsilon}$};
\node at (18:1.8) {\small $\epsilon$};
\node at (-18:1.35) {\small $\epsilon$};
\node at (-18:1.8) {\small $\bar{\epsilon}$};

\end{scope}
}

\begin{scope}[xshift=2.5cm, rotate=-36]

\node at (90:1.35) {\small $\bar{\alpha}$};
\node at (90:1.8) {\small $\alpha$};
\node at (-90:1.35) {\small $\alpha$};
\node at (-90:1.8) {\small $\bar{\alpha}$};

\end{scope}

\end{scope}

\end{tikzpicture}
\caption{Flip modification of $E_{\pentagon}2$.}
\label{emt2mod1}
\end{figure}

The flip modification in Figure \ref{emt2mod1} is not the viewpoint adopted in the main theorem. We will take an alternative viewpoint, which is not only the viewpoint in the main theorem, and is also necessary for understanding the other modifications in case $\alpha=\beta$. 

The alternative viewpoint is based on a different choice of the timezone, indicated by the shaded region in Figure \ref{emt2mod2}. The meridian (boundary between timezones) consists of six edges of the same length $a$. Then $E_{\pentagon}2$ is divided into two identical half earth map tilings ${\mc P}^{q+1}_q$, each again with $\gamma^{q+1}$ at one end and $\gamma^q$ at the other end. The second row of Figure \ref{emt2mod2} shows the angles and the angle values along the boundary between the two halves, with $\bar{\theta}=2\pi-\theta$. Then we may rotate the inner half by $120^{\circ}$ or $240^{\circ}$, and still get a tiling. We remark that the two rotation modifications are changed to each other by exchanging the inner and outer halves. Therefore we only get one new tiling $RE_{\pentagon}2$. This is the same as the flip modification in Figure \ref{emt2mod1}.

\begin{figure}[htp]
\centering
\begin{tikzpicture}[>=latex]


\fill[gray!50,xshift=2.1 cm]
	(0,0.7) -- (0,0.4) -- (0.25,0.2) -- (0.25,-0.2) -- (-0.25,-0.2) -- (-0.5,-0.4) -- (-0.5,-0.7) -- (0.5,-0.7) -- (0.5,-0.4) -- (0.75,-0.2) -- (1.25,-0.2) -- (1.25,0.2) -- (1,0.4) -- (1,0.7);
	
\foreach \b in {1,-1}
{
\begin{scope}[xshift=-2.1*\b cm, yscale=\b]

\draw[<->]
	(-0.5,0.9) -- node[fill=white] {\small $\gamma^{q+1}$} ++(4,0);

\draw[<->]
	(0,-0.9) -- node[fill=white] {\small $\gamma^q$} ++(3,0);
	
\foreach \a in {0,...,3}
{
\begin{scope}[xshift=\a cm]

\draw
	(-0.5,0.7) -- (-0.5,0.4) -- (-0.25,0.2) -- (0.25,0.2) -- (0.5,0.4) -- (0.5,0.7)
	(0.25,0.2) -- (0.25,-0.2) -- (0.75,-0.2) -- (0.75,0.2) -- (0.5,0.4)
	(0.25,-0.2) -- (0,-0.4);
	
\draw[line width=1.2]	
	(0.25,0.2) -- (0.5,0.4);

\end{scope}
}
			
\foreach \a in {0,1,2}
{
\begin{scope}[xshift=\a cm]

\draw
	(0,-0.7) -- (0,-0.4) -- (0.25,-0.2) -- (0.75,-0.2) -- (1,-0.4) -- (1,-0.7);
	
\draw[line width=1.2]
	(0.75,-0.2) -- (1,-0.4);

\end{scope}
}

\end{scope}
}

\node at (6.4,0) {\small ${\mc P}^{q+1}_q$};
\node at (-3,0) {\small ${\mc P}^{q+1}_q$};

\foreach \a in {0,...,3}
{
\node at (2.6+\a, -0.05) {\small $-$};
\node at (2.1+\a, -0.45) {\small $-$};
}

\foreach \a in {0,...,2}
{
\node at (-1.1+\a, -0.05) {\small $-$};
\node at (-1.6+\a, -0.45) {\small $-$};
}


\begin{scope}[shift={(1.7 cm, -3.3cm)}]

\foreach \a in {0,...,11}
\foreach \x in {1,-1}
\draw[xshift=2.5*\x cm, rotate=30*\a]
	(0:1.6) -- (30:1.6);

\begin{scope}[xshift=-2.5 cm]

\draw
	(-30:1.6) -- (-30:2)
	(210:1.6) -- (210:2);
	
\draw[line width=1.2]
	(120:1.6) -- (120:2)
	(240:1.6) -- (240:2)
	(0:1.6) -- (0:2);

\node at (115:1.8) {\small $\epsilon$};
\node at (126:1.8) {\small $\delta$};
\node at (150:1.8) {\small $\beta$};
\node at (180:1.8) {\small $\alpha$};
\node at (204:1.8) {\small $\gamma$};
\node at (217:1.8) {\small $\beta$};
\node at (245:1.8) {\small $\epsilon$};
\node at (235:1.8) {\small $\delta$};

\node at (60:1.8) {\small $\alpha$};
\node at (30:1.8) {\small $\beta$};
\node at (-5:1.8) {\small $\epsilon$};
\node at (6:1.8) {\small $\delta$};
\node at (-25:1.8) {\small $\gamma$};
\node at (-37:1.8) {\small $\beta$};
\node at (-60:1.8) {\small $\alpha$};

\node at (90:1.8) {\small $\gamma^{q+1}$};
\node at (-90:1.8) {\small $\gamma^q$};

\draw
	(30:1.6) -- (30:1.2)
	(150:1.6) -- (150:1.2);
	
\draw[line width=1.2]
	(180:1.6) -- (180:1.2)
	(60:1.6) -- (60:1.2)
	(-60:1.6) -- (-60:1.2);

\node at (120:1.35) {\small $\alpha$};
\node at (141:1.35) {\small $\beta$};
\node at (158:1.35) {\small $\gamma$};
\node at (173:1.35) {\small $\epsilon$};
\node at (188:1.35) {\small $\delta$};
\node at (210:1.35) {\small $\beta$};
\node at (240:1.35) {\small $\alpha$};

\node at (67:1.35) {\small $\epsilon$};
\node at (53:1.35) {\small $\delta$};
\node at (37:1.35) {\small $\beta$};
\node at (22:1.35) {\small $\gamma$};
\node at (0:1.35) {\small $\alpha$};
\node at (-30:1.35) {\small $\beta$};
\node at (-67:1.35) {\small $\epsilon$};
\node at (-52:1.35) {\small $\delta$};

\node at (-90:1.3) {\small $\gamma^{q+1}$};
\node at (90:1.3) {\small $\gamma^q$};

\node at (0,0) {\small ${\mc P}^{q+1}_q$};
\node at (-40:2.4) {\small ${\mc P}^{q+1}_q$};

\end{scope}


\foreach \b in {0,1,2}
{
\begin{scope}[xshift=2.5 cm, rotate=120*\b]

\node at (0:1.35) {\small $\alpha$};
\node at (0:1.8) {\small $\bar{\alpha}$};
\node at (30:1.35) {\small $\bar{\beta}$};
\node at (30:1.8) {\small $\beta$};
\node at (-30:1.35) {\small $\beta$};
\node at (-30:1.8) {\small $\bar{\beta}$};
\node at (60:1.35) {\small $\bar{\alpha}$};
\node at (60:1.8) {\small $\alpha$};

\end{scope}
}

\begin{scope}[xshift=2.5 cm]

\draw[gray, ->]
	(-90:0.9) arc (-90:30:0.9);
\draw[gray, ->]
	(-90:0.6) arc (-90:150:0.6);

\end{scope}

\end{scope}

\end{tikzpicture}
\caption{Rotation modification $RE_{\pentagon}2$ of $E_{\pentagon}2$.}
\label{emt2mod2}
\end{figure}

In the earlier classification work, the earth map tiling appears only for the edge combination $a^5$ (equilateral pentagon) \cite{awy}. They are the $a=b$ reductions of $E_{\pentagon}2$ and $RE_{\pentagon}2$ (called flip modification in \cite{awy}, referring to Figure \ref{emt2mod1}). However, the reduction happens only for $f=16,20,24,24$, with two different equilateral pentagons for $f=24$. Therefore there are only four equilateral $E_{\pentagon}2$. Among these, only $f=20$ satisfies $f=4$ mod $8$. Therefore there is only one equilateral $RE_{\pentagon}2$.

\subsubsection*{Modifications of the Second Earth Map Tiling in Case $\alpha=\gamma$}

We consider $E_{\pentagon}2$ as the union of the upper half tiling consisting of positively oriented (i.e., unlabelled) tiles, and the lower half tiling consisting of negatively oriented (i.e., labelled with $-$) tiles. The first of the second row of Figure \ref{emt2mod7} indicates the angles along the boundary between the two halves. If $\alpha=\gamma$, then we get the angle values in the second of the second row, with $\bar{\alpha}=2\pi-\alpha$. Then the horizontal flip of the lower half still gives a tiling, which is the one on the right of the first row.

\begin{figure}[htp]
\centering
\begin{tikzpicture}[>=latex]


\foreach \a in {0,...,4}
{
\begin{scope}[xshift=\a cm]

\draw
	(-0.5,-0.7) -- (-0.5,-0.4) -- (-0.25,-0.2) -- (-0.25,0.2) -- (0,0.4) -- (0,0.7);

\draw[line width=1.2]	
	(-0.25,0.2) -- ++(0.25,0.2);
	
\end{scope}
}
	
\foreach \a in {0,...,3}
{
\begin{scope}[xshift=\a cm]

\draw
	(-0.25,-0.2) -- (0.25,-0.2) -- (0.5,-0.4) 
	(0,0.4) -- (0.25,0.2) -- (0.75,0.2)
	(0.25,0.2) -- (0.25,-0.2);	

\draw[line width=1.2]
	(0.25,-0.2) -- ++(0.25,-0.2);
	
\node at (0.5, -0.05) {\small $-$};
\node at (0, -0.45) {\small $-$};
	
\end{scope}
}


\draw[very thick,->]
	(4.2,0) -- ++(1.5,0);
	
\node at (4.9,0.3) {\small flip};
\node at (4.9,-0.3) {\small lower half};


\begin{scope}[xshift=6.5cm]

\foreach \a in {0,...,4}
{
\begin{scope}[xshift=\a cm]

\draw
	(-0.5,-0.7) -- (-0.5,-0.4) -- (-0.25,-0.2) -- (-0.25,0.2) -- (0,0.4) -- (0,0.7);

\draw[line width=1.2]	
	(-0.25,0.2) -- ++(0.25,0.2);
	
\end{scope}
}
	
\foreach \a in {0,...,3}
{
\begin{scope}[xshift=\a cm]

\draw
	(-0.25,-0.2) -- (0.25,-0.2) -- (0.5,-0.4) 
	(0,0.4) -- (0.25,0.2) -- (0.75,0.2)
	(0.25,0.2) -- (0.25,-0.2);	

\draw[line width=1.2]
	(-0.25,-0.2) -- ++(-0.25,-0.2);
	
\end{scope}
}

\end{scope}


\begin{scope}[shift={(0.2cm, -2cm)}]


\foreach \a in {0,1}
{
\begin{scope}[xshift=2*\a cm]

\draw
	(0,-0.8) -- (0,0.8)
	(1,-0.8) -- (1,0.8)
	(0,-0.4) -- (1,-0.4)
	(1,0.4) -- (2,0.4);
	
\draw[line width=1.2]	
	(0,0.4) -- ++(0,0.4)
	(1,-0.4) -- ++(0,-0.4);

\node at (0.8,-0.2) {\small $\alpha$};
\node at (0.8,-0.6) {\small $\delta$};
\node at (1.2,-0.4) {\small $\epsilon$};

\node at (0.8,0.4) {\small $\beta$};
\node at (1.2,0.6) {\small $\beta$};
\node at (1.2,0.2) {\small $\gamma$};

\node at (0.2,-0.2) {\small $\gamma$};
\node at (0.2,-0.6) {\small $\beta$};
\node at (-0.2,-0.4) {\small $\beta$};

\node at (0.2,0.4) {\small $\epsilon$};
\node at (-0.2,0.6) {\small $\delta$};
\node at (-0.2,0.2) {\small $\alpha$};

\end{scope}
}

\draw
	(-1,0.4) -- (0,0.4);


\begin{scope}[xshift=6.5cm]

\foreach \a in {0,1}
{
\begin{scope}[xshift=2*\a cm]

\draw
	(0,0.4) -- (0,-0.4) -- (1,-0.4) -- (1,0.4) -- (2,0.4);

\node at (0.8,-0.2) {\small $\alpha$};
\node at (0.2,-0.2) {\small $\alpha$};
\node at (1.1,-0.6) {\small $\bar{\alpha}$};
\node at (-0.1,-0.6) {\small $\bar{\alpha}$};

\node at (-0.2,0.2) {\small $\alpha$};
\node at (1.2,0.2) {\small $\alpha$};
\node at (0.1,0.6) {\small $\bar{\alpha}$};
\node at (0.9,0.6) {\small $\bar{\alpha}$};

\end{scope}
}

\draw
	(-1,0.4) -- (0,0.4);

\end{scope}

\end{scope}

\end{tikzpicture}
\caption{Flip modification $S_{16\pentagon}$ of $E_{\pentagon}2$ in case $\alpha=\gamma$.}
\label{emt2mod7}
\end{figure}

There is a catch in the flip modification above. In the case ``$\alpha\delta\epsilon,\alpha\gamma^2$ are vertices'' of Proposition \ref{4ade}, where $(\beta,\delta)$ and $(\gamma,\epsilon)$ are exchanged, the calculation using spherical trigonometry shows that the pentagon exists if and only if $f=16$. This means the flip modified tiling has only four timezones. We consider the tiling to be sporadic and denote it by $S_{16\pentagon}$.

\subsubsection*{Modifications of the Second Earth Map Tiling in Case $\alpha=\beta$}

For each $f=0$ mod $4$, there is a unique pentagon in $E_{\pentagon}2$ satisfying $\alpha=\beta$. If $f=4$ mod $8$ and $\alpha=\beta$, then $E_{\pentagon}2$ can be modified in many ways to give new tilings. 

By $\alpha=\beta$, the angle values in the second of the second row of Figure \ref{emt2mod2} becomes the first of Figure \ref{emt2mod3}. Then we see three possible flip modifications with respect to the three gray lines. 

\begin{figure}[htp]
\centering
\begin{tikzpicture}[>=latex]

\foreach \a in {0,...,11}
\foreach \x in {0,1,2}
\draw[xshift=4.5*\x cm, rotate=30*\a]
	(0:1.6) -- (30:1.6);


\foreach \b in {1,2,3}
{
\begin{scope}[rotate=120*\b]

\draw[gray!50]
	(45:2) -- (45:-2);
	
\node at (0:1.35) {\small $\alpha$};
\node at (0:1.8) {\small $\bar{\alpha}$};
\node at (30:1.35) {\small $\bar{\alpha}$};
\node at (30:1.8) {\small $\alpha$};
\node at (-30:1.35) {\small $\alpha$};
\node at (-30:1.8) {\small $\bar{\alpha}$};
\node at (60:1.35) {\small $\bar{\alpha}$};
\node at (60:1.8) {\small $\alpha$};

\end{scope}

\node at (-15+60*\b:2.2)	{\small $\b$};
}


\foreach \a in {1,2}
{
\begin{scope}[xshift=4.5*\a cm]

\draw
	(-30:1.6) -- (-30:2.3) -- (-15:2.5) -- 
	(0:2.3) to[out=90, in=90] (90:2) -- (90:1.6) 
	(210:1.6) -- (210:2)
	(-60:1.6) -- 
	(-60:1.2) to[out=120, in=255] (-15:0.7)
	(60:1.6) -- (60:0) -- (-15:0.7);
	
\draw[line width=1.2]
	(120:1.6) -- (120:2)
	(240:1.6) -- (240:2)
	(0:1.6) -- (0:2.3)
	(30:1.6) -- 
	(30:1.2) to[out=210, in=75] (-15:0.7);

\node at (115:1.8) {\small $\epsilon$};
\node at (126:1.8) {\small $\delta$};
\node at (150:1.8) {\small $\beta$};
\node at (180:1.8) {\small $\alpha$};
\node at (204:1.8) {\small $\gamma$};
\node at (217:1.8) {\small $\beta$};
\node at (245:1.8) {\small $\epsilon$};
\node at (235:1.8) {\small $\delta$};

\node at (60:1.8) {\small $\alpha$};
\node at (30:1.8) {\small $\beta$};
\node at (-5:1.8) {\small $\epsilon$};
\node at (6:1.8) {\small $\delta$};
\node at (-25:1.8) {\small $\gamma$};
\node at (-37:1.8) {\small $\beta$};
\node at (-60:1.8) {\small $\alpha$};

\node at (85:1.8) {\small $\gamma$};
\node at (97:1.8) {\small $\gamma^q$};
\node at (-90:1.8) {\small $\gamma^q$};

\node at (-51:1.35) {\small $\gamma$};
\node at (52:1.35) {\small $\gamma$};
\node at (-15:0.85) {\small $\epsilon$};
\node at (4:0.55) {\small $\delta$};
\node at (30:0.3) {\small $\beta$};
\node at (45:0.7) {\small $\alpha$};

\node at (-15:2.3) {\small $\beta$};
\node at (5:2.1) {\small $\epsilon$};
\node at (-5:2.15) {\small $\delta$};
\node at (-25:2.15) {\small $\alpha$};

\node[inner sep=0.5,draw,shape=circle] at (45:1.85) {\small 1};
\node[inner sep=0.5,draw,shape=circle] at (-15:1.15) {\small 2};
\node[inner sep=0.5,draw,shape=circle] at (45:1.05) {\small 3};
\node[inner sep=0.5,draw,shape=circle] at (-15:1.95) {\small 4};

\end{scope}
}

\foreach \a/\b in {2/150,1/30}
{
\begin{scope}[xshift=4.5*\a cm, rotate=\b]

\draw
	(30:1.6) -- (30:1.2)
	(150:1.6) -- (150:1.2);
	
\draw[line width=1.2]
	(0:1.6) -- (0:1.2)
	(120:1.6) -- (120:1.2)
	(240:1.6) -- (240:1.2);

\node at (60:1.35) {\small $\alpha$};
\node at (39:1.35) {\small $\beta$};
\node at (22:1.35) {\small $\gamma$};
\node at (7:1.35) {\small $\epsilon$};
\node at (-8:1.35) {\small $\delta$};
\node at (-30:1.35) {\small $\beta$};
\node at (-60:1.35) {\small $\alpha$};

\node at (112:1.35) {\small $\epsilon$};
\node at (128:1.35) {\small $\delta$};
\node at (142:1.35) {\small $\beta$};
\node at (159:1.35) {\small $\gamma$};
\node at (180:1.35) {\small $\alpha$};
\node at (210:1.35) {\small $\beta$};
\node at (247:1.35) {\small $\epsilon$};
\node at (232:1.35) {\small $\delta$};

\node at (90:1.3) {\small $\gamma^q$};

\end{scope}
}

\begin{scope}[xshift=4.5cm]

\node at (-71:1.35) {\small $\gamma^q$};
\node at (225:0.6) {flip 2};

\node[inner sep=0.5,draw,shape=circle] at (90:0.8) {\small 5};
\node[inner sep=0.5,draw,shape=circle] at (-45:2.1) {\small 6};

\end{scope}

\begin{scope}[xshift=9cm]

\node at (70:1.3) {\small $\gamma^q$};
\node at (225:0.6) {flip 3};

\end{scope}

\end{tikzpicture}
\caption{Flip modifications $F_1E_{\pentagon}2,F_2E_{\pentagon}2,F_3E_{\pentagon}2$ of $E_{\pentagon}2$ in case $\alpha=\beta$.}
\label{emt2mod3}
\end{figure}

\begin{figure}[htp]
\centering
\begin{tikzpicture}[>=latex,scale=1]

\foreach \a in {1,-1}
{

\begin{scope}[xshift=-1.7 cm, scale=\a]

\draw
	(-0.5,0.7) -- (0,1.1) -- (0.5,0.7) -- (0.5,-0.7)
	(0,0) -- (0.5,0)
	(0.5,0.7) -- (1.3,0.7) -- (1.3,-0.7) -- (0.5,-0.7);

\draw[line width=1.2]
	(-0.5,0) -- (-0.5,-0.7);	

\node at (0.3,0.2) {\small $\delta$};	
\node at (-0.3,0.2) {\small $\beta$}; 
\node at (0.3,0.6) {\small $\epsilon$};
\node at (-0.3,0.6) {\small $\alpha$};
\node at (0,0.85) {\small $\gamma$};

\node at (1.1,0.5) {\small $\beta$};
\node at (0.7,0.5) {\small $\delta$};
\node at (1.1,-0.5) {\small $\alpha$};
\node at (0.7,0) {\small $\epsilon$};
\node at (0.7,-0.5) {\small $\gamma$};

\end{scope}

\begin{scope}[xshift=1.7 cm, scale=\a]

\draw[gray!50]
	(0,0) -- (0,1.4);

\draw
	(-1.3,0.7) -- (-0.5,0.7) -- (0,1.1) -- (0.5,0.7) -- (1.3,0.7) -- (1.3,-0.7);

\node at (0,0.85) {\small $\gamma$};
\node at (0.4,0.5) {\small $\bar{\alpha}$};
\node at (-0.4,0.5) {\small $\bar{\alpha}$};
\node at (1.1,0.5) {\small $\alpha$};
\node at (1.1,-0.5) {\small $\alpha$};

\end{scope}

}

\begin{scope}[xshift=-1.7 cm]

\node[inner sep=0.5,draw,shape=circle] at (0,0.4) {\small $1$};
\node[inner sep=0.5,draw,shape=circle] at (0,-0.4) {\small $2$};
\node[inner sep=0.5,draw,shape=circle] at (-1,0) {\small $3$};
\node[inner sep=0.5,draw,shape=circle] at (1,0) {\small $4$};

\end{scope}

\end{tikzpicture}
\caption{Patch ${\mc A}$ and its flip.}
\label{emt2mod8}
\end{figure}
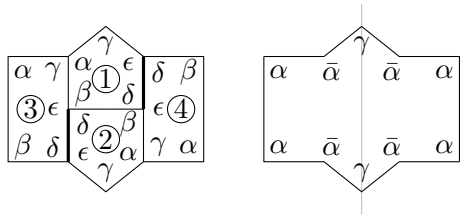

The second and third flip modifications $F_2E_{\pentagon}2,F_3E_{\pentagon}2$ in Figure \ref{emt2mod3} contain four tiles $T_1,T_2,T_3,T_4$ along the boundary between the two halves ${\mc P}^{q+1}_q$. They form the patch ${\mc A}$ in the first of Figure \ref{emt2mod8}. The second of Figure \ref{emt2mod8} shows the angle values along the boundary of ${\mc A}$, which implies that we may flip ${\mc A}$ with respect to the gray line to get a new tiling. The operation can be applied to $F_2E_{\pentagon}2$ and $F_3E_{\pentagon}2$, but cannot be applied to $F_1E_{\pentagon}2$ because it does not contain ${\mc A}$.

The flip modifications $F_2E_{\pentagon}2,F_3E_{\pentagon}2$ have further modifications, and are related by these modifications. Note that the two tilings are ${\mc P}^q_q\cup {\mc A}\cup {\mc P}^q_q$, where the left ${\mc P}^q_q$ (consisting of $q$ timezones given by the shaded region in Figure \ref{emt2mod2}) is the inner half minus $T_2,T_3$, and the right ${\mc P}^q_q$ is the outer half minus $T_1,T_4$. 

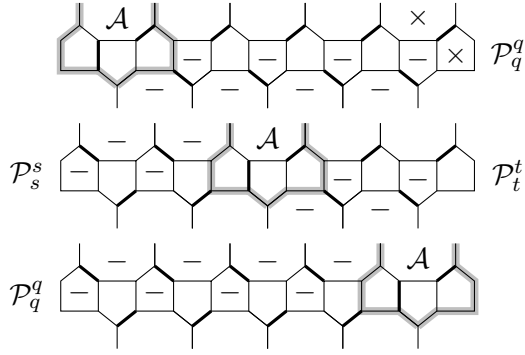
\begin{figure}[htp]
\centering
\begin{tikzpicture}[>=latex]

\foreach \a/\b in {0/0, 2/1, 4/2}
{
\begin{scope}[shift={(\a cm,-1.6*\b cm)}]

\foreach \c in {1,-1}
\draw[xscale=\c, gray!50, line width=3]
	(0,-0.4) -- (0.25,-0.2) -- (0.75,-0.2) -- (0.75,0.2) -- (0.5,0.4) -- (0.5,0.7);

\node at (0,0.5) {\small ${\mc A}$};

\end{scope}
}

\foreach \b in {0,1,2}
{
\begin{scope}[yshift=-1.6*\b cm]

\foreach \a in {0,...,4}
\draw[xshift=\a cm]
	(-0.5,0.7) -- (-0.5,0.4) -- (-0.25,0.2) -- (0.25,0.2) -- (0.5,0.4) -- (0.5,0.7)
	(-0.25,0.2) -- (-0.25,-0.2) -- (0,-0.4) -- (0.25,-0.2) -- (0.25,0.2);

\foreach \a in {0,...,3}
\draw[xshift=\a cm]
	(0,-0.7) -- (0,-0.4) -- (0.25,-0.2) -- (0.75,-0.2) -- (1,-0.4) -- (1,-0.7);

\draw
	(-0.5,0.4) -- (-0.75,0.2) -- (-0.75,-0.2) -- (-0.25,-0.2)
	(4.5,0.4) -- (4.75,0.2) -- (4.75,-0.2) -- (4.25,-0.2);

\end{scope}	
}


\foreach \a in {0,...,4}
\draw[line width=1.2, xshift=\a cm]	
	(0.25,0.2) -- ++(0.25,0.2);

\foreach \a in {1,...,4}
\draw[line width=1.2, xshift=\a cm]		
	(-0.25,-0.2) -- ++(0.25,-0.2);

\draw[line width=1.2]
	(-0.25,-0.2) -- ++(0,0.4);
	
\foreach \a in {0,...,3}
{
\node at (1+\a, -0.05) {\small $-$};
\node at (0.5+\a, -0.45) {\small $-$};
}

\node at (5.2,0) {\small ${\mc P}^q_q$};

\node at (4,0.5) {\small $\times$};
\node at (4.5,0) {\small $\times$};


\begin{scope}[yshift=-1.6 cm]
	
\foreach \a in {0,1,2}
\draw[line width=1.2, xshift=\a cm]		
	(2.25,0.2) -- ++(0.25,0.2);
	
\foreach \a in {0,1}
\draw[line width=1.2, xshift=\a cm]
	(0.25,-0.2) -- ++(-0.25,-0.2)
	(-0.25,0.2) -- ++(-0.25,0.2)
	(2.75,-0.2) -- ++(0.25,-0.2);

\draw[line width=1.2]
	(1.75,-0.2) -- ++(0,0.4);

\foreach \a in {0,1}
{
\node at (-0.5+\a, 0.05) {\small $-$};
\node at (\a, 0.45) {\small $-$};
\node at (3+\a, -0.05) {\small $-$};
\node at (2.5+\a, -0.45) {\small $-$};
}

\node at (-1.2,0) {\small ${\mc P}^s_s$};
\node at (5.2,0) {\small ${\mc P}^t_t$};

\end{scope}


\begin{scope}[yshift=-3.2 cm]
	
\foreach \a in {0,1,2,3}
\draw[line width=1.2, xshift=\a cm]
	(0.25,-0.2) -- ++(-0.25,-0.2)
	(-0.25,0.2) -- ++(-0.25,0.2);

\draw[line width=1.2]
	(3.75,-0.2) -- ++(0,0.4)
	(4.25,0.2) -- ++(0.25,0.2);

\foreach \a in {0,...,3}
{
\node at (-0.5+\a, 0.05) {\small $-$};
\node at (\a, 0.45) {\small $-$};
}

\node at (-1.2,0) {\small ${\mc P}^q_q$};

\end{scope}

\end{tikzpicture}
\caption{Shift ${\mc A}$ inside ${\mc P}^s_s\cup {\mc A}\cup {\mc P}^t_t$.}
\label{emt2mod4}
\end{figure}

The first of Figure \ref{emt2mod4} shows the {\em partial tiling} ${\mc A}\cup {\mc P}^q_q$ (this is the outer half plus $T_2,T_3$) in $F_2E_{\pentagon}2,F_3E_{\pentagon}2$, with ${\mc A}$ outlined by shaded edges. The partial tiling can be modified to ${\mc P}^s_s\cup {\mc A}\cup {\mc P}^t_t$ ($s+t=q$) in the second picture. This does not change the angle values on the boundary, and therefore gives new tilings. Visually this is the shift of ${\mc A}$ from one side of the partial tiling to the interior. We may shift ${\mc A}$ all the way to the other side and get ${\mc P}^q_q\cup{\mc A}$ in the third of Figure \ref{emt2mod4}. In addition to the shift, we may also flip ${\mc A}$.

We denote by $F'_2E_{\pentagon}2$ all the further modifications of $F_2E_{\pentagon}2$ by the shift in Figure \ref{emt2mod4} and the flip of ${\mc A}$. We remark that $F_3E_{\pentagon}2$ is obtained by shifting ${\mc A}$ in $F_2E_{\pentagon}2$ from one side to the other side. Therefore $F_3E_{\pentagon}2$ is a special case of $F'_2E_{\pentagon}2$, and we no longer use the notation $F_3E_{\pentagon}2$.

\begin{figure}[htp]
\centering
\begin{tikzpicture}[>=latex,scale=1]

\foreach \a in {1,-1}
{
\begin{scope}[scale=\a]

\draw
	(0.5,1) -- (0.5,0) -- (-0.5,0) -- (-1.5,1)
	(2.5,0) -- (3.5,0) -- (3.5,1) -- (-2.5,1) -- (-2.5,0) -- (-1.5,-1)
	(-1,0.5) -- (-2,-0.5);

\draw[line width=1.2]
	(-1,0.5) -- (-0.5,0)
	(-2,-0.5) -- (-2.5,0);

\foreach \b in {0,-1}
{
\begin{scope}[xshift=3*\b cm]

\node at (-0.3,-0.8) {\small $\alpha$};
\node at (-0.3,-0.2) {\small $\beta$};
\node at (1.05,-0.8) {\small $\gamma$};
\node at (0.4,-0.2) {\small $\delta$};
\node at (0.8,-0.55) {\small $\epsilon$};

\node at (1.6,-0.8) {\small $\alpha$};
\node at (1.3,-0.45) {\small $\beta$};
\node at (2.3,-0.8) {\small $\gamma$};
\node at (1.95,0.2) {\small $\delta$};
\node at (2.3,-0.1) {\small $\epsilon$};

\end{scope}
}

\end{scope}
}

\node[inner sep=0.5, draw, shape=circle] at (-0.1,0.5) {\small 1};
\node[inner sep=0.5, draw, shape=circle] at (0.1,-0.5) {\small 2};
\node[inner sep=0.5, draw, shape=circle] at (-1.1,-0.4) {\small 3};
\node[inner sep=0.5, draw, shape=circle] at (1.1,0.4) {\small 4};
\node[inner sep=0.5, draw, shape=circle] at (-1.9,0.4) {\small 5};
\node[inner sep=0.5, draw, shape=circle] at (1.9,-0.4) {\small 6};
\node[inner sep=0.5, draw, shape=circle] at (-2.9,-0.5) {\small 7};
\node[inner sep=0.5, draw, shape=circle] at (2.9,0.5) {\small 8};

\end{tikzpicture}
\caption{Interlocked ${\mc A}$.}
\label{emt2mod5}
\end{figure}
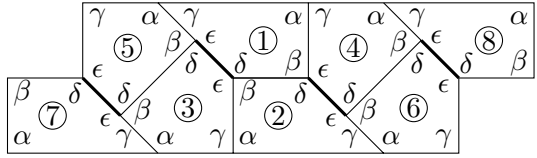

The second flip modification $F_2E_{\pentagon}2$ has another further modification. In the second of Figure \ref{emt2mod3}, we further find tiles $T_5,T_6$ and their companions, i.e., tiles $T_7,T_8$ sharing the $b$-edges with $T_5,T_6$. If we flip the patch ${\mc A}$ consisting of $T_1,T_2,T_3,T_4$ inside the eight tile patch $T_1,\dots,T_8$, then we get the patch in Figure \ref{emt2mod5}. The patch contains three interlocked ${\mc A}$:
\begin{align*}
{\mc A} &\colon T_1,T_2,T_3,T_4; \\
{\mc A}_+ &\colon T_1,T_3,T_5,T_7; \\
{\mc A}_- &\colon T_2,T_4,T_6,T_8.
\end{align*}

If we take away $T_2,T_3,T_5,T_7$ from the inner half of $F_2E_{\pentagon}2$, then we get ${\mc P}^q_{q-1}$. If we take away $T_1,T_4,T_6,T_8$ from the outer half of $F_2E_{\pentagon}2$, then we get ${\mc P}^q_{q-1}$. Therefore after flipping the patch ${\mc A}$ consisting of $T_1,T_2,T_3,T_4$, the tiling we get is the union of two half sphere tilings ${\mc P}^q_{q-1}\cup {\mc A}_-$ and ${\mc P}^q_{q-1}\cup {\mc A}_+$. The half sphere tilings are obtained by deleting the two tiles labelled $\times$ in the first of Figure \ref{emt2mod4}. Therefore we may shift and flip ${\mc A}_-$ and ${\mc A}_+$ independently in the respective half sphere tilings. We denote the tilings obtained this way by $F''_2E_{\pentagon}2$.

We conclude the following modifications of $E_{\pentagon}2$, and their characterisations:
\begin{enumerate}
\item $F_1E_{\pentagon}2$: Flip modification that does not contain ${\mc A}$.
\item $F_2E_{\pentagon}2$: Flip modification that contains ${\mc A}$. 
\item $F'_2E_{\pentagon}2$: Further modification by the shift and flip of ${\mc A}$, that contains a single ${\mc A}$. 
\item $F''_2E_{\pentagon}2$: Further modification by the shift and flip of ${\mc A}$, that contains two disjoint ${\mc A}$. 
\end{enumerate}
In $F'_2E_{\pentagon}2$, there are actually two ${\mc A}$ sharing two common tiles. See $T_1,T_3,T_5,T_6$ and $T_5,T_6,T_7,T_8$ in Figure \ref{4adeD2}. However, the tiling does not contain two disjoint ${\mc A}$.  

\subsubsection*{Combinatorially Possible Tilings}

In the classification proof, we obtain several tilings that are combinatorially possible, but fail the geometrical requirements that the pentagon must be simple \cite[Lemma 1]{gsy}, and the equality \eqref{coolsaet_eq1} must be satisfied. These ``tilings'' are not included in the context of this paper, but may be realisable in some other more general context.

The first example is the non-symmetric version of $E_{\pentagon}1$. Part of this combinatorially possible tiling is given by the first of Figure \ref{b2d_c2e_abcA}. The tiling even has a rich collection of flip modifications. 

The second example is the flip modification of $E_{\pentagon}2$ in Figure \ref{emt2mod7}. Combinatorially, the tiling can have any number of timezones. However, only four timezones is geometrically possible. 

The third example is the case $\alpha\delta^2,\gamma\epsilon^2,\alpha\beta\gamma$ are vertices, which is the last case of Proposition \ref{a2d_c2e}. Combinatorially, it is possible to have an earth map tiling with the second of Figure \ref{4a2eC} as the timezone. In fact, the tiling is combinatorially the same as the second example above, but with different angle combinations at vertices.

The fourth example is two tilings obtained by glueing two copies of any one of the half sphere tilings in Figure \ref{4a2eB}. The tilings have 28 tiles. 

The final example is that, in case $\alpha=\beta$ and $f=16,20$, there is an earth map tiling obtained by mixing the timezones in the first and third examples. See Figure \ref{4a2eC}.

\section{Technique}
\label{basic}

By a {\em tiling}, we mean an edge-to-edge tiling of the sphere by congruent almost equilateral pentagons. In particular, all vertices have degree $\ge 3$. Moreover, as explained in Section \ref{companion}, we will assume the number of tiles $f\ge 16$. Unless otherwise stated (mostly in Section \ref{geom}), the assumptions will always be implicit in the discussions, lemmas, and propositions. 

The pentagon is given by Figure \ref{pentagon}, with indicated edge lengths $a,b$ and angles $\alpha,\beta,\gamma,\delta,\epsilon$. The edges are straight, in the sense that they are parts of big circles. The edge lengths $a,b$ have distinct values. But some of the five angles may have the same values. 

All edge lengths and angle values are strictly between $0$ and $2\pi$, with Lemma \ref{geometry7} as the only exception. 
By \cite[Lemma 1]{gsy}, the pentagon is simple, in the sense the boundary is a simple closed curve. This implies $a<\pi$. 

Our problem has the exchange symmetry $(\beta,\delta)\leftrightarrow(\gamma,\epsilon)$. For example, if a result holds for vertices $\alpha\delta^2,\beta\delta\epsilon$, then similar result holds for vertices $\alpha\epsilon^2,\gamma\delta\epsilon$.

\subsection{Notations and Conventions}

\subsubsection*{Angles at a vertex}

We denote a vertex by $\alpha^a\beta^b\gamma^c\delta^d\epsilon^e$, which means that the vertex consists of $a$ copies of $\alpha$, $b$ copies of $\beta$, etc. The {\em angle sum of the vertex} is
\[
a\alpha+b\beta+c\gamma+d\delta+e\epsilon=2\pi. 
\]
Since angle sums are used so often, we will simply say ``by $\alpha\beta\gamma$'' to mean ``by the angle sum $\alpha+\beta+\gamma=2\pi$ of $\alpha\beta\gamma$''.

We use $\alpha\beta^2\cdots$ to mean a vertex with at least one $\alpha$ and at least two $\beta$. In other words, $\alpha\beta^2\cdots=\alpha^a\beta^b\gamma^c\delta^d\epsilon^e$, with $a\ge 1$ and $b\ge 2$. The angle combination in $\cdots$ is the {\em remainder}, which we denote by  $R(\alpha\beta^2)=\alpha^{a-1}\beta^{b-2}\gamma^c\delta^d\epsilon^e$. We also use $R$ to denote the total value of the remainder. Therefore $R(\alpha\beta^2)=2\pi-\alpha-2\beta=(a-1)\alpha+(b-2)\beta+c\gamma+d\delta+e\epsilon$.

Let $f$ be the number of tiles. By \cite[Lemma 4]{wy1}, the {\em angle sum for pentagon} is 
\begin{equation}\label{psum}
\alpha+\beta+\gamma+\delta+\epsilon=(3+\tfrac{4}{f})\pi.
\end{equation}
The equality means that the area of the pentagon is the area $4\pi$ of the sphere divided by the number $f$ of tiles. 

The {\em anglewise vertex combination}, abbreviated as AVC, is the collection of all possible vertices in a tiling. The following is an example from Proposition \ref{ade_2bc},
\begin{equation}\label{avc_example}
\text{AVC}=\{\alpha\delta\epsilon,\beta^2\gamma,\beta\gamma^k,\gamma^k\}.
\end{equation}
Strictly speaking, we should use $\beta\gamma^{k_1},\gamma^{k_2}$, because the numbers of $\gamma$ in these vertices are not the same. We always use $k,l$ as generic notations, that may take different values in different vertices. In fact, we also allow $k$ or $l$ to be $0$. For example, we have $\alpha^l\beta\gamma^k=\beta\gamma^k$ in case $l=0$. 

Some vertices in an AVC actually may not appear in the tiling. For example, if $\gamma^k$ is a vertex, then we derive the earth map tiling from the AVC, in which $\beta\gamma^k$ is not a vertex. On the other hand, if $\gamma^k$ is not a vertex, then $\beta\gamma^k$ is a vertex, and we derive the rotation modification of the earth map tiling.

\subsubsection*{$b$-Edge}

The most distinguished feature of almost equilateral pentagon is the $b$-edge. By considering the two bounding edges, the angles $\delta,\epsilon$ are {\em $b$-angles}, and $\alpha,\beta,\gamma$ are {\em $\hat{b}$-angles} (non-$b$-angles). Moreover, a {\em $b$-vertex} has $\delta,\epsilon$, and a {\em $\hat{b}$-vertex} has only $\alpha,\beta,
\gamma$ (and has no $\delta,\epsilon$). 

We indicate the {\em arrangement} of angles at a vertex by inserting the $a$-edge $\thin$ and the $b$-edge $\thick$ between the angles. For example, we may denote the vertex $\alpha\beta\delta\epsilon$ in Figure \ref{adjacent_deduce} by $\thick\delta\thin\beta\thin\alpha\thin\epsilon\thick$. The notation can be rotated or reversed. For example, the vertex can also be denoted as $\thin\alpha\thin\epsilon\thick\delta\thin\beta\thin$ (rotation) or $\thick\epsilon\thin\alpha\thin\beta\thin\delta\thick$ (reversion).

We use the same notation for the arrangement of consecutive angles at a vertex. For example, we have $\thick\delta\thin\beta\thin\alpha\thin$ as part of $\thick\delta\thin\beta\thin\alpha\thin\epsilon\thick$. We can also denote the part by its reversion $\thin\alpha\thin\beta\thin\delta\thick$ (but not rotation).

We will be very flexible in using the arrangement notation. For example, the vertex in Figure \ref{adjacent_deduce} is also $\alpha\thin\beta\cdots$ or $\delta\thick\epsilon\cdots$. 

At a $b$-vertex, the $b$-edges divide the vertex into several segments of consecutive $\hat{b}$-angles bounded by a pair of $b$-angles. Such a segment is a {\em fan}. Specifically, we have $\delta^2$-fan $\thick\delta\thin\cdots\thin\delta\thick$, $\epsilon^2$-fan $\thick\epsilon\thin\cdots\thin\epsilon\thick$, and $\delta\epsilon$-fan $\thick\delta\thin\cdots\thin\epsilon\thick$, where the interior $\cdots$ of the fans are filled by $\alpha,\beta,\gamma$. The {\em empty} fans are $\thick\delta\thin\delta\thick,\thick\epsilon\thin\epsilon\thick,\thick\delta\thin\epsilon\thick$.

\begin{figure}[htp]
\centering
\begin{tikzpicture}[>=latex,scale=1]


\foreach \a in {20,50,90,140,180,230}
\draw
	(0,0) -- (\a:1.1);

\foreach \b in {0,70,160,250}
\draw[line width=1.2]
	(0,0) -- (\b:1.1);

\foreach \x in {10,60,240}
\node at (\x:0.8) {\footnotesize $\delta$};

\foreach \y in {80,150,170}
\node at (\y:0.8) {\footnotesize $\epsilon$};

\draw[dotted]
	(25:0.8) arc (25:45:0.8)
	(95:0.8) arc (95:135:0.8)
	(185:0.8) arc (185:225:0.8);

\node[rotate=-55] at (35:1.25) {\small $\delta^2$-fan};
\node[rotate=25] at (115:1.25) {\small $\epsilon^2$-fan};
\node[rotate=115] at (205:1.25) {\small $\delta\epsilon$-fan};
		
\end{tikzpicture}
\caption{Fans at a $b$-vertex.}
\label{vertex_bedge}
\end{figure}
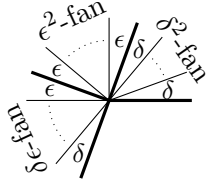

The number of fans at a vertex is the same as the number of $b$-edges. Since each fan has two $b$-angles, we get the following basic fact. 

\begin{lemma}[Parity Lemma]\label{beven}
The total number of $\delta$ and $\epsilon$ at any vertex is even. 
\end{lemma}

The lemma is the adaptation of \cite[Lemma 2]{cly} and \cite[Lemma 10]{wy2} to almost equilateral pentagon. Here are some typical applications of the parity lemma:
\begin{itemize}
\item A vertex $\alpha\delta\cdots$ is $\alpha\delta^2\cdots$ or $\alpha\delta\epsilon\cdots$. 
\item If $\alpha\delta\epsilon$ is a vertex, and $\delta<\epsilon$, then $\alpha\epsilon\cdots=\alpha\delta\epsilon$.
\item If the degree of $\alpha\delta\epsilon^k$ is $\ge 4$, then $k\ge 3$.
\item If $\delta>\epsilon>\frac{1}{2}\pi$, then there is at most one $b$-edge at any vertex, and $\delta\thin\delta\cdots,\delta\thin\epsilon\cdots,\epsilon\thin\epsilon\cdots$ are not vertices.
\item If $\delta+\epsilon>\pi$ and $\delta<\epsilon$, then $R(\epsilon^2)$ has no $\delta,\epsilon$, and $R(\delta\epsilon)$ has no $\epsilon$, and $\epsilon\thin\epsilon\cdots$ is not a vertex.
\end{itemize}
We use the parity lemma so often, that we will not mention the lemma in the proofs. Moreover, we will simply state the facts above without further explanation.

\subsubsection*{Adjacent angle deduction}

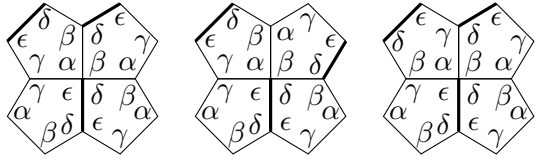
\begin{figure}[htp]
\centering
\begin{tikzpicture}[>=latex,scale=1]

\foreach \b in {0,1,2}
{
\begin{scope}[xshift=2.5*\b cm]

\foreach \a in {0,...,3}
\draw[rotate=90*\a]
	(0,0) -- (0.7,0) -- (1,0.5) -- (0.5,1) -- (0,0.7) -- (0,0);

\draw[line width=1.2]
	(0,0) -- (0,-0.7);

\node at (0.6,-0.25) {\small $\beta$}; 
\node at (0.2,-0.2) {\small $\delta$};
\node at (0.8,-0.45) {\small $\alpha$};
\node at (0.2,-0.6) {\small $\epsilon$};
\node at (0.5,-0.8) {\small $\gamma$}; 

\node at (-0.6,-0.2) {\small $\gamma$}; 
\node at (-0.2,-0.2) {\small $\epsilon$};
\node at (-0.8,-0.45) {\small $\alpha$};
\node at (-0.2,-0.6) {\small $\delta$};
\node at (-0.45,-0.75) {\small $\beta$}; 

\end{scope}
}


\draw[line width=1.2]
	(0,0.7) -- (0.5,1)
	(-1,0.5) -- (-0.5,1);

\node at (0.6,0.2) {\small $\alpha$}; 
\node at (0.2,0.2) {\small $\beta$};
\node at (0.8,0.45) {\small $\gamma$};
\node at (0.2,0.6) {\small $\delta$};
\node at (0.5,0.8) {\small $\epsilon$}; 

\node at (-0.6,0.2) {\small $\gamma$}; 
\node at (-0.2,0.2) {\small $\alpha$};
\node at (-0.8,0.5) {\small $\epsilon$};
\node at (-0.2,0.55) {\small $\beta$};
\node at (-0.5,0.8) {\small $\delta$};


\begin{scope}[xshift=2.5cm]

\draw[line width=1.2]
	(0.7,0) -- (1,0.5)
	(-1,0.5) -- (-0.5,1);

\node at (0.6,0.2) {\small $\delta$}; 
\node at (0.2,0.2) {\small $\beta$};
\node at (0.8,0.5) {\small $\epsilon$};
\node at (0.2,0.6) {\small $\alpha$};
\node at (0.5,0.75) {\small $\gamma$}; 

\node at (-0.6,0.2) {\small $\gamma$}; 
\node at (-0.2,0.2) {\small $\alpha$};
\node at (-0.8,0.5) {\small $\epsilon$};
\node at (-0.2,0.55) {\small $\beta$};
\node at (-0.5,0.8) {\small $\delta$};

\end{scope}


\begin{scope}[xshift=5cm]

\draw[line width=1.2]
	(0,0.7) -- (0.5,1)
	(-1,0.5) -- (-0.5,1);
	
\node at (0.6,0.2) {\small $\alpha$}; 
\node at (0.2,0.2) {\small $\beta$};
\node at (0.8,0.45) {\small $\gamma$};
\node at (0.2,0.6) {\small $\delta$};
\node at (0.5,0.8) {\small $\epsilon$}; 

\node at (-0.55,0.2) {\small $\beta$}; 
\node at (-0.2,0.2) {\small $\alpha$};
\node at (-0.8,0.5) {\small $\delta$};
\node at (-0.2,0.55) {\small $\gamma$};
\node at (-0.5,0.8) {\small $\epsilon$};

\end{scope}

\end{tikzpicture}
\caption{Adjacent angle deduction.}
\label{adjacent_deduce}
\end{figure}

The {\em adjacent angle deduction}, abbreviated as AAD, is the notation that keeps track of the arrangement of tiles around a vertex. For example, all three pictures in Figure \ref{adjacent_deduce} are AADs of the arrangement $\thick\delta\thin\beta\thin\alpha\thin\epsilon\thick$ of the vertex $\alpha\beta\delta\epsilon$. We denote them respectively by $\thick^{\epsilon}\delta^{\beta}\thin^{\alpha}\beta^{\delta}\thin^{\beta}\alpha^{\gamma}\thin^{\gamma}\epsilon^{\delta}\thick$, $\thick^{\epsilon}\delta^{\beta}\thin^{\delta}\beta^{\alpha}\thin^{\beta}\alpha^{\gamma}\thin^{\gamma}\epsilon^{\delta}\thick$, $\thick^{\epsilon}\delta^{\beta}\thin^{\alpha}\beta^{\delta}\thin^{\gamma}\alpha^{\beta}\thin^{\gamma}\epsilon^{\delta}\thick$.

Similar to the arrangement of angles at a vertex, the AAD can be rotated and reversed. For example, the first AAD in Figure \ref{adjacent_deduce} can also be denoted as $\thin^{\beta}\alpha^{\gamma}\thin^{\gamma}\epsilon^{\delta}\thick^{\epsilon}\delta^{\beta}\thin^{\delta}\beta^{\alpha}\thin$ (rotation) or $\thick^{\delta}\epsilon^{\gamma}\thin^{\gamma}\alpha^{\beta}\thin^{\delta}\beta^{\alpha}\thin^{\beta}\delta^{\epsilon}\thick$ (reversion). 

We also use AAD for consecutive angles at a vertex. For example, in the first picture, we have the AADs $\thick^{\epsilon}\delta^{\beta}\thin^{\alpha}\beta^{\delta}\thin^{\beta}\alpha^{\gamma}\thin$ and $\thin^{\beta}\alpha^{\gamma}\thin^{\gamma}\epsilon^{\delta}\thick$ at $\alpha\beta\delta\epsilon$. The AAD of consecutive angles can be reversed but not rotated.

We will be very flexible in using the AAD notation. For example, $\thick^{\epsilon}\delta^{\beta}\thin\beta\thin$ means $\thick^{\epsilon}\delta^{\beta}\thin^{\alpha}\beta^{\delta}\thin$ or $\thick^{\epsilon}\delta^{\beta}\thin^{\delta}\beta^{\alpha}\thin$. Moreover, ${}^{\beta}\thin^{\alpha}\beta^{\delta}\thin$ can be $\thin^{\gamma}\alpha^{\beta}\thin^{\alpha}\beta^{\delta}\thin$ or $\thick^{\epsilon}\delta^{\beta}\thin^{\alpha}\beta^{\delta}\thin$. We may even use ${}^{\beta}\thin\beta\thin$ to denote ${}^{\beta}\thin^{\alpha}\beta^{\delta}\thin$ or ${}^{\beta}\thin^{\delta}\beta^{\alpha}\thin$. 

The AAD is a convenient substitution for drawing some pictures. When we use AAD to make an argument, we actually have a picture in mind. Without AAD, this paper would have many more pictures. 

An AAD implies new AAD or new vertex. For example, the right half of the first of Figure \ref{adjacent_deduce} shows the AAD $\thick^{\epsilon}\delta^{\beta}\thin^{\alpha}\beta^{\delta}\thin$ (at a vertex $\beta\delta\cdots$) implies an AAD $\thin^{\gamma}\alpha^{\beta}\thin^{\delta}\beta^{\alpha}\thin$ (at a vertex $\alpha\beta\cdots$). In general, we have the following {\em reciprocity property}: An AAD $\lambda^{\theta}\thin^{\rho}\mu$ implies an AAD $\theta^{\lambda}\thin^{\mu}\rho$.

Here is a typical AAD argument using the reciprocity property. The AADs of $\thin\beta\thin\beta\thin$ are $\thin^{\alpha}\beta^{\delta}\thin^{\alpha}\beta^{\delta}\thin,\thin^{\alpha}\beta^{\delta}\thin^{\delta}\beta^{\alpha}\thin,\thin^{\delta}\beta^{\alpha}\thin^{\alpha}\beta^{\delta}\thin$. If we know $\alpha\delta\cdots,\delta\thin\delta\cdots$ are not vertices, then we know the AAD of $\thin\beta\thin\beta\thin$ is $\thin^{\delta}\beta^{\alpha}\thin^{\alpha}\beta^{\delta}\thin$.

Another example of the AAD argument is that a vertex $\alpha^3$ (or more generally, $\alpha^k$ with odd $k$) implies a vertex $\beta\gamma\cdots$. See \cite[Lemma 10]{wy1}. 

\subsection{Special Tile and Companion Pair}
\label{companion}

\begin{figure}[htp]
\centering
\begin{tikzpicture}[>=latex,scale=1]

\foreach \a in {-1,1}
\foreach \b in {0,1}
{
\begin{scope}[xshift=4*\b cm - 0.647*\a cm,xscale=\a]

\foreach \c in {0,...,4}
\draw[rotate=72*\c]
	(36:0.8) -- (-36:0.8);

\draw[line width=1.2]
	(36:0.8) -- (-36:0.8);

\end{scope}
}

\foreach \b in {0,1}
{
\begin{scope}[xshift=4*\b cm-0.647 cm]

\node at (180:0.5) {$\alpha$};
\node at (-108:0.5) {$\beta$};
\node at (108:0.5) {$\gamma$};
\node at (-36:0.5) {$\delta$};
\node at (36:0.5) {$\epsilon$};

\end{scope}
}

\foreach \b in {-1,1}
{
\begin{scope}[xshift=2*\b cm+2.647 cm,yscale=-\b]

\node at (0:0.5) {$\alpha$};
\node at (72:0.5) {$\beta$};
\node at (144:0.5) {$\delta$};
\node at (-144:0.5) {$\epsilon$};
\node at (-72:0.5) {$\gamma$};

\end{scope}
}

\node at (0,-1) {twisted};	
\node at (4,-1) {matched};	

\end{tikzpicture}
\caption{Twisted and matched companion pairs.}
\label{twist_match}
\end{figure}
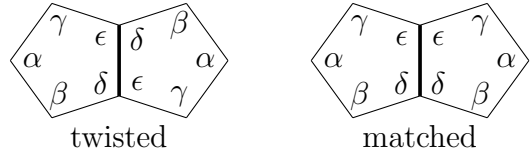

In an edge-to-edge tiling by congruent almost equilateral pentagons, any tile has a {\em companion tile} sharing the common $b$-edge. Figure \ref{twist_match} shows that the two tiles can be {\em twisted} or {\em matched}. A {\em companion pair} has two shared vertices (the two ends of the common $b$-edge), and six non-shared vertices. The tiling can be regarded as a tiling of companion pairs.

Let $v,e,f$ be the numbers of vertices, edges, and tiles. Let $v_k$ be the number of vertices of degree $k$. Then the Euler formula and the pentagonal tiling imply (see \cite[Section 2]{wy1}, for example)
\begin{align}
f
&=12+2\ssum_{k\ge 4}(k-3)v_k=12+2v_4+4v_5+6v_6+\cdots, \label{fnumber} \\
v_3
&=20+\ssum_{k\ge 4}(3k-10)v_k=20+2v_4+5v_5+8v_6+\cdots. \label{v3}
\end{align}
The equality \eqref{fnumber} implies that $f$ is an even integer $\ge 12$. Gao, Shi and Yan \cite{ay,gsy} classified edge-to-edge tilings of the sphere by $12$ congruent pentagons. Moreover, Yan \cite[Theorem 1]{yan1} showed that $f\ne 14$. Like \cite{awy,wy1,wy2}, therefore, we will assume $f\ge 16$ throughout the argument in this paper. By \eqref{fnumber}, this implies not all vertices have degree $3$.

The equality \eqref{v3} implies that most vertices have degree $3$. A degree $3$ $b$-vertex is $\theta\delta\epsilon, \theta\delta^2,\theta\epsilon^2$, where $\theta=\alpha,\beta,\gamma$. In view of Figure \ref{twist_match}, we call $\theta\delta\epsilon$ {\em twisted} and call $\theta\delta^2,\theta\epsilon^2$ {\em matched}. We denote by $\text{AVC}_3$ all the degree $3$ vertices in a tiling. We also call vertices of degree $>3$ {\em high degree vertices}. 

For each $k>3$, we assign a weight $w_k>0$ to degree $k$ vertices. Then we define the weight of a tile $T$ to be the sum of weights of all high degree vertices of $T$ 
\[
w(T)
=\ssum\{w_{\deg v}\colon \text{$v$ is a high degree vertex of $T$}\}.
\]
We also define the weight of a companion pair $P$ to be the sum of the weights of its two tiles (the weight of shared high degree vertex is counted twice). Then the total sum of the weights is
\[
\ssum_Tw(T)=\ssum_Pw(P)=\ssum_{k\ge 4}kw_kv_k.
\]

Now we fix the values of the weights
\begin{equation}\label{weight}
kw_k=2(k-3),\quad
w_k=2-\tfrac{6}{k}
=\tfrac{1}{2},\;\tfrac{4}{5},\;1,\;\tfrac{8}{7},\;\tfrac{5}{4},\;\dots.
\end{equation}
By \eqref{fnumber}, we get
\begin{equation}\label{weightsum}
f=12+\sum_{k\ge 4}kw_kv_k=12+\sum_Pw(P)>\sum_Pw(P).
\end{equation}
Since the number of terms in the sum on the right is the number $\frac{f}{2}$ of companion pairs, the inequality implies the following. 

\begin{lemma}\label{bb_pair}
Any tiling has a companion pair $P$ satisfying $w(P)<2$. 
\end{lemma}

We call a companion pair $P$ {\em special} if $w(P)<2$. 

\begin{lemma}\label{special_pair}
If all companion pairs have weight $\ge 2w_4=1$, then $f\ge 24$. If all companion pairs have weight $\ge w_4+w_5=\frac{13}{10}$, then $f\ge 36$. If all companion pairs have weight $\ge 3w_4=\frac{3}{2}$, then $f\ge 48$. If all companion pairs have weight $\ge 2w_5=\frac{8}{5}$, then $f\ge 60$.  
\end{lemma}

\begin{proof}
The sum $\sum_Pw(P)$ in \eqref{weightsum} has $\frac{f}{2}$ terms. If all $w(P)\ge 1$, then we get $f\ge 12+\frac{f}{2}$. This means $f\ge 24$. The proof for the other claims are similar.   
\end{proof}

We also have $f=12+\sum_Tw(T)>\sum_Tw(T)$. The argument above for companion pairs can be repeated for tiles. Then we find a tiling has a tile $T$ satisfying $w(T)<1$, which we call a {\em special tile}. This means four vertices of $T$ have degree $3$, and the fifth vertex has degree $3,4,5$. We call the tile respectively $3^5$-, $3^44$-, $3^45$-tile. Among the two tiles in a special companion pair, at least one is a special tile. 

The following is Lemmas 1, 2, 3 of \cite{wy1}. Also see \cite[Lemma 4]{wy2}. 

\begin{lemma}\label{special_tile}
Any tiling has a special tile. Moreover, 
\begin{enumerate}
\item If there is no $3^5$-tile, then $f\ge 24$. And $f=24$ implies that each tile is a $3^44$-tile, and the tiling is the pentagonal subdivision of the octahedron.
\item If there is no $3^5$-tile and no $3^44$-tile, then $f\ge 60$. And $f=60$ implies that each tile is a $3^45$-tile, and the tiling is the pentagonal subdivision of the icosahedron.
\end{enumerate}
\end{lemma}

\subsection{Counting}
\label{count}

We use $\#$ to denote various total numbers. For example, $\#\alpha$ is the total number of $\alpha$ in the tiling. If $\alpha$ appears once in the pentagon, then we have $\#\alpha=f$. For the AVC \eqref{avc_example}, we have
\begin{align*}
f=\#\alpha=\#\delta=\#\epsilon
&=\#\alpha\delta\epsilon, \\
f=\#\beta
&=2\#\beta^2\gamma
+\#\beta\gamma^k, \\
f=\#\gamma
&=\#\beta^2\gamma
+k\#\beta\gamma^k
+k\#\gamma^k.
\end{align*}
Here $\#\alpha\delta\epsilon$ is the total number of the vertex $\alpha\delta\epsilon$, etc. 

A subtle point about counting is the meaning of {\em distinct} angles, i.e., the criterion used for distinguishing angles. We may certainly use angle values as the criterion. When some angle values are the same, we may further combine the $b$-edge to distinguish angles. By Lemma \ref{geometry11} (also see \cite[Lemma 21]{gsy}), we know $\beta=\gamma$ if and only if $\delta=\epsilon$. Therefore we will divide the classification into three parts:
\begin{enumerate}
\item Symmetric: $\beta=\gamma$ and $\delta=\epsilon$. In this case, the pentagon has three angle values $\alpha,\beta,\delta$, and we may have $\alpha=\beta$. Moreover, $\delta$ is distinguished from $\alpha,\beta$ by being a $b$-angle.
\item $\alpha,\beta,\gamma$ have distinct values. In this case, we also have $\delta\ne\epsilon$. Then we may use angle values and $b$-edge to distinguish all five angles. 
\item $\alpha=\beta\ne \gamma$. In this case, we also have $\delta\ne\epsilon$. Then we may use angle values and $b$-edge to distinguish $\alpha,\gamma,\delta,\epsilon$. However, we need more information to distinguish $\alpha,\beta$. 
\end{enumerate}
We remark that, due to the exchange symmetry $(\beta,\delta)\leftrightarrow(\gamma,\epsilon)$, the case $\alpha=\gamma\ne \beta$ is included in the case $\alpha=\beta\ne \gamma$.

The following is \cite[Lemma 4]{cly}. We give a special name because it is frequently used. 

\begin{lemma}[Counting Lemma]\label{counting}
Suppose $\theta,\rho$ appear the same number of times in the pentagon. If all vertices are of the form $\theta^k\rho^l\cdots$, with $k\le l$ and no $\theta,\rho$ in the remainder, then all vertices are of the form $\theta^k\rho^k\cdots$, with no $\theta,\rho$ in the remainder. 
\end{lemma}

\begin{proof}
Suppose all vertices are $\theta^{k_i}\rho^{l_i}\cdots$, with no $\theta,\rho$ in the remainders. Then by $k_i\le l_i$ for all $i$, we get
\[
\#\theta=\sum_ik_i\#(\theta^{k_i}\rho^{l_i}\cdots)
\le \sum_il_i\#(\theta^{k_i}\rho^{l_i}\cdots)=\#\rho.
\]
Since $\theta,\rho$ appear the same number of times in the pentagon, we get $\#\theta=\#\rho$. Then the inequality above implies $k_i=l_i$ or $\#(\theta^{k_i}\rho^{l_i}\cdots)=0$, for each $i$. Of course $\#(\theta^{k_i}\rho^{l_i}\cdots)=0$ means that $\theta^{k_i}\rho^{l_i}\cdots$ is actually not a vertex. 
\end{proof}

Here is a typical application of the counting lemma. Suppose we know $\alpha\delta\epsilon$ is the only $b$-vertex. Then in a vertex $\alpha^l\delta^k\cdots$ (no $\alpha,\delta$ in the remainder), we have $k=l=1$ at $\alpha\delta\epsilon$, and $k=0\le l$ if the vertex is not $\alpha\delta\epsilon$. Then by applying the counting lemma to $\theta=\delta$ and $\rho=\alpha$, we get $k=l=0$ at any vertex other than $\alpha\delta\epsilon$. In other words, if a vertex is not $\alpha\delta\epsilon$, then it has no $\alpha,\delta$. By the same reason, the vertex also has no $\alpha,\epsilon$. Therefore the only vertex besides $\alpha\delta\epsilon$ is $\beta^k\gamma^l$.

The following is the adaptation of \cite[Lemma 11]{wy2} to almost equilateral pentagon. Again the lemma is used so often, that we give the same special name as in the earlier paper. 

\begin{lemma}[Balance Lemma]\label{balance}
$\delta^2\cdots$ is a vertex if and only if $\epsilon^2\cdots$ is a vertex. If both are not vertices, then a $b$-vertex is $\delta\epsilon\cdots$, with no $\delta,\epsilon$ in the remainder. 
\end{lemma}

\begin{proof}
If $\delta^2\cdots$ is not a vertex, then $\delta\cdots=\delta\epsilon^l\cdots$, with no $\delta,\epsilon$ in the remainder. By the parity lemma (Lemma \ref{beven}), we know $1+l$ is even. Therefore $1\le l$. Then by applying the counting lemma (Lemma \ref{counting}) to $\delta,\epsilon$, we get $l=1$. In particular, we know $\epsilon^2\cdots$ is not a vertex.
\end{proof}

The following is a useful application of the balance lemma.

\begin{lemma}\label{square}
If $\alpha\delta\epsilon$ is a vertex, then $\alpha^2\cdots$ is a $\hat{b}$-vertex. More generally, if $\lambda$ is a combination of $\alpha,\beta,\gamma$ with value $\ge 2\alpha$, then $\lambda\cdots$ is a $\hat{b}$-vertex.
\end{lemma}

The lemma is still true with $\beta$ or $\gamma$ in place of $\alpha$. 

\begin{proof}
If $\alpha^2\delta^2\cdots$ is a vertex, then the angle sum of $\alpha^2\delta^2\cdots$ implies $\alpha+\delta\le \pi$. By (the angle sum of) $\alpha\delta\epsilon$, this implies $\epsilon\ge \pi$. Therefore $\epsilon^2\cdots$ is not a vertex. By the balance lemma (Lemma \ref{balance}), however, this contradicts the assumption that $\alpha^2\delta^2\cdots$ is a vertex. 

Therefore $\alpha^2\delta^2\cdots$ is not a vertex. By the same reason, $\alpha^2\epsilon^2\cdots$ is not a vertex. By $\alpha\delta\epsilon$, the angle sum implies $\alpha^2\delta\epsilon\cdots$ is not a vertex. Then by the parity lemma, we conclude $\alpha^2\cdots$ has no $\delta,\epsilon$. 

The argument is still valid if we replace $\alpha^2$ by $\lambda$.
\end{proof}

The following is the more refined version of the balance lemma.

\begin{lemma}\label{fbalance}
There is no $\delta^2$-fan if and only if there is no $\epsilon^2$-fan. In this case, the only fans are $\delta\epsilon$-fans, and every vertex is of the form $\delta^k\epsilon^k\cdots$, with no $\delta,\epsilon$ in the remainder.
\end{lemma}

\begin{proof}
Suppose a vertex has no $\delta^2$-fan, and $k$ is the number of $\delta\epsilon$-fans at the vertex, and $l$ is the number of $\epsilon^2$-fans at the vertex. Then the vertex is $\delta^k\epsilon^{k+2l}\cdots$, with no $\delta,\epsilon$ in the remainder. If all vertices are like this, then we may apply the counting lemma to $\delta,\epsilon$, with $k\le  k+2l$. We conclude $k=k+2l$, i.e., $l=0$ at all vertices.
\end{proof}

Since most vertices have degree $3$, an angle appearing at every degree $3$ vertex means the abundance of the angle in the whole tiling. This leads to the subsequent statistical results.

The following is \cite[Lemmas 6, 7]{wy1}.

\begin{lemma}\label{degree3}
If an angle appears at all degree $3$ vertices, then it appears at least twice in the pentagon. If an angle (or two angles together) appears at least twice at all degree $3$ vertices, then the angle (or the two together) appears at least three times in the pentagon. 
\end{lemma}

The following is \cite[Lemma 8]{wy1}. 

\begin{lemma}\label{ndegree3}
There is at most one angle not appearing at all degree $3$ vertices. Moreover, if $\theta$ is such an angle, then one of $\rho\theta^3$, $\theta^4$, $\theta^5$ is a vertex, where $\rho\ne\theta$.
\end{lemma}

\subsection{Geometry}
\label{geom}

In this section, by almost equilateral pentagon, we mean the pentagon in Figure \ref{pentagon}. The pentagon is {\em simple} if the boundary does not intersect itself. By \cite[Lemma 1]{gsy}, the simple property is required for tiling. 

The following is the adaptation of \cite[Lemma 2]{wy2} to almost equilateral pentagon. 

\begin{lemma}\label{geometry1}
If an almost equilateral pentagon is simple, then $\beta>\gamma$ if and only if $\delta<\epsilon$.
\end{lemma}

The following is \cite[Lemma 21]{gsy}, which can be regarded as the equality version of the lemma above.

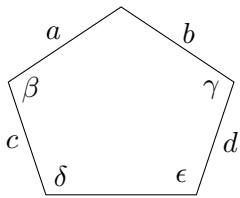
\begin{figure}[htp]
\centering
\begin{tikzpicture}

\draw 
	(-1,-1.5) -- (-1.5,0) -- (0,1) -- (1.5,0) -- (1,-1.5) -- cycle
	;
    
\node at (-1.2,-0.1) {\small $\beta$};
\node at (1.2,-0.1) {\small $\gamma$};
\node at (-0.8,-1.25) {\small $\delta$};
\node at (0.8,-1.25) {\small $\epsilon$};

\node at (-0.9,0.65) {\small $a$};
\node at (0.9,0.65) {\small $b$};
\node at (-1.45,-0.8) {\small $c$};
\node at (1.45,-0.8) {\small $d$};

\end{tikzpicture}
\caption{Four equalities.}
\label{pentcombo}
\end{figure}

\begin{lemma}\label{geometry11}
In the pentagon in Figure \ref{pentcombo}, three of the following equalities imply the fourth
\[
a=b,\quad
c=d,\quad
\beta=\gamma,\quad
\delta=\epsilon.
\]
\end{lemma}

For almost equilateral pentagon, the lemma implies $\beta=\gamma$ if and only if $\delta=\epsilon$. This means the pentagon is symmetric. Therefore we call the pentagon {\em non-symmetric} if $\beta\ne\gamma$ and $\delta\ne\epsilon$. The consideration leads to the discussion about the three parts in the classification argument in Section \ref{count}.

Another way of using the lemma is that, if $\alpha=\delta$ and $\gamma=\epsilon$ (or $\alpha=\epsilon$ and $\beta=\delta$) in an almost equilateral pentagon, then $a=b$. This is considered a contradiction in this paper.

The following is the adaptation of \cite[Lemma 5]{wy2} to almost equilateral pentagon.

\begin{lemma}\label{geometry2}
If an almost equilateral pentagon is simple and non-symmetric, then the following are the only possible combinations of two degree $3$ $b$-vertices without $\alpha$.
\begin{enumerate}
\item $\beta\delta\epsilon,\gamma\epsilon^2$.
\item $\gamma\delta\epsilon,\beta\delta^2$.
\item $\beta\delta^2,\gamma\epsilon^2$.
\end{enumerate}
\end{lemma}

\begin{proof}
By $\beta\ne\gamma$ and $\delta\ne\epsilon$, if the pair has $\beta\delta\epsilon$, then the other vertex in the pair is $\gamma\delta^2,\gamma\epsilon^2$. The angle sums of $\beta\delta\epsilon,\gamma\delta^2$ imply $\beta+\epsilon=\gamma+\delta$. Then we know $\beta>\gamma$ implies $\delta>\epsilon$, contradicting Lemma \ref{geometry1}. Therefore the other vertex can only be $\gamma\epsilon^2$.

By the similar argument, if the pair has $\gamma\delta\epsilon$, then the other vertex in the pair is $\beta\delta^2$. 

It remains to consider the pair $\beta\delta^2,\gamma\epsilon^2$ and the pair $\beta\epsilon^2,\gamma\delta^2$. The angle sums of $\beta\epsilon^2,\gamma\delta^2$ imply $\beta+2\epsilon=2\pi=\gamma+2\delta$. Then we know $\beta>\gamma$ implies $\delta>\epsilon$, contradicting Lemma \ref{geometry1}. 
\end{proof}

\begin{lemma}\label{geometry5}
Suppose $AP=AB=BQ=a<\pi$ in the first of Figure \ref{geom_proof2}, and $\alpha=\angle PAB$ and $\beta=\angle QBA$. If the three edges do not cross each other, then $\alpha+2\beta>\pi$ or $2\alpha+\beta>\pi$.
\end{lemma}

\begin{proof}
We have two isosceles triangles $\triangle PAB$ and $\triangle QBA$. If $\alpha+2\beta\le \pi$, then the isosceles triangle $\triangle PAB$ implies $\angle PBA>\frac{1}{2}(\pi-\alpha)\ge \beta$. Therefore $BQ$ lies inside $\angle PBA$. By the same reason, if $2\alpha+\beta\le \pi$, then $AP$ lies inside $\angle QAB$. If both $\alpha+2\beta\le \pi$ and $2\alpha+\beta\le \pi$ hold, then the two inside positions imply that $AP$ and $BQ$ must intersect.
\end{proof}

\begin{figure}[htp]
\centering
\begin{tikzpicture}[>=latex]

\begin{scope}[xshift=-3.8cm]

\draw
	(0.3,1.8) -- (-1.2,0) -- (1.2,0) -- (-0.6,1.5);

\draw[dashed]
	(0.3,1.8) -- (1.2,0)
	(-1.2,0) -- (-0.6,1.5);

\node at (-0.8,0.2) {\small $\alpha$};	
\node at (0.6,0.2) {\small $\beta$};

\node at (-1.4,0) {\small $A$};
\node at (1.4,0) {\small $B$};
\node at (0.3,2) {\small $P$};
\node at (-0.6,1.7) {\small $Q$};

\node at (0,-0.15) {\small $a$};
\node at (-0.55,1) {\small $a$};
\node at (0.4,0.9) {\small $a$};

\end{scope}

\draw
	(-1.5,0) -- (0,2) -- (1.5,0)
	(-1.5,0) -- ++(1,0.2) -- ++(0.6,0.4) -- ++(0.1,0.6);
	
\draw[dashed]
	(-1.5,0) -- (1.5,0)
	(-0.5,0.2) -- ++(1.6,0.32)
	(0.1,0.6) -- ++(0.63,0.42);

\node at (0.6,1.5) {\small $a$};
\node at (-0.6,1.5) {\small $a$};

\node at (1.7,0) {\small $A$};
\node at (-1.7,0) {\small $C$};
\node at (0,2.2) {\small $B$};

\node at (-0.05,1.65) {\small $\beta$};
\node at (-1.1,0.25) {\small $\gamma$};

\node at (-0.4,0) {\small $X_1$};
\node at (0.3,0.45) {\small $X_2$};

\node at (-0.6,0.4) {\small $\xi_1$};
\node at (-0.1,0.75) {\small $\xi_2$};

\node at (1.35,0.6) {\small $Y_1$};
\node at (1,1.1) {\small $Y_2$};

\begin{scope}[xshift=3.8cm]

\draw
	(-1.3,0.5) -- (-0.8,2) -- (0.8,2) -- (1.3,0.5);

\draw[densely dotted]
	(-1.3,0.5) -- (0,0) -- (1.3,0.5);

\draw[dashed]
	(-1.3,0.5) -- (1.3,0.5)
	(-1.3,0.5) -- (0.8,2)
	(1.3,0.5) -- (-0.8,2);

\node at (0,2.15) {\small $a$};
\node at (1.2,1.3) {\small $a$};
\node at (-1.2,1.3) {\small $a$};
\node at (0,0.65) {\small $x$};

\node at (-0.9,2.2) {\small $A$};
\node at (0.9,2.2) {\small $B$};

\node at (-1.5,0.5) {\small $C$};
\node at (1.5,0.5) {\small $D$};

\node at (0,-0.2) {\small $E$};

\node at (-0.3,1.85) {\small $\alpha_2$};
\node at (-0.65,1.65) {\small $\alpha_1$};
\node at (0.6,1.8) {\small $\beta$};
\node at (-1,0.75) {\small $\gamma'$};

\end{scope}

\begin{scope}[shift={(7cm,1cm)}]

\foreach \a in {0,...,4}
\draw[rotate=72*\a]
	(18:1.2) -- (90:1.2);

\draw[line width=1.2]
	(234:1.2) -- (-54:1.2);

\draw[dashed]
	(234:1.2) -- (18:1.2);
	
\node at (90:0.95) {\small $\alpha$};
\node at (162:0.95) {\small $\beta$};
\node at (-54:0.95) {\small $\epsilon$};

\node at (225:0.8) {\small $\delta_1$};	
\node at (265:0.75) {\small $\delta_2$};

\node at (26:0.85) {\small $\gamma_1$};	
\node at (-12:0.75) {\small $\gamma_2$};	

\foreach \a in {0,...,3}
\node at (-18+72*\a:1.15) {\small $a$};

\node at (-90:1.2) {\small $b$};

\end{scope}

\end{tikzpicture}
\caption{Lemmas \ref{geometry5}, \ref{geometry6}, \ref{geometry10}, \ref{geometry4}.}
\label{geom_proof2}
\end{figure}

\begin{lemma}\label{geometry6}
Suppose $A,B,C,X_1,X_2,\dots$ are successive vertices of a simple polygon, with corresponding angles $\alpha,\beta,\gamma,\xi_1,\xi_2,\dots$. If $AB=BC=a<\pi$ and $\xi_i\le\pi$, then $\beta+2\alpha>\pi$ and $\beta+2\gamma>\pi$.
\end{lemma}

Applying the lemma to a simple almost equilateral pentagon, we get the following:
\begin{enumerate}
\item If $\delta,\epsilon<\pi$, then $\alpha+2\beta>\pi$ and $\alpha+2\gamma>\pi$.
\item If $\beta,\delta<\pi$, then $2\alpha+\gamma>\pi$ and $\gamma+2\epsilon>\pi$.
\item If $\gamma,\epsilon<\pi$, then $2\alpha+\beta>\pi$ and $\beta+2\delta>\pi$.
\end{enumerate}

\begin{proof}
The conclusion holds for $\beta\ge\pi$. Therefore we may assume $\beta<\pi$. Then the angle $\beta$ and two edges $AB,BC$ form a simple and convex isosceles triangle $\triangle ABC$. We prove $\gamma\ge \angle ACB$, the base angle of the isosceles triangle. Then $\beta+2\gamma$ is no less than the sum of three angles in $\triangle ABC$, which we know is $>\pi$.

Suppose $\gamma<\angle ACB$. Then in the second of Figure \ref{geom_proof2}, $CX_1$ lies inside $\triangle ABC$, and the extension of $CX_1$ intersects $BA$ at $Y_1$. Then $\triangle BCY_1$ is convex and contained in $\triangle ABC$. By $\xi_1\le\pi$, we know $X_1X_2$ lies inside $\triangle BCY_1$, and the extension of $X_1X_2$ intersects $BY_1$ at $Y_2$. Then $\square BCX_1Y_2$ is convex and contained in $\triangle BCY_1$. Keep going, we get smaller and smaller convex polygons, and $A$ lies in the smallest polygon. However, we know $A$ actually lies outside $\triangle BCY_1$, a contradiction. 
\end{proof}

\begin{lemma}\label{geometry10}
Suppose the pentagon in the third of Figure \ref{geom_proof2} is convex and satisfies $AB=AC=BD=a<\pi$. If the angles at $A,B,C$ are $\alpha,\beta,\gamma$, and $\beta\ge\gamma$, then $2\alpha+\beta+\gamma>2\pi$.
\end{lemma}

If $CE=a$, then the pentagon is almost equilateral. In this case, we always have $2\alpha+\beta+\gamma>2\pi$, whether $\beta\ge\gamma$ or not.

\begin{proof}
By the convexity, we have $\gamma'=\angle ACD\le \angle ACE=\gamma<\beta=\angle ABD$. Then by the isosceles triangle $\triangle ABC$, we get $\angle CBD\ge\angle BCD$. This implies $CD \ge BD=AC$, and further implies $\alpha_1=\angle CAD\ge\angle ADC$. On the other hand, the isosceles triangle $\triangle ABD$ implies $2\alpha_2+\beta>\pi$, where $\alpha_2=\angle BAD$. Moreover, the triangle $\triangle ACD$ implies $2\alpha_1+\gamma\ge \alpha_1+\angle ADC+\gamma'>\pi$. Adding the two inequalities together, we get $2\alpha+\beta+\gamma>2\pi$.
\end{proof}

\begin{lemma}\label{geometry4}
If an almost equilateral pentagon is strictly convex, then any one of the following implies $a=b$:
\begin{enumerate}
\item $\beta\ge\alpha\ge\gamma$, $\beta\ge\epsilon$, $\delta\le\gamma$.
\item $\beta\le\alpha\le\gamma$, $\beta\le\epsilon$, $\delta\ge\gamma$.
\end{enumerate}
\end{lemma}

Since we assume $a\ne b$ in this paper, the lemma means we cannot have either set of inequalities happening at the same time.

\begin{proof}
Suppose $\beta\ge\alpha$ and $\delta\ge\gamma$. In the fourth of Figure \ref{geom_proof2}, we connect $\gamma,\delta$, and divide the angles $\gamma,\delta$ into $\gamma_1,\gamma_2,\delta_1,\delta_2$. By the strict convexity, we know the quadrilateral and the triangle are proper subsets of the pentagon. Therefore $\gamma_1,\gamma_2,\delta_1,\delta_2$ are positive angles, and $\gamma=\gamma_1+\gamma_2$, and $\delta=\delta_1+\delta_2$. By \cite[Lemma 2]{wy2}, we know $\beta\ge\alpha$ implies $\delta_1\le \gamma_1$. Then by $\delta\ge \gamma$, we get $\delta_2=\delta-\delta_1\ge \gamma_2=\gamma-\gamma_1$. By $\epsilon<\pi$, this further implies $a\ge b$.

Suppose $\gamma\le \alpha$ and $\epsilon\le \beta$. Then we may exchange $\beta,\delta,\ge$ with $\gamma,\epsilon,\le$, and carry out the same argument. The conclusion is $a\le b$. This completes the proof that the first set of inequalities imply $a=b$.

If we reverse the direction of all the inequalities, then the same proof shows that the second set of inequalities also imply $a=b$.
\end{proof}

\subsection{Existence of Pentagon}
\label{exist}

The focus of this section is Lemma \ref{geometry7}. This is the pentagonal analogue of \cite[Lemma 18]{cly} and \cite[Theorem 2.1]{coolsaet}.

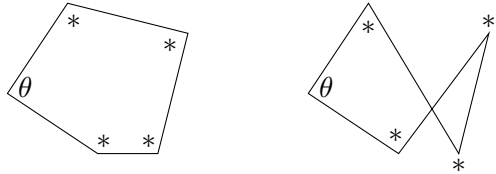
\begin{figure}[htp]
\centering
\begin{tikzpicture}[scale=0.8]

\draw
	(0,0) -- (1,1.5) -- (3,1) -- (2.5,-1) -- (1.5,-1) -- cycle;

\node at (0.3,0.1) {\small $\theta$};
\node at (1.1,1.2) {\small $*$};
\node at (2.7,0.8) {\small $*$};
\node at (2.35,-0.8) {\small $*$};
\node at (1.6,-0.8) {\small $*$};

\begin{scope}[xshift=5cm]

\draw
	(0,0) -- (1,1.5) -- (2.5,-1) -- (3,1) -- (1.5,-1) -- cycle;

\node at (0.3,0.1) {\small $\theta$};
\node at (1,1.1) {\small $*$};
\node at (3,1.2) {\small $*$};
\node at (2.5,-1.2) {\small $*$};
\node at (1.45,-0.7) {\small $*$};
	
\end{scope}

\end{tikzpicture}
\caption{Polygon determined by closed path and one angle.}
\label{polygon}
\end{figure}

The lemma is actually valid for more general kind of almost equilateral pentagon. We define (for this section only) a {\em polygon} to be a sequence of arcs (called {\em edges}) connected together (at {\em vertices}) to form a closed path, together with the choice of a {\em side}. The side means picking one angle $\theta$ at one vertex determines the choice of all the angles $*$ at the other vertices, by the criterion of being ``on the same side'' of $\theta$ along the path. See Figure \ref{polygon}. Picking the complement $2\pi-\theta$ of $\theta$ gives another side, and therefore another polygon. 

We remark that Lemma \ref{geometry11} can be extended to general almost equilateral pentagon. Specifically, if $D\ne E$ and $\beta=\gamma$, then $\delta=\epsilon$. Similarly, if $B\ne C$ and $\delta=\epsilon$, then $\beta=\gamma$. For the proof, see \cite[Lemma 21]{gsy}.

A polygon is {\em simple} if the closed path is simple, i.e., not crossing itself. The first of Figure \ref{polygon} is simple, and the second is not simple. The following holds for general almost equilateral pentagon defined above, and the pentagon is not required to be simple. However, we do assume edge lengths $a,b>0$.

\begin{lemma}\label{geometry7}
If the angles $\alpha,\beta,\gamma,\delta,\epsilon$ and the edge length $a$ of an almost equilateral pentagon satisfy $a\not\in {\bb Z}\pi$, then
\begin{align}
&
[
((1-\cos\beta)\sin(\delta-\tfrac{1}{2}\alpha)
-(1-\cos\gamma)\sin(\epsilon-\tfrac{1}{2}\alpha))\sin\tfrac{1}{2}(\delta-\epsilon) \nonumber \\
&-(1-\cos(\beta-\gamma))\sin\tfrac{1}{2}\alpha\sin\tfrac{1}{2}(\delta+\epsilon)]\cos\tfrac{1}{2}(\delta+\epsilon-\alpha)=0, \label{coolsaet_eq1} \\
& \sin\tfrac{1}{2}\alpha\sin\tfrac{1}{2}(\beta-\gamma)
(\sin\tfrac{1}{2}\beta\sin\delta\cos a
-\cos\tfrac{1}{2}\beta\cos\delta) \nonumber \\
& +\sin\tfrac{1}{2}\gamma
\sin\tfrac{1}{2}(\delta-\epsilon)
\cos\tfrac{1}{2}(\delta+\epsilon-\alpha)=0, \label{coolsaet_eq2} \\
& \sin\tfrac{1}{2}\alpha\sin\tfrac{1}{2}(\beta-\gamma)
(\sin\tfrac{1}{2}\gamma\sin\epsilon\cos a
-\cos\tfrac{1}{2}\gamma\cos\epsilon) \nonumber \\
& +\sin\tfrac{1}{2}\beta
\sin\tfrac{1}{2}(\delta-\epsilon)
\cos\tfrac{1}{2}(\delta+\epsilon-\alpha)=0. \label{coolsaet_eq3}
\end{align}
Conversely, if $\alpha,\beta,\gamma,\delta,\epsilon,a$ satisfy the three equalities, and $\alpha,\beta-\gamma\not\in 2{\bb Z}\pi$, and one of $\delta,\epsilon\not\in {\bb Z}\pi$, then there is an almost equilateral pentagon with the given $\alpha,\beta,\gamma,\delta,\epsilon,a$.
\end{lemma}

The first equality \eqref{coolsaet_eq1} in Lemma \ref{geometry7} means either $\delta+\epsilon-\alpha\in (2{\bb Z}+1)\pi$, or the following equality holds
\begin{align}
&
((1-\cos\beta)\sin(\delta-\tfrac{1}{2}\alpha)
-(1-\cos\gamma)\sin(\epsilon-\tfrac{1}{2}\alpha))\sin\tfrac{1}{2}(\delta-\epsilon) \nonumber \\
=\; & (1-\cos(\beta-\gamma))\sin\tfrac{1}{2}\alpha\sin\tfrac{1}{2}(\delta+\epsilon).  \label{coolsaet_eq15}
\end{align}

Suppose $\delta+\epsilon-\alpha\in (2{\bb Z}+1)\pi$. Then the three equalities become the following two
\begin{align*}
\sin\tfrac{1}{2}\alpha\sin\tfrac{1}{2}(\beta-\gamma)
(\sin\tfrac{1}{2}\beta\sin\delta\cos a
-\cos\tfrac{1}{2}\beta\cos\delta) &=0, \\
\sin\tfrac{1}{2}\alpha\sin\tfrac{1}{2}(\beta-\gamma)
(\sin\tfrac{1}{2}\gamma\sin\epsilon\cos a
-\cos\tfrac{1}{2}\gamma\cos\epsilon) &=0.
\end{align*}
If $\alpha,\beta-\gamma\not\in 2{\bb Z}\pi$, then we get
\begin{align*}
\sin\tfrac{1}{2}\beta\sin\delta\cos a
-\cos\tfrac{1}{2}\beta\cos\delta &=0, \\
\sin\tfrac{1}{2}\gamma\sin\epsilon\cos a
-\cos\tfrac{1}{2}\gamma\cos\epsilon &=0.
\end{align*}
This implies
\begin{equation}\label{coolsaet_eq10}
\cos\tfrac{1}{2}\beta\cos\delta\sin\tfrac{1}{2}\gamma\sin\epsilon
=\sin\tfrac{1}{2}\beta\sin\delta\cos\tfrac{1}{2}\gamma\cos\epsilon.
\end{equation}

The proof of Lemma \ref{geometry7} follows the same idea of the proof of the similar lemma for quadrilateral in \cite{cly}. 

\begin{proof}
Let $Y(\theta),Z(\theta)$ be the rotations of ${\bb R}^3$ with with respect to the $Y$- and $Z$-axes by angle $\theta$. Then the almost equilateral pentagon can be interpreted as a sequence of rotations satisfying
\begin{equation}\label{coolsaet_eq4}
Y(a)Z(\pi-\alpha)Y(a)Z(\pi-\beta)Y(a)Z(\pi-\delta)Y(b)Z(\pi-\epsilon)Y(a)Z(\pi-\gamma)=I.
\end{equation}
This is the same as 
\[
Y(b)Z(\pi-\epsilon)Y(a)Z(\pi-\gamma)Y(a)Z(\pi-\alpha)
=(Y(a)Z(\pi-\beta)Y(a)Z(\pi-\delta))^T.
\]
The $(2,3)$-entry of the matrix equality is $\sin a$ multiplied to the equality
\begin{equation}\label{coolsaet_eq5}
((1-\cos\beta)\sin\delta-(1-\cos\gamma)\sin\epsilon)\cos a
=\sin\beta\cos\delta-\sin\gamma\cos\epsilon.
\end{equation}
By $a\not\in {\bb Z}\pi$, the equality holds.

The matrix equation \eqref{coolsaet_eq4} is also the same as 
\[
Y(a)Z(\pi-\alpha)Y(a)Z(\pi-\beta)Y(a)Z(\pi-\delta)
=(Y(b)Z(\pi-\epsilon)Y(a)Z(\pi-\gamma))^T.
\]
The $(2,2)$- and $(3,2)$-entries of the matrix equality give
\begin{align}
& \sin\alpha(1-\cos\beta)\sin\delta\cos^2a \nonumber \\
&-(\sin\alpha\sin\beta\cos\delta+\cos\alpha\sin\beta\sin\delta+\sin\gamma\sin\epsilon)\cos a \nonumber  \\
& -\sin\alpha\sin\delta+\cos\alpha\cos\beta\cos\delta+\cos\gamma\cos\epsilon = 0, \label{coolsaet_eq6} \\
& (1-\cos\alpha)(1-\cos\beta)\sin\delta\cos^2a \nonumber  \\
& -((1-\cos\alpha)\cos\delta+\sin\alpha\sin\delta)\sin\beta\cos a \nonumber  \\
& +\cos\alpha\sin\delta+\sin\alpha\cos\beta\cos\delta-\sin\epsilon = 0. \label{coolsaet_eq7}
\end{align}
We actually get $\sin a$ multiplied to \eqref{coolsaet_eq7}, which can be dropped by $a\not\in {\bb Z}\pi$. 

We may cancel $\cos^2a$ from \eqref{coolsaet_eq6} and \eqref{coolsaet_eq7}, and get
\begin{align*}
&(1-\cos\alpha)(\sin\beta\sin\delta-\sin\gamma\sin\epsilon)\cos a \\
=\; &
\sin\alpha(\sin\delta-\sin\epsilon)
+(1-\cos\alpha)(\cos\beta\cos\delta-\cos\gamma\cos\epsilon).
\end{align*}
If $\alpha\not\in 2{\bb Z}\pi$, then we may divide $2\sin\tfrac{1}{2}\alpha$ and get
\begin{align*}
&\sin\tfrac{1}{2}\alpha(\sin\beta\sin\delta-\sin\gamma\sin\epsilon)\cos a  \\
=\; & 
\cos\tfrac{1}{2}\alpha(\sin\delta-\sin\epsilon)
+\sin\tfrac{1}{2}\alpha(\cos\beta\cos\delta-\cos\gamma\cos\epsilon) \nonumber \\
=\; &
\cos\tfrac{1}{2}\alpha(\sin\delta-\sin\epsilon)
-\sin\tfrac{1}{2}\alpha(\cos\delta-\cos\epsilon)  \\
&+\sin\tfrac{1}{2}\alpha((1+\cos\beta)\cos\delta 
-(1+\cos\gamma)\cos\epsilon) \\
=\; &
2\sin\tfrac{1}{2}(\delta-\epsilon)\cos\tfrac{1}{2}(\delta+\epsilon-\alpha)
+\sin\tfrac{1}{2}\alpha((1+\cos\beta)\cos\delta 
-(1+\cos\gamma)\cos\epsilon).
\end{align*}
We take the $\sin\tfrac{1}{2}\gamma$ multiple of the above and subtract the $\sin\tfrac{1}{2}\alpha\cos\tfrac{1}{2}\gamma$ multiple of \eqref{coolsaet_eq5}, and use the qualities
\[
\cos\tfrac{1}{2}\gamma(1-\cos\gamma)
=\sin\tfrac{1}{2}\gamma\sin\gamma,
\quad
\sin\tfrac{1}{2}\gamma(1+\cos\gamma)
=\cos\tfrac{1}{2}\gamma\sin\gamma,
\]
to get
\begin{align*}
& \sin\tfrac{1}{2}\alpha(\sin\beta\sin\tfrac{1}{2}\gamma-(1-\cos\beta)\cos\tfrac{1}{2}\gamma)\sin\delta\cos a \\
=\; & 
2\sin\tfrac{1}{2}(\delta-\epsilon)\cos\tfrac{1}{2}(\delta+\epsilon-\alpha)\sin\tfrac{1}{2}\gamma \\
& +\sin\tfrac{1}{2}\alpha\cos\delta((1+\cos\beta)\sin\tfrac{1}{2}\gamma-\sin\beta\cos\tfrac{1}{2}\gamma).
\end{align*}
Then by 
\begin{align*}
\sin\beta\sin\tfrac{1}{2}\gamma-(1-\cos\beta)\cos\tfrac{1}{2}\gamma
&=-2\sin\tfrac{1}{2}\beta\sin\tfrac{1}{2}(\beta-\gamma), \\
(1+\cos\beta)\sin\tfrac{1}{2}\gamma-\sin\beta\cos\tfrac{1}{2}\gamma
&=-2\cos\tfrac{1}{2}\beta\sin\tfrac{1}{2}(\beta-\gamma),
\end{align*}
we get \eqref{coolsaet_eq2}. We get \eqref{coolsaet_eq3} in the similar way.

Next, we derive \eqref{coolsaet_eq1} by substituting \eqref{coolsaet_eq2} into \eqref{coolsaet_eq7}. Denote \eqref{coolsaet_eq7} by $A_2\cos^2a+A_1\cos a+A_0=0$, and denote 
\[
X=\sin\tfrac{1}{2}(\delta-\epsilon)\cos\tfrac{1}{2}(\delta+\epsilon-\alpha)
=\tfrac{1}{2}(\sin(\delta-\tfrac{1}{2}\alpha)-\sin(\epsilon-\tfrac{1}{2}\alpha)).
\]
Then
\begin{align*}
&(A_2\cos a+A_1)\sin\tfrac{1}{2}(\beta-\gamma) \\
=\; & 
4\sin^2\tfrac{1}{2}\alpha\sin^2\tfrac{1}{2}\beta\sin\tfrac{1}{2}(\beta-\gamma)\sin\delta\cos a +A_1\sin\tfrac{1}{2}(\beta-\gamma) \\
=\; & 
4\sin\tfrac{1}{2}\alpha\sin\tfrac{1}{2}\beta(\sin\tfrac{1}{2}\alpha\cos\tfrac{1}{2}\beta\sin\tfrac{1}{2}(\beta-\gamma)\cos\delta-X\sin\tfrac{1}{2}\gamma) \\
&-((1-\cos\alpha)\cos\delta+\sin\alpha\sin\delta)\sin\beta\sin\tfrac{1}{2}(\beta-\gamma) \\
=\; &
-4X\sin\tfrac{1}{2}\alpha\sin\tfrac{1}{2}\beta\sin\tfrac{1}{2}\gamma
-\sin\alpha\sin\beta\sin\tfrac{1}{2}(\beta-\gamma)\sin\delta \\
=\; &
-4\sin\tfrac{1}{2}\alpha\sin\tfrac{1}{2}\beta
(X\sin\tfrac{1}{2}\gamma
+\cos\tfrac{1}{2}\alpha\cos\tfrac{1}{2}\beta\sin\tfrac{1}{2}(\beta-\gamma)\sin\delta).
\end{align*}
Here \eqref{coolsaet_eq2} is used in the second equality. Then
\begin{align*}
&(A_2\cos^2a+A_1\cos a)\sin^2\tfrac{1}{2}(\beta-\gamma)\sin\delta \\
=\; &
-4\sin\tfrac{1}{2}\alpha\sin\tfrac{1}{2}\beta\sin\tfrac{1}{2}(\beta-\gamma)\sin\delta\cos a \\
& \times (X\sin\tfrac{1}{2}\gamma
+\cos\tfrac{1}{2}\alpha\cos\tfrac{1}{2}\beta\sin\tfrac{1}{2}(\beta-\gamma)\sin\delta)  \\
=\; &
4(X\sin\tfrac{1}{2}\gamma
-\sin\tfrac{1}{2}\alpha\cos\tfrac{1}{2}\beta\sin\tfrac{1}{2}(\beta-\gamma)\cos\delta) \\
& \times (X\sin\tfrac{1}{2}\gamma
+\cos\tfrac{1}{2}\alpha\cos\tfrac{1}{2}\beta\sin\tfrac{1}{2}(\beta-\gamma)\sin\delta) \\
=\; &
4X^2\sin^2\tfrac{1}{2}\gamma
+4X\cos\tfrac{1}{2}\beta\sin\tfrac{1}{2}\gamma\sin\tfrac{1}{2}(\beta-\gamma)\sin(\delta-\tfrac{1}{2}\alpha) \\
&-\sin\alpha(1+\cos\beta)\sin^2\tfrac{1}{2}(\beta-\gamma)\sin\delta\cos\delta. 
\end{align*}
Here \eqref{coolsaet_eq2} is again used in the second equality. Then
\begin{align*}
& (A_2\cos^2a+A_1\cos a+A_0)\sin^2\tfrac{1}{2}(\beta-\gamma)\sin\delta \\
=\; & 
4X^2\sin^2\tfrac{1}{2}\gamma
+4X\cos\tfrac{1}{2}\beta\sin\tfrac{1}{2}\gamma\sin\tfrac{1}{2}(\beta-\gamma)\sin(\delta-\tfrac{1}{2}\alpha) \\
&-\sin\alpha(1+\cos\beta)\sin^2\tfrac{1}{2}(\beta-\gamma)\sin\delta\cos\delta \\
&+(\cos\alpha\sin\delta+\sin\alpha\cos\beta\cos\delta-\sin\epsilon)\sin^2\tfrac{1}{2}(\beta-\gamma)\sin\delta \\
=\; & 
4X^2\sin^2\tfrac{1}{2}\gamma
+4X\cos\tfrac{1}{2}\beta\sin\tfrac{1}{2}\gamma\sin\tfrac{1}{2}(\beta-\gamma)\sin(\delta-\tfrac{1}{2}\alpha) \\
&+2\cos\tfrac{1}{2}(\delta+\epsilon-\alpha)\sin\tfrac{1}{2}(\delta-\epsilon-\alpha)\sin^2\tfrac{1}{2}(\beta-\gamma)\sin\delta \\
=\; & 
W\cos\tfrac{1}{2}(\delta+\epsilon-\alpha).
\end{align*}
Here
\begin{align*}
W=\; &
2\sin\tfrac{1}{2}(\delta-\epsilon)(\sin(\delta-\tfrac{1}{2}\alpha)-\sin(\epsilon-\tfrac{1}{2}\alpha))\sin^2\tfrac{1}{2}\gamma \\
& +4\sin\tfrac{1}{2}(\delta-\epsilon)\cos\tfrac{1}{2}\beta\sin\tfrac{1}{2}\gamma\sin\tfrac{1}{2}(\beta-\gamma)\sin(\delta-\tfrac{1}{2}\alpha) \\
&+2\sin\tfrac{1}{2}(\delta-\epsilon-\alpha)\sin^2\tfrac{1}{2}(\beta-\gamma)\sin\delta \\
=\; & 
2\sin\tfrac{1}{2}(\delta-\epsilon)\sin(\delta-\tfrac{1}{2}\alpha)(\sin^2\tfrac{1}{2}\gamma+2\cos\tfrac{1}{2}\beta\sin\tfrac{1}{2}\gamma\sin\tfrac{1}{2}(\beta-\gamma)) \\
& -2\sin\tfrac{1}{2}(\delta-\epsilon)\sin(\epsilon-\tfrac{1}{2}\alpha)\sin^2\tfrac{1}{2}\gamma \\
& +2\sin\tfrac{1}{2}(\delta-\epsilon-\alpha)\sin^2\tfrac{1}{2}(\beta-\gamma)\sin\delta\\ 
=\; & 
2\sin\tfrac{1}{2}(\delta-\epsilon)\sin(\delta-\tfrac{1}{2}\alpha)(\sin^2\tfrac{1}{2}\beta-\sin^2\tfrac{1}{2}(\beta-\gamma)) \\
& -2\sin\tfrac{1}{2}(\delta-\epsilon)\sin(\epsilon-\tfrac{1}{2}\alpha)\sin^2\tfrac{1}{2}\gamma \\
& +2\sin\tfrac{1}{2}(\delta-\epsilon-\alpha)\sin^2\tfrac{1}{2}(\beta-\gamma)\sin\delta\\
=\; & 
2\sin^2\tfrac{1}{2}\beta\sin\tfrac{1}{2}(\delta-\epsilon)\sin(\delta-\tfrac{1}{2}\alpha) 
-2\sin^2\tfrac{1}{2}\gamma\sin\tfrac{1}{2}(\delta-\epsilon)\sin(\epsilon-\tfrac{1}{2}\alpha) \\
& +2\sin^2\tfrac{1}{2}(\beta-\gamma)(\sin\tfrac{1}{2}(\delta-\epsilon-\alpha)\sin\delta-\sin\tfrac{1}{2}(\delta-\epsilon)\sin(\delta-\tfrac{1}{2}\alpha))
\\
=\; & 
((1-\cos\beta)\sin(\delta-\tfrac{1}{2}\alpha)
-(1-\cos\gamma)\sin(\epsilon-\tfrac{1}{2}\alpha))\sin\tfrac{1}{2}(\delta-\epsilon) \\
& -(1-\cos(\beta-\gamma))\sin\tfrac{1}{2}\alpha\sin\tfrac{1}{2}(\delta+\epsilon).
\end{align*}
By $A_2\cos^2a+A_1\cos a+A_0=0$, we get \eqref{coolsaet_eq1}.

The argument assumes $\alpha\not\in 2{\bb Z}\pi$. If $\alpha\in 2{\bb Z}\pi$, then $B=C$, and $\triangle BDE$ is an isosceles triangle with top angle $\beta+\gamma$ mod $2\pi$ and base angle $\delta=\epsilon$ mod $2\pi$. Then \eqref{coolsaet_eq1}, \eqref{coolsaet_eq2}, \eqref{coolsaet_eq3} remain true. 

Conversely, suppose \eqref{coolsaet_eq1}, \eqref{coolsaet_eq2}, \eqref{coolsaet_eq3} hold. We will show the $(2,2)$-entry of the matrix
\begin{equation}\label{coolsaet_eq8}
K=Z(\pi-\epsilon)Y(a)Z(\pi-\gamma)Y(a)Z(\pi-\alpha)Y(a)Z(\pi-\beta)Y(a)Z(\pi-\delta)
\end{equation}
is $1$. This implies $K=Y(b)^T$ for a unique $b$ mod $2\pi$. Then \eqref{coolsaet_eq4} holds, and we have the corresponding pentagon.

To show the $(2,2)$-entry of $K$ is $1$, we introduce
\[
L=Y(a)Z(\pi-\alpha)Y(a)Z(\pi-\beta)Y(a)Z(\pi-\delta)-Z(\pi-\gamma)^TY(a)^TZ(\pi-\epsilon)^T.
\]
We will prove $L\vec{e}_2=\vec{0}$, where $\vec{e}_2=(0,1,0)^T$ is the second standard basis vector. This implies $K\vec{e}_2=\vec{e}_2$, and further implies the $(2,2)$-entry of $K$ is $1$.

We first argue that \eqref{coolsaet_eq1}, \eqref{coolsaet_eq2}, \eqref{coolsaet_eq3} imply \eqref{coolsaet_eq5}, \eqref{coolsaet_eq6}, \eqref{coolsaet_eq7}, under the assumption $\alpha,\gamma,\beta-\gamma\not\in 2{\bb Z}\pi$ and $\delta\not\in {\bb Z}\pi$.

We multiply \eqref{coolsaet_eq2} by $\sin\tfrac{1}{2}\beta$, and multiply \eqref{coolsaet_eq3} by $\sin\tfrac{1}{2}\gamma$. Then we take the difference and use $\alpha,\beta-\gamma\not\in 2{\bb Z}\pi$ to get
\[
\sin^2\tfrac{1}{2}\beta\sin\delta\cos a
-\sin\tfrac{1}{2}\beta\cos\tfrac{1}{2}\beta\cos\delta
=\sin^2\tfrac{1}{2}\gamma\sin\epsilon\cos a
-\sin\tfrac{1}{2}\gamma\cos\tfrac{1}{2}\gamma\cos\epsilon.
\]
This is \eqref{coolsaet_eq5}. The earlier argument (specifically, the calculation of $A_2\cos^2a+A_1\cos a+A_0$, with the help of \eqref{coolsaet_eq2}) shows that, by $\beta-\gamma\not\in 2{\bb Z}\pi$ and $\delta\not\in {\bb Z}\pi$, the equalities \eqref{coolsaet_eq1} and \eqref{coolsaet_eq2} imply \eqref{coolsaet_eq7}. Then the earlier argument also shows that, by $\alpha,\gamma\not\in 2{\bb Z}\pi$, the equalities \eqref{coolsaet_eq2}, \eqref{coolsaet_eq5}, \eqref{coolsaet_eq7} imply the equality \eqref{coolsaet_eq6}.

Next, we use \eqref{coolsaet_eq5}, \eqref{coolsaet_eq6}, \eqref{coolsaet_eq7} to prove $L\vec{e}_2=\vec{0}$. 

The second coordinate of $L\vec{e}_2$ is exactly the negative of the left of \eqref{coolsaet_eq6}. The third coordinate of $L\vec{e}_2$ is exactly the left of \eqref{coolsaet_eq7}. Therefore the two coordinates are $0$.

The first coordinate of $L\vec{e}_2$ is $B_3\cos^3a+B_2\cos^2a+B_1\cos a+B_0$, where
\begin{align*}
B_3 &
=-(1-\cos\alpha)(1-\cos\beta)\sin\delta, \\
B_2 &
=((1-\cos\alpha)\sin\beta\cos\delta+\sin\alpha\sin\beta\sin\delta, \\
B_1 &
=-\sin\alpha\cos\beta\cos\delta+(1-\cos\alpha-\cos\beta)\sin\delta+\cos\gamma\sin\epsilon, \\
B_0 &
=-\sin\beta\cos\delta+\sin\gamma\cos\epsilon.
\end{align*}
Then
\begin{align*}
& B_3\cos^3a+B_2\cos^2a+B_1\cos a+B_0 \\
=\; & (B_3\cos^2a+B_2\cos a)\cos a+B_1\cos a+B_0 \\
=\; & (\cos\alpha\sin\delta+\sin\alpha\cos\beta\cos\delta-\sin\epsilon)\cos a+B_1\cos a+B_0 \\
=\; & ((1-\cos\beta)\sin\delta-(1-\cos\gamma)\sin\epsilon)\cos a-\sin\beta\cos\delta+\sin\gamma\cos\epsilon \\
=\; & 
0.
\end{align*}
Here \eqref{coolsaet_eq7} is used in the second equality, and \eqref{coolsaet_eq5} is used in the fourth equality.

We proved the converse when $\alpha,\gamma,\beta-\gamma\not\in 2{\bb Z}\pi$ and $\delta\not\in{\bb Z}\pi$. Since the condition $\delta\not\in{\bb Z}\pi$ was used only in proving that \eqref{coolsaet_eq1} and \eqref{coolsaet_eq2} imply \eqref{coolsaet_eq7}, we still need to prove that, in case $\alpha,\gamma,\beta-\gamma\not\in 2{\bb Z}\pi$ and $\delta\in{\bb Z}\pi$, the equalities \eqref{coolsaet_eq1}, \eqref{coolsaet_eq2}, \eqref{coolsaet_eq3} still imply \eqref{coolsaet_eq7}. We remark that, by the assumption of the lemma and $\delta\in{\bb Z}\pi$, we know $\epsilon\not\in{\bb Z}\pi$.

By $\delta\in {\bb Z}\pi$, the equality \eqref{coolsaet_eq7} becomes
\begin{equation}\label{coolsaet_eq12}
(1-\cos\alpha)\sin\beta\cos a
-\sin\alpha\cos\beta
+(-1)^d\sin\epsilon=0,\quad
d=\tfrac{\delta}{\pi}.
\end{equation}
Presumably, the property $\delta\in {\bb Z}\pi$ means the pentagon is reduced to the almost equilateral quadrilateral in Figure \ref{coolsaet_f3}. Then \eqref{coolsaet_eq7} is the third equality in \cite[Lemma 18]{cly}, for this quadrilateral. However, we should derive \eqref{coolsaet_eq12} solely from \eqref{coolsaet_eq1}, \eqref{coolsaet_eq2}, \eqref{coolsaet_eq3}.

\begin{figure}[htp]
\centering
\begin{tikzpicture}[>=latex]


\draw
	(1.3,0) -- (0,0) -- (0.5,1.2) -- (1.7,1) -- (2.5,0);
	
\draw[line width=1.2]
	(1.3,0) -- (2.5,0);

\node at (0.6,1) {\small $\alpha$};
\node at (0.25,0.2) {\small $\beta$};
\node at (1.6,0.8) {\small $\gamma$};
\node at (1.3,0.2) {\small $\delta$};
\node at (2.1,0.2) {\small $\epsilon$};

\node at (1.3,-0.6) {\small $\delta\in (2{\bb Z}+1)\pi$};


\begin{scope}[xshift=-3.5cm]

\draw
	(1,0.1) -- (0.5,1.2) -- (1.7,1) -- (2.5,0)
	(1,0.1) -- (-0.3,0.1) arc (90:270:0.05)
	(-0.1,0.1) arc (40:-40:0.08)
	(0.8,0.1) arc (-180:115:0.2);
	
\draw[line width=1.2]
	(-0.3,0) -- (2.5,0);
	
\node at (0.8,1) {\small $\alpha$};
\node at (1.3,0.3) {\small $\beta$};
\node at (1.6,0.8) {\small $\gamma$};
\node at (-0.1,0.3) {\small $\delta$};
\node at (2.1,0.2) {\small $\epsilon$};

\node at (1.1,-0.6) {\small $\delta\in 2{\bb Z}\pi$};

\end{scope}

\end{tikzpicture}
\caption{Lemma \ref{geometry7}: Case $\delta\in {\bb Z}\pi$.}
\label{coolsaet_f3}
\end{figure}
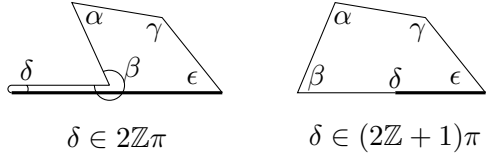

By $\delta\in {\bb Z}\pi$, the equality \eqref{coolsaet_eq2} becomes 
\[
2\sin\tfrac{1}{2}\alpha\sin\tfrac{1}{2}(\beta-\gamma)\cos\tfrac{1}{2}\beta
+\sin\tfrac{1}{2}\gamma(\sin\tfrac{1}{2}\alpha+(-1)^d\sin(\epsilon-\tfrac{1}{2}\alpha))=0.
\]
By
\begin{align*}
& 2\sin\tfrac{1}{2}\alpha\sin\tfrac{1}{2}(\beta-\gamma)\cos\tfrac{1}{2}\beta+\sin\tfrac{1}{2}\gamma\sin\tfrac{1}{2}\alpha \\
=\; &
\sin\tfrac{1}{2}\alpha(\sin(\beta-\tfrac{1}{2}\gamma)-\sin\tfrac{1}{2}\gamma)+\sin\tfrac{1}{2}\gamma\sin\tfrac{1}{2}\alpha \\
=\; &
\sin\tfrac{1}{2}\alpha\sin(\beta-\tfrac{1}{2}\gamma),
\end{align*}
we get
\begin{equation}\label{coolsaet_eq13}
\sin\tfrac{1}{2}\alpha\sin(\beta-\tfrac{1}{2}\gamma)
+(-1)^d\sin\tfrac{1}{2}\gamma\sin(\epsilon-\tfrac{1}{2}\alpha)
=0.
\end{equation}
This is exactly the first equality in \cite[Lemma 18]{cly} (also \cite[Theorem 2.1]{coolsaet}).

For $\delta\in {\bb Z}\pi$, we multiply $\sin\frac{1}{2}\gamma$ to \eqref{coolsaet_eq3} and get 
\begin{align*}
& \sin\tfrac{1}{2}\alpha\sin\tfrac{1}{2}(\beta-\gamma)
((1-\cos\gamma)\sin\epsilon\cos a
-\sin\gamma\cos\epsilon)  \\
& -\sin\tfrac{1}{2}\beta\sin\tfrac{1}{2}\gamma((-1)^d\sin\tfrac{1}{2}\alpha+\sin(\epsilon-\tfrac{1}{2}\alpha))=0.
\end{align*}
Then
\begin{align*}
& \sin\tfrac{1}{2}\alpha\sin\tfrac{1}{2}(\beta-\gamma)
((1-\cos\gamma)\sin\epsilon\cos a
-\sin\gamma\cos\epsilon+(-1)^d\sin\beta)  \\
=\; & 
(-1)^d\sin\tfrac{1}{2}\alpha\sin\tfrac{1}{2}(\beta-\gamma)\sin\beta
+\sin\tfrac{1}{2}\beta\sin\tfrac{1}{2}\gamma((-1)^d\sin\tfrac{1}{2}\alpha+\sin(\epsilon-\tfrac{1}{2}\alpha)) \\
=\; & 
(-1)^d\sin\tfrac{1}{2}\alpha\sin\tfrac{1}{2}\beta
(2\sin\tfrac{1}{2}(\beta-\gamma)\cos\tfrac{1}{2}\beta
+\sin\tfrac{1}{2}\gamma)
+\sin\tfrac{1}{2}\beta\sin\tfrac{1}{2}\gamma\sin(\epsilon-\tfrac{1}{2}\alpha) \\
=\; & 
(-1)^d\sin\tfrac{1}{2}\alpha\sin\tfrac{1}{2}\beta\sin(\beta-\tfrac{1}{2}\gamma)
+\sin\tfrac{1}{2}\beta\sin\tfrac{1}{2}\gamma\sin(\epsilon-\tfrac{1}{2}\alpha) 
=0.
\end{align*}
Here the equality \eqref{coolsaet_eq13} is used in the last equality. By $\alpha,\beta-\gamma\not\in 2{\bb Z}\pi$, we get 
\begin{equation}\label{coolsaet_eq14}
(1-\cos\gamma)\sin\epsilon\cos a
-\sin\gamma\cos\epsilon+(-1)^d\sin\beta
=0.
\end{equation}
This is exactly the second equality in \cite[Lemma 18]{cly}.

By $\gamma\not\in 2{\bb Z}\pi$ and $\epsilon\not\in {\bb Z}\pi$, and the remark after \cite[Lemma 18]{cly}, we know \eqref{coolsaet_eq13} and \eqref{coolsaet_eq14} imply \eqref{coolsaet_eq12}. This establishes \eqref{coolsaet_eq12} solely from \eqref{coolsaet_eq2}, \eqref{coolsaet_eq3}. In fact, unlike the case $\delta\not\in {\bb Z}\pi$, \eqref{coolsaet_eq1} is not used in the argument. 

This finishes the proof of the converse for the case $\alpha,\gamma,\beta-\gamma\not\in 2{\bb Z}\pi$, and one of $\delta,\epsilon\not\in {\bb Z}\pi$. Similarly, we may prove the converse for the case $\alpha,\beta,\beta-\gamma\not\in 2{\bb Z}\pi$, and one of $\delta,\epsilon\not\in {\bb Z}\pi$. Since $\beta-\gamma\not\in 2{\bb Z}\pi$ implies one of $\beta,\gamma\not\in 2{\bb Z}\pi$, we conclude that the converse holds under the assumption $\alpha,\beta-\gamma\not\in 2{\bb Z}\pi$, and one of $\delta,\epsilon\not\in {\bb Z}\pi$. 
\end{proof}

We explain that \eqref{coolsaet_eq1} is not satisfied if $a\in {\bb Z}\pi$.

\begin{figure}[htp]
\centering
\begin{tikzpicture}[>=latex]


\begin{scope}[xshift=-3.5cm]

\draw[line width=1.2]
	(0,0) circle (1.2);
	
\coordinate (AA) at (-1.2,0);
\coordinate (A) at (1.2,0);

\coordinate (X) at (0,0.7);
\coordinate (Y) at (0,-0.7);
\coordinate (Z) at (0,0.3);
\coordinate (W) at (0,-0.2);

\arcThroughThreePoints{A}{X}{AA};
\arcThroughThreePoints{AA}{Y}{A};
\arcThroughThreePoints{A}{Z}{AA};
\arcThroughThreePoints{AA}{W}{A};

\draw[->] (0,0.7) -- ++(-0.02,0);
\draw[->] (0,-0.2) -- ++(-0.02,0);
\draw[->] (0,-0.7) -- ++(0.02,0);
\draw[->] (0,0.3) -- ++(0.02,0);
\draw[->] (0,1.2) -- ++(-0.02,0);
\draw[->] (0,-1.2) -- ++(0.02,0);

\node at (0,0.9) {\small $BA$};
\node at (0,-0.4) {\small $CE$};
\node at (0,-0.9) {\small $CA$};
\node at (0,0.5) {\small $BD$};

\node at (0.8,0.65) {\small $\delta$};
\node at (0.8,-0.62) {\small $\epsilon$};
\node at (-0.6,0.38) {\small $\beta$};
\node at (-0.65,-0.35) {\small $\gamma$};
\node at (0.55,0.05) {\small $\alpha$};

\draw[xshift=1.2cm]
	(130:0.5) arc (130:230:0.5);
\draw[xshift=1.2cm]
	(106:0.6) arc (106:159:0.6);
\draw[xshift=1.2cm]
	(-106:0.6) arc (-106:-165:0.6);

\draw[xshift=-1.2cm]
	(50:0.5) arc (50:22:0.5);
\draw[xshift=-1.2cm]
	(-50:0.5) arc (-50:-16:0.5);

\node at (0,-1.6) {\small $a\in (2{\bb Z}+1)\pi$};
	
\end{scope}


\begin{scope}[xshift=0cm]

\draw
	(0:1) -- (0:-1)
	(40:1) -- (40:-1)
	(100:0.9) -- (100:-0.9)
	(160:1) -- (160:-1);

\draw[line width=1.2]
	(60:1) -- (60:-1);

\draw[->] (0:0.7) -- ++(0:-0.02);
\draw[->] (100:-0.7) -- ++(100:-0.02);
\draw[->] (100:0.7) -- ++(100:-0.02);
\draw[->] (40:0.7) -- ++(40:0.02);
\draw[->] (40:-0.7) -- ++(40:0.02);
\draw[->] (60:-0.7) -- ++(60:-0.02);
\draw[->] (60:0.7) -- ++(60:-0.02);
\draw[->] (0:-0.7) -- ++(0:-0.02);
\draw[->] (160:0.7) -- ++(160:0.02);
\draw[->] (160:-0.7) -- ++(160:0.02);

\node at (0:1.2) {\small $A_i$};	
\node at (100:-1.1) {\small $A_o$};	
\node at (100:1.05) {\small $B_i$};		
\node at (35:1.2) {\small $B_o$};
\node at (160:-1.25) {\small $C_i$};
\node at (0:-1.2) {\small $C_o$};	
\node at (35:-1.2) {\small $D_i$};	
\node at (60:-1.2) {\small $D_o$};
\node at (60:1.2) {\small $E_i$};
\node at (163:1.2) {\small $E_o$};

\draw
	(0:0.46) arc (0:-80:0.46)
	(-20:0.31) arc (-20:-180:0.31)
	(40:0.31) arc (40:100:0.31)
	(220:0.22) arc (220:-120:0.22)
	(60:0.4) arc (60:-200:0.4);
	
\node at (-40:0.58) {\small $\alpha$};
\node at (80:0.5) {\small $\beta$};

\node at (0,-1.6) {\small $a\in 2{\bb Z}\pi$};

\end{scope}

\end{tikzpicture}
\caption{Lemma \ref{geometry7}: Case $a\in {\bb Z}\pi$.}
\label{coolsaet_f1}
\end{figure}
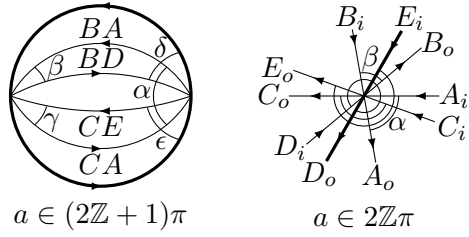

Suppose $a\in(2{\bb Z}+1)\pi$. Then we get the first of Figure \ref{coolsaet_f1}, where $A=D=E$ is the right vertex, and $B=C$ is the left vertex. Moreover, we have $b\in 2{\bb Z}\pi$, and the only relation between angles is $\alpha+\delta+\epsilon=\beta+\gamma+\pi$ mod $2\pi$. The relation does not imply \eqref{coolsaet_eq1}.

Suppose $a\in 2{\bb Z}\pi$. Then we get the second of Figure \ref{coolsaet_f1}, where $A=B=C=D=E$ is the only vertex, and $b\in 2{\bb Z}\pi$. Moreover, an edge $XY$ leaves $X$ through $X_o$ (for \underline{o}ut of $X$) and arrives at $Y$ through $Y_i$ (for \underline{i}nto $Y$). If we also denote by $X_i,X_o$ the angles of the direction of $X_i,X_o$, then we have $Y_i=X_o+\pi$ mod $2\pi$. Moreover, we have $\alpha=A_i-A_o=C_o-A_o+\pi$ mod $2\pi$, and similar formulae for the other angles. Then the only relation between the angles is $\alpha+\beta+\gamma+\delta+\epsilon=\pi$ mod $2\pi$. The relation does not imply \eqref{coolsaet_eq1}.

The converse in Lemma \ref{geometry7} requires extra conditions on the angle values. We investigate the necessity of these conditions. 

Suppose $\alpha\in 2{\bb Z}\pi$ or $\beta-\gamma\in 2{\bb Z}\pi$. Then \eqref{coolsaet_eq1}, \eqref{coolsaet_eq2}, \eqref{coolsaet_eq3} become
\begin{align*}
&
((1-\cos\beta)\sin(\delta-\tfrac{1}{2}\alpha)
-(1-\cos\gamma)\sin(\epsilon-\tfrac{1}{2}\alpha)) \\
& \quad\sin\tfrac{1}{2}(\delta-\epsilon)\cos\tfrac{1}{2}(\delta+\epsilon-\alpha)=0, \\
& \sin\tfrac{1}{2}\gamma
\sin\tfrac{1}{2}(\delta-\epsilon)
\cos\tfrac{1}{2}(\delta+\epsilon-\alpha)=0, \\
& \sin\tfrac{1}{2}\beta
\sin\tfrac{1}{2}(\delta-\epsilon)
\cos\tfrac{1}{2}(\delta+\epsilon-\alpha)=0.
\end{align*}
Since $\sin\tfrac{1}{2}\beta,\sin\tfrac{1}{2}\gamma$ are factors of $1-\cos\beta,1-\cos\gamma$, the second and third equalities imply the first. Therefore the three equalities is equivalent to the second and third equalities, which means either $\beta,\gamma\in 2{\bb Z}\pi$, or $\delta-\epsilon\in 2{\bb Z}\pi$, or $\delta+\epsilon-\alpha\in (2{\bb Z}+1)\pi$.

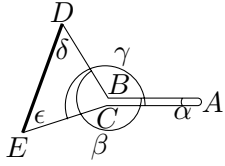
\begin{figure}[htp]
\centering
\begin{tikzpicture}[>=latex]


\draw
	(200:1.2) -- (-0.01,-0.05) -- (1.2,-0.05) arc (-90:90:0.05) -- (0.01,0.05) -- (120:1.2)
	(1,0.05) arc (140:220:0.08)
	(0.45,0.05) arc (0:-238:0.43)
	(0.5,-0.05) arc (0:200:0.53);
	
\draw[line width=1.2]
	(120:1.2) -- (200:1.2);

\node at (1.4,0) {\small $A$};
\node at (0.15,0.23) {\small $B$};
\node at (0,-0.22) {\small $C$};
\node at (-0.6,1.2) {\small $D$};
\node at (-1.2,-0.6) {\small $E$};

\node at (1,-0.15) {\small $\alpha$};
\node at (-0.1,-0.6) {\small $\beta$};
\node at (0.2,0.6) {\small $\gamma$};
\node at (-0.6,0.7) {\small $\delta$};
\node at (-0.9,-0.2) {\small $\epsilon$};

\end{tikzpicture}
\caption{Lemma \ref{geometry7}: Cases $\alpha\in 2{\bb Z}\pi$.}
\label{coolsaet_f2}
\end{figure}

Geometrically, if $\alpha\in 2{\bb Z}\pi$, then the pentagon is given by Figure \ref{coolsaet_f2}. We have $B=C$ and the triangle $\triangle BDE$ is isosceles, with top angle $\beta+\gamma$ mod $2\pi$ and the base angle $\delta=\epsilon$ mod $2\pi$. If $\beta-\gamma\in 2{\bb Z}\pi$, then for $b\not\in 2{\bb Z}\pi$, by Lemma \ref{geometry11}, we again get $\delta-\epsilon\in 2{\bb Z}\pi$, and the pentagon is symmetric. However, we may have $\beta,\gamma\in 2{\bb Z}\pi$, or $\delta+\epsilon-\alpha\in (2{\bb Z}+1)\pi$, for \eqref{coolsaet_eq1}, \eqref{coolsaet_eq2}, \eqref{coolsaet_eq3} to be satisfied. Together with $\alpha\in 2{\bb Z}\pi$ or $\beta-\gamma\in 2{\bb Z}\pi$, these angle conditions do not imply $\delta-\epsilon\in 2{\bb Z}\pi$, and therefore do not imply the existence of the pentagon.

We remark there is a subtle difference between Lemma \ref{geometry7} and its quadrilateral counterpart \cite[Lemma 18]{cly}. In the quadrilateral lemma, the first equality (relating the four angles) is almost the necessary and sufficient condition for the existence of the quadrilateral in the following sense: If the first equality is valid, then we may use either the second or the third equality to calculate $\cos a$. In fact, by the remark after \cite[Lemma 18]{cly}, the first equality implies the second and third equalities are almost equivalent. The only thing we need to verify is that the expected value of $\cos a$ has absolute value $\le 1$. In practise, this means that for given angle values of an expected quadrilateral, we only need to verify the first equality.

In the pentagonal lemma, however, there are cases that the first equality becomes trivial, and non-trivial conditions can be derived from the second and third equalities. 

Finally, we discuss the existence of pentagon in a tiling, especially the simple property is satisfied. In the earlier works on pentagonal tilings \cite{awy,wy1,wy2}, after tilings were obtained by combinatorial argument, further geometrical argument was made to show the existence of the pentagon. 

In the subsequent discussion, we always assume $0<a<\pi$, $0<b<2\pi$, $0<\alpha,\beta,\gamma,\delta<2\pi$. The recipe for arguing the existence is the following. First, we make sure the equality \eqref{coolsaet_eq1} holds. Then for the case $\beta\ne\gamma$ and one of $\delta,\epsilon\ne\pi$, we use \eqref{coolsaet_eq2} or \eqref{coolsaet_eq3} to uniquely determine $a<\pi$. Of course getting $a$ this way requires us to verify that the expected value of $\cos a$ has absolute value $<1$. Then by the converse of Lemma \ref{geometry7}, we already know a pentagon with the given angles and the deduced $a$ exists. The edge length $b$ is uniquely determined by $Y(b)=K^T$, where $K$ is given by \eqref{coolsaet_eq8}. If angles and edges satisfy the following, then we know the pentagon is simple.

\begin{lemma}\label{geometry8}
Suppose $\alpha,\beta,\gamma<\pi$. If $\alpha$ is bigger than the base angles of the isosceles triangles $\triangle ABD$ and $\triangle ACE$, and $(\beta-\gamma)(\delta-\epsilon)<0$, then the pentagon is simple.
\end{lemma}

The condition $(\beta-\gamma)(\delta-\epsilon)<0$ means exactly that the property in Lemma \ref{geometry1} holds.

\begin{proof}
In Figure \ref{geom_proof3}, we form the $2$-gon with angle $\alpha$. Then $B,C$ lie on the two boundary arcs of the $2$-gon. Since $\beta,\gamma<\pi$, and $\alpha$ is bigger than the base angles of the isosceles triangles $\triangle ABD$ and $\triangle ACE$, we know $D,E$ lie inside the $2$-gon. 

We always use black thick line to indicate the minor arc connecting $D,E$, and use gray thick line to indicate the other arc connecting $D,E$. Moreover, the angles $\delta,\epsilon$ are for the pentagon with the black thick line.

Assume $\beta>\gamma$. Then the condition $(\beta-\gamma)(\delta-\epsilon)<0$ implies $\delta<\epsilon$, and we get three possibilities: $D,E$ are outside $\triangle ABC$; $D$ is outside $\triangle ABC$ and $E$ is inside $\triangle ABC$; and $D,E$ are inside $\triangle ABC$. 

Suppose $D,E$ are outside $\triangle ABC$. Then we consider further three possibilities: Let $X$ be the intersection of (the extensions of) $BD$ and $CE$. The first picture is $X$ is not in $BD$ and $CE$; the second is $X$ is in one of $BD$ and $CE$; and the third is $X$ is in both $BD$ and $CE$. 

The first of Figure \ref{geom_proof3} is the case $D,E$ are outside $\triangle ABC$, and $X$ is not in $BD$ and $CE$. In this case, the black $DE$ is inside $\triangle A^*BC$.  We get $\delta,\epsilon<\pi$, and the pentagon is the union of simple $\triangle ABC$ and simple $\square BCED$ that do not overlap. Therefore the pentagon is simple. 

If we use grey $DE$ instead of the black one, then by $\alpha,\beta,\gamma,\delta,\epsilon<\pi$, the pentagon with black $DE$ lies in one hemisphere bounded by $\bigcirc DE$. Moreover, the pentagon with grey $DE$ is the union of the pentagon with black $DE$ with the other hemisphere. Therefore the pentagon with grey $DE$ is also simple.

\begin{figure}[htp]
\centering
\begin{tikzpicture}[>=latex, scale=1]

\foreach \b in {0,1,2,3}
{
\begin{scope}[xshift=2.5*\b cm]

\foreach \a in {-1,1}
{
\begin{scope}[xscale=\a]
	
\draw[gray!50]
	(0.75,0.6) to[out=-50, in=50] (0.75,-0.6) -- (0,-1.5);
\draw
	(0,1.5) -- (0.75,0.6);
	
\end{scope}
}

\draw[dashed]
	(0.75,0.6) -- (-0.75,0.6);
	
\node at (0,1.25) {\small $\alpha$};

\node at (0,1.7) {\small $A$};
\node at (-0.9,0.7) {\small $B$};
\node at (0.95,0.7) {\small $C$};
\node at (0,-1.7) {\small $A^*$};

\end{scope}
}


\draw[line width=1.2, gray!50]
	(0.3,-0.5) -- ++(0.7,0.1)
	(-0.4,-0.6) -- ++(-0.7,-0.1);

\draw[line width=1.2]
	(0.3,-0.5) -- (-0.4,-0.6);

\draw
	(0.75,0.6) -- (0.3,-0.5)
	(-0.75,0.6) -- (-0.4,-0.6)
	;

\node at (-0.5,0.6) {\small $\beta$};
\node at (0.5,0.6) {\small $\gamma$};	
\node at (-0.25,-0.35) {\small $\delta$};
\node at (0.2,-0.3) {\small $\epsilon$};

\node at (-0.4,-0.8) {\small $D$};
\node at (0.3,-0.7) {\small $E$};	
	

\begin{scope}[xshift=2.5cm]

\draw[gray!50]
	(-0.1,0.3) -- ++(18:-0.37);
	
\draw[line width=1.2, gray!50]
	(0.1,-0.5) -- ++(0.2,-0.8)
	(-0.1,0.3) -- ++(-0.2,0.8);

\draw[line width=1.2]
	(0.1,-0.5) -- (-0.1,0.3);

\draw
	(0.75,0.6) -- (-0.1,0.3) 
	(-0.75,0.6) -- (0.1,-0.5);
	
\node at (-0.5,0.65) {\small $\beta$};
\node at (0.45,0.65) {\small $\gamma$};	
\node at (-0.15,0.03) {\small $\delta$};
\node at (-0.2,0.35) {\small $\epsilon$};

\node at (0.25,-0.55) {\small $D$};
\node at (0.15,0.2) {\small $E$};	

\node at (-0.65,0.1) {\small $X$};
	
\end{scope}


\begin{scope}[xshift=5cm]

\draw 
	(0.75,0.6) -- (-0.5,-0.2)
	(-0.75,0.6) -- (0.3,-0.3);

\draw[line width=1.2, gray!50]
	(-0.5,-0.2) -- ++(-0.8,0.1)
	(0.3,-0.3) -- ++(0.8,-0.1);

\draw[line width=1.2]
	(-0.5,-0.2) -- (0.3,-0.3);

\node at (-0.5,0.65) {\small $\beta$};
\node at (0.45,0.65) {\small $\gamma$};
\node at (0.4,-0.3) {\small $\delta$};
\node at (-0.6,-0.2) {\small $\epsilon$};

\node at (-0.05,0.3) {\small $X$};
	
\end{scope}


\begin{scope}[xshift=7.5cm]

\draw[line width=1.2, gray!50]
	(0.1,-0.5) -- ++(0.3,-1.3)
	(-0.2,0.8) -- ++(-0.15,0.65);
	
\draw[line width=1.2]
	(0.1,-0.5) -- (-0.2,0.8);
	
\draw
	(0.75,0.6) -- (-0.2,0.8)
	(-0.75,0.6) -- (0.1,-0.5);
	
\node at (-0.5,0.65) {\small $\beta$};
\node at (0.35,0.83) {\small $\gamma$};
\node at (-0.15,0.03) {\small $\delta$};
\node at (-0.25,0.9) {\small $\epsilon$};

\node at (0.25,-0.55) {\small $D$};
\node at (0.05,0.55) {\small $E$};	

\end{scope}


\begin{scope}[xshift=10cm]

\foreach \a in {-1,1}
{
\begin{scope}[xscale=\a]
	
\draw[gray!50]
	(0.75,-0.6) -- (0,-1.5);
\draw
	(0,1.5) -- (0.75,0.6) to[out=-50, in=50] (0.75,-0.6);
	
\end{scope}
}

\draw[dashed]
	(0.75,-0.6) -- (-0.75,-0.6);
	
\node at (0,1.25) {\small $\alpha$};

\node at (0,1.7) {\small $A$};
\node at (-0.95,-0.7) {\small $B$};
\node at (0.9,-0.7) {\small $C$};
\node at (0,-1.7) {\small $A^*$};

\draw[line width=1.2, gray!50]
	(-0.3,0.5) -- ++(-0.7,-0.1)
	(0.4,0.6) -- ++(0.7,0.1);

\draw[line width=1.2]
	(-0.3,0.5) -- (0.4,0.6);

\draw
	(-0.75,-0.6) -- (-0.3,0.5)
	(0.75,-0.6) -- (0.4,0.6)
	;

\node at (-0.75,-0.3) {\small $\beta$};
\node at (0.75,-0.3) {\small $\gamma$};	
\node at (-0.35,0.65) {\small $\delta$};
\node at (0.45,0.75) {\small $\epsilon$};

\node at (-0.2,0.3) {\small $D$};
\node at (0.25,0.4) {\small $E$};

\end{scope}
	
\end{tikzpicture}
\caption{Lemma \ref{geometry8}.}
\label{geom_proof3}
\end{figure}

The second of Figure \ref{geom_proof3} is the case $D,E$ are outside $\triangle ABC$, and $X$ is in $BD$. Then we get $BX<BD=a=CE<CX$. By the isosceles triangle $\triangle ABC$, this implies $\beta>\gamma$. This actually proves that, if $\pi>\beta>\gamma$, then $X$ cannot be on $CE$ and cannot be on $BD$. In other words, the second case must have $X$ in $BD$ as in the picture. This implies the pentagon with black $CE$ is simple. Moreover, the pentagon with grey $CE$ has the $\delta$-angle $\delta+\pi$ and the $\epsilon$-angle $\epsilon-\pi$. By $\delta+\pi>\pi>\epsilon-\pi$, the pentagon with grey $CE$ fails the condition $\delta<\epsilon$.

The third of Figure \ref{geom_proof3} is the case $D,E$ are outside $\triangle ABC$, and $X$ is in both $BD$ and $CE$. By $\beta>\gamma$ and the isosceles triangle $\triangle ABC$, we get $\angle XBC>\angle XCB$. This implies $BX<CX$. Then by $BD=a=CE$, we get $DX>EX$. This implies $2\pi-\delta=\angle XDE<\angle XED=2\pi-\epsilon$. Therefore $\delta>\epsilon$ for the pentagon with black $DE$, and $\delta-\pi>\epsilon-\pi$ for the pentagon with grey $DE$. Therefore both pentagons fail the condition $\delta<\epsilon$. 

The fourth of Figure \ref{geom_proof3} is the case $D$ is outside $\triangle ABC$ and $E$ is inside $\triangle ABC$. In this case, the black $DE$ does not intersect the other four edges, and the pentagon is simple. We also have $\delta+\pi>\pi>\epsilon-\pi$. This means the pentagon with grey $DE$ fails the condition. 

The fifth of Figure \ref{geom_proof3} is the case both $D,E$ are inside $\triangle ABC$, and $X$ is not in $BD,CE$. This is the same as both $D,E$ being outside $\triangle A^*BC$, and we may use the exchange $(A,B)\leftrightarrow (A^*,C)$ to convert the situation to the first three pictures of Figure \ref{geom_proof3}. The conversion changes $\beta,\gamma,\delta,\epsilon$ to $\pi-\gamma,\pi-\beta,2\pi-\epsilon,2\pi-\delta$. Then the argument for the case $D,E$ are outside $\triangle ABC$ can be converted to the argument for the current case.
\end{proof}

\section{Symmetric Pentagon}
\label{symmetric_tiling}

As explained in Section \ref{companion}, we may assume $f\ge 16$. This implies not all vertices have degree $3$. 

We will mention $f\ge 16$ and the parity lemma the first time we use them. Then we will omit mentioning them subsequently. 

By Lemma \ref{geometry11}, we know $\beta=\gamma$ if and only if $\delta=\epsilon$. This means the pentagon is symmetric, with angles $\alpha,\beta,\beta,\delta,\delta$. 

\begin{proposition}\label{symmetric}
Tilings of the sphere by congruent symmetric almost equilateral pentagons are the first earth map tiling {\rm $E_{\pentagon}1$} in Figure \ref{emt}, and the two flip modifications {\rm $F_1E_{\pentagon}1,F_2E_{\pentagon}1$} in Figure \ref{emt1flip}.
\end{proposition}

\begin{proof}
The angle $\delta$ appears twice in the pentagon. By Lemma \ref{ndegree3}, this implies $\delta$ appears in a degree $3$ vertex. Then (by the parity lemma) one of $\alpha\delta^2,\beta\delta^2$ is a vertex. 

If $\alpha=\beta$, then the angle sum of $\alpha\delta^2$ and the angle sum for pentagon \eqref{psum} imply
\[
\alpha=(\tfrac{1}{2}+\tfrac{2}{f})\pi,\;
\delta=(\tfrac{3}{4}-\tfrac{1}{f})\pi.
\]
We have $\alpha,\delta>\frac{1}{2}\pi$ (by $f\ge 16$). This implies (by angle sums of vertices) all vertices have degree $3$, a contradiction. 

Therefore $\alpha\ne\beta$. Then $\alpha$ appears once in the pentagon. By applying Lemma \ref{degree3} to $\alpha$, we know there is a degree $3$ vertex without $\alpha$. This means one of $\beta^3,\beta\delta^2$ is a vertex. 

If $\alpha\delta^2$ is a vertex, then by $\alpha\ne\beta$, we know $\beta\delta^2$ is not a vertex. Therefore $\alpha\delta^2,\beta^3$ are vertices. The angle sums of $\alpha\delta^2,\beta^3$ and the angle sum for pentagon imply $f=12$, a contradiction. 

Therefore $\alpha\delta^2$ is not a vertex. Since one of $\alpha\delta^2,\beta\delta^2$ is a vertex, we know $\beta\delta^2$ is a vertex. This is the only degree $3$ $b$-vertex. Since $\delta$ appears twice in the pentagon, by applying Lemma \ref{degree3} to $\delta$, we know there is a degree $3$ $\hat{b}$-vertex. Therefore one of $\alpha^3,\alpha^2\beta,\alpha\beta^2,\beta^3$ is a vertex.

\subsubsection*{Case. $\alpha^3$ is a vertex}

The angle sum of $\beta\delta^2,\alpha^3$ and the angle sum for pentagon imply
\[
\alpha=\tfrac{2}{3}\pi,\;
\beta=(\tfrac{1}{3}+\tfrac{4}{f})\pi,\;
\delta=(\tfrac{5}{6}-\tfrac{2}{f})\pi.
\]
We have $\beta<\alpha,\delta$. 

By $\beta\delta^2$, and $\alpha>\beta$, we know $\beta\delta^2$ is the only $b$-vertex. This implies $\alpha\delta\cdots,\delta\thin\delta\cdots$ are not vertices. Therefore the AAD of $\thin\beta\thin\beta\thin$ is $\thin^{\delta}\beta^{\alpha}\thin^{\alpha}\beta^{\delta}\thin$. This implies no consecutive $\beta\beta\beta$.

By the values of $\alpha,\beta$, we know $\alpha^3,\alpha\beta^3,\beta^4,\beta^5$ are all the $\hat{b}$-vertices. By no consecutive $\beta\beta\beta$, we know $\alpha\beta^3,\beta^4,\beta^5$ are not vertices. Then $\beta\delta^2,\alpha^3$ are the only vertices. Then all vertices have degree $3$, a contradiction. 
	
\subsubsection*{Case. $\beta^3$ is a vertex}

The angle sums of $\beta\delta^2,\beta^3$ and the angle sum for pentagon imply
\[
\alpha=(\tfrac{1}{3}+\tfrac{4}{f})\pi,\;
\beta=\delta=\tfrac{2}{3}\pi.
\]
By $\alpha<R(\delta^2)=\beta<2\alpha<2\delta$, we know $\beta\delta^2$ is the only $b$-vertex. This implies  $\delta\thin\delta\cdots$ is not a vertex. Moreover, by the values of $\alpha,\beta$, we know $\beta^3,\alpha^3\beta,\alpha^4,\alpha^5$ are all the $\hat{b}$-vertices. Then we get $\beta^2\cdots=\beta^3$.

The AAD $\thin^{\beta}\alpha^{\beta}\thin^{\beta}\alpha^{\beta}\thin$ of $\thin\alpha\thin\alpha\thin$ implies a vertex $\thin^{\delta}\beta^{\alpha}\thin^{\alpha}\beta^{\delta}\thin\cdots=\beta^3=\thin^{\alpha}\beta^{\delta}\thin^{\alpha}\beta^{\delta}\thin^{\delta}\beta^{\alpha}\thin$, contradicting no $\delta\thin\delta\cdots$. Therefore $\alpha\thin\alpha\cdots$ is not a vertex. This implies $\alpha^3\beta,\alpha^4,\alpha^5$ are not vertices. Then $\alpha\cdots$ is not a vertex, a contradiction. 
	
\subsubsection*{Case. $\alpha^2\beta$ is a vertex}

The angle sums of $\beta\delta^2,\alpha^2\beta$ and the angle sum for pentagon imply
\[
\alpha=\delta=(1-\tfrac{4}{f})\pi,\;
\beta=\tfrac{8}{f}\pi.
\]
By $R(\alpha^2)=R(\delta^2)=\beta<\alpha,\delta$, we get $\alpha^2\cdots=\alpha^2\beta$, and $\beta\delta^2$ is the only $b$-vertex. This implies $\alpha\delta\cdots,\delta\thin\delta\cdots$ are not vertices. Therefore the AAD of $\thin\beta\thin\beta\thin$ is $\thin^{\delta}\beta^{\alpha}\thin^{\alpha}\beta^{\delta}\thin$.  This implies no consecutive $\beta\beta\beta$.

By the values of $\alpha,\beta$, we know $\alpha^2\beta,\alpha\beta^k(k\ge 3),\beta^k(k\ge 4)$ are all the $\hat{b}$-vertices. By no $\beta\beta\beta$, we know $\alpha\beta^k,\beta^k$ are not vertices. Then $\beta\delta^2,\alpha^2\beta$ are the only vertices. Then all vertices have degree $3$, a contradiction. 

\subsubsection*{Case. $\alpha\beta^2$ is a vertex}

The angle sums of $\beta\delta^2,\alpha\beta^2$ and the angle sum for pentagon imply
\[
\alpha=\tfrac{8}{f}\pi,\;
\beta=(1-\tfrac{4}{f})\pi,\;
\delta=(\tfrac{1}{2}+\tfrac{2}{f})\pi.
\]
We have $\alpha<\delta<\beta<2\delta$. By $\beta\delta^2$, and $\beta<2\delta$, we get $\delta^2\cdots=\beta\delta^2,\alpha^k\delta^2$. These are all the $b$-vertices. By the values of $\alpha,\beta$, we further get all the $\hat{b}$-vertices, and then get all vertices
\[
\text{AVC}
=\{\alpha\beta^2,\beta\delta^2,\alpha^k,\alpha^k\beta,\alpha^k\delta^2\}.
\]

For a tiling by congruent symmetric pentagons, a tile determines (the edges and angles of) its companion tile.

A consecutive pair $\thin\alpha\thin\alpha\thin$ determines $T_1,T_2$ in the first of Figure \ref{symmetricA}. Then $T_1,T_2$ share $\beta^2\cdots=\alpha\beta^2$. This determines $T_3$. Then $T_1,T_2,T_3$ determine their companions $T_4,T_5,T_6$. 

If $\alpha^k$ is a vertex, then we may apply the argument to all the pairs $\thin\alpha\thin\alpha\thin$ in $\alpha^k=\thin\alpha\thin\cdots\thin\alpha\thin$, and get the earth map tiling with $k$ timezones. Here  each timezone consists of $T_2,T_3,T_5,T_6$ in the first of Figure \ref{symmetricA}. This is the first earth map tiling {\rm $E_{\pentagon}1$} in Figure \ref{emt}.

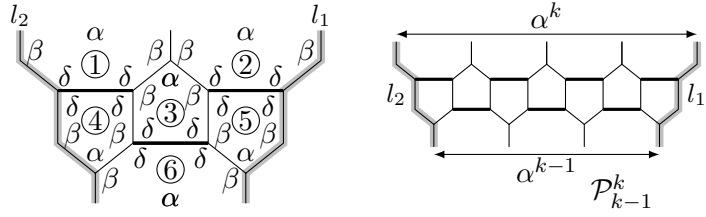
\begin{figure}[htp]
\centering
\begin{tikzpicture}[>=latex]

\foreach \a in {-1,1}
\draw[gray!50, xscale=\a, line width=3]
	(2,1.15) -- (2,0.75) -- (1.5,0.35) -- (1.5,-0.35) -- (1,-0.75) -- (1,-1.15);

\foreach \a in {-1,1}
{
\begin{scope}[xscale=\a]

\draw
	(0,1.15) -- (0,0.75) -- (0.5,0.35) -- (0.5,-0.35) -- (1,-0.75) -- (1,-1.15) 
	(2,1.15) -- (2,0.75) -- (1.5,0.35) -- (1.5,-0.35) -- (1,-0.75);
	
\draw[line width=1.2]
	(0.5,-0.35) -- (-0.5,-0.35)
	(0.5,0.35) -- (1.5,0.35);	

\node at (1.8,0.85) {\small $\beta$}; 
\node at (0.2,0.85) {\small $\beta$};
\node at (1.4,0.55) {\small $\delta$};
\node at (0.6,0.55) {\small $\delta$};
\node at (1,1.1) {\small $\alpha$};

\node at (0.3,0.25) {\small $\beta$}; 
\node at (0.3,-0.15) {\small $\delta$};	
\node at (0,0.5) {\small $\alpha$};	

\node at (1.3,-0.25) {\small $\beta$};
\node at (0.7,-0.25) {\small $\beta$};	
\node at (1.3,0.15) {\small $\delta$}; 
\node at (0.7,0.15) {\small $\delta$};
\node at (1,-0.5) {\small $\alpha$};

\node at (0.8,-0.85) {\small $\beta$}; 
\node at (0.4,-0.55) {\small $\delta$};
\node at (0,-1.1) {\small $\alpha$};

\end{scope}
}

\node at (-2,1.35) {\footnotesize $l_2$}; 
\node at (2,1.35) {\footnotesize $l_1$};

\node[inner sep=0.5,draw,shape=circle] at (-1,0.7) {\small 1};
\node[inner sep=0.5,draw,shape=circle] at (1,0.7) {\small 2};
\node[inner sep=0.5,draw,shape=circle] at (0,0.05) {\small 3};
\node[inner sep=0.5,draw,shape=circle] at (-1,-0.05) {\small 4};
\node[inner sep=0.5,draw,shape=circle] at (1,-0.05) {\small 5};
\node[inner sep=0.5,draw,shape=circle] at (0,-0.7) {\small 6};


\begin{scope}[shift={(3.5cm,0.3cm)}]

\draw[gray!50, line width=3]
	(-0.5,0.7) -- (-0.5,0.4) -- (-0.25,0.2) -- (-0.25,-0.2) -- (0,-0.4) -- (0,-0.7)
	(3.5,0.7) -- (3.5,0.4) -- (3.25,0.2) -- (3.25,-0.2) -- (3,-0.4) -- (3,-0.7);

\node at (-0.5,0) {\footnotesize $l_2$}; 
\node at (3.5,0) {\footnotesize $l_1$};
	
\foreach \a in {0,...,3}
\draw[xshift=\a cm]
	(-0.5,0.7) -- (-0.5,0.4) -- (-0.25,0.2) -- (0.25,0.2) -- (0.5,0.4) -- (0.5,0.7)
	(-0.25,0.2) -- (-0.25,-0.2) -- (0,-0.4) -- (0.25,-0.2) -- (0.25,0.2)
	(0,-0.4) -- (0,-0.7);

\foreach \a in {0,...,3}
\draw[line width=1.2, xshift=\a cm]	
	(-0.25,0.2) -- (0.25,0.2);

\foreach \a in {0,...,2}
\draw[line width=1.2, xshift=\a cm]
	(0.25,-0.2) -- (0.75,-0.2);

\draw[<->]
	(-0.5,0.8) -- ++(4,0);

\node at (1.5,1.05) {\small $\alpha^k$};
	
\draw[<->]
	(0,-0.8) -- ++(3,0);

\node at (1.5,-1) {\small $\alpha^{k-1}$};

\node at (2.5,-1.3) {\small ${\mc P}^k_{k-1}$};
	
\end{scope}

\end{tikzpicture}
\caption{Proposition \ref{symmetric}: Earth map tiling, and ${\mc P}^k_{k-1}$.}
\label{symmetricA}
\end{figure}

If there is consecutive $\alpha^k=\thin\alpha\thin\cdots\thin\alpha\thin$ at a vertex, then we may apply the argument to all the pairs $\thin\alpha\thin\alpha\thin$ in $\alpha^k$, and get a {\em partial earth map tiling} ${\mc P}^k_{k-1}$, with $\alpha^k,\alpha^{k-1}$ at the two ends. See the second of Figure \ref{symmetricA}. The boundary of ${\mc P}^k_{k-1}$ is given by the shaded edges $l_1,l_2$. 

Next, we assume $\alpha^k$ is not a vertex. Then we get updated list of vertices
\[
\text{AVC}
=\{
\alpha\beta^2,
\beta\delta^2,
\alpha^{q+1}\beta,
\alpha^q\delta^2\},\quad
f=8q+4. 
\]
In constructing the tiling, we may need to distinguish between two $\beta$ (or two $\delta$) in the same tile. We indicate the ditinction by $\beta,\dot{\beta}$ (or $\delta,\dot{\delta}$). 

Suppose $\alpha^q\delta^2$ is a vertex. The vertex is $\thick\delta\thin\alpha\thin\cdots\thin\alpha\thin\delta\thick$. Applying the argument for $\thin\alpha\thin\alpha\thin$ to all the pairs in the $\alpha^q=\thin\alpha\thin\cdots\thin\alpha\thin$ part of the vertex, we get a partial map tiling ${\mc P}^q_{q-1}$. This is the left of $l_1$ and the right of $l_2$ in the first of Figure \ref{symmetricB}. We also indicate the angles along the boundary of ${\mc P}^q_{q-1}$. By the first of Figure \ref{symmetricA}, these are $\beta,\delta\thick\delta,\beta,\alpha\thin\beta$ along the left of $l_1$ and along the right of $l_2$.

In the first of Figure \ref{symmetricB}, we also get $T_1,T_2$ containing the $\thin\delta\thick\delta\thin$ part of $\alpha^q\delta^2$. Then $\beta_1\cdots=\beta^2\cdots=\alpha\beta^2$ (note that the two $\beta$ in $T_1$ are $\beta_1$ and $\dot{\beta}_1$). This determines $T_3$. Then $T_3$ determines its companion $T_4$. By the same argument, we determine $T_5,T_6$. By adding $T_3,T_4,T_5,T_6$ to the partial earth map tiling ${\mc P}^q_{q-1}$ on the left of $l_1$ and the right of $l_2$, we actually get a partial earth map tiling ${\mc P}^{q+1}_q$.

Both $\alpha_1\cdots$ and $\alpha_2\cdots$ are $\alpha\beta\cdots=\alpha\beta^2,\alpha^{q+1}\beta$. Due to the symmetry of the constructed part of the tiling so far, we may independently assume $\alpha_1\cdots=\alpha\beta^2$ and $\alpha_2\cdots=\alpha^{q+1}\beta$, and then make independent argument. 

Suppose $\alpha_1\cdots=\alpha\beta^2$. Then by two $\beta$ not adjacent in a tile, this implies $\dot{\delta}_3\dot{\delta}_4\cdots=\alpha^q\delta^2$, and $\dot{\beta}_1\cdots=\beta\delta\cdots=\beta\delta^2$. This determines $T_7$.

Suppose $\alpha_2\cdots=\alpha^{q+1}\beta$. Then the $\alpha^{q+1}$ part of the vertex induces a partial earth map tiling ${\mc P}^{q+1}_q$, which consists of $P,T_1,T_2,T_7,T_8$, and is the complement of the earlier ${\mc P}^{q+1}_q$. In particular, this also determines the same $T_7$. 

We find either $\alpha_1\cdots=\alpha\beta^2$ or $\alpha_2\cdots=\alpha^{q+1}\beta$ implies the same $T_7$ as indicated. Since $T_7$ is not symmetric with respect to the initial $\alpha^q\delta^2$, we conclude one of $\alpha_1\cdots,\alpha_2\cdots$ is $\alpha\beta^2$, and the other is $\alpha^{q+1}\beta$. By the symmetry of the initial picture, we may assume $\alpha_2\cdots=\alpha^{q+1}\beta$. As argued earlier, this implies the tiling is obtained by glueing two ${\mc P}^{q+1}_q$ together. The tiling is the first flip modification $F_1E_{\pentagon}1$ in  Figure \ref{emt1flip}.

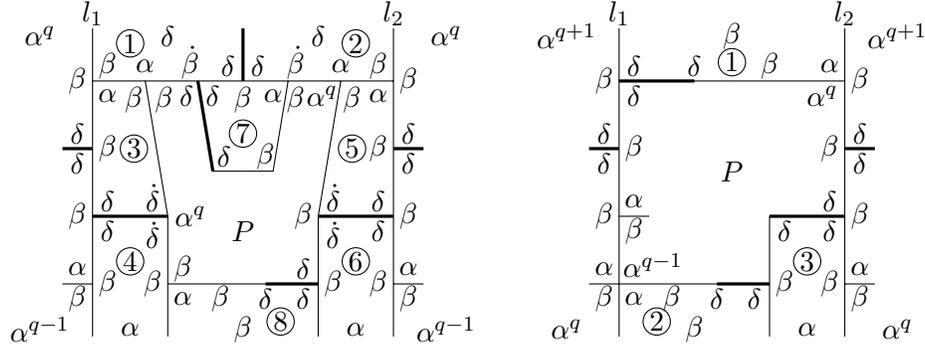
\begin{figure}[htp]
\centering
\begin{tikzpicture}[>=latex,scale=1]

	
\foreach \a in {1,-1}
{
\begin{scope}[xscale=\a]

\draw
	(2,2.5) -- (2,-1.6)
	(0,1.8) -- (2,1.8)
	(2,-0.9) -- ++(0.4,0)
	(1.3,1.8) -- (1,0) -- (1,-0.9) -- (1,-1.6)
	(0.6,1.8) -- (0.4,0.6) -- (0,0.6)
	(0,-0.9) -- (1,-0.9);

\draw[line width=1.2] 
	(2,0.9) -- ++(0.4,0)
	(2,0) -- (1,0);

\node at (1.3,2) {\small $\alpha$};
\node at (0.7,2.05) {\small $\dot{\beta}$};
\node at (1.8,2) {\small $\beta$};
\node at (0.2,2) {\small $\delta$};
\node at (1,2.4) {\small $\delta$};
	
\node at (2.7,2.4) {\small $\alpha^q$};
\node at (2.2,1.8) {\small $\beta$};
\node at (2.2,1.1) {\small $\delta$};
\node at (2.2,0.7) {\small $\delta$};
\node at (2.2,0) {\small $\beta$};
\node at (2.2,-0.7) {\small $\alpha$};
\node at (2.2,-1.1) {\small $\beta$};
\node at (2.7,-1.5) {\small $\alpha^{q-1}$};

\node at (1.8,1.6) {\small $\alpha$};
\node at (1.8,0.9) {\small $\beta$};
\node at (1.45,1.55) {\small $\beta$};
\node at (1.8,0.2) {\small $\delta$};
\node at (1.2,0.25) {\small $\dot{\delta}$};

\node at (1.5,-1.5) {\small $\alpha$};
\node at (1.8,-0.9) {\small $\beta$};
\node at (1.2,-0.9) {\small $\beta$};
\node at (1.8,-0.2) {\small $\delta$};
\node at (1.2,-0.25) {\small $\dot{\delta}$};

\end{scope}
}

\draw[line width=1.2] 
	(0,1.8) -- (0,2.5)
	(-0.6,1.8) -- (-0.4,0.6)
	(0.3,-0.9) -- (1,-0.9);	

\node at (0.4,1.6) {\small $\alpha$};
\node at (0,1.55) {\small $\beta$};
\node at (0.3,0.8) {\small $\beta$};
\node at (-0.4,1.6) {\small $\delta$};
\node at (-0.25,0.8) {\small $\delta$};

\node at (-0.8,-1.1) {\small $\alpha$};
\node at (0,-1.5) {\small $\beta$};
\node at (-0.3,-1.1) {\small $\beta$};
\node at (0.8,-1.1) {\small $\delta$};
\node at (0.3,-1.1) {\small $\delta$};

\node at (1.05,1.6) {\small $\alpha^q$};
\node at (0.7,1.55) {\small $\beta$};
\node at (0.8,0) {\small $\beta$};
\node at (0.8,-0.7) {\small $\delta$};

\node at (-1.05,1.55) {\small $\beta$};
\node at (-0.75,1.6) {\small $\delta$};
\node at (-0.7,0) {\small $\alpha^q$};
\node at (-0.8,-0.7) {\small $\beta$};

\node at (2,2.7) {\small $l_2$};
\node at (-2,2.7) {\small $l_1$};
\node at (0,-0.2) {\small $P$};

\node[inner sep=0.5,draw,shape=circle] at (-1.5,2.3) {\small 1};
\node[inner sep=0.5,draw,shape=circle] at (1.5,2.3) {\small 2};
\node[inner sep=0.5,draw,shape=circle] at (-1.45,0.9) {\small 3};
\node[inner sep=0.5,draw,shape=circle] at (-1.5,-0.6) {\small 4};
\node[inner sep=0.5,draw,shape=circle] at (1.45,0.9) {\small 5};
\node[inner sep=0.5,draw,shape=circle] at (1.5,-0.6) {\small 6};
\node[inner sep=0.5,draw,shape=circle] at (0,1.1) {\small 7};
\node[inner sep=0.5,draw,shape=circle] at (0.5,-1.4) {\small 8};


\begin{scope}[xshift=6.5cm]
	
\foreach \a in {1,-1}
{
\begin{scope}[xscale=\a]

\draw
	(1.5,2.5) -- (1.5,-1.6)
	(0,1.8) -- (1.5,1.8)
	(1.5,-0.9) -- ++(0.4,0);

\draw[line width=1.2] 
	(1.5,0.9) -- ++(0.4,0);
	
\node at (2.2,2.4) {\small $\alpha^{q+1}$};
\node at (1.7,1.8) {\small $\beta$};
\node at (1.7,1.1) {\small $\delta$};
\node at (1.7,0.7) {\small $\delta$};
\node at (1.7,0) {\small $\beta$};
\node at (1.7,-0.7) {\small $\alpha$};
\node at (1.7,-1.1) {\small $\beta$};
\node at (2.2,-1.5) {\small $\alpha^q$};

\end{scope}
}

\draw
	(-1.5,-0.9) -- (0.5,-0.9)
	(1.5,0) -- (0.5,0) -- (0.5,-1.6)
	(-1.5,0) -- ++(0.4,0);

\draw[line width=1.2] 
	(-1.5,1.8) -- (-0.5,1.8)
	(0.5,-0.9) -- (-0.2,-0.9)
	(0.5,0) -- (1.5,0);
	
\node at (1.5,2.7) {\small $l_2$};
\node at (-1.5,2.7) {\small $l_1$};
\node at (0,0.6) {\small $P$};

\node at (1.3,2) {\small $\alpha$};
\node at (0.5,2) {\small $\beta$};
\node at (0,2.4) {\small $\beta$};
\node at (-1.3,2) {\small $\delta$};
\node at (-0.5,2) {\small $\delta$};

\node at (1.2,1.6) {\small $\alpha^q$};
\node at (1.3,0.9) {\small $\beta$};
\node at (1.3,0.2) {\small $\delta$};

\node at (-1.05,-0.65) {\small $\alpha^{q-1}$};
\node at (-1.3,0.2) {\small $\alpha$};
\node at (-1.3,-0.2) {\small $\beta$};
\node at (-1.3,0.9) {\small $\beta$};
\node at (-1.3,1.6) {\small $\delta$};

\node at (-1.3,-1.1) {\small $\alpha$};
\node at (-0.8,-1.1) {\small $\beta$};
\node at (-0.5,-1.5) {\small $\beta$};
\node at (-0.2,-1.1) {\small $\delta$};
\node at (0.3,-1.1) {\small $\delta$};

\node at (1,-1.5) {\small $\alpha$};
\node at (0.7,-0.9) {\small $\beta$};
\node at (1.3,-0.9) {\small $\beta$};
\node at (1.3,-0.2) {\small $\delta$};
\node at (0.7,-0.2) {\small $\delta$};

\node[inner sep=0.5,draw,shape=circle] at (0,2.05) {\small 1};
\node[inner sep=0.5,draw,shape=circle] at (-1,-1.4) {\small 2};
\node[inner sep=0.5,draw,shape=circle] at (1,-0.6) {\small 3};

\end{scope}

\end{tikzpicture}
\caption{Proposition \ref{symmetric}: Tilings for $\alpha^q\delta^2,\alpha^{q+1}\beta$.}
\label{symmetricB}
\end{figure}

Finally, we assume $\alpha^q\delta^2$ is also not a vertex. Then we get further updated list of vertices
\[
\text{AVC}=\{
\alpha\beta^2,
\beta\delta^2,
\alpha^{q+1}\beta\}.
\]
The $\alpha^{q+1}$ part of the vertex $\alpha^{q+1}\beta$ induces a partial earth map tiling ${\mc P}^{q+1}_q$, which is the left of $l_1$ and the right of $l_2$ in the second of Figure \ref{symmetricB}. We also get $T_1$ containing $\beta$ at the vertex. 

The vertex $\delta\thick\delta\cdots$ on the right of $l_2$ is $\beta\delta^2$. Since two $\beta$ in a tile are not adjacent, this implies $\alpha_1\cdots=\alpha\beta\cdots$ is not $\alpha\beta^2$, and must be $\alpha^{q+1}\beta$. Then the $\alpha^{q+1}$ part of the vertex induces a partial earth map tiling ${\mc P}^{q+1}_q$, which consists of $P,T_1,T_2,T_3$. The whole tiling is then obtained by gluing two ${\mc P}^{q+1}_q$ together. The tiling is the second flip modification $F_2E_{\pentagon}1$ in Figure \ref{emt1flip}.

\medskip

\noindent{\em Geometry of Pentagon}

\medskip

We need to justify the existence of the symmetric pentagon in the earth map tiling. In Figure \ref{symmetricC}, we divide the symmetric pentagon into two congruent quadrilaterals. Then the existence of the pentagon is the same as the existence of the quadrilateral $\square ABDX$, with angles
\[
\bar{\alpha}=\tfrac{1}{2}\alpha=\tfrac{4}{f}\pi,\;
\beta=(1-\tfrac{4}{f})\pi,\;
\delta=(\tfrac{1}{2}+\tfrac{2}{f})\pi,
\]
and the angle of the quadrilateral at $X$ is $\frac{1}{2}\pi$.

The triangle $\triangle ABD$ is isosceles, with base angle $\theta$. For any $\frac{2}{f}\pi<\theta\le \frac{4}{f}\pi=\bar{\alpha}$, we have $\beta+2\theta>\pi$. Therefore we may draw the isosceles triangle $\triangle ABD$ with top angle $\beta$ and base angle $\theta$. Then we draw $AX$ satisfying $\angle BAX=\bar{\alpha}=\frac{4}{f}\pi$, and draw $DX$ satisfying $\angle BDX=\delta=(\tfrac{1}{2}+\tfrac{2}{f})\pi$. By $0\le \angle DAX=\bar{\alpha}-\theta<\frac{2}{f}\pi$ and $(\tfrac{1}{2}-\tfrac{2}{f})\pi\le \angle ADX=\delta-\theta<\frac{1}{2}\pi$, we get a triangle $\triangle ADX$ inside one side of the great circle $\bigcirc AD$. The triangle $\triangle ABD$ is on the other side of $\bigcirc AD$. This implies the quadrilateral $\square ABDX$ is simple.

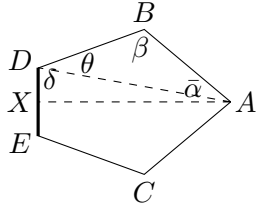
\begin{figure}[htp]
\centering
\begin{tikzpicture}[scale=1]

\foreach \a in {1,-1}
\draw[yscale=\a]
	(3,0) -- ++(140:1.5) -- ++(200:1.5) -- ++(0,-0.455);

\draw[dashed]
	(3,0) -- (0.44,0)
	(3,0) -- (0.44,0.455);
	
\draw[line width=1.2]
	(0.44,-0.455) -- (0.44,0.455);

\node at (3.2,0) {\small $A$};
\node at (1.85,1.2) {\small $B$};
\node at (1.85,-1.2) {\small $C$};
\node at (0.2,0.55) {\small $D$};
\node at (0.2,-0.55) {\small $E$};
\node at (0.2,0) {\small $X$};

\node at (1.8,0.7) {\small $\beta$};
\node at (2.5,0.17) {\small $\bar{\alpha}$};
\node at (1.1,0.52) {\small $\theta$};
\node at (0.6,0.3) {\small $\delta$};

\end{tikzpicture}
\caption{Proposition \ref{symmetric}: Existence of symmetric pentagon.}
\label{symmetricC}
\end{figure}

For the existence of the quadrilateral $\square ABDX$, it remains to justify that the angle at $X$ can be $\frac{1}{2}\pi$. This is equivalent to that the area of the quadrilateral is $\bar{\alpha}+\beta+\delta+\frac{1}{2}\pi-2\pi=\frac{2}{f}\pi$. The area is a continuous function $A(\theta)$ of $\theta$. As $\theta$ approaches $\frac{2}{f}\pi$, we see $\beta+2\theta$ approaches $\pi$. Therefore the triangle $\triangle ABD$ approaches to a point, and $\lim_{\theta\to\frac{2}{f}\pi^+} A(\theta)=0$. On the other hand, at $\theta=\frac{4}{f}\pi=\bar{\alpha}$, the quadrilateral is reduced to the triangle $\triangle ABD$. Therefore $A(\frac{4}{f}\pi)=\beta+2\cdot \frac{4}{f}\pi-\pi=\frac{4}{f}\pi$. By $0<\frac{2}{f}\pi<\bar{\alpha}$ and the intermediate value theorem, there is $\theta$ satisfying $\frac{2}{f}\pi<\theta\le \frac{4}{f}\pi$ and $A(\theta)=\frac{2}{f}\pi$. 

We may actually calculate the edge length, and get
\[
\cos a
=1+\tfrac{\sqrt{5}-3}{4}\sec^2\tfrac{2}{f}\pi.
\]
The details can be found in \cite{cl}. By $\tfrac{2}{f}\pi<\tfrac{1}{4}\pi$, we get $\cos a>0$. Therefore $a<\frac{1}{2}\pi$. This implies $\square ABDX$ lies in one side of the great circle $\bigcirc AX$. Therefore $\square ACEX$ lies in the other side of the great circle $\bigcirc AX$. Since the quadrilateral is simple, this implies the pentagon is simple.
\end{proof}

\begin{proposition}\label{b2d_c2e_abc}
If $\beta\delta^2,\gamma\epsilon^2,\alpha\beta\gamma$ are vertices in a tiling, then the pentagon is symmetric.
\end{proposition}

\begin{proof}
The angle sums of $\beta\delta^2,\gamma\epsilon^2,\alpha\beta\gamma$ imply
\[
\alpha=2(\delta+\epsilon-\pi),\;
\beta=2\pi-2\delta,\;
\gamma=2\pi-2\epsilon.
\]
Substituting into \eqref{coolsaet_eq1}, we get
\[
(\sin\delta-\sin\epsilon)^2
(\sin^2\delta+\sin^2\epsilon+3\sin\delta\sin\epsilon)=0.
\]

By $\beta\delta^2,\gamma\epsilon^2$, we get $\delta,\epsilon<\pi$. Therefore $\sin\delta,\sin\epsilon>0$. This implies the equation is the same as $\sin\delta=\sin\epsilon$. By the angle sum of $\alpha\beta\gamma$ and the angle sum for pentagon, we get $\pi<\delta+\epsilon=(1+\frac{4}{f})\pi<2\pi$. Therefore $\sin\delta=\sin\epsilon$ implies $\delta=\epsilon$. By Lemma \ref{geometry11}, this implies $\beta=\gamma$, and the pentagon is symmetric. 
\end{proof}

If we do not use \eqref{coolsaet_eq1}, then we may start with vertices $\beta\delta^2,\gamma\epsilon^2,\alpha\beta\gamma$, and obtain a tiling that is combinatorially possible, but geometrically impossible.

Specifically, the angle sums of $\beta\delta^2,\gamma\epsilon^2,\alpha\beta\gamma$ and the angle sum for pentagon imply
\[
\alpha=\tfrac{8}{f}\pi,\;
\beta+\gamma=(2-\tfrac{8}{f})\pi,\;
\delta+\epsilon=(1+\tfrac{4}{f})\pi.
\]
If $\beta>\gamma$, then we may derive the list of vertices (if $\beta<\gamma$, then we get another AVC by exchanging $(\beta,\delta)$ with $(\gamma,\epsilon)$)
\begin{equation}\label{b2d_c2eAVC}
\text{AVC}
=\{\alpha\beta\gamma,\beta\delta^2,\gamma\epsilon^2,\alpha^k\beta^l,\alpha^k\gamma,\alpha^k\delta^2\}.
\end{equation}
Then we may derive tilings purely based on the AVC.

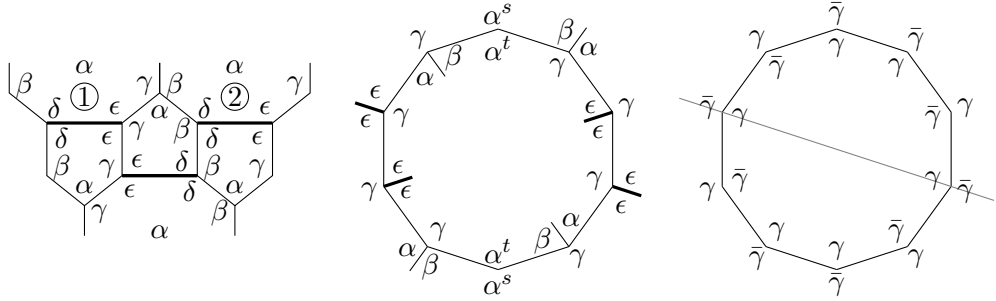
\begin{figure}[htp]
\centering
\begin{tikzpicture}[>=latex]

\foreach \a in {-1,1}
{
\begin{scope}[xscale=\a]

\draw
	(0,1.15) -- (0,0.75) -- (0.5,0.35) -- (0.5,-0.35) -- (1,-0.75) -- (1,-1.15) 
	(2,1.15) -- (2,0.75) -- (1.5,0.35) -- (1.5,-0.35) -- (1,-0.75);
	
\draw[line width=1.2]
	(0.5,-0.35) -- (-0.5,-0.35)
	(0.5,0.35) -- (1.5,0.35);	

\end{scope}
}

\foreach \a in {0,-1}
{
\begin{scope}[xshift=2*\a cm]

\node at (1,-0.5) {\small $\alpha$};
\node at (1.3,-0.25) {\small $\gamma$};
\node at (0.7,-0.25) {\small $\beta$};	
\node at (1.3,0.15) {\small $\epsilon$}; 
\node at (0.7,0.15) {\small $\delta$};

\node at (1.8,0.85) {\small $\gamma$}; 
\node at (0.2,0.85) {\small $\beta$};
\node at (1.4,0.55) {\small $\epsilon$};
\node at (0.6,0.55) {\small $\delta$};
\node at (1,1.1) {\small $\alpha$};

\end{scope}
}

\node at (0,0.5) {\small $\alpha$};	
\node at (-0.3,0.25) {\small $\gamma$}; 
\node at (0.3,0.25) {\small $\beta$}; 
\node at (-0.3,-0.15) {\small $\epsilon$};	
\node at (0.3,-0.15) {\small $\delta$};	

\node at (0,-1.1) {\small $\alpha$};
\node at (0.8,-0.85) {\small $\beta$}; 
\node at (-0.8,-0.85) {\small $\gamma$}; 
\node at (0.4,-0.55) {\small $\delta$};
\node at (-0.4,-0.55) {\small $\epsilon$};

\node[inner sep=0.5,draw,shape=circle] at (-1,0.7) {\small 1};
\node[inner sep=0.5,draw,shape=circle] at (1,0.7) {\small 2};

\begin{scope}[xshift=4.5cm]

\foreach \a in {0,...,9}
\foreach \b in {0,1}
\draw[xshift=4.5*\b cm, rotate=36*\a]
	(-18:1.6) -- (18:1.6);

\foreach \a in {0,1}
{
\begin{scope}[rotate=180*\a]

\draw
	(126:1.6) -- (126:1.2)
	(54:1.6) -- (54:2);
	
\draw[line width=1.2]
	(-18:1.6) -- (-18:2)
	(198:1.6) -- (198:1.2);

\node at (90:1.35) {\small $\alpha^t$};		
\node at (116:1.35) {\small $\beta$};
\node at (136:1.35) {\small $\alpha$};
\node at (162:1.35) {\small $\gamma$};
\node at (190:1.35) {\small $\epsilon$};
\node at (205:1.35) {\small $\epsilon$};
\node at (234:1.35) {\small $\gamma$};

\node at (60:1.8) {\small $\beta$};
\node at (48:1.8) {\small $\alpha$};
\node at (18:1.8) {\small $\gamma$};
\node at (-12:1.8) {\small $\epsilon$};
\node at (-25:1.8) {\small $\epsilon$};
\node at (-54:1.8) {\small $\gamma$};
\node at (90:1.8) {\small $\alpha^s$};

\end{scope}
}

\begin{scope}[xshift=4.5cm]

\draw[gray]
	(-18:2.2) -- (-18:-2.2);
	
\foreach \a in {0,1}
{

\foreach \b in {-18,54,90}
{
\node at (180*\a+\b:1.35) {\small $\gamma$};
\node at (180*\a+\b:1.8) {\small $\bar{\gamma}$};
}

\foreach \b in {18,-54}
{
\node at (180*\a+\b:1.35) {\small $\bar{\gamma}$};
\node at (180*\a+\b:1.8) {\small $\gamma$};
}

}

\end{scope}

\end{scope}

\end{tikzpicture}
\caption{Proposition \ref{b2d_c2e_abc}: Combinatorially possible tiling.}
\label{b2d_c2e_abcA}
\end{figure}

The AAD $\thin^{\beta}\alpha^{\gamma}\thin^{\beta}\alpha^{\gamma}\thin$ determines $T_1,T_2$ in the first of Figure \ref{b2d_c2e_abcA}. Then we can derive all the other tiles similar to the first of Figure \ref{symmetricA}. If $\alpha^k$ is a vertex, then we get an earth map tiling similar to $E_{\pentagon}1$ in Figure \ref{emt}. We note that $k=\frac{f}{4}$ in $\alpha^k$.

Next, suppose $\gamma=t\alpha$ for some integer $t$. In the second of Figure \ref{b2d_c2e_abcA}, we divide the earth map tiling into two parts, one with $t$ timezones, and the other with $s=\frac{f}{4}-t$ timezones. We indicate the angles along the boundary between the two parts. The third of Figure \ref{b2d_c2e_abcA} shows the angle values along the boundary between. Then we see that we may flip the inner part and still get a tiling. There is another obvious flip, but it gives the same result.

The flip modification can even be repeated. Let us call an $\alpha^t$-timezone to be $t$ consecutive timezones in the earth map tiling. We may pick $l$ disjoint $\alpha^t$-timezones in the earth map tiling, and flip each of them. Then we still get a tiling, with the vertex $\alpha^k$ in the earth map tiling changed to the vertex $\alpha^{k-lt}\gamma^l$. 

It turns out the earth map tiling and these flip modifications are all the tilings for the AVC \eqref{b2d_c2eAVC}. Of course, once we introduce geometrical considerations and find the pentagon must be symmetric, then we do not have so many flip modifications.

\section{$\alpha,\beta,\gamma$ Are Distinct}
\label{distinct}

We discuss tilings under the assumption that $\alpha,\beta,\gamma$ have distinct values. By $\beta\ne\gamma$ and Lemma \ref{geometry1}, we know $\delta\ne\epsilon$. Then by using the angle values and the $b$-edge, we may distinguish all five angles. In particular, each angle appears $f$ times in the tiling, and we may use all the counting results in Section \ref{count}. On the other hand, we remark that the assumption still allows some from $\alpha,\beta,\gamma$ has the same value as some from $\delta,\epsilon$. 

We may still assume $f\ge 16$, which implies not all vertices have degree $3$. We will omit mentioning $f\ge 16$ and parity lemma in the proofs. 

By applying Lemma \ref{degree3} and Lemma \ref{ndegree3} to $\delta,\epsilon$, we know a tiling has degree $3$ $b$-vertex, and has degree $3$ $\hat{b}$-vertex. The degree $3$ $\hat{b}$-vertices are $\alpha^3,\beta^3,\gamma^3,\alpha\beta\gamma,\alpha^2\beta,\alpha^2\gamma,\alpha\beta^2,\alpha\gamma^2,\beta^2\gamma,\beta\gamma^2$. 

We divide the classification according to the types of degree $3$ $b$-vertices. First, we consider twisted degree $3$ $b$-vertex. Up to the exchange symmetry $(\beta,\delta)\leftrightarrow(\gamma,\epsilon)$, this means either $\alpha\delta\epsilon$ or $\beta\delta\epsilon$ is a vertex. Sections \ref{ade} and \ref{bde} classify tilings for the two cases. Next, we consider matched degree $3$ $b$-vertices. These are $\theta\delta^2,\theta\epsilon^2$, where $\theta=\alpha,\beta,\gamma$. We further divide the classification into two matched degree $3$ $b$-vertices in Section \ref{2add}, and single matched degree $3$ $b$-vertex in Section \ref{singleb}.

\subsection{$\alpha\delta\epsilon$ and One Degree $3$ $\hat{b}$-Vertex}
\label{ade}

We assume $\alpha\delta\epsilon$ is a vertex. The vertex is symmetric with respect to the exchange $(\beta,\delta)\leftrightarrow(\gamma,\epsilon)$. The tiling also has a degree $3$ $\hat{b}$-vertex. Up to the exchange symmetry, this means one of $\alpha^3,\beta^3,\alpha\beta\gamma,\alpha^2\gamma,\alpha\beta^2,\beta^2\gamma$ is a vertex. 

\begin{proposition}\label{ade_abc}
There is no tiling, such that $\alpha,\beta,\gamma$ have distinct values, and $\alpha\delta\epsilon,\alpha\beta\gamma$ are vertices.
\end{proposition}

\begin{proof}
The angle sums of $\alpha\delta\epsilon,\alpha\beta\gamma$ and the angle sum for pentagon imply
\[
\alpha=(1-\tfrac{4}{f})\pi,\,
\beta+\gamma=\delta+\epsilon=(1+\tfrac{4}{f})\pi.
\] 
The pair of vertices $\alpha\delta\epsilon,\alpha\beta\gamma$ is symmetric with respect to the exchange $(\beta,\delta)\leftrightarrow(\gamma,\epsilon)$. By the exchange symmetry, and Lemma \ref{geometry1}, we may assume $\beta>\gamma$ and $\delta<\epsilon$. This implies $\beta,\epsilon>(\frac{1}{2}+\tfrac{2}{f})\pi$ and $\gamma,\delta<(\frac{1}{2}+\tfrac{2}{f})\pi$. In particular, we get $\alpha>\frac{2}{3}\pi>\gamma$ and $\beta+\epsilon>\pi$.

By $\alpha\delta\epsilon$ and Lemma \ref{square}, we know $\alpha^2\cdots$ is a $\hat{b}$-vertex. Then by $\alpha>\frac{2}{3}\pi$ and $2\alpha+\beta>\alpha+\beta+\gamma=2\pi$, we get $\alpha^2\cdots=\alpha^2\gamma^k$. By $\alpha\beta\gamma$ and $\alpha\ne\beta$, we know $k\ge 2$ in $\alpha^2\gamma^k$.

By $\delta+\epsilon>\pi$ and $\delta<\epsilon$, we know $R(\epsilon^2)$ has no $\delta,\epsilon$, and $R(\delta\epsilon)$ has no $\epsilon$. By $\alpha\delta\epsilon$ and $\delta<\epsilon$, we get $R(\epsilon^2)<R(\delta\epsilon)=\alpha<\beta+\gamma$. Then by $\beta>\gamma$, we get $\epsilon\cdots=\alpha\delta\epsilon,\beta\epsilon^2,\gamma^k\epsilon^2,\beta\delta^l\epsilon(l\ge 3),\gamma^k\delta^l\epsilon$. This implies $\alpha\epsilon\cdots=\alpha\delta\epsilon$, and $\epsilon\thin\epsilon\cdots$ is not a vertex. 

\subsubsection*{Case. $\alpha<\beta$} 

By $\alpha\delta\epsilon$, and $\alpha<\beta$, and $\delta<\epsilon$, we get $\epsilon\cdots=\alpha\delta\epsilon,\gamma^k\epsilon^2,\gamma^k\delta^l\epsilon$. This implies $\beta\epsilon\cdots$ is not a vertex. 

By $\alpha\delta\epsilon$, and $\alpha<\beta$, and Lemma \ref{square}, we know $\alpha\beta\cdots,\beta^2\cdots$ are $\hat{b}$-vertices. By $\alpha\beta\gamma$ and $\alpha<\beta$, we get $R(\beta^2)<R(\alpha\beta)=\gamma<\alpha,\beta$. This implies $\alpha\beta\cdots=\alpha\beta\gamma$, and $\beta^2\cdots$ is not a vertex. Combined with no $\beta\epsilon\cdots$, we get $\beta\cdots=\alpha\beta\gamma,\beta\gamma^k\delta^l$. 

Let $\theta=\alpha,\gamma$, and $\theta^{\rho}=\alpha^{\beta},\gamma^{\epsilon}$. By no $\beta^2\cdots,\beta\epsilon\cdots,\epsilon\thin\epsilon\cdots$, we know the AAD of ${}^{\rho}\thin\theta\thin\cdots\thin\theta\thin$ is ${}^{\rho}\thin\theta^{\rho}\thin\cdots\thin\theta^{\rho}\thin$. This implies the AAD of $\thick\delta\thin\theta\thin\theta\thin\cdots\thin\theta\thin$ is $\thick^{\epsilon}\delta^{\beta}\thin\theta^{\rho}\thin\theta^{\rho}\thin\cdots\thin\theta^{\rho}\thin$, and further implies a $\delta^2$-fan cannot have only $\theta$. Therefore a $\delta^2$-fan must have $\beta$. This implies $\beta\gamma^k\delta^l=\beta\gamma^k,\beta\gamma^k\delta^2$, and $\gamma^k\delta^l\epsilon=\gamma^k\delta\epsilon$. By $\alpha\delta\epsilon$ and $\alpha>\gamma$, we get $k\ge 2$ in $\gamma^k\delta\epsilon$.

The AAD $\thin^{\epsilon}\gamma^{\alpha}\thin^{\epsilon}\gamma^{\alpha}\thin$ implies a vertex $\thin^{\beta}\alpha^{\gamma}\thin^{\gamma}\epsilon^{\delta}\thick\cdots=\alpha\delta\epsilon=\thick^{\epsilon}\delta^{\beta}\thin^{\beta}\alpha^{\gamma}\thin^{\gamma}\epsilon^{\delta}\thick$, contradicting no $\beta^2\cdots$. The AAD $\thin^{\epsilon}\gamma^{\alpha}\thin^{\alpha}\gamma^{\epsilon}\thin$ implies a vertex $\thin^{\beta}\alpha^{\gamma}\thin^{\gamma}\alpha^{\beta}\thin\cdots=\alpha^2\gamma^k=\thin^{\gamma}\alpha^{\beta}\thin\gamma\thin\cdots\thin\gamma\thin^{\beta}\alpha^{\gamma}\thin$, contradicting the AAD ${}^{\beta}\thin^{\alpha}\gamma^{\epsilon}\thin\cdots\thin^{\alpha}\gamma^{\epsilon}\thin$ of ${}^{\beta}\thin\gamma\thin\cdots\thin\gamma\thin$. By no $\epsilon\thin\epsilon\cdots$, we do not have the AAD $\thin^{\alpha}\gamma^{\epsilon}\thin^{\epsilon}\gamma^{\alpha}\thin$. Therefore $\gamma\thin\gamma\cdots$ is not a vertex. Then $\beta\gamma^k$ and $\gamma^k\delta\epsilon$ are not vertices, and $\gamma^k\epsilon^2=\gamma\epsilon^2$. Therefore $\epsilon\cdots=\alpha\delta\epsilon,\gamma\epsilon^2$, and $\beta\cdots=\alpha\beta\gamma,\beta\gamma^k\delta^2$.

By $\alpha\beta\gamma$, and $\epsilon\cdots=\alpha\delta\epsilon,\gamma\epsilon^2$, and applying the counting lemma to $\alpha,\epsilon$, we know $\gamma\epsilon^2$ is a vertex. By $\gamma\epsilon^2$, and $\beta\cdots=\alpha\beta\gamma,\beta\gamma^k\delta^2$, and applying the counting lemma to $\beta,\gamma$, we know $\beta\gamma^k\delta^2=\beta\delta^2$, and is a vertex. By $\beta\delta^2$ and $\beta+\gamma=\delta+\epsilon=(1+\tfrac{4}{f})\pi$, we get $\gamma+2\epsilon=(1+\tfrac{12}{f})\pi<2\pi$. This contradicts the vertex $\gamma\epsilon^2$.

\subsubsection*{Case. $\alpha>\beta$}

By $\alpha\beta\gamma$, and $\alpha>\beta>\gamma$, and $\beta+\gamma>\pi$, we know $\alpha\beta\gamma,\beta^3,\alpha\gamma^k,\beta\gamma^k,\gamma^k$ are all the $\hat{b}$-vertices. In particular, we know $\alpha^2\cdots=\alpha^2\gamma^k$ is not a vertex. 

By no $\alpha^2\cdots,\epsilon\thin\epsilon\cdots$, we know the AAD of $\thin\gamma\thin\gamma\thin$ is $\thin^{\alpha}\gamma^{\epsilon}\thin^{\alpha}\gamma^{\epsilon}\thin$. The AAD implies a vertex $\thin^{\beta}\alpha^{\gamma}\thin^{\gamma}\epsilon^{\delta}\thick\cdots=\alpha\delta\epsilon=\thick^{\epsilon}\delta^{\beta}\thin^{\beta}\alpha^{\gamma}\thin^{\gamma}\epsilon^{\delta}\thick$. Therefore a vertex $\gamma\thin\gamma\cdots$ implies a vertex $\beta\thin\beta\cdots$. 

We divide the argument into further cases by considering $\epsilon^2\cdots=\beta\epsilon^2,\gamma^k\epsilon^2$.

\subsubsection*{Subcase. $\beta\epsilon^2$ is a vertex} 

By $\beta\epsilon^2$ and $\beta+\gamma=\delta+\epsilon=(1+\tfrac{4}{f})\pi$, we get $\gamma+2\delta=(1+\tfrac{12}{f})\pi>\beta+\gamma>\pi$. This implies $\beta<2\delta$. By $\beta+2\delta>\gamma+2\delta>\pi$, and $\delta+\epsilon>\pi$, we get the updated list $\epsilon\cdots=\alpha\delta\epsilon,\beta\epsilon^2,\gamma^k\epsilon^2(k\ge 2),\gamma^k\delta\epsilon(k\ge 2),\delta^l\epsilon$. Here by $\beta\epsilon^2$, we get $k\ge 2$ in $\gamma^k\epsilon^2$. 

Suppose $\beta^2\cdots$ is not a vertex. Then $\delta\thin\delta\cdots=\thick^{\epsilon}\delta^{\beta}\thin^{\beta}\delta^{\epsilon}\thick\cdots$ is not a vertex. This implies $\delta^l\epsilon$ is not a vertex. Moreover, by $\gamma\thin\gamma\cdots$ implying $\beta\thin\beta\cdots$, we also know $\gamma\thin\gamma\cdots$ is not a vertex. This implies $\alpha\gamma^k,\beta\gamma^k,\gamma^k,\gamma^k\delta\epsilon,\gamma^k\epsilon^2$ are not vertices. Therefore $\alpha\beta\gamma,\beta^3$ are all the $\hat{b}$-vertices, and $\epsilon\cdots=\alpha\delta\epsilon,\beta\epsilon^2$. In particular, $\delta\thin\epsilon\cdots$ is not a vertex.

The AAD $\thick^{\delta}\epsilon^{\gamma}\thin^{\alpha}\beta^{\delta}\thin^{\gamma}\epsilon^{\delta}\thick$ of $\beta\epsilon^2$ implies a vertex $\gamma\delta\cdots$. By no $\gamma\epsilon\cdots$, we know $\gamma\delta\cdots=\gamma\delta^k\cdots$, with no $\delta,\epsilon$ in the remainder. If $k\ge 4$, then by $\gamma+2\delta>\pi$, we get $R(\gamma\delta^k)<\gamma<\alpha,\beta$. This implies $\gamma\delta\cdots=\gamma\delta^k$, contradicting no $\delta\thin\delta\cdots$. Therefore $\gamma\delta\cdots=\gamma\delta^2\cdots$, with no $\delta,\epsilon$ in the remainder. By $\alpha\beta\gamma$ and $\beta<2\delta$, we get $R(\gamma\delta^2)<\alpha<\beta+\gamma$. Then by $\beta>\gamma$, we get $\gamma\delta\cdots=\beta\gamma\delta^2,\gamma^k\delta^2$. By no $\gamma\thin\gamma\cdots$, we get $k=1$ in $\gamma^k\delta^2$. However, by $\alpha\delta\epsilon$, and $\alpha>\gamma$, and $\delta<\epsilon$, we know $\gamma\delta^2$ is not a vertex. Therefore $\gamma\delta\cdots=\beta\gamma\delta^2$.

The angle sums of $\alpha\beta\gamma,\alpha\delta\epsilon,\beta\epsilon^2,\beta\gamma\delta^2$ and the angle sum for pentagon imply
\[
\alpha=(1-\tfrac{4}{f})\pi,\,
\beta=(1-\tfrac{12}{f})\pi,\,
\gamma=\tfrac{16}{f}\pi,\,
\delta=(\tfrac{1}{2}-\tfrac{2}{f})\pi,\,
\epsilon=(\tfrac{1}{2}+\tfrac{6}{f})\pi.
\]
By $\beta>\gamma$, we get $f>28$. By $\delta<\epsilon$ and $R(\alpha\gamma)=\beta<2\delta$, we know $\alpha\gamma\cdots$ is a $\hat{b}$-vertex. Then by all $\hat{b}$-vertices $\alpha\beta\gamma,\beta^3$, we get $\alpha\gamma\cdots=\alpha\beta\gamma$. 

The AAD $\thick^{\delta}\epsilon^{\gamma}\thin^{\alpha}\beta^{\delta}\thin^{\gamma}\epsilon^{\delta}\thick$ of $\beta\epsilon^2$ determines $T_1,T_2,T_3$ in the first of Figure \ref{ade_abcA}. Then $\alpha_1\gamma_3\cdots=\alpha\beta\gamma$ and no $\alpha^2\cdots$ determine $T_4$. Then $\alpha_4\gamma_1\cdots=\alpha\beta\gamma$ and no $\delta\thin\epsilon\cdots$ determine $T_5$. Then $\alpha\epsilon\cdots=\alpha\delta\epsilon$ determines $T_6$. Then we get $\gamma_2\delta_1\epsilon_6\cdots$, contradicting to $\epsilon\cdots=\alpha\delta\epsilon,\beta\epsilon^2$.

We conclude $\beta^2\cdots$ is a vertex. By the list of $\epsilon\cdots$, we know the vertex has no $\epsilon$. Therefore a $b$-vertex $\beta^2\cdots=\beta^2\delta^2\cdots$, with no $\epsilon$ in the remainder. Then by $\alpha>\beta>\gamma$, and $\beta+\gamma>\pi$, and $\beta+2\delta>\pi$, we find a $b$-vertex $\beta^2\cdots=\beta^2\delta^2$. Combined with the list of all $\hat{b}$-vertices, we get $\beta^2\cdots=\beta^3,\beta^2\delta^2$.

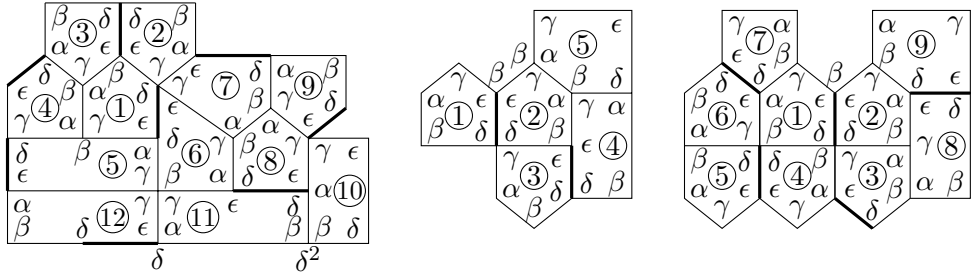
\begin{figure}[htp]
\centering
\begin{tikzpicture}[>=latex,scale=1]


\begin{scope}[shift={(-5.5cm,1cm)}]

\foreach \a in {0,1,2}
\draw[xshift=\a cm]
	(0,0.9) -- (0,0.2) -- (-0.5,-0.2) -- (-0.5,-0.9);

\foreach \a in {0,1}
\draw[xshift=\a cm]
	(0,0.2) -- (0.5,-0.2);
	
\draw
	(0,0.9) -- (2,0.9)
	(-0.5,-0.9) -- (1.5,-0.9)
	(-0.5,-1.6) -- (3.5,-1.6)
	(2.5,-1.6) -- (2.5,-0.9)
	(1.5,-0.9) -- (1.5,-2.3)
	(1.5,-0.2) -- (2.5,-0.9) -- (3,-0.5) -- (3.5,-0.9) -- (3.5,-2.3)
	(-0.5,-0.9) -- (-0.5,-2.3) -- (4.3,-2.3) -- (4.3,-0.9) -- (3.5,-0.9)
	(2,0.2) -- (4,0.2) -- (4,-0.5) -- (3.5,-0.9)
	(3,-0.5) -- (3,0.2);

\draw[line width=1.2]
	(1.5,-0.2) -- (1.5,-0.9)
	(1,0.2) -- (1,0.9)
	(-0.5,-1.6) -- (-0.5,-0.9)
	(0,0.2) -- (-0.5,-0.2)
	(0.5,-2.3) -- (1.5,-2.3)
	(2.5,-1.6) -- (3.5,-1.6)
	(2,0.2) -- (3,0.2)
	(4,-0.5) -- (3.5,-0.9);
		
\node at (1.3,-0.3) {\small $\delta$};	
\node at (0.95,-0.05) {\small $\beta$};
\node at (1.3,-0.7) {\small $\epsilon$}; 
\node at (0.7,-0.3) {\small $\alpha$};
\node at (0.7,-0.7) {\small $\gamma$};
		
\node at (-0.3,-0.3) {\small $\epsilon$};	
\node at (0,-0.05) {\small $\delta$};
\node at (-0.3,-0.7) {\small $\gamma$};
\node at (0.3,-0.3) {\small $\beta$};
\node at (0.3,-0.7) {\small $\alpha$}; 

\node at (0.2,0.3) {\small $\alpha$};	
\node at (0.5,0.05) {\small $\gamma$};
\node at (0.2,0.7) {\small $\beta$};	
\node at (0.8,0.3) {\small $\epsilon$};
\node at (0.8,0.7) {\small $\delta$}; 	
		
\node at (1.8,0.3) {\small $\alpha$};	
\node at (1.5,0.05) {\small $\gamma$};
\node at (1.8,0.7) {\small $\beta$}; 
\node at (1.2,0.3) {\small $\epsilon$};
\node at (1.2,0.7) {\small $\delta$};

\node at (-0.3,-1.1) {\small $\delta$}; 
\node at (-0.3,-1.4) {\small $\epsilon$};
\node at (0.5,-1.1) {\small $\beta$}; 
\node at (1.3,-1.4) {\small $\gamma$};
\node at (1.3,-1.1) {\small $\alpha$}; 

\node at (-0.3,-1.8) {\small $\alpha$}; 
\node at (-0.3,-2.1) {\small $\beta$}; 
\node at (1.3,-1.8) {\small $\gamma$};
\node at (0.5,-2.1) {\small $\delta$};
\node at (1.3,-2.1) {\small $\epsilon$}; 

\node at (1.7,-2.1) {\small $\alpha$};
\node at (3.3,-2.1) {\small $\beta$};
\node at (1.7,-1.8) {\small $\gamma$};
\node at (3.3,-1.8) {\small $\delta$};
\node at (2.5,-1.8) {\small $\epsilon$};

\node at (1.5,-2.5) {\small $\delta$};
\node at (3.5,-2.5) {\small $\delta^2$};

\node at (2.3,-1) {\small $\gamma$}; 
\node at (2.3,-1.4) {\small $\alpha$};
\node at (1.7,-0.5) {\small $\epsilon$}; 
\node at (1.7,-1.4) {\small $\beta$};
\node at (1.7,-0.9) {\small $\delta$}; 

\node at (3,-0.75) {\small $\alpha$};
\node at (2.7,-1) {\small $\beta$};
\node at (3.3,-1) {\small $\gamma$};
\node at (2.7,-1.4) {\small $\delta$};
\node at (3.3,-1.4) {\small $\epsilon$};

\node at (3.7,-1.6) {\small $\alpha$};
\node at (3.7,-2.1) {\small $\beta$};
\node at (3.7,-1.1) {\small $\gamma$};
\node at (4.1,-2.1) {\small $\delta$};
\node at (4.1,-1.1) {\small $\epsilon$};

\node at (2.5,-0.7) {\small $\alpha$};
\node at (2.8,-0.4) {\small $\beta$};
\node at (2.8,0) {\small $\delta$};  
\node at (1.8,-0.2) {\small $\gamma$}; 
\node at (2,0) {\small $\epsilon$}; 

\node at (3.2,0) {\small $\alpha$};
\node at (3.8,0) {\small $\beta$};
\node at (3.2,-0.4) {\small $\gamma$};
\node at (3.8,-0.4) {\small $\delta$};
\node at (3.5,-0.65) {\small $\epsilon$};

\node[inner sep=0.5,draw,shape=circle] at (1,-0.5) {\small $1$};
\node[inner sep=0.5,draw,shape=circle] at (0.5,0.5) {\small $3$};
\node[inner sep=0.5,draw,shape=circle] at (1.5,0.5) {\small $2$};
\node[inner sep=0.5,draw,shape=circle] at (0,-0.5) {\small $4$};
\node[inner sep=0.5,draw,shape=circle] at (0.9,-1.25) {\small $5$};
\node[inner sep=0.5,draw,shape=circle] at (2,-1.1) {\small $6$};
\node[inner sep=0.5,draw,shape=circle] at (2.4,-0.2) {\small $7$};
\node[inner sep=0.5,draw,shape=circle] at (3,-1.2) {\small $8$};
\node[inner sep=0.5,draw,shape=circle] at (3.5,-0.2) {\small $9$};
\node[inner sep=0,draw,shape=circle] at (4.05,-1.6) {\footnotesize $10$};
\node[inner sep=0,draw,shape=circle] at (2.1,-1.95) {\footnotesize $11$};
\node[inner sep=0,draw,shape=circle] at (0.9,-1.95) {\footnotesize $12$};

\end{scope}


\draw
	(1.5,0) -- (-0.5,0) -- (-0.5,0.7) -- (0,1.1) -- (0.5,0.7)
	(1.5,0.7) -- (2.3,0.7) 
	(1,1.1) -- (1,1.8) -- (2.3,1.8) -- (2.3,-0.7) -- (1.5,-0.7)	
	(0.5,-0.7) -- (0.5,0.7) -- (1,1.1) -- (1.5,0.7) -- (1.5,-0.7) -- (1,-1.1) -- (0.5,-0.7);

\draw[line width=1.2]
	(0.5,0) -- (0.5,0.7)
	(1.5,0) -- (1.5,-0.7);

\node at (-0.3,0.2) {\small $\beta$};
\node at (-0.3,0.6) {\small $\alpha$};
\node at (0.3,0.2) {\small $\delta$}; 
\node at (0,0.85) {\small $\gamma$};
\node at (0.3,0.6) {\small $\epsilon$};
	
\node at (0.7,-0.6) {\small $\alpha$};
\node at (1,-0.85) {\small $\beta$};
\node at (0.7,-0.2) {\small $\gamma$};
\node at (1.3,-0.6) {\small $\delta$};	
\node at (1.3,-0.2) {\small $\epsilon$}; 

\node at (1.3,0.6) {\small $\alpha$};	
\node at (1.3,0.2) {\small $\beta$};
\node at (1,0.85) {\small $\gamma$};
\node at (0.7,0.2) {\small $\delta$}; 
\node at (0.7,0.6) {\small $\epsilon$};

\node at (2.1,0.5) {\small $\alpha$};
\node at (2.1,-0.5) {\small $\beta$};
\node at (1.7,0.5) {\small $\gamma$};
\node at (1.7,-0.5) {\small $\delta$};
\node at (1.7,0) {\small $\epsilon$};

\node at (1.2,1.2) {\small $\alpha$};
\node at (1.6,0.9) {\small $\beta$};
\node at (1.2,1.6) {\small $\gamma$};
\node at (2.1,0.9) {\small $\delta$};
\node at (2.1,1.6) {\small $\epsilon$};

\node at (0.5,0.95) {\small $\beta$};
\node at (0.8,1.2) {\small $\beta$};

\node[inner sep=0.5,draw,shape=circle] at (0,0.4) {\small $1$};
\node[inner sep=0.5,draw,shape=circle] at (1,0.4) {\small $2$};
\node[inner sep=0.5,draw,shape=circle] at (1,-0.4) {\small $3$};
\node[inner sep=0.5,draw,shape=circle] at (2.05,0) {\small $4$};
\node[inner sep=0.5,draw,shape=circle] at (1.65,1.35) {\small $5$};


\begin{scope}[xshift=3.5cm]

\foreach \a in {0,1,2}
\draw[xshift=\a cm]
	(-0.5,-0.7) -- (-0.5,0.7) -- (0,1.1) -- (0.5,0.7) -- (0.5,-0.7) -- (0,-1.1) -- (-0.5,-0.7);

\draw
	(-0.5,0) -- (2.5,0)
	(0,1.1) -- (0,1.8) -- (1,1.8) -- (1,1.1)
	(2.5,-0.7) -- (3.3,-0.7) -- (3.3,1.8) -- (2,1.8) -- (2,1.1)
	(2.5,0.7) -- (3.3,0.7);

\draw[line width=1.2]
	(1.5,0) -- (1.5,0.7)
	(0.5,0) -- (0.5,-0.7)
	(0,1.1) -- (0.5,0.7)
	(2,-1.1) -- (1.5,-0.7)
	(2.5,0.7) -- (3.3,0.7);

\node at (-0.3,0.6) {\small $\beta$};
\node at (0,0.85) {\small $\delta$};
\node at (-0.3,0.2) {\small $\alpha$};	
\node at (0.3,0.2) {\small $\gamma$}; 
\node at (0.3,0.6) {\small $\epsilon$};

\node at (0.7,0.2) {\small $\beta$};
\node at (0.7,0.6) {\small $\alpha$};
\node at (1.3,0.2) {\small $\delta$}; 
\node at (1,0.85) {\small $\gamma$};
\node at (1.3,0.6) {\small $\epsilon$};

\node at (2.3,0.6) {\small $\alpha$};	
\node at (2.3,0.2) {\small $\beta$};
\node at (2,0.85) {\small $\gamma$};
\node at (1.7,0.2) {\small $\delta$}; 
\node at (1.7,0.6) {\small $\epsilon$};

\node at (1.5,0.95) {\small $\beta$};

\node at (-0.3,-0.6) {\small $\alpha$};
\node at (0,-0.85) {\small $\gamma$};
\node at (-0.3,-0.2) {\small $\beta$};	
\node at (0.3,-0.2) {\small $\delta$}; 
\node at (0.3,-0.6) {\small $\epsilon$};

\node at (1.3,-0.6) {\small $\alpha$};	
\node at (1.3,-0.2) {\small $\beta$}; 
\node at (1,-0.85) {\small $\gamma$};
\node at (0.7,-0.2) {\small $\delta$};
\node at (0.7,-0.6) {\small $\epsilon$};

\node at (1.7,-0.6) {\small $\epsilon$};
\node at (2,-0.85) {\small $\delta$};
\node at (1.7,-0.2) {\small $\gamma$};
\node at (2.3,-0.6) {\small $\beta$};	
\node at (2.3,-0.2) {\small $\alpha$}; 

\node at (0.8,1.6) {\small $\alpha$};
\node at (0.8,1.2) {\small $\beta$};
\node at (0.2,1.6) {\small $\gamma$};
\node at (0.5,0.95) {\small $\delta$};
\node at (0.2,1.2) {\small $\epsilon$};

\node at (3.1,0.5) {\small $\delta$};
\node at (3.1,-0.5) {\small $\beta$};
\node at (2.7,0.5) {\small $\epsilon$};
\node at (2.7,-0.5) {\small $\alpha$};
\node at (2.7,0) {\small $\gamma$};

\node at (2.2,1.2) {\small $\beta$};
\node at (2.6,0.9) {\small $\delta$};
\node at (2.2,1.6) {\small $\alpha$};
\node at (3.1,0.9) {\small $\epsilon$};
\node at (3.1,1.6) {\small $\gamma$};

\node[inner sep=0.5,draw,shape=circle] at (1,0.4) {\small $1$};
\node[inner sep=0.5,draw,shape=circle] at (2,0.4) {\small $2$};
\node[inner sep=0.5,draw,shape=circle] at (2,-0.4) {\small $3$};
\node[inner sep=0.5,draw,shape=circle] at (1,-0.4) {\small $4$};
\node[inner sep=0.5,draw,shape=circle] at (0,-0.4) {\small $5$};
\node[inner sep=0.5,draw,shape=circle] at (0,0.4) {\small $6$};
\node[inner sep=0.5,draw,shape=circle] at (0.5,1.4) {\small $7$};
\node[inner sep=0.5,draw,shape=circle] at (3.05,0) {\small $8$};
\node[inner sep=0.5,draw,shape=circle] at (2.65,1.35) {\small $9$};

\end{scope}

\end{tikzpicture}
\caption{Proposition \ref{ade_abc}: $\beta\epsilon^2$ is a vertex.}
\label{ade_abcA}
\end{figure}

\subsubsection*{Subsubcase. $\beta^3$ is a vertex} 

The angle sums of $\alpha\delta\epsilon,\alpha\beta\gamma,\beta\epsilon^2,\beta^3$ and the angle sum for pentagon imply
\[ 
\alpha=(1-\tfrac{4}{f})\pi,\,
\beta=\epsilon=\tfrac{2}{3}\pi,\,
\gamma=\delta=(\tfrac{1}{3}+\tfrac{4}{f})\pi.
\]
Then we get all the vertices besides $\alpha\delta\epsilon,\alpha\beta\gamma,\beta\epsilon^2,\beta^3$
\begin{align*}
f=24 &\colon \gamma^4,\gamma^2\delta^2,\delta^4. \\
f=36 &\colon \beta\gamma^3,\beta\gamma\delta^2,\gamma^2\delta\epsilon,\delta^3\epsilon. \\
f=60 &\colon \gamma^5,\gamma^3\delta^2,\gamma\delta^4. 
\end{align*}
By $\beta\epsilon^2$ and the balance lemma, we know $\delta^2\cdots$ is a vertex. 

By no $\alpha^2\cdots$, the AAD of $\beta^2\cdots=\beta^3$ is $\thin^{\alpha}\beta^{\delta}\thin^{\alpha}\beta^{\delta}\thin^{\alpha}\beta^{\delta}\thin$. This implies the AAD of $\thin\beta\thin\beta\thin$ is $\thin^{\alpha}\beta^{\delta}\thin^{\alpha}\beta^{\delta}\thin$, and further implies $\delta\thin\delta\cdots$ is not a vertex. Therefore $\delta^4,\delta^3\epsilon,\gamma\delta^4$ are not vertices. This implies $\delta\thin\epsilon\cdots=\delta^3\epsilon$ is not a vertex, and $\delta^2\cdots=\beta\gamma\delta^2,\gamma^2\delta^2,\gamma^3\delta^2$. Then by no $\alpha^2\cdots,\delta\thin\epsilon\cdots,\epsilon\thin\epsilon\cdots$, the AADs of $\delta^2\cdots$ are $\thick^{\epsilon}\delta^{\beta}\thin^{\delta}\beta^{\alpha}\thin^{\epsilon}\gamma^{\alpha}\thin^{\beta}\delta^{\epsilon}\thick$, $\thick^{\epsilon}\delta^{\beta}\thin^{\alpha}\beta^{\delta}\thin^{\alpha}\gamma^{\epsilon}\thin^{\beta}\delta^{\epsilon}\thick$, $\thick^{\epsilon}\delta^{\beta}\thin^{\alpha}\gamma^{\epsilon}\thin^{\alpha}\gamma^{\epsilon}\thin^{\beta}\delta^{\epsilon}\thick$, $\thick^{\epsilon}\delta^{\beta}\thin^{\alpha}\gamma^{\epsilon}\thin^{\alpha}\gamma^{\epsilon}\thin^{\alpha}\gamma^{\epsilon}\thin^{\beta}\delta^{\epsilon}\thick$. The last three AADs determine $T_1,T_2,T_3$ in the second of Figure \ref{ade_abcA}. Then $\beta_2\epsilon_3\cdots=\beta\epsilon^2$ determines $T_4$. Then $\alpha_2\gamma_4\cdots=\alpha\beta\gamma$ and no $\alpha^2\cdots$ determine $T_5$. Then $\epsilon_1\epsilon_2\cdots=\beta\epsilon^2$ and $\alpha_5\gamma_2\cdots=\alpha\beta\gamma$ imply two $\beta$ adjacent in a tile, a contradiction. 

The first AAD $\thick^{\epsilon}\delta^{\beta}\thin^{\delta}\beta^{\alpha}\thin^{\epsilon}\gamma^{\alpha}\thin^{\beta}\delta^{\epsilon}\thick$ determines $T_1,T_2,T_3,T_4$ in the third of Figure \ref{ade_abcA}. Then $\beta_1\delta_4\cdots=\beta\gamma\delta^2$ and no $\alpha^2\cdots$ determine $T_5,T_6$. Then $\alpha_1\epsilon_6\cdots=\alpha\delta\epsilon$ determines $T_7$. On the other hand, $\alpha_3\beta_2\cdots=\alpha\beta\gamma$ and no $\alpha^2\cdots$ determine $T_8$. Then $\alpha_2\epsilon_8\cdots=\alpha\delta\epsilon$ determines $T_9$. We also have $\epsilon_1\epsilon_2\cdots=\beta\epsilon^2$. This implies one of $\beta_7\gamma_1\cdots,\beta_9\gamma_2\cdots$ is $\thin^{\alpha}\beta^{\delta}\thin^{\alpha}\gamma^{\epsilon}\thin^{\beta}\delta^{\epsilon}\thick\cdots=\beta\gamma\delta^2=\thick^{\epsilon}\delta^{\beta}\thin^{\alpha}\beta^{\delta}\thin^{\alpha}\gamma^{\epsilon}\thin^{\beta}\delta^{\epsilon}\thick$. This is the second AAD above, which we know leads to contradiction.

\subsubsection*{Subsubcase. $\beta^2\delta^2$ is a vertex} 

The angle sums of $\alpha\beta\gamma,\alpha\delta\epsilon,\beta\epsilon^2,\beta^2\delta^2$ and the angle sum for pentagon imply 
\[
\alpha=(1-\tfrac{4}{f})\pi,\,
\beta=(\tfrac{2}{3}-\tfrac{8}{3f})\pi,\,
\gamma=(\tfrac{1}{3}+\tfrac{20}{3f})\pi,\,
\delta=(\tfrac{1}{3}+\tfrac{8}{3f})\pi,\,
\epsilon=(\tfrac{2}{3}+\tfrac{4}{3f})\pi.
\]
By $\beta>\gamma$, we get $f>28$. Then we get all the vertices besides $\alpha\delta\epsilon,\alpha\beta\gamma,\beta\epsilon^2,\beta^2\delta^2$
\begin{align*}
f=40 &\colon \gamma^4. \\
f=52 &\colon \beta\gamma^3,\gamma\delta^4,\gamma^2\delta\epsilon. \\
f=76 &\colon \gamma^3\delta^2. \\
f=100 &\colon \gamma^5. 
\end{align*}

The AAD $\thick^{\delta}\epsilon^{\gamma}\thin^{\alpha}\beta^{\delta}\thin^{\gamma}\epsilon^{\delta}\thick$ of $\beta\epsilon^2$ determines $T_1,T_2,T_3$ in the first of Figure \ref{ade_abcA}. By the list of vertices, we have $\alpha\gamma\cdots=\alpha\beta\gamma$, and $\alpha\epsilon\cdots=\alpha\delta\epsilon$, and no $\alpha^2\cdots,\delta\thin\epsilon\cdots$. This determines $T_4,T_5,T_6$, similar to the case $\beta^2\cdots$ is not a vertex. Then $\gamma_2\delta_1\epsilon_6\cdots=\gamma^2\delta\epsilon$ and no $\alpha^2\cdots$ determine $T_7$. Then $\alpha_7\gamma_6\cdots=\alpha\beta\gamma$ and no $\alpha^2\cdots$ determine $T_8$. Then $\alpha_8\beta_7\cdots=\alpha\beta\gamma$ and no $\delta\thin\epsilon\cdots$ determine $T_9$. Then $\gamma_8\epsilon_9\cdots=\gamma^2\delta\epsilon$ and no $\epsilon\thin\epsilon\cdots$ determine $T_{10}$. Then $\alpha_{10}\epsilon_8\cdots=\alpha\delta\epsilon$ determines $T_{11}$. Then $\beta_6\gamma_5\gamma_{11}\cdots=\beta\gamma^3$ and no $\alpha^2\cdots$ determine $T_{12}$. Then $\beta_{10}\beta_{11}\cdots=\beta^2\delta^2$ and $\alpha_{11}\epsilon_{12}\cdots=\alpha\delta\epsilon$ imply two $\delta$ in a tile, a contradiction.

\subsubsection*{Subcase. $\beta\epsilon^2$ is not a vertex}

Suppose $\gamma^k\epsilon^2$ is a vertex. Then $\gamma+2\epsilon\le 2\pi$. By $\beta+\gamma=\delta+\epsilon=(1+\tfrac{4}{f})\pi$, this implies $\beta+2\delta\ge (1+\frac{12}{f})\pi>\pi$. Then we get $\alpha+\beta+2\delta\ge (2+\frac{8}{f})\pi>2\pi$. 

By $\beta+2\delta>\pi$ and $\delta+\epsilon>\pi$, we know $\beta\delta^l\epsilon(l\ge 3)$ is not a vertex. Then by no $\beta\epsilon^2$, we get the updated list $\epsilon\cdots=\alpha\delta\epsilon,\gamma^k\epsilon^2,\gamma^k\delta^l\epsilon$. In particular, we know $\alpha\gamma\epsilon\cdots,\beta\epsilon\cdots$ are not vertices, and $\delta\thin\epsilon\cdots=\gamma^k\delta^l\epsilon$. Then by $\beta>\gamma$, and $\beta+\gamma>\pi>\alpha$, and $\beta+2\delta>\pi$, and $\alpha+\beta+2\delta>2\pi$, we know a $b$-vertex $\beta\cdots=\beta\delta^2\cdots=\beta^2\delta^2,\beta\gamma^k\delta^l$. Then by the list of all $\hat{b}$-vertices, we get $\beta\cdots=\alpha\beta\gamma,\beta^3,\beta^2\delta^2,\beta\gamma^k\delta^l$. In particular, we get $\alpha\beta\cdots=\alpha\beta\gamma$, and $\beta^2\cdots=\beta^3,\beta^2\delta^2$.

The angle sum of $\beta^2\delta^2$ implies $\beta+\delta=\pi<\beta+\gamma$. This implies $\gamma>\delta$, and $\gamma+\epsilon>\delta+\epsilon>\pi$. Then we get $k=1$ in $\gamma^k\epsilon^2$, and $\epsilon^2\cdots=\gamma\epsilon^2$. 

The vertex $\delta\thin\delta\cdots$ implies a vertex $\thin^{\alpha}\beta^{\delta}\thin^{\delta}\beta^{\alpha}\thin\cdots=\beta^3,\beta^2\delta^2$. If the vertex is $\beta^3$, then we get $\thin^{\alpha}\beta^{\delta}\thin^{\delta}\beta^{\alpha}\thin^{\alpha}\beta^{\delta}\thin$, contradicting no $\alpha^2\cdots$. Therefore the vertex is $\thin^{\alpha}\beta^{\delta}\thin^{\delta}\beta^{\alpha}\thin\cdots=\beta^2\delta^2=\thick^{\epsilon}\delta^{\beta}\thin^{\alpha}\beta^{\delta}\thin^{\delta}\beta^{\alpha}\thin^{\beta}\delta^{\epsilon}\thick$. This implies $\epsilon^2\cdots=\gamma\epsilon^2$.

The AAD $\thick^{\epsilon}\delta^{\beta}\thin^{\alpha}\beta^{\delta}\thin^{\delta}\beta^{\alpha}\thin^{\beta}\delta^{\epsilon}\thick$ determines $T_1,T_2,T_3,T_4$ in the first of Figure \ref{ade_abcB}. Then $\alpha_3\beta_2\cdots=\alpha\beta\gamma$ and no $\alpha^2\cdots$ determine $T_5$. Then $\alpha_2\epsilon_5\cdots=\alpha\delta\epsilon$ determines $T_6$. By the symmetry of horizontal flip, we also get $\beta$ just outside $\gamma_1$. Then $\epsilon_1\epsilon_2\cdots=\gamma\epsilon^2$ implies either $\gamma_1\cdots$ or $\gamma_2\cdots$ is $\beta\gamma\epsilon\cdots$, a contradiction. We conclude $\delta\thin\delta\cdots$ is not a vertex. 

\begin{figure}[htp]
\centering
\begin{tikzpicture}[>=latex,scale=1]


\foreach \a in {-1,1}
{
\begin{scope}[xscale=\a]

\draw
	(1,-0.7) -- (1,0.7) -- (0.5,1.1) -- (0,0.7) -- (0,-0.7) -- (0.5,-1.1) -- (1,-0.7)
	(0.5,1.1) -- (0.5,1.8);

\node at (0.2,-0.6) {\small $\delta$};
\node at (0.5,-0.85) {\small $\epsilon$};
\node at (0.2,-0.2) {\small $\beta$};
\node at (0.8,-0.6) {\small $\gamma$};	
\node at (0.8,-0.2) {\small $\alpha$}; 

\node at (0.8,0.6) {\small $\alpha$};	
\node at (0.8,0.2) {\small $\beta$};
\node at (0.5,0.85) {\small $\gamma$};
\node at (0.2,0.2) {\small $\delta$}; 
\node at (0.2,0.6) {\small $\epsilon$};

\node at (0.7,1.2) {\small $\beta$};

\end{scope}
}

\draw
	(1,-0.7) -- (1.8,-0.7) -- (1.8,1.8) -- (0.5,1.8) 
	(-1,0) -- (1,0);

\draw[line width=1.2]
	(0,0) -- (0,0.7)
	(-0.5,-1.1) -- (0,-0.7) -- (0.5,-1.1)
	(1,0.7) -- (1.8,0.7);

\node at (0,0.95) {\small $\gamma$};

\node at (1.6,0.5) {\small $\delta$};
\node at (1.6,-0.5) {\small $\beta$};
\node at (1.2,0.5) {\small $\epsilon$};
\node at (1.2,-0.5) {\small $\alpha$};
\node at (1.2,0) {\small $\gamma$};

\node at (0.7,1.6) {\small $\alpha$};
\node at (1.1,0.9) {\small $\delta$};
\node at (1.6,0.9) {\small $\epsilon$};
\node at (1.6,1.6) {\small $\gamma$};

\node[inner sep=0.5,draw,shape=circle] at (-0.5,0.4) {\small $1$};
\node[inner sep=0.5,draw,shape=circle] at (0.5,0.4) {\small $2$};
\node[inner sep=0.5,draw,shape=circle] at (0.5,-0.4) {\small $3$};
\node[inner sep=0.5,draw,shape=circle] at (-0.5,-0.4) {\small $4$};
\node[inner sep=0.5,draw,shape=circle] at (1.55,0) {\small $5$};
\node[inner sep=0.5,draw,shape=circle] at (1.15,1.35) {\small $6$};


\begin{scope}[xshift=3cm]

\foreach \a in {0,1,2}
\draw[xshift=\a cm]
	(-0.5,-0.7) -- (-0.5,0.7) -- (0,1.1) -- (0.5,0.7) -- (0.5,-0.7) -- (0,-1.1) -- (-0.5,-0.7);

\draw
	(-0.5,0) -- (2.5,0)
	(2,1.1) -- (2,1.8)
	(2.5,-0.7) -- (3.3,-0.7) -- (3.3,1.8) -- (1,1.8) -- (1,1.1)
	(2.5,0.7) -- (3.3,0.7);

\draw[line width=1.2]
	(1.5,0) -- (1.5,0.7)
	(0.5,0) -- (0.5,-0.7)
	(0,1.1) -- (0.5,0.7)
	(2,-1.1) -- (1.5,-0.7)
	(2.5,0.7) -- (3.3,0.7)
	(1,1.1) -- (1,1.8);

\node at (-0.3,0.6) {\small $\gamma$};
\node at (0,0.85) {\small $\epsilon$};
\node at (-0.3,0.2) {\small $\alpha$};	
\node at (0.3,0.2) {\small $\beta$}; 
\node at (0.3,0.6) {\small $\delta$};

\node at (0.7,0.2) {\small $\beta$};
\node at (0.7,0.6) {\small $\alpha$};
\node at (1.3,0.2) {\small $\delta$}; 
\node at (1,0.85) {\small $\gamma$};
\node at (1.3,0.6) {\small $\epsilon$};

\node at (2.3,0.6) {\small $\alpha$};	
\node at (2.3,0.2) {\small $\beta$};
\node at (2,0.85) {\small $\gamma$};
\node at (1.7,0.2) {\small $\delta$}; 
\node at (1.7,0.6) {\small $\epsilon$};

\node at (-0.3,-0.6) {\small $\alpha$};
\node at (0,-0.85) {\small $\gamma$};
\node at (-0.3,-0.2) {\small $\beta$};	
\node at (0.3,-0.2) {\small $\delta$}; 
\node at (0.3,-0.6) {\small $\epsilon$};

\node at (1.3,-0.6) {\small $\alpha$};	
\node at (1.3,-0.2) {\small $\beta$}; 
\node at (1,-0.85) {\small $\gamma$};
\node at (0.7,-0.2) {\small $\delta$};
\node at (0.7,-0.6) {\small $\epsilon$};

\node at (2.3,-0.2) {\small $\alpha$}; 
\node at (1.7,-0.2) {\small $\beta$};
\node at (2.3,-0.6) {\small $\gamma$};
\node at (1.7,-0.6) {\small $\delta$};
\node at (2,-0.85) {\small $\epsilon$};	

\node at (0.8,1.2) {\small $\epsilon$};
\node at (0.5,0.95) {\small $\epsilon$};

\node at (3.1,0.5) {\small $\delta$};
\node at (3.1,-0.5) {\small $\beta$};
\node at (2.7,0.5) {\small $\epsilon$};
\node at (2.7,-0.5) {\small $\alpha$};
\node at (2.7,0) {\small $\gamma$};

\node at (2.2,1.6) {\small $\alpha$};
\node at (2.2,1.2) {\small $\beta$};
\node at (2.6,0.9) {\small $\delta$};
\node at (3.1,0.9) {\small $\epsilon$};
\node at (3.1,1.6) {\small $\gamma$};

\node at (1.8,1.6) {\small $\beta$};
\node at (1.8,1.2) {\small $\alpha$};
\node at (1.5,0.95) {\small $\gamma$};
\node at (1.2,1.6) {\small $\delta$};
\node at (1.2,1.2) {\small $\epsilon$};

\node[inner sep=0.5,draw,shape=circle] at (1,0.4) {\small $1$};
\node[inner sep=0.5,draw,shape=circle] at (2,0.4) {\small $2$};
\node[inner sep=0.5,draw,shape=circle] at (2,-0.4) {\small $3$};
\node[inner sep=0.5,draw,shape=circle] at (1,-0.4) {\small $4$};
\node[inner sep=0.5,draw,shape=circle] at (3.05,0) {\small $5$};
\node[inner sep=0.5,draw,shape=circle] at (2.65,1.35) {\small $6$};
\node[inner sep=0.5,draw,shape=circle] at (1.5,1.4) {\small $7$};
\node[inner sep=0.5,draw,shape=circle] at (0,0.4) {\small $8$};
\node[inner sep=0.5,draw,shape=circle] at (0,-0.4) {\small $9$};

\end{scope}

\end{tikzpicture}
\caption{Proposition \ref{ade_abc}: $\beta\epsilon^2$ is not a vertex. }
\label{ade_abcB}
\end{figure}
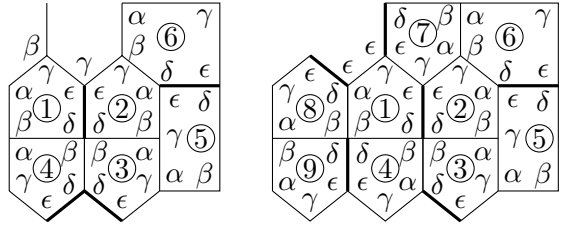

Suppose $\gamma^k\epsilon^2$ is not a vertex. Then $\epsilon^2\cdots$ is not a vertex. By the list of $\epsilon\cdots$ and the balance lemma, we know $\epsilon\cdots=\alpha\delta\epsilon,\gamma^k\delta\epsilon$ are all the $b$-vertices. In particular, we know $\beta\epsilon\cdots,\delta\thin\delta\cdots,\delta\thin\epsilon\cdots$ are not vertices. Moreover, by the list of all $\hat{b}$-vertices, we also get $\alpha\beta\cdots=\alpha\beta\gamma$.

We proved that, no matter $\gamma^k\epsilon^2$ is a vertex or not, we know $\beta\epsilon\cdots,\delta\thin\delta\cdots$ are not vertices, and $\alpha\beta\cdots=\alpha\beta\gamma$, and $\delta\thin\epsilon\cdots=\gamma^k\delta^l\epsilon$. 

By no $\beta\epsilon\cdots$, the AAD of $\thin\gamma\thin\delta\thick$ is $\thin^{\epsilon}\gamma^{\alpha}\thin^{\beta}\delta^{\epsilon}\thick$. By no $\alpha^2\cdots$, this implies a vertex $\thin^{\beta}\alpha^{\gamma}\thin^{\delta}\beta^{\alpha}\thin\cdots=\alpha\beta\gamma=\thin^{\delta}\beta^{\alpha}\thin^{\epsilon}\gamma^{\alpha}\thin^{\beta}\alpha^{\gamma}\thin$. By no $\beta\epsilon\cdots$, this further implies a vertex $\thin^{\beta}\alpha^{\gamma}\thin^{\alpha}\beta^{\delta}\thin\cdots=\alpha\beta\gamma=\thin^{\alpha}\beta^{\delta}\thin^{\epsilon}\gamma^{\alpha}\thin^{\beta}\alpha^{\gamma}\thin$. Then this implies a vertex $\delta\thin\epsilon\cdots=\gamma^k\delta^l\epsilon$. The vertex consists of one $\thick\delta\thin\epsilon\thick$ and several $\delta^2$-fans $\thick\delta\thin\gamma\thin\cdots\thin\gamma\thin\delta\thick$. By no $\thick\delta\thin\delta\thick$, the $\delta^2$-fan has at least one $\gamma$. Then the AAD of the $\delta^2$-fan implies one of $\alpha^2\cdots,\beta\epsilon\cdots,\epsilon\thin\epsilon\cdots$ is a vertex, a contradiction. Therefore $\gamma\thin\delta\cdots$ is not a vertex. 

By no $\gamma\thin\delta\cdots,\delta\thin\delta\cdots$, we know $\gamma^k\delta^l\epsilon$ is not a vertex. Then by no $\beta\epsilon\cdots$, we get updated list $\epsilon\cdots=\alpha\delta\epsilon,\gamma^k\epsilon^2$. Then by $\alpha\beta\gamma$, and applying the counting lemma to $\alpha,\epsilon$, we know $\gamma^k\epsilon^2$ is a vertex. 

By $\gamma^k\epsilon^2$ and the balance lemma, this implies $\delta^2\cdots$ is a vertex. Then by the earlier argument for the case $\gamma^k\epsilon^2$ is a vertex, we know $\beta\delta^2\cdots=\beta^2\delta^2,\beta\gamma^k\delta^l$. By no $\gamma\thin\delta\cdots$, this implies $\beta\delta^2\cdots=\beta^2\delta^2$. It remains to consider $\delta^2\cdots$ without $\beta$. By $\epsilon\cdots=\alpha\delta\epsilon,\gamma^k\epsilon^2$, the vertex has no $\epsilon$. Then by no $\alpha^2\cdots$, the vertex is $\alpha\gamma^k\delta^l,\gamma^k\delta^l$. Then by no $\gamma\thin\delta\cdots,\delta\thin\delta\cdots$, the vertex is $\alpha\delta^2$, contradicting $\alpha\delta\epsilon$. 

We conclude $\delta^2\cdots=\beta^2\delta^2$ is a vertex. By the earlier argument, the vertex implies $k=1$ in $\gamma\epsilon^2$, and $\epsilon^2\cdots=\gamma\epsilon^2$. By $\epsilon\cdots=\alpha\delta\epsilon,\gamma\epsilon^2$, we also get $\beta\delta\cdots=\beta\delta^2\cdots=\beta^2\delta^2$, and $\alpha\delta\cdots=\alpha\delta\epsilon$.

By no $\alpha^2\cdots,\delta\thin\delta\cdots$, we get the AAD $\thick^{\epsilon}\delta^{\beta}\thin^{\alpha}\beta^{\delta}\thin^{\alpha}\beta^{\delta}\thin^{\beta}\delta^{\epsilon}\thick$ of $\beta^2\delta^2$. This determines $T_1,T_2,T_3,T_4$ in the second of Figure \ref{ade_abcB}. Then we determine $T_5,T_6$ similar to the first picture. Then $\epsilon_1\epsilon_2\cdots=\gamma\epsilon^2$ and no $\beta\epsilon\cdots$ determine $T_7$, and the AAD of $\beta_1\delta_4\cdots=\beta^2\delta^2$ determines $T_8,T_9$. Then $\alpha_1\delta_8\cdots=\alpha\delta\epsilon$ and $\gamma_1\epsilon_7\cdots=\gamma\epsilon^2$ imply two adjacent $\epsilon$, a contradiction. 
\end{proof}

\begin{proposition}\label{ade_2bc}
Tilings with distinct $\alpha,\beta,\gamma$ values, and such that $\alpha\delta\epsilon,\beta^2\gamma$ are vertices, are the second earth map tiling {\rm $E_{\pentagon}2$} in Figure \ref{emt}, and the rotation modification {\rm $RE_{\pentagon}2$} in Figure \ref{emt2mod2}.
\end{proposition}

\begin{proof}
The angle sums of $\alpha\delta\epsilon,\beta^2\gamma$ and the angle sum for pentagon imply
\[
\alpha+\delta+\epsilon=2\pi,\;
\beta=(1-\tfrac{4}{f})\pi,\;
\gamma=\tfrac{8}{f}\pi.
\]
We have $\pi>\beta>\gamma$. By Lemma \ref{geometry1}, this implies $\delta<\epsilon$.

\subsubsection*{Case. $\alpha<\pi$, and $\epsilon^2\cdots$ is not a vertex}

By no $\epsilon^2\cdots$ and the balance lemma, we know a $b$-vertex is $\delta\epsilon\cdots$, with no $\delta,\epsilon$ in the remainder. By $\alpha\delta\epsilon$, we get $R(\delta\epsilon)=\alpha<\pi<\beta+\gamma$. By $\beta>\gamma$, this implies $\alpha\delta\epsilon,\gamma^k\delta\epsilon(k\ge 2)$ are all the $b$-vertices. In particular, we know $\beta\cdots$ is a $\hat{b}$-vertex, and $\delta\thin\delta\cdots, \delta\thin\epsilon\cdots,\epsilon\thin\epsilon\cdots$ are not vertices.  

Suppose $\gamma^k\delta\epsilon$ is a vertex. By $\alpha\delta\epsilon,\gamma^k\delta\epsilon$, we get $\alpha=k\gamma>\gamma$. By $\beta^2\gamma$, and $\alpha,\beta>\gamma$, and $\beta\cdots$ being $\hat{b}$-vertex, we get $\beta\cdots=\beta^2\gamma,\alpha^l\beta\gamma^k$. Combined with all the $b$-vertices, we know $\alpha\delta\epsilon,\beta^2\gamma,\alpha^l\beta\gamma^k,\alpha^l\gamma^k,\gamma^k\delta\epsilon$ are all the vertices.

By no $\beta\epsilon\cdots,\epsilon\thin\epsilon\cdots$, we know the AAD of $\gamma\delta\cdots=\gamma\epsilon\cdots=\gamma^k\delta\epsilon$ is $\thick^{\epsilon}\delta^{\beta}\thin^{\alpha}\gamma^{\epsilon}\thin\cdots\thin^{\alpha}\gamma^{\epsilon}\thin^{\gamma}\epsilon^{\delta}\thick$. This implies $\alpha\beta\cdots=\alpha^l\beta\gamma^k$ is a vertex, and the AAD of $\thin\gamma\thin\epsilon\thick$ is $\thin^{\alpha}\gamma^{\epsilon}\thin^{\gamma}\epsilon^{\delta}\thick$. 

By no $\beta\epsilon\cdots$, we know the AAD of $\thin\alpha\thin^{\epsilon}\gamma^{\alpha}\thin$ is $\thin^{\beta}\alpha^{\gamma}\thin^{\epsilon}\gamma^{\alpha}\thin$. This implies $\thin^{\epsilon}\gamma^{\alpha}\thin^{\gamma}\epsilon^{\delta}\thick$, contradicting the AAD $\thin^{\alpha}\gamma^{\epsilon}\thin^{\gamma}\epsilon^{\delta}\thick$ of $\thin\gamma\thin\epsilon\thick$. Therefore the AAD of $\thin\alpha\thin\gamma\thin$ is $\thin\alpha\thin^{\alpha}\gamma^{\epsilon}\thin$. Then by no $\epsilon\thin\epsilon\cdots$, the AAD of $\thin\alpha\thin\gamma\thin\cdots\thin\gamma\thin$ is $\thin\alpha\thin^{\alpha}\gamma^{\epsilon}\thin\cdots\thin^{\alpha}\gamma^{\epsilon}\thin$. 

A vertex $\alpha\gamma\cdots=\alpha^l\beta\gamma^k,\alpha^l\gamma^k$, with $k,l\ge 1$. The vertex contains $\thin\alpha\thin\gamma\thin\cdots\thin\gamma\thin=\thin\alpha\thin^{\alpha}\gamma^{\epsilon}\thin\cdots\thin^{\alpha}\gamma^{\epsilon}\thin$, and we cannot add more $\gamma$ to the right. By the AAD $\thin\alpha\thin^{\alpha}\gamma^{\epsilon}\thin$ of $\thin\alpha\thin\gamma\thin$, we cannot add $\alpha$ to the right. Then by no $\delta\thin\epsilon\cdots$, we get $\thin\alpha\thin^{\alpha}\gamma^{\epsilon}\thin\cdots\thin^{\alpha}\gamma^{\epsilon}\thin^{\alpha}\beta^{\delta}\thin$ in the vertex. Since $\alpha^l\beta\gamma^k$ has only one $\beta$, this implies no more $\thin\alpha\thin\gamma\thin\cdots\thin\gamma\thin$ in the rest of the vertex. Therefore the rest of the vertex consists of only $\alpha$, and $\alpha\gamma\cdots=\thin\alpha\thin\cdots\thin\alpha\thin^{\alpha}\gamma^{\epsilon}\thin\cdots\thin^{\alpha}\gamma^{\epsilon}\thin^{\alpha}\beta^{\delta}\thin$.

The AAD $\thin^{\gamma}\alpha^{\beta}\thin^{\beta}\alpha^{\gamma}\thin$ implies a vertex $\thin^{\delta}\beta^{\alpha}\thin^{\alpha}\beta^{\delta}\thin\cdots=\beta^2\gamma=\thin^{\alpha}\beta^{\delta}\thin^{\alpha}\gamma^{\epsilon}\thin^{\delta}\beta^{\alpha}\thin$, contradicting no $\delta\thin\epsilon\cdots$.

The AAD $\thin^{\beta}\alpha^{\gamma}\thin^{\beta}\alpha^{\gamma}\thin$ implies a vertex $\thin^{\delta}\beta^{\alpha}\thin^{\alpha}\gamma^{\epsilon}\thin\cdots=\beta^2\gamma,\alpha^l\beta\gamma^k$. If the vertex is $\beta^2\gamma$, then we get $\thin^{\alpha}\gamma^{\epsilon}\thin\beta\thin^{\delta}\beta^{\alpha}\thin$. This implies $\delta\thin\delta\cdots, \delta\thin\epsilon\cdots$, a contradiction. Therefore $\thin^{\delta}\beta^{\alpha}\thin^{\alpha}\gamma^{\epsilon}\thin\cdots=\alpha^l\beta\gamma^k$. Since $\thin^{\delta}\beta^{\alpha}\thin^{\alpha}\gamma^{\epsilon}\thin$ is incompatible with $\alpha\gamma\cdots=\thin\alpha\thin\cdots\thin\alpha\thin^{\alpha}\gamma^{\epsilon}\thin\cdots\thin^{\alpha}\gamma^{\epsilon}\thin^{\alpha}\beta^{\delta}\thin$, we get $\thin^{\delta}\beta^{\alpha}\thin^{\alpha}\gamma^{\epsilon}\thin\cdots=\beta\gamma^k$. Then the AAD of the vertex implies $\delta\thin\epsilon\cdots, \epsilon\thin\epsilon\cdots$, a contradiction.

We conclude the AAD of $\thin\alpha\thin\alpha\thin$ is $\thin^{\beta}\alpha^{\gamma}\thin^{\gamma}\alpha^{\beta}\thin$. This implies no consecutive $\alpha\alpha\alpha$. Then $\alpha\gamma\cdots=\thin\alpha\thin^{\alpha}\gamma^{\epsilon}\thin\cdots\thin^{\alpha}\gamma^{\epsilon}\thin^{\alpha}\beta^{\delta}\thin,\thin\alpha\thin\alpha\thin^{\alpha}\gamma^{\epsilon}\thin\cdots\thin^{\alpha}\gamma^{\epsilon}\thin^{\alpha}\beta^{\delta}\thin$. By no $\beta\delta\cdots$, and the AAD $\thin^{\beta}\alpha^{\gamma}\thin^{\gamma}\alpha^{\beta}\thin$ of $\thin\alpha\thin\alpha\thin$, we further get $\alpha\gamma\cdots=\thin^{\gamma}\alpha^{\beta}\thin^{\alpha}\gamma^{\epsilon}\thin\cdots\thin^{\alpha}\gamma^{\epsilon}\thin^{\alpha}\beta^{\delta}\thin$.

Recall that the AAD of $\gamma^k\delta\epsilon$ implies $\alpha\beta\cdots=\alpha^l\beta\gamma^k$ is a vertex. If $k=0$, then by no consecutive $\alpha\alpha\alpha$, we know the vertex is $\alpha^2\beta$. Then by the AAD $\thin^{\beta}\alpha^{\gamma}\thin^{\gamma}\alpha^{\beta}\thin$ of $\thin\alpha\thin\alpha\thin$, the vertex is $\alpha^2\beta=\thin^{\gamma}\alpha^{\beta}\thin^{\alpha}\beta^{\delta}\thin^{\beta}\alpha^{\gamma}\thin$, contradicting no $\beta\delta\cdots$. Therefore we have $k\ge 1$ in $\alpha\beta\cdots=\alpha^l\beta\gamma^k$. Then $\alpha\beta\cdots=\alpha\gamma\cdots=\thin^{\gamma}\alpha^{\beta}\thin^{\alpha}\gamma^{\epsilon}\thin\cdots\thin^{\alpha}\gamma^{\epsilon}\thin^{\alpha}\beta^{\delta}\thin$. This implies the AAD of $\thin\alpha\thin\beta\thin$ is $\thin^{\beta}\alpha^{\gamma}\thin^{\delta}\beta^{\alpha}\thin$. However, the AAD $\thin^{\gamma}\alpha^{\beta}\thin^{\alpha}\gamma^{\epsilon}\thin\cdots\thin^{\alpha}\gamma^{\epsilon}\thin^{\alpha}\beta^{\delta}\thin$ implies $\thin^{\beta}\alpha^{\gamma}\thin^{\alpha}\beta^{\delta}\thin$. The contradiction implies $\alpha\beta\cdots$ is not a vertex. Therefore $\gamma^k\delta\epsilon$ is not a vertex.

By no $\gamma^k\delta\epsilon$, we know $\alpha\delta\epsilon$ is the only $b$-vertex. By the counting lemma (especially the remark after the proof of Lemma \ref{counting}), this implies $\alpha\delta\epsilon,\beta^k\gamma^l$ are all the vertices. Then by $\beta^2\gamma$ and $\beta>\gamma$, we get all the vertices
\begin{equation}\label{ade_2bc_avc}
\text{AVC}=\{\alpha\delta\epsilon,\beta^2\gamma,\beta\gamma^k,\gamma^k\}.
\end{equation}

\subsubsection*{Case. $\alpha\ge\pi$, or $\epsilon^2\cdots$ is a vertex}

By $\alpha\delta\epsilon$, we know $\alpha\ge\pi$ implies $\epsilon<\pi$. Of course the vertex $\epsilon^2\cdots$ also implies $\epsilon<\pi$. Then by $\pi>\beta>\gamma$, and $\delta<\epsilon<\pi$, and Lemma \ref{geometry6}, we know $\beta+2\delta>\pi$, and $\gamma+2\epsilon>\pi$. Then $2\epsilon>\pi-\gamma\ge \gamma$.

If $\alpha\ge \pi$, then $\alpha+\beta+\gamma>2\pi$, and $\alpha+\beta+2\delta>2\pi$. By $\alpha>\beta>\gamma$, and $\delta<\epsilon$, this implies $\alpha^2\cdots,\alpha\beta\cdots$ are not vertices.

Next we assume $\alpha<\pi$ and $\epsilon^2\cdots$ is a vertex, and prove that $\alpha^2\cdots,\alpha\beta\cdots$ are also not vertices. 
 
By $\epsilon^2\cdots$ and the balance lemma, we know $\delta^2\cdots$ is a vertex. This implies $\delta,\epsilon<\pi$. Then by $\alpha,\beta,\gamma<\pi$, we know the pentagon is convex. 

By $\alpha\delta\epsilon$ and $\alpha<\pi$, we get $\delta+\epsilon>\pi$. By $\delta<\epsilon$, this implies $R(\epsilon^2)$ has no $\delta,\epsilon$. By $\alpha\delta\epsilon$ and $\delta<\epsilon$, we get $R(\epsilon^2)<\alpha<\pi<\beta+\gamma$. Then by $\beta>\gamma$, this implies $\epsilon^2\cdots=\beta\epsilon^2,\gamma^k\epsilon^2$. 

If $\beta\epsilon^2$ is a vertex, then by $\alpha\delta\epsilon,\beta^2\gamma$, we get $\beta+\epsilon=\alpha+\delta$ and $2\epsilon=\beta+\gamma$. By $\beta>\gamma$ and $\delta<\epsilon$, we get $\alpha>\beta>\epsilon>\gamma$. Then $\alpha+\epsilon<\alpha+\beta<\pi+(1-\frac{4}{f})\pi=(2-\frac{4}{f})\pi$. By $\alpha\delta\epsilon$, this implies $\delta>\frac{4}{f}\pi$, and $2\delta>\gamma$. Then by $\beta^2\gamma$ and $\alpha>\beta$, we get $R(\alpha^2)<R(\alpha\beta)<R(\beta^2)=\gamma<2\delta$. By $\alpha>\beta>\epsilon>\gamma$, this implies $\alpha^2\cdots,\alpha\beta\cdots$ are not vertices.

If $\gamma\epsilon^2$ is a vertex, then $\beta=\epsilon$ and $\alpha+\delta=\beta+\gamma=\gamma+\epsilon$. By $\delta<\epsilon$, we get $\alpha>\gamma$. If $\alpha<\beta$, then we get $\delta>\gamma$, and $\beta=\epsilon>\alpha>\gamma$, contradicting Lemma \ref{geometry4}. Therefore we have $\alpha>\beta$. Then $\delta<\gamma$, and we get $\alpha>\beta=\epsilon>\gamma>\delta$. Then we have $2\alpha+\delta>\alpha+\beta+\delta=\alpha+\delta+\epsilon=2\pi$. This implies $\alpha^2\cdots,\alpha\beta\cdots$ are not vertices.

If $\gamma^k\epsilon^2(k\ge 2)$ is a vertex, then $\gamma+\epsilon\le \pi<\delta+\epsilon,\beta+\gamma$. This implies $\gamma<\delta$ and $\beta>\epsilon$. Therefore $\pi>\beta>\epsilon>\delta>\gamma$. Then by $\alpha<\pi$, we know the pentagon is convex. By $\alpha\delta\epsilon,\gamma^k\epsilon^2$, and $\delta<\epsilon$, we also get $\alpha>\gamma$. Since the pentagon is convex, we may apply Lemma \ref{geometry4} to get $\alpha>\beta$. By $\beta^2\gamma$ and $\alpha>\beta>\epsilon>\delta>\gamma$, we know $\alpha^2\cdots,\alpha\beta\cdots$ are not vertices.

We conclude that, no matter $\alpha\ge \pi$ or $\alpha<\pi$, we always know that $\alpha^2\cdots,\alpha\beta\cdots$ are not vertices.

By $\beta^2\gamma$, and $\beta>\gamma$, and no $\alpha^2\cdots,\alpha\beta\cdots$, we know $\beta^2\gamma,\alpha\gamma^k,\beta\gamma^k,\gamma^k$ are all the $\hat{b}$-vertices. 

Let $\theta=\beta,\gamma$. By no $\alpha^2\cdots,\alpha\beta\cdots$, we know the AADs of ${}^{\alpha}\thin\theta\thin,{}^{\beta}\thin\theta\thin$ are ${}^{\alpha}\thin\theta^{\alpha}\thin,{}^{\beta}\thin\theta^{\alpha}\thin$. Then we know the AAD of $\alpha\gamma^k$ is $\thin^{\gamma}\alpha^{\beta}\thin^{\epsilon}\gamma^{\alpha}\thin\cdots\thin^{\epsilon}\gamma^{\alpha}\thin$, and the AAD of $\beta^l\gamma^k$ is $\thin\theta^{\alpha}\thin\cdots\thin\theta^{\alpha}\thin$. We also know the AAD of $\thick\delta\thin\theta\thin\cdots\thin\theta\thin$ is $\thick^{\epsilon}\delta^{\beta}\thin\theta^{\alpha}\thin\cdots\thin\theta^{\alpha}\thin$. This implies a $\delta^2$-fan is $\thick^{\epsilon}\delta^{\beta}\thin^{\beta}\delta^{\epsilon}\thick$ or has $\alpha$. Then by no $\alpha^2\cdots,\alpha\beta\cdots$, a $\delta^2$-fan with $\alpha$ has a single $\alpha$ and has no $\beta$. Then by the AAD $\thick^{\epsilon}\delta^{\beta}\thin\theta^{\alpha}\thin\cdots\thin\theta^{\alpha}\thin$ and no $\alpha\beta\cdots$, the fan is $\thick^{\epsilon}\delta^{\beta}\thin^{\epsilon}\gamma^{\alpha}\thin\cdots\thin^{\epsilon}\gamma^{\alpha}\thin^{\gamma}\alpha^{\beta}\thin^{\beta}\delta^{\epsilon}\thick$. We conclude a $\delta^2$-fan is $\thick\delta\thin\delta\thick$, $\thick\delta\thin\gamma\thin\cdots\thin\gamma\thin\alpha\thin\delta\thick$. 

By $\alpha\delta\epsilon$ and $\delta<\epsilon$, a $\delta\epsilon$-fan is the vertex $\alpha\delta\epsilon$ or has only $\theta$, and an $\epsilon^2$-fan has only $\theta$. Then by the AAD $\thick^{\epsilon}\delta^{\beta}\thin\theta^{\alpha}\thin\cdots\thin\theta^{\alpha}\thin$, we know a $\delta\epsilon$-fan without $\alpha$ is $\thick^{\epsilon}\delta^{\beta}\thin\theta^{\alpha}\thin\cdots\thin\theta^{\alpha}\thin^{\gamma}\epsilon^{\delta}\thick$. Therefore the only fan containing $\thin\alpha\thin\gamma\thin$ is the $\delta^2$-fan $\thick\delta\thin\gamma\thin\cdots\thin\gamma\thin\alpha\thin\delta\thick$.

By no $\alpha\beta\cdots$, we know $\beta^2\cdots$ has no $\alpha$. Then by $\beta^2\gamma$, and $\beta>\gamma$, and $\delta<\epsilon$, and $\beta+2\delta>\pi$, and $\gamma<2\epsilon$, we get $\beta^2\cdots=\beta^2\gamma,\beta^2\delta^2,\beta^2\delta\epsilon$. Since a $\delta^2$-fan is $\thick\delta\thin\delta\thick$, $\thick\delta\thin\gamma\thin\cdots\thin\gamma\thin\alpha\thin\delta\thick$, we know $\beta^2\delta^2$ is not a vertex. Moreover, the AAD $\thick^{\epsilon}\delta^{\beta}\thin^{\delta}\beta^{\alpha}\thin^{\delta}\beta^{\alpha}\thin^{\gamma}\epsilon^{\delta}\thick$ of $\beta^2\delta\epsilon$ implies a vertex $\alpha\gamma\cdots$. However, the angle sums of $\alpha\delta\epsilon,\beta^2\gamma,\beta^2\delta\epsilon$ imply $\alpha+\gamma=2\pi$, contradicting the vertex $\alpha\gamma\cdots$. Therefore $\beta^2\delta\epsilon$ is not a vertex. We conclude $\beta^2\cdots=\beta^2\gamma$. This implies a fan has at most one $\beta$. 

The AAD $\thick^{\epsilon}\delta^{\beta}\thin^{\beta}\delta^{\epsilon}\thick$ implies a vertex $\thin^{\alpha}\beta^{\delta}\thin^{\delta}\beta^{\alpha}\thin\cdots=\beta^2\gamma=\thin^{\delta}\beta^{\alpha}\thin^{\alpha}\gamma^{\epsilon}\thin^{\alpha}\beta^{\delta}\thin$, contradicting no $\alpha^2\cdots$. Therefore the only $\delta^2$-fan is $\thick\delta\thin\gamma\thin\cdots\thin\gamma\thin\alpha\thin\delta\thick$.

The AAD $\thin^{\epsilon}\gamma^{\alpha}\thin^{\gamma}\epsilon^{\delta}\thick$ implies a vertex $\thin^{\beta}\alpha^{\gamma}\thin^{\epsilon}\gamma^{\alpha}\thin\cdots$. The AAD $\thin^{\beta}\alpha^{\gamma}\thin^{\epsilon}\gamma^{\alpha}\thin$ is incompatible with the AAD of the $\hat{b}$-vertex $\alpha\gamma\cdots=\alpha\gamma^k$. Moreover, the only fan containing $\thin\alpha\thin\gamma\thin$ is $\thick\delta\thin\gamma\thin\cdots\thin\gamma\thin\alpha\thin\delta\thick$. However, $\thin^{\beta}\alpha^{\gamma}\thin^{\epsilon}\gamma^{\alpha}\thin$ is also incompatible with the AAD of the $\delta^2$-fan. Therefore we do not have $\thin^{\epsilon}\gamma^{\alpha}\thin^{\gamma}\epsilon^{\delta}\thick$, and the AAD of $\thin\gamma\thin\epsilon\thick$ is $\thin^{\alpha}\gamma^{\epsilon}\thin^{\gamma}\epsilon^{\delta}\thick$. 

We already know a $\delta^2$-fan is $\thick^{\epsilon}\delta^{\beta}\thin^{\epsilon}\gamma^{\alpha}\thin\cdots\thin^{\epsilon}\gamma^{\alpha}\thin^{\gamma}\alpha^{\beta}\thin^{\beta}\delta^{\epsilon}\thick$. We also know a $\delta\epsilon$-fan without $\alpha$ is $\thick^{\epsilon}\delta^{\beta}\thin\theta^{\alpha}\thin\cdots\thin\theta^{\alpha}\thin^{\gamma}\epsilon^{\delta}\thick$. Then by the AAD $\thin^{\alpha}\gamma^{\epsilon}\thin^{\gamma}\epsilon^{\delta}\thick$ of $\thin\gamma\thin\epsilon\thick$, and at most one $\beta$ in a fan, we know $\delta\epsilon$-fans are the vertex $\alpha\delta\epsilon$, and $\thick^{\epsilon}\delta^{\beta}\thin^{\gamma}\epsilon^{\delta}\thick,\thick^{\epsilon}\delta^{\beta}\thin^{\epsilon}\gamma^{\alpha}\thin\cdots\thin^{\epsilon}\gamma^{\alpha}\thin^{\delta}\beta^{\alpha}\thin^{\gamma}\epsilon^{\delta}\thick$. We also know an $\epsilon^2$-fan has no $\alpha$ and has at most one $\beta$. Then by the AAD $\thin^{\alpha}\gamma^{\epsilon}\thin^{\gamma}\epsilon^{\delta}\thick$ of $\thin\gamma\thin\epsilon\thick$, and no $\alpha^2\cdots$, we know $\epsilon^2$-fans are $\thick^{\delta}\epsilon^{\gamma}\thin^{\gamma}\epsilon^{\delta}\thick,\thick^{\delta}\epsilon^{\gamma}\thin^{\epsilon}\gamma^{\alpha}\thin\cdots\thin^{\epsilon}\gamma^{\alpha}\thin^{\delta}\beta^{\alpha}\thin^{\gamma}\epsilon^{\delta}\thick$.

The AADs $\thick^{\epsilon}\delta^{\beta}\thin^{\gamma}\epsilon^{\delta}\thick$ and $\thick^{\delta}\epsilon^{\gamma}\thin^{\gamma}\epsilon^{\delta}\thick$ imply vertices $\thin^{\alpha}\beta^{\delta}\thin^{\epsilon}\gamma^{\alpha}\thin\cdots$ and $\thin^{\alpha}\gamma^{\epsilon}\thin^{\epsilon}\gamma^{\alpha}\thin\cdots$. The AADs $\thin^{\alpha}\beta^{\delta}\thin^{\epsilon}\gamma^{\alpha}\thin$ and $\thin^{\alpha}\gamma^{\epsilon}\thin^{\epsilon}\gamma^{\alpha}\thin$ are incompatible with the AADs $\thin\theta^{\alpha}\thin\cdots\thin\theta^{\alpha}\thin$ and $\thin^{\gamma}\alpha^{\beta}\thin^{\epsilon}\gamma^{\alpha}\thin\cdots\thin^{\epsilon}\gamma^{\alpha}\thin$ of $\hat{b}$-vertices $\beta^2\gamma,\alpha\gamma^k,\beta\gamma^k,\gamma^k$. They are also incompatible with the AADs of all the fans. 

The fan $\thick^{\epsilon}\delta^{\beta}\thin^{\epsilon}\gamma^{\alpha}\thin\cdots\thin^{\epsilon}\gamma^{\alpha}\thin^{\delta}\beta^{\alpha}\thin^{\gamma}\epsilon^{\delta}\thick$ implies $\thin^{\alpha}\beta^{\delta}\thin^{\gamma}\epsilon^{\delta}\thick$. The AAD is not compatible with the AAD of any fan. 

Therefore the fans are the vertex $\alpha\delta\epsilon$, and $\thick\delta\thin\gamma\thin\cdots\thin\gamma\thin\alpha\thin\delta\thick,\thick\epsilon\thin\gamma\thin\cdots\thin\gamma\thin\beta\thin\epsilon\thick$. By $\beta^2\cdots=\beta^2\gamma$, and no $\alpha^2\cdots,\alpha\beta\cdots$, we know the fans cannot be combined. Therefore all fans are vertices, and $\alpha\delta\epsilon,\alpha\gamma^k\delta^2(k\ge 1),\beta\gamma^k\epsilon^2$ are all the $b$-vertices.

The vertex $\epsilon^2\cdots=\beta\gamma^k\epsilon^2=\thick^{\delta}\epsilon^{\gamma}\thin^{\epsilon}\gamma^{\alpha}\thin\cdots\thin^{\epsilon}\gamma^{\alpha}\thin^{\delta}\beta^{\alpha}\thin^{\gamma}\epsilon^{\delta}\thick=\thin^{\delta}\beta^{\alpha}\thin^{\gamma}\epsilon^{\delta}\thick^{\delta}\epsilon^{\gamma}\thin\cdots$ determines $T_1,T_2,T_3$ in the first of Figure \ref{ade_2bcA}. Since the AVC $\thin^{\gamma}\alpha_3^{\beta}\thin^{\epsilon}\gamma_1^{\alpha}\thin$ is incompatible with the AAD of $\alpha\gamma^k\delta^2$, we get $\thin^{\gamma}\alpha_3^{\beta}\thin^{\epsilon}\gamma_1^{\alpha}\thin\cdots=\alpha\gamma^k$. Then by no $\alpha^2\cdots$, we determine $T_4$. Then $\alpha_1\epsilon_4\cdots=\alpha\delta\epsilon$ determines $T_5$. Then $\beta_1\beta_5\cdots=\beta^2\gamma$ and no $\alpha^2\cdots$ determine $T_6$. Then $\alpha_6\delta_1\delta_2\cdots=\alpha\gamma^k\delta^2$. However, the AAD $\thick^{\epsilon}\delta_1^{\beta}\thin^{\gamma}\alpha_6^{\beta}\thin$ is incompatible with the AAD of $\alpha\gamma^k\delta^2$.

Therefore $\epsilon^2\cdots$ is not a vertex. By the balance lemma, this implies $\alpha\gamma^k\delta^2$ is not a vertex, and $\alpha\delta\epsilon$ is the only $b$-vertex. Then by the counting lemma, $\alpha\delta\epsilon,\beta^k\gamma^l$ are all the vertices. By $\beta^2\cdots=\beta^2\gamma$ and $\beta>\gamma$, we get the AVC \eqref{ade_2bc_avc}.

\begin{figure}[htp]
\centering
\begin{tikzpicture}[>=latex,scale=1]

\begin{scope}[shift={(-5cm,0.8cm)}, yscale=-1]

\draw	
	(0.5,-0.7) -- (0,-1.1) -- (-0.5,-0.7) -- (-0.5,0.7) -- (0,1.1) -- (0.5,0.7) -- (0.5,-0.7)
	(-0.5,0.7) -- (-1,1.1) -- (-1,1.8) -- (1,1.8) -- (1,1.1) -- (0.5,0.7) -- (2,0.7) -- (2,0) -- (-1.5,0) -- (-1.5,0.7) -- (-1,1.1)
	(0,1.8) -- (0,1.1) ;

\draw[line width=1.2]
	(0.5,0.7) -- (0.5,0)
	(0.5,-0.7) -- (0,-1.1)
	(-0.5,0.7) -- (-1,1.1)	
	(1,1.8) -- (0,1.8);

\node at (0.2,1.6) {\small $\epsilon$};
\node at (0.2,1.2) {\small $\gamma$};
\node at (0.8,1.6) {\small $\delta$};
\node at (0.5,0.95) {\small $\alpha$};
\node at (0.8,1.2) {\small $\beta$}; 

\node at (-1.3,0.6) {\small $\beta$};
\node at (-1,0.85) {\small $\delta$};
\node at (-1.3,0.2) {\small $\alpha$};
\node at (-0.7,0.6) {\small $\epsilon$};	
\node at (-0.7,0.2) {\small $\gamma$};

\node at (-0.3,-0.6) {\small $\gamma$};
\node at (0,-0.85) {\small $\epsilon$};
\node at (-0.3,-0.2) {\small $\alpha$};
\node at (0.3,-0.6) {\small $\delta$};	
\node at (0.3,-0.2) {\small $\beta$}; 

\node at (-0.3,0.6) {\small $\alpha$};
\node at (0,0.85) {\small $\beta$};
\node at (-0.3,0.2) {\small $\gamma$};
\node at (0.3,0.6) {\small $\delta$};	
\node at (0.3,0.2) {\small $\epsilon$}; 

\node at (1.8,0.2) {\small $\alpha$};
\node at (1.8,0.5) {\small $\beta$};  
\node at (1.4,0.2) {\small $\gamma$}; 
\node at (0.7,0.5) {\small $\delta$}; 
\node at (0.7,0.2) {\small $\epsilon$};

\node at (-0.8,1.6) {\small $\gamma$};
\node at (-0.8,1.2) {\small $\epsilon$};
\node at (-0.2,1.6) {\small $\alpha$};
\node at (-0.5,0.95) {\small $\delta$};
\node at (-0.2,1.2) {\small $\beta$};

\node[inner sep=0.5,draw,shape=circle] at (0,0.4) {\small $1$};
\node[inner sep=0.5,draw,shape=circle] at (-1,0.4) {\small $4$};
\node[inner sep=0.5,draw,shape=circle] at (1.05,0.35) {\small $2$};
\node[inner sep=0.5,draw,shape=circle] at (0,-0.4) {\small $3$};
\node[inner sep=0.5,draw,shape=circle] at (-0.5,1.4) {\small $5$};
\node[inner sep=0.5,draw,shape=circle] at (0.5,1.4) {\small $6$};

\end{scope}


\draw[gray!50, xshift=-2 cm, line width=3]
	(0,1.5) -- (0,1.1) -- (0.5,0.7) -- (1.5,0.7) -- (1.5,0) -- (1,-0.4) -- (1,-0.8);
	
\draw[gray!50, xshift=2 cm, line width=3]
	(0,1.5) -- (0,1.1) -- (0.5,0.7) -- (0.5,0) -- (-0.5,0) -- (-1,-0.4) -- (-1,-0.8);

\draw
	(1,-0.8) -- (1,-0.4) -- (0.5,0) -- (-0.5,0) -- (-1,-0.4) -- (-1,-0.8)
	(2,1.1) -- (2.5,0.7) -- (2.5,0) -- (1.5,0) -- (1,-0.4)
	(0.5,0) -- (0.5,0.7);
		
\draw[line width=1.2]
	(1,-0.4) -- (0.5,0);

\foreach \a in {0,-1}
{
\begin{scope}[xshift=2*\a cm]

\draw
	(0,1.5) -- (0,1.1) -- (0.5,0.7) -- (1.5,0.7) -- (2,1.1) -- (2,1.5)
	(1.5,0) -- (1.5,0.7);

\draw[line width=1.2]
	(2,1.1) -- (1.5,0.7);
	
\node at (1,1.4) {\small $\gamma$};
\node at (1.8,1.2) {\small $\epsilon$}; 
\node at (1.4,0.9) {\small $\delta$};
\node at (0.6,0.9) {\small $\beta$};
\node at (0.2,1.2) {\small $\alpha$};
	
\node at (1.7,0.6) {\small $\epsilon$};
\node at (2.3,0.6) {\small $\beta$};
\node at (1.7,0.2) {\small $\gamma$}; 
\node at (2.3,0.2) {\small $\alpha$};	
\node at (2,0.85) {\small $\delta$};	

\end{scope}
}

\node at (1.3,0.5) {\small $\alpha$}; 
\node at (0.7,0.5) {\small $\gamma$};
\node at (0.7,0.1) {\small $\epsilon$};	
\node at (1,-0.15) {\small $\delta$};	
\node at (1.3,0.1) {\small $\beta$};

\node at (0,-0.7) {\small $\gamma$};
\node at (0.8,-0.5) {\small $\epsilon$}; 
\node at (0.4,-0.2) {\small $\delta$};
\node at (-0.4,-0.2) {\small $\beta$};
\node at (-0.8,-0.5) {\small $\alpha$}; 

\node at (-2,1.7) {\footnotesize $l_2$}; 
\node at (2,1.7) {\footnotesize $l_1$};

\node[inner sep=0.5,draw,shape=circle] at (1,1) {\small 1};
\node[inner sep=0.5,draw,shape=circle] at (-1,1) {\small 2};
\node[inner sep=0.5,draw,shape=circle] at (2,0.4) {\small 3};
\node[inner sep=0.5,draw,shape=circle] at (0,0.4) {\small 4};
\node[inner sep=0.5,draw,shape=circle] at (1,0.3) {\small 5};
\node[inner sep=0.5,draw,shape=circle] at (0,-0.3) {\small 6};

\end{tikzpicture}
\caption{Proposition \ref{ade_2bc}: No $\epsilon^2\cdots$, and ${\mc P}^{q+1}_q$.}
\label{ade_2bcA}
\end{figure}
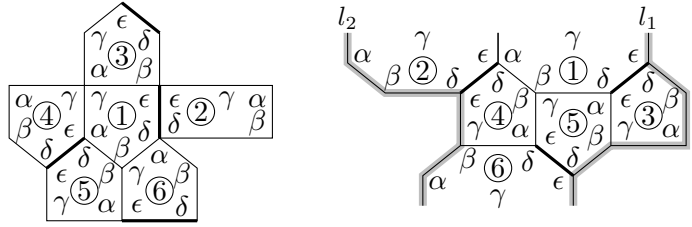

\medskip

\noindent{\em Tiling}

\medskip

Based on the AVC \eqref{ade_2bc_avc}, we construct the tiling. 

Since $\alpha\delta\epsilon$ is the only $b$-vertex, a tile determines its companion. Moreover, by no $\alpha^2\cdots$, the AADs of $\beta\gamma^k,\gamma^k$ are $\thin^{\alpha}\beta^{\delta}\thin^{\alpha}\gamma^{\epsilon}\thin\cdots\thin^{\alpha}\gamma^{\epsilon}\thin,\thin^{\alpha}\gamma^{\epsilon}\thin\cdots\thin^{\alpha}\gamma^{\epsilon}\thin$.

In the second of Figure \ref{ade_2bcA}, the AAD $\thin^{\alpha}\gamma^{\epsilon}\thin^{\alpha}\gamma^{\epsilon}\thin$ determines $T_1,T_2$. The tiles determine their companions $T_3,T_4$. Then $\delta_1\epsilon_3\cdots=\alpha\delta\epsilon$ and $\beta_1\beta_4\cdots=\beta^2\gamma$ determine $T_5$. Then $T_5$ determine its companion $T_6$.

If $\gamma^k$ is a vertex, then we may apply the argument to the $k$ pairs of $\thin^{\alpha}\gamma^{\epsilon}\thin^{\alpha}\gamma^{\epsilon}\thin$ in $\gamma^k$. We get an earth map tiling, with $T_1,T_3,T_5,T_6$ as a timezone. Another choice of the timezone is $T_1,T_4,T_5,T_6$. The second choice gives the earth map tiling $E_{\pentagon}2$ in Figure \ref{emt}.

If $\gamma^k$ is not a vertex, then the AVC is 
\[
\text{AVC}=\{\alpha\delta\epsilon,\beta^2\gamma,\beta\gamma^{q+1}\},\quad
f=8q+4.
\]
The $\gamma^{q+1}$ part of $\beta\gamma^{q+1}$ consists of $q$ pairs of $\thin^{\alpha}\gamma^{\epsilon}\thin^{\alpha}\gamma^{\epsilon}\thin$. Applying the second of Figure \ref{ade_2bcA} to these pairs, we get a partial earth map tiling ${\mc P}^{q+1}_q$, with $\gamma^{q+1}$ and $\gamma^q$ as the two ends, and with the shaded edges $l_1,l_2$ as the boundary (the picture has $q=1$). This gives the left of $l_1$ and the right of $l_2$ in Figure \ref{ade_2bcB}. Moreover, the AAD of $\beta\gamma^{q+1}$ determines $T_1$. Then $T_1$ determines its companion $T_2$.

We have $\gamma_2\cdots=\beta\gamma\cdots=\beta^2\gamma,\beta\gamma^{q+1}$. The case $\gamma_2\cdots=\beta^2\gamma$ is the first picture. Since the vertex $\delta\epsilon\cdots$ on the right of $l_2$ is $\alpha\delta\epsilon$, we determine $T_3$. Then $\gamma_3\cdots=\beta\gamma^2\cdots=\beta\gamma^{q+1}$. The $\gamma^{q+1}$ part of the vertex induces a partial earth map tiling ${\mc P}^{q+1}_q$. If we follow the angles along the boundary of this ${\mc P}^{q+1}_q$, we find the other end $\gamma^q$ of ${\mc P}^{q+1}_q$ is at $\gamma_1\cdots$. Therefore $\gamma_1\cdots=\beta\gamma^q\cdots=\beta\gamma^{q+1}$. The $\gamma^{q+1}$ part of this vertex induces another partial earth map tiling ${\mc P}^{q+1}_q$, which consists of $P,T_1,T_2,T_3$, and is exactly the part between $l_1,l_2$. Therefore the whole tiling is obtained by gluing two copies of ${\mc P}^{q+1}_q$ together. The tiling is the rotation modification (by $120^{\circ}$) $RE_{\pentagon}2$ in Figure \ref{emt2mod2}.

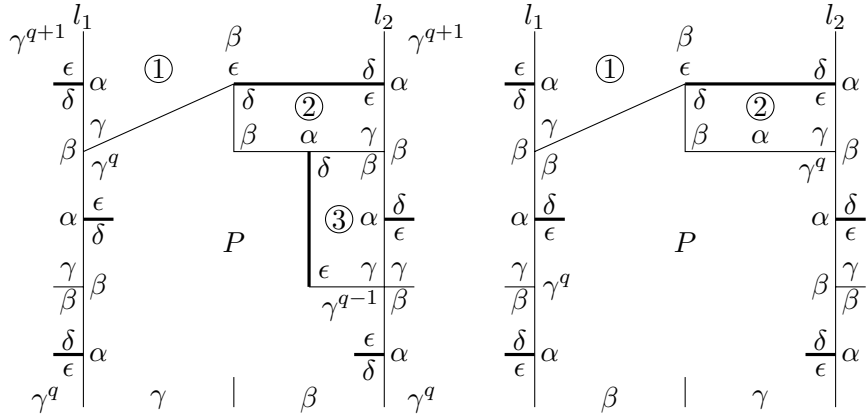
\begin{figure}[htp]
\centering
\begin{tikzpicture}[>=latex,scale=1]

\foreach \a in {0,1}
{
\begin{scope}[xshift=6*\a cm]

\draw
	(2,0.9) -- (0,0.9) -- (0,1.8)
	(-2,0.9) -- (0,1.8)
	(-2,2.5) -- (-2,-2.5)
	(2,2.5) -- (2,-2.5)
	(-2,-0.9) -- ++(-0.4,0)
	(2,-0.9) -- ++(0.4,0);

\draw[line width=1.2]
	(2,1.8) -- (0,1.8)
	(2,0) -- ++(0.4,0) 
	(-2,1.8) -- ++(-0.4,0)
	(-2,-1.8) -- ++(-0.4,0);

\node at (-2,2.7) {\small $l_1$};
\node at (2,2.7) {\small $l_2$};

\node at (2.2,1.8) {\small $\alpha$};
\node at (2.2,0.9) {\small $\beta$};
\node at (2.2,0.2) {\small $\delta$};
\node at (2.2,-0.2) {\small $\epsilon$};
\node at (2.2,-0.7) {\small $\gamma$};
\node at (2.2,-1.1) {\small $\beta$};
\node at (2.2,-1.8) {\small $\alpha$};

\node at (-2.2,2) {\small $\epsilon$};
\node at (-2.2,1.6) {\small $\delta$};
\node at (-2.2,0.9) {\small $\beta$};
\node at (-2.2,0) {\small $\alpha$};
\node at (-2.2,-0.7) {\small $\gamma$};
\node at (-2.2,-1.1) {\small $\beta$};
\node at (-2.2,-1.6) {\small $\delta$};
\node at (-2.2,-2) {\small $\epsilon$};

\node at (0,2.4) {\small $\beta$};
\node at (1.8,2) {\small $\delta$};
\node at (-1.8,1.8) {\small $\alpha$};
\node at (0,2) {\small $\epsilon$};
\node at (-1.8,1.2) {\small $\gamma$};

\node at (1,1.1) {\small $\alpha$};
\node at (1.8,1.1) {\small $\gamma$};
\node at (0.2,1.1) {\small $\beta$};
\node at (1.8,1.6) {\small $\epsilon$};
\node at (0.2,1.6) {\small $\delta$};

\node at (0,-0.3) {\small $P$};

\node[inner sep=0.5,draw,shape=circle] at (-1,2) {\small 1};
\node[inner sep=0.5,draw,shape=circle] at (1,1.5) {\small 2};

\end{scope}
}


\draw
	(2,-0.9) -- (1,-0.9)
	(0,-2.5) -- ++(0,0.4);

\draw[line width=1.2]
	(1,-0.9) -- (1,0.9)
	(2,-1.8) -- ++(-0.4,0)
	(-2,0) -- ++(0.4,0);

\node at (-2.6,2.4) {\small $\gamma^{q+1}$};
\node at (2.7,2.4) {\small $\gamma^{q+1}$};
\node at (-2.5,-2.4) {\small $\gamma^q$};
\node at (2.5,-2.4) {\small $\gamma^q$};
	
\node at (1.8,0) {\small $\alpha$};
\node at (1.8,0.7) {\small $\beta$};
\node at (1.8,-0.7) {\small $\gamma$};
\node at (1.2,0.7) {\small $\delta$};
\node at (1.2,-0.7) {\small $\epsilon$};

\node at (1.55,-1.15) {\small $\gamma^{q-1}$};
\node at (1.8,-1.6) {\small $\epsilon$};
\node at (1.8,-2) {\small $\delta$};

\node at (-1.7,0.7) {\small $\gamma^q$};
\node at (-1.8,0.2) {\small $\epsilon$};
\node at (-1.8,-0.2) {\small $\delta$};
\node at (-1.8,-0.9) {\small $\beta$};
\node at (-1.8,-1.8) {\small $\alpha$};

\node at (1,-2.4) {\small $\beta$};
\node at (-1,-2.4) {\small $\gamma$};

\node[inner sep=0.5,draw,shape=circle] at (1.4,0) {\small 3};


\begin{scope}[xshift=6cm]

\draw
	(0,-2.5) -- ++(0,0.4);

\draw[line width=1.2]
	(2,-1.8) -- ++(-0.4,0)
	(-2,0) -- ++(0.4,0);

\node at (1.7,0.65) {\small $\gamma^q$};
\node at (1.8,0) {\small $\alpha$};
\node at (1.8,-0.9) {\small $\beta$};
\node at (1.8,-1.6) {\small $\delta$};
\node at (1.8,-2) {\small $\epsilon$};

\node at (-1.8,0.7) {\small $\beta$};
\node at (-1.8,0.2) {\small $\delta$};
\node at (-1.8,-0.2) {\small $\epsilon$};
\node at (-1.7,-0.9) {\small $\gamma^q$};
\node at (-1.8,-1.8) {\small $\alpha$};

\node at (-1,-2.4) {\small $\beta$};
\node at (1,-2.4) {\small $\gamma$};

\end{scope}

\end{tikzpicture}
\caption{Proposition \ref{ade_2bc}: Rotation modification.}
\label{ade_2bcB}
\end{figure}

The case $\gamma_2\cdots=\beta\gamma^{q+1}$ is the second picture. The $\gamma^{q+1}$ part of the vertex induces a partial earth map tiling ${\mc P}^{q+1}_q$, which consists of $P,T_1,T_2$, and is exactly the part between $l_1,l_2$. Therefore the whole tiling is again obtained by gluing two copies of ${\mc P}^{q+1}_q$ together. The tiling is the rotation modification (by $240^{\circ}$) $RE_{\pentagon}2$ in Figure \ref{emt2mod2}, and is actually the same as the previous rotation modification.

\medskip

\noindent{\em Geometry of Pentagon}

\medskip

We discuss the geometrical existence of simple pentagon for the tilings. Figure \ref{ade_2bcC} gives two extreme cases. In the first case, the pentagon is the union of two isosceles triangles, with respective top angles $\beta=\pi-\theta$ and $\gamma=2\theta$, where $\theta=\frac{4}{f}\pi$. Given any length $a\in (0,\frac{1}{2}\pi]$, we construct the two isosceles triangles. The area of the two isosceles triangles is strictly increasing in $a$. As $a\to 0^+$, the area converges to $0$. For $a=\frac{1}{2}\pi$, the area is $(\pi-\theta)+2\theta=\pi+\theta>\theta$. Therefore there is unique $a_1(f)$, such that the corresponding area is exactly $\theta=\frac{4}{f}\pi$. Now for $a$ slightly smaller than $a_1(f)$, the area is slight smaller than $\theta$. Then we may add a small third triangle to the bottoms of the two isosceles triangles, so that the total area of the simple pentagon is exactly $\theta$. The pentagon is suitable for tiling.

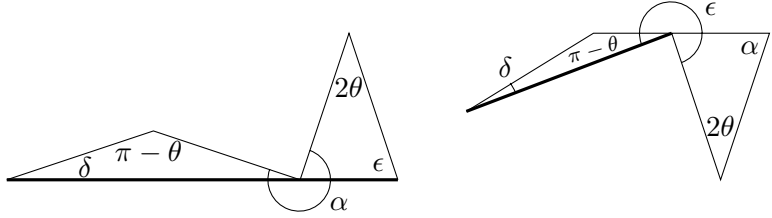
\begin{figure}[htp]
\centering
\begin{tikzpicture}[>=latex, scale=1.3]

\begin{scope}[shift={(-6cm,-1.5cm)}]

\draw
	(0,0) -- (1.5,0.5) -- (3,0) -- (3.5,1.5) -- (4,0)
	(3.1,0.3) arc (70:-200:0.32);

\draw[line width=1.2]
	(0,0) -- (4,0);

\node at (1.45,0.28) {\small $\pi-\theta$};
\node at (3.5,0.95) {\small $2\theta$};
\node at (3.4,-0.25) {\small $\alpha$};
\node at (3.8,0.15) {\small $\epsilon$};
\node at (0.8,0.13) {\small $\delta$};

\end{scope}

\draw
	(-1.3,-0.8) -- (0,0) -- (1.8,0) -- (1.3,-1.5) -- (0.8,0)
	(0.9,-0.3) arc (-70:200:0.32)
	(-0.8,-0.6) arc (25:35:0.6);
	
\draw[line width=1.2]
	(-1.3,-0.8) -- (0.8,0);
	
\node[rotate=15] at (0,-0.17) {\scriptsize $\pi-\theta$};	
\node at (1.3,-0.95) {\small $2\theta$};
\node at (1.6,-0.15) {\small $\alpha$};
\node at (-0.9,-0.35) {\small $\delta$};
\node at (1.2,0.25) {\small $\epsilon$};

\end{tikzpicture}
\caption{Proposition \ref{ade_2bc}: Two extreme cases, where $\theta=\frac{4}{f}\pi$.}
\label{ade_2bcC}
\end{figure}

The lower bound on $a$ is given by another extreme case, in the second of Figure \ref{ade_2bcC}. Again there is a unique $a_0(f)$, such that the $\alpha+\delta+\epsilon=2\pi$, and the area of the two triangles is $\theta=\frac{4}{f}\pi$. Then the tiling is parameterised by $a\in (a_0(f),a_1(f))$.\end{proof}

\begin{proposition}\label{ade_a2b}
There is no tiling, such that $\alpha,\beta,\gamma$ have distinct values, and $\alpha\delta\epsilon,\alpha\beta^2$ are vertices.
\end{proposition}

\begin{proof}
The angle sums of $\alpha\delta\epsilon,\alpha\beta^2$ and the angle sum for pentagon imply
\[
\beta+\gamma=(1+\tfrac{4}{f})\pi,\,
\delta+\epsilon=2\beta.
\]
We have $2\gamma+\delta+\epsilon=2(\beta+\gamma)>2\pi$. Then by $\alpha\delta\epsilon$, this implies $\alpha<2\gamma$.

By $\alpha\delta\epsilon$ and Proposition \ref{ade_2bc}, we know $\beta^2\gamma,\beta\gamma^2$ are not vertices. Then by $\alpha\delta\epsilon,\alpha\beta^2$ and distinct $\alpha,\beta,\gamma$ values, degree $3$ vertices besides $\alpha\delta\epsilon,\alpha\beta^2$ are $\alpha^2\gamma,\gamma^3,\beta\delta^2,\beta\epsilon^2,\gamma\delta^2,\gamma\epsilon^2$. 

By Lemma \ref{geometry1}, we divide the proof into the case $\beta>\gamma$ and $\delta<\epsilon$, and the case $\beta<\gamma$ and $\delta>\epsilon$.

\subsubsection*{Case. $\beta>\gamma$ and $\delta<\epsilon$}

By $\delta+\epsilon=2\beta$ and $\delta<\epsilon$, we get $\delta<\beta<\epsilon$. We also have $\delta+\epsilon=2\beta>\beta+\gamma>\pi$.

By $\delta+\epsilon>\pi$ and $\delta<\epsilon$, we know $R(\epsilon^2)$ has no $\delta,\epsilon$, and $R(\delta\epsilon)$ has no $\epsilon$. By $\alpha\delta\epsilon$ and $\delta<\epsilon$, we get $R(\epsilon^2)<R(\delta\epsilon)=\alpha<2\gamma$. Then by $\beta>\gamma$ and distinct $\alpha,\beta,\gamma$ values, we get $\epsilon\cdots=\alpha\delta\epsilon,\beta\epsilon^2,\gamma\epsilon^2,\delta^l\epsilon$. This implies $\epsilon\thin\epsilon\cdots$ is not a vertex. 

If $\epsilon^2\cdots$ is not a vertex, then $\epsilon\cdots=\alpha\delta\epsilon,\delta^l\epsilon$. Applying the counting lemma to $\delta,\epsilon$, we get $\epsilon\cdots=\alpha\delta\epsilon$. Then by $\alpha\beta^2$ and applying the counting lemma to $\alpha,\epsilon$, we get a contradiction. 

We conclude $\epsilon^2\cdots=\beta\epsilon^2,\gamma\epsilon^2$ is a vertex. By the balance lemma, this implies $\delta^2\cdots$ is a vertex. By the list $\epsilon\cdots$, we know $\delta^2\cdots$  has no $\epsilon$.

\subsubsection*{Subcase. $\beta\epsilon^2$ is a vertex}

By $\beta\epsilon^2$ and $\beta<\epsilon$, we get $\beta<\frac{2}{3}\pi<\epsilon$. By $\alpha\beta^2$ and $\beta+\gamma=(1+\tfrac{4}{f})\pi$, this implies $\alpha>\frac{2}{3}\pi$, and $\gamma>(\frac{1}{3}+\frac{4}{f})\pi$. Then we get $\beta<\alpha,2\gamma$. 

The AAD $\thick^{\delta}\epsilon^{\gamma}\thin^{\alpha}\beta^{\delta}\thin^{\gamma}\epsilon^{\delta}\thick$ of $\beta\epsilon^2$ implies a vertex $\alpha\gamma\cdots$. By $\alpha\delta\epsilon,\beta\epsilon^2$, we get $\alpha+\delta=\beta+\epsilon$. Then $\alpha+\gamma+2\delta=\beta+\gamma+\delta+\epsilon>2\gamma+\delta+\epsilon>2\pi$. By $\delta<\epsilon$, this implies $\alpha\gamma\cdots$ is a $\hat{b}$-vertex. Then by $\alpha\beta^2$, and $\alpha+\gamma>\beta+\gamma>\pi$, and $\gamma<\beta<2\gamma$, we get $\alpha\gamma\cdots=\alpha^2\gamma,\alpha\gamma^3$. The  angle sums of $\alpha\delta\epsilon,\alpha\beta^2,\beta\epsilon^2$, and the angle sum of one of $\alpha^2\gamma,\alpha\gamma^3$, and the angle sum for pentagon imply
\begin{align*}
\alpha^2\gamma &\colon 
	\alpha=(\tfrac{4}{5}-\tfrac{8}{5f})\pi,\,
	\beta=(\tfrac{3}{5}+\tfrac{4}{5f})\pi,\,
	\gamma=(\tfrac{2}{5}+\tfrac{16}{5f})\pi,\,
	\delta=(\tfrac{1}{2}+\tfrac{2}{f})\pi,\,
	\epsilon=(\tfrac{7}{10}-\tfrac{2}{5f})\pi. \\
\alpha\gamma^3 &\colon 
	\alpha=(\tfrac{4}{5}-\tfrac{24}{5f})\pi,\,
	\beta=(\tfrac{3}{5}+\tfrac{12}{5f})\pi,\,
	\gamma=(\tfrac{2}{5}+\tfrac{8}{5f})\pi,\,
	\delta=(\tfrac{1}{2}+\tfrac{6}{f})\pi,\,
	\epsilon=(\tfrac{7}{10}-\tfrac{6}{5f})\pi. 
\end{align*}
We know $\delta^2\cdots$ is a vertex, and has no $\epsilon$. The angle values above imply  $\delta^2\cdots=\gamma^2\delta^2$. The angle sum of $\gamma^2\delta^2$ further implies
\begin{align*}
\alpha^2\gamma &\colon 
	\alpha=\tfrac{10}{13}\pi,\,
	\beta=\tfrac{8}{13}\pi,\,
	\gamma=\tfrac{6}{13}\pi,\,
	\delta=\tfrac{7}{13}\pi,\,
	\epsilon=\tfrac{9}{13}\pi,\,
	f=52. \\
\alpha\gamma^3 &\colon 
	\alpha=\tfrac{14}{19}\pi,\,
	\beta=\tfrac{12}{19}\pi,\,
	\gamma=\tfrac{8}{19}\pi,\,
	\delta=\tfrac{11}{19}\pi,\,
	\epsilon=\tfrac{13}{19}\pi,\,
	f=76. 
\end{align*}

For $\alpha^2\gamma$, by the angle values, we get $\alpha\beta\cdots=\alpha\beta^2$, and $\alpha\gamma\cdots=\alpha^2\gamma$, and $\beta\gamma\cdots=\beta\gamma^3$, and $\gamma^2\cdots=\beta\gamma^3,\gamma^2\delta^2$, and $\beta\epsilon\cdots=\epsilon^2\cdots=\beta\epsilon^2$. Moreover, by no $\epsilon\thin\epsilon\cdots$, we know the AAD of $\gamma^2\delta^2$ is $\thick^{\epsilon}\delta^{\beta}\thin^{\epsilon}\gamma^{\alpha}\thin\gamma\thin^{\beta}\delta^{\epsilon}\thick$. This determines $T_1,T_2,T_3,\gamma_4$ in the first of Figure \ref{ade_a2bA}. Then $\beta_1\epsilon_3\cdots=\beta\epsilon^2$ determines $T_5$. Then $\alpha_1\gamma_5\cdots=\alpha^2\gamma$ and $\epsilon_1\epsilon_2\cdots=\beta\epsilon^2$ give $\alpha_6,\beta_7$. Then by $\beta\gamma\cdots=\beta\gamma^3$ and $\alpha\gamma\cdots=\alpha^2\gamma$, we determine $T_6,T_7$. Then $\gamma_1\gamma_6\delta_7\cdots=\gamma^2\delta^2$ determines $T_8$. Then $\alpha_7\gamma_2\cdots=\alpha^2\gamma$ gives $\alpha_9$. If $T_9$ is not arranged as indicated, then $\gamma_7\cdots=\beta\gamma\cdots=\beta\gamma^3$ and $\epsilon_7\epsilon_8\cdots=\beta\epsilon^2$ imply $\beta,\gamma$ adjacent, a contradiction. Therefore $T_9$ is arranged as indicated. Then $\alpha_2\beta_9\cdots=\alpha\beta^2$ gives $\beta$ just outside $\alpha_2$. On the other hand, $\gamma_4$ implies $\beta_2\cdots=\alpha\beta\cdots=\alpha\beta^2$ or $\beta_2\cdots=\beta\epsilon\cdots=\beta\epsilon^2$. The first case implies two $\beta$ adjacent, a contradiction. The second case implies $\beta,\epsilon$ adjacent, also a contradiction. 

\begin{figure}[htp]
\centering
\begin{tikzpicture}[>=latex]


\foreach \a in {1,-1}
{
\begin{scope}[xscale=\a]

\draw
	(0,0.7) -- (0,-0.7) -- (0.5,-1.1) -- (1,-0.7) -- (1,0.7) -- (0.5,1.1) -- (0,0.7);

\node at (0.8,0.6) {\small $\alpha$};	
\node at (0.5,0.85) {\small $\gamma$};	
\node at (0.8,0.2) {\small $\beta$};
\node at (0.2,0.6) {\small $\epsilon$};	
\node at (0.2,0.2) {\small $\delta$};

\end{scope}
}

\draw
	(-1,0) -- (1,0)
	(-1,0.7) -- (-1.8,0.7) 
	(-0.5,1.1) -- (-1.3,1.5) -- (-1.3,2.5) -- (-0.5,2.5) -- (-0.5,1.8)
	(-1.3,1.5) -- (-1.8,1.5) -- (-1.8,-0.7) -- (-1,-0.7)
	(-0.5,1.1) -- (-0.5,1.8) -- (1.5,1.8) -- (1.5,1.1) -- (1,0.7)
	(0.5,1.1) -- (0.5,1.8)
	;

\draw[line width=1.2]
	(0,0.7) -- (0,0)
	(-1,-0.7) -- (-1,0)
	(-0.5,1.1) -- (-0.5,1.8)
	(1.5,1.1) -- (1.5,1.8);
	
\node at (0.2,-0.2) {\small $\gamma$};

\node at (-0.8,-0.6) {\small $\delta$};	
\node at (-0.5,-0.85) {\small $\beta$};	
\node at (-0.8,-0.2) {\small $\epsilon$};
\node at (-0.2,-0.6) {\small $\alpha$};	
\node at (-0.2,-0.2) {\small $\gamma$};

\node at (-1.6,0.5) {\small $\alpha$};
\node at (-1.6,-0.5) {\small $\beta$};
\node at (-1.2,0.5) {\small $\gamma$};
\node at (-1.2,-0.5) {\small $\delta$};
\node at (-1.2,0) {\small $\epsilon$};

\node at (0.3,1.2) {\small $\alpha$};
\node at (0,0.95) {\small $\beta$};
\node at (0.3,1.6) {\small $\gamma$};
\node at (-0.3,1.2) {\small $\delta$};
\node at (-0.3,1.6) {\small $\epsilon$};

\node at (-0.3,2) {\small $\beta$};

\node at (1.3,1.2) {\small $\delta$};
\node at (1,0.95) {\small $\beta$};
\node at (1.3,1.6) {\small $\epsilon$};
\node at (0.7,1.2) {\small $\alpha$};
\node at (0.7,1.6) {\small $\gamma$};

\node at (1.2,0.6) {\small $\beta$};

\node at (-1.1,0.85) {\small $\alpha$};
\node at (-1.6,0.9) {\small $\beta$};
\node at (-0.8,1.05) {\small $\gamma$};
\node at (-1.6,1.3) {\small $\delta$};
\node at (-1.3,1.35) {\small $\epsilon$};

\node at (-1.1,2.35) {\small $\alpha$};
\node at (-1.1,1.65) {\small $\beta$};
\node at (-0.7,2.3) {\small $\gamma$};
\node at (-0.7,1.4) {\small $\delta$};
\node at (-0.7,1.8) {\small $\epsilon$};

\node[inner sep=0.5,draw,shape=circle] at (-0.5,0.4) {\small 1};
\node[inner sep=0.5,draw,shape=circle] at (0.5,0.4) {\small 2};
\node[inner sep=0.5,draw,shape=circle] at (-0.5,-0.4) {\small 3};
\node[inner sep=0.5,draw,shape=circle] at (0.5,-0.4) {\small 4};
\node[inner sep=0.5,draw,shape=circle] at (-1.5,0) {\small 5};
\node[inner sep=0.5,draw,shape=circle] at (-1.3,1.05) {\small 6};
\node[inner sep=0.5,draw,shape=circle] at (-0.9,2.05) {\small 8};
\node[inner sep=0.5,draw,shape=circle] at (0,1.4) {\small 7};
\node[inner sep=0.5,draw,shape=circle] at (1,1.4) {\small 9};


\begin{scope}[xshift=4.2cm]

\foreach \a in {1,-1}
{
\begin{scope}[xscale=\a]

\draw
	(2,0) -- (0,0) -- (0,-0.7) -- (0.5,-1.1) -- (1,-0.7) -- (1.5,-1.1) -- (2,-0.7) -- (2,0.7) -- (1.5,1.1) -- (1,0.7) -- (1,-0.7)
	(0,0.7) -- (1,0.7);

\node at (0.8,-0.6) {\small $\alpha$};	
\node at (0.5,-0.85) {\small $\beta$};	
\node at (0.8,-0.2) {\small $\gamma$};
\node at (0.2,-0.6) {\small $\delta$};	
\node at (0.2,-0.2) {\small $\epsilon$};

\node at (1.8,-0.2) {\small $\alpha$};
\node at (1.8,-0.6) {\small $\beta$};	
\node at (1.2,-0.2) {\small $\gamma$};
\node at (1.5,-0.85) {\small $\delta$};
\node at (1.2,-0.6) {\small $\epsilon$};

\end{scope}
}

\draw[line width=1.2]
	(0,0) -- (0,-0.7)
	(-1.5,-1.1) -- (-1,-0.7)
	(1.5,-1.1) -- (1,-0.7)
	(1,0) -- (1,0.7)
	(-2,0) -- (-2,0.7);

\node at (-1.2,0.6) {\small $\alpha$};
\node at (-1.52,0.85) {\small $\beta$};
\node at (-1.2,0.2) {\small $\gamma$};
\node at (-1.8,0.6) {\small $\delta$};	
\node at (-1.8,0.2) {\small $\epsilon$};

\node at (1.8,0.6) {\small $\alpha$};
\node at (1.8,0.2) {\small $\beta$};
\node at (1.5,0.85) {\small $\gamma$};	
\node at (1.2,0.2) {\small $\delta$};
\node at (1.2,0.6) {\small $\epsilon$};	

\node at (1,0.9) {\small $\beta$};
\node at (-0.9,0.95) {\small $\gamma^2$};

\node at (-0.8,0.2) {\small $\alpha$};
\node at (0,0.2) {\small $\beta$};
\node at (-0.8,0.5) {\small $\gamma$};
\node at (0.8,0.2) {\small $\delta$};
\node at (0.8,0.5) {\small $\epsilon$};

\node[inner sep=0.5,draw,shape=circle] at (-0.5,-0.4) {\small 1};
\node[inner sep=0.5,draw,shape=circle] at (0.5,-0.4) {\small 2};
\node[inner sep=0.5,draw,shape=circle] at (0.4,0.35) {\small 3};
\node[inner sep=0.5,draw,shape=circle] at (-1.5,-0.4) {\small 4};
\node[inner sep=0.5,draw,shape=circle] at (-1.5,0.4) {\small 5};
\node[inner sep=0.5,draw,shape=circle] at (1.5,0.4) {\small 6};

\end{scope}

\end{tikzpicture}
\caption{Proposition \ref{ade_a2b}: $\beta\epsilon^2$ is a vertex.}
\label{ade_a2bA}
\end{figure}
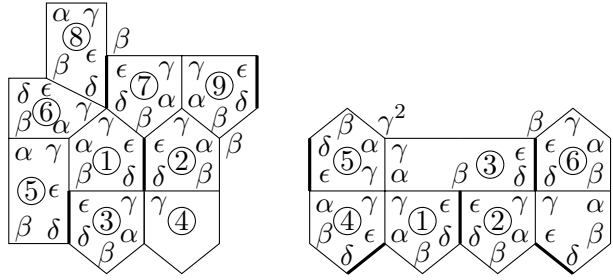

For $\alpha\gamma^3$, by the angle values, we get $\alpha\gamma\cdots=\alpha\gamma^3$, and $\gamma\delta\cdots=\gamma^2\delta^2$, and $\epsilon^2\cdots=\beta\epsilon^2$. We also know $\alpha^2\cdots$ is not a vertex. The AAD $\thick^{\delta}\epsilon^{\gamma}\thin^{\alpha}\beta^{\delta}\thin^{\gamma}\epsilon^{\delta}\thick$ of $\beta\epsilon^2$ determines $T_1,T_2,T_3$ in the second of Figure \ref{ade_a2bA}. Then $\alpha_3\gamma_1\cdots=\alpha\gamma^3$ and no $\alpha^2\cdots$ determine $T_4,T_5$, and $\gamma_2\delta_3\cdots=\gamma^2\delta^2$ determines $T_6$. Then $\alpha_5\gamma_3\cdots=\alpha\gamma^3$ and $\epsilon_3\epsilon_6\cdots=\beta\epsilon^2$ imply $\beta,\gamma$ adjacent, a contradiction. 

\subsubsection*{Subcase. $\gamma\epsilon^2$ is a vertex}

By $\gamma\epsilon^2$, and $\beta>\gamma$, and $\delta<\epsilon$, we know $\beta\epsilon^2,\gamma\delta^2$ are not vertices. In particular, $\beta\epsilon\cdots=\beta\epsilon^2$ is not a vertex. Moreover, the only degree $3$ vertices besides $\alpha\delta\epsilon,\alpha\beta^2,\gamma\epsilon^2$ are $\alpha^2\gamma,\beta\delta^2,\gamma^3$. The vertex $\gamma^3$ implies $\gamma=\tfrac{2}{3}\pi$. By $\beta+\gamma=(1+\tfrac{4}{f})\pi$, this implies $\beta<\gamma$, a contradiction. Therefore $\alpha^2\gamma,\beta\delta^2$ are the only possible extra degree $3$ vertices.

\subsubsection*{Subsubcase. $\alpha^2\gamma$ is a vertex}

The angle sums of $\alpha\delta\epsilon,\alpha\beta^2,\gamma\epsilon^2,\alpha^2\gamma$ and the angle sum for pentagon imply 
\[
\alpha=\epsilon=(\tfrac{4}{5}-\tfrac{8}{5f})\pi,\,
\beta=(\tfrac{3}{5}+\tfrac{4}{5f})\pi,\,
\gamma=\delta=(\tfrac{2}{5}+\tfrac{16}{5f})\pi.
\] 
By $\alpha\epsilon^2$ and the balance lemma, we know $\delta^2\cdots$ is a vertex. Then by the angle values, we get $\delta^2\cdots=\beta\gamma\delta^2,\gamma^2\delta^2,\delta^4$. The AAD of $\gamma^2\delta^2$ implies $\beta\epsilon\cdots,\epsilon\thin\epsilon\cdots$, a contradiction. 

The angle sum of $\beta\gamma\delta^2$ further implies
\[
\alpha=\epsilon=\tfrac{10}{13}\pi,\,
\beta=\tfrac{8}{13}\pi,\,
\gamma=\delta=\tfrac{6}{13}\pi,\,
f=52.
\]
By the angle values, we get $\beta\gamma\cdots=\beta\gamma^3,\beta\gamma\delta^2$, and $\beta\delta\cdots=\beta\gamma\delta^2$, and $\delta\thin\epsilon\cdots$ is not a vertex. By no $\beta\epsilon\cdots,\delta\thin\epsilon\cdots$, we know the AAD of $\beta\gamma\delta^2$ is $\thick^{\epsilon}\delta^{\beta}\thin^{\delta}\beta^{\alpha}\thin^{\epsilon}\gamma^{\alpha}\thin^{\beta}\delta^{\epsilon}\thick$. The AAD determines $T_1,T_2,T_3,T_4$ in the first of Figure \ref{ade_a2bB}. Then $\beta_1\delta_2\cdots=\beta\gamma\delta^2$, and the AAD of $\beta\gamma\delta^2$ determines $T_5,T_6$. Then $\alpha_1\epsilon_5\cdots=\alpha\delta\epsilon$ determines $T_7$. Then $\epsilon_1\epsilon_3\cdots=\gamma\epsilon^2$ and $\beta_7\gamma_1\cdots=\beta\gamma^3,\beta\gamma\delta^2$ imply either two $\gamma$ adjacent or $\gamma,\delta$ adjacent. Both are contradictions.

The angle sum of $\delta^4$ further implies
\[
\alpha=\epsilon=\tfrac{3}{4}\pi,\,
\beta=\tfrac{5}{8}\pi,\,
\gamma=\delta=\tfrac{1}{2}\pi,\,
f=32.
\]
By the angle values, we get $\alpha\gamma\cdots=\alpha\gamma^2$, and $\beta^2\cdots=\alpha\beta^2$, and $\beta\gamma\cdots$ is not a vertex. The AAD of $\delta^4$ implies a vertex $\thin^{\alpha}\beta^{\delta}\thin^{\delta}\beta^{\alpha}\thin\cdots=\alpha\beta^2=\thin^{\delta}\beta^{\alpha}\thin^{\beta}\alpha^{\gamma}\thin^{\alpha}\beta^{\delta}\thin$. This further implies a vertex $\thin^{\gamma}\alpha^{\beta}\thin^{\alpha}\gamma^{\epsilon}\thin\cdots=\alpha^2\gamma=\thin^{\alpha}\gamma^{\epsilon}\thin\alpha\thin^{\gamma}\alpha^{\beta}\thin$. Then we get $\beta\gamma\cdots,\beta\epsilon\cdots$, a contradiction.

\begin{figure}[htp]
\centering
\begin{tikzpicture}[>=latex,scale=1]


\begin{scope}[yshift=-0.9cm]

\foreach \a in {1,-1}
{
\begin{scope}[scale=\a]

\foreach \b in {0,1,-1}
\draw[xshift=\b cm]
	(-0.5,0) -- (-0.5,0.7) -- (0,1.1) -- (0.5,0.7) -- (0.5,0) -- (-0.5,0);
	
\draw[line width=1.2]
	(-0.5,0) -- (-0.5,0.7)
	(0.5,0.7) -- (1,1.1);

\node at (0.3,0.6) {\small $\alpha$};		
\node at (0.3,0.2) {\small $\beta$};
\node at (0,0.85) {\small $\gamma$};
\node at (-0.3,0.2) {\small $\delta$};
\node at (-0.3,0.6) {\small $\epsilon$};	

\node at (1.3,0.6) {\small $\beta$};		
\node at (1.3,0.2) {\small $\alpha$};
\node at (1,0.85) {\small $\delta$};
\node at (0.7,0.2) {\small $\gamma$};
\node at (0.7,0.6) {\small $\epsilon$};

\node at (-1.3,0.6) {\small $\alpha$};		
\node at (-1.3,0.2) {\small $\beta$};
\node at (-1,0.85) {\small $\gamma$};
\node at (-0.7,0.2) {\small $\delta$};
\node at (-0.7,0.6) {\small $\epsilon$};

\end{scope}
}

\draw
	(0,1.1) -- (0,1.8) -- (1,1.8) -- (1,1.1) -- (0.5,0.7); 

\node at (0.8,1.2) {\small $\epsilon$};
\node at (0.5,0.95) {\small $\delta$};
\node at (0.8,1.6) {\small $\gamma$};
\node at (0.2,1.2) {\small $\beta$};
\node at (0.2,1.6) {\small $\alpha$};

\node at (-0.5,0.95) {\small $\gamma$};

\node[inner sep=0.5,draw,shape=circle] at (0,0.4) {\small 1};
\node[inner sep=0.5,draw,shape=circle] at (0,-0.4) {\small 2};
\node[inner sep=0.5,draw,shape=circle] at (-1,0.4) {\small 3};
\node[inner sep=0.5,draw,shape=circle] at (-1,-0.4) {\small 4};
\node[inner sep=0.5,draw,shape=circle] at (1,0.4) {\small 5};
\node[inner sep=0.5,draw,shape=circle] at (1,-0.4) {\small 6};
\node[inner sep=0.5,draw,shape=circle] at (0.5,1.4) {\small 7};

\end{scope}

\begin{scope}[xshift=3.5cm]

\foreach \a in {-1,0,1,2}
\draw[xshift=\a cm]
	(0,0.9) -- (0,0.2) -- (-0.5,-0.2) -- (-0.5,-0.9);

\foreach \a in {-1,0,1}
\draw[xshift=\a cm]
	(0,0.2) -- (0.5,-0.2);
	
\draw
	(2,0.9) -- (-1,0.9) -- (-1,1.6) -- (1,1.6) -- (1,0.9)
	(-1.5,-0.9) -- (1.5,-0.9);

\draw[line width=1.2]
	(0,0.2) -- (0,0.9)
	(-0.5,-0.2) -- (-0.5,-0.9)
	(1,0.2) -- (1.5,-0.2)
	(-1,0.9) -- (-1,1.6);

\node at (0.2,0.3) {\small $\delta$};	
\node at (0.5,0.05) {\small $\beta$};
\node at (0.2,0.7) {\small $\epsilon$};	
\node at (0.8,0.3) {\small $\alpha$};
\node at (0.8,0.7) {\small $\gamma$}; 	

\node at (-0.2,0.3) {\small $\delta$};	
\node at (-0.5,0.05) {\small $\beta$};
\node at (-0.2,0.7) {\small $\epsilon$};	
\node at (-0.8,0.3) {\small $\alpha$};
\node at (-0.8,0.7) {\small $\gamma$}; 
	
\node at (-1.2,0.3) {\small $\beta$};
\node at (-1.2,0.9) {\small $\epsilon$};

\node at (1.2,0.3) {\small $\epsilon$};	
\node at (1.5,0.05) {\small $\delta$};
\node at (1.2,0.7) {\small $\gamma$};	
\node at (1.8,0.3) {\small $\beta$};
\node at (1.8,0.7) {\small $\alpha$};

\node at (1.3,-0.3) {\small $\epsilon$};	
\node at (1,-0.05) {\small $\delta$};
\node at (1.3,-0.7) {\small $\gamma$}; 
\node at (0.7,-0.3) {\small $\beta$};
\node at (0.7,-0.7) {\small $\alpha$};
		
\node at (-0.3,-0.3) {\small $\delta$};	
\node at (0,-0.05) {\small $\beta$};
\node at (-0.3,-0.7) {\small $\epsilon$};
\node at (0.3,-0.3) {\small $\alpha$};
\node at (0.3,-0.7) {\small $\gamma$}; 

\node at (-0.7,-0.3) {\small $\delta$};	
\node at (-1,-0.05) {\small $\beta$};
\node at (-0.7,-0.7) {\small $\epsilon$};
\node at (-1.3,-0.3) {\small $\alpha$};
\node at (-1.3,-0.7) {\small $\gamma$}; 

\node at (-0.8,1.4) {\small $\delta$};
\node at (0.8,1.4) {\small $\beta$};
\node at (-0.8,1.1) {\small $\epsilon$};
\node at (0.8,1.1) {\small $\alpha$};
\node at (0,1.1) {\small $\gamma$};

\node[inner sep=0.5,draw,shape=circle] at (1,-0.5) {\small $5$};
\node[inner sep=0.5,draw,shape=circle] at (0.5,0.5) {\small $2$};
\node[inner sep=0.5,draw,shape=circle] at (-0.5,0.5) {\small $1$};
\node[inner sep=0.5,draw,shape=circle] at (1.5,0.5) {\small $6$};
\node[inner sep=0.5,draw,shape=circle] at (0,-0.5) {\small $3$};
\node[inner sep=0.5,draw,shape=circle] at (-1,-0.5) {\small $4$};
\node[inner sep=0.5,draw,shape=circle] at (-0.4,1.25) {\small $7$};

\end{scope}


\begin{scope}[xshift=9.2cm]

\foreach \a in {-1,0,1,2}
\draw[xshift=\a cm]
	(0,0.9) -- (0,0.2) -- (-0.5,-0.2) -- (-0.5,-0.9);

\foreach \a in {-1,0,1}
\draw[xshift=\a cm]
	(0,0.2) -- (0.5,-0.2);
	
\draw
	(-1,0.9) -- (2,0.9) -- (2,1.6) -- (0,1.6) -- (0,0.9)
	(-1.5,-0.2) -- (-2,0.2) -- (-2.5,-0.2) -- (-2.5,-2.3) -- (2.3,-2.3) -- (2.3,-0.2) -- (1.5,-0.2)
	(-2.5,-0.9) -- (1.5,-0.9) -- (1.5,-1.6) -- (-3.3,-1.6) -- (-3.3,-0.2) -- (-2.5,-0.2)
	(-0.5,-0.9) -- (-0.5,-2.3);

\draw[line width=1.2]
	(0,0.2) -- (0.5,-0.2)
	(-1,0.2) -- (-0.5,-0.2)
	(1.5,-0.9) -- (0.5,-0.9)
	(-2.5,-0.9) -- (-1.5,-0.9)
	(-2.5,-2.3) -- (-1.5,-2.3)
	(2,0.9) -- (1,0.9)
	(1.5,-1.6) -- (2.3,-1.6)
	(-3.3,-1.6) -- (-3.3,-0.2);
		
\node at (1.3,-0.3) {\small $\gamma$};	
\node at (1,-0.05) {\small $\alpha$};
\node at (1.3,-0.7) {\small $\epsilon$}; 
\node at (0.7,-0.3) {\small $\beta$};
\node at (0.7,-0.7) {\small $\delta$};
		
\node at (-0.3,-0.3) {\small $\gamma$};	
\node at (0,-0.05) {\small $\epsilon$};
\node at (-0.3,-0.7) {\small $\alpha$};
\node at (0.3,-0.3) {\small $\delta$};
\node at (0.3,-0.7) {\small $\beta$}; 

\node at (-1.3,-0.3) {\small $\beta$};	
\node at (-1,-0.05) {\small $\delta$};
\node at (-1.3,-0.7) {\small $\alpha$};
\node at (-0.7,-0.3) {\small $\epsilon$};
\node at (-0.7,-0.7) {\small $\gamma$}; 

\node at (-2,-0.05) {\small $\alpha$};
\node at (-2.3,-0.3) {\small $\gamma$};
\node at (-1.7,-0.3) {\small $\beta$};
\node at (-2.3,-0.75) {\small $\epsilon$};
\node at (-1.7,-0.7) {\small $\delta$};

\node at (0.2,0.3) {\small $\epsilon$};	
\node at (0.5,0.05) {\small $\delta$};
\node at (0.2,0.7) {\small $\gamma$};	
\node at (0.8,0.3) {\small $\beta$};
\node at (0.8,0.7) {\small $\alpha$}; 	

\node at (-0.2,0.3) {\small $\gamma$};	
\node at (-0.5,0.05) {\small $\epsilon$};
\node at (-0.2,0.7) {\small $\alpha$};	
\node at (-0.8,0.3) {\small $\delta$};
\node at (-0.8,0.7) {\small $\beta$}; 
		
\node at (1.8,0.3) {\small $\gamma$};	
\node at (1.5,0.05) {\small $\alpha$};
\node at (1.8,0.7) {\small $\epsilon$}; 
\node at (1.2,0.3) {\small $\beta$};
\node at (1.2,0.7) {\small $\delta$};

\node at (0.2,1.4) {\small $\alpha$};
\node at (1.8,1.4) {\small $\beta$};
\node at (0.2,1.1) {\small $\gamma$};
\node at (1.8,1.1) {\small $\delta$};
\node at (1,1.1) {\small $\epsilon$};

\node at (-0.3,-1.1) {\small $\beta$}; 
\node at (-0.3,-1.4) {\small $\alpha$};
\node at (0.5,-1.1) {\small $\delta$}; 
\node at (1.3,-1.4) {\small $\gamma$};
\node at (1.3,-1.1) {\small $\epsilon$}; 

\node at (-0.7,-1.4) {\small $\alpha$};
\node at (-2.3,-1.4) {\small $\beta$};
\node at (-0.7,-1.1) {\small $\gamma$};
\node at (-2.3,-1.1) {\small $\delta$};
\node at (-1.5,-1.1) {\small $\epsilon$}; 

\node at (-2.65,-0.9) {\small $\alpha$}; 
\node at (-2.7,-0.4) {\small $\beta$}; 
\node at (-2.7,-1.4) {\small $\gamma$}; 
\node at (-3.1,-0.4) {\small $\delta$}; 
\node at (-3.1,-1.4) {\small $\epsilon$}; 

\node at (-0.7,-1.8) {\small $\alpha$};
\node at (-0.7,-2.1) {\small $\beta$};
\node at (-2.3,-1.8) {\small $\gamma$};
\node at (-1.5,-2.1) {\small $\delta$};
\node at (-2.3,-2.1) {\small $\epsilon$};

\node at (-2.7,-1.8) {\small $\alpha$};

\node at (-0.3,-2.1) {\small $\alpha$};
\node at (2.1,-2.1) {\small $\beta$};
\node at (-0.3,-1.8) {\small $\gamma$};
\node at (2.1,-1.8) {\small $\delta$};
\node at (1.5,-1.8) {\small $\epsilon$};

\node at (2.1,-0.4) {\small $\beta$}; 
\node at (2.1,-1.4) {\small $\delta$};
\node at (1.7,-0.4) {\small $\alpha$}; 
\node at (1.7,-1.4) {\small $\epsilon$};
\node at (1.7,-0.9) {\small $\gamma$}; 


\node[inner sep=0.5,draw,shape=circle] at (1,-0.5) {\small $1$};
\node[inner sep=0.5,draw,shape=circle] at (0.5,0.5) {\small $2$};
\node[inner sep=0.5,draw,shape=circle] at (1.5,0.5) {\small $3$};
\node[inner sep=0.5,draw,shape=circle] at (0,-0.5) {\small $4$};
\node[inner sep=0.5,draw,shape=circle] at (0.6,1.25) {\small $5$};
\node[inner sep=0.5,draw,shape=circle] at (0.9,-1.25) {\small $6$};
\node[inner sep=0.5,draw,shape=circle] at (-0.5,0.5) {\small $7$};
\node[inner sep=0.5,draw,shape=circle] at (2,-0.9) {\small $8$};
\node[inner sep=0.5,draw,shape=circle] at (-1,-0.5) {\small $9$};
\node[inner sep=0,draw,shape=circle] at (0.9,-1.95) {\footnotesize $10$};
\node[inner sep=0,draw,shape=circle] at (-1.1,-1.25) {\footnotesize $11$};
\node[inner sep=0,draw,shape=circle] at (-1.1,-1.95) {\footnotesize $13$};
\node[inner sep=0,draw,shape=circle] at (-2,-0.5) {\footnotesize $12$};
\node[inner sep=0,draw,shape=circle] at (-3,-0.9) {\footnotesize $14$};

\end{scope}

\end{tikzpicture}
\caption{Proposition \ref{ade_a2b}: $\gamma\epsilon^2$ is a vertex.}
\label{ade_a2bB}
\end{figure}
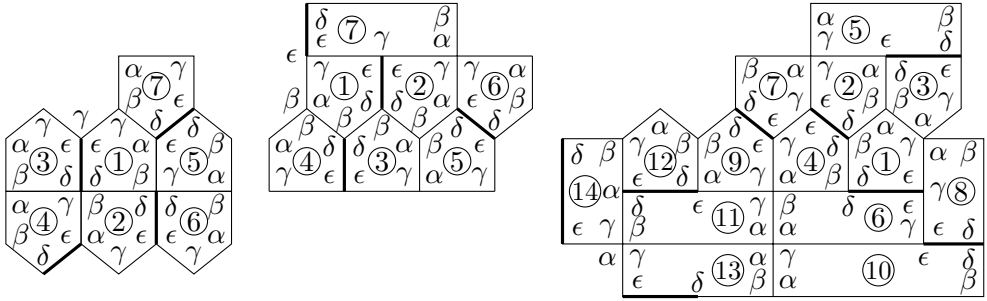

\subsubsection*{Subsubcase. $\beta\delta^2$ is a vertex}

The angle sums of $\alpha\delta\epsilon,\alpha\beta^2,\gamma\epsilon^2,\beta\delta^2$ and the angle sum for pentagon imply 
\[
\alpha=(\tfrac{1}{2}+\tfrac{2}{f})\pi,\,
\beta=(\tfrac{3}{4}-\tfrac{1}{f})\pi,\,
\gamma=(\tfrac{1}{4}+\tfrac{5}{f})\pi,\,
\delta=(\tfrac{5}{8}+\tfrac{1}{2f})\pi,\,
\epsilon=(\tfrac{7}{8}-\tfrac{5}{2f})\pi.
\]

By $\epsilon>\delta>\frac{1}{2}\pi$, we know a $b$-vertex has at most two from $\delta,\epsilon$. This implies that, if there is a $b$-vertex besides $\alpha\delta\epsilon,\gamma\epsilon^2,\beta\delta^2$, then one of $\alpha,\beta,\gamma$ is a combination of $\alpha,\beta,\gamma$. Guided by this, we find $\gamma^2\delta^2$ is the only possible such combination. However, the AAD of $\gamma^2\delta^2$ implies $\beta\epsilon\cdots,\epsilon\thin\epsilon\cdots$, a contradiction. Therefore $\alpha\delta\epsilon,\gamma\epsilon^2,\beta\delta^2$ are all the $b$-vertices.

We know $\alpha\beta\cdots$ is not a $b$-vertex. Then by the angle values, we get $\alpha\beta\cdots=\alpha\beta^2,\alpha\beta\gamma^2$. The angle sum of $\alpha\beta\gamma^2$ further implies 
\begin{equation}\label{ade_a2b_eq1}
\alpha=\tfrac{6}{11}\pi,\,
\beta=\tfrac{8}{11}\pi,\,
\gamma=\tfrac{4}{11}\pi,\,
\delta=\tfrac{7}{11}\pi,\,
\epsilon=\tfrac{9}{11}\pi,\,
f=44.
\end{equation}

Suppose $f\ne 44$. Then $\alpha\beta\cdots=\alpha\beta^2$. In the second of Figure \ref{ade_a2bB}, the AAD $\thick^{\epsilon}\delta^{\beta}\thin^{\alpha}\beta^{\delta}\thin^{\beta}\delta^{\epsilon}\thick$ of $\beta\delta^2$ determines $T_1,T_2,T_3$. Then $\beta_1\delta_3\cdots=\beta\delta^2$ determines $T_4$, and $\alpha_3\beta_2\cdots=\alpha\beta^2$ and no $\gamma\delta\cdots$ determine $T_5$. Then $\alpha_2\delta_5\cdots=\alpha\delta\epsilon$ determines $T_6$. Then $\epsilon_1\epsilon_2\cdots=\gamma\epsilon^2$ and no $\gamma^2\epsilon\cdots$ determine $T_7$. Then $\alpha_1\beta_4\cdots=\alpha\beta^2$ and $\gamma_1\epsilon_7\cdots=\gamma\epsilon^2$ imply $\beta,\epsilon$ adjacent, a contradiction.

It remains to consider the tiling for \eqref{ade_a2b_eq1}. The subsequent argument will be purely based on the angle values. 

The angle values imply the vertices are $\gamma^2\delta^2$ and the following
\[
\text{AVC}
=\{
\alpha\beta^2,\alpha\delta\epsilon,\beta\delta^2,\gamma\epsilon^2,\alpha\beta\gamma^2,\alpha^3\gamma,\alpha\gamma^4\}.
\]
In particular, we know $\beta\epsilon\cdots,\epsilon\thin\epsilon\cdots$ are not vertices. Then by the AAD argument as before, we know $\gamma^2\delta^2$ is not a vertex.

The AVC implies $\gamma\delta\cdots,\delta\thin\delta\cdots$ are not vertices. Therefore the AAD of $\beta^2\cdots=\alpha\beta^2$ is $\thin^{\beta}\alpha^{\gamma}\thin^{\alpha}\beta^{\delta}\thin^{\alpha}\beta^{\delta}\thin$. This implies the AAD of $\thin\beta\thin\beta\thin$ is $\thin^{\alpha}\beta^{\delta}\thin^{\alpha}\beta^{\delta}\thin$.

The AAD $\thin^{\beta}\alpha^{\gamma}\thin^{\alpha}\beta^{\delta}\thin^{\alpha}\beta^{\delta}\thin$ of $\alpha\beta^2$ determines $T_1,T_2,T_3$ in the third of Figure \ref{ade_a2bB}. Then $\beta_1\delta_2\cdots=\beta\delta^2$ determines $T_4$, and $\alpha_2\delta_3\cdots=\alpha\delta\epsilon$ determines $T_5$. Then $\beta_4\delta_1\cdots=\beta\delta^2$ determines $T_6$, and $\epsilon_2\epsilon_4\cdots=\gamma\epsilon^2$ and no $\gamma^2\epsilon\cdots$ determine $T_7$. Then $\epsilon_1\epsilon_6\cdots=\gamma\epsilon^2$ and no $\alpha\gamma\epsilon\cdots$ determine $T_8$, and $\gamma_4\epsilon_7\cdots=\gamma\epsilon^2$ determines $T_9$. Then $\gamma_6\epsilon_8\cdots=\gamma\epsilon^2$ determines $T_{10}$. Then $\alpha_4\beta_6\gamma_9\cdots=\alpha\beta\gamma^2$ and no $\alpha\gamma\epsilon\cdots$ determine $T_{11}$. Then $\alpha_9\epsilon_{11}\cdots=\alpha\delta\epsilon$ determines $T_{12}$, and $\alpha_6\alpha_{11}\gamma_{10}\cdots=\alpha^3\gamma$ and the AAD $\thin^{\alpha}\beta^{\delta}\thin^{\alpha}\beta^{\delta}\thin$ of $\thin\beta\thin\beta\thin$ determine $T_{13}$. Then $\delta_{11}\epsilon_{12}\cdots=\alpha\delta\epsilon$ and no $\beta^2\gamma\cdots$ determine $T_{14}$. Then $\beta_{11}\gamma_{13}\gamma_{14}\cdots=\alpha\beta\gamma^2$. This implies one of $\epsilon_{12}\cdots,\epsilon_{14}\cdots$ is $\beta\epsilon\cdots$, a contradiction. 

\subsubsection*{Subsubcase. $\alpha\delta\epsilon,\alpha\beta^2,\gamma\epsilon^2$ are all the degree $3$ vertices}

The angle sums of $\alpha\delta\epsilon,\alpha\beta^2,\gamma\epsilon^2$ and the angle sum for pentagon imply
\[
\alpha=2\gamma-\tfrac{8}{f}\pi,\,
\beta=(1+\tfrac{4}{f})\pi-\gamma,\,
\delta=(1+\tfrac{8}{f})\pi-\tfrac{3}{2}\gamma,\,
\epsilon=\pi-\tfrac{1}{2}\gamma.
\]
By $\alpha\delta\epsilon,\gamma\epsilon^2$, we get $\alpha+\delta=\gamma+\epsilon$. By $\delta<\epsilon$, this implies $\alpha>\gamma$. Moreover, we have $3\gamma+2\delta>\alpha+\gamma+2\delta=2\gamma+\delta+\epsilon>2\pi$. 

In a special tile, if $\gamma\cdots=\gamma\epsilon^2$, then $\alpha\cdots=\alpha\gamma\cdots$ has high degree. This implies $\delta\cdots,\epsilon\cdots$ have degree $3$, and $\delta\cdots=\epsilon\cdots=\alpha\delta\epsilon$. Then we get $\alpha,\epsilon$ adjacent, a contradiction. Therefore the vertex $\gamma\cdots$ has high degree. This implies $\delta\cdots=\epsilon\cdots=\alpha\delta\epsilon$. Then we get $\alpha\cdots=\beta\cdots=\alpha\beta^2$, as in the first of Figure \ref{ade_a2bC}. The AAD of the vertex $\beta\cdots$ is $\thin^{\beta}\alpha^{\gamma}\thin^{\delta}\beta^{\alpha}\thin^{\delta}\beta^{\alpha}\thin$, and implies a vertex $\gamma\delta\cdots$. Moreover, the vertex $H=\gamma\cdots$ of degree $4$ or $5$ is ${}^{\beta}\thin^{\alpha}\gamma^{\epsilon}\thin^{\alpha}\cdots=\thin^{\gamma}\alpha^{\beta}\thin^{\alpha}\gamma^{\epsilon}\thin^{\alpha}\cdots,\thick^{\epsilon}\delta^{\beta}\thin^{\alpha}\gamma^{\epsilon}\thin^{\alpha}\cdots$. The argument also shows there is no $3^5$-tile. By Lemma \ref{special_tile}, this implies $f\ge 24$. 

\begin{figure}[htp]
\centering
\begin{tikzpicture}[>=latex]


\draw
	(-0.5,-0.4) -- (-0.5,0.7) -- (0,1.1) -- (0.5,0.7) -- (0.5,-0.4)
	(0,1.1) -- ++(0,0.4)
	(-0.5,0.7) -- ++(-0.4,0);
		
\draw[line width=1.2]
	(-0.5,0) -- (0.5,0);

\node at (0,0.85) {\small $\alpha$}; 
\node at (-0.3,0.6) {\small $\beta$};
\node at (0.3,0.6) {\small $\gamma$};
\node at (-0.3,0.2) {\small $\delta$};
\node at (0.3,0.2) {\small $\epsilon$};	

\node at (0.7,0.7) {\small $H$};
\node at (-0.7,0.5) {\small $\beta$};
\node at (-0.6,0.9) {\small $\alpha$};
\node at (-0.7,0) {\small $\alpha$};
\node at (0.7,0) {\small $\alpha$};
\node at (0.3,-0.2) {\small $\delta$};
\node at (-0.3,-0.2) {\small $\epsilon$};	
\node at (-0.2,1.2) {\small $\beta$};
\node at (0.2,1.2) {\small $\beta$};


\begin{scope}[xshift=2.5cm]

\draw
	(0,0.7) -- (-0.5,1.1) -- (-1,0.7) -- (-1,-0.7) -- (-0.5,-1.1) -- (0,-0.7)
	(0,0.7) -- (1.6,0.7) -- (1.6,-0.7) -- (0,-0.7) -- cycle
	(0.8,0.7) -- (0.8,1.1)
	(-0.5,1.1) -- (-0.5,1.5);

\draw[line width=1.2]
	(0,0) -- (-1,0)
	(0.8,0.7) -- (0.8,-0.7)
	(0,0.7) -- (0.1,1.1);

\node at (-0.8,0.6) {\small $\beta$};
\node at (-0.8,0.2) {\small $\delta$};
\node at (-0.5,0.85) {\small $\alpha$};	
\node at (-0.2,0.2) {\small $\epsilon$};
\node at (-0.2,0.6) {\small $\gamma$};	

\node at (-0.2,-0.6) {\small $\beta$};
\node at (-0.2,-0.2) {\small $\delta$};
\node at (-0.5,-0.85) {\small $\alpha$};	
\node at (-0.8,-0.2) {\small $\epsilon$};
\node at (-0.8,-0.6) {\small $\gamma$};	

\node at (0.2,0) {\small $\alpha$};
\node at (0.2,0.5) {\small $\beta$};
\node at (0.2,-0.5) {\small $\gamma$};
\node at (0.6,0.5) {\small $\delta$};
\node at (0.6,-0.5) {\small $\epsilon$};

\node at (1.4,0) {\small $\alpha$};
\node at (1.4,0.5) {\small $\beta$};
\node at (1.4,-0.5) {\small $\gamma$};
\node at (1,0.5) {\small $\delta$};
\node at (1,-0.5) {\small $\epsilon$};

\node at (0.8,-0.9) {\small $\gamma$};

\node at (-0.1,1) {\small $\delta$};	
\node at (0.2,0.9) {\small $\delta$};	
\node at (0.6,0.9) {\small $\beta$};	
\node at (1,0.9) {\small $\gamma$};	

\node at (-0.3,1.2) {\small $\beta$};
\node at (-0.7,1.2) {\small $\beta$};

\node[inner sep=0.5,draw,shape=circle] at (-0.5,0.4) {\small 1};
\node[inner sep=0.5,draw,shape=circle] at (-0.5,-0.4) {\small 2};\node[inner sep=0.5,draw,shape=circle] at (0.55,0) {\small 3};\node[inner sep=0.5,draw,shape=circle] at (1.05,0) {\small 4};

\end{scope}

\end{tikzpicture}
\caption{Proposition \ref{ade_a2b}: $\alpha\delta\epsilon,\alpha\beta^2,\gamma\epsilon^2$ are all the degree $3$ vertices.}
\label{ade_a2bC}
\end{figure}
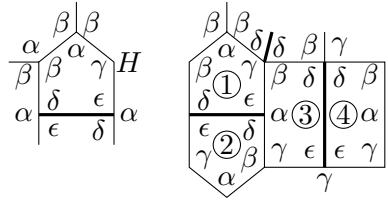

Suppose $H=\thin^{\gamma}\alpha^{\beta}\thin^{\alpha}\gamma^{\epsilon}\thin^{\alpha}\cdots=\alpha\beta\gamma\cdots,\alpha\gamma^2\cdots$. Then the AAD of the vertex $\alpha\cdots$ is $\thin^{\beta}\alpha^{\gamma}\thin^{\gamma}\beta^{\delta}\thin^{\delta}\beta^{\gamma}\thin$. This implies a vertex $\delta\thin\delta\cdots$.  By $\delta<\epsilon$, the vertex $\delta\thin\delta\cdots$ implies $\delta<\frac{1}{2}\pi$. This means $\gamma>(\frac{1}{3}+\frac{16}{3f})\pi$. Then $\alpha+\gamma=3\gamma-\frac{8}{f}\pi>\pi$, and $\alpha+4\gamma=6\gamma-\tfrac{8}{f}\pi>2\pi$. 

By $\alpha+\gamma+2\delta>2\pi$, and $\delta<\epsilon$, we know $H$ is a $\hat{b}$-vertex. By $\alpha,\beta>\gamma$, and $\alpha+\gamma>\pi$, and $\beta+\gamma>\pi$, and $\alpha+4\gamma>2\pi$, we get $H=\alpha\gamma^3$. The angle sum of $\alpha\gamma^3$ further implies
\[
\alpha=(\tfrac{4}{5}-\tfrac{24}{5f})\pi,\,
\beta=(\tfrac{3}{5}+\tfrac{12}{5f})\pi,\,
\gamma=(\tfrac{2}{5}+\tfrac{8}{5f})\pi,\,
\delta=(\tfrac{2}{5}+\tfrac{28}{5f})\pi,\,
\epsilon=(\tfrac{4}{5}-\tfrac{4}{5f})\pi.
\]

By $3\delta+\epsilon=(2+\tfrac{16}{f})\pi>2\pi$, and $\delta<\epsilon$, we get $\delta\thin\delta\cdots=\delta^4\cdots$. By $R(\delta^4)=(\tfrac{2}{5}-\tfrac{112}{5f})\pi<$ all angles, we get $\delta\thin\delta\cdots=\delta^4$. The angle sum of $\delta^4$ further implies
\[
\alpha=\tfrac{5}{7}\pi,\,
\beta=\tfrac{9}{14}\pi,\,
\gamma=\tfrac{3}{7}\pi,\,
\delta=\tfrac{1}{2}\pi,\,
\epsilon=\tfrac{11}{14}\pi,\,
f=56.
\]
Recall the AAD of the vertex $\beta\cdots$ in the special tile implies $\gamma\delta\cdots$ is a vertex. We have $R(\gamma\delta^2)=\frac{4}{7}\pi$ and $R(\gamma\delta\epsilon)=\frac{2}{7}\pi$. Since both are not sums of angles, we know $\gamma\delta\cdots$ is not a vertex, a contradiction. 

Suppose $H=\thick^{\epsilon}\delta^{\beta}\thin^{\alpha}\gamma^{\epsilon}\thin^{\alpha}\cdots=\beta\gamma\delta\cdots,\gamma^2\delta\cdots$. By $3\gamma+2\delta>2\gamma+\delta+\epsilon>2\pi$, and $\alpha,\beta>\gamma$, and $\delta<\epsilon$, we get $H=\beta\gamma\delta^2,\gamma^2\delta^2$. If $H=\gamma^2\delta^2$, then $H=\thick^{\epsilon}\delta^{\beta}\thin^{\alpha}\gamma^{\epsilon}\thin^{\alpha}\cdots=\thick^{\epsilon}\delta^{\beta}\thin^{\alpha}\gamma^{\epsilon}\thin^{\alpha}\gamma^{\epsilon}\thin^{\beta}\delta^{\epsilon}\thick$ implies a vertex $\beta\epsilon\cdots$, a contradiction. Therefore $H=\beta\gamma\delta^2$. The angle sum of $\beta\gamma\delta^2$ further implies
\[
	\alpha=(\tfrac{2}{3}+\tfrac{16}{3f})\pi,\,
	\beta=(\tfrac{2}{3}-\tfrac{8}{3f})\pi,\,
	\gamma=(\tfrac{1}{3}+\tfrac{20}{3f})\pi,\,
	\delta=(\tfrac{1}{2}-\tfrac{2}{f})\pi,\,
	\epsilon=(\tfrac{5}{6}-\tfrac{10}{3f})\pi. 
\]
We have $\alpha>\beta>\gamma$. Then by $\beta+\gamma,3\gamma>\pi$, we know a $\hat{b}$-vertex $\beta\gamma\cdots=\beta\gamma^3$.

By $f\ge 24$, we get $\gamma<2\delta,\epsilon$. Then by $\beta\gamma\delta^2$ and $\delta<\epsilon$, we get $R(\beta\delta\epsilon)<R(\beta\delta^2)=\gamma<\alpha,\beta,2\delta,\epsilon$. By no $\beta\epsilon\cdots$, this implies a $b$-vertex $\beta\cdots=\beta\gamma\delta^2$. Combined with $\hat{b}$-vertex $\beta\gamma\cdots=\beta\gamma^3$, we get $\beta\gamma\cdots=\beta\gamma^3,\beta\gamma\delta^2$.

By $f\ge 24$, we get $5\delta>2\pi$. Therefore $\delta^l\epsilon=\delta^3\epsilon$. The angle sum of $\delta^3\epsilon$ implies
\[
\alpha=\tfrac{6}{7}\pi,\,
\beta=\gamma=\tfrac{4}{7}\pi,\,
\delta=\tfrac{3}{7}\pi,\,
\epsilon=\tfrac{5}{7}\pi,\,
f=28.
\]
This contradicts Lemma \ref{geometry1}. Therefore $\delta^l\epsilon$ is not a vertex, and $\epsilon\cdots=\alpha\delta\epsilon,\gamma\epsilon^2$. In particular, we have $\delta\epsilon\cdots=\alpha\delta\epsilon$ and $\epsilon^2\cdots=\gamma\epsilon^2$.

The special tile in the first of Figure \ref{ade_a2bC} is the tile $T_1$ in the second picture, with $H=\beta\gamma\delta^2$. The special tile determines its companion $T_2$, and $H=\beta\gamma\delta^2$ and $\delta_2\epsilon_1\cdots=\alpha\delta\epsilon$ determine $T_3$. Then $H$ implies $\delta_3\cdots=\beta\delta\cdots=\beta\gamma\delta^2$. This determines $T_4$. Then $\epsilon_3\epsilon_4\cdots=\gamma\epsilon^2$ and $\beta_2\gamma_3\cdots=\beta\gamma^3,\beta\gamma\delta^2$ imply either two $\gamma$ adjacent or $\gamma,\delta$ adjacent. Both are contradictions.

\subsubsection*{Case. $\beta<\gamma$ and $\delta>\epsilon$}

By $\delta>\epsilon$ and $\delta+\epsilon=2\beta$, we get $\delta>\beta>\epsilon$.

We know a degree $3$ vertex $\gamma\cdots=\alpha^2\gamma,\gamma^3,\gamma\delta^2,\gamma\epsilon^2$. Then by Lemma \ref{ndegree3}, we know one of $\alpha^2\gamma,\gamma^3,\gamma\delta^2,\gamma\epsilon^2,\alpha\gamma^3,\beta\gamma^3,\gamma^4,\gamma^5$ is a vertex. The angle sums of $\alpha\delta\epsilon,\alpha\beta^2$, and the angle sum of one of $\alpha^2\gamma,\alpha\gamma^3,\beta\gamma^3,\gamma^4,\gamma^5$, and the angle sum for pentagon imply 
\begin{align*}
\alpha^2\gamma &\colon 
\alpha=(\tfrac{4}{5}-\tfrac{8}{5f})\pi,\,
\beta=(\tfrac{3}{5}+\tfrac{4}{5f})\pi,\,
\gamma=(\tfrac{2}{5}+\tfrac{16}{5f})\pi. \\
\alpha\gamma^3 &\colon 
\alpha=(\tfrac{4}{5}-\tfrac{24}{5f})\pi,\,
\beta=(\tfrac{3}{5}+\tfrac{12}{5f})\pi,\,
\gamma=(\tfrac{2}{5}+\tfrac{8}{5f})\pi. \\
\beta\gamma^3 &\colon 
\alpha=(1-\tfrac{12}{f})\pi,\,
\beta=(\tfrac{1}{2}+\tfrac{6}{f})\pi,\,
\gamma=(\tfrac{1}{2}-\tfrac{2}{f})\pi.   \\
\gamma^4 &\colon 
\alpha=(1-\tfrac{8}{f})\pi,\,
\beta=(\tfrac{1}{2}+\tfrac{4}{f})\pi,\,
\gamma=\tfrac{1}{2}\pi.  \\
\gamma^5 &\colon 
\alpha=(\tfrac{4}{5}-\tfrac{8}{f})\pi,\,
\beta=(\tfrac{3}{5}+\tfrac{4}{f})\pi,\,
\gamma=\tfrac{2}{5}\pi.
\end{align*}
Since $\beta>\gamma$ in all the cases above, we know $\alpha^2\gamma,\alpha\gamma^3,\beta\gamma^3,\gamma^4,\gamma^5$ are not vertices, and one of $\gamma^3,\gamma\delta^2,\gamma\epsilon^2$ is a vertex. 

The angle sums of $\alpha\delta\epsilon,\alpha\beta^2,\gamma\epsilon^2$ and the angle sum for pentagon imply
\[
\beta=(1+\tfrac{4}{f})\pi-\gamma,\,
\delta=(1+\tfrac{8}{f})\pi-\tfrac{3}{2}\gamma,\,
\epsilon=\pi-\tfrac{1}{2}\gamma.
\]
By $\beta<\gamma$, we get $\gamma>(\frac{1}{2}+\frac{2}{f})\pi$. By $\delta>\epsilon$, we get $\gamma<\frac{8}{f}\pi$. Then $\frac{1}{2}+\frac{2}{f}<\frac{8}{f}$, and we get $f<12$, a contradiction. 

\subsubsection*{Subcase. $\gamma^3$ is a vertex}

The angle sums of $\alpha\delta\epsilon,\alpha\beta^2,\gamma^3$ and the angle sum for pentagon imply
\[
\alpha=(\tfrac{4}{3}-\tfrac{8}{f})\pi,\,
\beta=(\tfrac{1}{3}+\tfrac{4}{f})\pi,\,
\gamma=\tfrac{2}{3}\pi,\,
\delta+\epsilon=(\tfrac{2}{3}+\tfrac{8}{f})\pi.
\]
We have $\alpha>\gamma>\beta$. By $\alpha\delta\epsilon$, and $\alpha>\beta$, and $\delta>\epsilon$, we know $\beta\epsilon^2$ is not a vertex. If $\gamma\epsilon^2$ is a vertex, then the angle sum of $\gamma\epsilon^2$ further implies $\epsilon=\gamma=\tfrac{2}{3}\pi$, and $\delta=\tfrac{8}{f}\pi$, contradicting $\delta>\epsilon$. Then by distinct $\alpha,\beta,\gamma$ values, the only degree $3$ vertices besides $\alpha\delta\epsilon,\alpha\beta^2,\gamma^3$ are $\beta\delta^2,\gamma\delta^2$. 

In a $3^5$-tile, by no degree $3$ vertex $\epsilon^2\cdots$, we know both vertices $\delta\cdots,\epsilon\cdots$ are $\alpha\delta\epsilon$. This implies $\gamma\cdots=\gamma^3$, and further implies $\alpha\cdots=\alpha\delta\epsilon$. Then $\beta\cdots=\beta^3\cdots,\beta^2\gamma\cdots$, which cannot be degree $3$. Therefore there is no $3^5$-tile. By Lemma \ref{special_tile}, we get $f\ge 24$.

By $f\ge 24$, we get $\alpha\ge \pi$, and $\beta\le\frac{1}{2}\pi$, and $\delta<\pi$. Then by Lemma \ref{geometry6}, we get $\gamma+2\epsilon>\pi$. This means $\epsilon>\frac{1}{6}\pi$, and further implies $\beta+2\epsilon>\gamma$. Then we have $R(\alpha\gamma)=\frac{8}{f}\pi\le\frac{1}{3}\pi<\alpha,\beta,\gamma,2\epsilon$, and $R(\gamma^2)=\gamma<\alpha,2\beta,\beta+2\epsilon,\delta+\epsilon$. By $\delta>\epsilon$, this implies $\alpha\gamma\cdots$ is not a vertex, and $\gamma^2\cdots=\gamma^3,\gamma^2\epsilon^k$. 

By $\alpha\ge \pi$, we know $\alpha^2\cdots$ is not a vertex. By no $\alpha^2\cdots,\alpha\gamma\cdots$, we know an $\epsilon^2$-fan is $\thick\epsilon\thin\epsilon\thick$ or has $\alpha$. This implies $\gamma^2\epsilon^k$ is not a vertex, and $\gamma^2\cdots=\gamma^3$. The AAD $\thick^{\delta}\epsilon^{\gamma}\thin^{\gamma}\epsilon^{\delta}\thick$ of $\thick\epsilon\thin\epsilon\thick$ implies a vertex $\thin^{\alpha}\gamma^{\epsilon}\thin^{\epsilon}\gamma^{\alpha}\thin\cdots=\gamma^3=\thin^{\epsilon}\gamma^{\alpha}\thin^{\alpha}\gamma^{\epsilon}\thin^{\alpha}\gamma^{\epsilon}\thin$, contradicting no $\alpha^2\cdots$. Therefore an $\epsilon^2$-fan has $\alpha$. 

By $\epsilon>\frac{1}{6}\pi$ and $f\ge 24$, we get $R(\alpha\epsilon^2)<(\frac{1}{3}+\frac{8}{f})\pi\le \alpha,2\beta,\gamma,\beta+\epsilon,\delta+\epsilon,4\epsilon$. By $\delta>\epsilon$, this implies $\alpha\epsilon^2\cdots=\alpha\beta\epsilon^2,\alpha\epsilon^4$. Since an $\epsilon^2$-fan has $\alpha$, we know $\alpha\epsilon^4$ is not a vertex, and the only $\epsilon^2$-fan is the vertex $\alpha\beta\epsilon^2$.  

By $\delta>\epsilon$, we get $R(\gamma\delta^2)<R(\gamma\delta\epsilon)=(\tfrac{2}{3}-\tfrac{8}{f})\pi<\alpha,2\beta,\gamma,\delta+\epsilon$. By $\alpha\delta\epsilon$, and an $\epsilon^2$-fan being the vertex $\alpha\beta\epsilon^2$, this implies $\gamma\delta\cdots=\gamma\delta^2,\beta\gamma\delta^2,\beta\gamma\delta\epsilon$.

The angle sum of $\alpha\beta\epsilon^2$ further implies
\[
\alpha=(\tfrac{4}{3}-\tfrac{8}{f})\pi,\,
\beta=(\tfrac{1}{3}+\tfrac{4}{f})\pi,\,
\gamma=\tfrac{2}{3}\pi,\,
\delta=(\tfrac{1}{2}+\tfrac{6}{f})\pi,\,
\epsilon=(\tfrac{1}{6}+\tfrac{2}{f})\pi.
\]
By $\beta+\gamma+2\delta=(2+\frac{16}{f})\pi>2\pi$, we get $\gamma\delta\cdots=\gamma\delta^2,\beta\gamma\delta\epsilon$. 

By no $\alpha\gamma\cdots$, we know the AAD of $\alpha\beta\epsilon^2$ is $\thick^{\delta}\epsilon^{\gamma}\thin^{\gamma}\alpha^{\beta}\thin^{\alpha}\beta^{\delta}\thin^{\gamma}\epsilon^{\delta}\thick$. This implies $\gamma\delta\cdots=\gamma\delta^2,\beta\gamma\delta\epsilon$ is a vertex. The angle sum of the vertex further imply
\[
	\alpha=\tfrac{10}{9}\pi,\,
	\beta=\tfrac{4}{9}\pi,\,
	\gamma=\delta=\tfrac{2}{3}\pi,\,
	\epsilon=\tfrac{2}{9}\pi,\,
	f=36. 
\]
The angle values imply $\alpha\beta\cdots=\alpha\beta^2,\alpha\beta\epsilon^2$, and $\delta^2\cdots=\gamma\delta^2$.

The AAD $\thick^{\delta}\epsilon^{\gamma}\thin^{\gamma}\alpha^{\beta}\thin^{\alpha}\beta^{\delta}\thin^{\gamma}\epsilon^{\delta}\thick$ of $\alpha\beta\epsilon^2$ determines $T_1,T_2,T_3,T_4$ in the first of Figure \ref{ade_a2bD}. Then $\gamma_1\gamma_3\cdots=\gamma^3$ and no $\alpha^2\cdots$ determine $T_5$. We have $\thin^{\alpha}\gamma_4^{\epsilon}\thin^{\beta}\delta_2^{\epsilon}\thick\cdots=\gamma\delta^2,\beta\gamma\delta\epsilon$. If the vertex is $\beta\gamma\delta\epsilon$, then we get $\thick^{\delta}\epsilon^{\gamma}\thin\beta\thin^{\alpha}\gamma^{\epsilon}\thin^{\beta}\delta^{\epsilon}\thick$. This implies $\alpha^2\cdots,\alpha\gamma\cdots$, a contradiction. Therefore $\thin^{\alpha}\gamma_4^{\epsilon}\thin^{\beta}\delta_2^{\epsilon}\thick\cdots=\gamma\delta^2$, and we determine $T_6$. Then $\thin^{\beta}\alpha_4^{\gamma}\thin^{\delta}\beta_6^{\alpha}\thin\cdots=\alpha\beta^2,\alpha\beta\epsilon^2$. Since $\thin^{\beta}\alpha_4^{\gamma}\thin^{\delta}\beta_6^{\alpha}\thin$ is incompatible with the AAD of $\alpha\beta\epsilon^2$, the vertex is $\alpha\beta^2$. Combined with no $\alpha^2\cdots$, we determine $T_7$. Then $\delta_3\delta_4\cdots=\gamma\delta^2$ and $\alpha_7\beta_4\cdots=\alpha\beta^2,\alpha\beta\epsilon^2$  determine $T_8,T_9$.

\begin{figure}[htp]
\centering
\begin{tikzpicture}[>=latex,scale=1]

\foreach \a/\b in {0/1, 0/0, 1.8/0.5, 1.8/1.5}
\draw[yshift=\a cm, xshift=\b cm]
	(0,-1.1) -- (0.5,-0.7) -- (0.5,0.7) -- (0,1.1) -- (-0.5,0.7) -- (-0.5,-0.7) -- (0,-1.1)
	(-0.5,0) -- (0.5,0);
	
\foreach \a/\b in {0/0, 1.8/0.5}
{
\begin{scope}[yshift=\a cm, xshift=\b cm]

\draw[line width=1.2]
	(0.5,0) -- (0.5,-0.7)
	(1.5,0) -- (1.5,0.7)
	(0,1.1) -- (-0.5,0.7);

\node at (1.3,-0.6) {\small $\alpha$};
\node at (1,-0.85) {\small $\beta$};	
\node at (1.3,-0.2) {\small $\gamma$};
\node at (0.7,-0.6) {\small $\delta$};
\node at (0.7,-0.2) {\small $\epsilon$}; 

\node at (-0.3,-0.6) {\small $\alpha$};
\node at (0,-0.85) {\small $\beta$};	
\node at (-0.3,-0.2) {\small $\gamma$};
\node at (0.3,-0.6) {\small $\delta$};
\node at (0.3,-0.2) {\small $\epsilon$}; 

\node at (0.7,0.6) {\small $\alpha$};
\node at (1,0.85) {\small $\gamma$};
\node at (0.7,0.2) {\small $\beta$};	
\node at (1.3,0.2) {\small $\delta$}; 
\node at (1.3,0.6) {\small $\epsilon$};

\node at (-0.3,0.6) {\small $\epsilon$};
\node at (0,0.85) {\small $\delta$};
\node at (-0.3,0.2) {\small $\gamma$};
\node at (0.3,0.6) {\small $\beta$};	
\node at (0.3,0.2) {\small $\alpha$}; 

\end{scope}
}

\draw
	(1.5,0.7) -- (2.3,0.3) -- (2.3,-0.7) -- (1.5,-0.7)
	(2,1.1) -- (2.8,1.1)
	(2,2.5) -- (2.8,2.5) -- (2.8,0.3) -- (2.3,0.3)
	(-0.5,0.7) -- (-1,1.3) -- (0,2.5)
	(-1,1.3) -- (-1.3,0) -- (-0.5,-0.7);

\draw[line width=1.2]
	(-1.3,0) -- (-0.5,-0.7);

\node at (2.1,-0.5) {\small $\alpha$};
\node at (1.7,-0.5) {\small $\beta$};
\node at (2.1,0.2) {\small $\gamma$};
\node at (1.7,0) {\small $\delta$};
\node at (1.7,0.4) {\small $\epsilon$};

\node at (1.8,0.75) {\small $\alpha$};
\node at (2.05,0.9) {\small $\beta$};
\node at (2.3,0.5) {\small $\gamma$};
\node at (2.65,0.9) {\small $\delta$};
\node at (2.65,0.5) {\small $\epsilon$};

\node at (2.6,1.3) {\small $\alpha$};
\node at (2.2,1.3) {\small $\beta$};
\node at (2.6,2.3) {\small $\gamma$};
\node at (2.2,1.8) {\small $\delta$};
\node at (2.2,2.3) {\small $\epsilon$};

\node at (-0.7,0.65) {\small $\alpha$};
\node at (-0.95,0.9) {\small $\beta$};
\node at (-0.65,0) {\small $\gamma$};
\node at (-1.15,0.15) {\small $\delta$};
\node at (-0.65,-0.35) {\small $\epsilon$};

\node at (-0.15,2.1) {\small $\alpha$};
\node at (-0.75,1.35) {\small $\beta$};
\node at (-0.2,1.8) {\small $\gamma$};
\node at (-0.5,1) {\small $\delta$};
\node at (-0.2,1.2) {\small $\epsilon$};

\node[inner sep=0.5,draw,shape=circle] at (0.5,2.2) {\small 1};
\node[inner sep=0.5,draw,shape=circle] at (1.5,2.2) {\small 2};
\node[inner sep=0.5,draw,shape=circle] at (0.5,1.4) {\small 3};
\node[inner sep=0.5,draw,shape=circle] at (1.5,1.4) {\small 4};
\node[inner sep=0.5,draw,shape=circle] at (-0.4,1.5) {\small 5};
\node[inner sep=0.5,draw,shape=circle] at (2.5,1.8) {\small 6};
\node[inner sep=0.5,draw,shape=circle] at (2.35,0.85) {\small 7};
\node[inner sep=0.5,draw,shape=circle] at (1,0.4) {\small 8};
\node[inner sep=0.5,draw,shape=circle] at (2,-0.15) {\small 9};
\node[inner sep=0,draw,shape=circle] at (0,0.4) {\footnotesize 10};
\node[inner sep=0,draw,shape=circle] at (0,-0.4) {\footnotesize 11};
\node[inner sep=0,draw,shape=circle] at (1,-0.4) {\footnotesize 12};
\node[inner sep=0,draw,shape=circle] at (-0.85,0.35) {\footnotesize 13};


\begin{scope}[xshift=4.5cm]

\foreach \a in {0,1}
{
\begin{scope}[xshift=2.5*\a cm]

\draw
	(-0.5,-0.4) -- (-0.5,0.7) -- (0,1.1) -- (0.5,0.7) -- (0.5,-0.4)
	(0,1.1) -- ++(0,0.4)
	(-0.5,0.7) -- ++(-0.4,0);
		
\draw[line width=1.2]
	(-0.5,0) -- (0.5,0);

\node at (0,0.85) {\small $\alpha$}; 
\node at (-0.3,0.6) {\small $\beta$};
\node at (0.3,0.6) {\small $\gamma$};
\node at (-0.3,0.2) {\small $\delta$};
\node at (0.3,0.2) {\small $\epsilon$};

\node at (-0.7,0.5) {\small $\beta$};
\node at (-0.6,0.9) {\small $\alpha$};
\node at (-0.7,0) {\small $\alpha$};
\node at (0.3,-0.2) {\small $\delta$};
\node at (-0.3,-0.2) {\small $\epsilon$};	
\node at (-0.2,1.2) {\small $\beta$};
\node at (0.2,1.2) {\small $\beta$};

\end{scope}
}

\draw[line width=1.2]
	(0.5,0.7) -- ++(0.4,0);

\node at (0.7,0.5) {\small $\delta$};
\node at (0.6,0.9) {\small $\delta$};
\node at (0.7,0) {\small $H$};

\begin{scope}[xshift=2.5cm]

\node at (0.7,0.7) {\small $H$};
\node at (0.7,0) {\small $\alpha$};

\end{scope}

\end{scope}

\end{tikzpicture}
\caption{Proposition \ref{ade_a2b}: $\beta<\gamma$ and $\delta>\epsilon$.}
\label{ade_a2bD}
\end{figure}
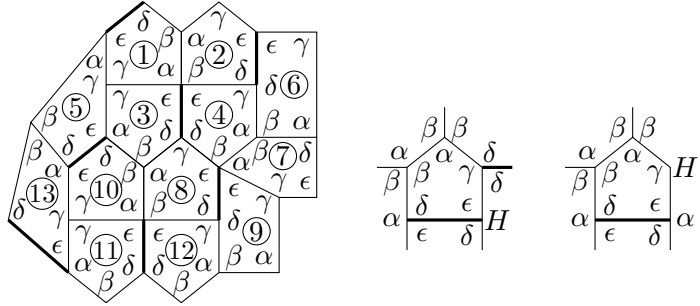

We derive $\alpha_7\beta_4\epsilon_8\epsilon_9=\alpha\beta\epsilon^2$ by starting from $\alpha_1\beta_2\epsilon_3\epsilon_4=\alpha\beta\epsilon^2$. If we repeat the argument with $\alpha_7\beta_4\epsilon_8\epsilon_9$ in place of $\alpha_1\beta_2\epsilon_3\epsilon_4$, then we derive the next $\alpha_{10}\beta_8\epsilon_{11}\epsilon_{12}=\alpha\beta\epsilon^2$. Similar to deriving $T_5$ from $\alpha_1\beta_2\epsilon_3\epsilon_4$, we also derive $T_{13}$ from $\alpha_{10}\beta_8\epsilon_{11}\epsilon_{12}$. The process repeats, and we derive more $\alpha\beta\epsilon^2$. As a result, we get more  $\beta$ at $\beta_5\beta_{13}\cdots$. This implies $\beta^k$ is a vertex, contradicting $\beta=\tfrac{4}{9}\pi$.

We conclude there is no $\epsilon^2$-fan. By Lemma \ref{fbalance}, there is only $\delta\epsilon$-fan. By $\gamma\delta\cdots=\gamma\delta^2,\beta\gamma\delta^2,\beta\gamma\delta\epsilon$, this implies a $b$-vertex $\gamma\cdots=\beta\gamma\delta\epsilon$. Moreover, by $R(\alpha\beta)=\beta<\alpha,\gamma,\delta+\epsilon$, this also implies $\alpha\beta\cdots=\alpha\beta^2$.

By no $\alpha^2\cdots,\alpha\gamma\cdots$, the AAD of $\beta\gamma\delta\epsilon$ is $\thick^{\epsilon}\delta^{\beta}\thin^{\alpha}\beta^{\delta}\thin^{\alpha}\gamma^{\epsilon}\thin^{\gamma}\epsilon^{\delta}\thick$ or $\thick^{\epsilon}\delta^{\beta}\thin^{\alpha}\gamma^{\epsilon}\thin^{\alpha}\beta^{\delta}\thin^{\gamma}\epsilon^{\delta}\thick$. The first AAD implies $\thin^{\gamma}\alpha^{\beta}\thin^{\delta}\beta^{\alpha}\thin\cdots=\alpha\beta^2=\thin^{\delta}\beta^{\alpha}\thin\beta\thin^{\gamma}\alpha^{\beta}\thin$, contradicting no $\alpha^2\cdots,\alpha\gamma\cdots$. The second AAD implies $\thick^{\epsilon}\delta^{\beta}\thin^{\epsilon}\gamma^{\alpha}\thin\cdots=\beta\gamma\delta\epsilon$. However, $\thick^{\epsilon}\delta^{\beta}\thin^{\epsilon}\gamma^{\alpha}\thin$ is incompatible with the second AAD of $\beta\gamma\delta\epsilon$. 

Therefore $\beta\gamma\delta\epsilon$ is not a vertex, and $\gamma\cdots$ is a $\hat{b}$-vertex. Then by no $\alpha\gamma\cdots,\gamma\delta\cdots$, the AAD shows $\beta\thin\epsilon\cdots$ is not a vertex.

By the values of $\alpha,\beta,\delta+\epsilon$, and only $\delta\epsilon$-fan, and $\gamma\cdots$ being $\hat{b}$-vertex, we know $\alpha\delta\epsilon,\beta^2\delta\epsilon,\beta^3\delta\epsilon,\beta\delta^2\epsilon^2,\delta^2\epsilon^2$ are all the $b$-vertices. Then by no $\beta\thin\epsilon\cdots$, we know $\alpha\delta\epsilon,\delta^2\epsilon^2$ are all the $b$-vertices. 

The $\delta\epsilon$-fan in $\delta^2\epsilon^2$ implies a vertex $\beta\gamma\cdots$. The angle sum of $\delta^2\epsilon^2$ further implies 
\[
\alpha=\delta+\epsilon=\pi,\,
\beta=\tfrac{1}{2}\pi,\, 
\gamma=\tfrac{2}{3}\pi,\,
f=24.
\]
Since $\beta\gamma\cdots$ is a $\hat{b}$-vertex, this implies $\beta\gamma\cdots$ is not a vertex. Therefore $\delta^2\epsilon^2$ is not a vertex, and $\alpha\delta\epsilon$ is the only $b$-vertex. By $\alpha\beta^2$ and applying the counting lemma to $\alpha,\delta$, we get a contradiction. 

\subsubsection*{Subcase. $\gamma^3$ is not a vertex, and $\gamma\delta^2$ is a vertex}

By Lemma \ref{geometry2}, we know $\alpha\delta\epsilon,\gamma\delta^2$ are the only degree $3$ $b$-vertices. Therefore $\alpha\delta\epsilon,\alpha\beta^2,\gamma\delta^2$ are all the degree $3$ vertices.

By $\alpha\beta^2,\gamma\delta^2$, we get $\beta,\delta<\pi$. Then by Lemma \ref{geometry6}, we get $\gamma+2\epsilon>\pi$. 

In the second and third of Figure \ref{ade_a2bD}, we consider a special tile. The second assumes $\gamma\cdots$ has degree $3$. Then $\gamma\cdots=\gamma\delta^2$. This implies $\epsilon\cdots=\beta\epsilon\cdots$ is the vertex $H$ of degree $4$ or $5$. The third assumes $\gamma\cdots$ is the vertex $H$ of degree $4$ or $5$. Then $\epsilon\cdots$ is the degree $3$ vertex $\alpha\delta\epsilon$. In both cases, we further get the degree $3$ vertices $\alpha\cdots=\beta\cdots=\alpha\beta^2$ and $\delta\cdots=\alpha\delta\epsilon$ as indicated. Moreover, both special tiles have $\delta\cdots=\thick^{\epsilon}\delta^{\beta}\thin^{\beta}\alpha^{\gamma}\thin^{\gamma}\epsilon^{\delta}\thick$. This implies a vertex $\gamma^2\cdots$. By Lemma \ref{special_tile}, the argument also implies $f\ge 24$.

By $\alpha\beta^2$ and $\beta<\gamma$, we get $R(\gamma^2)<\alpha$. By $\beta+\gamma>\pi$, we get $R(\gamma^2)<2\beta=\delta+\epsilon$. Then by $\gamma,\delta>\beta>\epsilon$, and the only degree $3$ $\hat{b}$-vertex $\alpha\beta^2$, we get $\gamma^2\cdots=\gamma^2\epsilon^2\cdots$, with no $\alpha,\delta$ in the remainder. Then by $\gamma+2\epsilon>\pi$, and $\beta+\gamma>\pi$, and $\beta<\gamma$, we get $\gamma^2\cdots=\gamma^2\epsilon^2$. The angle sum of $\gamma^2\epsilon^2$ further implies
\[
\alpha=\tfrac{24}{f}\pi,\,
\beta=(1-\tfrac{12}{f})\pi,\,
\gamma=\tfrac{16}{f}\pi,\,
\delta=(1-\tfrac{8}{f})\pi,\,
\epsilon=(1-\tfrac{16}{f})\pi.
\]

In a special tile, we have $H=\beta\delta\epsilon\cdots,\alpha\beta\gamma\cdots,\alpha\gamma^2\cdots,\beta\gamma\delta\cdots,\gamma^2\delta\cdots$. By $\alpha\beta^2$ and $\beta<\gamma$, we know $\alpha\beta\gamma\cdots,\alpha\gamma^2\cdots$ are not vertices. By $f\ge 24$, we get $2\beta+\delta+\epsilon=4\beta\ge 2\pi$. Then by $\beta<\gamma$ and $\delta>\epsilon$, this implies $H=\beta^2\delta\epsilon$. The angle sum of $\beta^2\delta\epsilon$ further implies
\[
\alpha=\pi,\,
\beta=\tfrac{1}{2}\pi,\,
\gamma=\delta=\tfrac{2}{3}\pi,\,
\epsilon=\tfrac{1}{3}\pi,\,
f=24.
\]
The angle values imply $\alpha^2\cdots,\alpha\gamma\cdots$ are not vertices. Then the AAD implies $\gamma^2\cdots=\gamma^2\epsilon^2$ is not a vertex.
\end{proof}

\begin{proposition}\label{ade_3a}
There is no tiling, such that $\alpha,\beta,\gamma$ have distinct values, and $\alpha\delta\epsilon,\alpha^3$ are vertices.
\end{proposition}

\begin{proof}
By $\alpha^3$ and distinct $\alpha,\beta,\gamma$ values, degree $3$ $\hat{b}$-vertices besides $\alpha^3$ are $\alpha\beta\gamma,\beta^2\gamma,\beta\gamma^2$. By $\alpha\delta\epsilon$ and Proposition \ref{ade_abc}, we know $\alpha\beta\gamma$ is not a vertex. By Proposition \ref{ade_2bc}, a tiling with $\alpha\delta\epsilon$ and one of $\beta^2\gamma,\beta\gamma^2$ does not have $\alpha^3$. Therefore $\beta^2\gamma,\beta\gamma^2$ are not vertices, and $\alpha^3$ is the only degree $3$ $\hat{b}$-vertex. 

If $\delta^2\cdots$ is not degree $3$, then all degree $3$ vertices are $\alpha\delta\epsilon,\alpha^3,\theta\epsilon^2$, with $\theta=\beta,\gamma$. Applying Lemma \ref{degree3} to $\alpha,\epsilon$, we get a contradiction. Therefore there is a degree $3$ vertex $\delta^2\cdots$. By the same reason, there is a degree $3$ vertex $\epsilon^2\cdots$. Then by $\alpha\delta\epsilon$ and Lemma \ref{geometry2}, we know $\beta\delta^2,\gamma\epsilon^2$ are vertices. The angle sums of $\alpha\delta\epsilon,\alpha^3,\beta\delta^2,\gamma\epsilon^2$ and the angle sum for pentagon imply $f=12$, a contradiction. 
\end{proof}

\begin{proposition}\label{ade_b2d_c2e} 
There is no tiling, such that $\alpha,\beta,\gamma$ have distinct values, and $\alpha\delta\epsilon,\beta\delta^2,\gamma\epsilon^2$ are vertices.
\end{proposition}

\begin{proof}
By distinct $\alpha,\beta,\gamma$ values, and by Propositions \ref{ade_abc}, \ref{ade_2bc}, \ref{ade_a2b}, \ref{ade_3a}, and by the exchange symmetry $(\beta,\delta)\leftrightarrow(\gamma,\epsilon)$, we know degree $3$ vertices besides $\alpha\delta\epsilon,\beta\delta^2,\gamma\epsilon^2$ are $\alpha^2\beta,\alpha^2\gamma,\beta^3,\gamma^3$. 

The angle sums of $\alpha\delta\epsilon,\beta\delta^2,\gamma\epsilon^2,\alpha^2\beta$ and the angle sum for pentagon imply
\[
	\alpha=\delta=(\tfrac{1}{2}+\tfrac{2}{f})\pi,\,
	\beta=\epsilon=(1-\tfrac{4}{f})\pi,\,
	\gamma=\tfrac{8}{f}\pi.
\]
We get $\pi>\beta=\epsilon>\alpha=\delta>\gamma$, contradicting Lemma \ref{geometry4}. Therefore $\alpha^2\beta$ is not a vertex. By the same reason, we know $\alpha^2\gamma$ is not a vertex. 

Therefore $\beta^3,\gamma^3$ are the only degree $3$ $\hat{b}$-vertices. By applying Lemma \ref{degree3} to $\delta,\epsilon$, and $\beta\ne\gamma$, we know exactly one of $\beta^3,\gamma^3$ is a vertex. Up to the exchange symmetry $(\beta,\delta)\leftrightarrow(\gamma,\epsilon)$, we may assume $\beta^3$ is the only degree $3$ $\hat{b}$-vertex.

The angle sums of $\alpha\delta\epsilon,\beta\delta^2,\gamma\epsilon^2,\beta^3$ and the angle sum for pentagon imply
\[
	\alpha=(\tfrac{1}{2}+\tfrac{2}{f})\pi,\,
	\beta=\delta=\tfrac{2}{3}\pi,\,
	\gamma=(\tfrac{1}{3}+\tfrac{4}{f})\pi,\,
	\epsilon=(\tfrac{5}{6}-\tfrac{2}{f})\pi.
\]
We have $\beta>\alpha>\gamma$.

The AAD $\thick^{\epsilon}\delta^{\beta}\thin^{\alpha}\beta^{\delta}\thin^{\beta}\delta^{\epsilon}\thick$ of $\beta\delta^2$ implies a vertex $\thin^{\gamma}\alpha^{\beta}\thin^{\delta}\beta^{\alpha}\thin\cdots$. By $\alpha\delta\epsilon$, and $\beta>\alpha$, and Lemma \ref{square}, we know $\alpha\beta\cdots$ is a $\hat{b}$-vertex. Then by $R(\alpha\beta)=(\frac{5}{6}-\frac{2}{f})\pi<\alpha+\gamma<3\gamma$, and $\beta>\alpha>\gamma$, and the only degree $3$ $\hat{b}$ vertex $\beta^3$, we get $\alpha\beta\cdots=\alpha\beta\gamma^2$. The angle sum of $\alpha\beta\gamma^2$ further implies
\[
\alpha=\tfrac{8}{15}\pi,\,
\beta=\delta=\tfrac{2}{3}\pi,\,
\gamma=\tfrac{2}{5}\pi,\,
\epsilon=\tfrac{4}{5}\pi,\,
f=60.
\]
The angle values imply $\beta\delta\cdots=\beta\delta^2$, and $\gamma\epsilon\cdots=\epsilon^2\cdots=\gamma\epsilon^2$, and $\alpha^2\cdots$ is not a vertex.

\begin{figure}[htp]
\centering
\begin{tikzpicture}[>=latex,scale=1]

\foreach \a in {-1,0}
\draw[xshift=\a cm]
	(-0.5,-0.7) -- (-0.5,0.7) -- (0,1.1) -- (0.5,0.7) -- (0.5,-0.7) -- (0,-1.1) -- (-0.5,-0.7);

\draw
	(0,1.1) -- (0,1.8) -- (1.3,1.8) -- (1.3,-0.7) -- (0.5,-0.7)
	(1.3,0.7) -- (0.5,0.7)
	(-1.5,0) -- (0.5,0);

\draw[line width=1.2]
	(0.5,0) -- (0.5,-0.7)
	(0,1.1) -- (0.5,0.7)
	(-0.5,-0.7) -- (-1,-1.1)
	(-1.5,0) -- (-1.5,0.7);

\node at (-0.7,0.6) {\small $\alpha$};
\node at (-1,0.85) {\small $\beta$};
\node at (-0.7,0.2) {\small $\gamma$};
\node at (-1.3,0.6) {\small $\delta$};	
\node at (-1.3,0.2) {\small $\epsilon$}; 

\node at (-0.7,-0.6) {\small $\epsilon$};
\node at (-1,-0.85) {\small $\delta$};
\node at (-0.7,-0.2) {\small $\gamma$};
\node at (-1.3,-0.6) {\small $\beta$};	
\node at (-1.3,-0.2) {\small $\alpha$}; 

\node at (0,0.85) {\small $\epsilon$};
\node at (-0.3,0.6) {\small $\gamma$};
\node at (0.3,0.6) {\small $\delta$};
\node at (-0.3,0.2) {\small $\alpha$}; 
\node at (0.3,0.2) {\small $\beta$};		

\node at (0,-0.85) {\small $\gamma$};
\node at (-0.3,-0.6) {\small $\alpha$};
\node at (0.3,-0.6) {\small $\epsilon$};
\node at (-0.3,-0.2) {\small $\beta$}; 
\node at (0.3,-0.2) {\small $\delta$};	

\node at (0.6,0.9) {\small $\delta$};
\node at (0.2,1.2) {\small $\epsilon$};
\node at (1.1,0.9) {\small $\beta$};
\node at (1.1,1.6) {\small $\gamma$};
\node at (0.2,1.6) {\small $\gamma$};

\node at (-0.2,1.2) {\small $\gamma$};

\node at (1.1,0.5) {\small $\alpha$};
\node at (1.1,-0.5) {\small $\gamma$};
\node at (0.7,0.5) {\small $\beta$};
\node at (0.7,-0.5) {\small $\epsilon$};
\node at (0.7,0) {\small $\delta$};

\node[inner sep=0.5,draw,shape=circle] at (0,0.4) {\small $1$};
\node[inner sep=0.5,draw,shape=circle] at (0.8,1.3) {\small $4$};
\node[inner sep=0.5,draw,shape=circle] at (0,-0.4) {\small $2$};
\node[inner sep=0.5,draw,shape=circle] at (1.05,0) {\small $3$};
\node[inner sep=0.5,draw,shape=circle] at (-1,0.4) {\small $5$};
\node[inner sep=0.5,draw,shape=circle] at (-1,-0.4) {\small $6$};

\end{tikzpicture}
\caption{Proposition \ref{ade_b2d_c2e}: $\alpha\delta\epsilon,\beta\delta^2,\gamma\epsilon^2$ are vertices.}
\label{ade_b2d_c2eA} 
\end{figure}

The AAD $\thick^{\epsilon}\delta^{\beta}\thin^{\alpha}\beta^{\delta}\thin^{\beta}\delta^{\epsilon}\thick$ of $\beta\delta^2$ determines $T_1,T_2,T_3$ in Figure \ref{ade_b2d_c2eA}. Then $\beta_3\delta_1\cdots=\beta\delta^2$ determines $T_4$, and $\alpha_1\beta_2\cdots=\alpha\beta\gamma^2$ and no $\alpha^2\cdots$ determine $T_5,T_6$. Then $\epsilon_1\epsilon_4\cdots=\gamma\epsilon^2$ implies $\alpha_5\gamma_1\cdots=\alpha^2\gamma\cdots,\alpha\gamma\epsilon\cdots$, contradicting no $\alpha^2\cdots$, and $\gamma\epsilon\cdots=\gamma\epsilon^2$. 
\end{proof}

\begin{proposition}\label{ade_2ac} 
There is no tiling, such that $\alpha,\beta,\gamma$ have distinct values, and $\alpha\delta\epsilon,\alpha^2\beta$ are vertices.
\end{proposition}

\begin{proof}
By distinct $\alpha,\beta,\gamma$ values, and by Propositions \ref{ade_abc}, \ref{ade_2bc}, \ref{ade_a2b}, \ref{ade_3a}, and by the exchange symmetry $(\beta,\delta)\leftrightarrow(\gamma,\epsilon)$, we know degree $3$ vertices besides $\alpha\delta\epsilon,\alpha^2\beta$ are $\gamma^3,\beta\delta^2,\beta\epsilon^2,\gamma\delta^2,\gamma\epsilon^2$. 

By Lemma \ref{geometry2} and Proposition \ref{ade_b2d_c2e}, we know $\beta\delta^2,\beta\epsilon^2,\gamma\delta^2,\gamma\epsilon^2$ are mutually exclusive. If $\gamma^3$ is not a vertex, then $\alpha\delta\epsilon,\alpha^2\beta$ and one of these are all the degree $3$ vertices. If $\alpha\delta\epsilon,\alpha^2\beta,\beta\delta^2$ are all the degree $3$ vertices, then by applying Lemma \ref{degree3} to $\alpha,\delta$, we get a contradiction. We get the same contradiction for all the other combinations.

Therefore $\gamma^3$ is a vertex. The angle sums of $\alpha\delta\epsilon,\alpha^2\beta,\gamma^3$ and the angle sum for pentagon imply
\[
\alpha=(\tfrac{5}{6}-\tfrac{2}{f})\pi,\;
\beta=(\tfrac{1}{3}+\tfrac{4}{f})\pi,\;
\gamma=\tfrac{2}{3}\pi,\;
\delta+\epsilon=(\tfrac{7}{6}+\tfrac{2}{f})\pi.
\]
We have $\alpha>\gamma>\beta$. By Lemma \ref{geometry1}, this implies $\delta>\epsilon$. 

Suppose $\delta^2\cdots,\epsilon^2\cdots$ are vertices. By $\delta>\epsilon$ and $\delta+\epsilon>\pi$, we know $R(\delta^2)$ has no $\delta,\epsilon$. By $\alpha\delta\epsilon$ and $\delta>\epsilon$, we get $R(\delta^2)<\alpha<\beta+\gamma,3\beta$. Therefore $\delta^2\cdots=\beta\delta^2,\gamma\delta^2,\beta^2\delta^2$. The angle sum of $\delta^2\cdots$ further implies
\begin{align*}
\beta\delta^2 & \colon
	\alpha=\delta=(\tfrac{5}{6}-\tfrac{2}{f})\pi,\,
	\beta=(\tfrac{1}{3}+\tfrac{4}{f})\pi,\,
	\gamma=\tfrac{2}{3}\pi,\,
	\epsilon=(\tfrac{1}{2}+\tfrac{2}{f})\pi. \\
\gamma\delta^2 & \colon
	\alpha=(\tfrac{5}{6}-\tfrac{2}{f})\pi,\,
	\beta=(\tfrac{1}{3}+\tfrac{4}{f})\pi,\,
	\gamma=\delta=\tfrac{2}{3}\pi,\,
	\epsilon=(\tfrac{1}{2}+\tfrac{2}{f})\pi. \\
\beta^2\delta^2 & \colon
	\alpha=(\tfrac{5}{6}-\tfrac{2}{f})\pi,\,
	\beta=(\tfrac{1}{3}+\tfrac{4}{f})\pi,\,
	\gamma=\tfrac{2}{3}\pi,\,
	\delta=(\tfrac{2}{3}-\tfrac{4}{f})\pi,\,
	\epsilon=(\tfrac{1}{2}+\tfrac{6}{f})\pi.
\end{align*}
We have $R(\epsilon^2)<\pi<\beta+\gamma,3\beta,2\epsilon$. By $\alpha\delta\epsilon$, and $\alpha>\gamma>\beta$, and $\delta>\epsilon$, this implies $\epsilon^2\cdots=\beta\epsilon^2,\gamma\epsilon^2,\beta^2\epsilon^2$. 

By $\gamma^3$ and $\alpha>\gamma$, we have $R(\alpha^2)<R(\alpha\gamma)<\gamma<\alpha,2\beta,2\epsilon$. By $\alpha^2\beta$ and $\delta>\epsilon$, this implies $\alpha^2\cdots=\alpha^2\beta$, and $\alpha\gamma\cdots$ is not a vertex. By no $\alpha\gamma\cdots$, the AAD shows $\beta\epsilon^2,\gamma\epsilon^2$ are not vertices. Moreover, the AAD of $\beta^2\epsilon^2$ is $\thick^{\delta}\epsilon^{\gamma}\thin^{\delta}\beta^{\alpha}\thin^{\alpha}\beta^{\delta}\thin^{\gamma}\epsilon^{\delta}\thick$. This implies $\thin^{\gamma}\alpha^{\beta}\thin^{\beta}\alpha^{\gamma}\thin\cdots=\alpha^2\beta=\thin^{\beta}\alpha^{\gamma}\thin^{\alpha}\beta^{\delta}\thin^{\gamma}\alpha^{\beta}\thin$, contradicting no $\alpha\gamma\cdots$. 

We conclude $\delta^2\cdots,\epsilon^2\cdots$ cannot both be vertices. By the balance lemma, this implies a $b$-vertex is $\delta\epsilon\cdots$, with no $\delta,\epsilon$ in the remainder. By $R(\delta\epsilon)=\alpha<\beta+\gamma,3\beta$, we know $\alpha\delta\epsilon,\beta^2\delta\epsilon$ are all the $b$-vertices. In particular, $\gamma\delta\cdots,\delta\thin\delta\cdots$ are not vertices, and $\alpha\gamma\cdots$ is a $\hat{b}$-vertex.

By $\alpha^2\beta$ and applying the counting lemma to $\alpha,\epsilon$, we know $\beta^2\delta\epsilon$ is a vertex. The angle sum of $\beta^2\delta\epsilon$ further implies
\[
\alpha=\tfrac{4}{5}\pi,\,
\beta=\tfrac{2}{5}\pi,\,
\gamma=\tfrac{2}{3}\pi,\,
\delta+\epsilon=\tfrac{6}{5}\pi,\,
f=60.
\]
By no $\gamma\delta\cdots,\delta\thin\delta\cdots$, the AAD of $\beta^2\delta
\epsilon$ is $\thick^{\delta}\epsilon^{\gamma}\thin^{\alpha}\beta^{\delta}\thin^{\alpha}\beta^{\delta}\thin^{\beta}\delta^{\epsilon}\thick$. This implies a vertex $\alpha\gamma\cdots$. We know $\alpha\gamma\cdots$ is a $\hat{b}$-vertex. Then by $\beta<R(\alpha\gamma)=\frac{8}{15}\pi<\alpha,2\beta,\gamma$, we know $\alpha\gamma\cdots$ is not a vertex.
\end{proof}

\begin{proposition}\label{ade_3b} 
Tilings with distinct $\alpha,\beta,\gamma$ values, and such that $\alpha\delta\epsilon,\beta^3$ are vertices, are the pentagonal subdivisions $P_{\pentagon}P_8,P_{\pentagon}P_{20}$ of the octahedron and the icosahedron.
\end{proposition}

The pentagonal subdivision $P_{\pentagon}P_4$ is not included because all propositions assume $f\ge 16$.

\begin{proof}
By $\beta^3$, and distinct $\alpha,\beta,\gamma$ values, and Propositions \ref{ade_abc}, \ref{ade_a2b}, \ref{ade_2ac}, we know $\beta^3$ is the only degree $3$ $\hat{b}$-vertex.

The angle sums of $\alpha\delta\epsilon,\beta^3$ and the angle sum for pentagon imply
\[
\alpha+\delta+\epsilon=2\pi,\,
\beta=\tfrac{2}{3}\pi,\,
\gamma=(\tfrac{1}{3}+\tfrac{4}{f})\pi.
\]
We have $2\gamma>\beta>\gamma$. By Lemma \ref{geometry1}, this implies $\delta<\epsilon$. 

\subsubsection*{Case. $\alpha<\pi$, and $\epsilon^2\cdots$ is not a vertex}

By no $\epsilon^2\cdots$ and the balance lemma, we know a $b$-vertex is $\delta\epsilon\cdots$, with no $\delta,\epsilon$ in the remainder. By $\alpha\delta\epsilon$, we get $R(\delta\epsilon)=\alpha<\pi<\beta+\gamma$. By $\beta>\gamma$, this implies $\alpha\delta\epsilon,\gamma^k\delta\epsilon(k\ge 2)$ are all the $b$-vertices. In particular, we know $\beta\cdots$ is a $\hat{b}$-vertex, and $\delta\thin\delta\cdots, \delta\thin\epsilon\cdots,\epsilon\thin\epsilon\cdots$ are not vertices. 

By no $\beta\epsilon\cdots,\epsilon\thin\epsilon\cdots$, we know the AAD of $\gamma^k\delta\epsilon$ is $\thick^{\epsilon}\delta^{\beta}\thin^{\alpha}\gamma^{\epsilon}\thin\cdots\thin^{\alpha}\gamma^{\epsilon}\thin^{\gamma}\epsilon^{\delta}\thick$. This implies $\alpha\beta\cdots$ is a vertex. By $\alpha\delta\epsilon,\gamma^k\delta\epsilon$, we get $\alpha=k\gamma\ge 2\gamma>\beta>\gamma$. Then by $\beta^3$, and $\alpha\ge 2\gamma>\beta$, and $\beta\cdots$ being a $\hat{b}$-vertex, and no $\alpha\beta\gamma$, we know $\alpha\beta\cdots$ is not a vertex. 

Therefore $\gamma^k\delta\epsilon$ is not a vertex, and $\alpha\delta\epsilon$ is the only $b$-vertex. Then by the counting lemma, $\alpha\delta\epsilon,\beta^k\gamma^l$ are all the vertices. By the angle values of $\beta,\gamma$, and $\beta^3$ being the only degree $3$ $\hat{b}$-vertex, we get $\beta^k\gamma^l=\beta^3,\beta\gamma^3,\gamma^4,\gamma^5$. By $\beta^3$, and applying the counting lemma to $\beta,\gamma$, we know one of $\beta\gamma^3,\gamma^4,\gamma^5$ is a vertex. The angle sum of $\beta\gamma^3,\gamma^4,\gamma^5$ further imply $f=36,24,60$, and we get the respective list of vertices
\begin{align}
f=36 &\colon
\text{AVC}
=\{\alpha\delta\epsilon,\beta^3,\beta\gamma^3\}.  \nonumber \\
f=24 &\colon
\text{AVC}
=\{\alpha\delta\epsilon,\beta^3,\gamma^4\}.   \label{ade_3b_avc}  \\
f=60 &\colon
\text{AVC}
=\{\alpha\delta\epsilon,\beta^3,\gamma^5\}. \nonumber 
\end{align}

\subsubsection*{Case. $\alpha\ge\pi$, or $\epsilon^2\cdots$ is a vertex}

By $\alpha\delta\epsilon$, we know $\alpha\ge\pi$ implies $\epsilon<\pi$. Of course the vertex $\epsilon^2\cdots$ also implies $\epsilon<\pi$. Then by $\pi>\beta>\gamma$, and Lemma \ref{geometry6}, we know $\beta+2\delta>\pi$. 

If $\alpha\ge \pi$, then $\alpha+\beta+\gamma>2\pi$, and $\alpha+\beta+2\delta>2\pi$. By $\alpha>\beta>\gamma$, and $\delta<\epsilon$, this implies $\alpha^2\cdots,\alpha\beta\cdots$ are not vertices.

Next we assume $\alpha<\pi$ and $\epsilon^2\cdots$ is a vertex, and prove that $\alpha^2\cdots,\alpha\beta\cdots$ are also not vertices. 

By $\epsilon^2\cdots$ and the balance lemma, we know $\delta^2\cdots$ is a vertex. This implies $\delta,\epsilon<\pi$. Then by $\alpha,\beta,\gamma<\pi$, we know the pentagon is convex. 

By $\alpha\delta\epsilon$ and $\alpha<\pi$, we get $\delta+\epsilon>\pi$. By $\delta<\epsilon$, this implies $R(\epsilon^2)$ has no $\delta,\epsilon$. By $\alpha\delta\epsilon$ and $\delta<\epsilon$, we get $R(\epsilon^2)<\alpha<\pi<\beta+\gamma$. Then by $\beta>\gamma$, this implies $\epsilon^2\cdots=\beta\epsilon^2,\gamma^k\epsilon^2$. 

If $\beta\epsilon^2$ is a vertex, then the angle sum of $\beta\epsilon^2$ further implies
\[
\alpha+\delta=\tfrac{4}{3}\pi,\,
\beta=\epsilon=\tfrac{2}{3}\pi,\,
\gamma=(\tfrac{1}{3}+\tfrac{4}{f})\pi.
\]
By $\delta<\epsilon$, we get $\alpha>\beta$. Then by $\alpha+\beta+\delta=2\pi$ and $\delta<\epsilon$, we know $\alpha^2\cdots,\alpha\beta\cdots$ are $\hat{b}$-vertices. Then by $\alpha>\beta=\frac{2}{3}\pi$, we get $R(\alpha^2)<R(\alpha\beta)<\alpha,\beta,2\gamma$. This implies $\alpha^2\cdots=\alpha^2\gamma$ and $\alpha\beta\cdots=\alpha\beta\gamma$. Since $\beta^3$ is the only degree 3 $\hat{b}$-vertex, we know $\alpha^2\cdots,\alpha\beta\cdots$ are not vertices.

If $\gamma\epsilon^2$ is a vertex, then the angle sum of $\gamma\epsilon^2$ further implies
\[
\alpha+\delta=(\tfrac{7}{6}+\tfrac{2}{f})\pi,\,
\beta=\tfrac{2}{3}\pi,\,
\gamma=(\tfrac{1}{3}+\tfrac{4}{f})\pi,\,
\epsilon=(\tfrac{5}{6}-\tfrac{2}{f})\pi.
\]
Assume $\alpha<\beta$. Then $\epsilon>\delta>(\tfrac{7}{6}+\tfrac{2}{f})\pi-\beta>\frac{1}{2}\pi$. This implies $R(\delta^2)$ has no $\delta,\epsilon$, and $\epsilon\thin\epsilon\cdots$ is not a vertex, and $\alpha+\gamma+2\delta>\alpha+\gamma+\delta+\frac{1}{2}\pi=(2+\tfrac{2}{f})\pi>2\pi$. Then by $\alpha<\beta$, and $\alpha+\delta>\pi$, we get $\alpha\delta^2\cdots=\alpha\delta^2$, and $\beta\delta^2\cdots=\beta\delta^2$. By $\alpha\delta\epsilon$ and $\delta<\epsilon$, we know $\alpha\delta^2$ is not a vertex. By $\alpha\delta\epsilon,\gamma\epsilon^2$, and Proposition \ref{ade_b2d_c2e}, we know $\beta\delta^2$ is not a vertex. Therefore $R(\delta^2)$ has only $\gamma$, and we get $\delta^2\cdots=\gamma^k\delta^2$. The AAD of $\gamma^k\delta^2$ implies $\beta\epsilon\cdots,\epsilon\thin\epsilon\cdots$. We already know $\epsilon\thin\epsilon\cdots$ is not a vertex. By $\alpha\delta\epsilon$, and $\alpha<\beta$, and $\delta<\epsilon$, we know $\beta\epsilon\cdots$ is also not a vertex. Therefore $\delta^2\cdots$ is not a vertex, contradicting $\gamma\epsilon^2$ and the balance lemma.

Therefore we have $\alpha>\beta$ in case $\gamma\epsilon^2$ is a vertex. By $\beta+2\delta>\pi$, we get $\delta>\frac{1}{6}\pi$, and $\beta+\delta>\frac{5}{6}\pi>\epsilon$. By $\alpha\delta\epsilon$, this implies $\alpha+\beta+2\delta>2\pi$. Then by $\alpha>\beta$, and $\delta<\epsilon$, this implies $\alpha^2\cdots,\alpha\beta\cdots$ are $\hat{b}$-vertices. Then by $\beta^3$, and $\beta<\alpha,2\gamma$, we get $\alpha^2\cdots=\alpha^2\gamma$ and $\alpha\beta\cdots=\alpha\beta\gamma$. Since $\beta^3$ is the only degree 3 $\hat{b}$-vertex, we know $\alpha^2\cdots,\alpha\beta\cdots$ are not vertices.

If $\gamma^k\epsilon^2(k\ge 2)$ is a vertex, then $\gamma+\epsilon\le \pi<\delta+\epsilon,\beta+\gamma$. This implies $\gamma<\delta$ and $\beta>\epsilon$. Therefore $\pi>\beta>\epsilon>\delta>\gamma$. Then by $\alpha<\pi$, we know the pentagon is convex. By $\alpha\delta\epsilon,\gamma^k\epsilon^2$, and $\delta<\epsilon$, we also get $\alpha>\gamma$. Since the pentagon is convex, we may apply Lemma \ref{geometry4} to get $\alpha>\beta$. Then the same argument for the $\alpha>\beta$ case of $\gamma\epsilon^2$ leads to a contradiction. 

We conclude that, no matter $\alpha\ge \pi$ or $\alpha<\pi$, we always know that $\alpha^2\cdots,\alpha\beta\cdots$ are not vertices.

By $\beta^3$, and $2\gamma>\beta>\gamma$, and no $\alpha^2\cdots,\alpha\beta\cdots$, we know $\beta^3,\alpha\gamma^k,\beta\gamma^k,\gamma^k$ are all the $\hat{b}$-vertices. 

Let $\theta=\beta,\gamma$. By no $\alpha^2\cdots,\alpha\beta\cdots$, we know the AADs of ${}^{\alpha}\thin\theta\thin,{}^{\beta}\thin\theta\thin$ are ${}^{\alpha}\thin\theta^{\alpha}\thin,{}^{\beta}\thin\theta^{\alpha}\thin$. By the AADs ${}^{\alpha}\thin\theta^{\alpha}\thin,{}^{\beta}\thin\theta^{\alpha}\thin$ of ${}^{\alpha}\thin\theta\thin,{}^{\beta}\thin\theta\thin$, we know the AAD of $\alpha\gamma^k$ is $\thin^{\gamma}\alpha^{\beta}\thin^{\epsilon}\gamma^{\alpha}\thin\cdots\thin^{\epsilon}\gamma^{\alpha}\thin$, and the AAD of $\beta^l\gamma^k$ is $\thin\theta^{\alpha}\thin\cdots\thin\theta^{\alpha}\thin$.

Let $\theta=\beta,\gamma$. By no $\alpha^2\cdots,\alpha\beta\cdots$, we know the AADs of ${}^{\alpha}\thin\theta\thin,{}^{\beta}\thin\theta\thin$ are ${}^{\alpha}\thin\theta^{\alpha}\thin,{}^{\beta}\thin\theta^{\alpha}\thin$. Then we know the AAD of $\alpha\gamma^k$ is $\thin^{\gamma}\alpha^{\beta}\thin^{\epsilon}\gamma^{\alpha}\thin\cdots\thin^{\epsilon}\gamma^{\alpha}\thin$, and the AAD of $\beta^l\gamma^k$ is $\thin\theta^{\alpha}\thin\cdots\thin\theta^{\alpha}\thin$. We also know the AAD of $\thick\delta\thin\theta\thin\cdots\thin\theta\thin$ is $\thick^{\epsilon}\delta^{\beta}\thin\theta^{\alpha}\thin\cdots\thin\theta^{\alpha}\thin$. This implies a $\delta^2$-fan is $\thick^{\epsilon}\delta^{\beta}\thin^{\beta}\delta^{\epsilon}\thick$ or has $\alpha$. Then by no $\alpha^2\cdots,\alpha\beta\cdots$, a $\delta^2$-fan with $\alpha$ has a single $\alpha$ and has no $\beta$. Then by the AAD $\thick^{\epsilon}\delta^{\beta}\thin\theta^{\alpha}\thin\cdots\thin\theta^{\alpha}\thin$ and no $\alpha\beta\cdots$, the fan is $\thick^{\epsilon}\delta^{\beta}\thin^{\epsilon}\gamma^{\alpha}\thin\cdots\thin^{\epsilon}\gamma^{\alpha}\thin^{\gamma}\alpha^{\beta}\thin^{\beta}\delta^{\epsilon}\thick$. We conclude a $\delta^2$-fan is $\thick\delta\thin\delta\thick$, $\thick\delta\thin\gamma\thin\cdots\thin\gamma\thin\alpha\thin\delta\thick$. 

By $\alpha\delta\epsilon$ and $\delta<\epsilon$, a $\delta\epsilon$-fan is the vertex $\alpha\delta\epsilon$ or has only $\theta$, and an $\epsilon^2$-fan has only $\theta$. Then by the AAD $\thick^{\epsilon}\delta^{\beta}\thin\theta^{\alpha}\thin\cdots\thin\theta^{\alpha}\thin$, we know a $\delta\epsilon$-fan without $\alpha$ is $\thick^{\epsilon}\delta^{\beta}\thin\theta^{\alpha}\thin\cdots\thin\theta^{\alpha}\thin^{\gamma}\epsilon^{\delta}\thick$. Therefore the only fan containing $\thin\alpha\thin\gamma\thin$ is the $\delta^2$-fan $\thick\delta\thin\gamma\thin\cdots\thin\gamma\thin\alpha\thin\delta\thick$.

By no $\alpha\beta\cdots$, we know $\beta^2\cdots$ has no $\alpha$. Then by $\beta^3$, and $\beta+\gamma>\pi$, and $\beta+2\delta>\pi$, and $\delta<\epsilon$, we get $\beta^2\cdots=\beta^3,\beta^2\delta^2,\beta^2\delta\epsilon$. By what we know about $\delta^2$-fans, we know $\beta^2\delta^2$ is not a vertex. The angle sum of $\beta^2\delta\epsilon$ further implies 
\[
\alpha=\tfrac{4}{3}\pi,\,
\beta=\tfrac{2}{3}\pi,\,
\gamma=(\tfrac{1}{3}+\tfrac{4}{f})\pi,\,
\delta+\epsilon=\tfrac{2}{3}\pi.
\]
By $\delta<\epsilon$, this implies $\alpha\cdots=\alpha\delta\epsilon,\alpha\gamma\delta^k,\alpha\delta^k$. By $\beta^2\delta\epsilon$, and applying the counting lemma to $\alpha,\delta$, we get a contradiction. Therefore $\beta^2\delta\epsilon$ is not a vertex. We conclude $\beta^2\cdots=\beta^3$. This implies a fan has at most one $\beta$. 

The AAD $\thick^{\epsilon}\delta^{\beta}\thin^{\beta}\delta^{\epsilon}\thick$ implies a vertex $\thin^{\alpha}\beta^{\delta}\thin^{\delta}\beta^{\alpha}\thin\cdots=\beta^3=\thin^{\delta}\beta^{\alpha}\thin^{\alpha}\beta^{\delta}\thin^{\alpha}\beta^{\delta}\thin$, contradicting no $\alpha^2\cdots$. Therefore the only $\delta^2$-fan is $\thick\delta\thin\gamma\thin\cdots\thin\gamma\thin\alpha\thin\delta\thick$.

The AAD $\thin^{\epsilon}\gamma^{\alpha}\thin^{\gamma}\epsilon^{\delta}\thick$ implies a vertex $\thin^{\beta}\alpha^{\gamma}\thin^{\epsilon}\gamma^{\alpha}\thin\cdots$. The AAD $\thin^{\beta}\alpha^{\gamma}\thin^{\epsilon}\gamma^{\alpha}\thin$ is incompatible with the AAD of the $\hat{b}$-vertex $\alpha\gamma\cdots=\alpha\gamma^k$. Moreover, the only fan containing $\thin\alpha\thin\gamma\thin$ is $\thick\delta\thin\gamma\thin\cdots\thin\gamma\thin\alpha\thin\delta\thick$. However, $\thin^{\beta}\alpha^{\gamma}\thin^{\epsilon}\gamma^{\alpha}\thin$ is also incompatible with the AAD of the $\delta^2$-fan. Therefore we do not have $\thin^{\epsilon}\gamma^{\alpha}\thin^{\gamma}\epsilon^{\delta}\thick$, and the AAD of $\thin\gamma\thin\epsilon\thick$ is $\thin^{\alpha}\gamma^{\epsilon}\thin^{\gamma}\epsilon^{\delta}\thick$. 

We already know a $\delta^2$-fan is $\thick^{\epsilon}\delta^{\beta}\thin^{\epsilon}\gamma^{\alpha}\thin\cdots\thin^{\epsilon}\gamma^{\alpha}\thin^{\gamma}\alpha^{\beta}\thin^{\beta}\delta^{\epsilon}\thick$. We also know a $\delta\epsilon$-fan without $\alpha$ is $\thick^{\epsilon}\delta^{\beta}\thin\theta^{\alpha}\thin\cdots\thin\theta^{\alpha}\thin^{\gamma}\epsilon^{\delta}\thick$. Then by the AAD $\thin^{\alpha}\gamma^{\epsilon}\thin^{\gamma}\epsilon^{\delta}\thick$ of $\thin\gamma\thin\epsilon\thick$, and at most one $\beta$ in a fan, we know $\delta\epsilon$-fans are the vertex $\alpha\delta\epsilon$, and $\thick^{\epsilon}\delta^{\beta}\thin^{\gamma}\epsilon^{\delta}\thick,\thick^{\epsilon}\delta^{\beta}\thin^{\epsilon}\gamma^{\alpha}\thin\cdots\thin^{\epsilon}\gamma^{\alpha}\thin^{\delta}\beta^{\alpha}\thin^{\gamma}\epsilon^{\delta}\thick$. We also know an $\epsilon^2$-fan has no $\alpha$ and has at most one $\beta$. Then by the AAD $\thin^{\alpha}\gamma^{\epsilon}\thin^{\gamma}\epsilon^{\delta}\thick$ of $\thin\gamma\thin\epsilon\thick$, and no $\alpha^2\cdots$, we know $\epsilon^2$-fans are $\thick^{\delta}\epsilon^{\gamma}\thin^{\gamma}\epsilon^{\delta}\thick,\thick^{\delta}\epsilon^{\gamma}\thin^{\epsilon}\gamma^{\alpha}\thin\cdots\thin^{\epsilon}\gamma^{\alpha}\thin^{\delta}\beta^{\alpha}\thin^{\gamma}\epsilon^{\delta}\thick$.

The AADs $\thick^{\epsilon}\delta^{\beta}\thin^{\gamma}\epsilon^{\delta}\thick$ and $\thick^{\delta}\epsilon^{\gamma}\thin^{\gamma}\epsilon^{\delta}\thick$ imply vertices $\thin^{\alpha}\beta^{\delta}\thin^{\epsilon}\gamma^{\alpha}\thin\cdots$ and $\thin^{\alpha}\gamma^{\epsilon}\thin^{\epsilon}\gamma^{\alpha}\thin\cdots$. The AADs $\thin^{\alpha}\beta^{\delta}\thin^{\epsilon}\gamma^{\alpha}\thin$ and $\thin^{\alpha}\gamma^{\epsilon}\thin^{\epsilon}\gamma^{\alpha}\thin$ are incompatible with the AADs $\thin\theta^{\alpha}\thin\cdots\thin\theta^{\alpha}\thin$ and $\thin^{\gamma}\alpha^{\beta}\thin^{\epsilon}\gamma^{\alpha}\thin\cdots\thin^{\epsilon}\gamma^{\alpha}\thin$ of $\hat{b}$-vertices $\beta^2\gamma,\alpha\gamma^k,\beta\gamma^k,\gamma^k$. They are also incompatible with the AADs of all the fans. 

The fan $\thick^{\epsilon}\delta^{\beta}\thin^{\epsilon}\gamma^{\alpha}\thin\cdots\thin^{\epsilon}\gamma^{\alpha}\thin^{\delta}\beta^{\alpha}\thin^{\gamma}\epsilon^{\delta}\thick$ implies $\thin^{\alpha}\beta^{\delta}\thin^{\gamma}\epsilon^{\delta}\thick$. The AAD is not compatible with the AAD of any fan. 

Therefore the fans are the vertex $\alpha\delta\epsilon$, and $\thick\delta\thin\gamma\thin\cdots\thin\gamma\thin\alpha\thin\delta\thick,\thick\epsilon\thin\gamma\thin\cdots\thin\gamma\thin\beta\thin\epsilon\thick$. By $\beta^2\cdots=\beta^3$, and no $\alpha^2\cdots,\alpha\beta\cdots$, we know the fans cannot be combined. Therefore all fans are vertices, and $\alpha\delta\epsilon,\alpha\gamma^k\delta^2(k\ge 1),\beta\gamma^k\epsilon^2$ are all the $b$-vertices.

The vertex $\epsilon^2\cdots=\beta\gamma^k\epsilon^2=\thick^{\delta}\epsilon^{\gamma}\thin^{\epsilon}\gamma^{\alpha}\thin\cdots\thin^{\epsilon}\gamma^{\alpha}\thin^{\delta}\beta^{\alpha}\thin^{\gamma}\epsilon^{\delta}\thick=\thin^{\delta}\beta^{\alpha}\thin^{\gamma}\epsilon^{\delta}\thick^{\delta}\epsilon^{\gamma}\thin\cdots$ determines $T_1,T_2,T_3$ in the first of Figure \ref{ade_3bA}. Since the AVC $\thin^{\gamma}\alpha_3^{\beta}\thin^{\epsilon}\gamma_1^{\alpha}\thin$ is incompatible with the AAD of $\alpha\gamma^k\delta^2$, we get $\thin^{\gamma}\alpha_3^{\beta}\thin^{\epsilon}\gamma_1^{\alpha}\thin\cdots=\alpha\gamma^k$. Then by no $\alpha^2\cdots$, we determine $T_4$. Then $\alpha_1\epsilon_4\cdots=\alpha\delta\epsilon$ determines $T_5$. Then $\beta_1\beta_5\cdots=\beta^3$ and no $\alpha^2\cdots$ determine $T_6$. Then $\alpha_6\delta_1\delta_2\cdots=\alpha\gamma^k\delta^2$. However, the AAD $\thick^{\epsilon}\delta_1^{\beta}\thin^{\gamma}\alpha_6^{\beta}\thin$ is incompatible with the AAD of $\alpha\gamma^k\delta^2$.

Therefore $\epsilon^2\cdots$ is not a vertex. By the balance lemma, this implies $\alpha\gamma^k\delta^2$ is not a vertex, and $\alpha\delta\epsilon$ is the only $b$-vertex. Then by the counting lemma, $\alpha\delta\epsilon,\beta^k\gamma^l$ are all the vertices. Then by $\beta=\tfrac{2}{3}\pi$ and $\gamma=(\tfrac{1}{3}+\tfrac{4}{f})\pi$, we get the AVC \eqref{ade_3b_avc}.

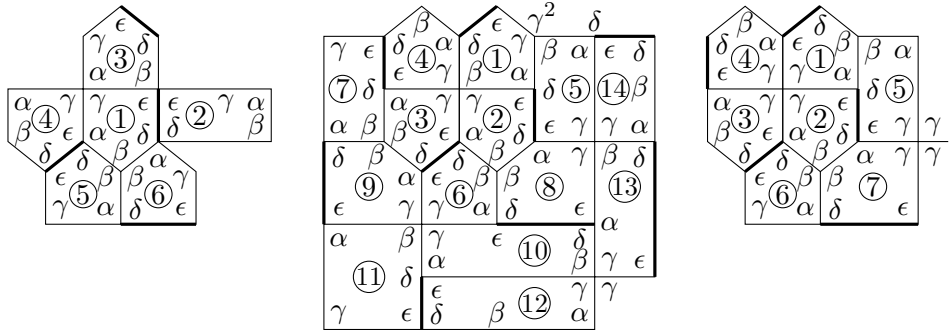
\begin{figure}[htp]
\centering
\begin{tikzpicture}[>=latex,scale=1]


\begin{scope}[xshift=-4.5cm, yscale=-1]

\draw	
	(0.5,-0.7) -- (0,-1.1) -- (-0.5,-0.7) -- (-0.5,0.7) -- (0,1.1) -- (0.5,0.7) -- (0.5,-0.7)
	(-0.5,0.7) -- (-1,1.1) -- (-1,1.8) -- (1,1.8) -- (1,1.1) -- (0.5,0.7) -- (2,0.7) -- (2,0) -- (-1.5,0) -- (-1.5,0.7) -- (-1,1.1)
	(0,1.8) -- (0,1.1) ;

\draw[line width=1.2]
	(0.5,0.7) -- (0.5,0)
	(0.5,-0.7) -- (0,-1.1)
	(-0.5,0.7) -- (-1,1.1)	
	(1,1.8) -- (0,1.8);

\node at (0.2,1.6) {\small $\delta$};
\node at (0.2,1.2) {\small $\beta$};
\node at (0.8,1.6) {\small $\epsilon$};
\node at (0.5,0.95) {\small $\alpha$};
\node at (0.8,1.2) {\small $\gamma$}; 

\node at (-1.3,0.6) {\small $\beta$};
\node at (-1,0.85) {\small $\delta$};
\node at (-1.3,0.2) {\small $\alpha$};
\node at (-0.7,0.6) {\small $\epsilon$};	
\node at (-0.7,0.2) {\small $\gamma$};

\node at (-0.3,-0.6) {\small $\gamma$};
\node at (0,-0.85) {\small $\epsilon$};
\node at (-0.3,-0.2) {\small $\alpha$};
\node at (0.3,-0.6) {\small $\delta$};	
\node at (0.3,-0.2) {\small $\beta$}; 

\node at (-0.3,0.6) {\small $\alpha$};
\node at (0,0.85) {\small $\beta$};
\node at (-0.3,0.2) {\small $\gamma$};
\node at (0.3,0.6) {\small $\delta$};	
\node at (0.3,0.2) {\small $\epsilon$}; 

\node at (1.8,0.2) {\small $\alpha$};
\node at (1.8,0.5) {\small $\beta$};  
\node at (1.4,0.2) {\small $\gamma$}; 
\node at (0.7,0.5) {\small $\delta$}; 
\node at (0.7,0.2) {\small $\epsilon$};

\node at (-0.8,1.6) {\small $\gamma$};
\node at (-0.8,1.2) {\small $\epsilon$};
\node at (-0.2,1.6) {\small $\alpha$};
\node at (-0.5,0.95) {\small $\delta$};
\node at (-0.2,1.2) {\small $\beta$};

\node[inner sep=0.5,draw,shape=circle] at (0,0.4) {\small $1$};
\node[inner sep=0.5,draw,shape=circle] at (-1,0.4) {\small $4$};
\node[inner sep=0.5,draw,shape=circle] at (1.05,0.35) {\small $2$};
\node[inner sep=0.5,draw,shape=circle] at (0,-0.4) {\small $3$};
\node[inner sep=0.5,draw,shape=circle] at (-0.5,1.4) {\small $5$};
\node[inner sep=0.5,draw,shape=circle] at (0.5,1.4) {\small $6$};

\end{scope}

\foreach \a in {1,-1}
\draw[xscale=\a]
	(0,0.7) -- (0.5,1.1) -- (1,0.7) -- (1,-0.7) -- (0.5,-1.1) -- (0,-0.7) -- (0,0.7);
	
\draw
	(-1,0) -- (1,0)
	(1,-0.7) -- (2.6,-0.7)
	(1,0.7) -- (2.6,0.7) -- (2.6,-2.5) -- (-0.5,-2.5)
	(1.8,0.7) -- (1.8,-3.2) -- (-1.8,-3.2) -- (-1.8,0.7) -- (-1,0.7)
	(-1.8,-0.7) -- (-1,-0.7)
	(-1.8,-1.8) -- (1.8,-1.8)
	(0.5,-1.1) -- (0.5,-1.8)
	(-0.5,-1.1) -- (-0.5,-3.2);

\draw[line width=1.2]
	(1,0) -- (1,-0.7)
	(-1,0) -- (-1,0.7)
	(0,0.7) -- (0.5,1.1)
	(0,-0.7) -- (-0.5,-1.1)
	(0.5,-1.8) -- (1.8,-1.8)
	(1.8,0.7) -- (2.6,0.7)
	(2.6,-0.7) -- (2.6,-2.5)
	(-0.5,-2.5) -- (-0.5,-3.2)
	(-1.8,-1.8) -- (-1.8,-0.7);
	
\node at (0.8,0.2) {\small $\alpha$};	
\node at (0.2,0.2) {\small $\beta$}; 
\node at (0.8,0.6) {\small $\gamma$};
\node at (0.2,0.6) {\small $\delta$};
\node at (0.5,0.85) {\small $\epsilon$};

\node at (-0.8,0.2) {\small $\epsilon$};	
\node at (-0.2,0.2) {\small $\gamma$}; 
\node at (-0.8,0.6) {\small $\delta$};
\node at (-0.2,0.6) {\small $\alpha$};
\node at (-0.5,0.85) {\small $\beta$};

\node at (0.8,-0.2) {\small $\epsilon$};	
\node at (0.2,-0.2) {\small $\gamma$}; 
\node at (0.8,-0.6) {\small $\delta$};
\node at (0.2,-0.6) {\small $\alpha$};
\node at (0.5,-0.85) {\small $\beta$};

\node at (-0.8,-0.2) {\small $\alpha$};	
\node at (-0.2,-0.2) {\small $\gamma$}; 
\node at (-0.8,-0.6) {\small $\beta$};
\node at (-0.2,-0.6) {\small $\epsilon$};
\node at (-0.5,-0.85) {\small $\delta$};

\node at (-0.3,-1.6) {\small $\gamma$};
\node at (-0.3,-1.2) {\small $\epsilon$};
\node at (0.3,-1.6) {\small $\alpha$};
\node at (0,-0.95) {\small $\delta$};
\node at (0.3,-1.2) {\small $\beta$};

\node at (1.6,0.5) {\small $\alpha$};
\node at (1.2,0.5) {\small $\beta$};
\node at (1.6,-0.5) {\small $\gamma$};
\node at (1.2,0) {\small $\delta$};
\node at (1.2,-0.5) {\small $\epsilon$};

\node at (-1.6,0.5) {\small $\gamma$};
\node at (-1.2,0.5) {\small $\epsilon$};
\node at (-1.6,-0.5) {\small $\alpha$};
\node at (-1.2,0) {\small $\delta$};
\node at (-1.2,-0.5) {\small $\beta$};

\node at (2,0.5) {\small $\epsilon$};
\node at (2,-0.5) {\small $\gamma$};
\node at (2.4,0.5) {\small $\delta$};
\node at (2.4,-0.5) {\small $\alpha$};
\node at (2.4,0) {\small $\beta$};

\node at (-1.1,-0.9) {\small $\beta$};
\node at (-0.7,-1.2) {\small $\alpha$};
\node at (-1.6,-0.9) {\small $\delta$};
\node at (-1.6,-1.6) {\small $\epsilon$};
\node at (-0.7,-1.6) {\small $\gamma$};

\node at (1.1,-0.9) {\small $\alpha$};
\node at (0.7,-1.2) {\small $\beta$};
\node at (1.6,-0.9) {\small $\gamma$};
\node at (1.6,-1.6) {\small $\epsilon$};
\node at (0.7,-1.6) {\small $\delta$};

\node at (1.1,0.95) {\small $\gamma^2$};
\node at (1.8,0.9) {\small $\delta$};
\node at (2,-2.7) {\small $\gamma$};

\node at (2,-1.8) {\small $\alpha$};
\node at (2,-0.9) {\small $\beta$};
\node at (2,-2.3) {\small $\gamma$};
\node at (2.4,-0.9) {\small $\delta$};
\node at (2.4,-2.3) {\small $\epsilon$};

\node at (-0.3,-2.3) {\small $\alpha$};
\node at (1.6,-2.3) {\small $\beta$};
\node at (-0.3,-2) {\small $\gamma$};
\node at (1.6,-2) {\small $\delta$};
\node at (0.5,-2) {\small $\epsilon$};

\node at (-0.3,-2.7) {\small $\epsilon$};
\node at (1.6,-2.7) {\small $\gamma$};
\node at (-0.3,-3) {\small $\delta$};
\node at (1.6,-3) {\small $\alpha$};
\node at (0.5,-3) {\small $\beta$};

\node at (-1.6,-2) {\small $\alpha$};
\node at (-0.7,-2) {\small $\beta$};
\node at (-1.6,-3) {\small $\gamma$};
\node at (-0.7,-2.5) {\small $\delta$};
\node at (-0.7,-3) {\small $\epsilon$};

\node[inner sep=0.5,draw,shape=circle] at (0.5,0.4) {\small $1$};
\node[inner sep=0.5,draw,shape=circle] at (0.5,-0.4) {\small $2$};
\node[inner sep=0.5,draw,shape=circle] at (-0.5,-0.4) {\small $3$};
\node[inner sep=0.5,draw,shape=circle] at (-0.5,0.4) {\small $4$};
\node[inner sep=0.5,draw,shape=circle] at (1.55,0) {\small $5$};
\node[inner sep=0.5,draw,shape=circle] at (0,-1.4) {\small $6$};
\node[inner sep=0.5,draw,shape=circle] at (-1.55,0) {\small $7$};
\node[inner sep=0.5,draw,shape=circle] at (1.2,-1.3) {\small $8$};
\node[inner sep=0.5,draw,shape=circle] at (-1.2,-1.3) {\small $9$};
\node[inner sep=0,draw,shape=circle] at (1,-2.15) {\footnotesize $10$};
\node[inner sep=0,draw,shape=circle] at (-1.2,-2.5) {\footnotesize $11$};
\node[inner sep=0,draw,shape=circle] at (1,-2.85) {\footnotesize $12$};
\node[inner sep=0,draw,shape=circle] at (2.2,-1.3) {\footnotesize $13$};
\node[inner sep=0,draw,shape=circle] at (2.05,0) {\footnotesize $14$};


\begin{scope}[xshift=4.3cm]

\foreach \a in {1,-1}
\draw[xscale=\a]
	(0,0.7) -- (0.5,1.1) -- (1,0.7) -- (1,-0.7) -- (0.5,-1.1) -- (0,-0.7) -- (0,0.7);

\draw
	(1,0.7) -- (1.8,0.7) -- (1.8,-1.8) -- (-0.5,-1.8) -- (-0.5,-1.1)
	(0.5,-1.8) -- (0.5,-1.1)
	(-1,0) -- (1,0)
	(1,-0.7) -- (2.2,-0.7);
	
\draw[line width=1.2]
	(1,0) -- (1,-0.7)
	(-1,0) -- (-1,0.7)
	(0,0.7) -- (0.5,1.1)
	(0,-0.7) -- (-0.5,-1.1)
	(0.5,-1.8) -- (1.8,-1.8);

\node at (0.8,0.2) {\small $\alpha$};	
\node at (0.2,0.2) {\small $\gamma$}; 
\node at (0.8,0.6) {\small $\beta$};
\node at (0.2,0.6) {\small $\epsilon$};
\node at (0.5,0.85) {\small $\delta$};

\node at (-0.8,0.2) {\small $\epsilon$};	
\node at (-0.2,0.2) {\small $\gamma$}; 
\node at (-0.8,0.6) {\small $\delta$};
\node at (-0.2,0.6) {\small $\alpha$};
\node at (-0.5,0.85) {\small $\beta$};

\node at (0.8,-0.2) {\small $\epsilon$};	
\node at (0.2,-0.2) {\small $\gamma$}; 
\node at (0.8,-0.6) {\small $\delta$};
\node at (0.2,-0.6) {\small $\alpha$};
\node at (0.5,-0.85) {\small $\beta$};

\node at (-0.8,-0.2) {\small $\alpha$};	
\node at (-0.2,-0.2) {\small $\gamma$}; 
\node at (-0.8,-0.6) {\small $\beta$};
\node at (-0.2,-0.6) {\small $\epsilon$};
\node at (-0.5,-0.85) {\small $\delta$};

\node at (1.1,-0.9) {\small $\alpha$};
\node at (0.7,-1.2) {\small $\beta$};
\node at (1.6,-0.9) {\small $\gamma$};
\node at (1.6,-1.6) {\small $\epsilon$};
\node at (0.7,-1.6) {\small $\delta$};

\node at (-0.3,-1.6) {\small $\gamma$};
\node at (-0.3,-1.2) {\small $\epsilon$};
\node at (0.3,-1.6) {\small $\alpha$};
\node at (0,-0.95) {\small $\delta$};
\node at (0.3,-1.2) {\small $\beta$};

\node at (1.6,0.5) {\small $\alpha$};
\node at (1.2,0.5) {\small $\beta$};
\node at (1.6,-0.5) {\small $\gamma$};
\node at (1.2,0) {\small $\delta$};
\node at (1.2,-0.5) {\small $\epsilon$};

\node at (2,-0.5) {\small $\gamma$};
\node at (2,-0.9) {\small $\gamma$};

\node[inner sep=0.5,draw,shape=circle] at (0.5,0.4) {\small $1$};
\node[inner sep=0.5,draw,shape=circle] at (0.5,-0.4) {\small $2$};
\node[inner sep=0.5,draw,shape=circle] at (-0.5,-0.4) {\small $3$};
\node[inner sep=0.5,draw,shape=circle] at (-0.5,0.4) {\small $4$};
\node[inner sep=0.5,draw,shape=circle] at (1.55,0) {\small $5$};
\node[inner sep=0.5,draw,shape=circle] at (0,-1.4) {\small $6$};
\node[inner sep=0.5,draw,shape=circle] at (1.2,-1.3) {\small $7$};

\end{scope}

\end{tikzpicture}
\caption{Proposition \ref{ade_3b}: No $\epsilon^2\cdots$, and pentagonal subdivision tiling}
\label{ade_3bA} 
\end{figure}

\medskip

\noindent{\em Tiling}

\medskip

Based on the AVC \eqref{ade_3b_avc}, we construct the tiling. 

For $f=36$, by no $\alpha^2\cdots$, we determine $T_1,T_2,T_3,T_4$ around a vertex $\beta\gamma^3$ in the first of Figure \ref{ade_3bA}. Then $\alpha_1\epsilon_2\cdots=\alpha_2\epsilon_3\cdots=\alpha_3\epsilon_4\cdots=\alpha\delta\epsilon$ determine $T_5,T_6,T_7$. Then $\beta_2\beta_6\cdots=\beta_3\beta_7\cdots=\beta^3$ and no $\alpha^2\cdots$ determine $T_8,T_9$. Then $\alpha_6\delta_8\cdots=\alpha\delta\epsilon$ determines $T_{10}$. Then $\gamma_6\gamma_9\gamma_{10}\cdots=\beta\gamma^3$ and $\alpha_{10}\cdots=\alpha\delta\epsilon$ determine $T_{11},T_{12}$. Then $\delta_{10}\epsilon_8\cdots=\alpha\delta\epsilon$ and $\beta_{10}\gamma_{12}\cdots=\beta\gamma^3$ determine $T_{13}$. Then $\beta_{13}\gamma_5\gamma_8\cdots=\beta\gamma^3$ and no $\alpha^2\cdots$ determine $T_{14}$. Then $\beta_5\gamma_1\cdots=\beta\gamma^3$ and $\alpha_5\epsilon_{14}\cdots=\alpha\delta\epsilon$ imply $\gamma,\delta$ adjacent in a tile, a contradiction. 

For $f=24$, by no $\alpha^2\cdots$, we determine $T_1,T_2,T_3,T_4$ around a vertex $\gamma^4$ in the second of Figure \ref{ade_3bA}. Then $\alpha_1\epsilon_2\cdots=\alpha_2\epsilon_3\cdots=\alpha\delta\epsilon$ determine $T_5,T_6$. Then $\beta_2\beta_6\cdots=\beta^3$ and no $\alpha^2\cdots$ determine $T_7$. The argument can be applied to any three consecutive $\gamma$ in $\gamma^4$, and actually determines the second layer of tiles around the initial four. Then we may repeat the argument for the new $\gamma_5\gamma_7\cdots=\gamma^4$ (there are four such $\gamma^4$ on the boundary of the second layer). More repetition of the argument gives the a pentagonal subdivision of the octahedron.

The argument for $f=24$ uses only three consecutive $\gamma$ in $\gamma^4$. Therefore it can also be applied to $\gamma^5$ for the case $f=60$. The result is the pentagonal subdivision of the icosahedron. 

Geometrically, the existence of simple pentagon for the pentagonal subdivision tiling is discussed in \cite{wy3}. The particular reduction $a=b$ that gives the pentagonal subdivision tiling by almost equilateral tiles is also discussed.
\end{proof}

\subsection{$\beta\delta\epsilon$ and One Degree $3$ $\hat{b}$-Vertex}
\label{bde}

We assume $\beta\delta\epsilon$ is a vertex. The tiling also has a degree $3$ $\hat{b}$-vertex. This means one of $\alpha^3,\beta^3,\gamma^3,\alpha\beta\gamma,\alpha^2\beta,\alpha^2\gamma,\alpha\beta^2,\alpha\gamma^2,\beta^2\gamma,\beta\gamma^2$ is a vertex. 

\begin{proposition}\label{bde_abc}
There is no tiling, such that $\alpha,\beta,\gamma$ have distinct values, and $\beta\delta\epsilon,\alpha\beta\gamma$ are vertices.
\end{proposition}

\begin{proof}
The angle sums of $\beta\delta\epsilon,\alpha\beta\gamma$ and the angle sum for pentagon imply
\[
\beta=(1-\tfrac{4}{f})\pi,\,
\alpha+\gamma=\delta+\epsilon=(1+\tfrac{4}{f})\pi.
\]  

By Lemma \ref{geometry1}, we have either $\beta<\gamma$ and $\delta>\epsilon$, or $\beta>\gamma$ and $\delta<\epsilon$. For the later case, we further consider $\alpha>\beta$ and $\alpha<\beta$.

Suppose $\beta<\gamma$ and $\delta>\epsilon$. Then by $\beta\delta\epsilon$, we know $\beta\delta\cdots=\beta\delta\epsilon$, and $\gamma\delta\cdots$ is not a vertex. By $\beta\delta\epsilon$, and $\beta<\gamma$, and Lemma \ref{square}, we know $\beta^2\cdots,\gamma^2\cdots,\beta\gamma\cdots$ are $\hat{b}$-vertices. By $\alpha\beta\gamma$ and $\gamma>\beta>\frac{2}{3}\pi$, we get $\alpha<\beta<\gamma$. Then $\beta^2\cdots=\alpha^k\beta^2$, and $\beta\gamma\cdots=\alpha\beta\gamma$, and $\gamma^2\cdots$ is not a vertex. Further by no $\gamma\delta\cdots$, this implies $\gamma\cdots=\alpha\beta\gamma,\alpha^k\gamma\epsilon^l$. By $\beta\delta\epsilon$ and applying the counting lemma to $\beta,\gamma$, we know $\alpha^k\gamma\epsilon^l$ is a vertex.

If $l\ge 2$ in $\alpha^k\gamma\epsilon^l$, then by the balance lemma, we know $\delta^2\cdots$ is a vertex. By $\beta\delta\epsilon$, and $\delta+\epsilon>\pi$, and $\beta<\gamma$, and $\delta>\epsilon$, we get $R(\delta^2)<\beta,\gamma,2\delta,2\epsilon$. This implies $\delta^2\cdots=\alpha^k\delta^2$, and $\delta\thin\delta\cdots$ is not a vertex. Moreover, the AAD of $\alpha^k\delta^2$ implies that $\beta^2\cdots=\alpha^k\beta^2$ is a vertex. However, by no $\gamma\delta\cdots,\delta\thin\delta\cdots$, the AAD of $\alpha^k\beta^2$ contains $\thin^{\gamma}\alpha^{\beta}\thin^{\delta}\beta^{\alpha}\thin$. This implies a vertex $\thick^{\epsilon}\delta^{\beta}\thin^{\alpha}\beta^{\delta}\thin\cdots=\beta\delta\epsilon=\thick^{\epsilon}\delta^{\beta}\thin^{\alpha}\beta^{\delta}\thin^{\gamma}\epsilon^{\delta}\thick$, contradicting no $\gamma\delta\cdots$.

Therefore $l=0$ in $\alpha^k\gamma\epsilon^l$, and $\gamma\cdots=\alpha\beta\gamma,\alpha^k\gamma(k\ge 2)$. By applying the counting lemma to $\alpha,\gamma$, this implies $\gamma\cdots=\alpha\beta\gamma$. Then by $\beta\delta\epsilon$ and applying the counting lemma to $\beta,\gamma$, we get a contradiction. 

Therefore we cannot have $\beta<\gamma$ and $\delta>\epsilon$. By Lemma \ref{geometry1}, we get  $\beta>\gamma$ and $\delta<\epsilon$. Then by $\beta\delta\epsilon$, we get $\beta\epsilon\cdots=\beta\delta\epsilon$. 

\subsubsection*{Case. $\alpha>\beta>\gamma$}

By $\beta\delta\epsilon$, and $\alpha>\beta$, and Lemma \ref{square}, we know $\alpha^2\cdots,\beta^2\cdots,\alpha\beta\cdots$ are $\hat{b}$-vertices. By $\alpha\beta\gamma$ and $\alpha>\beta>\frac{2}{3}\pi$, we get $\alpha>\beta>\gamma$. Then $\alpha\beta\cdots=\alpha\beta\gamma$, and $\beta^2\cdots=\beta^2\gamma^k$, and $\alpha^2\cdots$ is not a vertex. 

By $\beta\delta\epsilon$, and $\alpha>\beta$, and $\delta<\epsilon$, we know $\alpha\epsilon\cdots$ is not a vertex. By $\alpha\beta\cdots=\alpha\beta\gamma$ and no $\alpha^2\cdots,\alpha\epsilon\cdots$, we get $\alpha\cdots=\alpha\beta\gamma,\alpha\gamma^k\delta^l$. By $\beta\delta\epsilon$ and applying the counting lemma to $\alpha,\beta$, we know $\alpha\gamma^k\delta^l$ is a vertex. Then by applying the counting lemma to $\alpha,\gamma$, we also know that either $\alpha\delta^l$ is a vertex, or $\alpha\cdots=\alpha\beta\gamma,\alpha\gamma\delta^l$. 

In both cases, we know $\delta^2\cdots$ is a vertex. By the balance lemma, this implies $\epsilon^2\cdots$ is a vertex. By $\beta\delta\epsilon$, and $\delta+\epsilon>\pi$, and $\alpha>\beta$, and $\delta<\epsilon$, we get $R(\epsilon^2)<\alpha,\beta,2\delta,2\epsilon$. This implies $\epsilon^2\cdots=\gamma^k\epsilon^2$, and $\epsilon\thin\epsilon\cdots$ is not a vertex. Then by no $\alpha^2\cdots,\alpha\epsilon\cdots,\epsilon\thin\epsilon\cdots$, we know $\gamma\thin\gamma\cdots$ is not a vertex. Therefore $k=1$ in $\gamma^k\epsilon^2$, and $\epsilon^2\cdots=\gamma\epsilon^2$ is a vertex. by $\gamma\epsilon^2$ and applying the counting lemma to $\alpha,\gamma$, we cannot have $\alpha\cdots=\alpha\beta\gamma,\alpha\gamma\delta^l$. Therefore $\alpha\delta^l$ is a vertex. 

The angle sums of $\alpha\beta\gamma,\beta\delta\epsilon,\gamma\epsilon^2,\alpha\delta^2$ and the angle sum for pentagon imply $f=12$, a contradiction. Therefore $\alpha\delta^2$ is not a vertex, and $\alpha\delta^l(l\ge 4)$ is a vertex. The AAD of the vertex contains $\thick^{\epsilon}\delta^{\beta}\thin^{\beta}\delta^{\epsilon}\thick$. This implies a vertex $\thin^{\alpha}\beta^{\delta}\thin^{\delta}\beta^{\alpha}\thin\cdots=\beta^2\gamma^k=\thin^{\delta}\beta^{\alpha}\thin\gamma\thin\cdots\thin\gamma\thin^{\alpha}\beta^{\delta}\thin$. This further implies a vertex $\alpha^2\cdots$, a contradiction. 

\subsubsection*{Case. $\beta>\alpha,\gamma$}

By $\beta\delta\epsilon$ and Lemma \ref{square}, we know $\beta^2\cdots$ is a $\hat{b}$-vertex. Then by $\alpha\beta\gamma$ and $\beta>\alpha,\gamma$, we know $\beta^2\cdots$ is not a vertex. This implies no $\thick\delta\thin\alpha\thin\cdots\thin\alpha\thin\delta\thick$ (including no $\thick\delta\thin\delta\thick$). 
 
By $\beta\epsilon\cdots=\beta\delta\epsilon$, we know $\epsilon^2\cdots$ has no $\beta$. By $\delta+\epsilon>\pi$ and $\delta<\epsilon$, we know $R(\epsilon^2)$ has no $\delta,\epsilon$. Then by $\alpha+\gamma>\pi$, we get $\epsilon^2\cdots=\alpha\epsilon^2,\gamma^k\epsilon^2$ if $\alpha>\gamma$, and $\epsilon^2\cdots=\alpha^k\epsilon^2,\gamma\epsilon^2$ if $\alpha<\gamma$. In particular, $\epsilon\thin\epsilon\cdots$ is not a vertex.

By $\alpha\beta\gamma$, the vertex $\alpha^k\epsilon^2$ implies $\beta+\gamma\ge 2\epsilon$. Then by $\beta>\gamma$, we get $\beta>\epsilon$. By exchanging $\alpha,\gamma$, we know $\gamma^k\epsilon^2$ also implies $\beta>\epsilon$. Therefore a vertex $\epsilon^2\cdots$ implies $\beta>\epsilon$.

\subsubsection*{Subcase. $\beta>\alpha>\gamma$}

By $\alpha\beta\gamma$, and $\beta>\alpha>\gamma$, and $\alpha+\gamma>\pi$, we know $\alpha\beta\gamma,\alpha^3,\alpha\gamma^k,\beta\gamma^k,\gamma^k$ are all the $\hat{b}$-vertices. We also know $k\ge 3$ in $\gamma^k,\alpha\gamma^k,\beta\gamma^k$.

By $\beta\delta\epsilon$, we get $R(\delta\epsilon)=\beta<\pi<\alpha+\gamma,\delta+\epsilon$. Then by $\alpha>\gamma$ and $\delta<\epsilon$, this implies $\delta\epsilon\cdots=\beta\delta\epsilon,\alpha\delta^l\epsilon,\gamma^k\delta^l\epsilon$. By no $\thick\delta\thin\alpha\thin\cdots\thin\alpha\thin\delta\thick$, we get $l=1$ in $\alpha\delta^l\epsilon$, and $2k\ge l-1$ in $\gamma^k\delta^l\epsilon$. By $\beta\delta\epsilon$ and $\alpha<\beta$, we know $\alpha\delta\epsilon$ is not a vertex. Therefore $\delta\epsilon\cdots=\beta\delta\epsilon,\gamma^k\delta^l\epsilon(2k\ge l-1)$. 

If $\epsilon^2\cdots$ is not a vertex, then by the balance lemma, we know $b$-vertices are $\beta\delta\epsilon,\gamma^k\delta\epsilon$. In particular, $\alpha\cdots$ a $\hat{b}$-vertex. Then by the list of $\hat{b}$-vertices, we get $\alpha^2\cdots=\alpha^3$. Moreover, by $\alpha\beta\gamma$ and applying the counting lemma to $\beta,\delta$, we know $\gamma^k\delta\epsilon$ is a vertex. By $\beta\delta\epsilon$, we know $k\ge 2$ in $\gamma^k\delta\epsilon$. Then by no $\alpha\epsilon\cdots,\epsilon^2\cdots$, the AAD of $\gamma^k\delta\epsilon$ contains $\thin^{\epsilon}\gamma^{\alpha}\thin^{\alpha}\gamma^{\epsilon}\thin$. This implies a vertex $\thin^{\beta}\alpha^{\gamma}\thin^{\gamma}\alpha^{\beta}\thin\cdots=\alpha^3=\thin^{\gamma}\alpha^{\beta}\thin\alpha\thin^{\beta}\alpha^{\gamma}\thin$. This further implies a vertex $\beta^2\cdots$, a contradiction.

Therefore $\epsilon^2\cdots$ is a vertex. By $\alpha>\gamma$, we know $\epsilon^2\cdots=\alpha\epsilon^2,\gamma^k\epsilon^2$. 

\subsubsection*{Subsubcase. $\alpha\epsilon^2$ is a vertex}

The angle sums of $\beta\delta\epsilon,\alpha\beta\gamma,\alpha\epsilon^2,\alpha^3$ and the angle sum for pentagon imply
\[
\alpha=\epsilon=\tfrac{2}{3}\pi,\,
\beta=(1-\tfrac{4}{f})\pi,\,
\gamma=\delta=(\tfrac{1}{3}+\tfrac{4}{f})\pi.
\]
Then we have $\pi>\beta>\alpha=\epsilon>\gamma=\delta$, contradicting Lemma \ref{geometry4}. Therefore $\alpha^3$ is not a vertex, and $\alpha\beta\gamma,\alpha\gamma^k,\beta\gamma^k,\gamma^k$ are all the $\hat{b}$-vertices.

By $\alpha\epsilon^2$ and $\alpha+\gamma=\delta+\epsilon=(1+\tfrac{4}{f})\pi$, we get $\gamma+2\delta=(1+\tfrac{12}{f})\pi>\alpha+\gamma$. Then $\alpha<2\delta$. By $\beta+\gamma+2\delta=(2+\frac{8}{f})\pi>2\pi$, we get $R(\beta\delta^2)<\gamma<\alpha,\beta,2\delta$. Then by $\delta<\epsilon$, this implies $\beta\delta^2\cdots$ is not a vertex. 

By no $\beta\delta^2\cdots$, and no $\thick\delta\thin\alpha\thin\cdots\thin\alpha\thin\delta\thick$, we know a $\delta^2$-fan has only $\alpha,\gamma$, and has at least one $\gamma$. By $\gamma+2\delta>\pi$, a $\delta^2$-fan has value $>\pi$. Then by $2\epsilon>\delta+\epsilon>\pi$, all fans have value $>\pi$. This implies at most one $b$-edge at a vertex. In particular, we get $\delta^2\cdots=\alpha^l\gamma^k\delta^2(k\ge 1)$, and $\delta\epsilon\cdots=\beta\delta\epsilon,\gamma^k\delta\epsilon$. Then by $2\alpha>\alpha+\gamma>\pi$ and $2\pi=\beta+\delta+\epsilon>\alpha+2\delta>\gamma+2\delta>\pi$, we get $\delta^2\cdots=\alpha\gamma\delta^2,\gamma^k\delta^2(k\ge 2)$. We conclude $\beta\delta\epsilon,\alpha\epsilon^2,\alpha\gamma\delta^2,\gamma^k\delta^2,\gamma^k\delta\epsilon,\gamma^k\epsilon^2$ are all the $b$-vertices. The $\hat{b}$-vertex list and the $b$-vertex list imply $\alpha^2\cdots,\delta\thin\epsilon\cdots$ are not vertices, and $\alpha\beta\cdots=\alpha\beta\gamma$. 

By no $\beta^2\cdots$, the AAD of $\alpha\gamma\delta^2$ is $\thick^{\epsilon}\delta^{\beta}\thin^{\gamma}\alpha^{\beta}\thin^{\epsilon}\gamma^{\alpha}\thin^{\beta}\delta^{\epsilon}\thick$ or $\thick^{\epsilon}\delta^{\beta}\thin^{\gamma}\alpha^{\beta}\thin^{\alpha}\gamma^{\epsilon}\thin^{\beta}\delta^{\epsilon}\thick$, and determines $T_1,T_2,T_3,T_4$ in the two pictures in Figure \ref{bde_abcA}. In the first picture, $\beta_1\epsilon_2\cdots=\beta\delta\epsilon$ determines $T_5$, and $\alpha_2\beta_3\cdots=\alpha\beta\gamma$ and no $\alpha^2\cdots$ determine $T_6$. Then we have $\alpha_6\beta_2\cdots=\alpha\beta\gamma$. This implies $\delta_2\epsilon_5\cdots=\alpha\delta\epsilon\cdots,\delta\epsilon^2\cdots$, a contradiction. In the second picture, $\alpha_2\beta_1\cdots=\alpha\beta\gamma$ and no $\delta\thin\epsilon\cdots$ determine $T_5$. Then $\beta_2\epsilon_5\cdots=\beta\delta\epsilon$ determines $T_6$. Then $\alpha_1\delta_5\cdots=\alpha\gamma\delta^2$ and the AAD $\thick^{\epsilon}\delta^{\beta}\thin^{\gamma}\alpha^{\beta}\thin^{\alpha}\gamma^{\epsilon}\thin^{\beta}\delta^{\epsilon}\thick$ (we already dismissed the first AAD) of $\alpha\gamma\delta^2$ determine $T_7$. Then $\alpha_7\beta_5\cdots=\alpha\beta\gamma$ implies $\delta_5\epsilon_6\cdots=\alpha\delta\epsilon\cdots,\delta\epsilon^2\cdots$, a contradiction. 

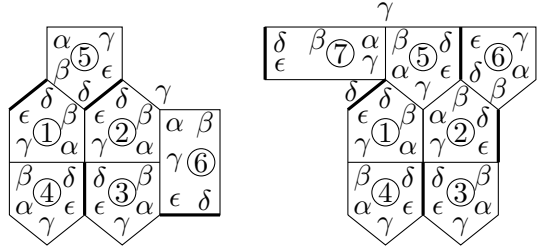
\begin{figure}[htp]
\centering
\begin{tikzpicture}[>=latex,scale=1]

\foreach \a in {0,1}
\draw[xshift=\a cm]
	(-0.5,-0.7) -- (-0.5,0.7) -- (0,1.1) -- (0.5,0.7) -- (0.5,-0.7) -- (0,-1.1) -- (-0.5,-0.7);

\draw
	(0,1.1) -- (0,1.8) -- (1,1.8) -- (1,1.1)
	(-0.5,0) -- (1.5,0)
	(1.5,0.7) -- (2.3,0.7) -- (2.3,-0.7) -- (1.5,-0.7);

\draw[line width=1.2]
	(0.5,0) -- (0.5,-0.7)
	(-0.5,0.7) -- (0,1.1)
	(0.5,0.7) -- (1,1.1)
	(2.3,-0.7) -- (1.5,-0.7);

\node at (-0.3,0.2) {\small $\gamma$};
\node at (-0.3,0.6) {\small $\epsilon$};
\node at (0.3,0.2) {\small $\alpha$}; 
\node at (0,0.85) {\small $\delta$};
\node at (0.3,0.6) {\small $\beta$};	

\node at (-0.3,-0.2) {\small $\beta$};
\node at (-0.3,-0.6) {\small $\alpha$};
\node at (0.3,-0.2) {\small $\delta$}; 
\node at (0,-0.85) {\small $\gamma$};
\node at (0.3,-0.6) {\small $\epsilon$};
	
\node at (0.7,0.6) {\small $\epsilon$};
\node at (1,0.85) {\small $\delta$};
\node at (0.7,0.2) {\small $\gamma$};
\node at (1.3,0.6) {\small $\beta$};	
\node at (1.3,0.2) {\small $\alpha$}; 
	
\node at (2.1,-0.5) {\small $\delta$};
\node at (2.1,0.5) {\small $\beta$};
\node at (1.7,-0.5) {\small $\epsilon$};
\node at (1.7,0.5) {\small $\alpha$};
\node at (1.7,0) {\small $\gamma$};

\node at (1.3,-0.6) {\small $\alpha$};	
\node at (1.3,-0.2) {\small $\beta$};
\node at (1,-0.85) {\small $\gamma$};
\node at (0.7,-0.2) {\small $\delta$}; 
\node at (0.7,-0.6) {\small $\epsilon$};

\node at (0.8,1.6) {\small $\gamma$};
\node at (0.8,1.2) {\small $\epsilon$};
\node at (0.2,1.6) {\small $\alpha$};
\node at (0.5,0.95) {\small $\delta$};
\node at (0.2,1.2) {\small $\beta$};

\node at (1.55,0.9) {\small $\gamma$};

\node[inner sep=0.5,draw,shape=circle] at (0,0.4) {\small $1$};
\node[inner sep=0.5,draw,shape=circle] at (1,0.4) {\small $2$};
\node[inner sep=0.5,draw,shape=circle] at (1,-0.4) {\small $3$};
\node[inner sep=0.5,draw,shape=circle] at (0,-0.4) {\small $4$};
\node[inner sep=0.5,draw,shape=circle] at (0.5,1.4) {\small $5$};
\node[inner sep=0.5,draw,shape=circle] at (2.05,0) {\small $6$};

\begin{scope}[xshift=4.5cm]

\foreach \a in {0,1}
\draw[xshift=\a cm]
	(-0.5,-0.7) -- (-0.5,0.7) -- (0,1.1) -- (0.5,0.7) -- (0.5,-0.7) -- (0,-1.1) -- (-0.5,-0.7);

\draw
	(0,1.1) -- (0,1.8) 
	(0,1.1) -- (-1.6,1.1) -- (-1.6,1.8) -- (2,1.8) -- (2,1.1) -- (1.5,0.7)
	(-0.5,0) -- (1.5,0);

\draw[line width=1.2]
	(1.5,0) -- (1.5,0.7)
	(0,1.1) -- (-0.5,0.7)
	(0.5,0) -- (0.5,-0.7)
	(1,1.1) -- (1,1.8)
	(-1.6,1.1) -- (-1.6,1.8);

\node at (-0.3,0.2) {\small $\gamma$};
\node at (-0.3,0.6) {\small $\epsilon$};
\node at (0.3,0.2) {\small $\alpha$}; 
\node at (0,0.85) {\small $\delta$};
\node at (0.3,0.6) {\small $\beta$};

\node at (0.7,0.6) {\small $\alpha$};
\node at (1,0.85) {\small $\beta$};
\node at (0.7,0.2) {\small $\gamma$};
\node at (1.3,0.6) {\small $\delta$};	
\node at (1.3,0.2) {\small $\epsilon$}; 

\node at (-0.3,-0.6) {\small $\alpha$};
\node at (0,-0.85) {\small $\gamma$};
\node at (-0.3,-0.2) {\small $\beta$};	
\node at (0.3,-0.2) {\small $\delta$}; 
\node at (0.3,-0.6) {\small $\epsilon$};

\node at (1.3,-0.6) {\small $\alpha$};	
\node at (1.3,-0.2) {\small $\beta$}; 
\node at (1,-0.85) {\small $\gamma$};
\node at (0.7,-0.2) {\small $\delta$};
\node at (0.7,-0.6) {\small $\epsilon$};

\node at (1.8,1.6) {\small $\gamma$};
\node at (1.8,1.2) {\small $\alpha$};
\node at (1.5,0.95) {\small $\beta$};
\node at (1.2,1.6) {\small $\epsilon$};
\node at (1.2,1.2) {\small $\delta$};

\node at (0.8,1.6) {\small $\delta$};
\node at (0.8,1.2) {\small $\epsilon$};
\node at (0.2,1.6) {\small $\beta$};
\node at (0.5,0.95) {\small $\gamma$};
\node at (0.2,1.2) {\small $\alpha$};

\node at (-0.2,1.6) {\small $\alpha$};
\node at (-0.2,1.3) {\small $\gamma$};
\node at (-0.9,1.6) {\small $\beta$};
\node at (-1.4,1.6) {\small $\delta$};
\node at (-1.4,1.3) {\small $\epsilon$};

\node at (0,2) {\small $\gamma$};
\node at (-0.45,0.95) {\small $\delta$};

\node[inner sep=0.5,draw,shape=circle] at (0,0.4) {\small $1$};
\node[inner sep=0.5,draw,shape=circle] at (1,0.4) {\small $2$};
\node[inner sep=0.5,draw,shape=circle] at (1,-0.4) {\small $3$};
\node[inner sep=0.5,draw,shape=circle] at (0,-0.4) {\small $4$};
\node[inner sep=0.5,draw,shape=circle] at (0.5,1.4) {\small $5$};
\node[inner sep=0.5,draw,shape=circle] at (1.5,1.4) {\small $6$};
\node[inner sep=0.5,draw,shape=circle] at (-0.6,1.45) {\small $7$};

\end{scope}

\end{tikzpicture}
\caption{Proposition \ref{bde_abc}: $\alpha\epsilon^2$ is a vertex.}
\label{bde_abcA}
\end{figure}

Therefore $\alpha\gamma\delta^2$ is not a vertex, and $\alpha\beta\gamma,\beta\delta\epsilon,\alpha\epsilon^2,\alpha\gamma^k,\beta\gamma^k,\gamma^k,\gamma^k\delta^2,\gamma^k\delta\epsilon,\gamma^k\epsilon^2$ are all the vertices. By $\alpha\epsilon^2$ and the balance lemma, we know $\delta^2\cdots=\gamma^k\delta^2$ is a vertex. We have
\begin{align*}
\#\delta &=
\#\beta\delta\epsilon+\#\gamma^k\delta\epsilon+
2\#\gamma^k\delta^2, \\
\#\epsilon &=
\#\beta\delta\epsilon+\#\gamma^k\delta\epsilon+
2\#\alpha\epsilon^2+2\#\gamma^k\epsilon^2.
\end{align*}
By $\#\delta=\#\epsilon$, this implies $\#\alpha\epsilon^2\le \#\gamma^k\delta^2$. Then
\begin{align*}
\#\alpha
&= \#\alpha\beta\gamma+\#\alpha\epsilon^2+\#\alpha\gamma^k \\
&\le \#\alpha\beta\gamma+\#\gamma^k\delta^2+\#\alpha\gamma^k \\
&< \#\alpha\beta\gamma+k\#\gamma^k\delta^2+k\#\alpha\gamma^k
\le \#\gamma.
\end{align*}
Here we use $\#\gamma^k\delta^2>0$ and $k\ge 2$ in $\gamma^k\delta^2$ in the second (strict) inequality. By $\#\alpha=\#\gamma$, we get a contradiction.

\subsubsection*{Subsubcase. $\alpha\epsilon^2$ is not a vertex, and $\gamma^k\epsilon^2$ is a vertex}

We know $\epsilon^2\cdots=\gamma^k\epsilon^2$. We also know $\delta\epsilon\cdots=\beta\delta\epsilon,\gamma^k\delta^l\epsilon$. Therefore $\alpha\epsilon\cdots,\epsilon\thin\epsilon\cdots$ are not vertices. This implies the AAD of $\thin\gamma\thin\gamma\thin$ is $\thin^{\epsilon}\gamma^{\alpha}\thin^{\alpha}\gamma^{\epsilon}\thin$, and further implies no consecutive $\gamma\gamma\gamma$. Therefore $\gamma^k,\alpha\gamma^k,\beta\gamma^k$ are not vertices, and $\alpha\beta\gamma,\alpha^3$ are all the $\hat{b}$-vertices.

If $k\ge 2$ in $\gamma^k\epsilon^2$, then $\gamma+\epsilon\le\pi<\alpha+\gamma=\delta+\epsilon$. This implies $\alpha>\epsilon$ and $\gamma<\delta$, contradicting to Lemma \ref{geometry4}. Therefore $\epsilon^2\cdots=\gamma\epsilon^2$ is a vertex. By $\gamma\epsilon^2$ and $\alpha+\gamma=\delta+\epsilon=(1+\frac{4}{f})\pi$, we get $\alpha+2\delta=(1+\frac{12}{f})\pi$. 

By $\alpha+2\delta>\pi$ and $\alpha+\gamma>\pi$, we get $R(\alpha^2\delta^2)<\gamma,2\delta$. By $\beta>\alpha>\gamma$, and $\delta<\epsilon$, this implies $\alpha^2\delta^2\cdots=\alpha^2\delta^2$. The AAD of $\alpha^2\delta^2$ implies a vertex $\beta^2\cdots$, a contradiction. Therefore $\alpha^2\delta^2\cdots$ is not a vertex.

By $\alpha+\beta+2\delta=(2+\frac{8}{f})\pi>2\pi$, and no $\alpha^2\delta^2\cdots,\alpha\epsilon\cdots$, we know $R(\alpha\delta^2)$ has no $\alpha,\beta,\epsilon$. Therefore $\alpha\delta^2\cdots=\alpha\gamma^k\delta^l$. By no $\thick\delta\thin\alpha\thin\cdots\thin\alpha\thin\delta\thick$, we get $2k\ge l$ in $\alpha\gamma^k\delta^l$. Combined with all $\hat{b}$-vertices and no $\alpha\epsilon\cdots$, we get $\alpha\cdots=\alpha\beta\gamma,\alpha^3,\alpha\gamma^k\delta^l(2k\ge l)$. By $\gamma\epsilon^2$ and applying the counting lemma to $\alpha,\gamma$, we know $\alpha^3$ is a vertex. The angle sums of $\beta\delta\epsilon,\alpha\beta\gamma,\gamma\epsilon^2,\alpha^3$ and the angle sum for pentagon imply
\[
\alpha=\tfrac{2}{3}\pi,\;
\beta=(1-\tfrac{4}{f})\pi,\;
\gamma=(\tfrac{1}{3}+\tfrac{4}{f})\pi,\;
\delta=(\tfrac{1}{6}+\tfrac{6}{f})\pi,\;
\epsilon=(\tfrac{5}{6}-\tfrac{2}{f})\pi.
\]

The AAD $\thin^{\epsilon}\gamma^{\alpha}\thin^{\alpha}\gamma^{\epsilon}\thin$ of $\thin\gamma\thin\gamma\thin$ implies a vertex $\thin^{\beta}\alpha^{\gamma}\thin^{\gamma}\alpha^{\beta}\thin\cdots=\alpha^3=\thin^{\gamma}\alpha^{\beta}\thin\alpha\thin^{\beta}\alpha^{\gamma}\thin$, contradicting no $\beta^2\cdots$. Therefore $\gamma\thin\gamma\cdots$ is not a vertex. 

By $\gamma\epsilon^2$ and Lemma \ref{fbalance}, there is $\delta^2$-fan. By no $\thick\delta\thin\alpha\thin\cdots\thin\alpha\thin\delta\thick$, a $\delta^2$-fan has $\beta$ or $\gamma$. By $R(\beta\delta^2)=(\frac{2}{3}-\frac{8}{f})\pi<\alpha,\beta,2\gamma$, a $\delta^2$-fan with $\beta$ is $\thick\delta\thin\beta\thin\delta\thick,\thick\delta\thin\beta\thin\gamma\thin\delta\thick$. By no $\gamma\thin\gamma\cdots$ and $R(\gamma\delta^2)=(\frac{4}{3}-\frac{12}{f})\pi<2\alpha,\alpha+2\gamma,4\gamma$, a $\delta^2$-fan with no $\beta$ is $\thick\delta\thin\gamma\thin\delta\thick,\thick\delta\thin\alpha\thin\gamma\thin\delta\thick,\thick\delta\thin\gamma\thin\alpha\thin\gamma\thin\delta\thick$. Altogether we get five possible $\delta^2$-fans $\thick\delta\thin\beta\thin\delta\thick,\thick\delta\thin\gamma\thin\delta\thick,\thick\delta\thin\alpha\thin\gamma\thin\delta\thick,\thick\delta\thin\beta\thin\gamma\thin\delta\thick,\thick\delta\thin\gamma\thin\alpha\thin\gamma\thin\delta\thick$, with respective values $(\frac{4}{3}+\frac{8}{f})\pi,(\frac{2}{3}+\frac{16}{f})\pi,(\frac{4}{3}+\frac{16}{f})\pi,(\frac{5}{3}+\frac{12}{f})\pi,(\frac{5}{3}+\frac{20}{f})\pi$. 

By $\delta\epsilon\cdots=\beta\delta\epsilon,\gamma^k\delta^l\epsilon$, if a vertex contains $\delta^2$-fan and $\epsilon$, then it is $\gamma^k\delta^l\epsilon(l\ge 3)$. By $2\gamma+3\delta+\epsilon=(2+\frac{24}{f})\pi>2\pi$, the vertex is $\gamma\delta^3\epsilon=\thick\delta\thin\epsilon\thick\delta\thin\gamma\thin\delta\thick$. 

A vertex contains $\delta^2$-fan and without $\epsilon$ is a combination of $\delta^2$-fans. By the values of $\delta^2$-fans, such vertices are $\alpha\gamma\delta^2,\beta\gamma\delta^2,\alpha\gamma^2\delta^2,\gamma^2\delta^4$. Moreover, we know $\gamma^2\delta^4=\thick\delta\thin\gamma\thin\delta\thick\delta\thin\gamma\thin\delta\thick$. 

We know one of $\alpha\gamma\delta^2,\beta\gamma\delta^2,\alpha\gamma^2\delta^2,\gamma^2\delta^4,\gamma\delta^3\epsilon$ is a vertex. The angle sum of one of the vertices implies the corresponding $f$ value, and we may further get all the vertices
\begin{align*}
f=24 &\colon
\text{AVC}=\{\beta\delta\epsilon,\alpha\beta\gamma,\alpha^3,\gamma\epsilon^2,\alpha\gamma\delta^2\}. \\
f=36 &\colon
\text{AVC}=\{\beta\delta\epsilon,\alpha\beta\gamma,\alpha^3,\gamma\epsilon^2,\beta\gamma\delta^2\}. \\
f=48 &\colon
\text{AVC}=\{\beta\delta\epsilon,\alpha\beta\gamma,\alpha^3,\gamma\epsilon^2,\gamma^2\delta^4\}. \\
f=60 &\colon
\text{AVC}=\{\beta\delta\epsilon,\alpha\beta\gamma,\alpha^3,\gamma\epsilon^2,\alpha\gamma^2\delta^2,\gamma\delta^3\epsilon\}.
\end{align*}

We first consider $f=36$. The AAD $\thick^{\delta}\epsilon^{\gamma}\thin^{\alpha}\gamma^{\epsilon}\thin^{\gamma}\epsilon^{\delta}\thick$ of $\gamma\epsilon^2$ determines $T_1,T_2,T_3$ in Figure \ref{bde_abcB}. Then $\gamma_2\epsilon_3\cdots=\gamma\epsilon^2$ determines $T_4$. Then $\alpha_2\gamma_4\cdots=\alpha_3\gamma_1\cdots=\alpha\beta\gamma$ and no $\alpha\delta\cdots$ determine $T_5,T_6$. Then $\delta_1\delta_2\cdots=\beta\gamma\delta^2$ and $\beta_2\delta_5\cdots=\beta\delta\epsilon,\beta\gamma\delta^2$  determine $T_7$. Then by no $\alpha\delta\cdots$, we further determine $T_8$. We note that this proves the AAD of $\beta\gamma\delta^2$ is $\thick^{\epsilon}\delta^{\beta}\thin^{\delta}\beta^{\alpha}\thin^{\alpha}\gamma^{\epsilon}\thin^{\beta}\delta^{\epsilon}\thick$.

By $\alpha_1\alpha_6\cdots=\alpha_4\alpha_5\cdots=\alpha_7\alpha_8\cdots=\alpha^3$ and no $\beta^2\cdots$, we determine $T_9,T_{10},T_{11}$. Then the AAD of $\beta_1\gamma_9\delta_8\cdots=\beta\gamma\delta^2$ determines $T_{12}$. Then $\beta_{11}\gamma_8\cdots=\alpha\beta\gamma,\beta\gamma\delta^2$ and $\epsilon_8\epsilon_{12}\cdots=\gamma\epsilon^2$ determine $T_{13}$. Since $\thin^{\delta}\beta_{13}{}^{\alpha}\thin^{\beta}\delta_{11}{}^{\epsilon}\thick$ is incompatible with the AAD of $\beta\gamma\delta^2$, we get $\beta_{13}\delta_{11}\cdots=\beta\delta\epsilon$. This determines $T_{14}$. Then $\delta_7\epsilon_5\cdots=\delta_{14}\epsilon_{11}\cdots=\beta\delta\epsilon$ and $\beta_7\gamma_{11}\cdots=\alpha\beta\gamma,\beta\gamma\delta^2$ determine $T_{15},T_{16}$. Then $\gamma_{16}\delta_{10}\cdots=\beta\gamma\delta^2$. However, $\thin^{\epsilon}\gamma_{16}{}^{\alpha}\thin^{\beta}\delta_{10}{}^{\epsilon}\thick$ is incompatible with the AAD of $\beta\gamma\delta^2$. We get a contradiction.  

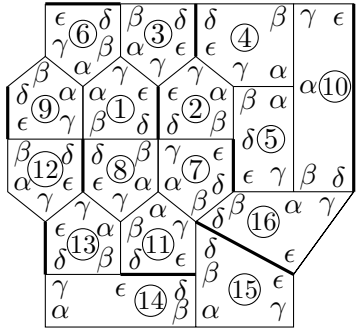
\begin{figure}[htp]
\centering
\begin{tikzpicture}[>=latex,scale=1]


\foreach \a in {1,-1}
{
\begin{scope}[scale=\a]

\foreach \b in {-1,0,1}
\draw[xshift=\b cm]
	(-0.5,0) -- (-0.5,0.7) -- (0,1.1) -- (0.5,0.7) -- (0.5,0)
	(0,1.1) -- (0,1.8);

\draw[line width=1.2]
	(0.5,0) -- (0.5,0.7)
	(-1.5,0) -- (-1.5,0.7)
	(1,1.1) -- (1,1.8)
	(0,1.8) -- (-1,1.8);	
	
\node at (-1.3,0.6) {\small $\delta$};
\node at (-1.05,0.85) {\small $\beta$};
\node at (-1.3,0.2) {\small $\epsilon$};	
\node at (-0.7,0.2) {\small $\gamma$}; 
\node at (-0.7,0.6) {\small $\alpha$};

\node at (-0.3,0.2) {\small $\beta$};
\node at (-0.3,0.6) {\small $\alpha$};
\node at (0.3,0.2) {\small $\delta$}; 
\node at (0,0.85) {\small $\gamma$};
\node at (0.3,0.6) {\small $\epsilon$};

\node at (1.3,0.6) {\small $\alpha$};	
\node at (1.3,0.2) {\small $\beta$};
\node at (1,0.85) {\small $\gamma$};
\node at (0.7,0.2) {\small $\delta$}; 
\node at (0.7,0.6) {\small $\epsilon$};

\node at (0.8,1.6) {\small $\delta$};
\node at (0.8,1.2) {\small $\epsilon$};
\node at (0.5,0.95) {\small $\gamma$};
\node at (0.2,1.6) {\small $\beta$};
\node at (0.2,1.2) {\small $\alpha$};

\node at (-0.2,1.6) {\small $\delta$};
\node at (-0.2,1.2) {\small $\beta$};
\node at (-0.8,1.6) {\small $\epsilon$};
\node at (-0.5,0.95) {\small $\alpha$};
\node at (-0.8,1.2) {\small $\gamma$};

\end{scope}
}

\draw
	(-1.5,0) -- (1.5,0)
	(1.5,0.7) -- (2.3,0.7)
	(2.3,-0.7) -- (2.3,1.8)
	(0,1.8) -- (3.1,1.8) -- (3.1,-0.7) -- (1.5,-0.7)
	(3.1,-0.7) -- (2.3,-1.8) -- (3.1,-0.7)
	(2.3,-1.8) -- (2.3,-2.5) -- (-1,-2.5) -- (-1,-1.8) -- (0,-1.8)
	(1,-2.5) -- (1,-1.8);

\draw[line width=1.2]
	(1,-1.1) -- (2.3,-1.8)
	(3.1,1.8) -- (3.1,-0.7);	
	
\node at (2.1,0.5) {\small $\alpha$};
\node at (2.1,-0.5) {\small $\gamma$};
\node at (1.7,0.5) {\small $\beta$};
\node at (1.7,-0.5) {\small $\epsilon$};
\node at (1.7,0) {\small $\delta$};

\node at (2.5,0.7) {\small $\alpha$};
\node at (2.5,-0.5) {\small $\beta$};
\node at (2.5,1.6) {\small $\gamma$};
\node at (2.9,-0.5) {\small $\delta$};
\node at (2.9,1.6) {\small $\epsilon$};

\node at (1.2,1.6) {\small $\delta$};
\node at (1.2,1.2) {\small $\epsilon$};
\node at (1.6,0.9) {\small $\gamma$};
\node at (2.1,0.9) {\small $\alpha$};
\node at (2.1,1.6) {\small $\beta$};

\node at (2.75,-0.9) {\small $\gamma$};
\node at (1.55,-0.9) {\small $\beta$};
\node at (2.3,-0.9) {\small $\alpha$};
\node at (1.3,-1.05) {\small $\delta$};
\node at (2.25,-1.55) {\small $\epsilon$};

\node at (1.2,-2.3) {\small $\alpha$};
\node at (1.2,-1.8) {\small $\beta$};
\node at (2.1,-2.3) {\small $\gamma$};
\node at (1.2,-1.45) {\small $\delta$};
\node at (2.1,-1.9) {\small $\epsilon$};

\node at (-0.8,-2.3) {\small $\alpha$};
\node at (0.8,-2.3) {\small $\beta$};
\node at (-0.8,-2) {\small $\gamma$};
\node at (0.8,-2) {\small $\delta$};
\node at (0,-2) {\small $\epsilon$};

\node[inner sep=0.5,draw,shape=circle] at (0,0.4) {\small $1$};
\node[inner sep=0.5,draw,shape=circle] at (1,0.4) {\small $2$};
\node[inner sep=0.5,draw,shape=circle] at (0.5,1.4) {\small $3$};
\node[inner sep=0.5,draw,shape=circle] at (1.65,1.35) {\small $4$};
\node[inner sep=0.5,draw,shape=circle] at (2,0) {\small $5$};
\node[inner sep=0.5,draw,shape=circle] at (-0.5,1.4) {\small $6$};
\node[inner sep=0.5,draw,shape=circle] at (1,-0.4) {\small $7$};
\node[inner sep=0.5,draw,shape=circle] at (0,-0.4) {\small $8$};
\node[inner sep=0.5,draw,shape=circle] at (-1,0.4) {\small $9$};
\node[inner sep=0,draw,shape=circle] at (2.85,0.7) {\footnotesize $10$};
\node[inner sep=0,draw,shape=circle] at (0.5,-1.4) {\footnotesize $11$};
\node[inner sep=0.5,draw,shape=circle] at (-1,-0.4) {\footnotesize $12$};
\node[inner sep=0,draw,shape=circle] at (-0.5,-1.4) {\footnotesize $13$};
\node[inner sep=0,draw,shape=circle] at (0.4,-2.15) {\footnotesize $14$};
\node[inner sep=0,draw,shape=circle] at (1.65,-2) {\footnotesize $15$};
\node[inner sep=0,draw,shape=circle] at (1.9,-1.1) {\footnotesize $16$};

\end{tikzpicture}
\caption{Proposition \ref{bde_abc}: $\gamma^k\epsilon^2$ is a vertex, $f=36$.}
\label{bde_abcB}
\end{figure}

Next we consider $\gamma\delta^3\epsilon$, which is the case $f=60$. By $\alpha\delta\cdots=\alpha\gamma^2\delta^2=\thick\delta\thin\gamma\thin\alpha\thin\gamma\thin\delta\thick$, we know $\alpha\thin\delta\cdots$ is not a vertex. By $\gamma\delta^3\epsilon=\thick\delta\thin\epsilon\thick\delta\thin\gamma\thin\delta\thick$, we get $T_1,T_2,T_3$ in the first, second and third of Figure \ref{bde_abcC}. Then $\beta_1\gamma_2\cdots=\alpha\beta\gamma$ gives $\alpha_4$. There are two possible arrangements of $T_4$. The first picture shows the first arrangement. By $\alpha_2\beta_4\cdots=\alpha\beta\gamma$ and no $\alpha\thin\delta\cdots$, we determine $T_5$. Then $\alpha_5\beta_2\cdots=\alpha\beta\gamma$ implies $\delta_2\epsilon_3\cdots=\alpha\delta\epsilon\cdots,\delta\epsilon^2\cdots$, a contradiction. 

The second and third pictures show the second arrangement of $T_4$. We have $\alpha_2\gamma_4\cdots=\alpha\beta\gamma,\alpha\gamma^2\delta^2$. If $\alpha_2\gamma_4\cdots=\alpha\beta\gamma$, then by no $\alpha\epsilon\cdots$, we determine $T_5$ in the second picture. Then we get the same contradiction as in the first picture. If $\alpha_2\gamma_4\cdots=\alpha\gamma^2\delta^2$, then $\alpha\gamma^2\delta^2=\thick\delta\thin\gamma\thin\alpha\thin\gamma\thin\delta\thick$ determines $T_5$. Then $\beta_5\epsilon_4\cdots=\beta\delta\epsilon$ determines $T_6$. Then $T_1,T_4,T_6$ is the same as $T_5,T_2,T_3$ in the second picture, and we get the same contradiction. 

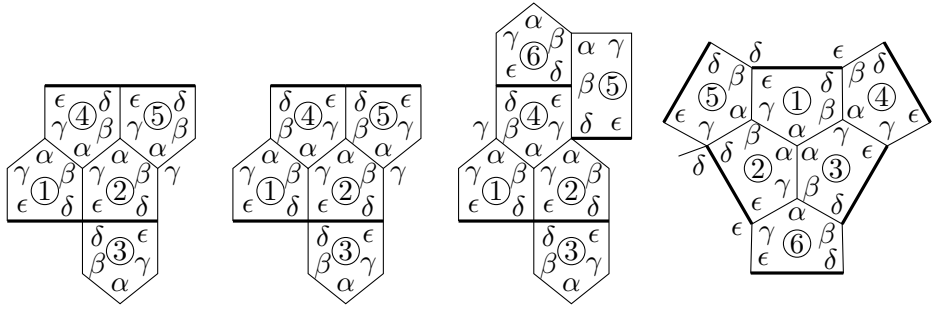
\begin{figure}[htp]
\centering
\begin{tikzpicture}[>=latex,scale=1]

\foreach \a in {0,1,2}
{
\begin{scope}[xshift=3*\a cm]

\draw
	(-1,0) -- (-1,0.7) -- (-0.5,1.1) -- (0,0.7) -- (0.5,1.1) -- (1,0.7) -- (1,-0.7) -- (0.5,-1.1) -- (0,-0.7) -- (0,0.7)
	(-0.5,1.1) -- (-0.5,1.8) -- (0.5,1.8) -- (0.5,1.1);

\draw[line width=1.2]
	(-1,0) -- (1,0)
	(-0.5,1.8) -- (0.5,1.8);

\node at (0.8,0.6) {\small $\beta$};		
\node at (0.8,0.2) {\small $\delta$};
\node at (0.5,0.85) {\small $\alpha$};
\node at (0.2,0.2) {\small $\epsilon$};
\node at (0.2,0.6) {\small $\gamma$};
	
\node at (-0.2,0.6) {\small $\beta$};		
\node at (-0.2,0.2) {\small $\delta$};
\node at (-0.5,0.85) {\small $\alpha$};
\node at (-0.8,0.2) {\small $\epsilon$};
\node at (-0.8,0.6) {\small $\gamma$};

\node at (0.8,-0.6) {\small $\gamma$};		
\node at (0.8,-0.2) {\small $\epsilon$};
\node at (0.5,-0.85) {\small $\alpha$};
\node at (0.2,-0.2) {\small $\delta$};
\node at (0.2,-0.6) {\small $\beta$};

\node[inner sep=0.5,draw,shape=circle] at (-0.5,0.4) {\small 1};
\node[inner sep=0.5,draw,shape=circle] at (0.5,0.4) {\small 2};
\node[inner sep=0.5,draw,shape=circle] at (0.5,-0.4) {\small 3};
\node[inner sep=0.5,draw,shape=circle] at (0,1.4) {\small 4};

\end{scope}
}

\foreach \a in {0,1}
{
\begin{scope}[xshift=\a cm]

\node at (0.3,1.2) {\small $\beta$};
\node at (0,0.95) {\small $\alpha$};
\node at (0.3,1.6) {\small $\delta$};
\node at (-0.3,1.2) {\small $\gamma$};
\node at (-0.3,1.6) {\small $\epsilon$};

\end{scope}
}

\foreach \a in {3,4,6}
{
\begin{scope}[xshift=\a cm]

\node at (-0.3,1.2) {\small $\beta$};
\node at (0,0.95) {\small $\alpha$};
\node at (-0.3,1.6) {\small $\delta$};
\node at (0.3,1.2) {\small $\gamma$};
\node at (0.3,1.6) {\small $\epsilon$};

\end{scope}
}


\foreach \a in {0,1}
{
\begin{scope}[xshift=3*\a cm]

\draw
	(0.5,1.8) -- (1.5,1.8) -- (1.5,1.1) -- (1,0.7);

\draw[line width=1.2]
	(0.5,1.8) -- (1.5,1.8);
	
\node at (1.2,0.6) {\small $\gamma$};

\node[inner sep=0.5,draw,shape=circle] at (1,1.4) {\small 5};

\end{scope}
}


\begin{scope}[xshift=6cm]

\draw
	(0.5,1.1) -- (1.3,1.1) -- (1.3,2.5) -- (0.5,2.5)	 		(0.5,1.8) -- (0.5,2.5) -- (0,2.9) -- (-0.5,2.5) -- (-0.5,1.8);

\draw[line width=1.2]
	(0.5,1.1) -- (1.3,1.1);

\node at (-0.3,2.4) {\small $\gamma$};
\node at (0,2.65) {\small $\alpha$};
\node at (-0.3,2) {\small $\epsilon$};
\node at (0.3,2.4) {\small $\beta$};
\node at (0.3,2) {\small $\delta$};

\node at (0.7,2.3) {\small $\alpha$};
\node at (0.7,1.8) {\small $\beta$};
\node at (1.1,2.3) {\small $\gamma$};
\node at (0.7,1.3) {\small $\delta$};
\node at (1.1,1.3) {\small $\epsilon$};

\node at (-0.7,1.2) {\small $\gamma$};

\node[inner sep=0.5,draw,shape=circle] at (1.05,1.8) {\small 5};
\node[inner sep=0.5,draw,shape=circle] at (0,2.2) {\small 6};

\end{scope}

\begin{scope}[shift={(9.5cm,1cm)}]

\foreach \a in {0,1,2}
{
\begin{scope}[rotate=120*\a]

\draw
	(0,0) -- (30:0.7) -- (60:1.2) -- (120:1.2) -- (150:0.7)
	(60:1.2) -- (50:1.8) -- (10:1.8) -- (0:1.2);

\draw[line width=1.2]
	(60:1.2) -- (120:1.2)
	(50:1.8) -- (10:1.8);
	
\node at (0,0.2) {\small $\alpha$};
\node at (0.4,0.45) {\small $\beta$};
\node at (-0.4,0.45) {\small $\gamma$};
\node at (0.4,0.85) {\small $\delta$};
\node at (-0.4,0.85) {\small $\epsilon$};

\end{scope}
}

\draw
	(180:1.2) -- ++(200:0.4);

\node at (30:0.9) {\small $\alpha$};
\node at (0.8,0.95) {\small $\beta$};
\node at (1.2,0.2) {\small $\gamma$};
\node at (1.1,1.1) {\small $\delta$};
\node at (1.55,0.4) {\small $\epsilon$};

\node at (150:0.9) {\small $\alpha$};
\node at (-0.8,0.9) {\small $\beta$};
\node at (-1.2,0.2) {\small $\gamma$};
\node at (-1.1,1.1) {\small $\delta$};
\node at (-1.55,0.4) {\small $\epsilon$};

\node at (0,-0.9) {\small $\alpha$};
\node at (0.4,-1.2) {\small $\beta$};
\node at (-0.4,-1.2) {\small $\gamma$};
\node at (0.45,-1.5) {\small $\delta$};
\node at (-0.45,-1.5) {\small $\epsilon$};

\node at (0.55,1.25) {\small $\epsilon$};
\node at (-0.55,1.25) {\small $\delta$};

\node at (-1.3,-0.25) {\small $\delta$};
\node at (-0.8,-1.1) {\small $\epsilon$};

\node[inner sep=0.5,draw,shape=circle] at (90:0.6) {\small 1};
\node[inner sep=0.5,draw,shape=circle] at (210:0.6) {\small 2};
\node[inner sep=0.5,draw,shape=circle] at (-30:0.6) {\small 3};
\node[inner sep=0.5,draw,shape=circle] at (30:1.3) {\small 4};
\node[inner sep=0.5,draw,shape=circle] at (150:1.3) {\small 5};
\node[inner sep=0.5,draw,shape=circle] at (-90:1.3) {\small 6};

\end{scope}
	
\end{tikzpicture}
\caption{Proposition \ref{bde_abc}: $\gamma^k\epsilon^2$ is a vertex, $f=24,48,60$.}
\label{bde_abcC}
\end{figure}

We conclude $\gamma\delta^3\epsilon$ is not a vertex. This implies $\delta\epsilon\cdots=\beta\delta\epsilon$, and we get updated list of vertices
\begin{align*}
f=24 &\colon
\text{AVC}=\{\beta\delta\epsilon,\alpha\beta\gamma,\alpha^3,\gamma\epsilon^2,\alpha\gamma\delta^2\}. \\
f=48 &\colon
\text{AVC}=\{\beta\delta\epsilon,\alpha\beta\gamma,\alpha^3,\gamma\epsilon^2,\gamma^2\delta^4\}. \\
f=60 &\colon
\text{AVC}=\{\beta\delta\epsilon,\alpha\beta\gamma,\alpha^3,\gamma\epsilon^2,\alpha\gamma^2\delta^2\}.
\end{align*}
We have $\gamma\delta\cdots=\alpha\gamma\delta^2,\alpha\gamma\delta^2,\alpha\gamma^2\delta^2$. This implies $\gamma\thin\delta\cdots=\thin\gamma\thin\delta\thick\delta\thin\cdots$.

Since $\gamma\epsilon^2$ is a vertex, by applying the counting lemma to $\alpha,\gamma$, we know $\alpha^3$ is a vertex. By no $\beta^2\cdots$, we know the AAD of $\alpha^3$ is $\thin^{\beta}\alpha^{\gamma}\thin^{\beta}\alpha^{\gamma}\thin^{\beta}\alpha^{\gamma}\thin$. This determines $T_1,T_2,T_3$ in the fourth Figure \ref{bde_abcC}. Then $\beta\gamma\cdots=\alpha\beta\gamma$ gives $\alpha_4,\alpha_5,\alpha_6$. Each of $T_4,T_5,T_6$ has two possible arrangements, clockwise or counterclockwise according to the orientation of $\alpha\to\beta$. 

Suppose $T_4$ has clockwise orientation, as indicated in the picture. Then by $\beta_4\delta_1\cdots\cdots=\beta\delta\epsilon$, we get $\epsilon_1\cdots=\delta\epsilon\cdots=\beta\delta\epsilon$. This implies $T_5$ has counterclockwise orientation. 

Suppose $T_5$ has counterclockwise orientation, as indicated in the picture. Then by $\gamma_5\thin\delta_2\cdots\cdots=\thin\gamma_5\thin\delta_2\thick\delta\thin\cdots$, we get $\epsilon_2\cdots=\epsilon^2\cdots=\gamma\epsilon^2$. This implies $T_6$ has counterclockwise orientation. 

By the rotation symmetry of the AAD of $\alpha^3$, we just proved that any two among $T_4,T_5,T_6$ have different orientation. Since the number of tiles is the odd number $3$, this is a contradiction. 

\subsubsection*{Subcase. $\beta>\gamma>\alpha$}

By $\alpha\beta\gamma$, and $\beta>\gamma>\alpha$, and $\alpha+\gamma>\pi$, we know $\alpha\beta\gamma,\gamma^3,\alpha^k,\alpha^k\beta,\alpha^k\gamma$ are all the $\hat{b}$-vertices. We also know $k\ge 3$ in $\alpha^k,\alpha^k\beta,\alpha^k\gamma$.

By $\beta\delta\epsilon$, we get $R(\delta\epsilon)=\beta<\pi<\alpha+\gamma,\delta+\epsilon$. Then by $\alpha<\gamma$ and $\delta<\epsilon$, this implies $\delta\epsilon\cdots=\beta\delta\epsilon,\alpha^k\delta^l\epsilon,\gamma\delta^l\epsilon$. By no  $\thick\delta\thin\alpha\thin\cdots\thin\alpha\thin\delta\thick$, we get $l=1$ in $\alpha^k\delta^l\epsilon$ and $l\le 3$ in $\gamma\delta^l\epsilon$. By $\beta\delta\epsilon$, we get $k\ge 2$ in $\alpha^k\delta\epsilon$ and $l\ge 3$ in $\gamma\delta^l\epsilon$. Therefore $\delta\epsilon\cdots=\beta\delta\epsilon,\gamma\delta^3\epsilon,\alpha^k\delta\epsilon (k\ge 2)$. 

If $\epsilon^2\cdots$ is not a vertex, then by the balance lemma, we know $\beta\delta\epsilon,\alpha^k\delta\epsilon$ are all the $b$-vertices. In particular, $\gamma\cdots$ is a $\hat{b}$-vertex. Then by the list of $\hat{b}$-vertices, we get $\beta\gamma\cdots=\alpha\beta\gamma$. Moreover, by $\alpha\beta\gamma$ and applying the counting lemma to $\beta,\delta$, we know $\alpha^k\delta\epsilon$ is a vertex. By no $\beta^2\cdots$, the AAD of $\alpha^k\delta\epsilon$ is $\thick^{\epsilon}\delta^{\beta}\thin^{\gamma}\alpha^{\beta}\thin\cdots\thin^{\gamma}\alpha^{\beta}\thin^{\gamma}\epsilon^{\delta}\thick$. Then by $k\ge 2$, this implies a vertex $\thin^{\delta}\beta^{\alpha}\thin^{\alpha}\gamma^{\epsilon}\thin\cdots=\alpha\beta\gamma=\thin^{\alpha}\gamma^{\epsilon}\thin\alpha\thin^{\delta}\beta^{\alpha}\thin$. This implies $\gamma\cdots$ is a $b$-vertex, a contradiction.

Therefore $\epsilon^2\cdots$ is a vertex. By $\alpha<\gamma$, we know $\epsilon^2\cdots=\alpha^k\epsilon^2,\gamma\epsilon^2$.

\subsubsection*{Subsubcase. $\alpha^k\epsilon^2$ is a vertex}

The AAD of $\alpha^k\epsilon^2$ implies $\gamma^2\cdots$ is a vertex. 

Suppose there is a $\hat{b}$-vertex $\gamma^2\cdots$. Then by the list of $\hat{b}$-vertices, the vertex is $\gamma^3$. The angle sum of $\gamma^3$ further implies $\gamma=\tfrac{2}{3}\pi$ and $\alpha=(\tfrac{1}{3}+\tfrac{4}{f})\pi$. Then by $2\epsilon>\delta+\epsilon=(1+\tfrac{4}{f})\pi$, we get $k=1,2$ in $\alpha^k\epsilon^2$. The angle sums of $\gamma^3,\alpha^k\epsilon^2$ further imply
\begin{align*}
\gamma^3,\alpha\epsilon^2 &\colon
	\alpha=(\tfrac{1}{3}+\tfrac{4}{f})\pi,\,
	\beta=(1-\tfrac{4}{f})\pi,\,
	\gamma=\tfrac{2}{3}\pi,\,
	\delta=(\tfrac{1}{6}+\tfrac{6}{f})\pi,\,
	\epsilon=(\tfrac{5}{6}-\tfrac{2}{f})\pi. \\
\gamma^3,\alpha^2\epsilon^2 &\colon
	\alpha=(\tfrac{1}{3}+\tfrac{4}{f})\pi,\,
	\beta=(1-\tfrac{4}{f})\pi,\,
	\gamma=\tfrac{2}{3}\pi,\,
	\delta=(\tfrac{1}{3}+\tfrac{8}{f})\pi,\,
	\epsilon=(\tfrac{2}{3}-\tfrac{4}{f})\pi.
\end{align*}
We get $\alpha<2\delta$. 

Suppose there is a $b$-vertex $\gamma^2\cdots$. Then by the lists of $\delta\epsilon\cdots,\epsilon^2\cdots$, the vertex is $\gamma^2\cdots=\gamma^2\delta^2\cdots$, with no $\epsilon$ in the remainder. The vertex implies $\gamma+\delta\le \pi<\delta+\epsilon$. Then $\gamma<\epsilon$, and $\alpha+\epsilon>\alpha+\gamma>\pi$. This implies $k=1$ in $\alpha^k\epsilon^2$. By $\alpha\epsilon^2$ and $\alpha+\gamma=\delta+\epsilon=(1+\tfrac{4}{f})\pi$, we get $\gamma+2\delta=(1+\tfrac{12}{f})\pi$. Then by $\alpha+\gamma>\pi$ and $\gamma+2\delta>\pi$, we get $R(\gamma^2\delta^2)<\alpha,2\delta$. By $\beta>\gamma>\alpha$, and no $\epsilon$ in $\gamma^2\delta^2\cdots$, this implies $\gamma^2\delta^2\cdots=\gamma^2\delta^2$. Then by $\alpha\beta\gamma,\gamma^2\delta^2$ and $\beta>\gamma$, we get $\alpha<2\delta$. 

Therefore we always have $\alpha<2\delta$. Then $\beta+\gamma+2\delta>\alpha+\beta+\gamma=2\pi$. By $\beta<\pi$, this implies $\gamma+2\delta>\pi$. Then by $\delta+\epsilon>\pi$, we know $\gamma\delta^3\epsilon$ is not a vertex.

We also proved that a $\hat{b}$-vertex $\gamma^2\cdots$ implies $k=1,2$ in $\alpha^k\epsilon^2$, and a $b$-vertex $\gamma^2\cdots$ implies $k=1$ in $\alpha^k\epsilon^2$. If $k=2$, we further know $\gamma^3,\alpha^2\epsilon^2$ are vertices, and obtained the corresponding angle values. The angle values imply $\gamma+2\epsilon=(2-\tfrac{8}{f})\pi<2\pi$. Therefore $\gamma\epsilon^2$ is not a vertex. If $k=1$, then $\alpha\epsilon^2$ and $\alpha<\gamma$ imply that $\gamma\epsilon^2$ is also not a vertex.

We proved $\gamma\epsilon^2,\gamma\delta^3\epsilon$ are not vertices. We also know $k=1,2$ in $\alpha^k\epsilon^2$. Then by the lists of $\delta\epsilon\cdots,\epsilon^2\cdots$, we get $\epsilon\cdots=\beta\delta\epsilon,\alpha\epsilon^2,\alpha^2\epsilon^2,\alpha^k\delta\epsilon(k\ge 2)$. This implies $\delta^2\cdots$ is a combination of $\delta^2$-fans. By no $\thick\delta\thin\alpha\thin\cdots\thin\alpha\thin\delta\thick$, a $\delta^2$-fan has $\beta,\gamma$. Then by $\gamma+2\delta>\pi$, and $\beta>\gamma$, we know a $\delta^2$-fan has value $>\pi$. Therefore $\delta^2\cdots$ is one $\delta^2$-fan. Then by $\gamma+2\delta>\pi$, and $\alpha+\gamma>\pi$, and $\beta+\gamma+2\delta>2\pi$ and $\beta>\gamma>\alpha$, we get $\delta^2\cdots=\gamma^2\delta^2,\alpha^k\beta\delta^2(k\ge 1),\alpha^k\gamma\delta^2(k\ge 1)$. 

We conclude $\beta\delta\epsilon,\alpha\epsilon^2,\alpha^2\epsilon^2,\gamma^2\delta^2,\alpha^k\delta\epsilon,\alpha^k\beta\delta^2,\alpha^k\gamma\delta^2$ are all the $b$-vertices. We also recall $\alpha\beta\gamma,\gamma^3,\alpha^k,\alpha^k\beta,\alpha^k\gamma$ are all the $\hat{b}$-vertices. Then we get $\beta\gamma\cdots=\alpha\beta\gamma$, and $\beta\epsilon\cdots=\beta\delta\epsilon$, and we know $\alpha\gamma^2\cdots,\beta^2\cdots,\gamma\epsilon\cdots,\delta\thin\epsilon\cdots,\epsilon\thin\epsilon\cdots$ are not vertices. By no $\epsilon\thin\epsilon\cdots$, the AAD of $\gamma^3$ is $\thin^{\alpha}\gamma^{\epsilon}\thin^{\alpha}\gamma^{\epsilon}\thin^{\alpha}\gamma^{\epsilon}\thin$. 

By no $\beta^2\cdots$, the AAD of ${}^{\beta}\thin\alpha\thin\alpha\thin$ is ${}^{\beta}\thin^{\gamma}\alpha^{\beta}\thin^{\gamma}\alpha^{\beta}\thin$. This determines $T_1,T_2$, and gives $\beta$ just outside $\gamma_1$ in the first of Figure \ref{bde_abcD}. Then $\beta\gamma_1\cdots=\beta_1\gamma_2\cdots=\alpha\beta\gamma$ and no $\gamma\epsilon\cdots$ determine $T_3,T_4$. Then $\beta_3\epsilon_1\cdots=\beta\delta\epsilon$ implies $\gamma_4\delta_1\cdots=\gamma\delta\epsilon\cdots$, a contradiction. Therefore there is no ${}^{\beta}\thin\alpha\thin\alpha\thin$. This implies no $\thick\delta\thin\alpha\thin\alpha\thin$. Then $\alpha^k\delta\epsilon$ is not a vertex, and we get $\alpha\epsilon\cdots=\epsilon^2\cdots=\alpha\epsilon^2,\alpha^2\epsilon^2$.

\begin{figure}[htp]
\centering
\begin{tikzpicture}[>=latex]


\draw
	(0.5,-0.7) -- (0,-1.1) -- (-0.5,-0.7) -- (-0.5,0.7) -- (0,1.1) -- (0.5,0.7) -- (0.5,-0.7)
	(-0.5,0) -- (0.5,0)
	(0.5,0.7) -- (1,1.1) -- (1,1.8) -- (0,1.8) -- (0,1.1)
	(-0.5,-0.7) -- (-1.3,-0.7) -- (-1.3,0.7) -- (-0.5,0.7);

\draw[line width=1.2]
	(-1.3,-0.7) -- (-1.3,0.7)
	(1,1.8) -- (0,1.8)
	(0,1.1) -- (-0.5,0.7)
	(0,-1.1) -- (-0.5,-0.7);

\node at (0.3,0.6) {\small $\gamma$};		
\node at (0.3,0.2) {\small $\alpha$};
\node at (0,0.85) {\small $\epsilon$};
\node at (-0.3,0.2) {\small $\beta$};
\node at (-0.3,0.6) {\small $\delta$};	

\node at (0.3,-0.6) {\small $\beta$};		
\node at (0.3,-0.2) {\small $\alpha$};
\node at (0,-0.85) {\small $\delta$};
\node at (-0.3,-0.2) {\small $\gamma$};
\node at (-0.3,-0.6) {\small $\epsilon$};	

\node at (0.7,0.6) {\small $\beta$};
\node at (-0.2,1.2) {\small $\delta$};	

\node at (0.8,1.2) {\small $\gamma$};
\node at (0.5,0.95) {\small $\alpha$};
\node at (0.8,1.6) {\small $\epsilon$};
\node at (0.2,1.2) {\small $\beta$};
\node at (0.2,1.6) {\small $\delta$};

\node at (-0.7,0) {\small $\alpha$};
\node at (-1.1,0.5) {\small $\epsilon$};
\node at (-1.1,-0.5) {\small $\delta$};
\node at (-0.7,0.5) {\small $\gamma$};
\node at (-0.7,-0.5) {\small $\beta$};

\node[inner sep=0.5,draw,shape=circle] at (0,0.4) {\small 1};
\node[inner sep=0.5,draw,shape=circle] at (0,-0.4) {\small 2};
\node[inner sep=0.5,draw,shape=circle] at (0.5,1.4) {\small 3};
\node[inner sep=0.5,draw,shape=circle] at (-1.05,0) {\small 4};


\begin{scope}[xshift=3.2cm]

\draw
	(-1,2.5) -- (-1,1.1) -- (-0.5,0.7) -- (0,1.1) -- (0.5,0.7) -- (1,1.1) -- (1.5,0.7) -- (2,1.1) -- (2,1.8) -- 	(-1,1.8)	
	(1,1.8) -- (2,2) -- (2,3) -- (1,3) -- (1,2.5) -- (1,1.1)
	(1,2.5) -- (-1.8,2.5) -- (-1.8,0.3) -- (-1,0.3)
	(-1.8,1.1) -- (-1,1.1)
	(0,1.8) -- (0,1.1)
	(-0.5,0.7) -- (-0.5,-0.7) -- (0,-1.1) -- (0.5,-0.7)  -- (0.5,0.7)
	(1.5,0.7) -- (1.5,0) -- (-0.5,0)
	(-0.5,0.7) -- (-1,0.3) -- (-1.8,-0.7) -- (-0.7,-0.7) -- (-0.5,0);

\draw[line width=1.2]
	(-0.5,0) -- (0.5,0)
	(1.5,0.7) -- (1,1.1)
	(-1,0.3) -- (-1.8,-0.7)
	(-0.5,0.7) -- (-1,1.1)
	(0,1.8) -- (1,1.8)
	(2,3) -- (1,3)
	(-1,2.5) -- (-1.8,2.5);

\node at (0.3,0.6) {\small $\beta$};		
\node at (0.3,0.2) {\small $\delta$};
\node at (0,0.85) {\small $\alpha$};
\node at (-0.3,0.2) {\small $\epsilon$};
\node at (-0.3,0.6) {\small $\gamma$};	

\node at (1.3,0.6) {\small $\delta$};		
\node at (1.3,0.2) {\small $\beta$};
\node at (1,0.85) {\small $\epsilon$};
\node at (0.7,0.2) {\small $\alpha$};
\node at (0.7,0.6) {\small $\gamma$};	

\node at (0.3,-0.6) {\small $\beta$};		
\node at (0.3,-0.2) {\small $\delta$};
\node at (0,-0.85) {\small $\epsilon$};
\node at (-0.3,-0.2) {\small $\epsilon$};
\node at (-0.3,-0.6) {\small $\gamma$};		

\node at (-0.7,0) {\small $\alpha$};
\node at (-0.9,0.2) {\small $\epsilon$};
\node at (-1.45,-0.5) {\small $\delta$};
\node at (-0.7,0.3) {\small $\gamma$};
\node at (-0.8,-0.5) {\small $\beta$};

\node at (-0.8,1.2) {\small $\epsilon$};
\node at (-0.5,0.95) {\small $\delta$};
\node at (-0.8,1.6) {\small $\gamma$};
\node at (-0.2,1.2) {\small $\beta$};
\node at (-0.2,1.6) {\small $\alpha$};

\node at (-1.6,0.9) {\small $\gamma$};
\node at (-1.6,0.5) {\small $\alpha$};
\node at (-0.8,0.7) {\small $\delta$};
\node at (-1.05,0.9) {\small $\epsilon$};
\node at (-1.05,0.5) {\small $\beta$};

\node at (0.8,1.2) {\small $\beta$};
\node at (0.5,0.95) {\small $\alpha$};
\node at (0.8,1.6) {\small $\delta$};
\node at (0.2,1.2) {\small $\gamma$};
\node at (0.2,1.6) {\small $\epsilon$};

\node at (1.8,1.2) {\small $\gamma$};
\node at (1.5,0.95) {\small $\epsilon$};
\node at (1.8,1.6) {\small $\alpha$};
\node at (1.2,1.2) {\small $\delta$};
\node at (1.2,1.6) {\small $\beta$};

\node at (1.8,2.8) {\small $\delta$};
\node at (1.2,2.8) {\small $\epsilon$};
\node at (1.8,2.15) {\small $\beta$};
\node at (1.2,2.5) {\small $\gamma$};
\node at (1.2,2) {\small $\alpha$};

\node at (0.8,2.3) {\small $\beta$};
\node at (-0.8,2.3) {\small $\alpha$};
\node at (0.8,2) {\small $\delta$};
\node at (-0.8,2) {\small $\gamma$};
\node at (0,1.95) {\small $\epsilon$};

\node at (-1,2.7) {\small $\epsilon$};
\node at (0.8,2.7) {\small $\alpha$};

\node at (-1.2,1.8) {\small $\gamma$};
\node at (-1.6,2.3) {\small $\delta$};
\node at (-1.6,1.3) {\small $\beta$};
\node at (-1.2,2.3) {\small $\epsilon$};
\node at (-1.2,1.3) {\small $\alpha$};

\node[inner sep=0.5,draw,shape=circle] at (1,0.4) {\small 1};
\node[inner sep=0.5,draw,shape=circle] at (0,0.4) {\small 2};
\node[inner sep=0.5,draw,shape=circle] at (0,-0.4) {\small 3};
\node[inner sep=0.5,draw,shape=circle] at (0.5,1.4) {\small 4};
\node[inner sep=0.5,draw,shape=circle] at (-1.05,-0.25) {\small 5};
\node[inner sep=0.5,draw,shape=circle] at (-0.5,1.4) {\small 6};
\node[inner sep=0.5,draw,shape=circle] at (0,2.25) {\small 8};
\node[inner sep=0.5,draw,shape=circle] at (-1.35,0.7) {\small 7};
\node[inner sep=0.5,draw,shape=circle] at (-1.55,1.8) {\small 9};
\node[inner sep=0,draw,shape=circle] at (1.5,1.4) {\footnotesize 10};
\node[inner sep=0,draw,shape=circle] at (1.5,2.4) {\footnotesize 11};

\end{scope}


\begin{scope}[xshift=7cm]

\foreach \a in {0,-1}
\draw[xshift=\a cm]
	(0.5,-0.7) -- (0,-1.1) -- (-0.5,-0.7) -- (-0.5,0.7) -- (0,1.1) -- (0.5,0.7) -- (0.5,-0.7);

\draw	
	(-1,1.1) -- (-1,1.8) -- (0,1.8) -- (0,1.1) 
	(0.5,0.7) -- (1.3,0.7) -- (1.3,-0.7) -- (0.5,-0.7) 
	(-1.5,0) -- (-0.5,0);

\draw[line width=1.2]
	(-0.5,0) -- (0.5,0)
	(-0.5,0.7) -- (-1,1.1)	
	(-0.5,-0.7) -- (-1,-1.1)
	(1.3,0.7) -- (1.3,-0.7);

\node at (-1.3,0.6) {\small $\beta$};
\node at (-1,0.85) {\small $\delta$};
\node at (-1.3,0.2) {\small $\alpha$};
\node at (-0.7,0.6) {\small $\epsilon$};	
\node at (-0.7,0.2) {\small $\gamma$};

\node at (-1.3,-0.6) {\small $\beta$};
\node at (-1,-0.85) {\small $\delta$};
\node at (-1.3,-0.2) {\small $\alpha$};
\node at (-0.7,-0.6) {\small $\epsilon$};	
\node at (-0.7,-0.2) {\small $\gamma$};

\node at (-0.3,-0.6) {\small $\beta$};
\node at (0,-0.85) {\small $\alpha$};
\node at (-0.3,-0.2) {\small $\delta$};
\node at (0.3,-0.6) {\small $\gamma$};	
\node at (0.3,-0.2) {\small $\epsilon$}; 

\node at (-0.3,0.6) {\small $\beta$};
\node at (0,0.85) {\small $\alpha$};
\node at (-0.3,0.2) {\small $\delta$};
\node at (0.3,0.6) {\small $\gamma$};	
\node at (0.3,0.2) {\small $\epsilon$}; 

\node at (1.1,-0.5) {\small $\delta$};
\node at (1.1,0.5) {\small $\epsilon$}; 
\node at (0.7,-0.5) {\small $\beta$};
\node at (0.7,0.5) {\small $\gamma$}; 
\node at (0.7,0) {\small $\alpha$}; 

\node at (-0.8,1.6) {\small $\gamma$};
\node at (-0.8,1.2) {\small $\epsilon$};
\node at (-0.2,1.6) {\small $\alpha$};
\node at (-0.5,0.95) {\small $\delta$};
\node at (-0.2,1.2) {\small $\beta$};

\node[inner sep=0.5,draw,shape=circle] at (0,-0.4) {\small $4$};
\node[inner sep=0.5,draw,shape=circle] at (0,0.4) {\small $3$};
\node[inner sep=0.5,draw,shape=circle] at (-1,0.4) {\small $1$};
\node[inner sep=0.5,draw,shape=circle] at (-1,-0.4) {\small $2$};
\node[inner sep=0.5,draw,shape=circle] at (-0.5,1.4) {\small $6$};
\node[inner sep=0.5,draw,shape=circle] at (1.05,0) {\small $5$};

\end{scope}


\begin{scope}[xshift=11cm]

\foreach \a in {0,-1}
\draw[xshift=\a cm]
	(0.5,-0.7) -- (0,-1.1) -- (-0.5,-0.7) -- (-0.5,0.7) -- (0,1.1) -- (0.5,0.7) -- (0.5,-0.7);

\draw	
	(-1,1.1) -- (-1,1.8) -- (0,1.8) -- (0,1.1) 
	(-1.5,0.7) -- (-2.3,0.7) -- (-2.3,-0.7) -- (-1.5,-0.7) 
	(-1.5,0) -- (-0.5,0);

\draw[line width=1.2]
	(-0.5,0) -- (0.5,0)
	(-1.5,0.7) -- (-1.5,0)	
	(-0.5,-0.7) -- (-1,-1.1)
	(0,1.1) -- (0,1.8);

\node at (-1.6,0.95) {\small $\gamma^2$}; 

\node at (-1.3,0.6) {\small $\delta$};
\node at (-1,0.85) {\small $\beta$};
\node at (-1.3,0.2) {\small $\epsilon$};
\node at (-0.7,0.6) {\small $\alpha$};	
\node at (-0.7,0.2) {\small $\gamma$};

\node at (-1.3,-0.6) {\small $\beta$};
\node at (-1,-0.85) {\small $\delta$};
\node at (-1.3,-0.2) {\small $\alpha$};
\node at (-0.7,-0.6) {\small $\epsilon$};	
\node at (-0.7,-0.2) {\small $\gamma$};

\node at (-0.3,-0.6) {\small $\beta$};
\node at (0,-0.85) {\small $\alpha$};
\node at (-0.3,-0.2) {\small $\delta$};
\node at (0.3,-0.6) {\small $\gamma$};	
\node at (0.3,-0.2) {\small $\epsilon$}; 

\node at (-0.3,0.6) {\small $\beta$};
\node at (0,0.85) {\small $\alpha$};
\node at (-0.3,0.2) {\small $\delta$};
\node at (0.3,0.6) {\small $\gamma$};	
\node at (0.3,0.2) {\small $\epsilon$}; 

\node at (-0.8,1.6) {\small $\beta$};
\node at (-0.8,1.2) {\small $\alpha$};
\node at (-0.2,1.6) {\small $\delta$};
\node at (-0.5,0.95) {\small $\gamma$};
\node at (-0.2,1.2) {\small $\epsilon$};

\node at (-1.2,1.2) {\small $\gamma$};

\node at (-2.1,-0.5) {\small $\alpha$};
\node at (-2.1,0.5) {\small $\beta$}; 
\node at (-1.7,-0.5) {\small $\gamma$};
\node at (-1.7,0.5) {\small $\delta$}; 
\node at (-1.7,0) {\small $\epsilon$}; 

\node[inner sep=0.5,draw,shape=circle] at (0,-0.4) {\small $4$};
\node[inner sep=0.5,draw,shape=circle] at (0,0.4) {\small $3$};
\node[inner sep=0.5,draw,shape=circle] at (-1,0.4) {\small $1$};
\node[inner sep=0.5,draw,shape=circle] at (-1,-0.4) {\small $2$};
\node[inner sep=0.5,draw,shape=circle] at (-0.5,1.4) {\small $6$};
\node[inner sep=0.5,draw,shape=circle] at (-2.05,0) {\small $5$};

\end{scope}

\end{tikzpicture}
\caption{Proposition \ref{bde_abc}: $\alpha^k\epsilon^2$ is a vertex.}
\label{bde_abcD}
\end{figure}
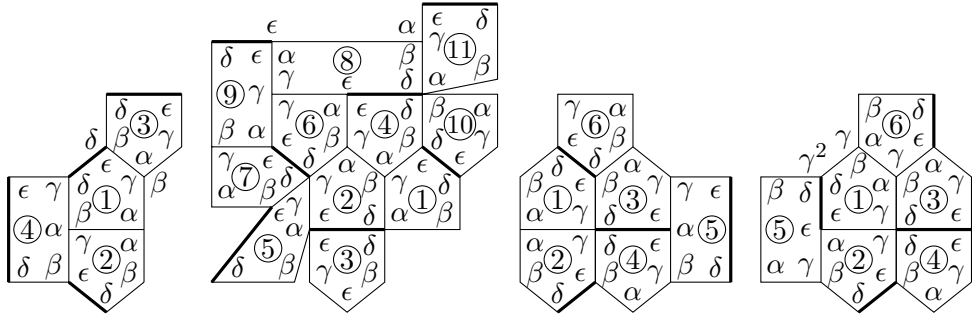

By no $\beta^2\cdots$, the AAD of $\thin\alpha\thin\delta\thick\delta\thin$ is $\thin^{\beta}\alpha^{\gamma}\thin^{\beta}\delta^{\epsilon}\thick^{\epsilon}\delta^{\beta}\thin$. This determines $T_1,T_2,T_3$ in the second of Figure \ref{bde_abcD}. Then $\beta_2\gamma_1\cdots=\alpha\beta\gamma$ and no $\gamma\epsilon\cdots$ determine $T_4$, and $\epsilon_2\epsilon_3\cdots=\alpha\epsilon^2,\alpha^2\epsilon^2$ gives $\alpha_5$. If $T_5$ is not arranged as indicated, then $\gamma_2\cdots=\beta\gamma\cdots=\alpha\beta\gamma$. By no $\alpha\gamma^2\cdots$, this implies $\alpha_2\gamma_4\cdots=\alpha\beta\gamma$, and further  implies $\epsilon_4\cdots=\delta\thin\epsilon\cdots$, a contradiction. Therefore $T_5$ is arranged as indicated. Then $\gamma_2\gamma_5\cdots=\gamma^3,\gamma^2\delta^2$. If $\gamma_2\gamma_5\cdots=\gamma^3$, then the AAD $\thin^{\alpha}\gamma^{\epsilon}\thin^{\alpha}\gamma^{\epsilon}\thin^{\alpha}\gamma^{\epsilon}\thin$ of $\gamma^3$ implies $\alpha_2\gamma_4\cdots=\alpha\gamma\epsilon\cdots$, a contradiction. Therefore $\gamma_2\gamma_5\cdots=\gamma^2\delta^2$. This determines $T_6,T_7$. We also note that the vertex $\gamma^2\delta^2$ implies $k=1$ in $\alpha^k\epsilon^2$, and $\alpha\epsilon\cdots=\epsilon^2\cdots=\alpha\epsilon^2$. Therefore $\alpha_6\epsilon_4\cdots=\alpha\epsilon^2$, and this determines $T_8$. Then $\epsilon_6\epsilon_7\cdots=\alpha\epsilon^2$ and $\gamma_6\gamma_8\cdots=\gamma^3,\gamma^2\delta^2$ determine $T_9$, and $\beta_4\epsilon_1\cdots=\beta\delta\epsilon$ determine $T_{10}$. Then $\beta_{10}\delta_4\delta_8\cdots=\alpha^k\beta\delta^2$ and no $\beta^2\cdots$ determine $T_{11}$. Then $\alpha_8\epsilon_9\cdots=\alpha\epsilon^2$ and $\beta_8\gamma_{11}\cdots=\alpha\beta\gamma$ imply $\alpha,\epsilon$ adjacent, a contradiction. 

Therefore there is no $\thin\alpha\thin\delta\thick\delta\thin$. This implies $\alpha^k\beta\delta^2,\alpha^k\gamma\delta^2$ are not vertices. We also know $\alpha^k\delta\epsilon$ is not a vertex. Therefore $\beta\delta\epsilon,\alpha\epsilon^2,\alpha^2\epsilon^2,\gamma^2\delta^2$ are all the $b$-vertices. Since $\alpha^k\epsilon^2$ is a vertex, by the balance lemma, we know $\gamma^2\delta^2$ is a vertex. This implies $k=1$ in $\alpha\epsilon^2$, and $\beta\delta\epsilon,\alpha\epsilon^2,\gamma^2\delta^2$ are all the $b$-vertices.

By the AAD of $\gamma^3$, we know $\thin^{\epsilon}\gamma^{\alpha}\thin^{\alpha}\gamma^{\epsilon}\thin\cdots$ is not $\gamma^3$. Therefore $\thin^{\epsilon}\gamma^{\alpha}\thin^{\alpha}\gamma^{\epsilon}\thin\cdots=\gamma^2\delta^2=\thick^{\epsilon}\delta^{\beta}\thin^{\epsilon}\gamma^{\alpha}\thin^{\alpha}\gamma^{\epsilon}\thin^{\beta}\delta^{\epsilon}\thick$. This determines $T_1,T_2,T_3,T_4$ in the third of Figure \ref{bde_abcD}. Then $\epsilon_3\epsilon_4\cdots=\alpha\epsilon^2$ gives $\alpha_5$. Up to the symmetry of vertical flip, we may assume $T_5$ is arranged as indicated. Then $\beta_3\epsilon_1\cdots=\beta\delta\epsilon$ determines $T_6$, and $\gamma_3\gamma_5\cdots=\gamma^3,\gamma^2\delta^2$. If $\gamma_3\gamma_5\cdots=\gamma^3$, then by the AAD of $\gamma^3$, we get $\alpha_3\beta_6\cdots=\alpha\beta\epsilon\cdots$, a contradiction. If $\gamma_3\gamma_5\cdots=\gamma^2\delta^2$, then $\alpha_3\beta_6\cdots=\alpha\beta^2\cdots$, also a contradiction. Therefore $\thin^{\epsilon}\gamma^{\alpha}\thin^{\alpha}\gamma^{\epsilon}\thin\cdots$ is not a vertex. This implies no AAD $\thin^{\beta}\alpha^{\gamma}\thin^{\gamma}\alpha^{\beta}\thin$. Combined with no ${}^{\beta}\thin\alpha\thin\alpha\thin$, we know there is no consecutive $\alpha\alpha\alpha$. Therefore $\alpha^k,\alpha^k\beta,\alpha^k\gamma$ are not vertices, and $\beta\delta\epsilon,\alpha\beta\gamma,\gamma^3,\alpha\epsilon^2,\gamma^2\delta^2$ are all the vertices. We note that $\alpha\beta\cdots=\alpha\beta\gamma$, and $\alpha^2\cdots$ is not a vertex.

We know $\gamma^2\delta^2$ is a vertex. By no $\alpha^2\cdots,\epsilon\thin\epsilon\cdots$, the AAD of the vertex is $\thick^{\epsilon}\delta^{\beta}\thin^{\alpha}\gamma^{\epsilon}\thin^{\alpha}\gamma^{\epsilon}\thin^{\beta}\delta^{\epsilon}\thick$. This determines $T_1,T_2,T_3,T_4$ in the fourth of Figure \ref{bde_abcD}. Then $\alpha_2\epsilon_1\cdots=\alpha\epsilon^2$ determines $T_5$, and $\alpha_1\beta_3\cdots=\alpha\beta\gamma$ and no $\alpha^2\cdots$ determine $T_6$. Then $\delta_1\delta_5\cdots=\gamma^2\delta^2$ and $\alpha_6\beta_1\cdots=\alpha\beta\gamma$ imply two $\gamma$ adjacent, a contradiction. 
 
\subsubsection*{Subsubcase. $\alpha^k\epsilon^2$ is not a vertex, and $\gamma\epsilon^2$ is a vertex}

By $\gamma\epsilon^2$ and $\alpha+\gamma=\delta+\epsilon=(1+\tfrac{4}{f})\pi$, we get $\alpha+2\delta=(1+\tfrac{12}{f})\pi>\alpha+\gamma$. Then $\gamma<2\delta$. By $\alpha+\beta+2\delta=(2+\frac{8}{f})\pi>2\pi$, we get $R(\beta\delta^2)<\alpha<\beta,\gamma,2\delta,2\epsilon$. This implies $\beta\delta^2\cdots$ is not a vertex. 

By no $\beta\delta^2\cdots$ and no $\thick\delta\thin\alpha\thin\cdots\thin\alpha\thin\delta\thick$, we know a $\delta^2$-fan has only $\alpha,\gamma$, and has at least one $\gamma$. By $\gamma+2\delta>\alpha+2\delta>\pi$, a $\delta^2$-fan has value $>\pi$. Then by $2\epsilon>\delta+\epsilon>\pi$, all fans have value $>\pi$. This implies at most one $b$-edge at a vertex. In particular, $\delta^2\cdots$ is one $\delta^2$-fan, which we know has only $\alpha,\gamma$, and has at least one $\gamma$. Then by $\alpha<\gamma$, and $\alpha+\gamma>\pi$, and $\alpha+2\delta>\pi$, we get $\delta^2\cdots=\alpha\gamma\delta^2,\gamma^2\delta^2$. Here by $\gamma\epsilon^2$, we know $\gamma\delta^2$ is not a vertex. Therefore $\beta\delta\epsilon,\gamma\epsilon^2,\alpha\gamma\delta^2,\gamma^2\delta^2,\alpha^k\delta\epsilon$ are all the $b$-vertices. This implies $\epsilon\thin\epsilon\cdots$ is not a vertex.  

By no $\beta^2\cdots$, the AAD of $\alpha\epsilon\cdots=\alpha^k\delta\epsilon$ is $\thick^{\epsilon}\delta^{\beta}\thin^{\gamma}\alpha^{\beta}\thin\cdots\thin^{\gamma}\alpha^{\beta}\thin^{\gamma}\epsilon^{\delta}\thick$. This implies the AAD of $\thin\alpha\thin\epsilon\thick$ is $\thin^{\gamma}\alpha^{\beta}\thin^{\gamma}\epsilon^{\delta}\thick$. This further implies no AAD $\thin^{\alpha}\gamma^{\epsilon}\thin^{\alpha}\gamma^{\epsilon}\thin$. Then by no $\epsilon\thin\epsilon\cdots$, the AAD of $\thin\gamma\thin\gamma\thin$ is $\thin^{\epsilon}\gamma^{\alpha}\thin^{\alpha}\gamma^{\epsilon}\thin$. This implies $\gamma^3$ is not a vertex.

Suppose $\gamma^2\delta^2$ is a vertex. The angle sums of $\gamma\epsilon^2,\gamma^2\delta^2$ further imply
\[
\alpha=(\tfrac{1}{3}+\tfrac{20}{3f})\pi,\,
\beta=(1-\tfrac{4}{f})\pi,\,
\gamma=(\tfrac{2}{3}-\tfrac{8}{3f})\pi,\,
\delta=(\tfrac{1}{3}+\tfrac{8}{3f})\pi,\,
\epsilon=(\tfrac{2}{3}+\tfrac{4}{3f})\pi.
\]
By the AAD $\thin^{\epsilon}\gamma^{\alpha}\thin^{\alpha}\gamma^{\epsilon}\thin$ of $\thin\gamma\thin\gamma\thin$, we get the AAD $\thick^{\epsilon}\delta^{\beta}\thin^{\epsilon}\gamma^{\alpha}\thin^{\alpha}\gamma^{\epsilon}\thin^{\beta}\delta^{\epsilon}\thick$ of $\gamma^2\delta^2$. This implies a vertex $\thin^{\beta}\alpha^{\gamma}\thin^{\gamma}\alpha^{\beta}\thin\cdots$. By no $\beta^2\cdots$, we know the vertex is not $\alpha^k,\alpha^k\delta\epsilon$. By the angle values, we know $\alpha^k\beta$ is not a vertex. Therefore the vertex is $\alpha^k\gamma$. By the angle values, we get $k=3$ in $\alpha^k\gamma$, and the angle sum of $\alpha^3\gamma$ further implies
\[
\alpha=\tfrac{6}{13}\pi,\,
\beta=\tfrac{12}{13}\pi,\,
\gamma=\tfrac{8}{13}\pi,\,
\delta=\tfrac{5}{13}\pi,\,
\epsilon=\tfrac{9}{13}\pi,\,
f=52.
\]
Then we find $\beta\delta\epsilon,\alpha\beta\gamma,\gamma\epsilon^2,\gamma^2\delta^2,\alpha^3\gamma,\alpha^2\delta\epsilon$ are all the vertices.

The AAD $\thick^{\epsilon}\delta^{\beta}\thin^{\epsilon}\gamma^{\alpha}\thin^{\alpha}\gamma^{\epsilon}\thin^{\beta}\delta^{\epsilon}\thick$ of $\gamma^2\delta^2$ determines $T_1,T_2,T_3,T_4$ in Figure \ref{bde_abcE}. By $\epsilon_2\epsilon_3\cdots=\gamma\epsilon^2$ and the symmetry of vertical flip, we may get $T_5$ as indicated. Then $\beta_2\epsilon_1\cdots=\beta\delta\epsilon$ determines $T_6$. Then $\alpha_2\beta_6\cdots=\alpha\beta\gamma$ and no $\alpha_5\gamma_3\epsilon\cdots$ determine $T_7$. Then $\alpha_5\alpha_7\gamma_2\cdots=\alpha^3\gamma$. This implies either $\beta_5\cdots$ or $\beta_7\cdots$ is $\beta^2\cdots$, a contradiction.

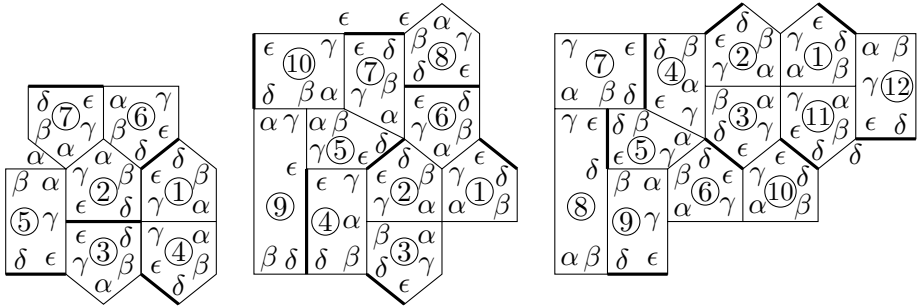
\begin{figure}[htp]
\centering
\begin{tikzpicture}[>=latex,scale=1]

\foreach \a in {0,1}
\draw[xshift=\a cm]
	(-0.5,-0.7) -- (-0.5,0.7) -- (0,1.1) -- (0.5,0.7) -- (0.5,-0.7) -- (0,-1.1) -- (-0.5,-0.7);
	
\draw
	(-0.5,-0.7) -- (-1.3,-0.7) -- (-1.3,0.7) -- (-0.5,0.7) -- (-1,1.1) -- (-1,1.8) -- (1,1.8) -- (1,1.1) -- (0.5,0.7)
	(0,1.1) -- (0,1.8)
	(0.5,0) -- (1.5,0);

\draw[line width=1.2]
	(1,1.1) -- (0.5,0.7)
	(1,-1.1) -- (0.5,-0.7)
	(0.5,0) -- (-0.5,0)
	(-1,1.8) -- (0,1.8)
	(-0.5,-0.7) -- (-1.3,-0.7);

\node at (0.7,0.6) {\small $\epsilon$};
\node at (1,0.85) {\small $\delta$};
\node at (0.7,0.2) {\small $\gamma$};
\node at (1.3,0.6) {\small $\beta$};	
\node at (1.3,0.2) {\small $\alpha$}; 

\node at (0.7,-0.6) {\small $\epsilon$};
\node at (1,-0.85) {\small $\delta$};
\node at (0.7,-0.2) {\small $\gamma$};
\node at (1.3,-0.6) {\small $\beta$};	
\node at (1.3,-0.2) {\small $\alpha$}; 

\node at (-0.3,-0.6) {\small $\gamma$};
\node at (0,-0.85) {\small $\alpha$};
\node at (-0.3,-0.2) {\small $\epsilon$};	
\node at (0.3,-0.2) {\small $\delta$}; 
\node at (0.3,-0.6) {\small $\beta$};

\node at (-0.3,0.6) {\small $\gamma$};
\node at (0,0.85) {\small $\alpha$};
\node at (-0.3,0.2) {\small $\epsilon$};	
\node at (0.3,0.2) {\small $\delta$}; 
\node at (0.3,0.6) {\small $\beta$};

\node at (0.8,1.6) {\small $\gamma$};
\node at (0.8,1.2) {\small $\epsilon$};
\node at (0.2,1.6) {\small $\alpha$};
\node at (0.5,0.95) {\small $\delta$};
\node at (0.2,1.2) {\small $\beta$};

\node at (-0.8,1.6) {\small $\delta$};
\node at (-0.8,1.2) {\small $\beta$};
\node at (-0.2,1.6) {\small $\epsilon$};
\node at (-0.5,0.95) {\small $\alpha$};
\node at (-0.2,1.2) {\small $\gamma$};

\node at (-1.1,0.5) {\small $\beta$};
\node at (-1.1,-0.5) {\small $\delta$};
\node at (-0.7,0.5) {\small $\alpha$};
\node at (-0.7,-0.5) {\small $\epsilon$};
\node at (-0.7,0) {\small $\gamma$};

\node at (-0.9,0.85) {\small $\alpha$};

\node[inner sep=0.5,draw,shape=circle] at (0,-0.4) {\small $3$};
\node[inner sep=0.5,draw,shape=circle] at (0,0.4) {\small $2$};
\node[inner sep=0.5,draw,shape=circle] at (1,0.4) {\small $1$};
\node[inner sep=0.5,draw,shape=circle] at (1,-0.4) {\small $4$};
\node[inner sep=0.5,draw,shape=circle] at (-1.05,0) {\small $5$};
\node[inner sep=0.5,draw,shape=circle] at (0.5,1.4) {\small $6$};
\node[inner sep=0.5,draw,shape=circle] at (-0.5,1.4) {\small $7$};


\begin{scope}[xshift=4.5cm]

\draw
	(0,0.7) -- (-0.5,1.1) -- (-1,0.7) -- (-1,-0.7) -- (-0.5,-1.1) -- (0,-0.7) -- (0,0.7) -- (0.5,1.1) -- (1,0.7) -- (1,0) -- (-1,0) 
	(0.5,1.1) -- (0.5,2.5) -- (0,2.9) -- (-0.5,2.5) -- (-0.5,1.1)
	(-0.5,2.5) -- (-2.5,2.5) -- (-2.5,-0.7) -- (-1,-0.7)
	(-1.3,1.5) -- (-0.5,1.1)
	(-1.3,2.5) -- (-1.3,1.5) -- (-2.5,1.5) 
	(-1.8,1.5) -- (-1.8,0.7) -- (-1,0.7);
	
\draw[line width=1.2]
	(1,0.7) -- (0.5,1.1)
	(-1,0.7) -- (-0.5,1.1)
	(-1,-0.7) -- (-0.5,-1.1)
	(-1.8,-0.7) -- (-1.8,0.7)
	(0.5,1.8) -- (-0.5,1.8)
	(-0.5,2.5) -- (-1.3,2.5)
	(-2.5,2.5) -- (-2.5,1.5);
	
\node at (0.8,0.6) {\small $\delta$};	
\node at (0.5,0.85) {\small $\epsilon$};	
\node at (0.8,0.2) {\small $\beta$};
\node at (0.2,0.6) {\small $\gamma$};	
\node at (0.2,0.2) {\small $\alpha$};

\node at (-0.8,0.6) {\small $\epsilon$};	
\node at (-0.5,0.85) {\small $\delta$};	
\node at (-0.8,0.2) {\small $\gamma$};
\node at (-0.2,0.6) {\small $\beta$};	
\node at (-0.2,0.2) {\small $\alpha$};

\node at (-0.8,-0.6) {\small $\delta$};	
\node at (-0.5,-0.85) {\small $\epsilon$};	
\node at (-0.8,-0.2) {\small $\beta$};
\node at (-0.2,-0.6) {\small $\gamma$};	
\node at (-0.2,-0.2) {\small $\alpha$};

\node at (-1.6,0.5) {\small $\epsilon$};
\node at (-1.6,-0.5) {\small $\delta$};
\node at (-1.2,0.5) {\small $\gamma$};
\node at (-1.2,-0.5) {\small $\beta$};
\node at (-1.2,0) {\small $\alpha$};

\node at (0,0.95) {\small $\alpha$};
\node at (0.3,1.2) {\small $\beta$};
\node at (-0.3,1.2) {\small $\gamma$};
\node at (0.3,1.6) {\small $\delta$};
\node at (-0.3,1.6) {\small $\epsilon$};

\node at (0,2.65) {\small $\alpha$};
\node at (-0.3,2.4) {\small $\beta$};
\node at (0.3,2.4) {\small $\gamma$};
\node at (-0.3,2) {\small $\delta$};
\node at (0.3,2) {\small $\epsilon$};

\node at (-1.05,0.85) {\small $\epsilon$};
\node at (-1.65,0.9) {\small $\gamma$};
\node at (-0.8,1.05) {\small $\delta$};
\node at (-1.65,1.3) {\small $\alpha$};
\node at (-1.35,1.3) {\small $\beta$};

\node at (-2.3,1.3) {\small $\alpha$};
\node at (-2.3,-0.5) {\small $\beta$};
\node at (-2,1.3) {\small $\gamma$};
\node at (-2,-0.5) {\small $\delta$};
\node at (-2,0.7) {\small $\epsilon$};

\node at (-1.1,2.3) {\small $\epsilon$};
\node at (-1.1,1.65) {\small $\gamma$};
\node at (-0.7,2.3) {\small $\delta$};
\node at (-0.7,1.4) {\small $\alpha$};
\node at (-0.7,1.8) {\small $\beta$};

\node at (-1.5,1.7) {\small $\alpha$};
\node at (-1.8,1.7) {\small $\beta$};
\node at (-1.5,2.3) {\small $\gamma$};
\node at (-2.3,1.7) {\small $\delta$};
\node at (-2.3,2.3) {\small $\epsilon$};

\node at (-1.3,2.7) {\small $\epsilon$};
\node at (-0.5,2.7) {\small $\epsilon$};

\node[inner sep=0.5,draw,shape=circle] at (0.5,0.4) {\small 1};
\node[inner sep=0.5,draw,shape=circle] at (-0.5,0.4) {\small 2};
\node[inner sep=0.5,draw,shape=circle] at (-0.5,-0.4) {\small 3};
\node[inner sep=0.5,draw,shape=circle] at (-1.55,0) {\small 4};
\node[inner sep=0.5,draw,shape=circle] at (-1.35,0.95) {\small 5};
\node[inner sep=0.5,draw,shape=circle] at (-2.15,0.2) {\small 9};
\node[inner sep=0.5,draw,shape=circle] at (-0.95,2) {\small 7};
\node[inner sep=0.5,draw,shape=circle] at (0,1.4) {\small 6};
\node[inner sep=0.5,draw,shape=circle] at (0,2.2) {\small 8};
\node[inner sep=0,draw,shape=circle] at (-1.9,2.1) {\footnotesize 10};

\end{scope}


\begin{scope}[shift={(9cm,1.8cm)}]

\draw
	(0,-0.7) -- (0,0.7) -- (0.5,1.1) -- (1,0.7) -- (1,-0.7) -- (0.5,-1.1) -- (0,-0.7)
	(0,0.7) -- (-0.5,1.1) -- (-1,0.7) -- (-1,-0.7) -- (-0.5,-1.1) -- (0,-0.7)
	(-1,0) -- (1,0)
	(0.5,-1.1) -- (0.5,-1.8) -- (-1.5,-1.8)
	(-2.3,-1.1) -- (-1.5,-1.1) -- (-1,-0.7) -- (-1.8,-0.3) -- (-3,-0.3)
	(-1,0.7) -- (-3,0.7) -- (-3,-2.5) -- (-1.5,-2.5) -- (-1.5,-1.1)
	(-2.3,-0.3) -- (-2.3,-2.5)
	(-0.5,-1.1) -- (-0.5,-1.8) 
	(1,0.7) -- (1.8,0.7) -- (1.8,-0.7) -- (1,-0.7) ;

\draw[line width=1.2]
	(1,0.7) -- (0.5,1.1)
	(-1,0.7) -- (-0.5,1.1)
	(-1,-0.7) -- (-0.5,-1.1)
	(-1.8,-0.3) -- (-1.8,0.7)
	(0,-0.7) -- (0.5,-1.1)
	(-2.3,-0.3) -- (-2.3,-1.1)
	(-2.3,-2.5) -- (-1.5,-2.5)
	(1.8,-0.7) -- (1,-0.7);

\node at (0.8,0.6) {\small $\delta$};	
\node at (0.5,0.85) {\small $\epsilon$};	
\node at (0.8,0.2) {\small $\beta$};
\node at (0.2,0.6) {\small $\gamma$};	
\node at (0.2,0.2) {\small $\alpha$};

\node at (-0.8,0.6) {\small $\epsilon$};	
\node at (-0.5,0.85) {\small $\delta$};	
\node at (-0.8,0.2) {\small $\gamma$};
\node at (-0.2,0.6) {\small $\beta$};	
\node at (-0.2,0.2) {\small $\alpha$};

\node at (-0.8,-0.6) {\small $\delta$};	
\node at (-0.5,-0.85) {\small $\epsilon$};	
\node at (-0.8,-0.2) {\small $\beta$};
\node at (-0.2,-0.6) {\small $\gamma$};	
\node at (-0.2,-0.2) {\small $\alpha$};

\node at (-1.6,0.5) {\small $\delta$};
\node at (-1.6,-0.2) {\small $\epsilon$};
\node at (-1.2,0.5) {\small $\beta$};
\node at (-1.2,-0.4) {\small $\gamma$};
\node at (-1.2,0) {\small $\alpha$};

\node at (1.6,0.5) {\small $\beta$};
\node at (1.6,-0.5) {\small $\delta$};
\node at (1.2,0.5) {\small $\alpha$};
\node at (1.2,-0.5) {\small $\epsilon$};
\node at (1.2,0) {\small $\gamma$};

\node at (1,-0.9) {\small $\delta$};

\node at (0.8,-0.6) {\small $\beta$};	
\node at (0.5,-0.85) {\small $\delta$};	
\node at (0.8,-0.2) {\small $\alpha$};
\node at (0.2,-0.6) {\small $\epsilon$};	
\node at (0.2,-0.2) {\small $\gamma$};

\node at (0,-0.95) {\small $\epsilon$};
\node at (0.3,-1.2) {\small $\delta$};
\node at (-0.3,-1.2) {\small $\gamma$};
\node at (0.3,-1.6) {\small $\beta$};
\node at (-0.3,-1.6) {\small $\alpha$};

\node at (-1,-0.95) {\small $\delta$};
\node at (-0.7,-1.2) {\small $\epsilon$};
\node at (-1.3,-1.2) {\small $\beta$};
\node at (-0.7,-1.6) {\small $\gamma$};
\node at (-1.3,-1.65) {\small $\alpha$};

\node at (-1.3,-0.7) {\small $\alpha$};
\node at (-1.85,-0.5) {\small $\beta$};
\node at (-1.5,-0.95) {\small $\gamma$};
\node at (-2.15,-0.5) {\small $\delta$};
\node at (-2.15,-0.95) {\small $\epsilon$};

\node at (-1.7,-1.3) {\small $\alpha$};
\node at (-1.7,-1.8) {\small $\gamma$};
\node at (-2.1,-1.3) {\small $\beta$};
\node at (-2.1,-2.3) {\small $\delta$};
\node at (-1.7,-2.3) {\small $\epsilon$};

\node at (-2.8,-0.1) {\small $\alpha$};
\node at (-2.3,-0.1) {\small $\beta$};
\node at (-2.8,0.5) {\small $\gamma$};
\node at (-2,-0.1) {\small $\delta$};
\node at (-2,0.5) {\small $\epsilon$};

\node at (-2.8,-2.3) {\small $\alpha$};
\node at (-2.5,-2.3) {\small $\beta$};
\node at (-2.8,-0.5) {\small $\gamma$};
\node at (-2.5,-1.1) {\small $\delta$};
\node at (-2.5,-0.5) {\small $\epsilon$};

\node[inner sep=0.5,draw,shape=circle] at (0.5,0.4) {\small 1};
\node[inner sep=0.5,draw,shape=circle] at (-0.5,0.4) {\small 2};
\node[inner sep=0.5,draw,shape=circle] at (-0.5,-0.4) {\small 3};
\node[inner sep=0.5,draw,shape=circle] at (-1.45,0.2) {\small 4};
\node[inner sep=0.5,draw,shape=circle] at (-1.85,-0.85) {\small 5};
\node[inner sep=0.5,draw,shape=circle] at (-1,-1.4) {\small 6};
\node[inner sep=0.5,draw,shape=circle] at (-2.4,0.3) {\small 7};
\node[inner sep=0.5,draw,shape=circle] at (-2.65,-1.6) {\small 8};
\node[inner sep=0.5,draw,shape=circle] at (-2.05,-1.8) {\small 9};
\node[inner sep=0,draw,shape=circle] at (0,-1.4) {\footnotesize 10};
\node[inner sep=0,draw,shape=circle] at (0.5,-0.4) {\footnotesize 11};
\node[inner sep=0,draw,shape=circle] at (1.55,0) {\footnotesize 12};

\end{scope}

\end{tikzpicture}
\caption{Proposition \ref{bde_abc}: $\gamma\epsilon^2$ is a vertex.}
\label{bde_abcE}
\end{figure}

Therefore $\gamma^2\delta^2$ is not a vertex. Combined with no $\alpha^k\epsilon^2,\gamma^3$, we find $\beta\delta\epsilon,\alpha\beta\gamma,\gamma\epsilon^2,\alpha\gamma\delta^2,\alpha^k,\alpha^k\gamma,\alpha^k\beta,\alpha^k\delta\epsilon$ are all the vertices. 

By no $\beta^2\cdots,\gamma^2\cdots$, we know the AAD of consecutive $\alpha\alpha\alpha$ is $\thin^{\beta}\alpha^{\gamma}\thin^{\beta}\alpha^{\gamma}\thin^{\beta}\alpha^{\gamma}\thin$. This determines $T_1,T_2,T_3$ in the second and third of Figure \ref{bde_abcE}. Then $\beta_3\gamma_2\cdots=\alpha\beta\gamma$ gives $\alpha_4$. There are two ways of arranging $T_4$, given respectively by the second and third pictures. In the second picture, $\gamma_4\epsilon_2\cdots=\gamma\epsilon^2$ determines $T_5$. Then $\delta_2\delta_5\cdots=\alpha\gamma\delta^2$, and $\beta_2\gamma_1\cdots=\alpha\beta\gamma$, and no $\beta^2\cdots$ determine $T_6,T_7$. Then $\beta_7\epsilon_6\cdots$ determines $T_8$. On the other hand, $\gamma_5\epsilon_4\cdots=\gamma\epsilon^2$ determines $T_9$. Then $\beta_5\gamma_7\cdots=\alpha\beta\gamma$ and $\alpha_5\gamma_9\cdots=\alpha\beta\gamma,\alpha\gamma\delta^2,\alpha^k\gamma$ determine $T_{10}$. Then  $\beta_8\delta_7\cdots=\beta\delta\epsilon$ and $\gamma_{10}\epsilon_7\cdots=\gamma\epsilon^2$ imply two $\epsilon$ adjacent, a contradiction. 

In the third of Figure \ref{bde_abcE}, $\gamma_4\delta_3\cdots=\alpha\gamma\delta^2$ and no $\beta^2\cdots$ determine $T_5,T_6$. Then $\beta_5\epsilon_4\cdots=\beta\delta\epsilon$ determines $T_7$. Then $\beta_7\delta_5\cdots=\beta\delta\epsilon$ determines $T_8$. Then $\beta_6\gamma_5\cdots=\alpha\beta\gamma$ and $\delta_8\epsilon_5\cdots=\beta\delta\epsilon,\alpha^k\delta\epsilon$ determine $T_9$. On the other hand, $\epsilon_3\epsilon_6\cdots=\gamma\epsilon^2$ gives $\gamma_{10}$. If $T_{10}$ is not arranged as indicated, then $\gamma_6\cdots=\gamma\epsilon\cdots=\gamma\epsilon^2$. This implies $\alpha_6\gamma_9\cdots=\alpha\gamma^2\cdots$, a contradiction. Therefore $T_{10}$ is arranged as indicated. Then $\gamma_3\epsilon_{10}\cdots=\gamma\epsilon^2$ determines $T_{11}$. 

Note that $T_{11}$ is attached to $T_3$, and we do not yet know it is attached to $T_1$. We have actually proved that the AAD $\thin^{\beta}\alpha^{\gamma}\thin^{\beta}\alpha^{\gamma}\thin^{\beta}\alpha^{\gamma}\thin$ of consecutive $\alpha\alpha\alpha$ must be extended to $\thin^{\beta}\alpha^{\gamma}\thin^{\beta}\alpha^{\gamma}\thin^{\beta}\alpha^{\gamma}\thin^{\epsilon}\gamma^{\alpha}\thin$. This implies $\alpha\thin\alpha\thin\alpha\cdots=\alpha^3\gamma$. In particular, $\alpha^k,\alpha^k\beta$ are not vertices, and $k=3$ in $\alpha^k\gamma$. Then $\alpha\beta\cdots=\alpha\beta\gamma$, and we know $T_{11}$ is attached to $T_1$. 

By $\alpha_{11}\beta_1\cdots=\alpha\beta\gamma$ and no $\delta\thin\epsilon\cdots$, we determine $T_{12}$. Then $\beta_{11}\epsilon_{12}\cdots=\beta\delta\epsilon$ implies $\delta_{10}\delta_{11}\cdots=\beta\delta^2\cdots$, a contradiction.  

Therefore there is no consecutive $\alpha\alpha\alpha$. Then $\alpha^k,\alpha^k\beta,\alpha^k\gamma$ are not vertices. By $\beta\delta\epsilon$, this also implies $k=2$ in $\alpha^k\delta\epsilon$. Then we get the updated list $\beta\delta\epsilon,\alpha\beta\gamma,\gamma\epsilon^2,\alpha\gamma\delta^2,\alpha^2\delta\epsilon$ of all the vertices. By $\gamma\epsilon^2$ and the balance lemma, we know $\alpha\gamma\delta^2$ is a vertex. By $\gamma\epsilon^2$ and applying the counting lemma to $\alpha,\gamma$, we know $\alpha^2\delta\epsilon$ is a vertex. The angle sums of the five vertices and the angle sum for pentagon imply
\[
\alpha=\delta=\tfrac{4}{9}\pi,\;
\beta=\tfrac{8}{9}\pi,\;
\gamma=\epsilon=\tfrac{2}{3}\pi,\;
f=36.
\]
By Lemma \ref{geometry11}, this implies $a=b$, a contradiction. 
\end{proof}

\begin{proposition}\label{bde_2ac}
There is no tiling, such that $\alpha,\beta,\gamma$ have distinct values, and $\beta\delta\epsilon$ is a vertex, and $2\alpha+\gamma=2\pi$.
\end{proposition}

\begin{proof}
The angle sum of $\beta\delta\epsilon$, and the equality $2\alpha+\gamma=2\pi$, and the angle sum for pentagon imply
\[
\alpha=(1-\tfrac{4}{f})\pi,\;
\gamma=\tfrac{8}{f}\pi,\;
\beta+\delta+\epsilon=2\pi.
\]
We have $\pi>\alpha>\gamma$ and $2\alpha>\pi$. 

Suppose $\beta<\gamma$ and $\delta>\epsilon$. Then by $\beta\delta\epsilon$, and $\alpha>\gamma>\beta$, and $\delta>\epsilon$, we get $\beta\delta\cdots=\beta\delta\epsilon$, and $\alpha\delta\cdots,\gamma\delta\cdots$ are not vertices. This implies $\delta^2\cdots$ has no $\alpha,\beta,\gamma$. By $\beta\delta\epsilon$, and $2\alpha+\gamma=2\pi$, and $\beta<\gamma$, we get $\delta+\epsilon>2\alpha>\pi$. Then by $\delta>\epsilon$, we know $R(\delta^2)$ has no $\delta,\epsilon$. Therefore $\delta^2\cdots$ is not a vertex. By the balance lemma, this implies a $b$-vertex is $\delta\epsilon\cdots$, with no $\delta,\epsilon$ in the remainder. By $\beta\delta\epsilon$ and $\beta<\gamma<\alpha$, we find $\beta\delta\epsilon$ is the only $b$-vertex. Then by the counting lemma, we know $\beta\delta\epsilon,\alpha^k\gamma^l$ are all the vertices. By the values of $\alpha,\gamma$, we get $\alpha^k\gamma^l=\alpha^2\gamma,\alpha\gamma^l(l\ge 3),\gamma^l(l\ge 3)$. By no $\alpha\epsilon\cdots,\epsilon^2\cdots$, we know the AAD of $\thin\gamma\thin\gamma\thin$ is $\thin^{\epsilon}\gamma^{\alpha}\thin^{\alpha}\gamma^{\epsilon}\thin$. This implies no consecutive $\gamma\gamma\gamma$. Then $\alpha\gamma^l,\gamma^l$ are not vertices, and $\beta\delta\epsilon,\alpha^2\gamma$ are all the vertices. Applying the counting lemma to $\alpha,\gamma$, we get a contradiction.

By Lemma \ref{geometry1}, we get $\beta>\gamma$ and $\delta<\epsilon$. Then by $\beta\delta\epsilon$, we get $\beta\epsilon\cdots=\beta\delta\epsilon$. 

\subsubsection*{Case. $\alpha<\beta$}

By $2\alpha+\gamma=2\pi$, and $\alpha<\beta$, we get $R(\beta^2)<R(\alpha\beta)<\gamma<\alpha<\beta$. By $\beta\delta\epsilon$ and $\delta<\epsilon$, we get $R(\alpha\beta)<\delta+\epsilon<2\epsilon$. By $\beta\delta\epsilon$ and Lemma \ref{square}, we know $\beta^2\cdots$ has no $\delta,\epsilon$. Therefore $\alpha\beta\cdots=\alpha\beta\delta^k$, and $\beta^2\cdots$ is not a vertex. By no $\beta^2\cdots$, the AAD implies $\delta\thin\delta\cdots$ is not a vertex. 

If $\alpha\beta\cdots=\alpha\beta\delta^k$ is a vertex, then by $\beta\delta\epsilon$ and $k$ even, we get $\epsilon=\alpha+(k-1)\delta>\alpha$. Then by $2\alpha+\gamma=2\pi$, we get $R(\epsilon^2)<\gamma<\alpha<\beta,\epsilon$. Therefore $\epsilon^2\cdots=\delta^k\epsilon^2$. This implies $\epsilon\cdots=\delta^k\epsilon\cdots (k\ge 1),\delta^k\epsilon^2(k\ge 2)$, with no $\epsilon$ in the remainder. Then by $\alpha\beta\delta^k$ and applying the counting lemma to $\delta,\epsilon$, we get a contradiction. Therefore $\alpha\beta\cdots$ is not a vertex.

By $2\alpha+\gamma=2\pi$, we get $R(\alpha^2)=\gamma<\alpha<\beta$. Therefore $\alpha^2\cdots=\alpha^2\gamma,\alpha^2\delta^k\epsilon^l$. By no $\beta^2\cdots$, the AAD of $\alpha^2\delta^k\epsilon^l$ implies $l\ge 1$, and $\beta\gamma\cdots$ is a vertex. By $2\alpha+\gamma=2\pi$ and $\delta<\epsilon$, the angle sum of $\alpha^2\delta^k\epsilon^l(l\ge 1)$ implies $\gamma\ge \delta+\epsilon$. Then $\beta+\gamma\ge \beta+\delta+\epsilon=2\pi$. This implies $\beta\gamma\cdots$ is not a vertex. Therefore $\alpha^2\delta^k\epsilon^l$ is not a vertex, and $\alpha^2\cdots=\alpha^2\gamma$.

The AAD $\thin^{\epsilon}\gamma^{\alpha}\thin^{\alpha}\gamma^{\epsilon}\thin$ implies a vertex $\thin^{\beta}\alpha^{\gamma}\thin^{\gamma}\alpha^{\beta}\thin\cdots=\alpha^2\gamma=\thin^{\gamma}\alpha^{\beta}\thin^{\alpha}\gamma^{\epsilon}\thin^{\beta}\alpha^{\gamma}\thin$, contradicting no $\alpha\beta\cdots$. Therefore we do not have $\thin^{\epsilon}\gamma^{\alpha}\thin^{\alpha}\gamma^{\epsilon}\thin$.

By no $\alpha\beta\cdots$, the AAD of $\thin\gamma\thin\delta\thick$ is $\thin^{\alpha}\gamma^{\epsilon}\thin^{\beta}\delta^{\epsilon}\thick$. This implies the a vertex $\thin^{\alpha}\beta^{\delta}\thin^{\gamma}\epsilon^{\delta}\thick\cdots=\beta\delta\epsilon=\thick^{\epsilon}\delta^{\beta}\thin^{\alpha}\beta^{\delta}\thin^{\gamma}\epsilon^{\delta}\thick$, contradicting no $\alpha\beta\cdots$. Therefore $\gamma\thin\delta\cdots$ is not a vertex.

By $\beta\epsilon\cdots=\beta\delta\epsilon$ and no $\alpha\beta\cdots,\beta^2\cdots,\gamma\thin\delta\cdots$, we get $\beta\cdots=\beta\delta\epsilon,\beta\gamma^k,\beta\delta^k$. The AAD of $\beta\delta^k$ contains $\thick^{\epsilon}\delta^{\beta}\thin^{\alpha}\beta^{\delta}\thin^{\beta}\delta^{\epsilon}\thick$, contradicting no $\alpha\beta\cdots$. Therefore $\beta\cdots=\beta\delta\epsilon,\beta\gamma^k$. This implies a $\delta^2$-fan has only $\alpha,\gamma$. By no $\alpha\beta\cdots,\beta^2\cdots$, the AAD of such a $\delta^2$-fan implies $\thin^{\epsilon}\gamma^{\alpha}\thin^{\alpha}\gamma^{\epsilon}\thin$, a contradiction. Therefore there is no $\delta^2$-fan. By Lemma \ref{fbalance}, we know all fans are $\delta\epsilon$-fans. In particular, $\epsilon\thin\epsilon\cdots$ is not a vertex.

By no $\thin^{\epsilon}\gamma^{\alpha}\thin^{\alpha}\gamma^{\epsilon}\thin$, and no $\epsilon\thin\epsilon\cdots$, we know the AAD of $\thin\gamma\thin\gamma\thin$ is $\thin^{\epsilon}\gamma^{\alpha}\thin^{\epsilon}\gamma^{\alpha}\thin$. The AAD implies $\thin^{\beta}\alpha^{\gamma}\thin^{\gamma}\epsilon^{\delta}\thick$. This is one side of a fan. Since all fans are $\delta\epsilon$-fans, and $\beta\epsilon\cdots=\beta\delta\epsilon$, the fan is $\thick^{\epsilon}\delta^{\beta}\thin\cdots\thin^{\beta}\alpha^{\gamma}\thin^{\gamma}\epsilon^{\delta}\thick$, where $\cdots$ consists of $\alpha,\gamma$. The AAD of this fan implies $\alpha\beta\cdots,\beta^2\cdots,\thin^{\epsilon}\gamma^{\alpha}\thin^{\alpha}\gamma^{\epsilon}\thin$, a contradiction. Therefore $\gamma\thin\gamma\cdots$ is not a vertex. 

By no $\gamma\thin\gamma\cdots$, we know $\beta\gamma^k$ is not a vertex, and $\beta\cdots=\beta\delta\epsilon$. Then by the counting lemma, we know $\beta\delta\epsilon,\alpha^k\gamma^l$ are all the vertices. By $2\alpha+\gamma=2\pi$, and $\alpha>\gamma$, and no $\gamma\thin\gamma\cdots$, we get $\alpha^k\gamma^l=\alpha^2\gamma$. Therefore $\beta\delta\epsilon,\alpha^2\gamma$ are all the vertices, contradicting the counting lemma.

\subsubsection*{Case. $\alpha>\beta$}

By $\beta\delta\epsilon$, and $\alpha>\beta$, and Lemma \ref{square}, we know $\alpha^2\cdots,\alpha\beta\cdots$ have no $\delta,\epsilon$. By $2\alpha+\gamma=2\pi$ and $\alpha>\beta>\gamma$, we get $R(\alpha^2)=\gamma<\beta<\alpha$ and $R(\alpha\beta)<\alpha$. Therefore $\alpha^2\cdots=\alpha^2\gamma$, and $\alpha\beta\cdots=\alpha\beta^k\gamma^l$.

By $\beta\delta\epsilon$ and $\beta<\alpha<\pi$, we get $\delta+\epsilon>\pi>\alpha$. By $\beta\delta\epsilon$ and $\delta<\epsilon$, we get $R(\epsilon^2)<R(\delta\epsilon)=\beta<\alpha<\delta+\epsilon<2\epsilon$. This implies $\epsilon\cdots=\beta\delta\epsilon,\gamma^k\epsilon^2,\gamma^k\delta^l\epsilon$. Then $\alpha\epsilon\cdots,\epsilon\thin\epsilon\cdots$ are not vertices, and $\delta\thin\epsilon\cdots=\gamma^k\delta^l\epsilon$. By no $\alpha\epsilon\cdots,\epsilon\thin\epsilon\cdots$, the AAD of $\thin\gamma\thin\gamma\thin$ is $\thin^{\epsilon}\gamma^{\alpha}\thin^{\alpha}\gamma^{\epsilon}\thin$. This implies no consecutive $\gamma\gamma\gamma$, and  we get $\epsilon\cdots=\beta\delta\epsilon,\gamma\epsilon^2,\gamma^2\epsilon^2,\gamma^k\delta^l\epsilon$. 

The AAD $\thin^{\epsilon}\gamma^{\alpha}\thin^{\alpha}\gamma^{\epsilon}\thin$ implies a vertex $\thin^{\beta}\alpha^{\gamma}\thin^{\gamma}\alpha^{\beta}\thin\cdots=\alpha^2\gamma=\thin^{\gamma}\alpha^{\beta}\thin^{\alpha}\gamma^{\epsilon}\thin^{\beta}\alpha^{\gamma}\thin$. This again implies a vertex $\thin^{\beta}\alpha^{\gamma}\thin^{\alpha}\beta^{\delta}\thin\cdots=\alpha\beta^k\gamma^l=\thin^{\beta}\alpha^{\gamma}\thin^{\alpha}\beta^{\delta}\thin\beta\thin\cdots,\thin^{\beta}\alpha^{\gamma}\thin^{\alpha}\beta^{\delta}\thin\gamma\thin\cdots$. This further implies a vertex $\alpha\delta\cdots,\delta\thin\delta\cdots,\delta\thin\epsilon\cdots$. By $\delta<\epsilon$, we know $\alpha\delta\cdots$ implies $\alpha+2\delta\le 2\pi$, and we know $\delta\thin\delta\cdots,\delta\thin\epsilon\cdots$ imply $4\delta<2\pi$. By $\alpha=(1-\frac{4}{f})\pi$, we get $\delta\le (\frac{1}{2}+\frac{2}{f})\pi$ in all the cases. 

If $\gamma\epsilon^2,\gamma^2\epsilon^2$ are not vertices, then we get $\epsilon\cdots=\beta\delta\epsilon,\gamma^k\delta^l\epsilon$. Applying the counting lemma to $\delta,\epsilon$, we know $b$-vertices are $\beta\delta\epsilon,\gamma^k\delta\epsilon$. By $\beta\delta\epsilon$ and no consecutive $\gamma\gamma\gamma$, we further know $b$-vertices are $\beta\delta\epsilon,\gamma^2\delta\epsilon$. 

The AAD of $\gamma^2\delta\epsilon$ contains the AAD $\thin^{\epsilon}\gamma^{\alpha}\thin^{\alpha}\gamma^{\epsilon}\thin$ of $\thin\gamma\thin\gamma\thin$. By the argument above, this implies a vertex $\alpha\delta\cdots,\delta\thin\delta\cdots,\delta\thin\epsilon\cdots$, contradicting the only $b$-vertices $\beta\delta\epsilon,\gamma^2\delta\epsilon$. Therefore $\beta\delta\epsilon$ is the only $b$-vertex. Then by the counting lemma, $\beta\delta\epsilon,\alpha^k\gamma^l$ are all the vertices. By $2\alpha+\gamma=2\pi$, and $\alpha>\gamma$, and no consecutive $\gamma\gamma\gamma$, we get $\alpha^k\gamma^l=\alpha^2\gamma$. Therefore $\beta\delta\epsilon,\alpha^2\gamma$ are all the vertices, contradicting the counting lemma.

We conclude either $\gamma\epsilon^2$ or $\gamma^2\epsilon^2$ is a vertex.

\subsubsection*{Subcase. $\gamma^2\epsilon^2$ is a vertex}

The AAD of $\gamma^2\epsilon^2$ contains the AAD $\thin^{\epsilon}\gamma^{\alpha}\thin^{\alpha}\gamma^{\epsilon}\thin$ of $\thin\gamma\thin\gamma\thin$. We know this implies $\delta\le (\frac{1}{2}+\frac{2}{f})\pi$ and a vertex $\thin^{\beta}\alpha^{\gamma}\thin^{\alpha}\beta^{\delta}\thin\cdots$. Moreover, the angle sum of $\gamma^2\epsilon^2$ implies $\epsilon=\pi-\gamma=(1-\frac{8}{f})\pi$. Then by $\beta\delta\epsilon$, we get $\beta\ge (\frac{1}{2}+\frac{6}{f})\pi$, and $\alpha+2\beta\ge (2+\frac{8}{f})\pi>2\pi$. Then $\alpha\beta\cdots=\alpha\beta^k\gamma^l$ becomes $\alpha\beta\cdots=\alpha\beta\gamma^l$. By $2\alpha+\gamma=2\pi$ and $\alpha\ne\beta$, we get $l\ge 2$ in $\alpha\beta\gamma^l$. Then by no consecutive $\gamma\gamma\gamma$, we get $\alpha\thin\beta\cdots=\alpha\beta\gamma^2$. 

By the AAD $\thin^{\epsilon}\gamma^{\alpha}\thin^{\alpha}\gamma^{\epsilon}\thin$ of $\thin\gamma\thin\gamma\thin$, we get $\thin^{\beta}\alpha^{\gamma}\thin^{\alpha}\beta^{\delta}\thin\cdots=\alpha\beta\gamma^2=\thin^{\beta}\alpha^{\gamma}\thin^{\alpha}\beta^{\delta}\thin^{\epsilon}\gamma^{\alpha}\thin^{\alpha}\gamma^{\epsilon}\thin$. This determines $T_1,T_2,T_3,T_4$ in Figure \ref{bde_2acA}. Then $\beta_2\epsilon_1\cdots=\beta\delta\epsilon$ determines $T_5$, and $\delta_4\thin\epsilon_3\cdots=\gamma^k\delta^l\epsilon$ determines $T_6$. Then $\alpha_1\alpha_3\cdots=\alpha^2\gamma$ implies either $\beta_1\cdots$ or $\beta_2\cdots$ is $\alpha\thin\beta\cdots=\alpha\beta\gamma^2$. This implies either $\delta_1\epsilon_5\cdots$ or $\delta_3\epsilon_6\cdots$ is $\alpha\delta\epsilon\cdots,\delta\epsilon^2\cdots$, a contradiction. 

\begin{figure}[htp]
\centering
\begin{tikzpicture}[>=latex,scale=1]

\foreach \a in {-1,1}
{
\begin{scope}[xscale=\a]

\draw
	(1,-0.7) -- (1,0.7) -- (0.5,1.1) -- (0,0.7) -- (0,-0.7) -- (0.5,-1.1) -- (1,-0.7)
	(0.5,1.1) -- (0.5,1.5);

\node at (0.8,0.6) {\small $\delta$};	
\node at (0.8,0.2) {\small $\epsilon$};
\node at (0.5,0.85) {\small $\beta$};
\node at (0.2,0.2) {\small $\gamma$}; 
\node at (0.2,0.6) {\small $\alpha$};

\node at (0.3,1.2) {\small $\alpha$};
\node at (0.75,1.25) {\small $\gamma^2$};

\end{scope}
}

\draw
	(-1,0.7) -- (-1.8,0.7) -- (-1.8,-0.7) -- (-1,-0.7) 
	(-1,0) -- (2,0) -- (2,0.7) -- (1.5,1.1) -- (1,0.7);

\draw[line width=1.2]
	(-0.5,-1.1) -- (-1,-0.7)
	(1,-0.7) -- (1,0.7)
	(-1,0) -- (-1,0.7);

\node at (0,0.95) {\small $\gamma$};

\node at (0.2,-0.6) {\small $\alpha$};
\node at (0.5,-0.85) {\small $\gamma$};
\node at (0.2,-0.2) {\small $\beta$};
\node at (0.8,-0.6) {\small $\epsilon$};	
\node at (0.8,-0.2) {\small $\delta$}; 

\node at (-0.2,-0.6) {\small $\gamma$};
\node at (-0.5,-0.85) {\small $\epsilon$};
\node at (-0.2,-0.2) {\small $\alpha$};
\node at (-0.8,-0.6) {\small $\delta$};	
\node at (-0.8,-0.2) {\small $\beta$}; 

\node at (1.8,0.6) {\small $\alpha$};	
\node at (1.8,0.2) {\small $\beta$};
\node at (1.5,0.85) {\small $\gamma$};
\node at (1.2,0.2) {\small $\delta$}; 
\node at (1.2,0.6) {\small $\epsilon$};

\node at (-1.6,0.5) {\small $\gamma$};
\node at (-1.6,-0.5) {\small $\alpha$};
\node at (-1.2,0.5) {\small $\epsilon$};
\node at (-1.2,-0.5) {\small $\beta$};
\node at (-1.2,0) {\small $\delta$};

\node[inner sep=0.5,draw,shape=circle] at (-0.5,0.4) {\small $1$};
\node[inner sep=0.5,draw,shape=circle] at (0.5,0.4) {\small $3$};
\node[inner sep=0.5,draw,shape=circle] at (0.5,-0.4) {\small $4$};
\node[inner sep=0.5,draw,shape=circle] at (-0.5,-0.4) {\small $2$};
\node[inner sep=0.5,draw,shape=circle] at (-1.55,0) {\small $5$};
\node[inner sep=0.5,draw,shape=circle] at (1.5,0.4) {\small $6$};

\end{tikzpicture}
\caption{Proposition \ref{bde_2ac}: $\gamma^2\epsilon^2$ is a vertex.}
\label{bde_2acA}
\end{figure}
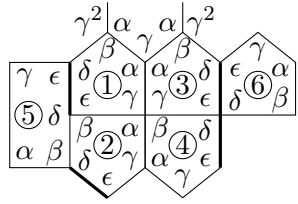

\subsubsection*{Subcase. $\gamma\epsilon^2$ is a vertex}

By $\gamma\epsilon^2$, we know $\gamma^2\epsilon^2$ is not a vertex, and $\epsilon\cdots=\beta\delta\epsilon,\gamma\epsilon^2,\gamma^k\delta^l\epsilon$.

The angle sums of $\beta\delta\epsilon,\gamma\epsilon^2$, and the equality $2\alpha+\gamma=2\pi$, and the angle sum for pentagon imply
\[
\alpha=\epsilon=(1-\tfrac{4}{f})\pi,\;
\beta+\delta=(1+\tfrac{4}{f})\pi,\;
\gamma=\tfrac{8}{f}\pi.
\]

By $\alpha+\gamma=\beta+\delta$ and $\alpha>\beta$, we get $\gamma<\delta$. Then $\alpha+\beta+\gamma<\alpha+\beta+\delta=2\pi$. This implies $\alpha\beta\gamma$ is not a vertex.

If $\beta=\delta$, then by Lemma \ref{geometry11}, we get $a=b$, a contradiction. Therefore $\beta\ne \delta$. By $\alpha+\beta+\delta=\pi$, this implies $\alpha\beta^2$ is not a vertex. 

We already have $\alpha=\epsilon>\beta,\delta$, and $\gamma<\beta,\delta$. Then by $\beta\ne\delta$, we have either $\alpha=\epsilon>\beta>\delta>\gamma$, or $\alpha=\epsilon>\delta>\beta>\gamma$.  

\subsubsection*{Subsubcase. $\alpha=\epsilon>\beta>\delta>\gamma$}

By $\beta>\delta$, we get $\alpha+2\beta>\alpha+\beta+\delta=2\pi$. By no $\alpha\beta\gamma$, this implies $\alpha\beta\cdots=\alpha\beta^k\gamma^l=\alpha\beta\gamma^l(l\ge 2)$. Then by no $\alpha\epsilon\cdots$, no consecutive $\gamma\gamma\gamma$, and the AAD $\thin^{\epsilon}\gamma^{\alpha}\thin^{\alpha}\gamma^{\epsilon}\thin$ of $\thin\gamma\thin\gamma\thin$, we get $\alpha\thin\beta\cdots=\thin\alpha\thin\gamma\thin\gamma\thin\beta\thin=\thin\alpha\thin^{\epsilon}\gamma^{\alpha}\thin^{\alpha}\gamma^{\epsilon}\thin^{\delta}\beta^{\alpha}\thin$. Therefore the AAD of $\thin\alpha\thin\beta\thin$ is $\thin\alpha\thin^{\alpha}\beta^{\delta}\thin$. 

The AAD $\thin^{\epsilon}\gamma^{\alpha}\thin^{\beta}\delta^{\epsilon}\thick$ implies $\thin^{\beta}\alpha^{\gamma}\thin^{\delta}\beta^{\alpha}\thin$, contradicting the AAD $\thin\alpha\thin^{\alpha}\beta^{\delta}\thin$ of $\thin\alpha\thin\beta\thin$. The AAD $\thin^{\alpha}\gamma^{\epsilon}\thin^{\beta}\delta^{\epsilon}\thick$ implies a vertex $\thin^{\alpha}\beta^{\delta}\thin^{\gamma}\epsilon^{\delta}\thick\cdots=\beta\delta\epsilon=\thick^{\epsilon}\delta^{\beta}\thin^{\alpha}\beta^{\delta}\thin^{\gamma}\epsilon^{\delta}\thick$. This further implies $\thin^{\gamma}\alpha^{\beta}\thin^{\delta}\beta^{\alpha}\thin$, also contradicting the AAD $\thin\alpha\thin^{\alpha}\beta^{\delta}\thin$ of $\thin\alpha\thin\beta\thin$. Therefore $\gamma\thin\delta\cdots$ is not a vertex. This implies $\delta\thin\epsilon\cdots=\gamma^k\delta^l\epsilon=\delta^l\epsilon(l\ge 3)$. 

By $R(\alpha\delta)=\beta<\alpha$ and no $\alpha\epsilon\cdots$, we get $\alpha\delta\cdots=\alpha\gamma^k\delta^l$. Then by no $\gamma\thin\delta\cdots$, we get $\alpha\delta\cdots=\alpha\delta^l$. By $\alpha+2\delta<\alpha+\beta+\delta=2\pi$, we get $l\ge 4$ in $\alpha\delta^l$.

By $\beta\delta\epsilon$ and Lemma \ref{square}, we know $\beta^2\cdots$ has no $\delta,\epsilon$. Then by $\alpha+2\beta>\alpha+\beta+\delta=2\pi$ and $2\beta>\beta+\delta>\pi$, we get $\beta^2\cdots=\beta^3\gamma^k,\beta^2\gamma^k$. By $\gamma\epsilon^2$ and $\beta<\epsilon$, we know $k\ge 2$ in $\beta^2\gamma^k$. 

The AAD $\thick^{\epsilon}\delta^{\beta}\thin^{\beta}\delta^{\epsilon}\thick$ of $\thick\delta\thin\delta\thick$ implies a vertex $\thin^{\alpha}\beta^{\delta}\thin^{\delta}\beta^{\alpha}\thin\cdots=\beta^3\gamma^k,\beta^2\gamma^k(k\ge 2)$. By the AAD $\thin^{\epsilon}\gamma^{\alpha}\thin^{\alpha}\gamma^{\epsilon}\thin$ of $\thin\gamma\thin\gamma\thin$, and no $\alpha\epsilon\cdots$, the vertex is $\beta^3,\beta^3\gamma$. The vertex implies $\beta\le \frac{2}{3}\pi$. Then by $\beta\delta\epsilon$, we get $\delta+\epsilon\ge \frac{4}{3}\pi$. Further by $\gamma\epsilon^2$, we get $\epsilon<\pi$. Then $\alpha+3\delta=3\delta+\epsilon=3(\delta+\epsilon)-2\epsilon>2\pi$. This implies $\alpha\delta\cdots=\alpha\delta^l(l\ge 4)$ and $\gamma^k\delta^l\epsilon=\delta^l\epsilon(l\ge 3)$ are not vertices. Then we get $\epsilon\cdots=\beta\delta\epsilon,\gamma\epsilon^2$, and $\delta\thin\epsilon\cdots,\epsilon\thin\epsilon\cdots$ are not vertices. By no $\alpha\delta\cdots,\alpha\epsilon\cdots,\delta\thin\epsilon\cdots$, the AAD implies $\thin^{\alpha}\beta^{\delta}\thin^{\delta}\beta^{\alpha}\thin\cdots=\beta^3,\beta^3\gamma$ is not a vertex. 

Therefore $\delta\thin\delta\cdots$ is not a vertex. This further implies $\delta\thin\epsilon\cdots=\delta^l\epsilon$ is not a vertex. 

Let $\theta=\beta,\gamma$. By no $\alpha\delta\cdots,\alpha\epsilon\cdots,\delta\thin\delta\cdots,\delta\thin\epsilon\cdots,\epsilon\thin\epsilon\cdots$, we know the AAD of $\thin\theta\thin\theta\thin$ is $\thin\theta^{\alpha}\thin^{\alpha}\theta\thin$. This implies no consecutive $\theta\theta\theta$. Then $\alpha\beta\cdots=\alpha\beta\gamma^l(l\ge 2)$ is not a vertex.

By $\gamma\epsilon^2$ and Lemma \ref{fbalance}, we know there is $\delta^2$-fan. By no $\alpha\delta\cdots,\delta\thin\delta\cdots$, we know a $\delta^2$-fan is $\thick\delta\thin\theta\thin\cdots\thin\theta\thin\delta\thick$, where the number $k$ of $\theta$ is at least $1$. By no $\alpha\beta\cdots$, the AAD implies $k\ne 1$. By no consecutive $\theta\theta\theta$, we have $k\le 2$. By $\beta+\delta>\pi$, we cannot have two $\beta$. Then by the AAD $\thin\theta^{\alpha}\thin^{\alpha}\theta\thin$ of $\thin\theta\thin\theta\thin$, we know a $\delta^2$-fan is $\thick^{\epsilon}\delta^{\beta}\thin^{\delta}\beta^{\alpha}\thin^{\alpha}\gamma^{\epsilon}\thin^{\beta}\delta^{\epsilon}\thick,\thick^{\epsilon}\delta^{\beta}\thin^{\epsilon}\gamma^{\alpha}\thin^{\alpha}\gamma^{\epsilon}\thin^{\beta}\delta^{\epsilon}\thick$. Then we get a vertex $\thin^{\alpha}\beta^{\delta}\thin^{\gamma}\epsilon^{\delta}\thick\cdots=\beta\delta\epsilon=\thick^{\epsilon}\delta^{\beta}\thin^{\alpha}\beta^{\delta}\thin^{\gamma}\epsilon^{\delta}\thick$, contradicting no $\alpha\beta\cdots$.

\subsubsection*{Subsubcase. $\alpha=\epsilon>\delta>\beta>\gamma$}

By $\beta+\delta=(1+\tfrac{4}{f})\pi$ and $\beta\le \delta$, we get $\epsilon>\delta>(\frac{1}{2}+\frac{2}{f})\pi$. This implies at most one $b$-edge at a vertex. Then $\epsilon\cdots=\beta\delta\epsilon,\gamma\epsilon^2,
\gamma^k\delta\epsilon$, and $\delta\thin\delta\cdots,\delta\thin\epsilon\cdots,\epsilon\thin\epsilon\cdots$ are not vertices. Moreover, by $\alpha+2\delta>\alpha+\beta+\delta=2\pi$ and $\delta<\epsilon$, we know $\alpha\delta\cdots,\alpha\epsilon\cdots$ are not vertices.

Let $\theta=\beta,\gamma$. By no $\alpha\delta\cdots,\alpha\epsilon\cdots,\delta\thin\delta\cdots,\delta\thin\epsilon\cdots,\epsilon\thin\epsilon\cdots$, we know the AAD of $\thin\theta\thin\theta\thin$ is $\thin\theta^{\alpha}\thin^{\alpha}\theta\thin$. This implies no consecutive $\theta\theta\theta$. Then by no $\alpha\beta^2,\alpha\beta\gamma$, we know $\alpha\beta\cdots=\alpha\beta^k\gamma^l$ is not a vertex.

By $\gamma\epsilon^2$ and Lemma \ref{fbalance}, we know there is $\delta^2$-fan. The rest of the proof is the same as the previous subsubcase $\alpha=\epsilon>\beta>\delta>\gamma$.
\end{proof}

\begin{proposition}\label{bde_a2c}
There is no tiling, such that $\alpha,\beta,\gamma$ have distinct values, and $\beta\delta\epsilon$ is a vertex, and $\alpha+2\gamma=2\pi$.
\end{proposition}

\begin{proof}
The angle sum of $\beta\delta\epsilon$, and the equality $\alpha+2\gamma=2\pi$, and the angle sum for pentagon imply
\[
\alpha=\tfrac{8}{f}\pi,\;
\gamma=(1-\tfrac{4}{f})\pi,\;
\beta+\delta+\epsilon=2\pi.
\]
We have $\alpha<\frac{2}{3}\pi<\gamma$. 

By Lemma \ref{geometry1}, we divide the proof into the case $\beta>\gamma$ and $\delta<\epsilon$, and the case $\beta<\gamma$ and $\delta>\epsilon$.

\subsubsection*{Case. $\beta>\gamma$ and $\delta<\epsilon$}

By $\beta\delta\epsilon$ and $\delta<\epsilon$, we get $\beta\epsilon\cdots=\beta\delta\epsilon$.

By $\alpha+2\gamma=2\pi$ and $\beta>\gamma$, we get $R(\beta^2)<R(\beta\gamma)<\alpha<\gamma<\beta$. By $\beta\delta\epsilon$ and Lemma \ref{square}, we know $\beta^2\cdots$ is a $\hat{b}$-vertex. Then by $\beta\epsilon\cdots=\beta\delta\epsilon$, we get $\beta\gamma\cdots=\beta\gamma\delta^k$, and $\beta^2\cdots$ is not a vertex. By no $\beta^2\cdots$, the AAD implies $\delta\thin\delta\cdots$ is not a vertex. 

If $\beta\gamma\cdots=\beta\gamma\delta^k$ is a vertex, then by $\beta\delta\epsilon$ and $k$ even, we get $\epsilon=\gamma+(k-1)\delta>\gamma$. Then by $\alpha+2\gamma=2\pi$, we get $R(\epsilon^2)<\alpha<\gamma<\beta,\epsilon$. Therefore $\epsilon^2\cdots=\delta^k\epsilon^2$. This implies $\epsilon\cdots=\delta^k\epsilon\cdots (k\ge 1),\delta^k\epsilon^2 (k\ge 2)$, with no $\epsilon$ in the remainder. Then by $\beta\gamma\delta^k$ and applying the counting lemma to $\delta,\epsilon$, we get a contradiction. Therefore $\beta\gamma\cdots$ is not a vertex.

By no $\beta^2\cdots,\beta\gamma\cdots$, we know $\alpha\thin\delta\cdots,\delta\thin\delta\cdots,\delta\thin\epsilon\cdots$ are not vertices, and the AAD of $\thin\alpha\thin\alpha\thin$ is $\thin^{\beta}\alpha^{\gamma}\thin^{\gamma}\alpha^{\beta}\thin$. This implies no consecutive $\alpha\alpha\alpha$. 

By $\beta\epsilon\cdots=\beta\delta\epsilon$ and no $\alpha\thin\delta\cdots,\beta^2\cdots,\beta\gamma\cdots$, we get $\beta\cdots=\beta\delta\epsilon,\alpha^k\beta,\beta\delta^k$. By $\beta\delta\epsilon$ and $\delta<\epsilon$, we get $k\ge 4$ in $\beta\delta^k$, contradicting no $\delta\thin\delta\cdots$. Moreover, by no consecutive $\alpha\alpha\alpha$, we get $k=2$ in $\alpha^k\beta$. Therefore $\beta\cdots=\beta\delta\epsilon,\alpha^2\beta$. This implies $\alpha\beta\cdots=\alpha^2\beta$. 

Suppose $\alpha^2\beta$ is a vertex. Then the angle sum further implies
\[
\alpha=\tfrac{8}{f}\pi,\;
\beta=(2-\tfrac{16}{f})\pi,\;
\gamma=(1-\tfrac{4}{f})\pi,\;
\delta+\epsilon=\tfrac{16}{f}\pi.
\]
By $\delta<\epsilon$, we get $\delta<\alpha=\tfrac{8}{f}\pi<\epsilon$. Then by $R(\gamma^2)=\alpha<\beta,\gamma,\epsilon$, we get $\gamma^2\cdots=\alpha\gamma^2,\gamma^2\delta^k$.

By the AAD $\thin^{\beta}\alpha^{\gamma}\thin^{\gamma}\alpha^{\beta}\thin$ of $\thin\alpha\thin\alpha\thin$, we get $\alpha\beta\cdots=\alpha^2\beta=\thin^{\beta}\alpha^{\gamma}\thin^{\gamma}\alpha^{\beta}\thin\beta\thin=\thin^{\gamma}\alpha^{\beta}\thin\beta\thin^{\beta}\alpha^{\gamma}\thin$. Therefore the AAD of $\thin\alpha\thin\beta\thin$ is $\thin^{\gamma}\alpha^{\beta}\thin\beta\thin$. Since the AAD $\thin^{\epsilon}\gamma^{\alpha}\thin^{\beta}\delta^{\epsilon}\thick$ implies $\thin^{\beta}\alpha^{\gamma}\thin^{\delta}\beta^{\alpha}\thin$, contradicting $\thin^{\gamma}\alpha^{\beta}\thin\beta\thin$, we know the AAD of $\thin\gamma\thin\delta\thick$ is $\thin^{\alpha}\gamma^{\epsilon}\thin^{\beta}\delta^{\epsilon}\thick$. 

By $\beta\cdots=\beta\delta\epsilon,\alpha^2\beta$, we know a $\delta^2$-fan has no $\beta$. Then by no $\alpha\thin\delta\cdots,\delta\thin\delta\cdots$, and the AAD  $\thin^{\alpha}\gamma^{\epsilon}\thin^{\beta}\delta^{\epsilon}\thick$ of $\thin\gamma\thin\delta\thick$, we know a $\delta^2$-fan is $\thick^{\epsilon}\delta^{\beta}\thin^{\epsilon}\gamma^{\alpha}\thin\cdots\thin^{\alpha}\gamma^{\epsilon}\thin^{\beta}\delta^{\epsilon}\thick$, where $\cdots$ consists of $\alpha,\gamma$. Then by $\gamma^2\cdots=\alpha\gamma^2,\gamma^2\delta^k$, and no $\delta\thin\delta\cdots$, we know a $\delta^2$-fan is the vertex $\thick^{\epsilon}\delta^{\beta}\thin^{\epsilon}\gamma^{\alpha}\thin^{\alpha}\gamma^{\epsilon}\thin^{\beta}\delta^{\epsilon}\thick$. This determines $T_1,T_2,T_3,T_4$ in the first of Figure \ref{bde_a2cA}. Then $\beta_1\epsilon_2\cdots=\beta\delta\epsilon$ determines $T_5$. Then $\alpha_1\beta_5\cdots=\alpha^2\beta$ and no $\beta\gamma\cdots$ determine $T_6$. Then $\gamma_1\gamma_6\cdots=\alpha\gamma^2,\gamma^2\delta^k$. By the same argument, we get $\gamma_3\cdots=\alpha\gamma^2,\gamma^2\delta^k$. Then we get $\epsilon_1\epsilon_3\cdots=\beta\epsilon^2\cdots,\gamma^2\epsilon^2\cdots$. By $\beta\delta\epsilon, \gamma^2\delta^2$, and $\delta<\epsilon$, both are not vertices. 

Therefore there is no $\delta^2$-fan. This implies $\gamma^2\cdots=\alpha\gamma^2$. Moreover, by Lemma \ref{fbalance}, we know all fans are $\delta\epsilon$-fans. By $\beta\delta\epsilon$, a $\delta\epsilon$-fan with $\beta$ is the vertex $\beta\delta\epsilon$. By $\gamma^2\cdots=\alpha\gamma^2$, and no $\alpha\thin\delta\cdots,\delta\thin\epsilon\cdots$, a $\delta\epsilon$-fan without $\beta$ is $\thick\delta\thin\gamma\thin\alpha\thin\cdots\thin\alpha\thin\epsilon\thick$. By the AAD $\thin^{\beta}\alpha^{\gamma}\thin^{\gamma}\alpha^{\beta}\thin$ of $\thin\alpha\thin\alpha\thin$, and the AAD $\thin^{\alpha}\gamma^{\epsilon}\thin^{\beta}\delta^{\epsilon}\thick$ of $\thin\gamma\thin\delta\thick$, and no $\beta\gamma\cdots$, we know the fan is $\thick^{\epsilon}\delta^{\beta}\thin^{\epsilon}\gamma^{\alpha}\thin^{\beta}\alpha^{\gamma}\thin^{\gamma}\epsilon^{\delta}\thick$. However, this implies $\thin^{\beta}\alpha^{\gamma}\thin^{\alpha}\beta^{\delta}\thin$, contradicting the AAD $\thin^{\gamma}\alpha^{\beta}\thin\beta\thin$ of $\thin\alpha\thin\beta\thin$. 

We conclude that, if $\alpha^2\beta$ is a vertex, then $\beta\delta\epsilon$ is the only $b$-vertex. If $\alpha^2\beta$ is not a vertex, then by $\beta\cdots=\beta\delta\epsilon,\alpha^2\beta$, we get $\beta\cdots=\beta\delta\epsilon$. In either case, by the counting lemma, this implies $\beta\delta\epsilon,\alpha^k\gamma^l$ are all the vertices. By $\alpha+2\gamma=2\pi$, and $\alpha<\gamma$, and no consecutive $\alpha\alpha\alpha$, we get $\alpha^k\gamma^l=\alpha\gamma^2$. Therefore $\beta\delta\epsilon,\alpha^2\gamma$ are all the vertices, contradicting the counting lemma.

\begin{figure}[htp]
\centering
\begin{tikzpicture}[>=latex,scale=1]


\foreach \a in {-1,1}
{
\begin{scope}[xscale=\a]

\draw
	(1,-0.7) -- (1,0.7) -- (0.5,1.1) -- (0,0.7) -- (0,-0.7) -- (0.5,-1.1) -- (1,-0.7)
	(0.5,1.1) -- (0.5,1.8);

\node at (0.8,0.6) {\small $\alpha$};	
\node at (0.8,0.2) {\small $\beta$};
\node at (0.5,0.85) {\small $\gamma$};
\node at (0.2,0.2) {\small $\delta$}; 
\node at (0.2,0.6) {\small $\epsilon$};

\node at (0.2,-0.6) {\small $\alpha$};
\node at (0.5,-0.85) {\small $\beta$};
\node at (0.2,-0.2) {\small $\gamma$};
\node at (0.8,-0.6) {\small $\delta$};	
\node at (0.8,-0.2) {\small $\epsilon$}; 

\node at (0.7,1.2) {\small $\gamma$};

\end{scope}
}

\draw
	(1,-0.7) -- (1.8,-0.7) -- (1.8,1.8) -- (0.5,1.8) 
	(-1,0) -- (1,0)
	(1,0.7) -- (1.8,0.7);

\draw[line width=1.2]
	(0,0) -- (0,0.7)
	(1,0) -- (1,-0.7) 
	(-1,0) -- (-1,-0.7) 
	(0.5,1.8) -- (1.8,1.8);

\node at (1.6,0.5) {\small $\alpha$};
\node at (1.6,-0.5) {\small $\gamma$};
\node at (1.2,0.5) {\small $\beta$};
\node at (1.2,-0.5) {\small $\epsilon$};
\node at (1.2,0) {\small $\delta$};

\node at (0.7,1.6) {\small $\epsilon$};
\node at (1.1,0.9) {\small $\alpha$};
\node at (1.6,0.9) {\small $\beta$};
\node at (1.6,1.6) {\small $\delta$};

\node[inner sep=0.5,draw,shape=circle] at (-0.5,0.4) {\small $3$};
\node[inner sep=0.5,draw,shape=circle] at (0.5,0.4) {\small $1$};
\node[inner sep=0.5,draw,shape=circle] at (0.5,-0.4) {\small $2$};
\node[inner sep=0.5,draw,shape=circle] at (-0.5,-0.4) {\small $4$};
\node[inner sep=0.5,draw,shape=circle] at (1.55,0) {\small $5$};
\node[inner sep=0.5,draw,shape=circle] at (1.15,1.35) {\small $6$};


\begin{scope}[xshift=4cm]

\foreach \a in {-1,1}
{
\begin{scope}[xscale=\a]

\draw
	(1,-0.7) -- (1,0.7) -- (0.5,1.1) -- (0,0.7) -- (0,-0.7) -- (0.5,-1.1) -- (1,-0.7)
	(0.5,1.1) -- (0.5,1.8);

\node at (0.8,0.6) {\small $\alpha$};	
\node at (0.8,0.2) {\small $\beta$};
\node at (0.5,0.85) {\small $\gamma$};
\node at (0.2,0.2) {\small $\delta$}; 
\node at (0.2,0.6) {\small $\epsilon$};

\node at (0.2,-0.6) {\small $\gamma$};
\node at (0.5,-0.85) {\small $\epsilon$};
\node at (0.2,-0.2) {\small $\alpha$};
\node at (0.8,-0.6) {\small $\delta$};	
\node at (0.8,-0.2) {\small $\beta$}; 

\node at (0.7,1.2) {\small $\gamma$};
\node at (0.3,1.2) {\small $\alpha$};

\end{scope}
}

\draw
	(1,-0.7) -- (1.8,-0.7) -- (1.8,0.7) -- (1,0.7) -- (1.8,1) -- (1.8,1.8) -- (0.5,1.8) 
	(-1,0) -- (1,0);

\draw[line width=1.2]
	(0,0) -- (0,0.7)
	(0.5,-1.1) -- (1,-0.7) 
	(-0.5,-1.1) -- (-1,-0.7)
	(1.8,0.7) -- (1,0.7);

\node at (1.6,0.5) {\small $\epsilon$};
\node at (1.6,-0.5) {\small $\gamma$};
\node at (1.2,0.5) {\small $\delta$};
\node at (1.2,-0.5) {\small $\alpha$};
\node at (1.2,0) {\small $\beta$};

\node at (1.05,0.9) {\small $\theta$};

\node[inner sep=0.5,draw,shape=circle] at (-0.5,0.4) {\small $3$};
\node[inner sep=0.5,draw,shape=circle] at (0.5,0.4) {\small $1$};
\node[inner sep=0.5,draw,shape=circle] at (0.5,-0.4) {\small $2$};
\node[inner sep=0.5,draw,shape=circle] at (-0.5,-0.4) {\small $4$};
\node[inner sep=0.5,draw,shape=circle] at (1.55,0) {\small $5$};
\node[inner sep=0.5,draw,shape=circle] at (1.15,1.35) {\small $6$};

\end{scope}

\end{tikzpicture}
\caption{Proposition \ref{bde_a2c}: $\gamma^2\delta^2$ and $\alpha^2\delta^2$.}
\label{bde_a2cA}
\end{figure}
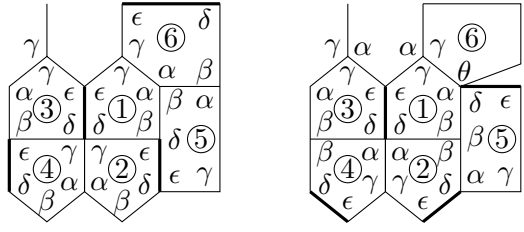

\subsubsection*{Case. $\beta<\gamma$ and $\delta>\epsilon$}

By $\beta\delta\epsilon$, and $\beta<\gamma$, and $\delta>\epsilon$, we get $\beta\delta\cdots=\beta\delta\epsilon$, and $\gamma\delta\cdots$ is not a vertex. Then the AAD of $\thin\alpha\thin^{\delta}\beta^{\alpha}\thin$ is $\thin^{\gamma}\alpha^{\beta}\thin^{\delta}\beta^{\alpha}\thin$. This implies a vertex $\thick^{\epsilon}\delta^{\beta}\thin^{\alpha}\beta^{\delta}\thin\cdots=\beta\delta\epsilon=\thick^{\epsilon}\delta^{\beta}\thin^{\alpha}\beta^{\delta}\thin^{\gamma}\epsilon^{\delta}\thick$, contradicting no $\gamma\delta\cdots$. Therefore the AAD of $\thin\alpha\thin\beta\thin$ is $\thin\alpha\thin^{\alpha}\beta^{\delta}\thin$. 

By $\beta\delta\epsilon$, and $\beta<\gamma$, and Lemma \ref{square}, we know $\beta^2\cdots, \gamma^2\cdots,\beta\gamma\cdots$ are $\hat{b}$-vertices. 

By $\beta\delta\epsilon$ and $\beta<\gamma<\pi$, we get $\delta+\epsilon>\pi$. By $\delta>\epsilon$, this implies $R(\delta^2)$ has no $\delta,\epsilon$. By $\beta\delta\epsilon$ and $\delta>\epsilon$, we get $R(\delta^2)<\beta<\gamma$. This implies $\delta^2\cdots=\alpha^k\delta^2$. In particular, $\delta\thin\delta\cdots$ is not a vertex.

\subsubsection*{Subcase. $\delta^2\cdots=\alpha\delta^2$}

The angle sum of $\alpha\delta^2$ further implies
\[
\alpha=\tfrac{8}{f}\pi,\,
\beta+\epsilon=(1+\tfrac{4}{f})\pi,\,
\gamma=\delta=(1-\tfrac{4}{f})\pi.
\]
By $\delta>\epsilon$, we get $\beta>\tfrac{8}{f}\pi=\alpha$.

By $\pi>\gamma=\delta>\alpha,\beta,\epsilon$, the pentagon is convex. By Lemma \ref{geometry10}, we get $2\alpha+\beta+\gamma>2\pi$. 

The AAD $\thick^{\epsilon}\delta^{\beta}\thin^{\beta}\alpha^{\gamma}\thin^{\beta}\delta^{\gamma}\thick$ of $\alpha\delta^2$ implies a vertex $\thin^{\alpha}\beta^{\delta}\thin^{\alpha}\gamma^{\epsilon}\thin\cdots$. We know $\beta\gamma\cdots$ is a $\hat{b}$-vertex. By $2\alpha+\beta+\gamma>2\pi$, we get $R(\beta\gamma)<2\alpha$. By $\alpha<\beta,\gamma$, this implies $\beta\gamma\cdots$ has degree $3$. Therefore $\beta\gamma\cdots=\alpha\beta\gamma,\beta^2\gamma,\beta\gamma^2$. By $\alpha+2\gamma=2\pi$, the vertices $\alpha\beta\gamma,\beta\gamma^2$ respectively imply $\beta=\gamma$, $\alpha=\beta$, both contradictions. The angle sum of $\beta^2\gamma$ further implies
\begin{equation}\label{bde_a2c_eq1}
\alpha=\tfrac{8}{f}\pi,\,
\beta=\epsilon=(\tfrac{1}{2}+\tfrac{2}{f})\pi,\,
\gamma=\delta=(1-\tfrac{4}{f})\pi.
\end{equation}
The inequality $2\alpha+\beta+\gamma>2\pi$ means $f<28$. By $\delta>\epsilon>\frac{1}{2}\pi$, we know $\delta\thin\epsilon\cdots$ is not a vertex. 

By $f<28$ and the angle values, we know $\alpha\delta^2,\beta\delta\epsilon,\gamma\epsilon^2,\alpha\gamma^2,\beta^2\gamma$ are all the possible degree $3$ vertices. In a $3^5$-tile, we have $\delta\cdots=\alpha\delta^2,\epsilon\cdots=\gamma\epsilon^2$ in the matched case, or $\delta\cdots=\epsilon\cdots=\beta\delta\epsilon$ in the twisted case. We also have $\alpha\cdots=\alpha\delta^2,\alpha\gamma^2$. If $\alpha\cdots=\alpha\delta^2$, then $\beta\cdots=\gamma\cdots=\beta^2\gamma$. No matter whether $\delta\cdots,\epsilon\cdots$ are matched or twisted, this implies $\beta,\gamma$ adjacent in a tile, a contradiction. If $\alpha\cdots=\alpha\gamma^2$, then $\beta\cdots=\beta\delta\epsilon$, and $\delta\cdots=\beta\delta\epsilon$. Therefore $\delta\cdots,\epsilon\cdots$ are twisted, and $\epsilon\cdots=\beta\delta\epsilon$. Combined with $\alpha\cdots=\alpha\gamma^2$, we find $\gamma\cdots$ cannot have degree $3$. 

Therefore there is no $3^5$-tile. By Lemma \ref{special_tile}, we get $f\ge 24$. Combined with $f<28$, we get
\begin{align*}
f=24 &\colon
	\alpha=\tfrac{1}{3}\pi,\,
	\beta=\epsilon=\tfrac{7}{12}\pi,\,
	\gamma=\delta=\tfrac{5}{6}\pi. \\
f=26 &\colon
	\alpha=\tfrac{4}{13}\pi,\,
	\beta=\epsilon=\tfrac{15}{26}\pi,\,
	\gamma=\delta=\tfrac{11}{13}\pi.
\end{align*}
By no $\delta\thin\epsilon\cdots$, we know the vertex $\thin^{\alpha}\beta^{\delta}\thin^{\alpha}\gamma^{\epsilon}\thin\cdots$ implied by the AAD of $\alpha\delta^2$ is $\beta^2\gamma=\thin^{\alpha}\beta^{\delta}\thin^{\alpha}\gamma^{\epsilon}\thin^{\alpha}\beta^{\delta}\thin$. This implies a vertex $\alpha\epsilon\cdots$. However, by the angle values, we see that $R(\alpha\delta\epsilon)=\beta-\alpha$ and $R(\alpha\epsilon^2)=\gamma-\alpha$ are not integral combinations of angle values. Therefore $\alpha\epsilon\cdots$ cannot be a vertex.

\subsubsection*{Subcase. $\delta^2\cdots=\alpha^k\delta^2(k\ge 2)$}

By $\beta\delta\epsilon$ and $\delta>\epsilon$, we get $\beta>k\alpha\ge 2\alpha$. By $\alpha^k\delta^2$ and $\alpha+2\gamma=2\pi$, we get $\gamma>\delta$. Then all angles $\le\gamma<\pi$. By Lemma \ref{geometry10}, we get $2\alpha+\beta+\gamma>2\pi$. By $\alpha+2\gamma=2\pi$, this implies $\alpha+\beta>\gamma$. Then $3\alpha+2\beta>2\alpha+\beta+\gamma>2\pi$.

We know $\beta^2\cdots, \gamma^2\cdots,\beta\gamma\cdots$ are $\hat{b}$-vertices. Then by $2\alpha+\beta+\gamma>2\pi$, and $\alpha+2\gamma=2\pi$, and $2\alpha<\beta<\gamma$, we get $\gamma^2\cdots=\alpha\gamma^2$, and $\alpha\beta^2,\beta\gamma\cdots$ are not vertices. Then by $3\alpha+2\beta>2\pi$, and $2\alpha<\beta$, we get $\beta^2\cdots=\alpha^2\beta^2,\beta^3$. 

The AAD $\thin^{\gamma}\alpha^{\beta}\thin^{\beta}\alpha^{\gamma}\thin$ implies $\thin^{\delta}\beta^{\alpha}\thin^{\alpha}\beta^{\delta}\thin\cdots=\thin\alpha\thin^{\delta}\beta^{\alpha}\thin^{\alpha}\beta^{\delta}\thin\alpha\thin,\thin^{\alpha}\beta^{\delta}\thin^{\alpha}\beta^{\delta}\thin^{\delta}\beta^{\alpha}\thin$. The first contradicts the AAD $\thin\alpha\thin^{\alpha}\beta^{\delta}\thin$ of $\thin\alpha\thin\beta\thin$, and the second contradicts no $\delta\thin\delta\cdots$. Then by no $\beta\gamma\cdots$, we know the AAD of $\thin\alpha\thin\alpha\thin$ is $\thin^{\beta}\alpha^{\gamma}\thin^{\gamma}\alpha^{\beta}\thin$. This implies $k=2$ in $\alpha^k\delta^2$, and the AAD of $\alpha^2\delta^2$ is $\thick^{\epsilon}\delta^{\beta}\thin^{\beta}\alpha^{\gamma}\thin^{\gamma}\alpha^{\beta}\thin^{\beta}\delta^{\epsilon}\thick$. This determines $T_1,T_2,T_3,T_4$ in the second of Figure \ref{bde_a2cA}. Then $\thin^{\alpha}\beta_1^{\delta}\thin^{\alpha}\beta_2^{\delta}\thin\cdots=\thin\alpha\thin^{\alpha}\beta^{\delta}\thin^{\alpha}\beta^{\delta}\thin\alpha\thin,\thin^{\alpha}\beta^{\delta}\thin^{\alpha}\beta^{\delta}\thin\beta\thin$. Since the first contradicts the AAD $\thin\alpha\thin^{\alpha}\beta^{\delta}\thin$ of $\thin\alpha\thin\beta\thin$, and $\delta\thin\delta\cdots$ is not a vertex, we determine $T_5$. Then $\alpha_1\delta_5\cdots=\alpha\delta^2\cdots,\alpha\delta\epsilon\cdots$, and we consider the angle $\theta_6$ next to $\alpha_1$.

By $\beta\delta\epsilon$ and $\delta>\epsilon$, we get $R(\alpha\delta^2)<R(\alpha\delta\epsilon)<\beta<\gamma$. By  $\delta+\epsilon>\pi$ and $\delta>\epsilon$, the remainders have no $\delta$. Moreover, by $\beta\delta\epsilon,\alpha^2\delta^2$, and $\alpha\ne\beta$, both remainders are not empty. Therefore $\theta_6=\alpha,\epsilon$. Then by no $\beta\gamma\cdots$, we get $\gamma_1\cdots=\gamma^2\cdots=\alpha\gamma^2$. By the same argument, we get $\gamma_2\cdots=\alpha\gamma^2$. Then we find $\epsilon_1\epsilon_2\cdots=\beta^2\epsilon^2\cdots,\gamma^2\epsilon^2\cdots,\beta\gamma\epsilon^2\cdots$. By $\gamma+\epsilon>\beta+\epsilon=2\pi-\delta>\pi$, we get a contradiction.  
 
\subsubsection*{Subcase. $\delta^2\cdots$ is not a vertex}

By no $\delta^2\cdots$ and the balance lemma, a $b$-vertex is $\delta\epsilon\cdots$, with no $\delta,\epsilon$ in the remainder. Then by $\beta\delta\epsilon$ and $\beta<\gamma$, we know $\beta\delta\epsilon,\alpha^k\delta\epsilon(k\ge 2)$ are all the $b$-vertices. This implies $\gamma\epsilon\cdots,\delta\thin\delta\cdots,\delta\thin\epsilon\cdots,\epsilon\thin\epsilon\cdots$ are not vertices.

If $\alpha^k\delta\epsilon(k\ge 2)$ is a vertex. Then by $\beta\delta\epsilon$, we get $\alpha<\beta$. We know $\gamma^2\cdots$ is a $\hat{b}$-vertex. Then by $\alpha+2\gamma=2\pi$ and $\alpha<\beta<\gamma$, we get $\gamma^2\cdots=\alpha\gamma^2$. 

The AAD $\thin^{\beta}\alpha^{\gamma}\thin^{\gamma}\alpha^{\beta}\thin$ implies a vertex $\thin^{\epsilon}\gamma^{\alpha}\thin^{\alpha}\gamma^{\epsilon}\thin\cdots=\alpha\gamma^2=\thin^{\alpha}\gamma^{\epsilon}\thin^{\beta}\alpha^{\gamma}\thin^{\epsilon}\gamma^{\alpha}\thin$, contradicting no $\gamma\epsilon\cdots$. The AAD $\thin^{\gamma}\alpha^{\beta}\thin^{\beta}\alpha^{\gamma}\thin$ implies a vertex $\thin^{\delta}\beta^{\alpha}\thin^{\alpha}\beta^{\delta}\thin\cdots$. The AAD $\thin^{\beta}\alpha^{\gamma}\thin^{\beta}\alpha^{\gamma}\thin$ implies a vertex $\thin^{\delta}\beta^{\alpha}\thin^{\alpha}\gamma^{\epsilon}\thin\cdots$. Both are $\hat{b}$-vertices. By no $\delta\thin\delta\cdots,\delta\thin\epsilon\cdots,\epsilon\thin\epsilon\cdots$, the remainders cannot consist of only $\beta,\gamma$, and must have $\alpha$. By $\gamma^2\cdots=\alpha\gamma^2$, the vertices have at most one $\gamma$. Therefore the AADs $\thin^{\delta}\beta^{\alpha}\thin^{\alpha}\beta^{\delta}\thin,\thin^{\delta}\beta^{\alpha}\thin^{\alpha}\gamma^{\epsilon}\thin$ are extended to $\thin\alpha\thin\beta\thin\cdots\thin^{\delta}\beta^{\alpha}\thin^{\alpha}\theta\thin$, where $\theta=\beta,\gamma$, and $\cdots$ consists of $\beta$. Then by no $\delta\thin\delta\cdots$, we get $\thin\alpha\thin^{\delta}\beta^{\alpha}\thin\cdots\thin^{\delta}\beta^{\alpha}\thin^{\alpha}\gamma^{\epsilon}\thin$, contradicting the AAD $\thin\alpha\thin^{\alpha}\beta^{\delta}\thin$ of $\thin\alpha\thin\beta\thin$. We conclude $\alpha\thin\alpha\cdots$ is not a vertex. This implies $\alpha^k\delta\epsilon$ is not a vertex.

Therefore $\beta\delta\epsilon$ is the only $b$-vertex. By the counting lemma, this implies $\beta\delta\epsilon,\alpha^k\gamma^l$ are all the vertices. Then $\beta^2\cdots,\beta\gamma\cdots$ are not vertices. This implies the AAD of $\thin\alpha\thin\alpha\thin$ is $\thin^{\beta}\alpha^{\gamma}\thin^{\gamma}\alpha^{\beta}\thin$, and further implies no consecutive $\alpha\alpha\alpha$. Then by $\alpha+2\gamma=\pi$, and $\alpha<\gamma$, and no consecutive $\alpha\alpha\alpha$, we get $\alpha^k\gamma^l=\alpha\gamma^2$. Therefore $\beta\delta\epsilon,\alpha^2\gamma$ are all the vertices, contradicting the counting lemma.\end{proof}

\begin{proposition}\label{bde_b2c}
There is no tiling, such that $\alpha,\beta,\gamma$ have distinct values, and $\beta\delta\epsilon,\beta\gamma^2$ are vertices.
\end{proposition}

\begin{proof}
By $\beta\delta\epsilon$ and Proposition \ref{bde_2ac}, we know $\alpha^2\gamma$ is not a vertex. Then by Lemma \ref{geometry2} and distinct $\alpha,\beta,\gamma$ values, the only degree $3$ vertices besides $\beta\delta\epsilon,\beta\gamma^2$ are $\alpha^3,\alpha\beta^2,\alpha\delta^2,\alpha\epsilon^2,\gamma\epsilon^2$. 

The angle sums of $\beta\delta\epsilon,\beta\gamma^2$, and the angle sum of one of $\alpha^3,\alpha\beta^2$, and the angle sum for pentagon imply
\begin{align*}
\alpha^3 &\colon 
\alpha=\tfrac{2}{3}\pi,\,
\beta=(\tfrac{4}{3}-\tfrac{8}{f})\pi,\,
\gamma=(\tfrac{1}{3}+\tfrac{4}{f})\pi. \\
\alpha\beta^2 &\colon 
\alpha=(\tfrac{2}{5}+\tfrac{16}{5f})\pi,\,
\beta=(\tfrac{4}{5}-\tfrac{8}{5f})\pi,\,
\gamma=(\tfrac{3}{5}+\tfrac{4}{5f})\pi. 
\end{align*}
We have $\beta>\gamma$. By Lemma \ref{geometry1}, this implies $\delta<\epsilon$. Then by $\beta\delta\epsilon$, we know $\alpha\delta^2$ is not a vertex. Then by $\alpha\ne\beta$, we know $\alpha^3,\alpha\beta^2,\alpha\delta^2$ are mutually exclusive.

Suppose the vertices $\beta\cdots,\delta\cdots$ in a tile have degree $3$. Then $\delta\cdots=\beta\delta\epsilon,\alpha\delta^2$. If $\delta\cdots=\alpha\delta^2$, then by $\alpha^3,\alpha\beta^2,\alpha\delta^2$ mutually exclusive, and $\delta\ne\epsilon$, we know $\alpha^3,\alpha\beta^2,\alpha\epsilon^2$ are not vertices. This implies $\beta\cdots=\beta\gamma^2$, and $\alpha\cdots=\alpha^2\cdots,\alpha\epsilon\cdots$ has high degree. If $\delta\cdots=\beta\delta\epsilon$, then $\beta\cdots=\beta\delta\epsilon,\alpha\beta^2$. If $\beta\cdots=\beta\delta\epsilon$, then $\alpha\cdots=\alpha\gamma\cdots$ has high degree. If $\beta\cdots=\alpha\beta^2$, then by $\alpha^3,\alpha\beta^2,\alpha\delta^2$ mutually exclusive, we know $\alpha^3,\alpha\delta^2$ are not vertices. Then $\delta\cdots=\beta\delta\epsilon$ and $\beta\cdots=\alpha\beta^2$ imply $\alpha\cdots=\alpha^2\cdots,\alpha\delta\cdots$, which has high degree. 

We proved that, in any tile, one of $\alpha\cdots,\beta\cdots,\delta\cdots$ has high degree. By Lemma \ref{special_tile}, we get $f\ge 24$. Moreover, in a special tile, both $\gamma\cdots,\epsilon\cdots$ have degree $3$ in a special tile. Then $\gamma\cdots=\beta\gamma^2,\gamma\epsilon^2$. The two cases are described by the first and second of Figure \ref{bde_b2cA}. 

In the first picture, the vertex $\gamma\cdots=\beta\gamma^2$ implies the degree $3$ vertex $\epsilon\cdots=\alpha\epsilon^2$. This further implies $\alpha\delta^2$ is not a vertex, and $\delta\cdots=\delta^2\cdots$ is the vertex $H$ of degree $4$ or $5$. Therefore $\alpha\cdots,\beta\cdots$ have degree $3$. Then $\alpha\cdots=\alpha^3,\alpha\epsilon^2$. If $\alpha\cdots=\alpha^3$, then $\alpha\beta^2$ is not a vertex, and $\beta\cdots=\beta\gamma^2$. If $\alpha\cdots=\alpha\epsilon^2$, then we still get $\beta\cdots=\beta\gamma^2$.

In the second picture, the vertex $\gamma\cdots=\gamma\epsilon^2$ implies $\epsilon\cdots=\gamma\epsilon^2$, and $\alpha\cdots=\alpha\gamma\cdots$ is the vertex $H$ of degree $4$ or $5$. Therefore $\beta\cdots,\delta\cdots$ have degree $3$. Then $\delta\cdots=\delta^2\cdots=\alpha\delta^2$. This implies $\alpha\beta^2$ is not a vertex, and $\beta\cdots=\beta\gamma^2$. 

\begin{figure}[htp]
\centering
\begin{tikzpicture}[>=latex,scale=1]

\foreach \a in {0,1}
{
\begin{scope}[xshift=2.5*\a cm]


\draw
	(-0.5,-0.4) -- (-0.5,0.7) -- (0,1.1) -- (0.5,0.7) -- (0.5,-0.4)
	(-0.5,0.7) -- ++(-0.4,0);
		
\draw[line width=1.2]
	(-0.5,0) -- (0.5,0);

\node at (0,0.85) {\small $\alpha$}; 
\node at (-0.3,0.6) {\small $\beta$};
\node at (0.3,0.6) {\small $\gamma$};
\node at (-0.3,0.2) {\small $\delta$};
\node at (0.3,0.2) {\small $\epsilon$};	

\node at (-0.3,-0.2) {\small $\delta$};	
\node at (0.3,-0.2) {\small $\epsilon$};
\node at (-0.7,0.5) {\small $\gamma$};
\node at (-0.6,0.9) {\small $\gamma$};

\end{scope}
}

\draw
	(0,1.1) -- ++(0,0.5)
	(0.5,0.7) -- ++(0.6,0);

\node at (0.7,0.5) {\small $\beta$};
\node at (0.6,0.9) {\small $\gamma$};
\node at (0.95,0.5) {\small $\gamma$};
\node at (0.85,0.9) {\small $\beta$};
\node at (-0.2,1.15) {\small $\epsilon$};
\node at (0.2,1.15) {\small $\epsilon$};
\node at (-0.2,1.4) {\small $\alpha$};
\node at (0.2,1.4) {\small $\alpha$};
\node at (0.7,0) {\small $\alpha$};

\begin{scope}[xshift=2.5cm]

\draw[line width=1.2]
	(0.5,0.7) -- ++(0.4,0);

\node at (0,1.3) {\small $H$};
\node at (0.7,0.5) {\small $\epsilon$};
\node at (0.6,0.9) {\small $\epsilon$};
\node at (0.7,0) {\small $\gamma$};
	
\end{scope}


\begin{scope}[xshift=4.5cm]

\draw
	(0,-0.7) -- (0,0.7) -- (0.5,1.1) -- (1,0.7) -- (1,-0.7) -- (0.5,-1.1) -- (0,-0.7)
	(0,0.7) -- (-0.5,1.1) -- (-0.5,1.8) -- (1.8,1.8) -- (1.8,0.7)
	(0,0) -- (1,0)
	(0,-0.7) -- (-0.5,-1.1)
	(1,0.7) -- (2.5,0.7) -- (2.5,-1.5) -- (1,-1.5) -- (0.5,-1.1) -- (0,-1.5)
	(1,-0.7) -- (1.8,-0.7) -- (1.8,-1.5)
	(0,0) -- (-0.5,0);

\draw[line width=1.2]
	(0,0) -- (1,0)
	(0.5,1.1) -- (0.5,1.8)
	(1.8,0.7) -- (1.8,-0.7)
	(1,-1.5) -- (0.5,-1.1);

\node at (0.5,0.85) {\small $\alpha$};
\node at (0.2,0.6) {\small $\beta$};
\node at (0.8,0.6) {\small $\gamma$};	
\node at (0.2,0.2) {\small $\delta$};
\node at (0.8,0.2) {\small $\epsilon$}; 

\node at (-0.3,1.6) {\small $\beta$};
\node at (0,0.95) {\small $\gamma$};
\node at (-0.3,1.2) {\small $\alpha$};
\node at (0.3,1.6) {\small $\delta$};
\node at (0.3,1.2) {\small $\epsilon$};

\node at (0.5,-0.85) {\small $\alpha$};
\node at (0.2,-0.6) {\small $\beta$};
\node at (0.8,-0.6) {\small $\gamma$};	
\node at (0.2,-0.2) {\small $\delta$};
\node at (0.8,-0.2) {\small $\epsilon$}; 

\node at (1.6,-0.5) {\small $\epsilon$};
\node at (1.2,-0.5) {\small $\gamma$};
\node at (1.6,0.5) {\small $\delta$};
\node at (1.2,0) {\small $\alpha$};
\node at (1.2,0.5) {\small $\beta$};

\node at (1.6,0.9) {\small $\alpha$};
\node at (1.6,1.6) {\small $\beta$};
\node at (1,0.9) {\small $\gamma$};
\node at (0.7,1.6) {\small $\delta$};
\node at (0.7,1.2) {\small $\epsilon$};

\node at (1.6,-0.9) {\small $\alpha$};
\node at (1.05,-0.9) {\small $\beta$};
\node at (1.6,-1.3) {\small $\gamma$};
\node at (0.75,-1.07) {\small $\delta$};
\node at (1.05,-1.3) {\small $\epsilon$};

\node at (2.3,-1.3) {\small $\alpha$};
\node at (2.3,0.5) {\small $\beta$};
\node at (2,-1.3) {\small $\gamma$};
\node at (2,0.5) {\small $\delta$};
\node at (2,-0.7) {\small $\epsilon$};

\node at (2,0.9) {\small $\theta$};
\node at (0.5,-1.35) {\small $\delta$};
\node at (0.25,-1.12) {\small $\theta$};
\node at (-0.2,0.6) {\small $\gamma$};
\node at (-0.2,0.2) {\small $\alpha$}; 
\node at (-0.2,-0.2) {\small $\theta$}; 

\node[inner sep=0.5,draw,shape=circle] at (0.5,0.4) {\small $1$};
\node[inner sep=0.5,draw,shape=circle] at (0.5,-0.4) {\small $2$};
\node[inner sep=0.5,draw,shape=circle] at (0,1.4) {\small $3$};
\node[inner sep=0.5,draw,shape=circle] at (1.2,1.3) {\small $4$};
\node[inner sep=0.5,draw,shape=circle] at (1.5,0) {\small $5$};
\node[inner sep=0.5,draw,shape=circle] at (2.15,0) {\small $6$};
\node[inner sep=0.5,draw,shape=circle] at (1.3,-1.1) {\small $7$};

\end{scope}


\begin{scope}[xshift=8cm]

\draw
	(0,-0.7) -- (0,0.7) -- (0.5,1.1) -- (1,0.7) -- (1,-0.7) -- (0.5,-1.1) -- (0,-0.7)
	(0,0.7) -- ++(-0.4,0)
	(0,-0.7) -- ++(-0.4,0)
	(1,-0.7) -- ++(0.4,0)
	(0.5,1.1) -- (0.7,1.5);
	
\draw[line width=1.2]
	(0,0) -- (1,0)
	(1,0.7) -- ++(0.4,0)
	(0.5,-1.1) -- (0.5,-1.5)
	(0.5,1.1) -- (0.3,1.5);
	
\node at (0.5,0.85) {\small $\alpha$};
\node at (0.2,0.6) {\small $\beta$};
\node at (0.8,0.6) {\small $\gamma$};	
\node at (0.2,0.2) {\small $\delta$};
\node at (0.8,0.2) {\small $\epsilon$}; 

\node at (0.5,-0.85) {\small $\alpha$};
\node at (0.2,-0.6) {\small $\beta$};
\node at (0.8,-0.6) {\small $\gamma$};	
\node at (0.2,-0.2) {\small $\delta$};
\node at (0.8,-0.2) {\small $\epsilon$}; 

\node at (-0.2,0) {\small $\alpha$};
\node at (-0.2,0.5) {\small $\gamma$};
\node at (-0.2,-0.5) {\small $\beta$};
\node at (-0.1,0.9) {\small $\gamma$};

\node at (1.2,0) {\small $\gamma$};
\node at (1.2,0.5) {\small $\epsilon$};
\node at (1.1,0.9) {\small $\epsilon$};
\node at (1.2,-0.5) {\small $\alpha$};

\node at (0.7,-1.2) {\small $\delta$};
\node at (0.3,-1.2) {\small $\delta$};

\node at (0.75,1.15) {\small $\gamma$};
\node at (0.3,1.15) {\small $\epsilon$};

\node[inner sep=0.5,draw,shape=circle] at (0.5,0.4) {\small $1$};
\node[inner sep=0.5,draw,shape=circle] at (0.5,-0.4) {\small $2$};

\end{scope}

\end{tikzpicture}
\caption{Proposition \ref{bde_b2c}: Special tile and special companion pair.}
\label{bde_b2cA}
\end{figure}

By $\alpha^3,\alpha\beta^2,\alpha\delta^2$ mutually exclusive, the two special tiles are mutually exclusive.

Suppose a special tile is the first of Figure \ref{bde_b2cA}. Then $\alpha\epsilon^2$ is a vertex. This implies $\alpha\delta^2,\gamma\epsilon^2$ are not vertices, and $\beta\delta\epsilon,\beta\gamma^2,\alpha\epsilon^2,\alpha^3,\alpha\beta^2$ are all the degree $3$ vertices. 

 If $\alpha^3$ is also a vertex, then the angle sums of $\beta\delta\epsilon,\beta\gamma^2,\alpha^3,\alpha\epsilon^2$ and the angle sum for pentagon imply
\[
\alpha=\epsilon=\tfrac{2}{3}\pi,\,
\beta=(\tfrac{4}{3}-\tfrac{8}{f})\pi,\,
\gamma=(\tfrac{1}{3}+\tfrac{4}{f})\pi,\,
\delta=\tfrac{8}{f}\pi.
\]
By $f\ge 24$, we get $\beta\ge\pi>\alpha=\epsilon>\gamma>\delta$. By $\beta\delta\epsilon$ and $\delta<\epsilon$, we get $\beta\epsilon\cdots=\beta\delta\epsilon$.

By $R(\beta^2)<R(\alpha\beta)=\delta<\alpha,\beta,\gamma,\epsilon$, we know $\alpha\beta\cdots,\beta^2\cdots$ are not vertices. This implies $\delta\thin\delta\cdots$ is not a vertex, and the AAD of $\thin\gamma\thin\delta\thick$ is $\thin^{\alpha}\gamma^{\epsilon}\thin^{\beta}\delta^{\epsilon}\thick$. The AAD $\thin^{\alpha}\gamma^{\epsilon}\thin^{\beta}\delta^{\epsilon}\thick$ implies a vertex $\thin^{\alpha}\beta^{\delta}\thin^{\gamma}\epsilon^{\delta}\thick\cdots=\beta\delta\epsilon=\thick^{\epsilon}\delta^{\beta}\thin^{\alpha}\beta^{\delta}\thin^{\gamma}\epsilon^{\delta}\thick$, contradicting no $\alpha\beta\cdots$. Therefore $\gamma\thin\delta\cdots$ is not a vertex.

By $\alpha\epsilon^2$ and Lemma \ref{fbalance}, we know there is $\delta^2$-fan. If the $\delta^2$-fan has $\beta$, then by no $\alpha\beta\cdots,\beta^2\cdots$, the fan consists of one $\beta$ and several $\gamma$. Then by no $\gamma\thin\delta\cdots$, the $\delta^2$-fan is $\thick\delta\thin\beta\thin\delta\thick=\thick^{\epsilon}\delta^{\beta}\thin^{\alpha}\beta^{\delta}\thin^{\beta}\delta^{\epsilon}\thick$, contradicting no $\alpha\beta\cdots$. 

Therefore the $\delta^2$-fan has no $\beta$. By no $\gamma\thin\delta\cdots,\delta\thin\delta\cdots$, the fan has $\alpha$. By $\alpha^3$, the number of $\alpha$ in the fan is one or two. If the number of $\alpha$ is one, then by no $\gamma\thin\delta\cdots$, the $\delta^2$-fan is $\thick\delta\thin\alpha\thin\delta\thick=\thick^{\epsilon}\delta^{\beta}\thin^{\beta}\alpha^{\gamma}\thin^{\beta}\delta^{\epsilon}\thick$, contradicting no $\beta^2\cdots$. If the number of $\alpha$ is two, then by $\alpha+\gamma>\pi$, the fan has at most one $\gamma$. Then by no $\gamma\thin\delta\cdots$, the fan is $\thick\delta\thin\alpha\thin\alpha\thin\delta\thick,\thick\delta\thin\alpha\thin\gamma\thin\alpha\thin\delta\thick$. The AAD of the first implies $\beta^2\cdots$, and the AAD of the second implies $\alpha\beta\cdots,\beta^2\cdots$. Both are contradictions.

Therefore $\alpha^3$ is not a vertex. Then the vertex $\alpha\cdots=\alpha\epsilon^2$ in the first of Figure \ref{bde_b2cA}, and $\beta\delta\epsilon,\beta\gamma^2,\alpha\epsilon^2,\alpha\beta^2$ are all the degree $3$ vertices. The special tile is one tile in a special companion pair. The tiles $T_1,T_2$ in the third of Figure \ref{bde_b2cA} is a special companion pair, with $T_1$ being the special tile. Then $T_3,T_4,T_5$ are determined. Then one of $\alpha_2\cdots,\gamma_2\gamma_5\cdots$ has high degree. By Lemma \ref{bb_pair}, this implies $\deg H=4$, and $\beta_2\cdots$ has degree $3$. Therefore $H=\alpha^2\delta^2,\alpha\beta\delta^2,\alpha\gamma\delta^2,\delta^3\epsilon,\delta^2\epsilon^2$. The angle sums of $\beta\gamma^2,\alpha\epsilon^2,\alpha\beta\delta^2$ imply $\alpha+\beta+\gamma+\delta+\epsilon=3\pi$, contradicting the angle sum for pentagon. Therefore $H=\alpha^2\delta^2,\alpha\gamma\delta^2,\delta^3\epsilon,\delta^2\epsilon^2$. The angle sums of $\beta\delta\epsilon,\beta\gamma^2,\alpha\epsilon^2$, and the angle sum of $H$, and the angle sum for pentagon imply
\begin{align*}
\alpha^2\delta^2 & \colon
	\alpha=\tfrac{16}{f}\pi,\,
	\beta=\tfrac{24}{f}\pi,\,
	\gamma=(1-\tfrac{12}{f})\pi,\,
	\delta=(1-\tfrac{16}{f})\pi,\,
	\epsilon=(1-\tfrac{8}{f})\pi. \\
\alpha\gamma\delta^2 & \colon
	\alpha=(\tfrac{1}{3}+\tfrac{20}{3f})\pi,\,
	\beta=(\tfrac{2}{3}+\tfrac{16}{3f})\pi,\,
	\gamma=(\tfrac{2}{3}-\tfrac{8}{3f})\pi,\,
	\delta=(\tfrac{1}{2}-\tfrac{2}{f})\pi,\,
	\epsilon=(\tfrac{5}{6}-\tfrac{10}{3f})\pi. \\
\delta^3\epsilon & \colon
	\alpha=(\tfrac{2}{5}+\tfrac{24}{5f})\pi,\,
	\beta=(\tfrac{4}{5}+\tfrac{8}{5f})\pi,\,
	\gamma=(\tfrac{3}{5}-\tfrac{4}{5f})\pi,\,
	\delta=(\tfrac{2}{5}+\tfrac{4}{5f})\pi,\,
	\epsilon=(\tfrac{4}{5}-\tfrac{12}{5f})\pi. \\
\delta^2\epsilon^2 & \colon
	\alpha=(\tfrac{1}{2}+\tfrac{4}{f})\pi,\,
	\beta=\pi,\,
	\gamma=\tfrac{1}{2}\pi,\,
	\delta=(\tfrac{1}{4}+\tfrac{2}{f})\pi,\,
	\epsilon=(\tfrac{3}{4}-\tfrac{2}{f})\pi. 
\end{align*}

The AAD of $H=\delta^3\epsilon$ implies a vertex $\beta^2\cdots$. By $R(\beta^2)=(\tfrac{2}{5}-\tfrac{16}{5f})\pi<\alpha,\beta,\gamma,\delta,\epsilon$, we get a contradiction.

The AAD of $H=\delta^2\epsilon^2$ implies $\alpha_3\cdots=\alpha^2\cdots$ is a vertex. By $R(\alpha^2)=(1-\tfrac{8}{f})\pi<2\gamma=\delta+\epsilon<\gamma+2\delta<4\delta$, and $\gamma<\alpha<\beta=2\gamma$, and $\delta<\epsilon$, and $\alpha^2\cdots$ having high degree, we get $\alpha_3\cdots=\alpha^2\delta^2$. The AAD of $\alpha^2\delta^2$ implies $\beta^2\cdots$ is a vertex, contradicting $\beta=\pi$.

The AAD of $H=\alpha^2\delta^2$ implies $\beta^2\cdots$ is a vertex. By $\beta\delta\epsilon$ and Lemma \ref{square}, we know $\beta^2\cdots$ has no $\delta,\epsilon$. By $\delta<\epsilon$ and Lemma \ref{geometry1}, we get $\beta>\gamma$. Then by $\beta\gamma^2$, we know $R(\beta^2)$ has no $\beta,\gamma$. Then by $\alpha+\beta>\alpha+\gamma>\pi$, we get $\beta^2\cdots=\alpha\beta^2$. The angle sum of $\alpha\beta^2$ further implies
\[
\alpha=\delta=\tfrac{1}{2}\pi,\,
\beta=\epsilon=\tfrac{3}{4}\pi,\,
\gamma=\tfrac{5}{8}\pi,\,
f=32. 
\]
The angle values imply $\alpha\delta\cdots=\alpha^2\delta^2$, and $\epsilon^2\cdots=\alpha\epsilon^2$, and $\gamma^2\cdots=\beta\gamma^2$. In the third of Figure \ref{bde_b2cA} (all three $\theta$ are $\alpha$), $\alpha_4\delta\cdots=\alpha^2\delta^2$ determines $T_6$. Then $\gamma_2\gamma_5\cdots=\beta\gamma^2$ and $\epsilon_5\epsilon_6\cdots=\alpha\epsilon^2$ determine $T_7$. Then $\alpha_2\delta_7\cdots=\alpha^2\delta^2$. With $\theta=\alpha$, the degree $3$ vertex ${}^{\theta}\thin^{\alpha}\beta_2^{\delta}\thin^{\theta}\cdots=\beta\gamma^2=\thin^{\alpha}\beta^{\delta}\thin^{\alpha}\gamma^{\epsilon}\thin^{\epsilon}\gamma^{\alpha}\thin$. This implies a vertex $\epsilon\thin\epsilon\cdots$. However, the angle values imply $\epsilon\thin\epsilon\cdots$ is not a vertex.

The AAD of $H=\alpha\gamma\delta^2$ implies a vertex $\alpha\beta\cdots$. We have $R(\alpha\beta)=(1-\tfrac{12}{f})\pi<3\alpha,\alpha+\gamma,2\gamma,2\delta$. Then by $\beta>\alpha,\gamma$, and $\delta<\epsilon$, and degree $3$ vertex $\alpha\beta\cdots=\alpha\beta^2$, we get $\alpha\beta\cdots=\alpha\beta^2,\alpha^3\beta$. The angle sum of $\alpha\beta^2$ further implies $f=28$ and $\alpha=\gamma=\frac{4}{7}\pi$, a contradiction. The angle sum of $\alpha^3\beta$ further implies
\[
\alpha=\tfrac{8}{19}\pi,\,
\beta=\tfrac{14}{19}\pi,\,
\gamma=\tfrac{12}{19}\pi,\,
\delta=\tfrac{9}{19}\pi,\,
\epsilon=\tfrac{15}{19}\pi,\,
f=76.
\]
The angle values imply $\alpha\delta\cdots=\alpha\gamma\delta^2$, and $\epsilon^2\cdots=\alpha\epsilon^2$, and $\gamma^2\cdots=\beta\gamma^2$. Then we may repeat the argument for the case $H=\alpha^2\delta^2$ and $f=32$, using the third of Figure \ref{bde_b2cA}. The only difference is all three $\theta$ are $\gamma$. This implies $\beta_2\cdots\ne\beta\delta\epsilon,\beta\gamma^2$. Therefore $\beta_2\cdots$ has high degree, a contradiction. 

This completes the proof that there is no tiling with the first of Figure \ref{bde_b2cA} as a special tile. Then the special tile is the second picture, in which both $\alpha\delta^2,\gamma\epsilon^2$ are vertices. The angle sums of $\beta\delta\epsilon,\beta\gamma^2,\alpha\delta^2,\gamma\epsilon^2$ and the angle sum for pentagon imply
\[
\alpha=(\tfrac{1}{4}+\tfrac{5}{f})\pi,\,
\beta=(\tfrac{1}{2}+\tfrac{2}{f})\pi,\,
\gamma=(\tfrac{3}{4}-\tfrac{1}{f})\pi,\,
\delta=(\tfrac{7}{8}-\tfrac{5}{2f})\pi,\,
\epsilon=(\tfrac{5}{8}+\tfrac{1}{2f})\pi.
\]

The special tile is one tile in a special companion pair. This is $T_1$ in the fourth of Figure \ref{bde_b2cA}. Then $\beta_2\cdots=\beta^2\cdots$ and $\gamma_2\cdots=\alpha\gamma\cdots$ have high degree. Since we also know $\alpha_1\cdots$ has high degree, by Lemma \ref{bb_pair}, we know $\alpha_2\cdots$ has degree $3$. This implies $\alpha_2\cdots=\alpha\delta^2$, and further implies $\beta_2\cdots=\beta^3\cdots$. By $R(\beta^3)=(\tfrac{1}{2}-\tfrac{6}{f})\pi<2\alpha,\beta,\gamma,\delta,\epsilon$, we get $\beta_2\cdots=\alpha\beta^3$. The angle sum of $\alpha\beta^3$ further implies
\[
\alpha=\tfrac{4}{11}\pi,\,
\beta=\tfrac{6}{11}\pi,\,
\gamma=\tfrac{8}{11}\pi,\,
\delta=\tfrac{9}{11}\pi,\,
\epsilon=\tfrac{7}{11}\pi,\,
f=44.
\]

By $\delta>\epsilon>\frac{1}{2}\pi$, we know $\epsilon\thin\epsilon\cdots$ is not a vertex. Then the AAD of $\beta_1\cdots=\beta\gamma^2$ is $\thin^{\alpha}\beta^{\delta}\thin^{\alpha}\gamma^{\epsilon}\thin^{\alpha}\gamma^{\epsilon}\thin$. This implies $\alpha_1\cdots=\alpha\gamma\epsilon\cdots$. By $\alpha+\gamma+2\epsilon=\frac{24}{11}\pi>2\pi$ and $\delta>\epsilon$, we get a contradiction. 
\end{proof}

\begin{proposition}\label{bde_3b}
There is no tiling, such that $\alpha,\beta,\gamma$ have distinct values, and $\beta\delta\epsilon,\beta^3$ are vertices.
\end{proposition}

\begin{proof}
By $\beta\delta\epsilon$ and Proposition \ref{bde_abc}, \ref{bde_2ac}, \ref{bde_a2c}, we know $\alpha\beta\gamma,\alpha^2\gamma,\alpha\gamma^2$ are not vertices. Then by $\beta^3$ and distinct $\alpha,\beta,\gamma$ values, we know $\beta^3$ is the only degree $3$ $\hat{b}$-vertex. 

If $\delta^2\cdots$ is not degree $3$, then all degree $3$ vertices are $\beta\delta\epsilon,\beta^3,\theta\epsilon^2$, with $\theta=\alpha,\gamma$. Applying Lemma \ref{degree3} to $\beta,\epsilon$, we get a contradiction. Therefore there is a degree $3$ vertex $\delta^2\cdots$. Then by $\beta\delta\epsilon$ and Lemma \ref{geometry2}, we get a contradiction.
\end{proof}

\begin{proposition}\label{bde_a2d_c3e}
There is no tiling, such that $\alpha,\beta,\gamma$ have distinct values, and $\beta\delta\epsilon,\alpha\delta^2,\gamma\epsilon^2$ are vertices.
\end{proposition}

\begin{proof}
The angle sums of $\beta\delta\epsilon,\alpha\delta^2,\gamma\epsilon^2$ and the angle sum for pentagon imply
\[
\alpha+\gamma=(1+\tfrac{4}{f})\pi,\,
\beta=(\tfrac{1}{2}+\tfrac{2}{f})\pi,\,
\delta+\epsilon=(\tfrac{3}{2}-\tfrac{2}{f})\pi.
\]
By $\alpha+\gamma=2\beta$, we get $\alpha>\beta>\gamma$ or $\alpha<\beta<\gamma$. By Lemma \ref{geometry1}, we have $\delta<\epsilon$ in the first case, and $\delta>\epsilon$ in the second case.

\subsubsection*{Subcase. $\alpha>\beta>\gamma$ and $\delta<\epsilon$.}

By $\alpha\delta^2$ and $\delta<\epsilon$, we know $\alpha\epsilon\cdots$ is not a vertex, and $\alpha\delta\cdots=\alpha\delta^2$. 

By $\beta\delta\epsilon$, and $\alpha>\beta$, and Lemma \ref{square}, we know $\alpha^2\cdots,\beta^2\cdots$ are $\hat{b}$-vertices. Then by $\alpha+\gamma=2\beta>\pi$ and $\alpha>\beta>\gamma$, we get $\alpha^2\cdots=\alpha^3,\alpha^2\beta,\alpha^2\gamma$, and $\beta^2\cdots=\alpha\beta^2,\beta^2\gamma^k,\beta^3\gamma^k$. By $\alpha+\beta+\gamma=3\beta<2\pi$, we get $k\ge 2$ in $\beta^2\gamma^k$, and $k\ge 1$ in $\beta^3\gamma^k$.

By Proposition \ref{bde_2ac}, we know $\alpha^2\gamma$ is not a vertex. In fact, the angle sum of $\alpha\beta^2$ also implies $2\alpha+\gamma=2\pi$. Therefore $\alpha\beta^2$ is also not a vertex. Then we get $\alpha^2\cdots=\alpha^3,\alpha^2\beta$, and $\beta^2\cdots=\beta^2\gamma^k,\beta^3\gamma^k$.

The AAD $\thick^{\epsilon}\delta^{\beta}\thin^{\beta}\alpha^{\gamma}\thin^{\beta}\delta^{\epsilon}\thick$ of $\alpha\delta^2$ implies $\thin^{\alpha}\beta^{\delta}\thin^{\alpha}\beta^{\delta}\thin\cdots=\beta^2\gamma^k,\beta^3\gamma^k$ is a vertex. By no $\alpha\epsilon\cdots$ and $k\ge 1$, the AAD of the vertex implies $\alpha^2\cdots=\alpha^3,\alpha^2\beta$ is a vertex, and one of $\delta\thin\delta\cdots,\delta\thin\epsilon\cdots,\epsilon\thin\epsilon\cdots$ is a vertex. The angle sum of one of $\alpha^3,\alpha^2\beta$ further implies
\begin{align*}
\alpha^3 &\colon 
	\alpha=\delta=\tfrac{2}{3}\pi,\,
	\beta=(\tfrac{1}{2}+\tfrac{2}{f})\pi,\,
	\gamma=(\tfrac{1}{3}+\tfrac{4}{f})\pi,\,
	\epsilon=(\tfrac{5}{6}-\tfrac{2}{f})\pi. \\
\alpha^2\beta &\colon
	\alpha=(\tfrac{3}{4}-\tfrac{1}{f})\pi,\,
	\beta=(\tfrac{1}{2}+\tfrac{2}{f})\pi,\,
	\gamma=(\tfrac{1}{4}+\tfrac{5}{f})\pi,\,
	\delta=(\tfrac{5}{8}+\tfrac{1}{2f})\pi,\,
	\epsilon=(\tfrac{7}{8}-\tfrac{5}{2f})\pi. 
\end{align*}
We have $\epsilon>\delta>\frac{1}{2}\pi$. This implies $\delta\thin\delta\cdots,\delta\thin\epsilon\cdots,\epsilon\thin\epsilon\cdots$ are not vertices. We get a contradiction.

\subsubsection*{Subcase. $\alpha<\beta<\gamma$ and $\delta>\epsilon$.}

By $\beta\delta\epsilon$, and $\beta+\gamma>\alpha+\gamma=2\beta$, and Lemma \ref{square}, we know $\alpha\gamma\cdots,\beta\gamma\cdots$ are $\hat{b}$-vertices. Then by $\alpha+\gamma=2\beta>\pi$ and $\alpha<\beta<\gamma$, we get $\alpha\gamma\cdots=\alpha\gamma^2,\alpha^k\beta\gamma,\alpha^k\gamma$, and $\beta\gamma\cdots=\beta^2\gamma,\beta\gamma^2,\alpha^k\beta\gamma$. By $\alpha+\beta+\gamma=(\tfrac{3}{2}+\tfrac{6}{f})\pi<2\pi$, we get $k\ge 2$ in $\alpha^k\beta\gamma$.

By Propositions \ref{bde_a2c} and \ref{bde_b2c}, we know $\alpha\gamma^2,\beta\gamma^2$ are not vertices. In fact, the angle sum of $\beta^2\gamma$ also implies $\alpha+2\gamma=2\pi$. Therefore $\beta^2\gamma$ is also not a vertex. Then we get $\alpha\gamma\cdots=\alpha^k\beta\gamma,\alpha^k\gamma$, and $\beta\gamma\cdots=\alpha^k\beta\gamma$. 

If $\gamma^2\cdots$ is a vertex, then $\gamma<\pi$, and the pentagon is convex. By Lemma \ref{geometry10}, we get $2\alpha+\beta+\gamma>2\pi$. This implies $\beta\gamma\cdots=\alpha^k\beta\gamma$ is not a vertex, contradicting the AAD $\thick^{\epsilon}\delta^{\beta}\thin^{\beta}\alpha^{\gamma}\thin^{\beta}\delta^{\epsilon}\thick$ of $\alpha\delta^2$.

Therefore $\gamma^2\cdots$ is not a vertex. By $\gamma\epsilon^2$ and $\delta>\epsilon$, a $b$-vertex $\gamma\cdots=\gamma\epsilon^2$. By $\alpha\gamma\cdots=\alpha^k\beta\gamma,\alpha^k\gamma$, and $\beta\gamma\cdots=\alpha^k\beta\gamma$, and no $\gamma^2\cdots$, a $\hat{b}$-vertex $\gamma\cdots=\alpha^k\beta\gamma,\alpha^k\gamma$. Then we get $\gamma\cdots=\gamma\epsilon^2,\alpha^k\beta\gamma,\alpha^k\gamma$. This implies
\begin{align*}
f=\#\gamma
&=\#\gamma\epsilon^2+\#\alpha^k\beta\gamma+\#\alpha^k\gamma \\
&<\#\gamma\epsilon^2+\#\alpha^k\beta\gamma+\#\alpha^k\gamma+\tfrac{1}{2}\#\alpha\delta^2 \\
&\le \tfrac{1}{2}\#\epsilon+\tfrac{1}{2}\#\alpha=f.
\end{align*}
We get a contradiction.
\end{proof}

\begin{proposition}\label{bde_2ab}
There is no tiling, such that $\alpha,\beta,\gamma$ have distinct values, and $\beta\delta\epsilon,\alpha^2\beta$ are vertices.
\end{proposition}

\begin{proof}
The angle sums of $\beta\delta\epsilon,\alpha^2\beta$ and the angle sum for pentagon imply
\[
\alpha+\gamma=(1+\tfrac{4}{f})\pi,\,
\delta+\epsilon=2\alpha.
\]
We have $2\gamma+\delta+\epsilon=2(\alpha+\gamma)>2\pi$. By $\beta\delta\epsilon$, this implies $\beta<2\gamma$. 

By $\beta\delta\epsilon$ and Proposition \ref{bde_a2c}, we know $\alpha\gamma^2$ is not a vertex. Then by $\alpha^2\beta$, and Lemma \ref{geometry2}, and distinct $\alpha,\beta,\gamma$ values, the only degree $3$ vertices besides $\beta\delta\epsilon,\alpha^2\beta$ are $\beta^2\gamma,\gamma^3,\alpha\delta^2,\alpha\epsilon^2,\gamma\epsilon^2$. By Lemma \ref{degree3}, there is a degree $3$ vertex without $\beta$. Therefore one of $\gamma^3,\alpha\delta^2,\alpha\epsilon^2,\gamma\epsilon^2$ is a vertex.

\subsubsection*{Case. $\gamma^3$ is a vertex}

The angle sums of $\beta\delta\epsilon,\alpha^2\beta,\gamma^3$ and the angle sum for pentagon imply
\[
\alpha=(\tfrac{1}{3}+\tfrac{4}{f})\pi,\;
\beta=(\tfrac{4}{3}-\tfrac{8}{f})\pi,\;
\gamma=\tfrac{2}{3}\pi,\;
\delta+\epsilon=(\tfrac{2}{3}+\tfrac{8}{f})\pi.
\]
We have $\beta>\gamma>\alpha$. This implies $\delta<\epsilon$. Then by $\beta\delta\epsilon$, we get $\beta\epsilon\cdots=\beta\delta\epsilon$. Moreover, by $\delta+\epsilon=2\alpha$, we get $\delta<\alpha<\epsilon$. 

By $\alpha^2\beta$ and $\beta>\alpha>\delta$, we know $\alpha\delta^2$ is not a vertex. By $\gamma^3$ and $\beta\ne\gamma$, we know $\beta^2\gamma$ is not a vertex. Therefore $\beta\delta\epsilon,\alpha^2\beta,\gamma^3,\alpha\epsilon^2,\gamma\epsilon^2$ are all the degree $3$ vertices.

Consider a special tile in the first and second of Figure \ref{bde_2abA}. In the first picture, we assume $\alpha\cdots,\beta\cdots,\delta\cdots$ have degree $3$. Then $\delta\cdots=\beta\delta\epsilon$. This implies that, if $\beta\cdots=\beta\delta\epsilon$, then $\alpha\cdots=\alpha\gamma\cdots$ has high degree. Therefore $\beta\cdots=\alpha^2\beta$. This implies $\alpha\cdots=\alpha^2\beta$ as indicated, and we get $\beta\cdots=\thin^{\beta}\alpha^{\gamma}\thin^{\gamma}\alpha^{\beta}\thin^{\alpha}\beta^{\delta}\thin$. This further implies $\thin^{\epsilon}\gamma^{\alpha}\thin^{\alpha}\gamma^{\epsilon}\thin\cdots$ is a vertex.

In the second picture, we assume both $\gamma\cdots,\epsilon\cdots$ have degree $3$. If $\gamma\cdots=\gamma\epsilon^2$, then $\alpha\cdots=\alpha\gamma\cdots$ has high degree, and the degree $3$ vertex $\epsilon\cdots=\gamma\epsilon\cdots=\gamma\epsilon^2$. This implies $\delta\cdots=\delta^2\cdots$ has high degree, a contradiction. Therefore $\gamma\cdots=\gamma^3$. This implies $\epsilon\cdots=\alpha\epsilon^2$. Then $\delta\cdots=\delta^2\cdots$ has high degree, and $\alpha\cdots,\beta\cdots$ have degree $3$. Since $\alpha\cdots=\alpha\epsilon^2$ implies $\beta\cdots=\beta\gamma\cdots$ has high degree, we get $\alpha\cdots=\alpha^2\beta$ as indicated. Then  we get $\gamma\cdots=\thin^{\alpha}\gamma^{\epsilon}\thin^{\epsilon}\gamma^{\alpha}\thin^{\alpha}\gamma^{\epsilon}\thin$. This implies $\epsilon\thin\epsilon\cdots$ is a vertex.

\begin{figure}[htp]
\centering
\begin{tikzpicture}[>=latex,scale=1]

\foreach \a in {0,1}
{
\begin{scope}[xshift=2*\a cm]

\draw
	(-0.5,-0.4) -- (-0.5,0.7) -- (0,1.1) -- (0.5,0.7) -- (0.5,-0.4);
		
\draw[line width=1.2]
	(-0.5,0) -- (0.5,0);

\node at (0,0.85) {\small $\alpha$}; 
\node at (-0.3,0.6) {\small $\beta$};
\node at (0.3,0.6) {\small $\gamma$};
\node at (-0.3,0.2) {\small $\delta$};
\node at (0.3,0.2) {\small $\epsilon$};	

\end{scope}
}


\draw
	(0,1.1) -- ++(0,0.4)
	(-0.5,0.7) -- ++(-0.4,0);

\node at (-0.7,0.5) {\small $\alpha$};
\node at (-0.6,0.9) {\small $\alpha$};
\node at (-0.7,0) {\small $\beta$};
\node at (0.3,-0.2) {\small $\delta$};
\node at (-0.3,-0.2) {\small $\epsilon$};	
\node at (-0.2,1.2) {\small $\beta$};
\node at (0.2,1.2) {\small $\alpha$};


\begin{scope}[xshift=2 cm]

\draw
	(0,1.1) -- ++(0,0.4)
	(0.5,0.7) -- ++(0.4,0);

\node at (0.7,0.5) {\small $\gamma$};
\node at (0.6,0.9) {\small $\gamma$};
\node at (0.7,0) {\small $\alpha$};
\node at (-0.3,-0.2) {\small $\delta$};
\node at (0.3,-0.2) {\small $\epsilon$};	
\node at (-0.2,1.2) {\small $\beta$};
\node at (0.2,1.2) {\small $\alpha$};

\end{scope}


\begin{scope}[shift={(5.5 cm,0.5cm)}]

\foreach \a in {-1,1}
{
\begin{scope}[xscale=\a]
\draw
	(0,0) -- (0,0.7) -- (0.5,1.1) -- (1,0.7) -- (1.5,1.1) -- (2,0.7) -- (2,0)
	(1,0.7) -- (1,0);

\draw[line width=1.2]
	(1,0) -- (2,0);
	
\node at (0.8,0.6) {\small $\alpha$};	
\node at (0.8,0.2) {\small $\beta$}; 
\node at (0.5,0.85) {\small $\gamma$};
\node at (0.2,0.2) {\small $\delta$};
\node at (0.2,0.6) {\small $\epsilon$};

\node at (1.8,0.6) {\small $\gamma$};	
\node at (1.8,0.2) {\small $\epsilon$}; 
\node at (1.5,0.85) {\small $\alpha$};
\node at (1.2,0.2) {\small $\delta$};
\node at (1.2,0.6) {\small $\beta$};

\node at (1,0.9) {\small $\alpha$};

\end{scope}
}

\node at (0.8,-0.6) {\small $\epsilon$};	
\node at (0.8,-0.2) {\small $\gamma$}; 
\node at (0.5,-0.85) {\small $\delta$};
\node at (0.2,-0.2) {\small $\alpha$};
\node at (0.2,-0.6) {\small $\beta$};

\node at (-0.9,-0.55) {\small $\alpha$}; 
\node at (-1.8,-0.55) {\small $\beta$}; 
\node at (-0.45,-0.2) {\small $\gamma$}; 
\node at (-1.8,-0.2) {\small $\delta$}; 
\node at (-1,-0.2) {\small $\epsilon$}; 

\node at (-0.15,-0.35) {\small $\gamma$}; 
\node at (1.2,-0.2) {\small $\delta$}; 

\node at (0,0.9) {\small $\alpha$};
	
\draw
	(0,0) -- (-0.8,-0.7) -- (-2,-0.7) -- (-2,0) -- (2,0)
	(0,0) -- (0,-0.7) -- (0.5,-1.1) -- (1,-0.7) -- (1,0);
	
\draw[line width=1.2]
	(0,0) -- (0,0.7);
	
\node[inner sep=0.5,draw,shape=circle] at (-0.5,0.4) {\small 1};
\node[inner sep=0.5,draw,shape=circle] at (0.5,0.4) {\small 2};
\node[inner sep=0.5,draw,shape=circle] at (0.5,-0.4) {\small 3};
\node[inner sep=0.5,draw,shape=circle] at (-1.4,-0.35) {\small 4};
\node[inner sep=0.5,draw,shape=circle] at (-1.5,0.4) {\small 5};
\node[inner sep=0.5,draw,shape=circle] at (1.5,0.4) {\small 6};

\end{scope}


\begin{scope}[shift={(9.2 cm,0.5cm)}]

\foreach \a in {-1,1}
{
\begin{scope}[xscale=\a]
\draw
	(0,-0.7) -- (0,0.7) -- (0.5,1.1) -- (1,0.7) -- (1,-0.7) -- (0.5,-1.1) -- (0,-0.7) 
	(0,0) -- (1,0);

\draw[line width=1.2]
	(1,0) -- (1,-0.7);
	
\node at (0.8,0.6) {\small $\alpha$};	
\node at (0.8,0.2) {\small $\beta$}; 
\node at (0.5,0.85) {\small $\gamma$};
\node at (0.2,0.2) {\small $\delta$};
\node at (0.2,0.6) {\small $\epsilon$};

\node at (0.2,-0.6) {\small $\alpha$};
\node at (0.5,-0.85) {\small $\beta$};
\node at (0.2,-0.2) {\small $\gamma$};
\node at (0.8,-0.6) {\small $\delta$};	
\node at (0.8,-0.2) {\small $\epsilon$}; 

\end{scope}
}

\draw
	(1,0.7) -- (1.8,0.7) -- (1.8,-0.7) -- (1,-0.7) -- (1.5,-1.1) -- (1.5,-1.8) -- (-0.5,-1.8) -- (-0.5,-1.1);

\draw[line width=1.2]
	(0,0) -- (0,0.7)
	(0.5,-1.8) -- (0.5,-1.1);

\node at (1.3,-1.2) {\small $\alpha$};	
\node at (1.3,-1.6) {\small $\beta$}; 
\node at (1,-0.95) {\small $\gamma$};
\node at (0.7,-1.6) {\small $\delta$};
\node at (0.7,-1.2) {\small $\epsilon$};

\node at (0.3,-1.2) {\small $\delta$};	
\node at (0.3,-1.6) {\small $\epsilon$}; 
\node at (0,-0.95) {\small $\beta$};
\node at (-0.3,-1.6) {\small $\gamma$};
\node at (-0.3,-1.2) {\small $\alpha$};
	
\node at (1.6,0.5) {\small $\alpha$};
\node at (1.2,0.5) {\small $\beta$};
\node at (1.6,-0.5) {\small $\gamma$};
\node at (1.2,0) {\small $\delta$};
\node at (1.2,-0.5) {\small $\epsilon$};

\node[inner sep=0.5,draw,shape=circle] at (0.5,0.4) {\small 1};
\node[inner sep=0.5,draw,shape=circle] at (0.5,-0.4) {\small 2};
\node[inner sep=0.5,draw,shape=circle] at (-0.5,-0.4) {\small 3};
\node[inner sep=0.5,draw,shape=circle] at (-0.5,0.4) {\small 4};
\node[inner sep=0.5,draw,shape=circle] at (1.5,0) {\small 7};\node[inner sep=0.5,draw,shape=circle] at (0,-1.4) {\small 5};\node[inner sep=0.5,draw,shape=circle] at (1,-1.4) {\small 6};

\end{scope}

\end{tikzpicture}
\caption{Proposition \ref{bde_2ab}: $\gamma^3$ is a vertex.}
\label{bde_2abA}
\end{figure}

The argument above implies there is no $3^5$-tile. By Lemma \ref{special_tile}, we get $f\ge 24$. This implies $\beta\ge \pi$, and $\beta^2\cdots$ is not a vertex. This further implies no $\thick\delta\thin\alpha\thin\cdots\thin\alpha\thin\delta\thick$. In particular, $\delta\thin\delta\cdots$ is not a vertex.

By $\alpha^2\beta$ and $\alpha<\gamma$, we get $R(\beta\gamma)<R(\alpha\beta)=\alpha<\beta,\gamma,\epsilon$. This implies $\alpha\beta\cdots=\alpha^2\beta,\alpha\beta\delta^k$, and $\beta\gamma\cdots=\beta\gamma\delta^k$. By $\gamma^3$, we get $R(\gamma^2)=\gamma<2\alpha,\beta,\delta+\epsilon$. By $\delta<\epsilon$, this implies $\gamma^2\cdots=\gamma^3,\gamma^2\delta^k,\alpha\gamma^2\delta^k$. Then by no $\thick\delta\thin\alpha\thin\cdots\thin\alpha\thin\delta\thick$, we get $\alpha\beta\cdots=\alpha^2\beta,\alpha\beta\delta^2$, and $\beta\thin\gamma\cdots=\beta\gamma\delta^2$, and  $\gamma\thin\gamma\cdots=\gamma^3,\gamma^2\delta^2,\alpha\gamma^2\delta^2$. 

Suppose $\alpha\epsilon^2$ is a vertex. The AAD $\thick^{\delta}\epsilon^{\gamma}\thin^{\beta}\alpha^{\gamma}\thin^{\gamma}\epsilon^{\delta}\thick$ of $\alpha\epsilon^2$ implies $\beta\thin\gamma\cdots=\beta\gamma\delta^2$ is a vertex. Then by distinct $\alpha,\beta,\gamma$ values, we get $\alpha\beta\cdots=\alpha^2\beta$ and $\gamma\thin\gamma\cdots=\gamma^3,\alpha\gamma^2\delta^2$. The angle sums of $\alpha\epsilon^2,\beta\gamma\delta^2$ further imply
\[
\alpha=\tfrac{4}{9}\pi,\,
\beta=\tfrac{10}{9}\pi,\,
\gamma=\tfrac{2}{3}\pi,\,
\delta=\tfrac{1}{9}\pi,\,
\epsilon=\tfrac{7}{9}\pi,\,
f=36.
\]
By $R(\epsilon^2)=\alpha=4\delta<\beta,\gamma,\epsilon$ and no $\delta\thin\delta\cdots$, we get $\epsilon^2\cdots=\alpha\epsilon^2$. In particular, we know $\epsilon\thin\epsilon\cdots$ is not a vertex.

Since the second of Figure \ref{bde_2abA} implies a vertex $\epsilon\thin\epsilon\cdots$, we know the special tile is given by the first of Figure \ref{bde_2abA}. Then $\thin^{\epsilon}\gamma^{\alpha}\thin^{\alpha}\gamma^{\epsilon}\thin\cdots=\gamma^3,\alpha\gamma^2\delta^2$ is a vertex. If the vertex is $\gamma^3$, then we get $\thin^{\epsilon}\gamma^{\alpha}\thin^{\alpha}\gamma^{\epsilon}\thin\cdots=\thin^{\alpha}\gamma^{\epsilon}\thin^{\alpha}\gamma^{\epsilon}\thin^{\epsilon}\gamma^{\alpha}\thin$, contradicting no $\epsilon\thin\epsilon\cdots$. Therefore the vertex is $\alpha\gamma^2\delta^2$. By no $\beta^2\cdots$, we get $\thin^{\epsilon}\gamma^{\alpha}\thin^{\alpha}\gamma^{\epsilon}\thin\cdots=\alpha\gamma^2\delta^2=\thick^{\epsilon}\delta^{\beta}\thin^{\epsilon}\gamma^{\alpha}\thin^{\alpha}\gamma^{\epsilon}\thin^{\beta}\alpha^{\gamma}\thin^{\beta}\delta^{\epsilon}\thick$. This determines $T_1,T_2,T_3,T_4$ in the third of Figure \ref{bde_2abA}. Then $\beta_1\epsilon_4\cdots=\beta\delta\epsilon$ determines $T_5$, and $\beta_2\thin\gamma_3\cdots=\beta\gamma\delta^2$ determines $T_6$. Then $\alpha_1\beta_5\cdots=\alpha_2\beta_6\cdots=\alpha^2\beta$ and $\epsilon_1\epsilon_2\cdots=\alpha\epsilon^2$. Then either $\gamma_1\cdots$ or $\gamma_2\cdots$ is $\beta\thin\gamma\cdots=\beta\gamma\delta^2$. This implies $\alpha,\delta$ adjacent, a contradiction. 

We conclude $\alpha\epsilon^2$ is not a vertex. 

If $\beta\gamma\cdots=\beta\gamma\delta^k$ is a vertex, then by $\beta\delta\epsilon$, we get $\epsilon=\gamma+(k-1)\delta>\gamma=\tfrac{2}{3}\pi$. This implies $R(\epsilon^2)<2\alpha,\beta,\gamma,\epsilon$. Therefore $\epsilon^2\cdots=\alpha\delta^k\epsilon^2,\delta^k\epsilon^2$. By no $\alpha\epsilon^2$, we get $k\ge 2$ in both $\alpha\delta^k\epsilon^2,\delta^k\epsilon^2$. This implies the number of $\epsilon$ is always no more than the number of $\delta$ at any vertex. By the counting lemma, this implies the number of $\delta$ and $\epsilon$ are equal at every vertex. In particular, we know $\beta\gamma\cdots=\beta\gamma\delta^k$ is not a vertex.

By no $\beta^2\cdots,\beta\gamma\cdots$, we know $\alpha\thin\delta\cdots$ is not a vertex. Therefore $\alpha\beta\delta^k$ is not a vertex, and we get $\alpha\beta\cdots=\alpha^2\beta$. 

By no $\beta^2\cdots,\beta\gamma\cdots$, we also know the AAD of $\thin\alpha\thin\alpha\thin$ is $\thin^{\beta}\alpha^{\gamma}\thin^{\gamma}\alpha^{\beta}\thin$. This implies the AAD of $\alpha\beta\cdots=\alpha^2\beta$ is $\thin^{\beta}\alpha^{\gamma}\thin^{\gamma}\alpha^{\beta}\thin\beta\thin$. Therefore the AAD of $\thin\alpha\thin\beta\thin$ is $\thin^{\gamma}\alpha^{\beta}\thin\beta\thin$. 

By no $\alpha\thin\delta\cdots,\delta\thin\delta\cdots$, a $\delta^2$-fan has $\thin\beta\thin\delta\thick$ or $\thin\gamma\thin\delta\thick$ at both ends. Then by no $\beta^2\cdots,\beta\gamma\cdots$, a $\delta^2$-fan is $\thick\delta\thin\beta\thin\delta\thick$, $\thick\delta\thin\gamma\thin\cdots\thin\gamma\thin\delta\thick$, with only $\alpha,\gamma$ in $\cdots$. Further by $\gamma^3$ and $2\alpha>\gamma$, a $\delta^2$-fan is $\thick\delta\thin\beta\thin\delta\thick$, $\thick\delta\thin\gamma\thin\delta\thick$, $\thick\delta\thin\gamma\thin\gamma\thin\delta\thick$, $\thick\delta\thin\gamma\thin\alpha\thin\gamma\thin\delta\thick$. The fan $\thick\delta\thin\gamma\thin\delta\thick$ has the AAD $\thick^{\gamma}\delta^{\beta}\thin^{\alpha}\gamma^{\epsilon}\thin^{\beta}\delta^{\epsilon}\thick$, and the fan $\thick\delta\thin\gamma\thin\alpha\thin\gamma\thin\delta\thick$ gives the AAD $\thin^{\gamma}\alpha^{\beta}\thin\gamma\thin^{\beta}\delta^{\epsilon}\thick=\thin^{\gamma}\alpha^{\beta}\thin^{\epsilon}\gamma^{\alpha}\thin^{\beta}\delta^{\epsilon}\thick$ or $\thin^{\gamma}\alpha^{\beta}\thin^{\alpha}\gamma^{\epsilon}\thin^{\beta}\delta^{\epsilon}\thick$. Then we get $\thin^{\beta}\alpha^{\gamma}\thin^{\delta}\beta^{\alpha}\thin$ or $\thin^{\beta}\alpha^{\gamma}\thin^{\alpha}\beta^{\delta}\thin$, contradicting the AAD $\thin^{\gamma}\alpha^{\beta}\thin\beta\thin$ of $\thin\alpha\thin\beta\thin$. Therefore a $\delta^2$-fan is $\thick\delta\thin\beta\thin\delta\thick$, $\thick\delta\thin\gamma\thin\gamma\thin\delta\thick$. Both have values $>\beta\ge\pi$. By $\beta\delta\epsilon$ and $\delta<\epsilon$, $\delta\epsilon$-fan and $\epsilon^2$-fan have values $\ge 2\pi-\beta$. Therefore a $\delta^2$-fan cannot be combined with another fan, and the two $\delta^2$-fans are vertices $\thick\delta\thin\beta\thin\delta\thick=\beta\delta^2$, $\thick\delta\thin\gamma\thin\gamma\thin\delta\thick=\gamma^2\delta^2$. By $\beta\delta\epsilon$, we know $\beta\delta^2$ is not a vertex. This implies $\gamma^2\delta^2$ is the only $\delta^2$-fan. Then by $\beta\epsilon\cdots=\beta\delta\epsilon$, we get $\beta\delta\cdots=\beta\delta\epsilon$.

Suppose $\gamma^2\delta^2$ is a vertex. The angle sum of $\gamma^2\delta^2$ further implies
\[
\alpha=(\tfrac{1}{3}+\tfrac{4}{f})\pi,\;
\beta=(\tfrac{4}{3}-\tfrac{8}{f})\pi,\;
\gamma=\tfrac{2}{3}\pi,\;
\delta=\tfrac{1}{3}\pi,\;
\epsilon=(\tfrac{1}{3}+\tfrac{8}{f})\pi.
\]
Then $R(\gamma\delta\epsilon)<R(\gamma\delta^2)=R(\gamma^2)=\gamma=2\delta<2\alpha,\beta$. By $\beta>\gamma>\alpha$ and $\delta<\epsilon$, and $\gamma^2\delta^2$ being the only $\delta^2$-fan, this implies $\gamma\delta^2\cdots=\gamma^2\cdots=\gamma^2\delta^2$, and $\gamma\delta\epsilon\cdots=\alpha\gamma\delta\epsilon$. 

The angle sum of $\alpha\gamma\delta\epsilon$ further implies
\[
\alpha=\tfrac{4}{9}\pi,\;
\beta=\tfrac{10}{9}\pi,\;
\gamma=\tfrac{2}{3}\pi,\;
\delta=\tfrac{1}{3}\pi,\;
\epsilon=\tfrac{5}{9}\pi,\;
f=36.
\]
By $\gamma^2\delta^2$ and Lemma \ref{balance}, we know $\epsilon^2\cdots$ is a vertex. The angle values imply $\epsilon^2\cdots=\alpha^2\epsilon^2,\delta\epsilon^3$. The AADs of $\alpha^2\epsilon^2,\delta\epsilon^3$ imply $\beta^2\cdots,\beta\gamma\cdots$, a contradiction. Therefore $\gamma\delta\epsilon\cdots=\alpha\gamma\delta\epsilon$ is not a vertex, and we get $\gamma\delta\cdots=\gamma\delta^2\cdots=\gamma^2\delta^2$.

By the the AAD $\thin^{\gamma}\alpha^{\beta}\thin\beta\thin$ of $\thin\alpha\thin\beta\thin$, the AAD of $\thin\gamma\thin\delta\thick$ is $\thin^{\alpha}\gamma^{\epsilon}\thin^{\beta}\delta^{\epsilon}\thick$. This implies the AAD of $\gamma^2\delta^2$ is $\thick^{\epsilon}\delta^{\beta}\thin^{\epsilon}\gamma^{\alpha}\thin^{\alpha}\gamma^{\epsilon}\thin^{\beta}\delta^{\epsilon}\thick$, and determines $T_1,T_2,T_3,T_4$ in the fourth of Figure \ref{bde_2abA}. If $\alpha_2\alpha_3\cdots\ne \alpha^2\beta$, then the vertex is $\thin^{\gamma}\alpha^{\beta}\thin\cdots\thin^{\beta}\alpha^{\gamma}\thin$, where $\cdots$ has no $\beta$. By no $\beta^2\cdots,\beta\gamma\cdots$, we do not have $\thin^{\gamma}\alpha^{\beta}\thin\alpha\thin,\thin^{\gamma}\alpha^{\beta}\thin\delta\thick,\thin^{\gamma}\alpha^{\beta}\thin\epsilon\thick$. Therefore $\thin^{\gamma}\alpha^{\beta}\thin\cdots\thin^{\beta}\alpha^{\gamma}\thin=\thin^{\gamma}\alpha^{\beta}\thin\gamma\thin\cdots\thin\gamma\thin^{\beta}\alpha^{\gamma}\thin$. By $\alpha+\gamma=(1+\frac{4}{f})\pi>\pi$ and no $\alpha^2\gamma$, we get a contradiction. Therefore $\alpha_2\alpha_3\cdots=\alpha^2\beta$. By the symmetry of horizontal flip, we may assume $T_5$ is arranged as indicated. Then $\beta_2\delta_5\cdots=\beta\delta\epsilon$ determines $T_6$, and $\beta_1\epsilon_2\cdots=\beta\delta\epsilon$ determines $T_7$. Then we get $\gamma_6\delta_2\epsilon_7\cdots$, contradicting $\gamma\delta\cdots=\gamma^2\delta^2$. 

We conclude that the only $\delta^2$-fan $\gamma^2\delta^2$ is not a vertex. Then by Lemma \ref{fbalance}, all fans are $\delta\epsilon$-fans. By $\beta\delta\epsilon$ and no $\alpha\thin\delta\cdots,\beta\gamma\cdots$, a fan is either the vertex $\beta\delta\epsilon$ or $\thick\delta\thin\gamma\thin\cdots\thin\epsilon\thick$, where $\cdots$ consists of $\alpha,\gamma$. By $R(\gamma\delta\epsilon)=(\frac{2}{3}-\frac{8}{f})\pi<2\alpha,\beta,\gamma,\delta+\epsilon$, we know the fan $\thick\delta\thin\gamma\thin\cdots\thin\epsilon\thick$ cannot be combined with another fan, and is actually a vertex $\thick\delta\thin\gamma\thin\alpha\thin\epsilon\thick$. By no $\beta\gamma\cdots$, the AAD of $\thick\delta\thin\gamma\thin\alpha\thin\epsilon\thick$ is $\thick^{\epsilon}\delta^{\beta}\thin^{\alpha}\gamma^{\epsilon}\thin^{\beta}\alpha^{\gamma}\thin^{\gamma}\epsilon^{\delta}\thick$ or $\thick^{\epsilon}\delta^{\beta}\thin^{\epsilon}\gamma^{\alpha}\thin^{\beta}\alpha^{\gamma}\thin^{\gamma}\epsilon^{\delta}\thick$. The AAD implies either $\thin^{\beta}\alpha^{\gamma}\thin^{\delta}\beta^{\alpha}\thin$ or $\thin^{\beta}\alpha^{\gamma}\thin^{\alpha}\beta^{\delta}\thin$, contradicting the AAD $\thin^{\gamma}\alpha^{\beta}\thin\beta\thin$ of $\thin\alpha\thin\beta\thin$. 

Therefore $\beta\delta\epsilon$ is the only $b$-vertex. By $\alpha^2\beta$ and applying the counting lemma to $\beta,\delta$, we get a contradiction.

\subsubsection*{Case. $\gamma\epsilon^2$ is a vertex}

The angle sums of $\beta\delta\epsilon,\alpha^2\beta,\gamma\epsilon^2$ and the angle sum for pentagon imply
\[
\alpha=(1+\tfrac{4}{f})\pi-\gamma,\;
\beta=2\gamma-\tfrac{8}{f}\pi,\;
\delta=(1+\tfrac{8}{f})\pi-\tfrac{3}{2}\gamma,\;
\epsilon=\pi-\tfrac{1}{2}\gamma.
\]
We have $\beta+\gamma+2\delta=(2+\frac{8}{f})\pi>2\pi$. 

If $\beta<\gamma$ and $\delta>\epsilon$, then $\gamma<\tfrac{8}{f}\pi$. This implies $\alpha>\gamma$. By $\beta\delta\epsilon$ and $\delta>\epsilon$, we get $R(\delta^2)<\beta<\gamma<\alpha$. By $\beta\delta\epsilon$ and $\beta<\gamma<\pi$, we get $\delta+\epsilon>\pi$. Then by $\delta>\epsilon$, we know $R(\delta^2)$ has no $\delta,\epsilon$. Therefore $\delta^2\cdots$ is not a vertex, contradicting $\gamma\epsilon^2$ and the balance lemma. 

By Lemma \ref{geometry1}, therefore, we get $\beta>\gamma$ and $\delta<\epsilon$. Then by $\beta\delta\epsilon$, we get $\beta\epsilon\cdots=\beta\delta\epsilon$. 

We know there is no tiling with $\gamma^3$. By $\gamma\epsilon^2$ and $\alpha\ne\gamma$, we know $\alpha\epsilon^2$ is not a vertex. By Proposition \ref{bde_a2d_c3e}, we know $\alpha\delta^2$ is not a vertex. The angle sums of $\beta^2\gamma$ further implies
\[
\alpha=(\tfrac{3}{5}+\tfrac{4}{5f})\pi,\,
\beta=\epsilon=(\tfrac{4}{5}-\tfrac{8}{5f})\pi,\,
\gamma=\delta=(\tfrac{2}{5}+\tfrac{16}{5f})\pi,\,
\]
We have $\pi>\beta=\epsilon>\alpha>\gamma=\delta$, contradicting Lemma \ref{geometry4}. 

Therefore $\beta\delta\epsilon,\alpha^2\beta,\gamma\epsilon^2$ are all the degree $3$ vertices. This implies that, in any tile, one of $\alpha\cdots,\gamma\cdots$ has high degree. By Lemma \ref{special_tile}, we get $f\ge 24$. Moreover, in a special tile, we know $\beta\cdots,\delta\cdots,\epsilon\cdots$ have degree $3$. This implies $\delta\cdots=\epsilon\cdots=\beta\delta\epsilon$. Then $\gamma\cdots=\alpha\gamma\cdots,\gamma\delta\cdots$ has high degree. Therefore $\alpha\cdots,\beta\cdots$ have degree $3$, and we get $\alpha\cdots=\beta\cdots=\alpha^2\beta$ as indicated by the first of Figure \ref{bde_2abB}. The vertex $H=\gamma\cdots$ has degree $4$ or $5$. 

\begin{figure}[htp]
\centering
\begin{tikzpicture}[>=latex,scale=1]

\begin{scope}[shift={(-4.5cm,0.3cm)}]

\draw
	(-0.5,-0.4) -- (-0.5,0.7) -- (0,1.1) -- (0.5,0.7) -- (0.5,-0.4)
	(-0.5,0.7) -- ++(-0.4,0)
	(0,1.1) -- ++(0,0.4);
		
\draw[line width=1.2]
	(-0.5,0) -- (0.5,0);

\node at (0,0.85) {\small $\alpha$}; 
\node at (-0.3,0.6) {\small $\beta$};
\node at (0.3,0.6) {\small $\gamma$};
\node at (-0.3,0.2) {\small $\delta$};
\node at (0.3,0.2) {\small $\epsilon$};	

\node at (-0.7,0) {\small $\beta$};
\node at (0.7,0) {\small $\beta$};
\node at (-0.3,-0.2) {\small $\epsilon$};
\node at (0.3,-0.2) {\small $\delta$};
\node at (-0.7,0.5) {\small $\alpha$};
\node at (-0.6,0.9) {\small $\alpha$};
\node at (-0.2,1.2) {\small $\beta$};
\node at (0.2,1.2) {\small $\alpha$};
\node at (0.7,0.7) {\small $H$};

\end{scope}


\begin{scope}[shift={(-2cm,-0.7cm)}]

\draw
	(-0.5,0.7) -- (-0.5,0) -- (0.5,0) -- (0.5,0.7) -- (0,1.1) -- (-0.5,0.7)
	(-0.5,0.7) -- (-1,1.1) -- (-1,2.5) -- (1,2.5) -- (1,1.1) -- (0.5,0.7)
	(0,1.1) -- (0,1.8)
	(-1,1.8) -- (1,1.8)
	(1,1.1) -- (1.4,1.1)
	(0.5,0) -- (0.5,-0.4);

\draw[line width=1.2]
	(-0.5,0) -- (0.5,0)
	(-1,1.8) -- (0,1.8)
	(1,1.1) -- (1,1.8);

\node at (-0.3,0.6) {\small $\beta$};
\node at (0,0.85) {\small $\alpha$};
\node at (-0.3,0.2) {\small $\delta$};
\node at (0.3,0.6) {\small $\gamma$};	
\node at (0.3,0.2) {\small $\epsilon$}; 

\node at (-0.5,0.95) {\small $\alpha$};
\node at (-0.2,1.2) {\small $\beta$};
\node at (-0.8,1.2) {\small $\gamma$};
\node at (-0.2,1.6) {\small $\delta$};
\node at (-0.8,1.6) {\small $\epsilon$};

\node at (0.5,0.95) {\small $\gamma$};
\node at (0.2,1.2) {\small $\alpha$};
\node at (0.8,1.2) {\small $\epsilon$};
\node at (0.2,1.6) {\small $\beta$};
\node at (0.8,1.6) {\small $\delta$};

\node at (0.8,2.3) {\small $\alpha$};
\node at (-0.8,2.3) {\small $\beta$};
\node at (0.8,2) {\small $\gamma$};
\node at (-0.8,2) {\small $\delta$};
\node at (0,1.95) {\small $\epsilon$};

\node at (1,0.85) {\small $\beta$};
\node at (0.7,0) {\small $\beta$};
\node at (1.2,1.3) {\small $\delta$};
\node at (0.3,-0.2) {\small $\delta$};

\node[inner sep=0.5,draw,shape=circle] at (0,0.4) {\small $1$};
\node[inner sep=0.5,draw,shape=circle] at (-0.5,1.4) {\small $2$};
\node[inner sep=0.5,draw,shape=circle] at (0.5,1.4) {\small $3$};
\node[inner sep=0.5,draw,shape=circle] at (0.4,2.15) {\small $4$};

\end{scope}


\foreach \a in {0,1}
{
\begin{scope}[shift={(3.5*\a cm,-0.7*\a cm)}]

\draw
	(0,-0.7) -- (0,0.7) -- (0.5,1.1) -- (1,0.7) -- (1,-0.7) -- (0.5,-1.1) -- (0,-0.7)
	(1,-0.7) -- (1.8,-0.7) -- (1.8,0.3)
	(1,0.7) -- (1.8,0.3) -- (2.5,0.3) -- (2.5,1.1) -- (1.5,1.1)
	(1,0.7) -- (1.5,1.1) -- (1.5,1.8) -- (0.5,1.8) -- (0.5,1.1);

\draw[line width=1.2]
	(0,0) -- (1,0)
	(1.5,1.8) -- (1.5,1.1);

\node at (0.5,0.85) {\small $\alpha$};	
\node at (0.2,0.6) {\small $\beta$};
\node at (0.8,0.6) {\small $\gamma$};
\node at (0.2,0.2) {\small $\delta$};
\node at (0.8,0.2) {\small $\epsilon$}; 

\node at (0.5,-0.85) {\small $\alpha$};	
\node at (0.2,-0.6) {\small $\gamma$};
\node at (0.8,-0.6) {\small $\beta$};
\node at (0.2,-0.2) {\small $\epsilon$};
\node at (0.8,-0.2) {\small $\delta$}; 

\node at (0.7,1.2) {\small $\alpha$};
\node at (1,0.95) {\small $\gamma$};
\node at (0.7,1.6) {\small $\beta$};
\node at (1.3,1.2) {\small $\epsilon$};
\node at (1.3,1.6) {\small $\delta$};

\node[inner sep=0.5,draw,shape=circle] at (0.5,0.4) {\small $1$};
\node[inner sep=0.5,draw,shape=circle] at (0.5,-0.4) {\small $2$};
\node[inner sep=0.5,draw,shape=circle] at (1,1.4) {\small $5$};
\node[inner sep=0.5,draw,shape=circle] at (1.95,0.8) {\small $6$};
\node[inner sep=0.5,draw,shape=circle] at (1.5,-0.15) {\small $7$};

\end{scope}
}


\draw
	(0.5,-1.1) -- (0.5,-1.8) -- (1.8,-1.8) -- (1.8,-0.7);
	
\draw[line width=1.2]
	(1,-0.7) -- (1.8,-0.7)
	(1.8,0.3) -- (2.5,0.3);

\node at (1.6,0.2) {\small $\gamma$};
\node at (1.6,-0.5) {\small $\epsilon$};
\node at (1.2,0.4) {\small $\alpha$};
\node at (1.2,-0.5) {\small $\delta$};
\node at (1.2,0) {\small $\beta$};

\node at (1.55,0.9) {\small $\alpha$};
\node at (2.3,0.9) {\small $\beta$};
\node at (1.3,0.7) {\small $\gamma$};
\node at (2.3,0.5) {\small $\delta$};
\node at (1.8,0.5) {\small $\epsilon$};

\node at (0.7,-1.2) {\small $\gamma$};
\node at (1,-0.9) {\small $\epsilon$};
\node at (0.7,-1.6) {\small $\alpha$};
\node at (1.6,-0.9) {\small $\delta$};
\node at (1.6,-1.6) {\small $\beta$};

\node at (2,-0.7) {\small $\beta$};

\node[inner sep=0.5,draw,shape=circle] at (1.15,-1.35) {\small $8$};


\begin{scope}[shift={(3.5cm,-0.7cm)}]

\draw
	(1.5,1.8) -- (1.5,2.5) -- (-0.5,2.5) -- (-0.5,1.1) -- (0,0.7);
	
\draw[line width=1.2]
	(1,0.7) -- (1.8,0.3)
	(-0.5,1.8) -- (0.5,1.8);
	
\node at (1.6,0.2) {\small $\epsilon$};
\node at (1.6,-0.5) {\small $\gamma$};
\node at (1.2,0.4) {\small $\delta$};
\node at (1.2,-0.5) {\small $\alpha$};
\node at (1.2,0) {\small $\beta$};

\node at (1.55,0.9) {\small $\beta$};
\node at (2.3,0.9) {\small $\alpha$};
\node at (1.3,0.75) {\small $\delta$};
\node at (2.3,0.5) {\small $\gamma$};
\node at (1.8,0.5) {\small $\epsilon$};

\node at (0.3,1.6) {\small $\delta$};
\node at (0.3,1.2) {\small $\beta$};
\node at (-0.3,1.6) {\small $\epsilon$};
\node at (0,0.95) {\small $\alpha$};
\node at (-0.3,1.2) {\small $\gamma$};

\node at (1.3,2.3) {\small $\alpha$};
\node at (-0.3,2.3) {\small $\beta$};
\node at (1.3,2) {\small $\gamma$};
\node at (-0.3,2) {\small $\delta$};
\node at (0.5,1.95) {\small $\epsilon$};

\node at (1.7,1.3) {\small $\delta$};

\node[inner sep=0.5,draw,shape=circle] at (0,1.4) {\small $8$};
\node[inner sep=0.5,draw,shape=circle] at (0.5,2.25) {\small $9$};

\end{scope}

\end{tikzpicture}
\caption{Proposition \ref{bde_2ab}: $\gamma\epsilon^2$ or $\alpha\epsilon^2$ is a vertex.}
\label{bde_2abB}
\end{figure}
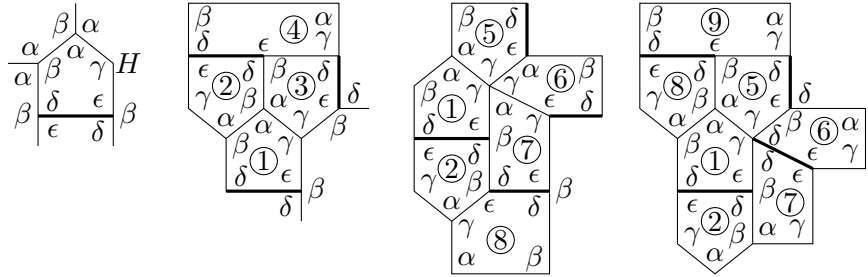

We divide the further discussion by whether the vertex $H$ in the first of Figure \ref{bde_2abB} is a $\hat{b}$-vertex or a $b$-vertex.

\subsubsection*{Subcase. $H$ is a $\hat{b}$-vertex}

In the first of Figure \ref{bde_2abB}, the $\hat{b}$-vertex $H={}^{\beta}\thin^{\epsilon}\gamma^{\alpha}\thin^{\alpha}\cdots=\thin^{\gamma}\alpha^{\beta}\thin^{\epsilon}\gamma^{\alpha}\thin^{\alpha}\beta^{\delta}\thin\cdots$ or $\thin^{\gamma}\alpha^{\beta}\thin^{\epsilon}\gamma^{\alpha}\thin^{\alpha}\gamma^{\epsilon}\thin\cdots$. By  $\alpha+\gamma>\pi$ and $\beta>\gamma$, the remainder of $H$ consists of one or two from $\beta,\gamma$. Then by $2\gamma>\beta>\gamma$, we know either $H=\alpha\gamma^3$, or the angle sum of $H$ implies $\alpha+\beta+2\gamma\le 2\pi$. 

Suppose $H\ne \alpha\gamma^3$. Then $\alpha+\beta+2\gamma=(1-\tfrac{4}{f})\pi+3\gamma\le 2\pi$. This means $\gamma\le \frac{1}{3}(1+\frac{4}{f})\pi$. Then $\alpha\ge 2\gamma>\beta>\gamma$, and $\epsilon>\delta>\frac{1}{2}\pi>\gamma$. This implies $\delta\thin\delta\cdots,\epsilon\thin\epsilon\cdots,\delta\thin\epsilon\cdots$ are not vertices. 

By $\alpha+\gamma+2\delta>\beta+\gamma+2\delta>2\pi$, and $\delta<\epsilon$, we get $R(\alpha\epsilon^2)<R(\alpha\delta\epsilon)<R(\alpha\delta^2)<\gamma<\alpha,\beta,\delta,\epsilon$. By all degree $3$ vertices $\beta\delta\epsilon,\alpha^2\beta,\gamma\epsilon^2$, this implies $\alpha\beta\cdots,\alpha\epsilon\cdots$ are not vertices. 

We have $H=\thin^{\gamma}\alpha^{\beta}\thin^{\epsilon}\gamma^{\alpha}\thin^{\alpha}\beta^{\delta}\thin\theta\thin\cdots$ or $\thin^{\gamma}\alpha^{\beta}\thin^{\epsilon}\gamma^{\alpha}\thin^{\alpha}\gamma^{\epsilon}\thin\theta\thin\cdots$, with $\theta=\beta,\gamma$. The AAD of $H$ implies $\alpha\delta\cdots,\alpha\epsilon\cdots,\delta\thin\delta\cdots,\epsilon\thin\epsilon\cdots,\delta\thin\epsilon\cdots$, a contradiction. 

Therefore $H=\alpha\gamma^3$. The angle sum of $\alpha\gamma^3$ further implies
\[
\alpha=(\tfrac{1}{2}+\tfrac{6}{f})\pi,\,
\beta=(1-\tfrac{12}{f})\pi,\,
\gamma=(\tfrac{1}{2}-\tfrac{2}{f})\pi,\,
\delta=(\tfrac{1}{4}+\tfrac{11}{f})\pi,\,
\epsilon=(\tfrac{3}{4}+\tfrac{1}{f})\pi.
\]
Then $H=\thin^{\gamma}\alpha^{\beta}\thin^{\epsilon}\gamma^{\alpha}\thin^{\alpha}\gamma^{\epsilon}\thin\gamma\thin$ implies $\alpha\epsilon\cdots,\epsilon\thin\epsilon\cdots$. By $\delta<\epsilon$ and $\delta+\epsilon>\pi$, we know $\epsilon\thin\epsilon\cdots$ is not a vertex. By $\delta<\epsilon$, we get $R(\alpha\epsilon^2)<R(\alpha\delta\epsilon)=(\tfrac{1}{2}-\tfrac{18}{f})\pi<\alpha,\beta,\gamma,2\delta,2\epsilon$. By all degree $3$ vertices $\beta\delta\epsilon,\alpha^2\beta,\gamma\epsilon^2$, this implies $\alpha\epsilon\cdots$ is not a vertex.

\subsubsection*{Subcase. $H$ is a $b$-vertex}

By $\beta+\gamma+2\delta>2\pi$ and $\delta<\epsilon$, we know $\beta\gamma\cdots$ is a $\hat{b}$-vertex. Therefore the $b$-vertex $H$ has no $\beta$, and $H=\thin^{\epsilon}\gamma^{\alpha}\thin^{\alpha}\gamma^{\epsilon}\thin^{\beta}\cdots$. By $\beta\delta\epsilon$, and $\beta<2\gamma$, and $\delta<\epsilon$, we know $H$ has no $\epsilon$. Then by $3\gamma+2\delta>\beta+\gamma+2\delta>2\pi$, the $b$-vertex $H=\gamma^2\delta^2,\alpha\gamma^2\delta^2$. 

Suppose $H=\gamma^2\delta^2$. The angle sum of $\gamma^2\delta^2$ further implies
\[
\alpha=(1-\tfrac{12}{f})\pi,\,
\beta=\tfrac{24}{f}\pi,\,
\gamma=\tfrac{16}{f}\pi,\,
\delta=(1-\tfrac{16}{f})\pi,\,
\epsilon=(1-\tfrac{8}{f})\pi.
\]
By $f\ge 24$, we get $R(\beta\delta^2)=R(\gamma\delta\epsilon)=\frac{8}{f}\pi<\alpha,\beta,\gamma,2\delta,2\epsilon$. This implies $\beta\delta^2\cdots,\gamma\delta\epsilon\cdots$ are not vertices. Then by $\beta\delta\epsilon$, we get $\beta\delta\cdots=\beta\delta\epsilon$. 

The tile $T_1$ in the second of Figure \ref{bde_2abB} is the special tile in the first picture. The special tile and $H=\thick^{\epsilon}\delta^{\beta}\thin^{\epsilon}\gamma^{\alpha}\thin^{\alpha}\gamma^{\epsilon}\thin^{\beta}\delta^{\epsilon}\thick$ determine $T_2,T_3$ and $\beta$ just outside $\epsilon_3$. Then $\beta_3\delta_2\cdots=\beta\delta\epsilon$ determines $T_4$. Then $\epsilon_3\cdots=\beta\epsilon\cdots=\beta\delta\epsilon$ implies $\gamma_4\delta_3\cdots=\gamma\delta\epsilon\cdots$, a contradiction. 

We are left with the special tile with $H=\alpha\gamma^2\delta^2$. Having exhausted all other special tiles, we know there is no special tile with $\deg H=3,4$. By Lemma \ref{special_tile}, this implies $f\ge 60$. 

The angle sum of $\alpha\gamma^2\delta^2$ further implies
\[
\alpha=(\tfrac{1}{2}-\tfrac{6}{f})\pi,\,
\beta=(1+\tfrac{12}{f})\pi,\,
\gamma=(\tfrac{1}{2}+\tfrac{10}{f})\pi,\,
\delta=(\tfrac{1}{4}-\tfrac{7}{f})\pi,\,
\epsilon=(\tfrac{3}{4}-\tfrac{5}{f})\pi.
\]
By $\beta>\pi$, we know $\beta^2\cdots$ is not a vertex. 

By $f\ge 60$, we get $2\delta<\alpha<R(\beta\delta^2)=R(\gamma\delta\epsilon)=(\tfrac{1}{2}+\tfrac{2}{f})\pi<2\alpha,\alpha+2\delta,\beta,\gamma,\epsilon$. This implies $\beta\delta^2\cdots=\beta\delta^k(k\ge 6)$ and $\gamma\delta\epsilon\cdots=\gamma\delta^k\epsilon(k\ge 5)$. The AADs of $\beta\delta^k,\gamma\delta^k\epsilon$ imply $\beta^2\cdots$, a contradiction. Therefore $\beta\delta^2\cdots, \gamma\delta\epsilon\cdots$ are not vertices, and we get $\beta\delta\cdots=\beta\delta\epsilon$. 

Now we may carry out the argument similar to the case $H=\gamma^2\delta^2$, using the second of Figure \ref{bde_2abB}. The key is that, by no $\beta^2\cdots$, we get $H=\thin^{\epsilon}\gamma^{\alpha}\thin^{\alpha}\gamma^{\epsilon}\thin\cdots=\alpha\gamma^2\delta^2=\thick^{\epsilon}\delta^{\beta}\thin^{\gamma}\alpha^{\beta}\thin^{\epsilon}\gamma^{\alpha}\thin^{\alpha}\gamma^{\epsilon}\thin^{\beta}\delta^{\epsilon}\thick$. This still gives $\beta$ just outside $\epsilon_3$. Then we get the same contradiction.

\subsubsection*{Case. $\alpha\delta^2$ or $\alpha\epsilon^2$ is a vertex}

The angle sums of $\beta\delta\epsilon,\alpha^2\beta,\beta^2\gamma$ and the angle sum for pentagon imply
\[
\alpha=(\tfrac{3}{5}+\tfrac{4}{5f})\pi,\,
\beta=(\tfrac{4}{5}-\tfrac{8}{5f})\pi,\,
\gamma=(\tfrac{2}{5}+\tfrac{16}{5f})\pi,\,
\delta+\epsilon=(\tfrac{6}{5}+\tfrac{8}{5f})\pi.
\]
We get $\beta>\alpha>\gamma$. By Lemma \ref{geometry1}, this implies $\delta<\epsilon$. Then by $\beta\delta\epsilon$, we know $\alpha\delta^2$ is not a vertex. If $\alpha\epsilon^2$ is a vertex, then its angle sum further implies
\[
\alpha=(\tfrac{3}{5}+\tfrac{4}{5f})\pi,\,
\beta=(\tfrac{4}{5}-\tfrac{8}{5f})\pi,\,
\gamma=(\tfrac{2}{5}+\tfrac{16}{5f})\pi,\,
\delta=(\tfrac{1}{2}+\tfrac{2}{f})\pi,\,
\epsilon=(\tfrac{7}{10}-\tfrac{2}{5f})\pi.
\]
We have $\pi>\beta>\alpha>\gamma$, and $\beta>\epsilon$, and $\pi>\delta>\gamma$, contradicting Lemma \ref{geometry4}.

We also know there is no tiling with $\gamma^3,\gamma\epsilon^2$. Therefore $\beta\delta\epsilon,\alpha^2\beta,\alpha\delta^2,\alpha\epsilon^2$ are all the degree $3$ vertices. In particular, $\gamma\cdots$ always has high degree. By Lemma \ref{special_tile}, this implies $f\ge 24$.  By $\beta\ne\gamma$ and Lemma \ref{geometry1}, we know $\alpha\delta^2,\alpha\epsilon^2$ are mutually exclusive. Then we get $\delta\cdots=\epsilon\cdots=\beta\delta\epsilon$ in a special tile, and a special tile is the first of Figure \ref{bde_2abB}. 

By Lemma \ref{ndegree3}, we know one of $\alpha\gamma^3,\beta\gamma^3,\gamma^4,\gamma^5$ is a vertex. The angle sums of $\beta\delta\epsilon,\alpha^2\beta$, and the angle sum of one of these, and the angle sum for pentagon imply
\begin{align*}
\alpha\gamma^3 &\colon 
	\alpha=(\tfrac{1}{2}+\tfrac{6}{f})\pi,\,
	\beta=(1-\tfrac{12}{f})\pi,\,
	\gamma=(\tfrac{1}{2}-\tfrac{2}{f})\pi.  \\
\beta\gamma^3 &\colon 
	\alpha=(\tfrac{3}{5}+\tfrac{12}{5f})\pi,\,
	\beta=(\tfrac{4}{5}-\tfrac{24}{5f})\pi,\,
	\gamma=(\tfrac{2}{5}+\tfrac{8}{5f})\pi. \\
\gamma^4 &\colon 
	\alpha=(\tfrac{1}{2}+\tfrac{4}{f})\pi,\,
	\beta=(1-\tfrac{8}{f})\pi,\,
	\gamma=\tfrac{1}{2}\pi.  \\
\gamma^5 &\colon 
	\alpha=(\tfrac{3}{5}+\tfrac{4}{f})\pi,\,
	\beta=(\tfrac{4}{5}-\tfrac{8}{f})\pi,\,
	\gamma=\tfrac{2}{5}\pi.
\end{align*}
By $f\ge 24$, we get $\pi>2\gamma>\alpha,\beta>\gamma$. Then by Lemma \ref{geometry1}, we get $\delta<\epsilon$. 

Suppose $\alpha\delta^2$ is a vertex. By the balance lemma, we know $\epsilon^2\cdots$ is a vertex. By $\alpha\delta^2$ and $\alpha<\pi$, we get $\epsilon>\delta>\frac{1}{2}\pi$. This implies $R(\epsilon^2)$ has no $\delta,\epsilon$. Then by $2\gamma+2\epsilon>\alpha+2\delta=2\pi$, and $\alpha,\beta>\gamma$, we get $\epsilon^2\cdots=\gamma\epsilon^2$, contradicting no degree $3$ vertex $\gamma\cdots$. 

Suppose $\alpha\epsilon^2$ is a vertex. The angle sums of $\alpha\epsilon^2$ further implies
\begin{align*}
\alpha\gamma^3 &\colon 
	\alpha=(\tfrac{1}{2}+\tfrac{6}{f})\pi,\,
	\beta=(1-\tfrac{12}{f})\pi,\,
	\gamma=(\tfrac{1}{2}-\tfrac{2}{f})\pi,\,
	\delta=(\tfrac{1}{4}+\tfrac{15}{f})\pi,\,
	\epsilon=(\tfrac{3}{4}-\tfrac{3}{f})\pi.  \\
\beta\gamma^3 &\colon 
	\alpha=(\tfrac{3}{5}+\tfrac{12}{5f})\pi,\,
	\beta=(\tfrac{4}{5}-\tfrac{24}{5f})\pi,\,
	\gamma=(\tfrac{2}{5}+\tfrac{8}{5f})\pi,\,
	\delta=(\tfrac{1}{2}+\tfrac{6}{f})\pi,\,
	\epsilon=(\tfrac{7}{10}-\tfrac{6}{5f})\pi. \\
\gamma^4 &\colon 
	\alpha=(\tfrac{1}{2}+\tfrac{4}{f})\pi,\,
	\beta=(1-\tfrac{8}{f})\pi,\,
	\gamma=\tfrac{1}{2}\pi,\,
	\delta=(\tfrac{1}{4}+\tfrac{10}{f})\pi,\,
	\epsilon=(\tfrac{3}{4}-\tfrac{2}{f})\pi.  \\
\gamma^5 &\colon 
	\alpha=(\tfrac{3}{5}+\tfrac{4}{f})\pi,\,
	\beta=(\tfrac{4}{5}-\tfrac{8}{f})\pi,\,
	\gamma=\tfrac{2}{5}\pi,\,
	\delta=(\tfrac{1}{2}+\tfrac{10}{f})\pi,\,
	\epsilon=(\tfrac{7}{10}-\tfrac{2}{f})\pi.
\end{align*}
Moreover, the AAD of $\alpha\epsilon^2$ implies $\beta\gamma\cdots$ is a vertex. 

For $\beta\gamma^3$ and $\gamma^5$, by $\delta<\epsilon$, we get $f>36$ and $f>60$, respectively. This implies $\pi>\beta>\epsilon>\alpha>\delta>\gamma$, contradicting Lemma \ref{geometry4}.

For $\alpha\gamma^3$ and $\gamma^4$, by $\delta<\epsilon$, we get $f>36$ and $f>24$, respectively. This implies $\beta>\alpha>\gamma$, and $R(\beta\gamma)=(\tfrac{1}{2}+\tfrac{14}{f})\pi$ or $(\tfrac{1}{2}+\tfrac{8}{f})\pi$ satisfies $R(\beta\gamma)<2\gamma,2\delta$. Then by $\beta>\alpha>\gamma$ and $\delta<\epsilon$, we know $\beta\gamma\cdots$ has degree $3$, contradicting all degree $3$ vertices $\beta\delta\epsilon,\alpha^2\beta,\alpha\epsilon^2$.
\end{proof}

\begin{proposition}\label{bde_2bc}
There is no tiling, such that $\alpha,\beta,\gamma$ have distinct values, and $\beta\delta\epsilon,\beta^2\gamma$ are vertices.
\end{proposition}

\begin{proof}
By $\beta^2\gamma$, and Lemma \ref{geometry2}, and distinct $\alpha,\beta,\gamma$ values, the only degree $3$ vertices besides $\beta\delta\epsilon,\beta^2\gamma$ are $\alpha^3,\alpha^2\beta,\alpha\beta\gamma,\alpha\gamma^2,\alpha\delta^2,\alpha\epsilon^2,\gamma\epsilon^2$. By $\beta\delta\epsilon$ and Propositions \ref{bde_abc}, \ref{bde_a2c}, \ref{bde_2ab}, we know $\alpha\beta\gamma,\alpha^2\beta,\alpha\gamma^2$ are not vertices. Therefore $\beta\delta\epsilon,\beta^2\gamma,\alpha^3,\alpha\delta^2,\alpha\epsilon^2,\gamma\epsilon^2$ are all the degree $3$ vertices. By $\alpha\ne\gamma$, and $\delta\ne\epsilon$, and Proposition \ref{bde_a2d_c3e}, we know $\alpha\delta^2,\alpha\epsilon^2,\gamma\epsilon^2$ are mutually exclusive. 

In Figure \ref{bde_2bcA}, we consider a special tile $T_1$. If $\delta_1\cdots=\epsilon_1\cdots=\delta\epsilon\cdots$ , and both have degree $3$, then they are $\beta\delta\epsilon$. This implies $\gamma_1\cdots=\alpha\gamma\cdots,\gamma\delta\cdots$ has high degree. Then $\alpha_1\cdots,\beta_1\cdots$ have degree $3$. Then $\beta_1\cdots=\beta\delta\epsilon$. This implies $\alpha_1\cdots=\alpha\gamma\cdots$ has high degree, a contradiction. By $\alpha\delta^2,\alpha\epsilon^2,\gamma\epsilon^2$ mutually exclusive, we also cannot have $\delta_1\cdots=\delta^2\cdots$, and $\epsilon_1\cdots=\epsilon^2\cdots$, and both have degree $3$. 

Therefore one of $\delta_1\cdots,\epsilon_1\cdots$ has high degree. This means $\alpha_1\cdots,\beta_1\cdots,\gamma_1\cdots$ have degree $3$. This implies $\gamma_1\cdots=\beta^2\gamma$, and $\alpha_1\cdots=\alpha^3,\alpha\delta^2$. In the first picture, we have $\alpha_1\cdots=\alpha^3$. This determines $T_3$, gives $\alpha_2$, and implies the degree $3$ vertex $\beta_1\cdots=\beta^2\gamma$. In the second picture, we have $\alpha_1\cdots=\alpha\delta^2$. This determines $T_2,T_3$, and implies the degree $3$ vertex $\beta_1\cdots=\beta^2\gamma$. 

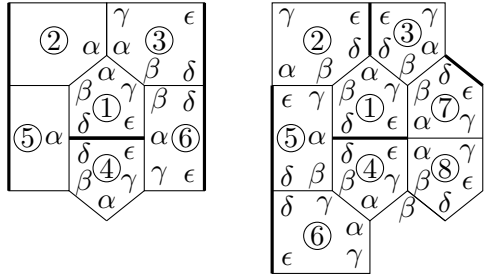
\begin{figure}[htp]
\centering
\begin{tikzpicture}[>=latex,scale=1]


\foreach \a in {-1,1}
{
\begin{scope}[xscale=\a]

\draw
	(0,-1.1) -- (0.5,-0.7) -- (0.5,0.7) -- (0,1.1) -- (0,1.8) -- (1.3,1.8) -- (1.3,-0.7) -- (0.5,-0.7)
	(1.3,0.7) -- (0.5,0.7);

\draw[line width=1.2]
	(-0.5,0) -- (0.5,0)
	(-1.3,1.8) -- (-1.3,-0.7);

\end{scope}
}

\node at (0,0.85) {\small $\alpha$};
\node at (-0.3,0.6) {\small $\beta$};
\node at (0.3,0.6) {\small $\gamma$};
\node at (-0.3,0.2) {\small $\delta$}; 
\node at (0.3,0.2) {\small $\epsilon$};		
	
\node at (0,-0.85) {\small $\alpha$};
\node at (-0.3,-0.6) {\small $\beta$};
\node at (0.3,-0.6) {\small $\gamma$};
\node at (-0.3,-0.2) {\small $\delta$}; 
\node at (0.3,-0.2) {\small $\epsilon$};	

\node at (0.7,0) {\small $\alpha$};
\node at (0.7,0.5) {\small $\beta$};
\node at (0.7,-0.5) {\small $\gamma$};
\node at (1.1,0.5) {\small $\delta$};
\node at (1.1,-0.5) {\small $\epsilon$};

\node at (0.6,0.9) {\small $\beta$};
\node at (0.2,1.2) {\small $\alpha$};
\node at (1.1,0.9) {\small $\delta$};
\node at (1.1,1.6) {\small $\epsilon$};
\node at (0.2,1.6) {\small $\gamma$};

\node at (-0.2,1.2) {\small $\alpha$};

\node at (-0.7,0) {\small $\alpha$};

\node[inner sep=0.5,draw,shape=circle] at (0,0.4) {\small $1$};
\node[inner sep=0.5,draw,shape=circle] at (-0.7,1.3) {\small $2$};
\node[inner sep=0.5,draw,shape=circle] at (0,-0.4) {\small $4$};
\node[inner sep=0.5,draw,shape=circle] at (0.7,1.3) {\small $3$};
\node[inner sep=0.5,draw,shape=circle] at (-1.05,0) {\small $5$};
\node[inner sep=0.5,draw,shape=circle] at (1.05,0) {\small $6$};


\begin{scope}[xshift=3.5cm]

\foreach \a in {0,1}
\draw[xshift=\a cm]
	(-0.5,-0.7) -- (-0.5,0.7) -- (0,1.1) -- (0.5,0.7) -- (0.5,-0.7) -- (0,-1.1) -- (-0.5,-0.7);

\draw
	(0,-1.1) -- (0,-1.8) -- (-1.3,-1.8) -- (-1.3,1.8) -- (1,1.8) -- (1,1.1)
	(0.5,0) -- (1.5,0)
	(0.5,-0.7) -- (1,-1.1)
	(-0.5,0.7) -- (-1.3,0.7)
	(-0.5,-0.7) -- (-1.3,-0.7) ;

\draw[line width=1.2]
	(0.5,0) -- (-0.5,0)
	(0,1.1) -- (0,1.8)
	(-1.3,0.7) -- (-1.3,-1.8)
	(1,1.1) -- (1.5,0.7);
	
\foreach \a in {1,-1}
{
\begin{scope}[yscale=\a]

\node at (0.7,0.6) {\small $\beta$};
\node at (1,0.85) {\small $\delta$};
\node at (0.7,0.2) {\small $\alpha$};
\node at (1.3,0.6) {\small $\epsilon$};	
\node at (1.3,0.2) {\small $\gamma$}; 

\node at (0,0.85) {\small $\alpha$};
\node at (-0.3,0.6) {\small $\beta$};
\node at (0.3,0.6) {\small $\gamma$};
\node at (-0.3,0.2) {\small $\delta$}; 
\node at (0.3,0.2) {\small $\epsilon$};		

\end{scope}
}

\node at (0.8,1.6) {\small $\gamma$};
\node at (0.8,1.2) {\small $\alpha$};
\node at (0.2,1.6) {\small $\epsilon$};
\node at (0.5,0.95) {\small $\beta$};
\node at (0.2,1.2) {\small $\delta$};

\node at (0.5,-0.95) {\small $\beta$};

\node at (-0.6,0.9) {\small $\beta$};
\node at (-0.2,1.2) {\small $\delta$};
\node at (-1.1,0.9) {\small $\alpha$};
\node at (-1.1,1.6) {\small $\gamma$};
\node at (-0.2,1.6) {\small $\epsilon$};

\node at (-0.6,-0.9) {\small $\gamma$};
\node at (-0.2,-1.2) {\small $\alpha$};
\node at (-1.1,-0.9) {\small $\delta$};
\node at (-1.1,-1.6) {\small $\epsilon$};
\node at (-0.2,-1.6) {\small $\gamma$};

\node at (-1.1,0.5) {\small $\epsilon$};
\node at (-1.1,-0.5) {\small $\delta$};
\node at (-0.7,0.5) {\small $\gamma$};
\node at (-0.7,-0.5) {\small $\beta$};
\node at (-0.7,0) {\small $\alpha$};

\node[inner sep=0.5,draw,shape=circle] at (0,0.4) {\small $1$};
\node[inner sep=0.5,draw,shape=circle] at (-0.7,1.3) {\small $2$};
\node[inner sep=0.5,draw,shape=circle] at (0.5,1.4) {\small $3$};
\node[inner sep=0.5,draw,shape=circle] at (0,-0.4) {\small $4$};
\node[inner sep=0.5,draw,shape=circle] at (-1.05,0) {\small $5$};
\node[inner sep=0.5,draw,shape=circle] at (-0.7,-1.3) {\small $6$};
\node[inner sep=0.5,draw,shape=circle] at (1,0.4) {\small $7$};
\node[inner sep=0.5,draw,shape=circle] at (1,-0.4) {\small $8$};

\end{scope}

\end{tikzpicture}
\caption{Proposition \ref{bde_2bc}: Special tile and special companion pair.}
\label{bde_2bcA}
\end{figure}

In both pictures, by $\beta_1\cdots=\gamma_1\cdots=\beta^2\gamma$, we cannot have $\delta_1\cdots=\beta\delta\epsilon$ or $\epsilon_1\cdots=\beta\delta\epsilon$. Therefore the companion pair is matched, and we get $T_4$. Moreover, we have $\delta_1\delta_4\cdots=\alpha\delta^2$ or $\epsilon_1\epsilon_4\cdots=\alpha\epsilon^2$. Therefore $\beta^2\gamma$ is a vertex, and one of $\alpha\delta^2,\alpha\epsilon^2$ is a vertex. 

In the first picture, we know $\alpha^3$ is also a vertex. The angle sums of $\beta\delta\epsilon,\beta^2\gamma,\alpha^3$, and the angle sum of one of $\alpha\delta^2,\alpha\epsilon^2$, and the angle sum for pentagon imply
\begin{align*}
\alpha\delta^2 &\colon
	\alpha=\delta=\tfrac{2}{3}\pi,\,
	\beta=(\tfrac{5}{6}-\tfrac{2}{f})\pi,\,
	\gamma=(\tfrac{1}{3}+\tfrac{4}{f})\pi,\,
	\epsilon=(\tfrac{1}{2}+\tfrac{2}{f})\pi. \\
\alpha\epsilon^2 &\colon
	\alpha=\epsilon=\tfrac{2}{3}\pi,\,
	\beta=(\tfrac{5}{6}-\tfrac{2}{f})\pi,\,
	\gamma=(\tfrac{1}{3}+\tfrac{4}{f})\pi,\,
	\delta=(\tfrac{1}{2}+\tfrac{2}{f})\pi.
\end{align*}
We have $\delta,\epsilon>\frac{1}{2}\pi$. This implies $\delta\thin\delta\cdots,\delta\thin\epsilon\cdots,\epsilon\thin\epsilon\cdots$ are not vertices. If $\delta_1\delta_4\cdots=\alpha\delta^2$. Then we get $\alpha_5$, and $T_2,T_5$ share a vertex $\delta\thin\delta\cdots,\delta\thin\epsilon\cdots,\epsilon\thin\epsilon\cdots$, a contradiction. If $\epsilon_1\epsilon_4\cdots=\alpha\epsilon^2$, then we determine $T_6$ and get $\delta_3\thin\delta_6\cdots$, a contradiction. 

In the second picture, by $\alpha_1\cdots=\alpha\delta^2$, we know $\beta\delta\epsilon,\beta^2\gamma,\alpha^3,\alpha\delta^2$ are all the degree $3$ vertices. Therefore $\delta_1\delta_4\cdots=\alpha\delta^2$, and $\epsilon_1\epsilon_4\cdots$ has high degree. This determines $T_5$. 

Since there is no tiling with the special tile in the first picture, a special tile must be $T_1$ in the second picture. Therefore we may assume this is part of a special companion pair $T_1,T_4$. If both $\alpha_4\cdots,\beta_4\cdots$ have degree $3$, then $\alpha_4\cdots=\alpha^3$ and $\beta_4\cdots=\beta_4\beta_5\cdots=\beta^2\gamma$. This determines $T_6$, and we have $\delta_5\thin\delta_6\cdots$. On the other hand, the vertices $\alpha_1\cdots=\alpha\delta^2$ and $\alpha_4\cdots=\alpha^3$ imply the angle values above, which imply $\delta\thin\delta\cdots$ is not a vertex. The contradiction implies one of $\alpha_4\cdots,\beta_4\cdots$ has high degree. By Lemma \ref{bb_pair}, this implies $\epsilon_1\epsilon_4\cdots$ has degree $4$, and $\gamma_4\cdots$ has degree $3$. By all degree $3$ vertices $\beta\delta\epsilon,\beta^2\gamma,\alpha^3,\alpha\delta^2$, we get $\gamma_4\cdots=\beta^2\gamma$. Then we get $\epsilon_1\epsilon_4\cdots=\alpha^2\epsilon^2,\delta^2\epsilon^2$. The angle sums of $\beta\delta\epsilon,\beta^2\gamma,\delta^2\epsilon^2$ imply $\gamma=0$, a contradiction. Therefore $\epsilon_1\epsilon_4\cdots=\alpha^2\epsilon^2$. This determines $T_7,T_8$.

The angle sums of $\beta\delta\epsilon,\beta^2\gamma,\alpha\delta^2,\alpha^2\epsilon^2$ and the angle sum for pentagon imply 
\[
	\alpha=(\tfrac{1}{2}-\tfrac{2}{f})\pi,\,
	\beta=(\tfrac{3}{4}-\tfrac{3}{f})\pi,\,
	\gamma=(\tfrac{1}{2}+\tfrac{6}{f})\pi,\,
	\delta=(\tfrac{3}{4}+\tfrac{1}{f})\pi,\,
	\epsilon=(\tfrac{1}{2}+\tfrac{2}{f})\pi. 
\]
We have $\alpha<\beta$ and $\delta>\epsilon$. By Lemma \ref{geometry1}, this implies $\beta<\gamma$. Then by $\alpha<\epsilon$, we get $\alpha<\beta,\gamma,\delta,\epsilon$. Then by $R(\gamma_7\gamma_8)=(1-\frac{12}{f})\pi<2\alpha$, and all degree $3$ vertices $\beta\delta\epsilon,\beta^2\gamma,\alpha^3,\alpha\delta^2$, we get a contradiction.
\end{proof}

\begin{proposition}\label{bde_all}
There is no tiling, such that $\alpha,\beta,\gamma$ have distinct values, and $\beta\delta\epsilon$ is a vertex, and one of $\alpha^3,\alpha\beta^2,\gamma^3$ is a vertex.
\end{proposition}

\begin{proof}
By $\beta\delta\epsilon$ and Lemma \ref{geometry2}, the only other degree $3$ $b$-vertices are $\alpha\delta^2,\alpha\epsilon^2,\gamma\epsilon^2$. By distinct $\alpha,\beta,\gamma$ values and Proposition \ref{bde_a2d_c3e}, we know they are mutually exclusive. By Propositions \ref{bde_abc}, \ref{bde_2ac}, \ref{bde_a2c}, \ref{bde_b2c}, \ref{bde_3b}, \ref{bde_2ab}, \ref{bde_2bc}, we also know the only degree $3$ $\hat{b}$-vertices are $\alpha^3,\alpha\beta^2,\gamma^3$. 

In a tile, if the vertices $\delta\cdots,\epsilon\cdots$ have degree $3$, then by $\alpha\delta^2,\alpha\epsilon^2,\gamma\epsilon^2$ mutually exclusive, we know $\delta\cdots=\epsilon\cdots=\beta\delta\epsilon$. This implies the vertex $\gamma\cdots=\alpha\gamma\cdots,\gamma\delta\cdots$ has high degree. 

Therefore in any tile, one of $\gamma\cdots,\delta\cdots,\epsilon\cdots$ has high degree. In particular, in a special tile, we know both $\alpha\cdots,\beta\cdots$ have degree $3$. Then $\beta\cdots=\beta\delta\epsilon,\alpha\beta^2$. The first of Figure \ref{bde_allA} is the case $\beta\cdots=\beta\delta\epsilon$. Then we get $\alpha\cdots=\alpha\beta^2$, and angles of $\beta\cdots=\beta\delta\epsilon$ are arranged as indicated. This implies $\gamma\cdots=\alpha\gamma\cdots,\gamma\delta\cdots$ and $\delta\cdots=\gamma\delta\cdots$ have high degrees, a contradiction. 

\begin{figure}[htp]
\centering
\begin{tikzpicture}[>=latex,scale=1]

\foreach \a in {0,1,2}
{
\begin{scope}[xshift=2.5*\a cm]

\draw
	(-0.5,-0.4) -- (-0.5,0.7) -- (0,1.1) -- (0.5,0.7) -- (0.5,-0.4)
	(-0.5,0.7) -- ++(-0.4,0)
	(0,1.1) -- ++(0,0.4);
		
\draw[line width=1.2]
	(-0.5,0) -- (0.5,0);

\node at (0,0.85) {\small $\alpha$}; 
\node at (-0.3,0.6) {\small $\beta$};
\node at (0.3,0.6) {\small $\gamma$};
\node at (-0.3,0.2) {\small $\delta$};
\node at (0.3,0.2) {\small $\epsilon$};	

\end{scope}
}


\draw[line width=1.2]
	(-0.5,0.7) -- ++(-0.4,0);

\node at (-0.7,0.5) {\small $\epsilon$};
\node at (-0.6,0.9) {\small $\delta$};
\node at (0.2,1.2) {\small $\beta$};
\node at (-0.2,1.2) {\small $\beta$};


\begin{scope}[xshift=2.5cm]
		
\node at (-0.7,0.5) {\small $\beta$};
\node at (-0.6,0.9) {\small $\alpha$};
\node at (0.2,1.2) {\small $\beta$};
\node at (-0.2,1.2) {\small $\beta$};

\node at (-0.3,-0.2) {\small $\epsilon$};
\node at (0.3,-0.2) {\small $\delta$};
\node at (0.7,0) {\small $\beta$};
\node at (-0.7,0) {\small $\beta$};

\end{scope}


\begin{scope}[xshift=5cm]
	
\draw[line width=1.2]
	(0,1.1) -- ++(0,0.4);
		
\node at (-0.7,0.5) {\small $\alpha$};
\node at (-0.6,0.9) {\small $\beta$};

\node at (0.2,1.2) {\small $\delta$};
\node at (-0.2,1.2) {\small $\delta$};

\end{scope}

\end{tikzpicture}
\caption{Proposition \ref{bde_all}: Special tile.}
\label{bde_allA}
\end{figure}
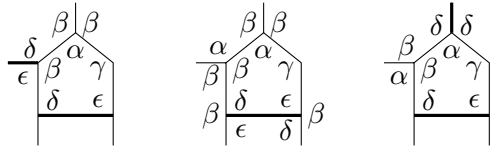

The second and third of Figure \ref{bde_allA} are the case $\beta\cdots=\alpha\beta^2$ in a special tile. By $\alpha\beta^2$ and $\alpha\ne\beta$, we know $\alpha^3$ is not a vertex. Then $\beta\cdots=\alpha\beta^2$ implies the degree $3$ vertex $\alpha\cdots=\alpha\beta^2,\alpha\delta^2$, with the angles at $\alpha\cdots,\beta\cdots$ arranged as indicated. In the second picture, the vertex $\gamma\cdots=\alpha\gamma\cdots,\gamma\delta\cdots$ has high degree. Then $\delta\cdots,\epsilon\cdots$ have degree $3$. As argued earlier, this implies $\delta\cdots=\epsilon\cdots=\beta\delta\epsilon$. Then we get two $\beta$ adjacent, a contradiction. 

In the third picture, we know $\beta\delta\epsilon,\alpha\beta^2,\alpha\delta^2$ are vertices. The angle sums of $\beta\delta\epsilon,\alpha\beta^2,\alpha\delta^2$ imply $\alpha=\epsilon$ and $\beta=\delta$. By Lemma \ref{geometry4}, this implies $a=b$, a contradiction. 
\end{proof}

\subsection{Two Matched Degree $3$ $b$-Vertices}
\label{2add}

We study tilings with two matched degree $3$ vertices, which are $\theta\delta^2,\theta\epsilon^2$, where $\theta=\alpha,\beta,\gamma$. By Lemma \ref{geometry2}, and up to the exchange symmetry $(\beta,\delta)\leftrightarrow(\gamma,\epsilon)$, we only need to consider three combinations: $\alpha\delta^2$ and $\beta\epsilon^2$; $\alpha\delta^2$ and $\gamma\epsilon^2$; $\beta\delta^2$ and $\gamma\epsilon^2$.

\begin{proposition}\label{a2d_b2e}
There is no tiling, such that $\alpha,\beta,\gamma$ have distinct values, and $\alpha\delta^2,\beta\epsilon^2$ are vertices.
\end{proposition}

\begin{proof}
The AAD $\thick^{\delta}\epsilon^{\gamma}\thin^{\alpha}\beta^{\delta}\thin^{\gamma}\epsilon^{\delta}\thick$ of $\beta\epsilon^2$ implies a vertex $\gamma\delta\cdots=\gamma\delta^2\cdots,\gamma\delta\epsilon\cdots$. If $\beta<\gamma$ and $\delta>\epsilon$, then $\gamma+2\delta>\gamma+\delta+\epsilon>\beta+2\epsilon=2\pi$, contradicting $\gamma\delta\cdots$. Therefore by Lemma \ref{geometry1}, we get $\beta>\gamma$ and $\delta<\epsilon$. Then by $\alpha\delta^2,\beta\epsilon^2$ and $\delta<\epsilon$, we get $\alpha>\beta>\gamma$. Moreover, by $\alpha\delta^2$ and $\delta<\epsilon$, we know a $b$-vertex $\alpha\cdots=\alpha\delta^2$. In particular, we know $\alpha\epsilon\cdots$ is not a vertex.

By Lemma \ref{geometry2} and $\alpha\ne\gamma$, we know $\alpha\delta^2,\beta\epsilon^2$ are all the degree $3$ $b$-vertices. This implies a special companion pair is matched, which determines $T_1,T_2$ in the first of Figure \ref{a2d_b2eA}. In the special companion pair, we know one of shared vertices have degree $3$.

If the shared vertex $\epsilon_1\epsilon_2\cdots$ has high degree, then by Lemma \ref{bb_pair}, the shared vertex $\delta_1\delta_2\cdots$ has degree $3$, and is $\alpha\delta^2$. Up to the vertical flip, we may assume $T_3$ is arranged as indicated. Since $\epsilon_1\epsilon_2\cdots$ has high degree, we know one of $\beta_1\beta_3\cdots,\beta_2\gamma_3\cdots$ has degree $3$. 

If the shared vertex $\delta_1\delta_2\cdots$ has high degree, then the shared vertex $\epsilon_1\epsilon_2\cdots$ has degree $3$, and is $\beta\epsilon^2$. Up to the vertical flip, we may assume $T_4$ is arranged as indicated. Since $\delta_1\delta_2\cdots$ has high degree, we know one of $\alpha_4\gamma_1\cdots,\gamma_2\delta_4\cdots$ has degree $3$. Since we already know $\gamma_2\delta_4\cdots$ has high degree, we know $\alpha_4\gamma_1\cdots$ has degree $3$.

If both shared vertices have degree $3$, then they are $\alpha\delta^2$ and $\beta\epsilon^2$. Then $\beta_1\cdots,\beta_2\cdots$ are $\beta^2\cdots,\beta\gamma\cdots$, and $\gamma_1\cdots,\gamma_2\cdots$ are $\alpha\gamma\cdots,\gamma\delta\cdots$. Since $\gamma\delta\cdots$ has high degree, we know one of $\beta^2\cdots,\beta\gamma\cdots,\alpha\gamma\cdots$ has degree $3$. 

\begin{figure}[htp]
\centering
\begin{tikzpicture}[>=latex,scale=1]


\draw
	(-0.5,-0.7) -- (-0.5,0.7) -- (0,1.1) -- (0.5,0.7) -- (0.5,-0.7) -- (0,-1.1) -- (-0.5,-0.7)
	(0.5,-0.7) -- (1.3,-0.7) -- (1.3,0.7) -- (0.5,0.7)
	(-0.5,-0.7) -- (-1.3,-0.7) -- (-1.3,0.7) -- (-0.5,0.7);

\draw[line width=1.2]
	(-0.5,0) -- (0.5,0)
	(-1.3,-0.7) -- (-1.3,0.7)
	(0.5,-0.7) -- (1.3,-0.7);

\node at (0.3,-0.6) {\small $\gamma$};
\node at (0,-0.85) {\small $\alpha$};
\node at (0.3,-0.2) {\small $\epsilon$};	
\node at (-0.3,-0.2) {\small $\delta$}; 
\node at (-0.3,-0.6) {\small $\beta$};

\node at (0.3,0.6) {\small $\gamma$};
\node at (0,0.85) {\small $\alpha$};
\node at (0.3,0.2) {\small $\epsilon$};	
\node at (-0.3,0.2) {\small $\delta$}; 
\node at (-0.3,0.6) {\small $\beta$};

\node at (1.1,0.5) {\small $\gamma$};
\node at (1.1,-0.5) {\small $\epsilon$};
\node at (0.7,0.5) {\small $\alpha$};
\node at (0.7,-0.5) {\small $\delta$};
\node at (0.7,0) {\small $\beta$};

\node at (-1.1,0.5) {\small $\delta$};
\node at (-1.1,-0.5) {\small $\epsilon$};
\node at (-0.7,0.5) {\small $\beta$};
\node at (-0.7,-0.5) {\small $\gamma$};
\node at (-0.7,0) {\small $\alpha$};

\node[inner sep=0.5,draw,shape=circle] at (0,0.4) {\small $1$};
\node[inner sep=0.5,draw,shape=circle] at (0,-0.4) {\small $2$};
\node[inner sep=0.5,draw,shape=circle] at (-1.05,0) {\small $3$};
\node[inner sep=0.5,draw,shape=circle] at (1.05,0) {\small $4$};


\begin{scope}[shift={(2.8cm,-0.5cm)}]

\draw
	(-0.5,-0.4) -- (-0.5,0.7) -- (0,1.1) -- (0.5,0.7) -- (0.5,-0.4)
	(0.5,0.7) -- ++(0.4,0)
	(0,1.1) -- ++(0,0.4);
		
\draw[line width=1.2]
	(-0.5,0) -- (0.5,0);

\node at (0,0.85) {\small $\alpha$}; 
\node at (-0.3,0.6) {\small $\beta$};
\node at (0.3,0.6) {\small $\gamma$};
\node at (-0.3,0.2) {\small $\delta$};
\node at (0.3,0.2) {\small $\epsilon$};

\node at (0.7,0) {\small $\beta$};
\node at (-0.3,-0.2) {\small $\delta$};
\node at (0.3,-0.2) {\small $\epsilon$};
\node at (0.2,1.2) {\small $\alpha$};
\node at (0.7,0.5) {\small $\alpha$};
\node at (0.6,0.9) {\small $\gamma$};

\end{scope}


\begin{scope}[xshift=5cm]

\foreach \a in {0,1,2,3}
\draw[xshift=\a cm]
	(0,0.7) -- (0,0) -- (-0.5,-0.4) -- (-0.5,-1.1);

\foreach \a in {0,1,2}
\draw[xshift=\a cm]
	(0,0) -- (0.5,-0.4);

\draw
	(0,0.7) -- (3,0.7)
	(-0.5,-1.1) -- (2.5,-1.1)
	(0,0.7) -- (0,1.4) -- (2,1.4) -- (2,0.7);

\draw[line width=1.2]
	(0.5,-0.4) -- (0,0)
	(1,0.7) -- (3,0.7)
	(0.5,-1.1) -- (2.5,-1.1);
		
\node at (0.35,-0.95) {\small $\gamma$}; 
\node at (0.35,-0.5) {\small $\epsilon$};	
\node at (0,-0.25) {\small $\delta$};
\node at (-0.35,-0.5) {\small $\beta$};
\node at (-0.35,-0.95) {\small $\alpha$};

\node at (0.5,-0.2) {\small $\epsilon$}; 
\node at (0.15,0.15) {\small $\delta$};	
\node at (0.15,0.5) {\small $\beta$};
\node at (0.85,0.5) {\small $\alpha$};	 
\node at (0.85,0.1) {\small $\gamma$};

\foreach \a in {0,1}
{
\begin{scope}[xshift=\a cm]

\node at (0.65,-0.9) {\small $\delta$}; 
\node at (1.35,-0.95) {\small $\epsilon$}; 		
\node at (1.35,-0.5) {\small $\gamma$};	
\node at (1,-0.25) {\small $\alpha$};
\node at (0.65,-0.5) {\small $\beta$};

\node at (1.5,-0.2) {\small $\alpha$}; 
\node at (1.15,0.15) {\small $\beta$};	
\node at (1.15,0.5) {\small $\delta$};
\node at (1.85,0.5) {\small $\epsilon$};	 
\node at (1.85,0.1) {\small $\gamma$};

\end{scope}
}

\node at (0.15,1.25) {\small $\alpha$};  
\node at (0.15,0.9) {\small $\beta$}; 	 	
\node at (1,0.9) {\small $\delta$};
\node at (1.85,0.9) {\small $\epsilon$};
\node at (1.85,1.25) {\small $\gamma$};

\node[inner sep=0.5,draw,shape=circle] at (0,-0.7) {\small $1$};
\node[inner sep=0.5,draw,shape=circle] at (1,-0.7) {\small $2$};
\node[inner sep=0.5,draw,shape=circle] at (0.5,0.3) {\small $3$};
\node[inner sep=0.5,draw,shape=circle] at (1.5,0.3) {\small $4$};
\node[inner sep=0.5,draw,shape=circle] at (1.5,1.05) {\small $5$};
\node[inner sep=0.5,draw,shape=circle] at (2,-0.7) {\small $6$};
\node[inner sep=0.5,draw,shape=circle] at (2.5,0.3) {\small $7$};

\end{scope}

\end{tikzpicture}
\caption{Proposition \ref{a2d_b2e}: Special companion pair, $\alpha\gamma^2$, $\alpha\beta\gamma$.}
\label{a2d_b2eA}
\end{figure}
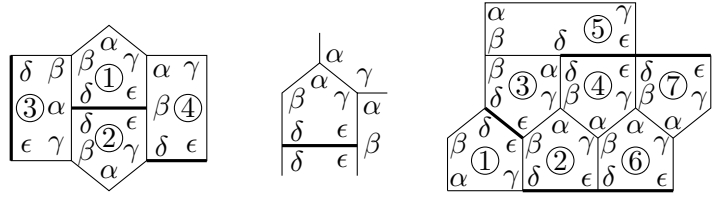

Combining all the cases, we conclude one of $\alpha\gamma\cdots,\beta^2\cdots,\beta\gamma\cdots$ has degree $3$. This means one of $\alpha\beta^2,\alpha\beta\gamma,\alpha^2\gamma,\alpha\gamma^2,\beta^3,\beta^2\gamma,\beta\gamma^2$ is a vertex.

\subsubsection*{Case. $\alpha^2\gamma$ is a vertex}

The angle sums of $\alpha\delta^2,\beta\epsilon^2,\alpha^2\gamma$ and the angle sum for pentagon imply
\[
\alpha=\pi-\tfrac{1}{2}\gamma,\;
\beta=(1+\tfrac{8}{f})\pi-\tfrac{3}{2}\gamma,\;
\delta=\tfrac{1}{2}\pi+\tfrac{1}{4}\gamma,\;
\epsilon=(\tfrac{1}{2}-\tfrac{4}{f})\pi+\tfrac{3}{4}\gamma.
\]
By $\alpha^2\gamma$ and $\alpha\ne\beta$, we know $\alpha\beta\gamma$ is not a vertex.

By $\epsilon>\delta>\frac{1}{2}\pi$, we know $\delta\thin\delta\cdots,\delta\thin\epsilon\cdots, \epsilon\thin\epsilon\cdots$ are not vertices. Then the AAD of $\beta^k\gamma^l$ is $\thin^{\alpha}\theta^{\rho}\thin^{\alpha}\theta^{\rho}\thin\cdots$, where $\theta^{\rho}=\beta^{\delta},\gamma^{\epsilon}$. Then by no $\alpha\epsilon\cdots$, we get $\beta^k\gamma^l=\beta^k$.

The AAD $\thick^{\epsilon}\delta^{\beta}\thin^{\beta}\alpha^{\gamma}\thin^{\beta}\delta^{\epsilon}\thick$ of $\alpha\delta^2$ implies a vertex $\beta\gamma\cdots$. By $\beta+\gamma+2\delta=(2+\tfrac{8}{f})\pi>2\pi$ and $\delta<\epsilon$, we know $\beta\gamma\cdots$ has no $\delta,\epsilon$. Then by $\beta^k\gamma^l=\beta^k$, we get $\beta\gamma\cdots=\alpha\beta\gamma\cdots$. Then by $\alpha>\beta>\gamma$, and $\alpha+\beta+2\gamma=(2+\tfrac{8}{f})\pi>2\pi$, and no $\alpha\beta\gamma$, we know $\alpha\beta\gamma\cdots$ is not a vertex.

\subsubsection*{Case. $\alpha\gamma^2$ is a vertex}

By distinct $\alpha,\beta,\gamma$ values, we know $\alpha\delta^2,\beta\epsilon^2,\alpha\gamma^2,\alpha^2\beta,\beta^3,\beta^2\gamma$ are all the degree $3$ vertices. By $\alpha>\beta>\gamma$ and $\alpha\gamma^2$, we know $\alpha^2\beta$ is not a vertex. The angle sums of $\alpha\delta^2,\beta\epsilon^2,\alpha\gamma^2,\beta^3$ and the angle sum for pentagon imply $f=12$, a contradiction. The angle sums of $\alpha\delta^2,\beta\epsilon^2,\alpha\gamma^2,\beta^2\gamma$ and the angle sum for pentagon imply $\alpha=(\tfrac{32}{f}-2)\pi$, contradicting $f\ge 16$. 

Therefore $\alpha\delta^2,\beta\epsilon^2,\alpha\gamma^2$ are all the degree $3$ vertices. This implies one of $\beta\cdots,\delta\cdots$ in a special tile has high degree. Then the vertices $\alpha\cdots,\gamma\cdots,\epsilon\cdots$ have degree $3$. This implies $\epsilon\cdots=\beta\epsilon^2$, and $\gamma\cdots=\alpha\gamma^2$ as indicated in the second of Figure \ref{a2d_b2eA}. Then by no $\alpha\epsilon\cdots$, we get $\alpha\cdots=\alpha^2\cdots$. This has high degree, a contradiction.

\subsubsection*{Case. $\alpha\beta\gamma$ is a vertex}

The angle sums of $\alpha\delta^2,\beta\epsilon^2,\alpha\beta\gamma$ and the angle sum for pentagon imply
\[
\alpha+\beta=(2-\tfrac{8}{f})\pi,\;
\gamma=\tfrac{8}{f}\pi,\;
\delta+\epsilon=(1+\tfrac{4}{f})\pi.
\]

By $\alpha\beta\gamma$ and $\alpha>\beta$, we get $R(\alpha^2)<R(\alpha\beta)=\gamma<\alpha,\beta$. Then by $b$-vertex $\alpha\cdots=\alpha\delta^2$, this implies $\alpha\beta\cdots=\alpha\beta\gamma$, and $\alpha^2\cdots$ is not a vertex.

By $\delta<\epsilon$ and $\delta+\epsilon>\pi$, we know $R(\epsilon^2)$ has no $\delta,\epsilon$. Then by $\beta\epsilon^2$ and $\alpha>\beta>\gamma$, we get $\epsilon^2\cdots=\beta\epsilon^2,\gamma^k\epsilon^2(k\ge 2)$. This implies no $\epsilon\thin\epsilon\cdots$. 

By no $\alpha^2\cdots,\alpha\epsilon\cdots,\epsilon\thin\epsilon\cdots$, the AAD implies no $\gamma\thin\gamma\cdots$. This implies $\gamma^k\epsilon^2$ is not a vertex. Therefore $\epsilon^2\cdots=\beta\epsilon^2$.

By no $\alpha^2\cdots,\gamma\thin\gamma\cdots$, we know $\alpha^k\gamma^l$ is not a vertex. Then by $b$-vertex $\alpha\cdots=\alpha\delta^2$, and $\alpha\beta\cdots=\alpha\beta\gamma$, we get $\alpha\cdots=\alpha\delta^2,\alpha\beta\gamma$.

The AAD $\thick^{\delta}\epsilon^{\gamma}\thin^{\alpha}\beta^{\delta}\thin^{\gamma}\epsilon^{\delta}\thick$ of $\beta\epsilon^2$ determines $T_1,T_2,T_3$ in the third of Figure \ref{a2d_b2eA}. Then $\alpha_2\gamma_3\cdots=\alpha\beta\gamma$ and no $\alpha^2\cdots$ determine $T_4$. Then $\alpha_3\delta_4\cdots=\alpha\delta^2$ determines $T_5$, and $\alpha_4\gamma_2\cdots=\alpha\beta\gamma$ and no $\alpha\epsilon\cdots$ determine $T_6$. Then $\alpha_6\gamma_4\cdots=\alpha\beta\gamma$ and no $\alpha\epsilon\cdots$ determine $T_7$. Then $\delta_7\epsilon_4\epsilon_5\cdots$ contradicts  $\epsilon^2\cdots=\beta\epsilon^2$.

\subsubsection*{Case. $\alpha\beta^2$ is a vertex}

By distinct $\alpha,\beta,\gamma$ values, we know $\alpha\delta^2,\beta\epsilon^2,\alpha\beta^2,\alpha^2\gamma,\beta\gamma^2,\gamma^3$ are all the degree $3$ vertices. We already proved there is no tiling with $\alpha^2\gamma$. By $\alpha>\beta>\gamma$ and $\alpha\beta^2$, we know $\beta\gamma^2,\gamma^3$ are not vertices. Therefore $\alpha\delta^2,\beta\epsilon^2,\alpha\beta^2$ are all the degree $3$ vertices. By Lemma \ref{ndegree3}, we know one of $\alpha\gamma^3,\beta\gamma^3,\gamma^4,\gamma^5$ is a vertex.

By $\alpha\delta^2,\beta\epsilon^2$, we get $\alpha+\beta+2(\delta+\epsilon)=4\pi$. By $\alpha\beta^2$, we know $\alpha+\beta<2\pi$. Therefore $\delta+\epsilon>\pi$. By $\delta<\epsilon$, this implies $\epsilon\thin\epsilon\cdots$ is not a vertex. By no $\alpha\epsilon\cdots,\epsilon\thin\epsilon\cdots$, we know the AAD of $\thin\gamma\thin\gamma\thin$ is $\thin^{\epsilon}\gamma^{\alpha}\thin^{\alpha}\gamma^{\epsilon}\thin$. This implies no consecutive $\gamma\gamma\gamma$. Therefore $\alpha\gamma^3,\beta\gamma^3,\gamma^4,\gamma^5$ are not vertices.

\subsubsection*{Case. $\beta^3$ is a vertex}

By distinct $\alpha,\beta,\gamma$ values, we know $\alpha\delta^2,\beta\epsilon^2,\beta^3,\alpha\beta\gamma,\alpha^2\gamma,\alpha\gamma^2$ are all the degree $3$ vertices. We already proved there is no tiling with any one of $\alpha\beta\gamma,\alpha^2\gamma,\alpha\gamma^2$. Therefore $\alpha\delta^2,\beta\epsilon^2,\beta^3$ are all the degree $3$ vertices. Therefore $\gamma\cdots$ has high degree. By Lemma \ref{special_tile}, we know $f\ge 24$. By Lemma \ref{ndegree3}, we also know one of $\alpha\gamma^3,\beta\gamma^3,\gamma^4,\gamma^5$ is a vertex. The angle sums of $\alpha\delta^2,\beta\epsilon^2,\beta^3$, and the angle sum of one of these, and the angle sum for pentagon imply 
\begin{align*}
\alpha\gamma^3 &\colon
	\alpha=\tfrac{24}{f}\pi,\,
	\beta=\epsilon=\tfrac{2}{3}\pi,\,
	\gamma=(\tfrac{2}{3}-\tfrac{8}{f})\pi,\,
	\delta=(1-\tfrac{12}{f})\pi. \\
\beta\gamma^3 &\colon
	\alpha=(\tfrac{4}{9}+\tfrac{8}{f})\pi,\,
	\beta=\epsilon=\tfrac{2}{3}\pi,\,
	\gamma=\tfrac{4}{9}\pi,\,
	\delta=(\tfrac{7}{9}-\tfrac{4}{f})\pi. \\
\gamma^4 &\colon
	\alpha=(\tfrac{1}{3}+\tfrac{8}{f})\pi,\,
	\beta=\epsilon=\tfrac{2}{3}\pi,\,
	\gamma=\tfrac{1}{2}\pi,\,
	\delta=(\tfrac{5}{6}-\tfrac{4}{f})\pi. \\
\gamma^5 &\colon
	\alpha=(\tfrac{8}{15}+\tfrac{8}{f})\pi,\,
	\beta=\epsilon=\tfrac{2}{3}\pi,\,
	\gamma=\tfrac{2}{5}\pi,\,
	\delta=(\tfrac{11}{15}-\tfrac{4}{f})\pi.
\end{align*}
By $f\ge 24$, we get $\delta+\epsilon>\pi$. Then similar to the case $\alpha\beta^2$ is a vertex, there is no consecutive $\gamma\gamma\gamma$. Therefore $\alpha\gamma^3,\beta\gamma^3,\gamma^4,\gamma^5$ are not vertices.

\subsubsection*{Case. $\beta\gamma^2$ is a vertex}

By distinct $\alpha,\beta,\gamma$ values, we know $\alpha\delta^2,\beta\epsilon^2,\beta\gamma^2,\alpha^3,\alpha\beta^2,\alpha^2\gamma$ are all the degree $3$ vertices. By $\alpha>\beta>\gamma$ and $\beta\gamma^2$, we know $\alpha^3,\alpha\beta^2$ are not vertices. We already proved there is no tiling with $\alpha^2\gamma$. Therefore $\alpha\delta^2,\beta\epsilon^2,\beta\gamma^2$ are all the degree $3$ vertices. This implies that, in any tile, one of $\alpha\cdots,\beta\cdots$ has high degree, and one of $\gamma\cdots,\epsilon\cdots$ has high degree. This implies no special tile, a contradiction.

\subsubsection*{Case. $\beta^2\gamma$ is a vertex}

By $\alpha\delta^2,\beta^2\gamma$, we get $\beta<\pi$, and $\beta+\delta+(\alpha+\beta+\gamma+\delta)=4\pi$. Then by the angle sum for pentagon, we get $\beta+\delta=(1-\frac{4}{f})\pi+\epsilon$. 

By $\beta\epsilon^2$ and $\beta<\pi$, we get $\epsilon>\frac{1}{2}\pi>\frac{4}{f}\pi$. Therefore $\alpha+\delta>\beta+\delta>\pi$. By $\delta<\epsilon$, this implies $\alpha^2\cdots,\beta^2\cdots$ have no $\delta,\epsilon$. Then by $\beta^2\gamma$ and $\alpha>\beta>\gamma$, we know $\alpha^2\cdots$ is not a vertex, and $\beta^2\cdots=\beta^2\gamma$.

The AAD $\thick^{\epsilon}\delta^{\beta}\thin^{\beta}\alpha^{\gamma}\thin^{\beta}\delta^{\epsilon}\thick$ of $\alpha\delta^2$ implies a vertex $\thin^{\alpha}\beta^{\delta}\thin^{\alpha}\beta^{\delta}\thin\cdots=\beta^2\gamma=\thin^{\alpha}\beta^{\delta}\thin\gamma\thin^{\alpha}\beta^{\delta}\thin$. This implies $\alpha^2\cdots,\alpha\epsilon\cdots$, a contradiction.  
\end{proof}

\begin{proposition}\label{a2d_c2e}
There is no tiling, such that $\alpha,\beta,\gamma$ have distinct values, and $\alpha\delta^2,\gamma\epsilon^2$ are vertices.
\end{proposition}

\begin{proof}
By Proposition \ref{bde_a2d_c3e}, we know $\beta\delta\epsilon$ is not a vertex. Then by distinct $\alpha,\beta,\gamma$ values, we know $\alpha\delta^2,\gamma\epsilon^2$ are all the degree $3$ $b$-vertices. This implies a special companion pair is matched, and a $b$-vertex $\beta\cdots$ has high degree. 

A special companion pair determines $T_1,T_2$ in Figure \ref{a2d_c2eA}. By Lemma \ref{bb_pair}, one of the shared vertices $\delta_1\delta_2\cdots, \epsilon_1\epsilon_2\cdots$ has degree $3$. If $\delta_1\delta_2\cdots$ has degree $3$, then it is $\alpha\delta^2$, and $\beta_1\cdots,\beta_2\cdots$ are $\beta^2\cdots,\beta\gamma\cdots$. If $\epsilon_1\epsilon_2\cdots$ has degree $3$, then it is $\gamma\epsilon^2$, and up to the vertical flip, we may assume $T_4$ is arranged as indicated. Then we get $\alpha_4\gamma_1\cdots,\gamma_2\epsilon_4\cdots$. Moreover, if $\gamma_2\epsilon_4\cdots$ has degree $3$, then it is $\gamma\epsilon^2$, and we get $\alpha_2\cdots=\alpha\gamma\cdots$. Then by the similar argument as in Proposition \ref{a2d_b2e}, we conclude one of $\beta^2\cdots,\beta\gamma\cdots,\alpha\gamma\cdots$ has degree $3$. This means one of $\alpha\beta^2,\alpha\beta\gamma,\alpha^2\gamma,\alpha\gamma^2,\beta^3,\beta^2\gamma,\beta\gamma^2$ is a vertex.

\begin{figure}[htp]
\centering
\begin{tikzpicture}[>=latex,scale=1]

\foreach \a in {1,-1}
\draw[xscale=\a]
	(0.5,0.7) -- (1.3,1) -- (1.3,1.8) -- (0,1.8) -- (0,1.1) -- (0.5,0.7) -- (0.5,-0.7) -- (0,-1.1) -- (0.3,-1.8) -- (1.3,-1.8) -- (1.3,-0.7)
	(0.5,0.7) -- (1.3,0.7) -- (1.3,-0.7) -- (0.5,-0.7);

\draw[line width=1.2]
	(-0.5,0) -- (0.5,0)
	(-1.3,-0.7) -- (-1.3,0.7)
	(0.5,-0.7) -- (1.3,-0.7)
	(0,1.8) -- (0,1.1)
	(0,-1.1) -- (-0.3,-1.8);

\node at (0.3,-0.6) {\small $\gamma$};
\node at (0,-0.85) {\small $\alpha$};
\node at (0.3,-0.2) {\small $\epsilon$};	
\node at (-0.3,-0.2) {\small $\delta$}; 
\node at (-0.3,-0.6) {\small $\beta$};

\node at (0.3,0.6) {\small $\gamma$};
\node at (0,0.85) {\small $\alpha$};
\node at (0.3,0.2) {\small $\epsilon$};	
\node at (-0.3,0.2) {\small $\delta$}; 
\node at (-0.3,0.6) {\small $\beta$};

\node at (1.1,0.5) {\small $\beta$};
\node at (1.1,-0.5) {\small $\delta$};
\node at (0.7,0.5) {\small $\alpha$};
\node at (0.7,-0.5) {\small $\epsilon$};
\node at (0.7,0) {\small $\gamma$};

\node at (-1.1,0.5) {\small $\delta$};
\node at (-1.1,-0.5) {\small $\epsilon$};
\node at (-0.7,0.5) {\small $\beta$};
\node at (-0.7,-0.5) {\small $\gamma$};
\node at (-0.7,0) {\small $\alpha$};

\node at (1.1,1.1) {\small $\alpha$};
\node at (0.6,0.9) {\small $\beta$};
\node at (1.1,1.6) {\small $\gamma$};
\node at (0.2,1.2) {\small $\delta$};
\node at (0.2,1.6) {\small $\epsilon$};

\node at (-1.1,1.1) {\small $\alpha$};
\node at (-0.5,0.9) {\small $\beta$};
\node at (-1.1,1.6) {\small $\gamma$};
\node at (-0.2,1.2) {\small $\delta$};
\node at (-0.2,1.6) {\small $\epsilon$};

\node at (1.1,-0.9) {\small $\delta$};
\node at (0.5,-0.9) {\small $\epsilon$};
\node at (1.1,-1.6) {\small $\beta$};
\node at (0.2,-1.2) {\small $\gamma$};
\node at (0.4,-1.6) {\small $\alpha$};

\node at (-1.1,-0.9) {\small $\alpha$};
\node at (-0.6,-0.9) {\small $\gamma$};
\node at (-1.1,-1.6) {\small $\beta$};
\node at (-0.2,-1.2) {\small $\epsilon$};
\node at (-0.4,-1.6) {\small $\delta$};

\node[inner sep=0.5,draw,shape=circle] at (0,0.4) {\small $1$};
\node[inner sep=0.5,draw,shape=circle] at (0,-0.4) {\small $2$};
\node[inner sep=0.5,draw,shape=circle] at (-1.05,0) {\small $3$};
\node[inner sep=0.5,draw,shape=circle] at (1.05,0) {\small $4$};
\node[inner sep=0.5,draw,shape=circle] at (-0.7,1.3) {\small $5$};
\node[inner sep=0.5,draw,shape=circle] at (0.7,1.3) {\small $6$};
\node[inner sep=0.5,draw,shape=circle] at (-0.7,-1.3) {\small $7$};
\node[inner sep=0.5,draw,shape=circle] at (0.7,-1.3) {\small $8$};

\end{tikzpicture}
\caption{Proposition \ref{a2d_c2e}: Special companion pair.}
\label{a2d_c2eA}
\end{figure}
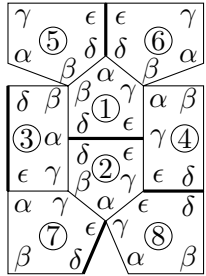

\subsubsection*{Case. $\alpha\beta^2$ is a vertex}

By $\alpha\delta^2,\alpha\beta^2$, we get $\beta=\delta$. Then by Lemma \ref{geometry1}, we know $\beta=\delta$ is strictly between $\gamma,\epsilon$. 

By distinct $\alpha,\beta,\gamma$ values, and all the degree $3$ $b$-vertices $\alpha\delta^2,\gamma\epsilon^2$, we know $\alpha\delta^2,\gamma\epsilon^2,\alpha\beta^2,\alpha^2\gamma,\beta\gamma^2,\gamma^3$ are all the degree $3$ vertices. If $\alpha^2\gamma$ is a vertex, then by $\gamma\epsilon^2$, we get $\alpha=\epsilon$. Then by $\beta=\delta$ and $\alpha=\epsilon$, and Lemma \ref{geometry11}, we get $a=b$, a contradiction. If $\beta\gamma^2$ is a vertex, then by $\gamma\epsilon^2$, we get $\beta+\gamma=2\epsilon$, contradicting $\beta$ being strictly between $\gamma,\epsilon$. If $\gamma^3$ is a vertex, then by $\gamma\epsilon^2$, we get $\gamma=\epsilon$, contradicting $\beta=\delta$ being strictly between $\gamma,\epsilon$.

Therefore $\alpha\delta^2,\gamma\epsilon^2,\alpha\beta^2$ are all the degree $3$ vertices. Consider the special companion pair in Figure \ref{a2d_c2eA}. We know one of $\alpha_1\cdots,\gamma_1\cdots$ has high degree, and one of $\alpha_2\cdots,\gamma_2\cdots$ has high degree. This implies $\delta_1\delta_2\cdots$ has degree $3$. Then $\delta_1\delta_2\cdots=\alpha\delta^2$. This implies one of one of $\beta_1\cdots,\beta_2\cdots$ is $\beta\gamma\cdots$, which has high degree. We find the special companion pair has at least three high degree vertices. Therefore the special companion pair has weight $\ge 3w_4$. By Lemma \ref{special_pair}, this implies $f\ge 48$.

The angle sums of $\alpha\delta^2,\gamma\epsilon^2,\alpha\beta^2$ and the angle sum for pentagon imply
\[
\alpha+2\beta=2\pi,\;
\beta=\delta,\,
\gamma=\tfrac{8}{f}\pi,\;
\epsilon=(1-\tfrac{4}{f})\pi.
\]
We have $\gamma<\epsilon$. By $\beta=\delta$ is strictly between $\gamma,\epsilon$, this implies $\gamma<\beta=\delta<\epsilon$. By $\alpha\delta^2,\gamma\epsilon^2$, this further implies $\alpha>\gamma$. 

By $\gamma\epsilon^2$, and $\gamma<\alpha,\beta,\delta,\epsilon$, we get $\epsilon^2\cdots=\gamma\epsilon^2$. Therefore $\epsilon\thin\epsilon\cdots$ is not a vertex. 

By $\alpha\delta^2$ and $\delta<\epsilon$, a $b$-vertex $\alpha\cdots=\alpha\delta^2$, and $\alpha\epsilon\cdots$ is not a vertex. By no $\alpha\epsilon\cdots,\epsilon\thin\epsilon\cdots$, we know the AAD of $\thin\gamma\thin\gamma\thin$ is $\thin^{\epsilon}\gamma^{\alpha}\thin^{\alpha}\gamma^{\epsilon}\thin$. This implies no consecutive $\gamma\gamma\gamma$.

We have $\gamma<\beta=\delta<\epsilon<\pi$. If $\alpha^2\cdots$ is a vertex, then $\alpha<\pi$, and the pentagon is convex. By Lemmas \ref{geometry6} and \ref{geometry10}, we get $\alpha+2\gamma>\pi$ and $2\alpha+\beta+\gamma>2\pi$. Then by $f\ge 48$, this implies $\alpha>(1-\frac{16}{f})\pi\ge\frac{2}{3}\pi$. Then we get $R(\alpha^2)<\alpha,\beta+\gamma,4\gamma$. By all the degree $3$ vertices $\alpha\delta^2,\gamma\epsilon^2,\alpha\beta^2$, and $b$-vertex $\alpha\cdots=\alpha\delta^2$, we also know $\alpha^2\cdots$ is a high degree $\hat{b}$-vertex. Then by $\beta>\gamma$, we get $\alpha^2\cdots=\alpha^2\gamma^2,\alpha^2\gamma^3$. 

Suppose $\alpha^2\cdots=\alpha^2\gamma^2$. The angle sum of $\alpha^2\gamma^2$ further implies
\[
	\alpha=(1-\tfrac{8}{f})\pi,\,
	\gamma=\tfrac{8}{f}\pi,\,
	\beta=\delta=(\tfrac{1}{2}+\tfrac{4}{f})\pi,\,
	\epsilon=(1-\tfrac{4}{f})\pi.
\]
By $\epsilon>\delta>\frac{1}{2}\pi$, there is at most one $b$-edge at a vertex, and $\delta\thin\delta\cdots,\delta\thin\epsilon\cdots,\epsilon\thin\epsilon\cdots$ are not vertices. 

By $\beta+\delta+\epsilon=(2+\tfrac{4}{f})\pi>2\pi$, and $\delta<\epsilon$, and at most one $b$-edge at a vertex, we know a $b$-vertex $\beta\cdots=\beta\delta^2\cdots$, with no $\delta,\epsilon$ in the remainder. Then by $\alpha\delta^2$ and $\beta+\delta=(1+\tfrac{8}{f})\pi>\pi$, we know $R(\beta\delta^2)$ has no $\alpha,\beta$. Therefore a $b$-vertex $\beta\cdots=\beta\gamma^k\delta^2$, and $\beta\epsilon\cdots$ is not a vertex. By $\alpha\delta^2$ and $\alpha\ne\beta$, we get $k\ge 1$ in $\beta\gamma^k\delta^2$. 

We know the AAD of $\thin\gamma\thin\gamma\thin$ is $\thin^{\epsilon}\gamma^{\alpha}\thin^{\alpha}\gamma^{\epsilon}\thin$. The AAD implies a vertex $\thin^{\beta}\alpha^{\gamma}\thin^{\gamma}\alpha^{\beta}\thin\cdots=\alpha^2\gamma^2=\thin^{\beta}\alpha^{\gamma}\thin^{\gamma}\alpha^{\beta}\thin^{\epsilon}\gamma^{\alpha}\thin^{\alpha}\gamma^{\epsilon}\thin$, contradicting no $\beta\epsilon\cdots$. Therefore $\gamma\thin\gamma\cdots$ is not a vertex. This implies $\beta\gamma^k\delta^2=\thick\delta\thin\beta\thin\gamma\thin\delta\thick,\thick\delta\thin\gamma\thin\beta\thin\gamma\thin\delta\thick$. Then the AAD of the vertex imply $\alpha\epsilon\cdots,\beta\epsilon\cdots,\delta\thin\epsilon\cdots$, a contradiction. Therefore $\beta\delta\cdots$ is not a vertex. 

By no $\gamma\thin\gamma\cdots$, we know $\alpha^2\gamma^2=\thin\alpha\thin\gamma\thin\alpha\thin\gamma\thin$. Then by no $\beta\epsilon\cdots$, the AAD of the vertex is $\thin^{\beta}\alpha^{\gamma}\thin^{\epsilon}\gamma^{\alpha}\thin^{\beta}\alpha^{\gamma}\thin^{\epsilon}\gamma^{\alpha}\thin$. This implies a vertex $\thin^{\beta}\alpha^{\gamma}\thin^{\alpha}\beta^{\delta}\thin\cdots$. By $\alpha\beta^2$, and $\alpha^2\cdots=\alpha^2\gamma^2$, and no $\beta\delta\cdots,\beta\epsilon\cdots$, we get $\thin^{\beta}\alpha^{\gamma}\thin^{\alpha}\beta^{\delta}\thin\cdots=\alpha\beta^2,\alpha\beta\gamma^k$. Then by no $\delta\thin\delta\cdots,\delta\thin\epsilon\cdots,\epsilon\thin\epsilon\cdots$, the AAD of the vertex is $\thin^{\beta}\alpha^{\gamma}\thin^{\alpha}\beta^{\delta}\thin^{\alpha}\beta^{\delta}\thin,\thin^{\beta}\alpha^{\gamma}\thin^{\alpha}\beta^{\delta}\thin^{\alpha}\gamma^{\epsilon}\thin
\cdots\thin^{\alpha}\gamma^{\epsilon}\thin$, contradicting no $\beta\delta\cdots,\beta\epsilon\cdots$.

We conclude $\alpha^2\gamma^2$ is not a vertex. This implies $\alpha^2\cdots=\alpha^2\gamma^3$. By no consecutive $\gamma\gamma\gamma$, this implies $\alpha\thin\alpha\cdots$ is not a vertex. Then the AAD $\thin^{\epsilon}\gamma^{\alpha}\thin^{\alpha}\gamma^{\epsilon}\thin$ of $\thin\gamma\thin\gamma\thin$ implies $\gamma\thin\gamma\cdots$ is not a vertex. This further implies $\alpha^2\gamma^3$ is not a vertex. 

We conclude $\alpha^2\cdots$ is not a vertex. By the AAD $\thin^{\epsilon}\gamma^{\alpha}\thin^{\alpha}\gamma^{\epsilon}\thin$ of $\thin\gamma\thin\gamma\thin$, we know $\gamma\thin\gamma\cdots$ is not a vertex. 

By $\alpha\beta^2$, and $b$-vertex $\alpha\cdots=\alpha\delta^2$, and no $\alpha^2\cdots$, we get $\alpha\cdots=\alpha\delta^2,\alpha\beta^2,\alpha\beta\gamma^k,\alpha\gamma^k$. By no $\gamma\thin\gamma\cdots$, we know $\alpha\gamma^k$ is not a vertex. By no $\alpha^2\cdots,\alpha\epsilon\cdots$, we do not have $\gamma\thin\beta\thin\gamma\cdots$. Then by no $\gamma\thin\gamma\cdots$, we get $k=1$ in $\alpha\beta\gamma^k$, contradicting $\alpha\beta^2$ and $\beta\ne\gamma$. Therefore $\alpha\cdots=\alpha\delta^2,\alpha\beta^2$. This implies
\[
f=\#\alpha=\#\alpha\delta^2+\#\alpha\beta^2
\le \tfrac{1}{2}\#\delta+\tfrac{1}{2}\#\beta=f.
\]
Since both sides are equal, this implies $\alpha\delta^2,\alpha\beta^2,\gamma^k\epsilon^l$ are all the vertices. By $\gamma\epsilon^2$, and $\gamma<\epsilon$, and no $\gamma\thin\gamma\cdots$, we get $\gamma^k\epsilon^l=\gamma\epsilon^2$. Then there is no high degree vertex, a contradiction.

\subsubsection*{Case. $\alpha^2\gamma$ is a vertex}

By distinct $\alpha,\beta,\gamma$ values, and all the degree $3$ $b$-vertices $\alpha\delta^2,\gamma\epsilon^2$, we know $\alpha\delta^2,\gamma\epsilon^2,\alpha^2\gamma,\alpha\beta^2,\beta^3,\beta\gamma^2$ are all the degree $3$ vertices. We already proved there is no tiling with $\alpha\beta^2$.

The angle sums of $\alpha\delta^2,\gamma\epsilon^2,\alpha^2\gamma,\beta\gamma^2$ and the angle sum for pentagon imply
\[
\alpha=\epsilon=(\tfrac{4}{7}+\tfrac{8}{7f})\pi,\,
\beta=(\tfrac{2}{7}+\tfrac{32}{7f})\pi,\,
\gamma=(\tfrac{6}{7}-\tfrac{16}{7f})\pi,\,
\delta=(\tfrac{5}{7}-\tfrac{4}{7f})\pi.
\]
We have $\beta<\alpha=\epsilon<\delta<\gamma<\pi$, contradicting Lemma \ref{geometry4}. 

Therefore $\alpha\delta^2,\gamma\epsilon^2,\alpha^2\gamma,\beta^3$ are all the degree $3$ vertices. The angle sums of the vertices and the angle sum for pentagon imply
\[
	\alpha=\epsilon=(\tfrac{4}{3}-\tfrac{8}{f})\pi,\,
	\beta=\tfrac{2}{3}\pi,\,
	\gamma=(\tfrac{16}{f}-\tfrac{2}{3})\pi,\,
	\delta=(\tfrac{1}{3}+\tfrac{4}{f})\pi.
\]
By $\gamma>0$, we get $f<24$. By Lemma \ref{special_tile}, there is a $3^5$-tile, given by Figure \ref{a2d_c2eB}. We have $\delta\cdots=\alpha\delta^2$ and $\epsilon\cdots=\gamma\epsilon^2$. This implies $\gamma\cdots=\alpha\gamma\cdots=\alpha^2\gamma$ or $\gamma\cdots=\gamma\epsilon\cdots=\gamma\epsilon^2$. In both cases, we get $\alpha\cdots=\alpha^2\gamma$, and the vertex is arranged as indicated. Then we get $\beta\cdots=\beta^3=\thin^{\delta}\beta^{\alpha}\thin^{\alpha}\beta^{\delta}\thin^{\alpha}\beta^{\delta}\thin$ (the middle $\beta$ is the $\beta$ of the special tile). This implies a vertex $\delta\thin\delta\cdots$. On the other hand, by $f<24$, we get $\epsilon>\delta>\frac{1}{2}$. This implies $\delta\thin\delta\cdots$ is not a vertex.

\begin{figure}[htp]
\centering
\begin{tikzpicture}[>=latex,scale=1]

\foreach \a in {0,1}
{
\begin{scope}[xshift=2.5*\a cm]

\draw
	(-0.5,-0.4) -- (-0.5,0.7) -- (0,1.1) -- (0.5,0.7) -- (0.5,-0.4)
	(0.5,0.7) -- ++(0.4,0)
	(0,1.1) -- ++(0,0.4);
		
\draw[line width=1.2]
	(-0.5,0) -- (0.5,0);

\node at (0,0.85) {\small $\alpha$}; 
\node at (-0.3,0.6) {\small $\beta$};
\node at (0.3,0.6) {\small $\gamma$};
\node at (-0.3,0.2) {\small $\delta$};
\node at (0.3,0.2) {\small $\epsilon$};	

\node at (-0.7,0) {\small $\alpha$};
\node at (0.7,0) {\small $\gamma$};
\node at (-0.3,-0.2) {\small $\delta$};
\node at (0.3,-0.2) {\small $\epsilon$};

\node at (0.2,1.2) {\small $\gamma$};
\node at (-0.2,1.2) {\small $\alpha$};

\node at (-0.7,0.7) {\small $H$};

\end{scope}
}

\node at (0.7,0.5) {\small $\alpha$};
\node at (0.6,0.9) {\small $\alpha$};

\begin{scope}[xshift=2.5 cm]

\draw[line width=1.2]
	(0.5,0.7) -- ++(0.4,0);
	
\node at (0.7,0.5) {\small $\epsilon$};
\node at (0.6,0.9) {\small $\epsilon$};

\end{scope}

\end{tikzpicture}
\caption{Proposition \ref{a2d_c2e}: $\alpha^2\gamma$.}
\label{a2d_c2eB}
\end{figure}
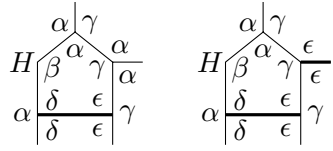

Therefore $\alpha\delta^2,\gamma\epsilon^2,\alpha^2\gamma$ are all the degree $3$ vertices. This implies $\beta\cdots$ has high degree. By Lemma \ref{special_tile}, this implies $f\ge 24$. Moreover, the other four vertices in a special tile have degree $3$, and the argument above shows the special tile is given by  Figure \ref{a2d_c2eB}. Then $H={}^{\alpha}\thin^{\alpha}\beta^{\delta}\thin^{\alpha}\cdots$ has degree $4$ or $5$. 

By Lemma \ref{ndegree3}, we know one of $\alpha\beta^3,\beta^3\gamma,\beta^4,\beta^5$ is a vertex.

\subsubsection*{Subcase. $\alpha\beta^3$ is a vertex}

The angle sums of $\alpha\delta^2,\gamma\epsilon^2,\alpha^2\gamma,\alpha\beta^3$ and the angle sum for pentagon imply
\[
	\alpha=\epsilon=(\tfrac{4}{5}-\tfrac{24}{5f})\pi,\,
	\beta=(\tfrac{2}{5}+\tfrac{8}{5f})\pi,\,
	\gamma=(\tfrac{2}{5}+\tfrac{48}{5f})\pi,\,
	\delta=(\tfrac{3}{5}+\tfrac{12}{5f})\pi.
\]	
By $f\ge 24$, we get $\beta<\gamma<\pi$. By Lemma \ref{geometry1}, this implies $\delta>\epsilon$. Then we get $f<36$, and $\beta<\alpha=\epsilon<\delta<\gamma<\pi$, contradicting Lemma \ref{geometry4}. 

\subsubsection*{Subcase. $\beta^3\gamma$ is a vertex}

The angle sums of $\alpha\delta^2,\gamma\epsilon^2,\alpha^2\gamma,\beta^3\gamma$ and the angle sum for pentagon imply
\[
	\alpha=\epsilon=\tfrac{24}{f}\pi,\,
	\beta=\tfrac{16}{f}\pi,\,
	\gamma=(2-\tfrac{48}{f})\pi,\,
	\delta=(1-\tfrac{12}{f})\pi.
\]	
	
Suppose $\beta<\gamma$ and $\delta>\epsilon$. Then we get $f>36$. This implies $\beta<\alpha=\epsilon<\delta<\gamma$. By Lemma \ref{geometry4}, this implies $\gamma\ge \pi$. By $\gamma\ge \pi$, and $\delta>\epsilon$, and $\delta+\epsilon>\pi$, we know $\gamma^2\cdots,\gamma\delta\cdots,\delta\thin\delta\cdots$ are not vertices. By the AAD, this implies an $\epsilon^2$-fan cannot have only $\alpha,\beta$. Then by $\gamma\epsilon^2$, an $\epsilon^2$-fan is the vertex $\gamma\epsilon^2$. Then by $\delta>\epsilon$ and $\delta+\epsilon>\pi$, this implies $\delta\thin\epsilon\cdots$ is not a vertex. 

If the vertex $H$ in Figure \ref{a2d_c2eB} is a $b$-vertex, then by $\beta+\delta>\pi$, we get $H=\thick^{\delta}\epsilon^{\gamma}\thin^{\delta}\beta^{\alpha}\thin\beta\thin^{\alpha}\beta^{\delta}\thin\theta\thick$, where $\theta=\delta,\epsilon$. By no $\gamma\delta\cdots$, we get a contradiction. 

Therefore $H$ is a $\hat{b}$-vertex. If $H$ has $\gamma$, then by $\beta^3\gamma$ and $\beta<\alpha,\gamma$, we get $H=\beta^3\gamma=\thin^{\delta}\beta^{\alpha}\thin\beta\thin^{\alpha}\beta^{\delta}\thin\gamma\thin$, $\thin^{\delta}\beta^{\alpha}\thin\beta\thin^{\alpha}\gamma^{\epsilon}\thin\beta\thin$. The AAD of $H$ implies $\delta\thin\delta\cdots,\delta\thin\epsilon\cdots$, a contradiction. Therefore $H$ has no $\gamma$. This means $H=\thin\beta\thin^{\alpha}\beta^{\delta}\thin\cdots\thin^{\delta}\beta^{\alpha}\thin$, with only $\alpha,\beta$ in $\cdots$. The AAD implies $\gamma^2\cdots,\gamma\delta\cdots,\delta\thin\delta\cdots$, a contradiction.  

By Lemma \ref{geometry1}, we conclude that $\beta>\gamma$ and $\delta<\epsilon$. This means $f<32$ and implies $\gamma<\delta$. 

By all degree the $3$ vertices $\alpha\delta^2,\gamma\epsilon^2,\alpha^2\gamma$, we know $\beta\epsilon\cdots$ has high degree. By $f\ge 24$, we get $\alpha\le\pi$. Then by the angle sum for pentagon, we get $\beta+\gamma+\delta+\epsilon>2\pi$. Then by $\delta<\epsilon$ and $\gamma<\alpha,\beta,\delta,\epsilon$, this implies $\beta\epsilon\cdots$ is not a vertex.

By $f\ge 24$, we get $\epsilon>\delta>\frac{1}{2}\pi$. This implies at most one $b$-edge at a vertex, and $\delta\thin\delta\cdots,\delta\thin\epsilon\cdots,\epsilon\thin\epsilon\cdots$ are not vertices. Then the remainder of $H=\thin\theta^{\alpha}\thin\beta\thin^{\alpha}\theta\thin\cdots$ cannot consist of only $\theta=\beta,\gamma$. 

If $H$ is a $\hat{b}$-vertex, then $H$ has $\alpha$. By $\beta^3\gamma$ and $\alpha>\beta>\gamma$, we get $H=\alpha\beta\gamma^2,\alpha\beta\gamma^3$. Then by no $\epsilon\thin\epsilon\cdots$, we get $H=\thin^{\epsilon}\gamma^{\alpha}\thin\beta\thin^{\alpha}\gamma^{\epsilon}\thin\alpha\thin,\thin^{\epsilon}\gamma^{\alpha}\thin\beta\thin^{\alpha}\gamma^{\epsilon}\thin^{\alpha}\gamma^{\epsilon}\thin\alpha\thin$, contradicting no $\beta\epsilon\cdots$.  

If $H$ is a $b$-vertex, then by no $\beta\epsilon\cdots$, we get $H=\thick\delta\thin\theta^{\alpha}\thin\beta\thin^{\alpha}\theta\thin\delta\thick$, with $\theta=\beta,\gamma$. By $f\ge 24$, we get $2\beta+\gamma+2\delta=(4-\frac{40}{f})\pi>2\pi$. Then by $\beta>\gamma$, this implies $H=\beta\gamma^2\delta^2=\thick^{\epsilon}\delta^{\beta}\thin^{\epsilon}\gamma^{\alpha}\thin\beta\thin^{\alpha}\gamma^{\epsilon}\thin^{\beta}\delta^{\epsilon}\thick$, contradicting no $\beta\epsilon\cdots$.  

\subsubsection*{Subcase. $\beta^4$ is a vertex}

The angle sums of $\alpha\delta^2,\gamma\epsilon^2,\alpha^2\gamma,\beta^4$ and the angle sum for pentagon imply
\[
	\alpha=\epsilon=(1-\tfrac{8}{f})\pi,\,
	\beta=\tfrac{1}{2}\pi,\,
	\gamma=\tfrac{16}{f}\pi,\,
	\delta=(\tfrac{1}{2}+\tfrac{4}{f})\pi.
\]	
By $f\ge 24$, we get $\delta<\epsilon$. By Lemma \ref{geometry1}, this implies $\beta>\gamma$. Then we get $\alpha=\epsilon>\delta>\beta>\gamma$. We also have $\alpha<\pi$ and $\delta>\frac{1}{2}\pi$. 

Similar to the case $\beta^3\gamma$, we know $\beta\epsilon\cdots$ is not a vertex, and the remainder of $H=\thin\theta^{\alpha}\thin\beta\thin^{\alpha}\theta\thin\cdots$ cannot consist of only $\theta=\beta,\gamma$. 

If $H$ is a $\hat{b}$-vertex, then $H$ has $\alpha$. By $\alpha>\beta>\gamma$ and $\alpha+\gamma=(1+\frac{8}{f})\pi>2\beta=\pi$, we still get $H=\thin^{\epsilon}\gamma^{\alpha}\thin\beta\thin^{\alpha}\gamma^{\epsilon}\thin\alpha\thin,\thin^{\epsilon}\gamma^{\alpha}\thin\beta\thin^{\alpha}\gamma^{\epsilon}\thin^{\alpha}\gamma^{\epsilon}\thin\alpha\thin$, contradicting no $\beta\epsilon\cdots$.  

If $H$ is a $b$-vertex. Then by no $\beta\epsilon\cdots$, we get $H=\thick\delta\thin\theta^{\alpha}\thin\beta\thin^{\alpha}\theta\thin\delta\thick$, with $\theta=\beta,\gamma$. By $\gamma<\beta=\frac{1}{2}\pi<\delta$, we still get $H=\beta\gamma^2\delta^2=\thick^{\epsilon}\delta^{\beta}\thin^{\epsilon}\gamma^{\alpha}\thin\beta\thin^{\alpha}\gamma^{\epsilon}\thin^{\beta}\delta^{\epsilon}\thick$, contradicting no $\beta\epsilon\cdots$. 

\subsubsection*{Subcase. $\beta^5$ is a vertex}

The angle sums of $\alpha\delta^2,\gamma\epsilon^2,\alpha^2\gamma,\beta^5$ and the angle sum for pentagon imply
\[
	\alpha=\epsilon=(\tfrac{4}{5}-\tfrac{8}{f})\pi,\,
	\beta=\tfrac{2}{5}\pi,\,
	\gamma=(\tfrac{2}{5}+\tfrac{16}{f})\pi,\,
	\delta=(\tfrac{3}{5}+\tfrac{4}{f})\pi.
\]	
By $\beta<\gamma$ and Lemma \ref{geometry1}, we get $\delta>\epsilon$. This implies $f<60$, and further implies $\delta<\gamma$. On the other hand, by $f\ge 24$, we get $\alpha=\epsilon>\beta$. Therefore $\beta<\alpha=\epsilon<\delta<\gamma<\pi$, contradicting Lemma \ref{geometry4}.

\subsubsection*{Case. $\alpha\gamma^2$ is a vertex}

By distinct $\alpha,\beta,\gamma$ values, and all the degree $3$ $b$-vertices $\alpha\delta^2,\gamma\epsilon^2$, we know $\alpha\delta^2,\gamma\epsilon^2,\alpha\gamma^2,\alpha^2\beta,\beta^3,\beta^2\gamma$ are all the degree $3$ vertices. 

The angle sums of $\alpha\delta^2,\gamma\epsilon^2,\alpha\gamma^2,\alpha^2\beta$ and the angle sum for pentagon imply 
\[
\alpha=(\tfrac{6}{7}-\tfrac{16}{7f})\pi,\,
\beta=(\tfrac{2}{7}+\tfrac{32}{7f})\pi,\,
\gamma=\delta=(\tfrac{4}{7}+\tfrac{8}{7f})\pi,\,
\epsilon=(\tfrac{5}{7}-\tfrac{4}{7f})\pi.
\]
We get $\beta<\gamma=\delta<\epsilon$, contradicting Lemma \ref{geometry1}.

The angle sums of $\alpha\delta^2,\gamma\epsilon^2,\alpha\gamma^2,\beta^2\gamma$ and the angle sum for pentagon imply 
\[
\alpha=\tfrac{8}{f}\pi,\,
\beta=\epsilon=(\tfrac{1}{2}+\tfrac{2}{f})\pi,\,
\gamma=\delta=(1-\tfrac{4}{f})\pi.
\] 
The angle values appeared in \eqref{bde_a2c_eq1}, for the subcase $\delta^2\cdots=\alpha\delta^2$ of the case $\beta<\gamma$ and $\delta>\epsilon$ in the proof of Proposition \ref{bde_a2c}. The same argument, including the use of Lemma \ref{geometry10} to get $f<28$, is still valid and leads to contradiction. 

Therefore $\alpha\delta^2,\gamma\epsilon^2,\alpha\gamma^2,\beta^3$ are all the degree $3$ vertices. The angle sums of the vertices and the angle sum for pentagon imply 
\[
\alpha=(\tfrac{16}{f}-\tfrac{2}{3})\pi,\,
\beta=\tfrac{2}{3}\pi,\,
\gamma=\delta=(\tfrac{4}{3}-\tfrac{8}{f})\pi,\,
\epsilon=(\tfrac{1}{3}+\tfrac{4}{f})\pi.
\]
If the vertices $\alpha\cdots,\beta\cdots$ in a tile have degree $3$, then $\alpha\cdots=\alpha\delta^2$ and $\beta\cdots=\beta^3$. This implies $\gamma\cdots=\beta\gamma\cdots$ has high degree. By Lemma \ref{special_tile}, this implies $f>24$. Then we get $\alpha\le 0$, a contradiction. 

Therefore $\alpha\delta^2,\gamma\epsilon^2,\alpha\gamma^2$ are all the degree $3$ vertices. Consider the special companion pair in Figure \ref{a2d_c2eA}. If $\delta_1\delta_2\cdots$ has degree $3$, then $\delta_1\delta_2\cdots=\alpha\delta^2$, and both $\beta_1\cdots,\beta_2\cdots$ have high degree. If $\epsilon_1\epsilon_2\cdots$ has degree $3$, then $\epsilon_1\epsilon_2\cdots=\gamma\epsilon^2$. Up to the vertical flip, we may assume $T_4$ is arranged as indicated. If $\alpha_4\gamma_1\cdots$ has degree $3$, then $\alpha_1\cdots=\alpha^2\cdots,\alpha\epsilon\cdots$ has high degree. We conclude the special companion pair has weight $\ge 3w_4$. By Lemma \ref{special_pair}, this implies $f\ge 48$. 

By all degree $3$ vertices $\alpha\delta^2,\gamma\epsilon^2,\alpha\gamma^2$ and Lemma \ref{ndegree3}, we know one of $\alpha\beta^3,\beta^3\gamma,\beta^4,\beta^5$ is a vertex. The angle sums of $\alpha\delta^2,\gamma\epsilon^2,\alpha\gamma^2$, and the angle sum of one of these, and the angle sum for pentagon imply 
\begin{align*}
\alpha\beta^3 &\colon
	\alpha=(2-\tfrac{48}{f})\pi,\,
	\beta=\tfrac{16}{f}\pi,\,
	\gamma=\delta=\tfrac{24}{f}\pi,\,
	\epsilon=(1-\tfrac{12}{f})\pi. \\
\beta^3\gamma &\colon
	\alpha=(\tfrac{2}{5}+\tfrac{48}{5f})\pi,\,
	\beta=(\tfrac{2}{5}+\tfrac{8}{5f})\pi,\,
	\gamma=\delta=(\tfrac{4}{5}-\tfrac{24}{5f})\pi,\,
	\epsilon=(\tfrac{3}{5}+\tfrac{12}{5f})\pi. \\
\beta^4 &\colon
	\alpha=\tfrac{16}{f}\pi,\,
	\beta=\tfrac{1}{2}\pi,\,
	\gamma=\delta=(1-\tfrac{8}{f})\pi,\,
	\epsilon=(\tfrac{1}{2}+\tfrac{4}{f})\pi. \\
\beta^5 &\colon
	\alpha=(\tfrac{2}{5}+\tfrac{16}{f})\pi,\,
	\beta=\tfrac{2}{5}\pi,\,
	\gamma=\delta=(\tfrac{4}{5}-\tfrac{8}{f})\pi,\,
	\epsilon=(\tfrac{3}{5}+\tfrac{4}{f})\pi.
\end{align*}

For $\alpha\beta^3$, by $f\ge 48$, we get $\beta<\gamma=\delta<\epsilon$, contradicting Lemma \ref{geometry1}.

For $\beta^3\gamma$, by $f\ge 48$, we get $\pi>\gamma=\delta>\epsilon>\alpha>\beta$, contradicting Lemma \ref{geometry4}.

For $\beta^4$, the pentagon is convex. By Lemma \ref{geometry10}, we get $2\alpha+\beta+\gamma=(\tfrac{3}{2}\pi+\frac{24}{f})\pi>2\pi$, contradicting $f\ge 48$.  

For $\beta^5$, by $f\ge 48$, we get $\beta<\gamma$. By Lemma \ref{geometry1}, this implies $\delta>\epsilon$, which means $f>60$. Then we get $\beta<\alpha<\epsilon<\gamma=\delta<\pi$, contradicting Lemma \ref{geometry4}.

\subsubsection*{Case. $\beta^2\gamma$ is a vertex}

By distinct $\alpha,\beta,\gamma$ values, and all the degree $3$ $b$-vertices $\alpha\delta^2,\gamma\epsilon^2$, we know $\alpha\delta^2,\gamma\epsilon^2,\beta^2\gamma,\alpha^3,\alpha^2\beta,\alpha\gamma^2$ are all the degree $3$ vertices. The angle sums of $\alpha\delta^2,\gamma\epsilon^2,\beta^2\gamma,\alpha^3$ imply $f=12$, a contradiction. The angle sums of $\alpha\delta^2,\gamma\epsilon^2,\beta^2\gamma,\alpha^2\beta$ and the angle sum for pentagon imply $\gamma=(\tfrac{32}{f}-2)\pi\le 0$, a contradiction. We already proved there is no tiling with $\alpha\gamma^2$. 

Therefore $\alpha\delta^2,\gamma\epsilon^2,\beta^2\gamma$ are all the degree $3$ vertices. If the vertex $\gamma\cdots$ in a tile has degree $3$, then it is $\gamma\epsilon^2,\beta^2\gamma$. If the vertex is $\gamma\epsilon^2$, then $\alpha\cdots=\alpha\gamma\cdots$ has high degree. If the vertex is $\beta^2\gamma$, then $\epsilon\cdots=\alpha\epsilon\cdots,\delta\epsilon\cdots$ has high degree. Therefore there is no $3^5$-tile. By Lemma \ref{special_tile}, we get $f\ge 24$. 

The angle sums of $\alpha\delta^2,\gamma\epsilon^2,\beta^2\gamma$ and the angle sum for pentagon imply
\[
\alpha=\tfrac{8}{f}\pi,\,
\beta=\epsilon,\,
2\beta+\gamma=2\pi,\,
\delta=(1-\tfrac{4}{f})\pi.
\]
We have $\delta,\epsilon<\pi$. Then by Lemma \ref{geometry6}, we get $\alpha+2\beta>\pi$ and $\alpha+2\gamma>\pi$. By $f\ge 24$, this means $\beta,\gamma,\epsilon>(\frac{1}{2}-\frac{4}{f})\pi>\alpha$. We also have $\delta>\alpha$. 

By $\alpha\delta^2,\gamma\epsilon^2$, and $\alpha<\gamma$, we get $\delta>\epsilon$. By Lemma \ref{geometry1}, this implies $\beta<\gamma$. Then by $\beta^2\gamma$, we get $\beta<\frac{2}{3}\pi<\gamma$. We conclude $\alpha<\beta=\epsilon<\gamma,\delta$.

By $\alpha\delta^2$ and $\alpha<\beta,\gamma,\delta,\epsilon$, we get $\delta^2\cdots=\alpha\delta^2$. This implies $\delta\thin\delta\cdots$ is not a vertex, and the only $\delta^2$-fan is the vertex $\alpha\delta^2$. 

By $\gamma\epsilon^2,\beta^2\gamma$ and $\alpha<\beta=\epsilon<\gamma,\delta$, we get $\gamma\epsilon\cdots=\gamma\epsilon^2$ and $\beta\gamma\cdots=\beta^2\gamma,\alpha^k\beta\gamma(k\ge 2)$. 

If $\gamma<\delta$, then by $\delta<\pi$, the pentagon is convex. By Lemma \ref{geometry10}, we get $2\alpha+\beta+\gamma>2\pi$. By $\beta^2\gamma$, this implies $2\alpha>\beta=\epsilon$. By $\beta^2\gamma$ and $\beta<\gamma<\delta$, we get $R(\gamma\delta)<R(\gamma^2)<\beta=\epsilon<2\alpha,\gamma,\delta$. This implies $\gamma^2\cdots,\gamma\delta\cdots$ are not vertices.

If $\gamma>\delta$, then by $\alpha\delta^2$, we get  $R(\gamma^2)<R(\gamma\delta)<\alpha<\beta,\gamma,\delta,\epsilon$. This again implies $\gamma^2\cdots,\gamma\delta\cdots$ are not vertices. 

By no $\gamma^2\cdots,\gamma\delta\cdots$, and $\gamma\epsilon\cdots=\gamma\epsilon^2$, and $\beta\gamma\cdots=\beta^2\gamma,\alpha^k\beta\gamma(k\ge 2)$, we get $\gamma\cdots=\beta^2\gamma,\gamma\epsilon^2,\alpha^k\beta\gamma,\alpha^k\gamma$.

By no $\gamma^2\cdots,\gamma\delta\cdots,\delta\thin\delta\cdots$, we know an $\epsilon^2$-fan cannot consist of only $\alpha,\beta$. Then by $\gamma\epsilon^2$, we know the only $\epsilon^2$-fan is the vertex $\gamma\epsilon^2$.

Since the only $\delta^2$-fan is $\alpha\delta^2$, and the only $\epsilon^2$-fan is $\gamma\epsilon^2$, a $b$-vertex $\beta\cdots$ is a combination of $\delta\epsilon$-fans. By $\gamma\epsilon\cdots=\gamma\epsilon^2$, the vertex has no $\gamma$. By $\beta,\epsilon>(\frac{1}{2}-\frac{4}{f})\pi$ and $\delta=(1-\tfrac{4}{f})\pi$, we get $R(\beta\delta\epsilon)<\frac{12}{f}\pi\le 2\alpha,\delta$. Then by $\alpha<\beta$, and no $\beta\delta\epsilon$, and the vertex being combination of $\delta\epsilon$-fans, we know a $b$-vertex $\beta\cdots=\alpha\beta\delta\epsilon,\beta^2\delta\epsilon$. 

The only remaining $b$-vertex is a combination of $\delta\epsilon$-fans without $\beta$. By $\alpha\delta^2$, and $\gamma\epsilon\cdots=\gamma\epsilon^2$, and $\alpha<\beta,\gamma,\delta,\epsilon$, the vertex is $\alpha^k\delta\epsilon(k\ge 2)$. We conclude $\alpha\delta^2,\gamma\epsilon^2,\alpha\beta\delta\epsilon,\beta^2\delta\epsilon,\alpha^k\delta\epsilon$ are all the $b$-vertices. This implies $\delta\thin\epsilon\cdots$ is not a vertex. 

By the list of $b$-vertices, we get $\beta\thin\delta\cdots=\thick\delta\thin\beta\thin\alpha\thin\epsilon\thick,\thick\delta\thin\beta\thin\beta\thin\epsilon\thick$. By no $\gamma^2\cdots,\gamma\delta\cdots,\delta\thin\delta\cdots$, the AAD of the vertex is $\thick^{\epsilon}\delta^{\beta}\thin^{\delta}\beta^{\alpha}\thin^{\gamma}\alpha^{\beta}\thin^{\gamma}\epsilon^{\delta}\thick,\thick^{\epsilon}\delta^{\beta}\thin^{\delta}\beta^{\alpha}\thin^{\delta}\beta^{\alpha}\thin^{\gamma}\epsilon^{\delta}\thick$. Therefore the AAD of $\thin\beta\thin\delta\thick$ is $\thin^{\alpha}\beta^{\delta}\thin^{\beta}\delta^{\epsilon}\thick$. By no $\gamma\delta\cdots$, this implies the AAD of $\thin\alpha\thin\beta\thin$ cannot be $\thin\alpha\thin^{\delta}\beta^{\alpha}\thin$, and must be $\thin\alpha\thin^{\alpha}\beta^{\delta}\thin$. Then by no $\gamma^2\cdots,\gamma\delta\cdots,\delta\thin\delta\cdots$,  this implies $\alpha^k\beta^l=\alpha^k,\beta^l$, and the AAD of $\alpha^k\beta\gamma$ is $\thin\alpha\thin\cdots\thin\alpha\thin^{\alpha}\beta^{\delta}\thin\gamma\thin$. Then by $\beta\gamma\cdots=\beta^2\gamma,\alpha^k\beta\gamma$, and $b$-vertex $\beta\cdots=\alpha\beta\delta\epsilon,\beta^2\delta\epsilon$, we get $\beta\cdots=\beta^2\gamma,\alpha\beta\delta\epsilon,\beta^2\delta\epsilon,\alpha^k\beta\gamma,\beta^l$. 

The AAD $\thin^{\beta}\alpha^{\gamma}\thin^{\beta}\alpha^{\gamma}\thin$ implies a vertex $\thin^{\delta}\beta^{\alpha}\thin^{\alpha}\gamma^{\epsilon}\thin\cdots=\beta^2\gamma,\alpha^k\beta\gamma$. The AAD is incompatible with the AAD $\thin\alpha\thin\cdots\thin\alpha\thin^{\alpha}\beta^{\delta}\thin\gamma\thin$ of $\alpha^k\beta\gamma$. If the vertex is $\beta^2\gamma$, then we get $\thin^{\alpha}\gamma^{\epsilon}\thin\beta\thin^{\delta}\beta^{\alpha}\thin$, contradicting no $\delta\thin\delta\cdots,\delta\thin\epsilon\cdots$. 

The AAD $\thin^{\gamma}\alpha^{\beta}\thin^{\beta}\alpha^{\gamma}\thin$ implies a vertex $\thin^{\delta}\beta^{\alpha}\thin^{\alpha}\beta^{\delta}\thin\cdots=\beta^2\gamma,\beta^2\delta\epsilon,\beta^l$. The AAD is incompatible with the AAD $\thick^{\epsilon}\delta^{\beta}\thin^{\delta}\beta^{\alpha}\thin^{\delta}\beta^{\alpha}\thin^{\gamma}\epsilon^{\delta}\thick$ of $\beta^2\delta\epsilon$. If the vertex is $\beta^2\gamma$, then we get $\thin^{\alpha}\beta^{\delta}\thin\gamma\thin^{\delta}\beta^{\alpha}\thin$, contradicting no $\delta\thin\epsilon\cdots$. If the vertex is $\beta^l$, then we get $\thin^{\alpha}\beta^{\delta}\thin^{\delta}\beta^{\alpha}\thin$, contradicting no $\delta\thin\delta\cdots$.

By no $\gamma^2\cdots,\thin^{\beta}\alpha^{\gamma}\thin^{\beta}\alpha^{\gamma}\thin,\thin^{\gamma}\alpha^{\beta}\thin^{\beta}\alpha^{\gamma}\thin$, we know $\alpha\thin\alpha\cdots$ is not a vertex. This implies $\alpha^k,\alpha^k\gamma,\alpha^k\delta\epsilon$ are not vertices. By the AAD $\thin\alpha\thin\cdots\thin\alpha\thin^{\alpha}\beta^{\delta}\thin\gamma\thin$ of $\alpha^k\beta\gamma$, we also know $\alpha^k\beta\gamma$ is not a vertex. Then we conclude $\alpha\delta^2,\gamma\epsilon^2,\beta^2\gamma,\alpha\beta\delta\epsilon,\beta^2\delta\epsilon,\beta^l$ are all the vertices. This implies 
\[
f=\#\gamma=\#\gamma\epsilon^2+\#\beta^2\gamma\le \tfrac{1}{2}\#\epsilon+\tfrac{1}{2}\#\beta=f.
\]
Since both sides are equal, this implies the only vertex besides $\gamma\epsilon^2,\beta^2\gamma$ is $\alpha^k\delta^l=\alpha\delta^2$. Then there is no high degree vertex, a contradiction.

\subsubsection*{Case. $\beta\gamma^2$ is a vertex}

By distinct $\alpha,\beta,\gamma$ values, and all the degree $3$ $b$-vertices $\alpha\delta^2,\gamma\epsilon^2$, We know $\alpha\delta^2,\gamma\epsilon^2,\beta\gamma^2,\alpha^3,\alpha\beta^2,\alpha^2\gamma$ are all the degree $3$ vertices. We already proved there is no tiling with $\alpha\beta^2$ or $\alpha^2\gamma$. The angle sums of $\alpha\delta^2,\gamma\epsilon^2,\beta\gamma^2,\alpha^3$ and the angle sum for pentagon imply
\[
\alpha=\delta=\tfrac{2}{3}\pi,\,
\beta=(\tfrac{2}{9}+\tfrac{16}{3f})\pi,\,
\gamma=(\tfrac{8}{9}-\tfrac{8}{3f})\pi,\,
\epsilon=(\tfrac{5}{9}+\tfrac{4}{3f})\pi.
\]
We have $\beta<\epsilon<\alpha=\delta<\gamma<\pi$, contradicting Lemma \ref{geometry4}. Therefore $\alpha\delta^2,\gamma\epsilon^2,\beta\gamma^2$ are all the degree $3$ vertices. 

The angle sums of $\alpha\delta^2,\beta\gamma^2$ and the angle sum for pentagon imply $\gamma+\delta-\epsilon=(1-\frac{4}{f})\pi$. If $\beta>\gamma$, then by $\beta\gamma^2$, we get $\gamma<\frac{2}{3}\pi$. This implies $\delta-\epsilon>(\frac{1}{3}-\frac{4}{f})\pi>0$, contradicting Lemma \ref{geometry1}. By Lemma \ref{geometry1}, therefore, we get $\beta<\gamma$ and $\delta>\epsilon$. Then by $\alpha\delta^2,\gamma\epsilon^2$, we further get $\alpha<\gamma$. 

By $\beta\gamma^2$, we get $\gamma<\pi$. Then by $\gamma\epsilon^2$, we get $\epsilon>\frac{1}{2}\pi$. By $\delta>\epsilon>\frac{1}{2}\pi$, there is at most one $b$-edge at a vertex. In particular, $\epsilon\thin\epsilon\cdots$ is not a vertex.

Consider the special companion pair in Figure \ref{a2d_c2eA}. If $\epsilon_1\epsilon_2\cdots$ has degree $3$, then it is $\gamma\epsilon^2$, and up to the vertical flip, we may assume $T_4$ is arranged as indicated. This implies $\gamma_1\cdots=\alpha\gamma\cdots$ has high degree. Moreover, if $\gamma_2\cdots$ has degree $3$, then it is $\gamma\epsilon\cdots=\gamma\epsilon^2$, which implies $\alpha_2\cdots=\alpha\gamma\cdots$ has high degree. Therefore $\gamma_1\cdots$ has high degree, and one of $\alpha_2\cdots,\gamma_2\cdots$ has high degree. This implies $\delta_1\delta_2\cdots$ has degree $3$, and one of $\beta_1\cdots,\beta_2\cdots$ has degree $3$. If $\epsilon_1\epsilon_2\cdots$ has high degree, then $\delta_1\delta_2\cdots$ also has degree $3$, and one of $\beta_1\cdots,\beta_2\cdots$ has degree $3$. 

Therefore we always have $\delta_1\delta_2\cdots=\alpha\delta^2$, and up to the vertical flip, we may assume $T_3$ is arranged as indicated. Then $\beta_1\cdots=\beta^2\cdots$ has high degree. Since one of $\beta_1\cdots,\beta_2\cdots$ has degree $3$, this implies $\beta_2\cdots$ has degree $3$. Then $\beta_2\cdots=\beta\gamma^2$. By no $\epsilon\thin\epsilon\cdots$, this determines $T_7$. Then $\alpha_2\cdots=\alpha\epsilon\cdots$ has high degree. 

Given two high degree vertices $\beta_1\cdots,\alpha_2\cdots$, we know $\epsilon_1\epsilon_2\cdots$ has degree $3$, and at most one high degree among $\alpha_1\cdots,\gamma_1\cdots,\gamma_2\cdots$. Then we get $\epsilon_1\epsilon_2\cdots=\gamma\epsilon^2$. If $T_4$ is not arranged as indicated, then as argued earlier, we know $\gamma_2\cdots$ has high degree, and one of $\alpha_1\cdots,\gamma_1\cdots$ has high degree, a contradiction. Therefore $T_4$ is arranged as indicated. Then $\gamma\cdots=\alpha\gamma\cdots$ has high degree. This implies $\gamma_2\cdots$ has degree $3$. Then the vertex is $\gamma\epsilon\cdots=\gamma\epsilon^2$. This determines $T_8$. Then we get $\alpha_2\gamma_8\epsilon_7\cdots$. By $\gamma\epsilon^2$ and $\delta>\epsilon$, the angle sum of $\alpha\gamma\epsilon\cdots$ is $\ge \alpha+\gamma+2\epsilon>\gamma+2\epsilon=2\pi$, a contradiction.

\subsubsection*{Case. $\beta^3$ is a vertex}

By distinct $\alpha,\beta,\gamma$ values, and all the degree $3$ $b$-vertices $\alpha\delta^2,\gamma\epsilon^2$, we know $\alpha\delta^2,\gamma\epsilon^2,\beta^3,\alpha\beta\gamma,\alpha^2\gamma,\alpha\gamma^2$ are all the degree $3$ vertices. We already proved there is no tiling if $\alpha^2\gamma$ or $\alpha\gamma^2$ is a vertex. 

The angle sums of $\alpha\delta^2,\gamma\epsilon^2,\beta^3$ and the angle sum for pentagon imply
\[
\alpha+\gamma=(\tfrac{2}{3}+\tfrac{8}{f})\pi,\;
\beta=\tfrac{2}{3}\pi,\;
\delta+\epsilon=(\tfrac{5}{3}-\tfrac{4}{f})\pi.
\]
Then $\alpha+\beta+\gamma<2\pi$. Therefore $\alpha\beta\gamma$ is not a vertex, and $\alpha\delta^2,\gamma\epsilon^2,\beta^3$ are all the degree $3$ vertices. 

Consider the special companion pair in Figure \ref{a2d_c2eA}. We know one of $\alpha_1\cdots,\gamma_1\cdots$ has high degree, and one of $\alpha_2\cdots,\gamma_2\cdots$ has high degree. This implies $\delta_1\delta_2\cdots,\epsilon_1\epsilon_2\cdots$ have degree $3$. Then $\delta_1\delta_2\cdots=\alpha\delta^2$ and $\epsilon_1\epsilon_2\cdots=\gamma\epsilon^2$. This implies one of $\beta_1\cdots,\beta_2\cdots$ is $\beta\gamma\cdots$, which has high degree. Then the the special companion pair has at least three high degree vertices, and has weight $\ge 3w_4$. By Lemma \ref{special_pair}, this implies $f\ge 48$.

By $\epsilon_1\epsilon_2\cdots=\gamma\epsilon^2$, and up to the vertical flip, we may assume $T_4$ is arranged as indicated. This implies $\gamma_1\cdots=\alpha\gamma\cdots$ has high degree. On the other hand, a special companion pair can have at most three high degree vertices. Therefore exactly one of $\alpha_1\cdots,\gamma_1\cdots$ has high degree. Then $\alpha_1\cdots$ has degree $3$, and $\alpha_1\cdots=\alpha\delta^2$. Then we get a vertex $H=\alpha_4\beta_6\gamma_1\cdots=\thin^{\beta}\alpha^{\gamma}\thin^{\epsilon}\gamma^{\alpha}\thin^{\delta}\beta^{\alpha}\thin\cdots$. 

By $\alpha\delta^2,\gamma\epsilon^2$, we get $\delta,\epsilon<\pi$. Then by $\delta+\epsilon=(\tfrac{5}{3}-\tfrac{4}{f})\pi$ and $f\ge 48$, we get $\delta,\epsilon>(\tfrac{2}{3}-\tfrac{4}{f})\pi>\frac{1}{2}\pi$. This implies $2\delta,2\epsilon>\pi>\alpha+\gamma,\beta$, and $\delta\thin\delta\cdots,\delta\thin\epsilon\cdots,\epsilon\thin\epsilon\cdots$ are not vertices.

By $\alpha\delta^2,\gamma\epsilon^2$, a vertex $\alpha\beta\gamma\cdots$ has no $\delta,\epsilon$. By no $\alpha\beta\gamma$, we know the remainder is not empty. Then by $\alpha+\gamma>\beta=\frac{2}{3}\pi$, we get $\alpha\beta\gamma\cdots=\alpha^k\beta\gamma(k\ge 2),\alpha\beta\gamma^k(k\ge 2)$.

\subsubsection*{Subcase. $H=\alpha^2\beta\gamma$}

The angle sums of $\alpha\delta^2,\gamma\epsilon^2,\beta^3,\alpha^2\beta\gamma$ and the angle sum for pentagon imply
\[
	\alpha=(\tfrac{2}{3}-\tfrac{8}{f})\pi,\,
	\beta=\tfrac{2}{3}\pi,\,
	\gamma=\tfrac{16}{f}\pi,\,
	\delta=(\tfrac{2}{3}+\tfrac{4}{f})\pi,\,
	\epsilon=(1-\tfrac{8}{f})\pi.
\]
By $f\ge 48$, we get $\epsilon>\delta>\beta>\alpha>\gamma$. By $\alpha\delta^2,\gamma\epsilon^2$, and $\beta>\alpha>\gamma$, we know $\beta\cdots$ is a $\hat{b}$-vertex. 

The vertex $H=\thin^{\beta}\alpha^{\gamma}\thin^{\epsilon}\gamma^{\alpha}\thin^{\delta}\beta^{\alpha}\thin\cdots=\alpha^2\beta\gamma=\thin\alpha\thin^{\beta}\alpha^{\gamma}\thin^{\epsilon}\gamma^{\alpha}\thin^{\delta}\beta^{\alpha}\thin$ implies a vertex $\thin^{\delta}\beta^{\alpha}\thin^{\alpha}\beta^{\delta}\thin\cdots$ or $\thin^{\delta}\beta^{\alpha}\thin^{\alpha}\gamma^{\epsilon}\thin\cdots$. Since $\beta\cdots$ is a $\hat{b}$-vertex, the vertex has no $\delta,\epsilon$. By no $\delta\thin\delta\cdots,\delta\thin\epsilon\cdots,\epsilon\thin\epsilon\cdots$, the vertex is not $\beta^k\gamma^l$. Therefore the vertex has $\alpha$.

By $R(\alpha\beta^2)=\frac{8}{f}\pi<$ all angles, we know $\thin^{\delta}\beta^{\alpha}\thin^{\alpha}\beta^{\delta}\thin\cdots=\alpha\beta^2\cdots$ is not a vertex. Therefore by $\beta\cdots$ being $\hat{b}$-vertex, the vertex is $\alpha\thin^{\delta}\beta^{\alpha}\thin^{\alpha}\gamma^{\epsilon}\thin\cdots=\alpha^2\beta\gamma,\alpha\beta\gamma^k$. If the vertex is $\alpha^2\beta\gamma$, then by $\beta\cdots$ being $\hat{b}$-vertex, it is $\thin^{\gamma}\alpha^{\beta}\thin^{\beta}\alpha^{\gamma}\thin^{\delta}\beta^{\alpha}\thin^{\alpha}\gamma^{\epsilon}\thin$. However, this implies $\thin^{\delta}\beta^{\alpha}\thin^{\alpha}\beta^{\delta}\thin\cdots$, which we just proved cannot be a vertex. If the vertex is $\alpha\beta\gamma^k$, then by no $\beta\epsilon\cdots,\epsilon\thin\epsilon\cdots$, it is $\thin^{\delta}\beta^{\alpha}\thin^{\alpha}\gamma^{\epsilon}\thin\cdots\thin^{\alpha}\gamma^{\epsilon}\thin^{\gamma}\alpha^{\beta}\thin$ or $\thin^{\delta}\beta^{\alpha}\thin^{\alpha}\gamma^{\epsilon}\thin\cdots\thin^{\alpha}\gamma^{\epsilon}\thin^{\gamma}\alpha^{\beta}\thin^{\alpha}\gamma^{\epsilon}\thin\cdots\thin^{\alpha}\gamma^{\epsilon}\thin$. This implies $\beta\delta\cdots,\delta\thin\epsilon\cdots$, a contradiction.

\subsubsection*{Subcase. $\alpha^k\beta\gamma(k>2)$ is a vertex}

By $k>2$, the angle sum of $\alpha^k\beta\gamma$ implies $2\alpha\le R(\alpha\beta\gamma)=(\tfrac{2}{3}-\tfrac{8}{f})\pi<\alpha+\gamma$. By $\alpha\delta^2,\gamma\epsilon^2$, this implies $\delta>\epsilon$. By Lemma \ref{geometry1}, we get $\beta<\gamma$. This implies $\alpha<\frac{8}{f}\pi<\beta=\frac{2}{3}\pi<\gamma$.

By $\beta^3$ and $\beta<\gamma$, we get $R(\beta\gamma)<R(\beta^2)=\beta<\gamma,2\delta,2\epsilon$. This implies $\beta^2\cdots=\beta^3,\alpha^k\beta^2$, and $\beta\gamma\cdots=\alpha^k\beta\gamma$. 

The vertex $H=\thin^{\beta}\alpha^{\gamma}\thin^{\epsilon}\gamma^{\alpha}\thin^{\delta}\beta^{\alpha}\thin\cdots=\alpha^k\beta\gamma=\thin\alpha\thin^{\beta}\alpha^{\gamma}\thin^{\epsilon}\gamma^{\alpha}\thin^{\delta}\beta^{\alpha}\thin\cdots$ implies a vertex $\thin^{\delta}\beta^{\alpha}\thin^{\alpha}\beta^{\delta}\thin\cdots=\beta^3,\alpha^k\beta^2$, or $\thin^{\delta}\beta^{\alpha}\thin^{\alpha}\gamma^{\epsilon}\thin\cdots=\alpha^k\beta\gamma$. If $\thin^{\delta}\beta^{\alpha}\thin^{\alpha}\beta^{\delta}\thin\cdots=\beta^3$, then the AAD of the vertex implies a vertex $\delta\thin\delta\cdots$, a contradiction. On the other hand, both $\thin^{\delta}\beta^{\alpha}\thin^{\alpha}\beta^{\delta}\thin\cdots=\alpha^k\beta^2$ and $\thin^{\delta}\beta^{\alpha}\thin^{\alpha}\gamma^{\epsilon}\thin\cdots=\alpha^k\beta\gamma$ are of the form $\thin\alpha\thin^{\delta}\beta^{\alpha}\thin\cdots$. This implies either $\beta\delta\cdots$ or $\gamma\delta\cdots$ is a vertex. However, by $\beta<\gamma$, and $\delta>\epsilon$, and $\beta+\delta+\epsilon=(\frac{7}{3}-\frac{4}{f})\pi>2\pi$, we know $\beta\delta\cdots,\gamma\delta\cdots$ are not vertices.

\subsubsection*{Subcase. $\alpha\beta\gamma^k(k>2)$ is a vertex}

By $k>2$, the angle sum of $\alpha\beta\gamma^k$ implies $2\gamma\le R(\alpha\beta\gamma)=(\tfrac{2}{3}-\tfrac{8}{f})\pi<\alpha+\gamma$. By $\alpha\delta^2,\gamma\epsilon^2$, this implies $\delta<\epsilon$. Then by $\alpha\delta^2$, we know $\alpha\epsilon\cdots$ is not a vertex. By no $\alpha\epsilon\cdots,\epsilon\thin\epsilon\cdots$, we know the AAD of $\thin\gamma\thin\gamma\thin$ is $\thin^{\epsilon}\gamma^{\alpha}\thin^{\alpha}\gamma^{\epsilon}\thin$. This implies no consecutive $\gamma\gamma\gamma$. 

We know $H\ne \alpha^2\beta\gamma$, and just proved that $\alpha^k\beta\gamma(k>2)$ is not a vertex. Then by no $\gamma\gamma\gamma$, we get $H=\thin^{\beta}\alpha^{\gamma}\thin^{\epsilon}\gamma^{\alpha}\thin^{\delta}\beta^{\alpha}\thin\cdots=\alpha\beta\gamma^k=\thin^{\beta}\alpha^{\gamma}\thin^{\epsilon}\gamma^{\alpha}\thin^{\delta}\beta^{\alpha}\thin^{\epsilon}\gamma^{\alpha}\thin^{\alpha}\gamma^{\epsilon}\thin$. This implies $\beta\epsilon\cdots$. By $\delta<\epsilon$ and $\beta+\delta+\epsilon=(\frac{7}{3}-\frac{4}{f})\pi>2\pi$, we know $\beta\epsilon\cdots$is not a vertex.

\subsubsection*{Subcase. $H=\alpha\beta\gamma^2$}

Since $\alpha\beta\gamma^2$ is a vertex, and $\alpha\ne\gamma$, we know $\alpha^2\beta\gamma$ is not a vertex. We also proved $\alpha^k\beta\gamma(k>2)$ is not a vertex. Therefore $\alpha\beta\gamma\cdots=\alpha\beta\gamma^k(k\ge 2)$. Again by $\alpha\beta\gamma^2$, we conclude $\alpha\beta\gamma\cdots=\alpha\beta\gamma^2$.

The angle sums of $\alpha\delta^2,\gamma\epsilon^2,\beta^3,\alpha\beta\gamma^2$ and the angle sum for pentagon imply
\[
	\alpha=\tfrac{16}{f}\pi,\,
	\beta=\tfrac{2}{3}\pi,\,
	\gamma=(\tfrac{2}{3}-\tfrac{8}{f})\pi,\,
	\delta=(1-\tfrac{8}{f})\pi,\,
	\epsilon=(\tfrac{2}{3}+\tfrac{4}{f})\pi.
\]
By $f\ge 48$, we get $\delta>\epsilon>\beta>\gamma>\alpha$. By $\alpha\delta^2,\gamma\epsilon^2$, and $\beta>\alpha>\gamma$, we know $\beta\cdots$ is a $\hat{b}$-vertex. 

By no $\beta\epsilon\cdots$, the vertex $H=\thin^{\beta}\alpha^{\gamma}\thin^{\epsilon}\gamma^{\alpha}\thin^{\delta}\beta^{\alpha}\thin\cdots=\alpha\beta\gamma^2=\thin^{\epsilon}\gamma^{\alpha}\thin^{\beta}\alpha^{\gamma}\thin^{\epsilon}\gamma^{\alpha}\thin^{\delta}\beta^{\alpha}\thin$. This implies a vertex $\thin^{\beta}\alpha^{\gamma}\thin^{\alpha}\beta^{\delta}\thin\cdots$. Since $\beta\cdots$ is a $\hat{b}$-vertex, and $\alpha\beta\gamma\cdots=\alpha\beta\gamma^2$, we get $\thin^{\beta}\alpha^{\gamma}\thin^{\alpha}\beta^{\delta}\thin\cdots=\alpha\beta\gamma^2,\alpha^k\beta^l$. 

If $\thin^{\beta}\alpha^{\gamma}\thin^{\alpha}\beta^{\delta}\thin\cdots=\alpha\beta\gamma^2$, then the AAD implies $\beta\epsilon\cdots,\delta\thin\epsilon\cdots,\epsilon\thin\epsilon\cdots$, a contradiction. If $\thin^{\beta}\alpha^{\gamma}\thin^{\alpha}\beta^{\delta}\thin\cdots=\alpha^k\beta^l$, then by no $\delta\thin\delta\cdots$, the vertex is $\thin^{\alpha}\beta^{\delta}\thin\cdots\thin^{\alpha}\beta^{\delta}\thin\alpha\thin\cdots$. By no $\beta\delta\cdots$, this implies a vertex $\gamma\delta\cdots$. By $\gamma\epsilon^2$ and $\delta>\epsilon$, however, $\gamma\delta\cdots$ is also not a vertex.

\subsubsection*{Case. $\alpha\beta\gamma$ is a vertex}

The angle sums of $\alpha\delta^2,\gamma\epsilon^2,\alpha\beta\gamma$ and the angle sum for pentagon imply
\[
\alpha=2\pi-2\delta,\,
\beta=\tfrac{8}{f}\pi,\,
\gamma=-\tfrac{8}{f}\pi+2\delta,\,
\epsilon=(1+\tfrac{4}{f})\pi-\delta.
\]

We use Lemma \ref{geometry7} to calculate the pentagon. First, we assume $\delta+\epsilon-\alpha\in (2{\bb Z}+1)\pi$. Then $\alpha=\frac{4}{f}\pi$, and we get
\[
\alpha=\tfrac{4}{f}\pi,\,
\beta=\tfrac{8}{f}\pi,\,
\gamma=(2-\tfrac{12}{f})\pi,\,
\delta=(1-\tfrac{2}{f})\pi,\,
\epsilon=\tfrac{6}{f}\pi.
\]
Substituting into \eqref{coolsaet_eq10}, we get
\[
16\cos^7 t-24\cos^5 t+7\cos^3 t+\cos t=0,\quad t=\tfrac{2}{f}\pi.
\]
The equation has no solution for even $f\ge 16$.

Next, we assume $\delta+\epsilon-\alpha\not\in (2{\bb Z}+1)\pi$, and substitute the angle values into \eqref{coolsaet_eq15} to get
\begin{align*}
0=&-4(64\cos^6 t-80\cos^4 t+16\cos^2 t+5)\cos^2 t\sin t\cos^3\delta\\
&+2(128 \cos^8 t-224\cos^6 t+96\cos^4 t+2\cos^2 t-1)\cos t\cos^2\delta\sin\delta\\
&+4(64 \cos^6 t-84\cos^4 t+20\cos^2 t+5)\cos^2 t\sin t\cos\delta\\
&-2(64\cos^8 t-120\cos^6 t+64\cos^4 t-6\cos^2 t-1)\cos t\sin\delta.
\end{align*}
This implicitly determines $\delta$ as a function of $\frac{2}{f}$. See Figure \ref{a2d_c2e_geometry}. The graph shows $\delta<(\frac{1}{2}+\frac{2}{f})\pi$. This implies $\delta<\epsilon$. Then by Lemma \ref{geometry1}, we get $\beta>\gamma$. This means $\delta<\frac{8}{f}\pi=\beta$.

\begin{figure}[htp]
\centering
\begin{tikzpicture}[>=latex,scale=1]

\node[scale=0.6] at (0,0) {\includegraphics{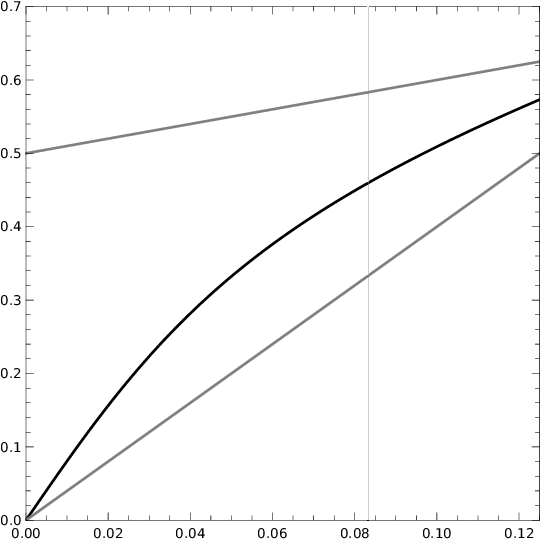}};

\node at (2.9,-2.5) {\scriptsize $\frac{2}{f}$};
\node at (0,0.6) {\scriptsize $\frac{\delta}{\pi}$};
\node at (0,1.95) {\scriptsize $\frac{1}{2}+\frac{2}{f}$};
\node at (0,-0.5) {\scriptsize $\frac{8}{f}$};

\end{tikzpicture}
\caption{Proposition \ref{a2d_c2e}: $\frac{\delta}{\pi}$ as function of $\frac{2}{f}$, for $\alpha\beta\gamma$.}
\label{a2d_c2e_geometry}
\end{figure}
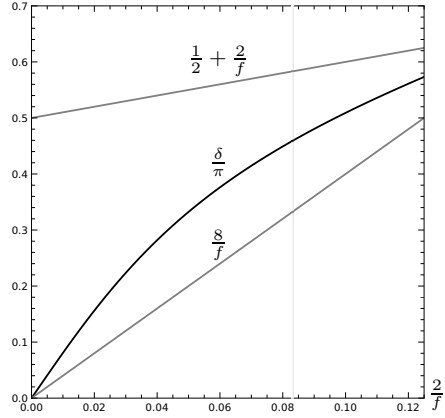

By $\gamma\epsilon^2$, we get $\epsilon<\pi$. We also have $\gamma<\beta<\pi$. Then by Lemma \ref{geometry6}, we get $\beta+2\delta>\pi$. By $\delta<\frac{8}{f}\pi=\beta$, this implies $16\le f<24$, the graph shows $\delta>\frac{8}{f}\pi$, a contradiction. 
\end{proof}

\begin{proposition}\label{b2d_c2e}
There is no tiling, such that $\alpha,\beta,\gamma$ have distinct values, and $\beta\delta^2,\gamma\epsilon^2$ are vertices.
\end{proposition}

\begin{proof}
By Proposition \ref{ade_b2d_c2e}, we know $\alpha\delta\epsilon$ is not a vertex. Then by distinct $\alpha,\beta,\gamma$ values, we know $\beta\delta^2,\gamma\epsilon^2$ are all the degree $3$ $b$-vertices. This implies a special companion pair is matched, and a $b$-vertex $\alpha\cdots$ has high degree.

By Proposition  \ref{b2d_c2e_abc}, we know the pentagon is symmetric in case $\alpha\beta\gamma$ is a vertex. Therefore $\alpha\beta\gamma$ is not a vertex.

By applying Lemma \ref{degree3} to $\delta,\epsilon$, we know there is a degree $3$ $\hat{b}$-vertex. Then by distinct $\alpha,\beta,\gamma$, and no $\alpha\beta\gamma$, and up to the exchange symmetry $(\beta,\delta)\leftrightarrow(\gamma,\epsilon)$, we may assume one of $\alpha^3,\alpha^2\beta,\alpha\beta^2,\beta^3,\beta^2\gamma$ is a vertex. 

By $\beta\delta^2,\gamma\epsilon^2$, we get $\delta,\epsilon<\pi$. Then by Lemma \ref{geometry6}, we get $\alpha+2\beta>\pi$ and $\alpha+2\gamma>\pi$.

\subsubsection*{Case. $\alpha^2\beta$ is a vertex}

By distinct $\alpha,\beta,\gamma$ values, and all the degree $3$ $b$-vertices $\beta\delta^2,\gamma\epsilon^2$, and no $\alpha\beta\gamma$, we know $\beta\delta^2,\gamma\epsilon^2,\alpha^2\beta,\alpha\gamma^2,\beta^2\gamma,\gamma^3$ are all the degree $3$ vertices. 

By $\beta\delta^2,\alpha^2\beta$, we get $\alpha=\delta$. If $\gamma^3$ is a vertex, then by $\gamma\epsilon^2$, we get $\gamma=\epsilon$. Then by Lemma \ref{geometry11}, we get $a=b$, a contradiction. 

The the angle sums of $\beta\delta^2,\gamma\epsilon^2,\alpha^2\beta$, and the angle sum of one of $\alpha\gamma^2,\beta^2\gamma$, and the angle sum for pentagon imply
\begin{align*}
\alpha\gamma^2 & \colon
	\alpha=\delta=(2-\tfrac{16}{f})\pi,\,
	\beta=(\tfrac{32}{f}-2)\pi,\,
	\gamma=\tfrac{8}{f}\pi,\,
	\epsilon=(1-\tfrac{4}{f})\pi. \\
\beta^2\gamma &\colon
	\alpha=\delta=(\tfrac{1}{2}+\tfrac{2}{f})\pi,\,
	\beta=\epsilon=(1-\tfrac{4}{f})\pi,\,
	\gamma=\tfrac{8}{f}\pi.
\end{align*}
For $\alpha\gamma^2$, we get $\beta\le 0$, a contradiction. For $\beta^2\gamma$, we get $\pi>\beta=\epsilon>\alpha=\delta>\gamma$, contradicting Lemma \ref{geometry4}.  

Therefore $\beta\delta^2,\gamma\epsilon^2,\alpha^2\beta$ are all the degree $3$ vertices. This implies that, in any tile, one of $\alpha,\gamma\cdots$ has high degree. By Lemma \ref{special_tile}, this implies $f\ge 24$.

The angle sums of $\beta\delta^2,\gamma\epsilon^2,\alpha^2\beta$ and the angle sum for pentagon imply
\[
\alpha=\delta=\pi-\tfrac{1}{2}\beta,\,
\gamma=\tfrac{8}{f}\pi,\,
\epsilon=(1-\tfrac{4}{f})\pi.
\]
We have $\gamma<\epsilon$.

By $\alpha+2\gamma>\pi$ and $f\ge 24$, we get $\alpha>(1-\tfrac{16}{f})\pi\ge \gamma$, and $\alpha+\epsilon=\delta+\epsilon>\gamma+\epsilon>\pi$. By $\delta+\epsilon>\pi$ and $2\epsilon>\pi$, we know $\epsilon\thin\epsilon\cdots$ is not a vertex.

\subsubsection*{Subcase. $\beta>\epsilon$}

By $\gamma<\epsilon$ and Lemma \ref{geometry1}, we get $\delta<\epsilon$. Then we get $\beta>\epsilon>\alpha=\delta>\gamma$. By Lemma \ref{geometry4}, this implies $\beta\ge \pi$. Then $\beta^2\cdots$ is not a vertex. This implies $\delta\thin\delta\cdots$ is not a vertex.

By $\gamma\epsilon^2$ and $\gamma<\alpha,\beta,\delta,\epsilon$, we get $\epsilon^2\cdots=\gamma\epsilon^2$.

By no $\delta\thin\delta\cdots$, and $\gamma+2\delta=2\alpha+\gamma>\alpha+2\gamma>\pi$, and $\gamma<\alpha,\beta,\delta,\epsilon$, we know a $\delta^2$-fan has value $>\pi$. Then by $2\epsilon>\delta+\epsilon>\pi$, we know a $b$-vertex has exactly one fan. In particular, we know $\delta\thin\delta\cdots,\delta\thin\epsilon\cdots,\epsilon\thin\epsilon\cdots$ are not vertices.

By $\beta\delta^2$ and $\delta<\epsilon$, we know a $b$-vertex $\beta\cdots=\beta\delta^2$. By $\alpha^2\beta$, and $\beta>\alpha>\gamma$, and no $\beta^2\cdots$, a $\hat{b}$-vertex $\beta\cdots=\alpha^2\beta,\alpha\beta\gamma^k(k\ge 2),\beta\gamma^k(k\ge 3)$. By $\alpha+2\gamma>\pi$ and $\beta\ge\pi$, we know $\alpha\beta\gamma^k$ is not a vertex. Therefore $\beta\cdots=\beta\delta^2,\alpha^2\beta,\beta\gamma^k(k\ge 3)$.

The vertex $\alpha\epsilon\cdots$ has high degree and has no $\beta$. By $\epsilon^2\cdots=\gamma\epsilon^2$, and $b$-vertex having exactly one fan, we get $\alpha\epsilon\cdots=\alpha\delta\epsilon\cdots$, with no $\delta,\epsilon$ in the remainder. By $\alpha+2\gamma>\pi$ and $\delta+\epsilon>\pi$, we get $R(\alpha\delta\epsilon)<2\gamma<2\alpha$. Then we get $\alpha\epsilon\cdots=\alpha\gamma\delta\epsilon,\alpha^2\delta\epsilon$. 

By no $\epsilon\thin\epsilon\cdots$, the AAD of consecutive $\gamma\gamma\gamma$ implies a vertex $\thin^{\beta}\alpha^{\gamma}\thin^{\gamma}\epsilon^{\delta}\thick\cdots$. By $\alpha\epsilon\cdots=\alpha\gamma\delta\epsilon,\alpha^2\delta\epsilon$, the vertex is $\thick^{\epsilon}\delta^{\beta}\thin\alpha\thin^{\beta}\alpha^{\gamma}\thin^{\gamma}\epsilon^{\delta}\thick,\thick^{\epsilon}\delta^{\beta}\thin\gamma\thin^{\beta}\alpha^{\gamma}\thin^{\gamma}\epsilon^{\delta}\thick$. Then we get $\beta^2\cdots,\beta\epsilon\cdots$, a contradiction. 

Therefore we do not have consecutive $\gamma\gamma\gamma$. This implies $\beta\gamma^k$ is not a vertex, and $\beta\cdots=\beta\delta^2,\alpha^2\beta$. Then we get
\[
f=\#\beta=\#\beta\delta^2+\#\alpha^2\beta
\le \tfrac{1}{2}\#\delta+\tfrac{1}{2}\#\alpha=f.
\]
Since both sides are equal, this implies the only vertex besides $\beta\delta^2,\alpha^2\beta$ is $\gamma^k\epsilon^l$. By $\gamma\epsilon^2$, and $\gamma<\epsilon$, and no consecutive $\gamma\gamma\gamma$, we get $\gamma^k\epsilon^l=\gamma\epsilon^2$. Then there is no high degree vertex, a contradiction.

\subsubsection*{Subcase. $\beta<\epsilon$}

Suppose $\beta<\gamma$ and $\delta>\epsilon$. Then we get $\alpha=\delta>\epsilon>\gamma>\beta$. By $\beta\delta^2,\gamma\epsilon^2$, and $\delta>\epsilon>\frac{1}{2}\pi$, and $\alpha>\gamma>\beta$, we know $\beta\delta^2,\gamma\epsilon^2,\beta^k\delta\epsilon(k\ge 2),\beta^k\epsilon^2(k\ge 2)$ are all the $b$-vertices. Then $\alpha\delta\cdots,\gamma\delta\cdots,\delta\thin\delta\cdots$ are not vertices, and the AAD implies $\beta^k\delta\epsilon(k\ge 2),\beta^k\epsilon^2(k\ge 2)$ are not vertices. Therefore $\beta\delta^2,\gamma\epsilon^2$ are the only $b$-vertices. This implies no consecutive $\theta\theta\theta$, for $\theta=\beta,\gamma$. On the other hand, the AAD of $\gamma\epsilon^2$ implies a vertex $\alpha\gamma\cdots$. We know this is a $\hat{b}$-vertex. Then by $\alpha^2\beta$ and $\beta<\gamma$, we get $\alpha\gamma\cdots=\alpha\beta^k\gamma^l(k+l\ge 3)$, contradicting no $\theta\theta\theta$.

By Lemma \ref{geometry1}, therefore, we get $\beta>\gamma$ and $\delta<\epsilon$. By $\alpha>\gamma$, we have two possibilities $\epsilon>\alpha=\delta>\beta>\gamma$ and $\epsilon>\beta>\alpha=\delta>\gamma$. 

By $\gamma\epsilon^2$ and $\gamma<\alpha,\beta,\gamma,\delta$, we get $\epsilon^2\cdots=\gamma\epsilon^2$.

By $\alpha^2\beta$ and $\beta<\epsilon=(1-\frac{4}{f})\pi$, we get $\alpha=\delta>(\frac{1}{2}+\frac{2}{f})\pi$. Then by $\epsilon>\delta>\frac{1}{2}\pi$, we know a vertex has at most one $b$-edge, and $\delta\thin\delta\cdots,\delta\thin\epsilon\cdots,\epsilon\thin\epsilon\cdots$ are not vertices.

\subsubsection*{Subsubcase. $\epsilon>\alpha=\delta>\beta>\gamma$}

By $\beta\delta^2$, and $\alpha>\beta>\gamma$, and $\delta<\epsilon$, and at most one $b$-edge, we get $\delta^2\cdots=\beta\delta^2,\gamma^k\delta^2(k\ge 2)$, and $\delta\epsilon\cdots=\gamma^k\delta\epsilon(k\ge 2)$. Combined with $\epsilon^2\cdots=\gamma\epsilon^2$, we know $\alpha\delta\cdots,\alpha\epsilon\cdots,\beta\epsilon\cdots$ are not vertices. By no $\alpha\epsilon\cdots,\epsilon\thin\epsilon\cdots$, we know the AAD of $\thin\gamma\thin\gamma\thin$ is $\thin^{\epsilon}\gamma^{\alpha}\thin^{\alpha}\gamma^{\epsilon}\thin$. This implies no consecutive $\gamma\gamma\gamma$. Then we get $k=2$ in $\gamma^k\delta^2,\gamma^k\delta\epsilon$. Then the vertices are $\thick^{\epsilon}\delta^{\beta}\thin^{\epsilon}\gamma^{\alpha}\thin^{\alpha}\gamma^{\epsilon}\thin\theta\thick$, with $\theta=\delta,\epsilon$, contradicting no $\beta\epsilon\cdots$. We conclude $\beta\delta^2,\gamma\epsilon^2$ are all the $b$-vertices. This implies no consecutive $\theta\theta\theta$, for $\theta=\beta,\gamma$. 

Therefore $\alpha\cdots$ is a $\hat{b}$-vertex. Then by $\alpha^2\beta$ and $\alpha>\beta$, a $\hat{b}$-vertex is $\alpha^2\beta,\alpha\beta^k\gamma^l(k+l\ge 3),\beta^k\gamma^l(k+l\ge 4)$. However, by no consecutive $\theta\theta\theta$, we know $\alpha\beta^k\gamma^l,\beta^k\gamma^l$ are not vertices. Therefore $\beta\delta^2,\gamma\epsilon^2,\alpha^2\beta$ are all the vertices. In particular, there is no high degree vertex, a contradiction. 

\subsubsection*{Subsubcase. $\epsilon>\beta>\alpha=\delta>\gamma$}

By $\beta\delta^2$ and $\delta<\epsilon$, a $b$-vertex $\beta\cdots=\beta\delta^2$, and $\beta\epsilon\cdots$ is not a vertex. Then by $\epsilon^2\cdots=\gamma\epsilon^2$, this implies $\alpha\epsilon\cdots=\alpha\delta\epsilon\cdots$, with no $\beta$ in the remainder. Then by $\alpha,\delta,\epsilon>\frac{1}{2}\pi$, and $\alpha+2\gamma>\pi$, and all degree $3$ vertices $\beta\delta^2,\gamma\epsilon^2,\alpha^2\beta$, we get $\alpha\epsilon\cdots=\alpha\gamma\delta\epsilon$. The angle sum of $\alpha\gamma\delta\epsilon$ and the angle sum for pentagon imply $\beta=(1+\frac{4}{f})\pi>\pi$, a contradiction. Therefore $\alpha\epsilon\cdots$ is not a vertex.

By $\alpha^2\beta$ and $\beta>\alpha>\gamma$, a $\hat{b}$-vertex $\beta\cdots=\alpha^2\beta,\alpha\beta\gamma^k,\beta^2\gamma^k,\beta\gamma^k$. By no $\alpha\epsilon\cdots,\epsilon\thin\epsilon\cdots$, we know the AAD of $\thin\gamma\thin\gamma\thin$ is $\thin^{\epsilon}\gamma^{\alpha}\thin^{\alpha}\gamma^{\epsilon}\thin$. By no $\alpha\epsilon\cdots,\delta\thin\epsilon\cdots,\epsilon\thin\epsilon\cdots$, this implies no consecutive $\beta\gamma\gamma,\gamma\gamma\gamma$. Therefore $\beta\gamma^k$ is not a vertex, and $\alpha\beta\gamma^k=\thin\alpha\thin\gamma\thin\beta\thin\gamma\thin$, and $\beta^2\gamma^k=\thin\beta\thin\gamma\thin\beta\thin\gamma\thin$. By no $\alpha\epsilon\cdots,\delta\thin\epsilon\cdots$, we know the AAD of $\thin\gamma\thin\beta\thin\gamma\thin$ is  $\thin^{\epsilon}\gamma^{\alpha}\thin^{\alpha}\beta^{\delta}\thin^{\alpha}\gamma^{\epsilon}\thin$. Then the AADs of $\thin\alpha\thin\gamma\thin\beta\thin\gamma\thin,\thin\beta\thin\gamma\thin\beta\thin\gamma\thin$ imply $\alpha\epsilon\cdots,\beta\epsilon\cdots$, a contradiction. 

We conclude $\beta\cdots=\beta\delta^2,\alpha^2\beta$. By the same argument for the subcase $\beta>\epsilon$, we get a contradiction.

\subsubsection*{Case. $\beta^2\gamma$ is a vertex}

Suppose $\beta\delta^2,\gamma\epsilon^2,\beta^2\gamma$ are all the degree $3$ vertices. Then in a special companion pair, both vertices $\alpha\cdots$ have high degree. This implies both $\delta\cdots,\epsilon\cdots$ have degree $3$, and they are $\beta\delta^2,\gamma\epsilon^2$. This implies one of the two $\beta\cdots$ is $\alpha\beta\cdots$, and one of the two $\gamma\cdots$ is $\alpha\gamma\cdots$. Both have high degree. Then the special companion pair has at least four high degree vertices, a contradiction. 

Therefore there is a degree $3$ vertex besides $\beta\delta^2,\gamma\epsilon^2,\beta^2\gamma$. By distinct $\alpha,\beta,\gamma$ values, and all the degree $3$ $b$-vertices $\beta\delta^2,\gamma\epsilon^2$, and no $\alpha\beta\gamma$, these degree $3$ vertices are $\alpha^3,\alpha^2\beta,\alpha\gamma^2$. We already proved there is no tiling with $\alpha^2\beta$. 

If $\alpha^3$ is a vertex, then by distinct $\alpha,\beta,\gamma$ values, and no $\alpha\beta\gamma$, we know $\alpha\beta\cdots$ has high degree. If $\beta\cdots,\delta\cdots$ in a tile have degree $3$, then by all the degree $3$ $b$-vertices $\beta\delta^2,\gamma\epsilon^2$, and $\alpha\beta\cdots$ having high degree, we get $\delta\cdots=\beta\cdots=\beta\delta^2$. This implies $\alpha\cdots=\alpha\beta\cdots$, which has high degree. Then by Lemma \ref{special_tile}, we get $f\ge 24$. On the other hand, the angle sums of $\beta\delta^2,\gamma\epsilon^2,\beta^2\gamma,\alpha^3$ and the angle sum for pentagon imply $\gamma=(\tfrac{16}{f}-\tfrac{2}{3})\pi$, contradicting $f\ge 24$. 

Therefore $\alpha\gamma^2$ is a vertex, and $\beta\delta^2,\gamma\epsilon^2,\beta^2\gamma,\alpha\gamma^2$ are all the degree $3$ vertices. The angle sums of the vertices and the angle sum for pentagon imply
\[
	\alpha=(\tfrac{2}{7}+\tfrac{32}{7f})\pi,\,
	\beta=\epsilon=(\tfrac{4}{7}+\tfrac{8}{7f})\pi,\,
	\gamma=(\tfrac{6}{7}-\tfrac{16}{7f})\pi,\,
	\delta=(\tfrac{5}{7}-\tfrac{4}{7f})\pi.
\]
We have $\alpha<\beta=\epsilon<\delta<\gamma$. By $\delta>\epsilon>\frac{1}{2}\pi$, there is at most one $b$-edge at a vertex, and $\delta\thin\delta\cdots,\delta\thin\epsilon\cdots,\epsilon\thin\epsilon\cdots$ are not vertices.

By $\alpha\gamma^2$ and $\alpha<\beta,\gamma,\delta,\epsilon$, we get $\gamma^2\cdots=\alpha\gamma^2$. By $\gamma\epsilon^2$ and $\beta=\epsilon<\delta$, we get $R(\gamma\delta)<R(\beta\gamma)=R(\gamma\epsilon)=\beta=\epsilon<2\alpha,\gamma,\delta$. This implies $\beta\gamma\cdots=\beta^2\gamma$, and $\gamma\epsilon\cdots=\gamma\epsilon^2$, and $\gamma\delta\cdots$ is not a vertex. We conclude $\gamma\cdots=\alpha\gamma^2,\beta^2\gamma,\gamma\epsilon^2,\alpha^k\gamma$. 

By $\beta=\epsilon<\delta$, we get $R(\beta\delta\epsilon)<R(\beta\epsilon^2)=R(\beta^3)=(\frac{2}{7}-\frac{24}{f})\pi<\alpha<\beta,\gamma,\delta,\epsilon$. This implies $\beta\epsilon\cdots,\beta^3\cdots$ are not vertices. Then by $\beta\delta^2$, we know a $b$-vertex $\beta\cdots=\beta\delta^2$. Then $\beta^2\cdots$ is a $\hat{b}$-vertex, with no $\beta$ in the remainder. Then by $\beta^2\gamma$ and $\gamma<3\alpha$, we get $\beta^2\cdots=\beta^2\gamma,\alpha^2\beta^2$. It remains to consider $\hat{b}$-vertex $\beta\cdots$, with no $\beta,\gamma$ in the remainder. The vertex is $\alpha^k\beta$. We conclude $\beta\cdots=\beta\delta^2,\beta^2\gamma,\alpha^2\beta^2,\alpha^k\beta$.

The AAD $\thin^{\beta}\alpha^{\gamma}\thin^{\gamma}\alpha^{\beta}\thin$ implies a vertex $\thin^{\epsilon}\gamma^{\alpha}\thin^{\alpha}\gamma^{\epsilon}\thin\cdots=\alpha\gamma^2=\thin^{\alpha}\gamma^{\epsilon}\thin\alpha\thin^{\epsilon}\gamma^{\alpha}\thin$, contradicting no $\beta\epsilon\cdots$. The AAD $\thin^{\gamma}\alpha^{\beta}\thin^{\gamma}\alpha^{\beta}\thin$ implies a vertex $\thin^{\delta}\beta^{\alpha}\thin^{\alpha}\gamma^{\epsilon}\thin\cdots=\beta^2\gamma=\thin^{\alpha}\gamma^{\epsilon}\thin\beta\thin^{\delta}\beta^{\alpha}\thin$, contradicting no $\delta\thin\delta\cdots,\delta\thin\epsilon\cdots$. Therefore the AAD of $\thin\alpha\thin\alpha\thin$ is $\thin^{\gamma}\alpha^{\beta}\thin^{\beta}\alpha^{\gamma}\thin$. The AAD implies a vertex $\thin^{\delta}\beta^{\alpha}\thin^{\alpha}\beta^{\delta}\thin\cdots=\beta^2\gamma,\alpha^2\beta^2$. If the vertex is $\beta^2\gamma$, then it is $\thin^{\alpha}\beta^{\delta}\thin\gamma\thin^{\delta}\beta^{\alpha}\thin$, contradicting no $\delta\thin\epsilon\cdots$. If the vertex is $\alpha^2\beta^2$, then by the AAD $\thin^{\gamma}\alpha^{\beta}\thin^{\beta}\alpha^{\gamma}\thin$ of $\thin\alpha\thin\alpha\thin$, it is $\thin^{\alpha}\beta^{\delta}\thin^{\gamma}\alpha^{\beta}\thin^{\beta}\alpha^{\gamma}\thin^{\delta}\beta^{\alpha}\thin$, contradicting no $\gamma\delta\cdots$. 

Therefore $\alpha\thin\alpha\cdots$ is not a vertex. This implies $\beta\cdots=\beta\delta^2,\beta^2\gamma,\alpha^2\beta^2$, and $\gamma\cdots=\alpha\gamma^2,\beta^2\gamma,\gamma\epsilon^2$. 
Then a $b$-vertex $\alpha\cdots$ has no $\beta,\gamma$. Then by at most one $b$-edge at a vertex, and no $\alpha\thin\alpha\cdots$, and all the degree $3$ vertices $\beta\delta^2,\gamma\epsilon^2,\beta^2\gamma,\alpha\gamma^2$, we know $\alpha\cdots$ is a $\hat{b}$-vertex. Then by the lists of $\beta\cdots,\gamma\cdots$, and no $\alpha\thin\alpha\cdots$, we get $\alpha\cdots=\alpha\gamma^2,\alpha^2\beta^2$. 

By no $\alpha\thin\alpha\cdots,\gamma\delta\cdots$, we get $\alpha^2\beta^2=\thin\alpha\thin\beta\thin\alpha\thin\beta\thin=\thin^{\beta}\alpha^{\gamma}\thin^{\alpha}\beta^{\delta}\thin^{\beta}\alpha^{\gamma}\thin^{\alpha}\beta^{\delta}\thin$. This implies a vertex $\thin^{\gamma}\alpha^{\beta}\thin^{\alpha}\gamma^{\epsilon}\thin\cdots=\alpha\gamma^2=\thin^{\alpha}\gamma^{\epsilon}\thin\gamma\thin^{\gamma}\alpha^{\beta}\thin$, contradicting no $\alpha\epsilon\cdots,\epsilon\thin\epsilon\cdots$.

Therefore $\alpha^2\beta^2$ is not a vertex, and $\alpha\cdots=\alpha\gamma^2$. Applying the counting lemma to $\alpha,\gamma$, we get a contradiction.

\subsubsection*{Case. $\beta^3$ is a vertex}

By distinct $\alpha,\beta,\gamma$ values, and all the degree $3$ $b$-vertices $\beta\delta^2,\gamma\epsilon^2$, and no $\alpha\beta\gamma$, we know $\beta\delta^2,\gamma\epsilon^2,\beta^3,\alpha^2\gamma,\alpha\gamma^2$ are all the degree $3$ vertices.  

If $\alpha^2\gamma$ is a vertex, then we get $\beta=\delta$ and $\alpha=\epsilon$. By Lemma \ref{geometry11}, this implies $a=b$, a contradiction. 

Suppose $\beta\delta^2,\gamma\epsilon^2,\beta^3$ are all the degree $3$ vertices. Then in a special companion pair, both vertices $\alpha\cdots$ have high degree. This implies both $\delta\cdots,\epsilon\cdots$ have degree $3$, and they are $\beta\delta^2,\gamma\epsilon^2$. This implies one of the two $\beta\cdots$ is $\alpha\beta\cdots$, and one of the two $\gamma\cdots$ is $\alpha\gamma\cdots$. Both have high degree. Then the special companion pair has at least four high degree vertices, a contradiction. 

Therefore $\alpha\gamma^2$ is a vertex, and $\beta\delta^2,\gamma\epsilon^2,\beta^3,\alpha\gamma^2$ are all the degree $3$ vertices. The angle sums of the vertices and the angle sum for pentagon imply
\[
\alpha=(\tfrac{2}{9}+\tfrac{16}{3f})\pi,\,
\beta=\delta=\tfrac{2}{3}\pi,\,
\gamma=(\tfrac{8}{9}-\tfrac{8}{3f})\pi,\,
\epsilon=(\tfrac{5}{9}+\tfrac{4}{3f})\pi.
\]
We have $\alpha<\epsilon<\beta=\delta<\gamma$. By $\delta>\epsilon>\frac{1}{2}\pi$, there is at most one $b$-edge at a vertex, and $\delta\thin\delta\cdots,\delta\thin\epsilon\cdots,\epsilon\thin\epsilon\cdots$ are not vertices. 

By $\gamma\epsilon^2\cdots$ and $\delta>\epsilon$, we know $\gamma\delta\cdots$ is not a vertex. Then by at most one $b$-edge at a vertex, we get $\delta^2\cdots=\beta\delta^2,\alpha^k\delta^2(k\ge 2)$.

By $\alpha\gamma^2$ and $\alpha<\beta,\gamma,\delta,\epsilon$, we get $\gamma^2\cdots=\alpha\gamma^2$. By $R(\beta\gamma)=(\tfrac{4}{9}+\tfrac{8}{3f})\pi<2\alpha,\beta,\gamma,\delta,\epsilon$, and no $\alpha\beta\gamma$, we know $\beta\gamma\cdots$ is not a vertex.

By $R(\beta\delta\epsilon)<R(\beta\epsilon^2)=(\tfrac{2}{9}-\tfrac{8}{3f})\pi<\alpha<\beta,\gamma,\delta,\epsilon$, we know $\beta\epsilon\cdots$ is not a vertex, and a $b$-vertex $\beta\cdots=\beta\delta^2$. Then by $\beta^3$, and $\beta<3\alpha$, and $\alpha\ne\beta$, and no $\beta\gamma\cdots$, we get $\beta\cdots=\beta\delta^2,\beta^3,\alpha^2\beta^2,\alpha^k\beta$. 

By no $\beta\gamma\cdots$, we do not have $\thin^{\beta}\alpha^{\gamma}\thin^{\beta}\alpha^{\gamma}\thin$. The AAD $\thin^{\beta}\alpha^{\gamma}\thin^{\gamma}\alpha^{\beta}\thin$ implies a vertex $\thin^{\epsilon}\gamma^{\alpha}\thin^{\alpha}\gamma^{\epsilon}\thin\cdots=\alpha\gamma^2=\thin^{\alpha}\gamma^{\epsilon}\thin\alpha\thin^{\epsilon}\gamma^{\alpha}\thin$, contradicting no $\beta\epsilon\cdots$. Therefore the AAD of $\thin\alpha\thin\alpha\thin$ is $\thin^{\gamma}\alpha^{\beta}\thin^{\beta}\alpha^{\gamma}\thin$. This implies a vertex $\thin^{\delta}\beta^{\alpha}\thin^{\alpha}\beta^{\delta}\thin\cdots=\beta^3,\alpha^2\beta^2$. If the vertex is $\beta^3$, then it is $\thin^{\alpha}\beta^{\delta}\thin\beta\thin^{\delta}\beta^{\alpha}\thin$, contradicting no $\delta\thin\delta\cdots$. If the vertex is $\alpha^2\beta^2$, then by the AAD $\thin^{\gamma}\alpha^{\beta}\thin^{\beta}\alpha^{\gamma}\thin$ of $\thin\alpha\thin\alpha\thin$, it is $\thin^{\alpha}\beta^{\delta}\thin^{\gamma}\alpha^{\beta}\thin^{\beta}\alpha^{\gamma}\thin^{\delta}\beta^{\alpha}\thin$, contradicting no $\gamma\delta\cdots$. 

Therefore $\alpha\thin\alpha\cdots$ is not a vertex. This implies $\beta\cdots=\beta^3,\beta\delta^2,\alpha^2\beta^2$, and $\delta^2\cdots=\beta\delta^2$. Then by at most one $b$-edge at a vertex, a $b$-vertex $\alpha\cdots=\alpha\delta\epsilon\cdots,\alpha\epsilon^2\cdots$, with no $\beta,\delta,\epsilon$ in the remainder. By $\gamma\epsilon^2$ and $\delta>\epsilon$, the remainder has no $\gamma$. Then by no $\alpha\delta\epsilon,\alpha\epsilon^2$, we know a $b$-vertex $\alpha\cdots=\alpha^k\delta\epsilon,\alpha^k\epsilon^2$, with $k\ge 2$. Then by no $\alpha\thin\alpha\cdots$, we know $\alpha\cdots$ is a $\hat{b}$-vertex. Then by the list of $\beta\cdots$, we get $\alpha\cdots=\alpha^2\beta^2,\alpha^k\gamma^l(k\ge 1)$.

By $\gamma^2\cdots=\alpha\gamma^2$, we get $\alpha^k\gamma^l(k\ge 1)=\alpha\gamma^2,\alpha^k,\alpha^k\gamma$. By no $\alpha\thin\alpha\cdots$, we know $\alpha^k,\alpha^k\gamma$ are not vertices. Therefore $\alpha^k\gamma^l(k\ge 1)=\alpha\gamma^2$.

By no $\alpha\thin\alpha\cdots,\gamma\delta\cdots$, we get $\alpha^2\beta^2=\thin\alpha\thin\beta\thin\alpha\thin\beta\thin=\thin^{\beta}\alpha^{\gamma}\thin^{\alpha}\beta^{\delta}\thin^{\beta}\alpha^{\gamma}\thin^{\alpha}\beta^{\delta}\thin$. This implies a vertex $\thin^{\gamma}\alpha^{\beta}\thin^{\alpha}\gamma^{\epsilon}\thin\cdots=\alpha\gamma^2=\thin^{\alpha}\gamma^{\epsilon}\thin\gamma\thin^{\gamma}\alpha^{\beta}\thin$, contradicting no $\alpha\epsilon\cdots,\epsilon\thin\epsilon\cdots$.

We conclude $\alpha\cdots=\alpha\gamma^2$. Applying the counting lemma to $\alpha,\gamma$, we get a contradiction.

\subsubsection*{Case. $\alpha^3$ is a vertex}

By distinct $\alpha,\beta,\gamma$ values, and all the degree $3$ $b$-vertices $\beta\delta^2,\gamma\epsilon^2$, and no $\alpha\beta\gamma$, we know $\beta\delta^2,\gamma\epsilon^2,\alpha^3,\beta^2\gamma,\beta\gamma^2$ are all the degree $3$ vertices. We already proved there is no tiling with $\beta^2\gamma$. By the exchange symmetry $(\beta,\delta)\leftrightarrow(\gamma,\epsilon)$, there is also no tiling with $\beta\gamma^2$. Therefore $\beta\delta^2,\gamma\epsilon^2,\alpha^3$ are all the degree $3$ vertices.

A special companion pair is matched. This determines $T_1,T_2$ in Figure \ref{b2d_c2eD}. By all the degree $3$ vertices $\beta\delta^2,\gamma\epsilon^2,\alpha^3$, if $\alpha\cdots$ in a tile has degree $3$, then we get $\alpha\cdots=\alpha^3$, and both $\beta\cdots,\gamma\cdots$ have high degree. Therefore in the special  companion pair, one of $\alpha_1\cdots,\alpha_2\cdots$ has high degree. If $\delta_1\delta_2\cdots$ has degree $3$, then the vertex is $\beta\delta^2$. This implies one of $\beta_1\cdots,\beta_2\cdots$ is $\alpha\beta\cdots$, which has high degree. By the same reason, if $\epsilon_1\epsilon_2\cdots$ has degree $3$, then one of $\gamma_1\cdots,\gamma_2\cdots$ is the high degree vertex $\alpha\gamma\cdots$. These facts together imply that both shared vertices in the special companion pair have degree $3$. Moreover, since both $\alpha\beta\cdots,\alpha\gamma\cdots$ have high degree, we know that, up to the vertical flip symmetry, $T_3,T_4$ are arranged as indicated, and $\alpha_3\beta_1\cdots,\alpha_4\gamma_1\cdots$, and $\alpha_2\cdots=\thin^{\alpha}\beta^{\delta}\thin^{\beta}\alpha^{\gamma}\thin^{\epsilon}\gamma^{\alpha}\thin\cdots$ are all the high degree vertices in the special companion pair. 

\begin{figure}[htp]
\centering
\begin{tikzpicture}[>=latex,scale=1]

\foreach \a in {1,-1}
\draw[xscale=\a]
	(0,1.1) -- (0.5,0.7) -- (0.5,-0.7) -- (0,-1.1) -- (0.1,-1.4)
	(0.5,0.7) -- (1.3,0.7) -- (1.3,-0.7) -- (0.5,-0.7)
	(0,1.1) -- (0,1.5);

\draw[line width=1.2]
	(-0.5,0) -- (0.5,0)
	(0.5,-0.7) -- (1.3,-0.7)
	(-0.5,-0.7) -- (-1.3,-0.7);

\node at (0.2,1.2) {\small $\alpha$};
\node at (-0.2,1.2) {\small $\alpha$};

\node at (0.3,-0.6) {\small $\gamma$};
\node at (0,-0.85) {\small $\alpha$};
\node at (0.3,-0.2) {\small $\epsilon$};	
\node at (-0.3,-0.2) {\small $\delta$}; 
\node at (-0.3,-0.6) {\small $\beta$};

\node at (0.3,0.6) {\small $\gamma$};
\node at (0,0.85) {\small $\alpha$};
\node at (0.3,0.2) {\small $\epsilon$};	
\node at (-0.3,0.2) {\small $\delta$}; 
\node at (-0.3,0.6) {\small $\beta$};

\node at (1.1,0.5) {\small $\beta$};
\node at (1.1,-0.5) {\small $\delta$};
\node at (0.7,0.5) {\small $\alpha$};
\node at (0.7,-0.5) {\small $\epsilon$};
\node at (0.7,0) {\small $\gamma$};

\node at (-1.1,0.5) {\small $\gamma$};
\node at (-1.1,-0.5) {\small $\epsilon$};
\node at (-0.7,0.5) {\small $\alpha$};
\node at (-0.7,-0.5) {\small $\delta$};
\node at (-0.7,0) {\small $\beta$};

\node at (0.6,-0.9) {\small $\epsilon$};
\node at (0.2,-1.2) {\small $\gamma$};

\node at (-0.6,-0.9) {\small $\delta$};
\node at (-0.2,-1.2) {\small $\beta$};

\node[inner sep=0.5,draw,shape=circle] at (0,0.4) {\small $1$};
\node[inner sep=0.5,draw,shape=circle] at (0,-0.4) {\small $2$};
\node[inner sep=0.5,draw,shape=circle] at (-1.05,0) {\small $3$};
\node[inner sep=0.5,draw,shape=circle] at (1.05,0) {\small $4$};

\end{tikzpicture}
\caption{Proposition \ref{b2d_c2e}: $\alpha^3$ is a vertex.}
\label{b2d_c2eD}
\end{figure}
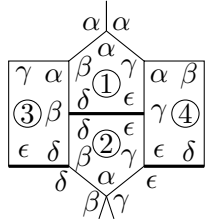

The angle sums of $\beta\delta^2,\gamma\epsilon^2,\alpha^3$ and the angle sum for pentagon imply
\[
\alpha=\tfrac{2}{3}\pi,\;
\beta+\gamma=(\tfrac{2}{3}+\tfrac{8}{f})\pi,\;
\delta+\epsilon=(\tfrac{5}{3}-\tfrac{4}{f})\pi.
\]

Up to the exchange of $(\beta,\delta)$ with $(\gamma,\delta)$, we may assume $\beta>\gamma$ and $\delta<\epsilon$. By $
\alpha+\delta+\epsilon>2\pi$, and $\delta+\epsilon>\pi$, and $\delta<\epsilon$, we know $\alpha\epsilon\cdots,\epsilon\thin\epsilon\cdots$ are not vertices. 

By $\beta\delta^2,\gamma\epsilon^2$, we know $\beta\gamma\cdots$ is a $\hat{b}$-vertex. Then by $R(\alpha\beta\gamma)=(\tfrac{2}{3}-\tfrac{8}{f})\pi<\alpha<\beta+\gamma$, and $\beta>\gamma$, we get $\alpha\beta\gamma\cdots=\alpha\beta^2\gamma,\alpha\beta\gamma^k$. Then $\alpha_2\cdots=\thin^{\beta}\alpha^{\gamma}\thin^{\epsilon}\gamma^{\alpha}\thin\beta\thin^{\alpha}\beta^{\delta}\thin,\thin^{\beta}\alpha^{\gamma}\thin^{\epsilon}\gamma^{\alpha}\thin\gamma\thin\cdots\thin\gamma\thin^{\alpha}\beta^{\delta}\thin$. The second case implies $\alpha\epsilon\cdots,\epsilon\thin\epsilon\cdots$, a contradiction. For the first case, the AAD implies a vertex $\alpha\delta\cdots$, and the angle sum of $\alpha\beta^2\gamma$ further implies 
\[
\alpha=\tfrac{2}{3}\pi,\,
\beta=(\tfrac{2}{3}-\tfrac{8}{f})\pi,\,
\gamma=\tfrac{16}{f}\pi,\,
\delta=(\tfrac{2}{3}+\tfrac{4}{f})\pi,\,
\epsilon=(1-\tfrac{8}{f})\pi.
\]
Then by $\alpha+2\delta=(2+\tfrac{8}{f})\pi>2\pi$ and no $\alpha\epsilon\cdots$, we get a contradiction.

\subsubsection*{Case. $\alpha\beta^2$ is a vertex}

By distinct $\alpha,\beta,\gamma$ values, and all the degree $3$ $b$-vertices $\beta\delta^2,\gamma\epsilon^2$, and no $\alpha\beta\gamma$, we know $\beta\delta^2,\gamma\epsilon^2,\alpha\beta^2,\alpha^2\gamma,\beta\gamma^2,\gamma^3$ are all the degree $3$ vertices. We already proved there is no tiling with $\alpha^2\beta$, $\beta^3$ or $\beta^2\gamma$. By the exchange symmetry $(\beta,\delta)\leftrightarrow(\gamma,\epsilon)$, there is also no tiling with $\alpha^2\gamma,\gamma^3$ or $\beta\gamma^2$. Therefore $\beta\delta^2,\gamma\epsilon^2,\alpha\beta^2$ are all the degree $3$ vertices. This implies that, in any tile, one of $\alpha\cdots,\gamma\cdots$ has high degree. Then by Lemma \ref{special_tile}, we get $f\ge 24$. 

The angle sums of $\beta\delta^2,\gamma\epsilon^2,\alpha\beta^2$ and the angle sum for pentagon imply
\[
\alpha=(\tfrac{2}{3}+\tfrac{16}{3f})\pi-\tfrac{2}{3}\gamma,\;
\beta=(\tfrac{2}{3}-\tfrac{8}{3f})\pi+\tfrac{1}{3}\gamma,\;
\delta=(\tfrac{2}{3}+\tfrac{4}{3f})\pi-\tfrac{1}{6}\gamma,\;
\epsilon=\pi-\tfrac{1}{2}\gamma.
\]
We have $2\alpha+\beta+\gamma=(2+\frac{4}{f})\pi>2\pi$, and $\alpha<\pi$, and $\alpha+\gamma>\beta>\frac{1}{2}\pi$. By $\beta\delta^2,\gamma\epsilon^2,\alpha\beta^2$, we also get $\beta,\delta,\epsilon<\pi$. 

By $\beta\delta^2,\gamma\epsilon^2$, we know $\beta\gamma\cdots$ is a $\hat{b}$-vertex. By $\beta\delta^2$ and $\alpha+\gamma>\beta$, we get $\alpha+\gamma+2\delta>2\pi$. By $\gamma\epsilon^2$, we get $\alpha+\gamma+2\epsilon>2\pi$. The two inequalities imply $\alpha\gamma\cdots$ is a $\hat{b}$-vertex. 

Suppose $\gamma\ge \pi$. Then $\gamma^2\cdots$ is not a vertex. Moreover, by $\beta<\pi\le\gamma$ and Lemma \ref{geometry1}, we get $\delta>\epsilon$. Then by $\gamma\epsilon^2$, we know $\gamma\delta\cdots$ is not a vertex. 

In a special tile, one of $\alpha\cdots,\gamma\cdots$ has degree $3$. This means $\alpha\cdots=\alpha\beta^2$ or $\gamma\cdots=\gamma\epsilon^2$. If $\alpha\cdots=\alpha\beta^2$, then by no $\gamma\delta\cdots$, we know $\gamma\cdots=\alpha\gamma\cdots$ and has degree $4$ or $5$. If $\gamma\cdots=\gamma\epsilon^2$, then $\alpha\cdots=\alpha\gamma\cdots$ and has degree $4$ or $5$. 

Therefore $\alpha\gamma\cdots$ is a vertex of degree $4$ or $5$. By $\alpha\beta^2$ and $\beta<\pi\le\gamma$, we get $\alpha\gamma\cdots=\alpha^k\gamma(k=3,4)$. By $f\ge 24$, the angle sums of $\alpha^3\gamma$ and $\alpha^4\gamma$ respectively imply $\gamma=\frac{16}{f}\pi<\pi$ and $\gamma=(\frac{2}{5}+\frac{64}{5f})\pi<\pi$, contradicting the assumption $\gamma\ge\pi$. 

Therefore $\gamma<\pi$. By $\beta\delta^2,\gamma\epsilon^2$, and $\beta,\gamma<\pi$, we get $\delta,\epsilon>\frac{1}{2}\pi$. This implies at most one $b$-edge at a vertex, and $\delta\thin\delta\cdots,\delta\thin\epsilon\cdots,\epsilon\thin\epsilon\cdots$ are not vertices.

By $\alpha\beta^2$, we get $\alpha+\beta>\pi$. Then by $\delta,\epsilon>\frac{1}{2}\pi$, this implies $\alpha\beta\cdots$ is a $\hat{b}$-vertex. Now we know $\alpha\beta\cdots,\alpha\gamma\cdots,\beta\gamma\cdots$ are all $\hat{b}$-vertices. 

By all the degree $3$ vertices $\beta\delta^2,\gamma\epsilon^2,\alpha\beta^2$, we know $\beta\epsilon\cdots,\gamma\delta\cdots$ have high degree. 

Since $\alpha\beta\cdots,\beta\gamma\cdots$ are $\hat{b}$-vertices, we know a $b$-vertex $\beta\cdots$ has no $\alpha,\gamma$. Then by $\beta,\delta,\epsilon>\frac{1}{2}\pi$, the high degree vertex $\beta\epsilon\cdots$ is actually not a vertex. Then by $\beta\delta^2$, we know a $b$-vertex $\beta\cdots=\beta\delta^2$. 

Since $\alpha\gamma\cdots,\beta\gamma\cdots$ are $\hat{b}$-vertices, we know the high degree vertex $\gamma\delta\cdots$ has no $\alpha,\beta$. Then by at most one $b$-edge at a vertex, we get $\gamma\delta\cdots=\gamma^k\delta^2(k\ge 2),\gamma^k\delta\epsilon(k\ge 2)$. The AAD of $\gamma^k\delta^2$ implies $\beta\epsilon\cdots,\epsilon\thin\epsilon\cdots$, a contradiction. The angle sum of $\gamma^k\delta\epsilon$ implies $2\pi\ge 2\gamma+\delta+\epsilon=(\frac{5}{3}+\frac{4}{3f})+\frac{4}{3}\gamma$. This implies $\gamma<\frac{1}{4}\pi$. Then $\beta>\frac{1}{2}\pi>\gamma$. By Lemma \ref{geometry1}, this implies $\delta<\epsilon$. Then $\alpha+2\epsilon>\alpha+\delta+\epsilon=(\frac{7}{3}+\frac{20}{3f})\pi-\frac{4}{3}\gamma>(\frac{7}{3}+\frac{20}{3f})\pi-\frac{1}{3}\pi>2\pi$. Therefore $\alpha\epsilon\cdots$ is not a vertex. Now $\gamma^k\delta\epsilon$ contains $\thick\delta\thin\gamma\thin\gamma\thin$, and the AAD implies $\alpha\epsilon\cdots,\beta\epsilon\cdots,\epsilon\thin\epsilon\cdots$, a contradiction.  Therefore $\gamma\delta\cdots$ is not a vertex. By $\gamma\epsilon^2$, this implies a $b$-vertex $\gamma\cdots=\gamma\epsilon^2$.

Next we prove $\alpha\beta\gamma\cdots$ is not a vertex. We know this is a high degree $\hat{b}$-vertex. By $\alpha\beta^2$, and $2\alpha+\beta+\gamma>2\pi$, we know the remainder has no $\alpha,\beta$. Therefore $\alpha\beta\gamma\cdots=\alpha\beta\gamma^k(k\ge 2)$. By $\alpha\beta^2$ and $\alpha\beta\gamma^k$, we get $\beta=k\gamma>\gamma$. This means $\gamma=\frac{2}{3k-1}(1-\frac{4}{f})\pi$. By Lemma \ref{geometry1}, we also get $\delta<\epsilon$. 

If $k\ge 3$ in $\alpha\beta\gamma^k$, then we get $\gamma<\frac{1}{4}\pi$. Then the earlier argument proves that $\alpha\epsilon\cdots$ is not a vertex. Now $\alpha\beta\gamma^k(k\ge 3)$ contains $\thin\beta\thin\gamma\thin\gamma\thin\cdots$, and the AAD implies $\alpha\epsilon\cdots,\delta\thin\epsilon\cdots,\epsilon\thin\epsilon\cdots$, a contradiction. 
Therefore $k=2$ in $\alpha\beta\gamma^k$. The angle sum of $\alpha\beta\gamma^2$ further implies 
\[
\alpha=(\tfrac{2}{5}+\tfrac{32}{5f})\pi,\,
\beta=(\tfrac{4}{5}-\tfrac{16}{5f})\pi,\,
\gamma=(\tfrac{2}{5}-\tfrac{8}{5f})\pi,\,
\delta=(\tfrac{3}{5}+\tfrac{8}{5f})\pi,\,
\epsilon=(\tfrac{4}{5}+\tfrac{4}{5f})\pi.
\]
By no $\gamma\delta\cdots,\delta\thin\delta\cdots$, the AAD of $\alpha\beta^2$ is $\thin^{\beta}\alpha^{\gamma}\thin^{\alpha}\beta^{\delta}\thin^{\alpha}\beta^{\delta}\thin$. This implies a vertex $\alpha\delta\cdots$. Since $\alpha\beta\cdots,\alpha\gamma\cdots$ are $\hat{b}$-vertices, we know $\alpha\delta\cdots$ has no $\beta,\gamma$. By at most one $b$-edge at a vertex, and $\alpha+\delta=(1+\frac{8}{f})\pi>\pi$, and $2\alpha+\delta+\epsilon=(\frac{11}{5}+\frac{76}{5f})\pi>2\pi$, we know $\alpha\delta\cdots$ has degree $3$, a contradiction. This completes the proof that $\alpha\beta\gamma\cdots$ is not a vertex.

By all the degree $3$ vertices $\beta\delta^2,\gamma\epsilon^2,\alpha\beta^2$, we know a special companion pair is matched. This determines $T_2,T_2$ in Figure \ref{b2d_c2eC}. Moreover, we know one of $\alpha_1\cdots,\gamma_1\cdots$ has high degree, and one of $\alpha_2\cdots,\gamma_2\cdots$ has high degree. This implies both shared vertices $\delta_1\delta_2\cdots,\epsilon_1\epsilon_2\cdots$ have degree $3$. Then $\delta_1\delta_2\cdots=\beta\delta^2$ and $\epsilon_1\epsilon_2\cdots=\gamma\epsilon^2$. Up to the vertical flip symmetry, we may assume $T_3$ is arranged as indicated. 

If $T_4$ is not arranged as indicated, then both $\beta_1\cdots,\gamma_1\cdots$ are $b$-vertices. This implies $\beta_1\cdots=\beta\delta^2$ and $\gamma_1\cdots=\gamma\epsilon^2$. Then $\alpha_1\cdots=\alpha\beta\gamma\cdots$, a contradiction. Therefore $T_4$ is arranged as indicated. Then the $b$-vertex $\gamma_2\cdots=\gamma\epsilon^2$ determine $T_5$.

\begin{figure}[htp]
\centering
\begin{tikzpicture}[>=latex,scale=1]

\foreach \a in {-1,1}
\draw[xscale=\a] 
	(0,1.1) -- (0.5,0.7) -- (0.5,-0.7) -- (0,-1.1) -- (0,-1.1)
	(0.5,0.7) -- (1.3,0.7) -- (1.3,-0.7) -- (0.5,-0.7) 
	(0.5,0) -- (0,0);

\draw[line width=1.2]
	(1.3,-0.7) -- (0.5,-0.7)
	(-1.3,0.7) -- (-0.5,0.7)
	(0.5,0) -- (-0.5,0)
	(-1.3,-1.8) -- (0,-1.8);

\draw
	(1.3,-0.7) -- (1.3,-1.8) -- (0.3,-1.8) -- (0,-1.1) -- (0,-1.8) -- (-1.3,-1.8) -- (-1.3,-1) -- (-0.5,-0.7);
	
\node at (0.3,-0.6) {\small $\gamma$};
\node at (0,-0.85) {\small $\alpha$};
\node at (0.3,-0.2) {\small $\epsilon$};	
\node at (-0.3,-0.2) {\small $\delta$}; 
\node at (-0.3,-0.6) {\small $\beta$};

\node at (0.3,0.6) {\small $\gamma$};
\node at (0,0.85) {\small $\alpha$};
\node at (0.3,0.2) {\small $\epsilon$};	
\node at (-0.3,0.2) {\small $\delta$}; 
\node at (-0.3,0.6) {\small $\beta$};

\node at (1.1,0.5) {\small $\beta$};
\node at (1.1,-0.5) {\small $\delta$};
\node at (0.7,0.5) {\small $\alpha$};
\node at (0.7,-0.5) {\small $\epsilon$};
\node at (0.7,0) {\small $\gamma$};

\node at (-1.1,0.5) {\small $\epsilon$};
\node at (-1.1,-0.5) {\small $\gamma$};
\node at (-0.7,0.5) {\small $\delta$};
\node at (-0.7,-0.5) {\small $\alpha$};
\node at (-0.7,0) {\small $\beta$};

\node at (-1.1,-1.1) {\small $\beta$};
\node at (-0.6,-0.9) {\small $\alpha$};
\node at (-1.1,-1.6) {\small $\delta$};
\node at (-0.2,-1.2) {\small $\gamma$};
\node at (-0.2,-1.6) {\small $\epsilon$};

\node at (1.1,-0.9) {\small $\delta$};
\node at (0.6,-0.9) {\small $\epsilon$};
\node at (1.1,-1.6) {\small $\beta$};
\node at (0.2,-1.2) {\small $\gamma$};
\node at (0.4,-1.6) {\small $\alpha$};

\node[inner sep=0.5,draw,shape=circle] at (0,0.4) {\small $1$};
\node[inner sep=0.5,draw,shape=circle] at (0,-0.4) {\small $2$};
\node[inner sep=0.5,draw,shape=circle] at (-1.05,0) {\small $3$};
\node[inner sep=0.5,draw,shape=circle] at (1.05,0) {\small $4$};
\node[inner sep=0.5,draw,shape=circle] at (0.7,-1.3) {\small $5$};
\node[inner sep=0.5,draw,shape=circle] at (-0.7,-1.3) {\small $6$};

\end{tikzpicture}
\caption{Proposition \ref{b2d_c2e}: $\alpha\beta^2$ is a vertex.}
\label{b2d_c2eC}
\end{figure}
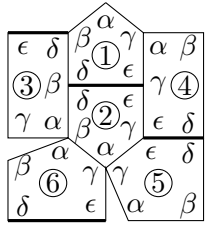

We know $\alpha_3\beta_2\cdots$ is a $\hat{b}$-vertex without $\gamma$. Therefore $\alpha_3\beta_2\cdots=\alpha\beta^2,\alpha^k\beta(k\ge 3)$. If $\alpha_3\beta_2\cdots=\alpha\beta^2$, then $\gamma_3\cdots=\gamma\delta\cdots$ or $\alpha_2\gamma_6\cdots=\gamma\delta\cdots$, both contradictions. Therefore $\alpha_3\beta_2\cdots=\alpha^k\beta(k\ge 3)$. This gives $\alpha_6$. Then by no $\alpha\beta\gamma\cdots$, we determine $T_6$. Since $\alpha_2,\beta_2\cdots,\gamma_1\cdots$ have high degree, by Lemma \ref{special_pair}, we get $f\ge 48$. 

We know $\alpha_2\gamma_6\gamma_7\cdots$ is a high degree $\hat{b}$-vertex without $\beta$. Then by $3\alpha+2\gamma=(2+\tfrac{16}{f})\pi>2\pi$, we get $\alpha_2\gamma_7\gamma_8\cdots=\alpha\gamma^k(k\ge 3),\alpha^2\gamma^k(k\ge 2)$. If $\alpha\gamma^k(k\ge 4)$ or $\alpha^2\gamma^k(k\ge 3)$ is a vertex, then $\alpha+4\gamma\le 2\pi$ or $2\alpha+3\gamma\le 2\pi$. By $3\alpha+2\gamma>2\pi$, this implies $\alpha>\gamma$. Then by $\alpha\beta^2$, the vertex $\alpha\gamma^k(k\ge 4)$ or $\alpha^2\gamma^k(k\ge 3)$ implies $\beta\ge 2\gamma$. Then we get $3\alpha+\beta\ge 3\alpha+2\gamma>2\pi$, contradicting the vertex $\alpha_3\beta_2\cdots=\alpha^k\beta$.

Therefore $\alpha_2\gamma_7\gamma_8\cdots=\alpha\gamma^3,\alpha^2\gamma^2$. The angle sum of the vertex further implies
\begin{align*}
\alpha\gamma^3 &\colon
\alpha=(\tfrac{2}{7}+\tfrac{48}{7f})\pi,\,
\beta=(\tfrac{6}{7}-\tfrac{24}{7f})\pi,\,
\gamma=(\tfrac{4}{7}-\tfrac{16}{7f})\pi,\,
\delta=(\tfrac{4}{7}+\tfrac{12}{7f})\pi,\,
\epsilon=(\tfrac{5}{7}+\tfrac{8}{7f})\pi.
\\
\alpha^2\gamma^2 &\colon
\alpha=\tfrac{16}{f}\pi,\,
\beta=(1-\tfrac{8}{f})\pi,\,
\gamma=(1-\tfrac{16}{f})\pi,\,
\delta=(\tfrac{1}{2}+\tfrac{4}{f})\pi,\,
\epsilon=(\tfrac{1}{2}+\tfrac{8}{f})\pi.
\end{align*}
We have $\beta>\gamma$ and $\delta<\epsilon$ in both cases. By $f\ge 48$, we get $\alpha<2\gamma$ and $\beta+3\gamma>2\pi$ in both cases.

By $b$-vertex $\beta\cdots=\beta\delta^2$, we know $\beta^2\cdots,\beta\gamma\cdots$ are $\hat{b}$-vertices. By $\alpha\beta^2$, and $\alpha<2\gamma$, and $\beta>\gamma$, we get $\beta^2\cdots=\alpha\beta^2$. By no $\alpha\beta\gamma\cdots$, we know $\beta\gamma\cdots$ has no $\alpha$. Then by $\beta+3\gamma>2\pi$ and $\beta>\gamma$, we know $\beta\gamma\cdots$ has degree $3$. By all the degree $3$ vertices $\beta\delta^2,\gamma\epsilon^2,\alpha\beta^2$, this implies $\beta\gamma\cdots$ is not a vertex.

By no $\beta\gamma\cdots$, the AAD of $\alpha_3\beta_2\cdots=\alpha^k\beta(k\ge 3)$ contains $\thin^{\gamma}\alpha^{\beta}\thin^{\beta}\alpha^{\gamma}\thin$. This implies a vertex $\thin^{\delta}\beta^{\alpha}\thin^{\alpha}\beta^{\delta}\thin\cdots=\alpha\beta^2=\thin^{\alpha}\beta^{\delta}\thin^{\beta}\alpha^{\gamma}\thin^{\delta}\beta^{\alpha}\thin$, contradicting no $\gamma\delta\cdots$. 
\end{proof}

\subsection{Single Matched Degree $3$ $b$-Vertex}
\label{singleb}

We study tilings with a single matched degree $3$ $b$-vertex. Up to the exchange symmetry $(\beta,\delta)\leftrightarrow(\gamma,\epsilon)$, we need to consider three vertices $\alpha\delta^2,\beta\delta^2,\gamma\delta^2$. Then $\epsilon$ does not appear at degree $3$ vertices. By Lemma \ref{ndegree3}, this implies $\delta\epsilon^3$ or $\epsilon^4$ is a vertex. Then the AADs of $\delta\epsilon^3,\epsilon^4$ imply $\gamma^2\cdots$ is a vertex. Moreover, in a special tile, the vertices $\alpha\cdots,\beta\cdots,\gamma\cdots,\delta\cdots$ have degree $3$, and $H=\epsilon\cdots$ has degree $4$ or $5$. By Lemma \ref{special_tile}, this implies $f\ge 24$. 

\begin{proposition}\label{a2d}
There is no tiling, such that $\alpha,\beta,\gamma$ have distinct values, and $\alpha\delta^2$ is the only degree $3$ $b$-vertex.
\end{proposition}

\begin{proof}
By applying Lemma \ref{degree3} to $\alpha$, there is a degree $3$ vertex without $\alpha$. By the only degree $3$ $b$-vertex $\alpha\delta^2$, degree $3$ vertices without $\alpha$ are $\beta^3,\beta^2\gamma,\beta\gamma^2,\gamma^3$. By $\beta\ne\gamma$, these are mutually exclusive. By applying Lemma \ref{degree3} to $\beta,\delta$ or $\gamma,\delta$, we know $\alpha\delta^2$ and one of $\beta^3,\beta^2\gamma,\beta\gamma^2,\gamma^3$ cannot be all the degree $3$ vertices. Therefore there is another degree $3$ $\hat{b}$-vertex with $\alpha$. This means there are at least three vertices, one from each of the following groups:
\begin{enumerate}
\item degree $3$ $b$-vertex: $\alpha\delta^2$.
\item degree $3$ $\hat{b}$-vertex without $\alpha$: $\beta^3,\beta^2\gamma,\beta\gamma^2,\gamma^3$.
\item degree $3$ $\hat{b}$-vertex with $\alpha$: $\alpha^3,\alpha^2\beta,\alpha^2\gamma,\alpha\beta^2,\alpha\gamma^2,\alpha\beta\gamma$. 
\end{enumerate}
By distinct $\alpha,\beta,\gamma$ values, we find there cannot be fourth degree $3$ vertex in the tiling. Then the collection $\text{AVC}_3$ of all degree $3$ vertices is exactly one of the following:
\[
\begin{array}{lll}
\{\alpha\delta^2,\beta^3,\alpha^2\gamma\}, &
\{\alpha\delta^2,\beta^3,\alpha\gamma^2\}, &
\{\alpha\delta^2,\beta^3,\alpha\beta\gamma\}, \\
\{\alpha\delta^2,\beta^2\gamma,\alpha^3\}, &
\{\alpha\delta^2,\beta^2\gamma,\alpha^2\beta\}, &
\{\alpha\delta^2,\beta^2\gamma,\alpha\gamma^2\}, \\
\{\alpha\delta^2,\beta\gamma^2,\alpha^3\}, &
\{\alpha\delta^2,\beta\gamma^2,\alpha\beta^2\}, &
\{\alpha\delta^2,\beta\gamma^2,\alpha^2\gamma\}, \\ 
\{\alpha\delta^2,\gamma^3,\alpha^2\beta\}, &
\{\alpha\delta^2,\gamma^3,\alpha\beta^2\}, &
\{\alpha\delta^2,\gamma^3,\alpha\beta\gamma\}.
\end{array}
\]
In a special tile, the degree $3$ vertex $\delta\cdots=\alpha\delta^2$, and the vertices $\alpha\cdots$, $\beta\cdots$, $\gamma\cdots$ have degree $3$.

For $\text{AVC}_3=\{\alpha\delta^2,\beta^3,\alpha\gamma^2\}$, in a special tile, the degree $3$ vertex $\beta\cdots=\beta^3$. This implies the degree $3$ vertex $\alpha\cdots=\alpha\delta^2$. Then $\gamma\cdots=\beta\gamma\cdots$ has high degree, a contradiction. 

For $\text{AVC}_3=\{\alpha\delta^2,\beta\gamma^2,\alpha\beta^2\}$, in a special tile, the degree $3$ vertex $\gamma\cdots=\beta\gamma^2$. This implies the degree $3$ vertex $\alpha\cdots=\alpha\delta^2$. Then by $\delta\cdots=\alpha\delta^2$, we get $\beta\cdots=\beta^3\cdots,\beta^2\gamma\cdots$. This has high degree, a contradiction. 

For $\text{AVC}_3=\{\alpha\delta^2,\gamma^3,\alpha^2\beta\}$,  in a special tile, by $\delta\cdots=\alpha\delta^2$, we get $\beta\cdots=\beta^2\cdots,\beta\gamma\cdots$. This has high degree, a contradiction. 

For $\text{AVC}_3=\{\alpha\delta^2,\gamma^3,\alpha\beta^2\}$, in a special tile, one of $\alpha\cdots,\gamma\cdots$ has high degree, a contradiction.  

We recall that $\delta\epsilon^3$ or $\epsilon^4$ is a vertex, and $\gamma^2\cdots$ is a vertex.

\subsubsection*{Case. $\text{AVC}_3=\{\alpha\delta^2,\beta^3,\alpha^2\gamma\}$ or $\{\alpha\delta^2,\beta^3,\alpha\beta\gamma\}$}

The angle sums of the degree $3$ vertices, and the angle sum of $\delta\epsilon^3$ or $\epsilon^4$, and the angle sum for pentagon imply
\begin{align*}
\alpha\delta^2,\beta^3,\alpha^2\gamma,\delta\epsilon^3 &\colon
	\alpha=(\tfrac{3}{4}-\tfrac{3}{f})\pi,\,
	\delta=(\tfrac{5}{8}+\tfrac{3}{2f})\pi,\,
	\epsilon=(\tfrac{11}{24}-\tfrac{1}{2f})\pi. \\
\alpha\delta^2,\beta^3,\alpha^2\gamma,\epsilon^4 &\colon
	\alpha=(\tfrac{7}{9}-\tfrac{8}{3f})\pi,\,
	\delta=(\tfrac{11}{18}+\tfrac{4}{3f})\pi,\,
	\epsilon=\tfrac{1}{2}\pi.  \\
\alpha\delta^2,\beta^3,\alpha\beta\gamma,\delta\epsilon^3 &\colon
	\alpha=(1-\tfrac{12}{f})\pi,\,
	\delta=(\tfrac{1}{2}+\tfrac{6}{f})\pi,\,
	\epsilon=(\tfrac{1}{2}-\tfrac{2}{f})\pi. \\
\alpha\delta^2,\beta^3,\alpha\beta\gamma,\epsilon^4 &\colon
	\alpha=(1-\tfrac{8}{f})\pi,\,
	\delta=(\tfrac{1}{2}+\tfrac{4}{f})\pi,\,
	\epsilon=\tfrac{1}{2}\pi.
\end{align*}
We have $\delta>\epsilon$ and $\alpha<2\epsilon$. By Lemma \ref{geometry1}, this implies $\beta<\gamma$. Then by $\beta^3,\alpha^2\gamma$ or $\beta^3,\alpha\beta\gamma$, we get $\alpha<\beta<\gamma$. 

By $\alpha^2\gamma$ or $\alpha\beta\gamma$, and $\alpha,\beta<\gamma$, we get $R(\gamma^2)<\alpha<\beta,\gamma,2\delta,2\epsilon$. This implies $\gamma^2\cdots$ is not a vertex, a contradiction.

\subsubsection*{Case. $\text{AVC}_3=
\{\alpha\delta^2,\beta^2\gamma,\alpha^3\}$ or $\{\alpha\delta^2,\beta^2\gamma,\alpha^2\beta\}$}

The angle sums of the degree $3$ vertices, and the angle sum of $\delta\epsilon^3$ or $\epsilon^4$, and the angle sum for pentagon imply
\begin{align*}
\alpha\delta^2,\beta^2\gamma,\alpha^3,\delta\epsilon^3 &\colon
	\beta=(\tfrac{7}{9}-\tfrac{4}{f})\pi,\,
	\delta=\tfrac{2}{3}\pi,\,
	\epsilon=\tfrac{4}{9}\pi. \\
\alpha\delta^2,\beta^2\gamma,\alpha^3,\epsilon^4 &\colon
	\beta=(\tfrac{5}{6}-\tfrac{4}{f})\pi,\,
	\delta=\tfrac{2}{3}\pi,\,
	\epsilon=\tfrac{1}{2}\pi. \\
\alpha\delta^2,\beta^2\gamma,\alpha^2\beta,\delta\epsilon^3 &\colon
	\beta=(\tfrac{3}{4}-\tfrac{3}{f})\pi,\,
	\delta=(\tfrac{11}{16}-\tfrac{3}{4f})\pi,\,
	\epsilon=(\tfrac{7}{16}+\tfrac{1}{4f})\pi. \\
\alpha\delta^2,\beta^2\gamma,\alpha^2\beta,\epsilon^4 &\colon
	\beta=(\tfrac{4}{5}-\tfrac{16}{5f})\pi,\,
	\delta=(\tfrac{7}{10}-\tfrac{4}{5f})\pi,\,
	\epsilon=\tfrac{1}{2}\pi.
\end{align*}
We have $\delta>\epsilon$ and $\beta<2\epsilon$. By Lemma \ref{geometry1}, this implies $\beta<\gamma$. Then by $\beta^2\gamma,\alpha^3$ or $\beta^2\gamma,\alpha^2\beta$, we get $\beta<\alpha<\gamma$. 

By $\beta^2\gamma$ and $\beta<\gamma$, we get $R(\gamma^2)<\beta<\alpha,\gamma,2\delta,2\epsilon$. This implies $\gamma^2\cdots$ is not a vertex, a contradiction.

\subsubsection*{Case. $\text{AVC}_3=\{\alpha\delta^2,\beta^2\gamma,\alpha\gamma^2\}$}

The angle sums of the degree $3$ vertices, and the angle sum of $\delta\epsilon^3$ or $\epsilon^4$, and the angle sum for pentagon imply
\begin{align*}
\delta\epsilon^3 &\colon
	\alpha=(\tfrac{2}{5}+\tfrac{48}{5f})\pi,\,
	\beta=(\tfrac{3}{5}+\tfrac{12}{5f})\pi,\,
	\gamma=\delta=(\tfrac{4}{5}-\tfrac{24}{5f})\pi,\,
	\epsilon=(\tfrac{2}{5}+\tfrac{8}{5f})\pi. \\
\epsilon^4 &\colon
	\alpha=\tfrac{16}{f}\pi,\,
	\beta=(\tfrac{1}{2}+\tfrac{4}{f})\pi,\,
	\gamma=\delta=(1-\tfrac{8}{f})\pi,\,
	\epsilon=\tfrac{1}{2}\pi.
\end{align*}
We have $\alpha+\gamma=2\beta>\pi$ and $\beta+\epsilon>\pi$. By $f\ge 24$, we get $\delta>\epsilon$ and $\alpha<2\epsilon$. By Lemma \ref{geometry1}, this implies $\beta<\gamma$. Then by $\alpha+\gamma=2\beta$, we get $\alpha<\beta<\gamma$. By $\delta>\epsilon$, and $\delta\epsilon^3$ or $\epsilon^4$, and $\alpha\delta^2$, we know $\delta\thin\delta\cdots,\delta\thin\alpha\thin\epsilon\cdots$ are not vertices, and $\delta\thin\epsilon\cdots=\delta\epsilon^3$.

By $R(\gamma^2)=\alpha<\beta,\gamma,\delta,2\epsilon$, we get $\gamma^2\cdots=\alpha\gamma^2$.

By $\alpha\delta^2$ and $\alpha<\beta,\gamma$, we know $\beta\delta^2\cdots,\gamma\delta^2\cdots$ are not vertices. By $R(\gamma\delta\epsilon)<R(\beta\delta\epsilon)=\gamma-\epsilon<\alpha<\beta,\gamma,2\delta,2\epsilon$, we know $\beta\delta\epsilon\cdots,\gamma\delta\epsilon\cdots$ are not vertices. Therefore $\beta\delta\cdots,\gamma\delta\cdots$ are not vertices.

By $\alpha+\gamma=\alpha+\delta=2\beta>\pi$, and $\beta+\epsilon>\pi$, we know $\alpha\epsilon^2\cdots=\alpha^k\epsilon^2,\alpha^k\beta\epsilon^2$. By no $\beta\delta\cdots,\gamma\delta\cdots$, the AAD implies $\alpha^k\beta\epsilon^2$ is not a vertex. Therefore $\alpha\epsilon^2\cdots=\alpha^k\epsilon^2$. Then by no $\beta\delta\cdots,\gamma\delta\cdots$, this further implies $\alpha\beta\epsilon\cdots,\alpha\gamma\epsilon\cdots$ are not vertices.

By no $\delta\thin\delta\cdots$, we know the AAD of $\thin\beta\thin\beta\thin$ is $\thin^{\delta}\beta^{\alpha}\thin^{\alpha}\beta^{\delta}\thin$ or $\thin^{\alpha}\beta^{\delta}\thin^{\alpha}\beta^{\delta}\thin$. By no $\beta\delta\cdots$ and $\beta+\epsilon>\pi$, we know $\beta^2\cdots$ is a $\hat{b}$-vertex. Then by no $\beta\delta\cdots,\gamma\delta\cdots$, we know $\thin^{\delta}\beta^{\alpha}\thin^{\alpha}\beta^{\delta}\thin\cdots,\thin^{\alpha}\beta^{\delta}\thin^{\alpha}\beta^{\delta}\thin\cdots$ cannot be $\thin^{\delta}\beta^{\alpha}\thin^{\alpha}\beta^{\delta}\thin\alpha\thin\cdots,\thin^{\alpha}\beta^{\delta}\thin^{\alpha}\beta^{\delta}\thin\alpha\thin\cdots$. Then by $\beta^2\gamma$, we get $\beta\thin\beta\cdots=\beta^2\gamma,\beta^k$. By $\beta^2\gamma$ and $\beta>\frac{1}{2}\pi$, we know $\beta^k$ is not a vertex. Therefore $\beta\thin\beta\cdots=\beta^2\gamma$.

In the special companion pair $T_1,T_2$ in the first of Figure \ref{a2dA}, we have $\delta_1\delta_2\cdots=\alpha\delta^2$, and $\epsilon_1\epsilon_2\cdots$ has degree $4$ or $5$. We may assume $T_1$ is a special tile. Then $\epsilon_1\epsilon_2\cdots$ is the only high degree vertex in $T_1$, and by all the degree $3$ vertices, we get $\alpha_1\cdots,\beta_1\cdots,\gamma_1\cdots$ as indicated. This determines $T_3,T_4,T_5$, and gives $\beta_6$. Then $\beta_2\thin\beta_3\cdots=\beta^2\gamma$ gives $\gamma_7$. 

If $T_7$ is not arranged as indicated, then $\delta_3\cdots=\delta\thin\epsilon\cdots=\delta\epsilon^3$. This implies $\alpha_4\epsilon_3\cdots=\alpha\delta\epsilon\cdots$. Since $\delta_3\cdots=\delta\epsilon^3$ is a vertex, we get $R(\alpha\delta\epsilon)=(\frac{2}{5}-\frac{32}{f})\pi<\alpha,\beta,\gamma,\delta,\epsilon$. This implies $\alpha\delta\epsilon\cdots=\alpha\delta\epsilon$, contradicting $\alpha\delta^2$. 

Therefore $T_7$ is arranged as indicated. Now $\beta_6$ implies $\epsilon_1\epsilon_2\cdots=\alpha\epsilon^2\cdots,\delta\epsilon^2\cdots$. If the vertex is $\delta\epsilon^2\cdots$, then it is $\delta\thin\epsilon\cdots=\delta\epsilon^3$. This implies $\gamma_2\cdots=\gamma^2\cdots=\alpha\gamma^2$, and further implies $\alpha_2\epsilon_7\cdots=\alpha\beta\epsilon\cdots,\alpha\gamma\epsilon\cdots$, a contradiction.

Therefore $\epsilon_1\epsilon_2\cdots=\alpha\epsilon^2\cdots=\alpha^k\epsilon^2$. This gives $\alpha_8$. If $T_8$ is not arranged as indicated, then the degree $3$ vertex $\gamma_2\cdots=\gamma^2\cdots=\alpha\gamma^2$. This implies $\alpha_2\epsilon_7\cdots=\alpha\beta\epsilon\cdots,\alpha\gamma\epsilon\cdots$, a contradiction. Therefore $T_8$ is arranged as indicated. Then the degree $3$ vertex $\gamma_2\cdots=\beta\gamma\cdots=\beta^2\gamma$. This implies either $\delta_8\cdots=\delta\thin\delta\cdots$, or $\alpha_2\epsilon_7\cdots=\thick\delta\thin\alpha\thin\epsilon\thick\cdots$, both contradictions. 

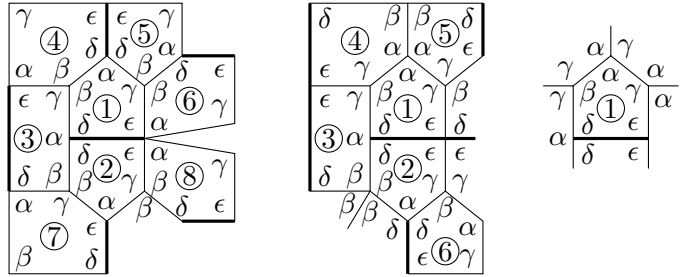
\begin{figure}[htp]
\centering
\begin{tikzpicture}[>=latex,scale=1]

\foreach \b in {1,-1}
{
\begin{scope}[yscale=\b]

\draw
	(-0.5,0) -- (-0.5,0.7) -- (0,1.1) -- (0.5,0.7) -- (0.5,0)
	(0.5,0) -- (1.7,0.2) -- (1.7,1.1) -- (1,1.1) -- (0.5,0.7)
	(-0.5,0.7) -- (-1.3,0.7);

\draw[line width=1.2]
	(-0.5,0) -- (0.5,0)
	(1,1.1) -- (1.7,1.1);
	
\node at (0,0.85) {\small $\alpha$}; 
\node at (-0.3,0.6) {\small $\beta$};
\node at (0.3,0.6) {\small $\gamma$};
\node at (-0.3,0.2) {\small $\delta$};
\node at (0.3,0.2) {\small $\epsilon$};	

\node at (1,0.9) {\small $\delta$}; 
\node at (0.7,0.6) {\small $\beta$};
\node at (1.5,0.9) {\small $\epsilon$};
\node at (0.7,0.2) {\small $\alpha$};
\node at (1.5,0.4) {\small $\gamma$};

\end{scope}
}

\draw
	(1,1.1) -- (1,1.8) -- (-1.3,1.8) -- (-1.3,-1.8) -- (0,-1.8) -- (0,-1.1);
	
\draw[line width=1.2]
	(-1.3,0.7) -- (-1.3,-0.7)
	(0,1.1) -- (0,1.8) 
	(0,-1.1) -- (0,-1.8) ;

\node at (0.5,-0.95) {\small $\beta$};

\node at (0.8,1.2) {\small $\alpha$};
\node at (0.5,0.95) {\small $\beta$};
\node at (0.8,1.6) {\small $\gamma$};
\node at (0.2,1.2) {\small $\delta$};
\node at (0.2,1.6) {\small $\epsilon$};
	
\node at (-1.1,0.5) {\small $\epsilon$};
\node at (-1.1,-0.5) {\small $\delta$};
\node at (-0.7,0.5) {\small $\gamma$};
\node at (-0.7,-0.5) {\small $\beta$};
\node at (-0.7,0) {\small $\alpha$};

\node at (-0.6,0.9) {\small $\beta$};
\node at (-0.2,1.2) {\small $\delta$};
\node at (-1.1,0.9) {\small $\alpha$};
\node at (-1.1,1.6) {\small $\gamma$};
\node at (-0.2,1.6) {\small $\epsilon$};

\node at (-0.6,-0.9) {\small $\gamma$};
\node at (-0.2,-1.2) {\small $\epsilon$};
\node at (-1.1,-0.9) {\small $\alpha$};
\node at (-1.1,-1.6) {\small $\beta$};
\node at (-0.2,-1.6) {\small $\delta$};

\node[inner sep=0.5,draw,shape=circle] at (0,0.4) {\small $1$};
\node[inner sep=0.5,draw,shape=circle] at (0,-0.4) {\small $2$};
\node[inner sep=0.5,draw,shape=circle] at (-1.05,0) {\small $3$};
\node[inner sep=0.5,draw,shape=circle] at (-0.7,1.3) {\small $4$};
\node[inner sep=0.5,draw,shape=circle] at (0.5,1.4) {\small $5$};
\node[inner sep=0.5,draw,shape=circle] at (1.1,0.5) {\small $6$};
\node[inner sep=0.5,draw,shape=circle] at (-0.7,-1.3) {\small $7$};
\node[inner sep=0.5,draw,shape=circle] at (1.1,-0.5) {\small $8$};


\begin{scope}[xshift=4cm]

\foreach \b in {1,-1}
{
\begin{scope}[yscale=\b]

\draw
	(-0.5,0) -- (-0.5,0.7) -- (0,1.1) -- (0.5,0.7) -- (0.5,0)
	(0.5,0.7) -- (1,1.1) -- (1,1.8) -- (0,1.8);

\node at (0,0.85) {\small $\alpha$}; 
\node at (-0.3,0.6) {\small $\beta$};
\node at (0.3,0.6) {\small $\gamma$};
\node at (-0.3,0.2) {\small $\delta$};
\node at (0.3,0.2) {\small $\epsilon$};	

\end{scope}
}

\draw
	(0,1.8) -- (-1.3,1.8) -- (-1.3,-0.7) -- (-0.5,-0.7) -- (-0.8,-1.2)
	(0,1.1) -- (0,1.8)
	(-0.5,0.7) -- (-1.3,0.7);
		
\draw[line width=1.2]
	(-0.5,0) -- (0.9,0)
	(0,-1.1) -- (0,-1.8)
	(1,1.1) -- (1,1.8)
	(-1.3,1.8) -- (-1.3,-0.7);

\node at (0.5,0.95) {\small $\gamma$};
\node at (0.2,1.2) {\small $\alpha$};
\node at (0.8,1.2) {\small $\epsilon$};
\node at (0.8,1.6) {\small $\delta$};
\node at (0.2,1.6) {\small $\beta$};

\node at (-0.6,0.9) {\small $\gamma$};
\node at (-0.2,1.2) {\small $\alpha$};
\node at (-1.1,0.9) {\small $\epsilon$};
\node at (-1.1,1.6) {\small $\delta$};
\node at (-0.2,1.6) {\small $\beta$};
	
\node at (-1.1,0.5) {\small $\epsilon$};
\node at (-1.1,-0.5) {\small $\delta$};
\node at (-0.7,0.5) {\small $\gamma$};
\node at (-0.7,-0.5) {\small $\beta$};
\node at (-0.7,0) {\small $\alpha$};

\node at (0.8,-1.2) {\small $\alpha$};
\node at (0.5,-0.95) {\small $\beta$};
\node at (0.8,-1.6) {\small $\gamma$};
\node at (0.2,-1.2) {\small $\delta$};
\node at (0.2,-1.6) {\small $\epsilon$};

\node at (0.7,0.2) {\small $\delta$};
\node at (0.7,0.6) {\small $\beta$};
\node at (0.7,-0.2) {\small $\epsilon$};
\node at (0.7,-0.6) {\small $\gamma$};
\node at (-0.2,-1.2) {\small $\delta$};
\node at (-0.8,-0.9) {\small $\beta$};
\node at (-0.5,-1) {\small $\beta$};

\node[inner sep=0.5,draw,shape=circle] at (0,0.4) {\small $1$};
\node[inner sep=0.5,draw,shape=circle] at (0,-0.4) {\small $2$};
\node[inner sep=0.5,draw,shape=circle] at (-1.05,0) {\small $3$};
\node[inner sep=0.5,draw,shape=circle] at (-0.7,1.3) {\small $4$};
\node[inner sep=0.5,draw,shape=circle] at (0.5,1.4) {\small $5$};
\node[inner sep=0.5,draw,shape=circle] at (0.5,-1.5) {\small $6$};

\end{scope}


\begin{scope}[xshift=6.7cm]

\draw
	(-0.5,-0.4) -- (-0.5,0.7) -- (0,1.1) -- (0.5,0.7) -- (0.5,-0.4)
	(0.5,0.7) -- ++(0.4,0)
	(-0.5,0.7) -- ++(-0.4,0)
	(0,1.1) -- ++(0,0.4);
		
\draw[line width=1.2]
	(-0.5,0) -- (0.5,0);

\node at (0,0.85) {\small $\alpha$}; 
\node at (-0.3,0.6) {\small $\beta$};
\node at (0.3,0.6) {\small $\gamma$};
\node at (-0.3,0.2) {\small $\delta$};
\node at (0.3,0.2) {\small $\epsilon$};	

\node at (-0.7,0) {\small $\alpha$};
\node at (-0.3,-0.2) {\small $\delta$};
\node at (0.3,-0.2) {\small $\epsilon$};
	
\node at (0.7,0.5) {\small $\alpha$};
\node at (-0.7,0.5) {\small $\gamma$};
\node at (0.6,0.9) {\small $\alpha$};
\node at (-0.6,0.9) {\small $\gamma$};
\node at (0.2,1.2) {\small $\gamma$};
\node at (-0.2,1.2) {\small $\alpha$};

\node[inner sep=0.5,draw,shape=circle] at (0,0.4) {\small $1$};

\end{scope}

\end{tikzpicture}
\caption{Proposition \ref{a2d}: Special companion pair and special tile.}
\label{a2dA}
\end{figure}

\subsubsection*{Case. $\text{AVC}_3=\{\alpha\delta^2,\beta\gamma^2,\alpha^3\}$}

The angle sums of the degree $3$ vertices, and the angle sum of $\delta\epsilon^3$ or $\epsilon^4$, and the angle sum for pentagon imply
\begin{align*}
\delta\epsilon^3 &\colon
	\alpha=\delta=\tfrac{2}{3}\pi,\,
	\beta=(\tfrac{4}{9}+\tfrac{8}{f})\pi,\,
	\gamma=(\tfrac{7}{9}-\tfrac{4}{f})\pi,\,
	\epsilon=\tfrac{4}{9}\pi. \\
\epsilon^4 &\colon
	\alpha=\delta=\tfrac{2}{3}\pi,\,
	\beta=(\tfrac{1}{3}+\tfrac{8}{f})\pi,\,
	\gamma=(\tfrac{5}{6}-\tfrac{4}{f})\pi,\,
	\epsilon=\tfrac{1}{2}\pi.
\end{align*}
We have $\delta>\epsilon$. By Lemma \ref{geometry1}, this implies $\beta<\gamma$. Then by $\beta\gamma^2,\alpha^3$, we get $\beta<\alpha<\gamma$.  

In the special companion pair $T_1,T_2$ in the second of Figure \ref{a2dA}, we have $\delta_2\delta_2\cdots=\alpha\delta^2$, and $\epsilon_1\epsilon_2\cdots$ has degree $4$ or $5$. We may assume $T_1$ is a special tile. Then $\epsilon_1\epsilon_2\cdots$ is the only high degree vertex in $T_1$, and by all the degree $3$ vertices, we get $\alpha_1\cdots,\beta_1\cdots,\gamma_1\cdots$ as indicated. This determines $T_3,T_4$ and gives $\alpha_5$. By all the degree $3$ vertices, we know $\beta_2\beta_3\cdots=\beta^2\cdots$ is a high degree vertex. Then $\alpha_2\cdots,\gamma_2\cdots$ have degree $3$. 

For $\epsilon^4$, by $\beta<\alpha=\delta<\gamma$ and Lemma \ref{geometry4}, we get $\beta>\epsilon$. This implies $\alpha,\beta,\gamma,\delta>\epsilon$. Then by $\epsilon^4$, this implies $\beta^2\cdots$ is not a (high degree) vertex. 

For $\delta\epsilon^3$, by $\alpha=\delta>\beta>\epsilon$, we get $\alpha+3\beta>2\beta+\delta+\epsilon>\alpha+\beta+2\epsilon>\delta+3\epsilon=2\pi$. Then by $\beta<\alpha,\gamma$, and $\delta>\epsilon$, and $\beta>\epsilon>\frac{2}{5}\pi$, we get $\beta^2\cdots=\beta^4,\beta^2\epsilon^2$. Since the degree $3$ vertex $\alpha_2\cdots=\alpha\delta^2,\alpha^3$, we know $\beta_2\beta_3\cdots\ne\beta^2\epsilon^2$. Therefore $\beta_2\beta_3\cdots=\beta^4$ is a vertex. This implies $\beta^2\epsilon^2$ is not a vertex, and $\beta^2\cdots=\beta^4$.

By $\delta\epsilon^3$ and $\delta>\epsilon$, we know $\delta\thin\delta\cdots$ is not a vertex. This implies the AAD of $\beta^4$ is $\thin^{\alpha}\beta^{\delta}\thin^{\alpha}\beta^{\delta}\thin^{\alpha}\beta^{\delta}\thin^{\alpha}\beta^{\delta}\thin$. Then $\beta_2\beta_3\cdots=\beta^4$ implies the degree $3$ vertex $\alpha_2\cdots=\alpha\delta\cdots=\alpha\delta^2$. This determines $T_6$. Then the degree $3$ vertex $\gamma_2\cdots=\beta\gamma^2$ as indicated. This implies $\epsilon_1\epsilon_2\cdots=\alpha\epsilon^2\cdots, \epsilon^3\cdots$. By $\delta\epsilon^3$ and $\alpha=\delta$, we get $R(\alpha\epsilon^2)=\epsilon$, and $R(\epsilon^3)=\delta<\beta+\epsilon$. Then by $\epsilon<\alpha,\beta,\gamma,\delta$ and $\beta<\alpha,\gamma,\delta,2\epsilon$, we get $\epsilon_1\epsilon_2\cdots=\delta\epsilon^3$. This implies the degree $3$ vertex $\gamma_1\cdots=\beta\gamma^2$ as indicated. Combined with $\alpha_5$, we determine $T_5$. Then we get $\thin^{\delta}\beta_4^{\alpha}\thin^{\alpha}\beta_5^{\delta}\thin\cdots=\beta^4$, contradicting the AAD $\thin^{\alpha}\beta^{\delta}\thin^{\alpha}\beta^{\delta}\thin^{\alpha}\beta^{\delta}\thin^{\alpha}\beta^{\delta}\thin$ of $\beta^4$.

\subsubsection*{Case. $\text{AVC}_3=\{\alpha\delta^2,\beta\gamma^2,\alpha^2\gamma\}$}

The angle sums of the degree $3$ vertices, and the angle sum of $\delta\epsilon^3$ or $\epsilon^4$, and the angle sum for pentagon imply
\begin{align*}
\delta\epsilon^3 &\colon
\alpha=(\tfrac{5}{8}+\tfrac{3}{2f})\pi,\,
\beta=(\tfrac{1}{2}+\tfrac{6}{f})\pi,\,
\gamma=(\tfrac{3}{4}-\tfrac{3}{f})\pi,\,
\delta=(\tfrac{11}{16}-\tfrac{3}{4f})\pi,\,
\epsilon=(\tfrac{7}{16}+\tfrac{1}{4f})\pi. \\
\epsilon^4 &\colon
\alpha=(\tfrac{3}{5}+\tfrac{8}{5f})\pi,\,
\beta=(\tfrac{2}{5}+\tfrac{32}{5f})\pi,\,
\gamma=(\tfrac{4}{5}-\tfrac{16}{5f})\pi,\,
\delta=(\tfrac{7}{10}-\tfrac{4}{5f})\pi,\,
\epsilon=\tfrac{1}{2}\pi.
\end{align*}
We have $\delta>\epsilon$. By Lemma \ref{geometry1}, this implies $\beta<\gamma$. Then by $\beta\gamma^2,\alpha^2\gamma$, we get $\beta<\alpha<\gamma$. 

By all the degree $3$ vertices, a special tile is the third of Figure \ref{a2dA}. Then $\epsilon_1\cdots=\beta\epsilon^2\cdots,\gamma\epsilon^2\cdots$. By $\beta+\gamma+2\epsilon>\alpha+\beta+2\epsilon=(2+\frac{8}{f})\pi>2\pi$, we get $R(\beta\epsilon^2)<\alpha<2\beta,\gamma,2\delta,2\epsilon$, and $R(\gamma\epsilon^2)<\beta<\alpha,\gamma,2\delta,2\epsilon$. This implies $\epsilon_1\cdots=\beta^2\epsilon^2$. 

The angle sum of $\beta^2\epsilon^2$ further implies
\begin{align*}
\delta\epsilon^3 &\colon
	\alpha=\tfrac{16}{25}\pi,\,
	\beta=\tfrac{14}{25}\pi,\,
	\gamma=\tfrac{18}{25}\pi,\,
	\delta=\tfrac{17}{25}\pi,\,
	\epsilon=\tfrac{11}{25}\pi,\,
	f=100. \\
\epsilon^4 &\colon
	\alpha=\tfrac{5}{8}\pi,\,
	\beta=\epsilon=\tfrac{1}{2}\pi,\,
	\gamma=\tfrac{3}{4}\pi,\,
	\delta=\tfrac{11}{16}\pi,\,
	f=64.
\end{align*}
The angle values imply $\beta^2\cdots=\beta^2\epsilon^2$, and $\gamma\delta\cdots$ is not a vertex. The AAD $\thick^{\epsilon}\delta^{\beta}\thin^{\beta}\alpha^{\gamma}\thin^{\beta}\delta^{\epsilon}\thick$ of $\alpha\delta^2$ implies a vertex $\thin^{\alpha}\beta^{\delta}\thin^{\alpha}\beta^{\delta}\thin\cdots=\beta^2\epsilon^2=\thick^{\delta}\epsilon^{\gamma}\thin^{\alpha}\beta^{\delta}\thin^{\alpha}\beta^{\delta}\thin^{\gamma}\epsilon^{\delta}\thick$. Then we get $\gamma\delta\cdots$, a contradiction. 

\subsubsection*{Case. $\text{AVC}_3\equiv\{\alpha\delta^2,\gamma^3,\alpha\beta\gamma\}$}

The angle sums of the degree $3$ vertices, and the angle sum of $\delta\epsilon^3$ or $\epsilon^4$, and the angle sum for pentagon imply
\begin{align*}
\delta\epsilon^3 &\colon
\alpha=(1-\tfrac{12}{f})\pi,\,
\beta=(\tfrac{1}{3}+\tfrac{12}{f})\pi,\,
\gamma=\tfrac{2}{3}\pi,\,
\delta=(\tfrac{1}{2}+\tfrac{6}{f})\pi,\,
\epsilon=(\tfrac{1}{2}-\tfrac{2}{f})\pi. \\
\epsilon^4 &\colon
\alpha=(1-\tfrac{8}{f})\pi,\,
\beta=(\tfrac{1}{3}+\tfrac{8}{f})\pi,\,
\gamma=\tfrac{2}{3}\pi,\,
\delta=(\tfrac{1}{2}+\tfrac{4}{f})\pi,\,
\epsilon=\tfrac{1}{2}\pi.
\end{align*}
We have $\delta>\epsilon$ and $\gamma<2\epsilon$. By Lemma \ref{geometry1}, this implies $\beta<\gamma$. Then by $\gamma^3,\alpha\beta\gamma$, we get $\beta<\gamma<\alpha$.

The AAD of $\gamma^3$ implies a vertex $\alpha\epsilon\cdots$. By $\alpha\beta\gamma$, and $\gamma<2\epsilon$, and $\delta>\epsilon$, we get $R(\alpha\delta\epsilon)<R(\alpha\epsilon^2)<\beta<\alpha,\gamma,2\delta,2\epsilon$. Then by $\alpha\delta^2$, this implies $\alpha\epsilon\cdots$ is not a vertex, a contradiction. 
\end{proof}

\begin{proposition}\label{b2d}
There is no tiling, such that $\alpha,\beta,\gamma$ have distinct values, and $\beta\delta^2$ is the only degree $3$ $b$-vertex.
\end{proposition}

\begin{proof}
By exchanging $\alpha$ and $\beta$, the initial discussion in the proof of Proposition \ref{a2d} can be adopted. We find that the collection $\text{AVC}_3$ of all degree $3$ vertices is exactly one of the following:
\[
\begin{array}{lll}
\{\beta\delta^2,\alpha^3,\alpha\beta\gamma\}, &
\{\beta\delta^2,\alpha^3,\beta^2\gamma\}, &
\{\beta\delta^2,\alpha^3,\beta\gamma^2\},  \\
\{\beta\delta^2,\alpha^2\gamma,\alpha\beta^2\}, &
\{\beta\delta^2,\alpha^2\gamma,\beta^3\}, &
\{\beta\delta^2,\alpha^2\gamma,\beta\gamma^2\}, \\
\{\beta\delta^2,\alpha\gamma^2,\alpha^2\beta\}, &
\{\beta\delta^2,\alpha\gamma^2,\beta^3\}, &
\{\beta\delta^2,\alpha\gamma^2,\beta^2\gamma\}, \\ 
\{\beta\delta^2,\gamma^3,\alpha\beta\gamma\}, &
\{\beta\delta^2,\gamma^3,\alpha^2\beta\}, &
\{\beta\delta^2,\gamma^3,\alpha\beta^2\}.
\end{array}
\]
In a special tile, the degree $3$ vertex $\delta\cdots=\beta\delta^2$, and the vertices $\alpha\cdots$, $\beta\cdots$, $\gamma\cdots$ have degree $3$. Note that $\delta\cdots=\beta\delta^2$ implies the degree $3$ vertex $\beta\cdots=\beta\delta^2$. This further implies the degree $3$ vertex $\alpha\cdots=\alpha\beta\cdots$. Therefore $\{\beta\delta^2,\alpha^3,\beta^2\gamma\}$, $\
\{\beta\delta^2,\alpha^3,\beta\gamma^2\}$, $\{\beta\delta^2,\alpha^2\gamma,\beta^3\}$, $
\{\beta\delta^2,\alpha^2\gamma,\beta\gamma^2\}$, $\{\beta\delta^2,\alpha\gamma^2,\beta^3\}$, $
\{\beta\delta^2,\alpha\gamma^2,\beta^2\gamma\}$ cannot be the $\text{AVC}_3$.

For $\text{AVC}_3=\{\beta\delta^2,\gamma^3,\alpha\beta^2\}$, in a special tile, one of $\alpha\cdots,\gamma\cdots$, has high degree, a contradiction.

We recall that $\delta\epsilon^3$ or $\epsilon^4$ is a vertex, and $\gamma^2\cdots$ is a vertex.

\subsubsection*{Case. $\text{AVC}_3= \{\beta\delta^2,\alpha^3,\alpha\beta\gamma\}$}

The angle sums of the degree $3$ vertices, and the angle sum of $\delta\epsilon^3$ or $\epsilon^4$, and the angle sum for pentagon imply
\begin{align*}
\delta\epsilon^3 &\colon
	\alpha=\tfrac{2}{3}\pi,\,
	\beta=(1-\tfrac{12}{f})\pi,\,
	\gamma=(\tfrac{1}{3}+\tfrac{12}{f})\pi,\,
	\delta=(\tfrac{1}{2}+\tfrac{6}{f})\pi,\,
	\epsilon=(\tfrac{1}{2}-\tfrac{2}{f})\pi. \\
\epsilon^4 &\colon
	\alpha=\tfrac{2}{3}\pi,\,
	\beta=(1-\tfrac{8}{f})\pi,\,
	\gamma=(\tfrac{1}{3}+\tfrac{8}{f})\pi,\,
	\delta=(\tfrac{1}{2}+\tfrac{4}{f})\pi,\,
	\epsilon=\tfrac{1}{2}\pi.
\end{align*}
We have $\delta>\epsilon$ and $\beta<2\epsilon$. By Lemma \ref{geometry1}, this implies $\beta<\gamma$. Then by $\alpha^3,\alpha\beta\gamma$, we get $\beta<\alpha<\gamma$. 

By $\alpha\beta\gamma$ and $\alpha<\gamma$, we get $R(\gamma^2)<\beta<\alpha,\gamma,2\delta,2\epsilon$. This implies $\gamma^2\cdots$ is not a vertex, a contradiction.

\subsubsection*{Case. $\text{AVC}_3= \{\beta\delta^2,\alpha^2\gamma,\alpha\beta^2\}$}

The angle sums of the degree $3$ vertices, and the angle sum of $\delta\epsilon^3$ or $\epsilon^4$, and the angle sum for pentagon imply
\begin{align*}
\delta\epsilon^3 &\colon
	\alpha=(\tfrac{3}{4}-\tfrac{3}{f})\pi,\,
	\beta=(\tfrac{5}{8}+\tfrac{3}{2f})\pi,\,
	\gamma=(\tfrac{1}{2}+\tfrac{6}{f})\pi,\,
	\delta=(\tfrac{11}{16}-\tfrac{3}{4f})\pi,\,
	\epsilon=(\tfrac{7}{16}+\tfrac{1}{4f})\pi. \\
\epsilon^4 &\colon
	\alpha=(\tfrac{4}{5}-\tfrac{16}{5f})\pi,\,
	\beta=(\tfrac{3}{5}+\tfrac{8}{5f})\pi,\,
	\gamma=(\tfrac{2}{5}+\tfrac{32}{5f})\pi,\,
	\delta=(\tfrac{7}{10}-\tfrac{4}{5f})\pi,\,
	\epsilon=\tfrac{1}{2}\pi.
\end{align*}
We have $\delta>\epsilon$ and $\alpha<2\epsilon$. By Lemma \ref{geometry1}, this implies $\beta<\gamma$. Then by $\alpha^2\gamma,\alpha\beta^2$, we get $\alpha<\beta<\gamma$. 

By $\alpha^2\gamma$ and $\alpha<\gamma$, we get $R(\gamma^2)<\alpha<\beta,\gamma,2\delta,2\epsilon$. This implies $\gamma^2\cdots$ is not a vertex, a contradiction.

\subsubsection*{Case. $\text{AVC}_3= \{\beta\delta^2,\alpha\gamma^2,\alpha^2\beta\}$}

The angle sums of the degree $3$ vertices, and the angle sum of $\delta\epsilon^3$ or $\epsilon^4$, and the angle sum for pentagon imply
\begin{align*}
\delta\epsilon^3 &\colon
	\alpha=\delta=(\tfrac{4}{5}-\tfrac{24}{5f})\pi,\,
	\beta=(\tfrac{2}{5}+\tfrac{48}{5f})\pi,\,
	\gamma=(\tfrac{3}{5}+\tfrac{12}{5f})\pi,\,
	\epsilon=(\tfrac{2}{5}+\tfrac{8}{5f})\pi. \\
\epsilon^4 &\colon 
	\alpha=\delta=(1-\tfrac{8}{f})\pi,\,
	\beta=\tfrac{16}{f}\pi,\,
	\gamma=(\tfrac{1}{2}+\tfrac{4}{f})\pi,\,
	\epsilon=\tfrac{1}{2}\pi.
\end{align*}
By $f\ge 24$, we get $\delta>\epsilon$ and $\alpha=\delta<2\epsilon$. By Lemma \ref{geometry1}, this implies $\beta<\gamma$. Then by $\alpha\gamma^2,\alpha^2\beta$, we get $\beta<\gamma<\alpha$.

By $\delta\epsilon^3$ or $\epsilon^4$, and $\delta>\epsilon$, we know $\delta\thin\delta\cdots$ is not a vertex, and $\delta\thin\epsilon\cdots=\delta\epsilon^3$. 

By $\alpha^2\beta$, we get $R(\alpha^2)=\beta<\alpha,\gamma,\delta,2\epsilon$. This implies $\alpha^2\cdots=\alpha^2\beta$.

By $\beta\delta^2$ and $\alpha>\beta$, we know a $b$-vertex $\alpha\cdots=\alpha\delta\epsilon\cdots,\alpha\epsilon^2\cdots$. By $\alpha^2\beta$, and $\alpha<2\epsilon$, and $\delta>\epsilon$, we get $R(\alpha\delta\epsilon)<R(\alpha\epsilon^2)<\beta<\alpha,\gamma,\delta,2\epsilon$. This implies $\alpha\cdots$ is a $\hat{b}$-vertex. 

By $\beta\delta^2$ and $\beta<\gamma$, we know $\gamma\delta\cdots=\gamma\delta\epsilon\cdots$. By $\alpha<\pi$ and the angle sum for pentagon, we get $\beta+\gamma+\delta+\epsilon>2\pi$. This implies $R(\gamma\delta\epsilon)<\beta<\alpha,\gamma,\delta,2\epsilon$. Therefore $\gamma\delta\cdots$ is not a vertex.

By no $\alpha\delta\cdots,\delta\thin\delta\cdots$, we know the AAD of $\thin\beta\thin\beta\thin$ is $\thin^{\delta}\beta^{\alpha}\thin^{\alpha}\beta^{\delta}\thin$. This implies a vertex $\thin^{\gamma}\alpha^{\beta}\thin^{\beta}\alpha^{\gamma}\thin\cdots=\alpha^2\beta=\thin^{\beta}\alpha^{\gamma}\thin^{\alpha}\beta^{\delta}\thin^{\gamma}\alpha^{\beta}\thin$, contradicting no $\gamma\delta\cdots$. Therefore $\beta\thin\beta\cdots$ is not a vertex. 

By $\beta\delta^2$ and $\delta<2\epsilon$, we know $R(\beta\delta\epsilon)$ has no $\delta,\epsilon$. By no $\alpha\delta\cdots,\gamma\delta\cdots$, we know $R(\beta\delta\epsilon)$ has no $\alpha,\gamma$. Therefore $\beta\delta\epsilon\cdots=\beta^k\delta\epsilon(k\ge 2)$. Then by no $\beta\thin\beta\cdots$, we know $\beta\delta\epsilon\cdots$ is not a vertex, and a $b$-vertex $\beta\cdots=\beta\delta^2$.

By $\alpha\gamma^2,\alpha^2\beta$, a $\hat{b}$ vertex $\beta\cdots=\alpha^2\beta,\alpha\beta^k,\alpha\beta^k\gamma,\beta^k\gamma^l$. By no $\beta\thin\beta\cdots$, we know $\alpha\beta^k$ is not a vertex, and we get $\thin\gamma\thin\beta\thin\gamma\thin$ in $\alpha\beta^k\gamma,\beta^k\gamma^l$. By no $\alpha\delta\cdots$, the AAD of $\thin\gamma\thin\beta\thin\gamma\thin$ contains $\thin^{\alpha}\beta^{\delta}\thin^{\epsilon}\gamma^{\alpha}\thin$. This implies a vertex $\delta\thin\epsilon\cdots=\delta\epsilon^3$. By the angle values for the case $\delta\epsilon^3$, we get $\beta+\gamma=(1+\frac{12}{f})\pi>\pi$. We know $\alpha\beta^k\gamma,\beta^k\gamma^l$ are high degree vertices, and $k\le l$ in $\beta^k\gamma^l$ (by no $\beta\thin\beta\cdots$). Then by $\beta+\gamma>\pi$ and $\gamma<\alpha$, we know $\alpha\beta^k\gamma,\beta^k\gamma^l$ are not vertices.

We conclude $\beta\cdots=\beta\delta^2,\alpha^2\beta$. Then we get
\[
f =\#\beta
=\#\beta\delta^2+\#\alpha^2\beta 
\le \tfrac{1}{2}\#\delta+\tfrac{1}{2}\#\alpha
=\tfrac{1}{2}f+\tfrac{1}{2}f=f.
\]
This implies $\tfrac{1}{2}\#\alpha=\#\alpha^2\beta$, and further implies $\alpha\cdots=\alpha^2\beta$, contradicting the vertex $\alpha\gamma^2$.

\subsubsection*{Case. $\text{AVC}_3= \{\beta\delta^2,\gamma^3,\alpha^2\beta\}$}

Suppose $\epsilon^4$ is a vertex. The angle sums of $\beta\delta^2,\gamma^3,\alpha^2\beta,\epsilon^4$ and the angle sum for pentagon imply $f=24$. Since there is no $3^5$-tile, by Lemma \ref{special_tile}, we know every tile has exactly one vertex of degree $4$. Since $\epsilon\cdots$ has high degree, this implies that, in every tile, the vertex $\alpha\cdots$ has degree $3$. Since any $\alpha$ belongs to some tile, we know any vertex $\alpha\cdots$ has degree $3$. On the other hand, the AAD of $\gamma^3$ implies $\alpha\epsilon\cdots$ is a vertex. This has high degree, a contradiction. 

Therefore $\delta\epsilon^3$ is a vertex. The angle sums of $\beta\delta^2,\gamma^3,\alpha^2\beta,\delta\epsilon^3$ and the angle sum for pentagon imply
\[
	\alpha=\delta=(1-\tfrac{12}{f})\pi,\,
	\beta=\tfrac{24}{f}\pi,\,
	\gamma=\tfrac{2}{3}\pi,\,
	\epsilon=(\tfrac{1}{3}+\tfrac{4}{f})\pi.
\]
By $f\ge 24$ and $\delta\ne\epsilon$, we get $\delta>\epsilon$. By Lemma \ref{geometry1}, this implies $\beta<\gamma<2\epsilon$. Then by $\gamma^3,\alpha^2\beta$, we get $\beta<\gamma<\alpha$. 

By $\delta\epsilon^3$ and $\delta>\epsilon$, we know $\delta\thin\delta\cdots$ is not a vertex. 

By $\alpha^2\beta$, we get $R(\alpha^2)=\beta<\alpha,\gamma,\delta,2\epsilon$. This implies $\alpha^2\cdots=\alpha^2\beta$.

By $\beta\delta^2$ and $\beta<\alpha,\gamma$, we know $\alpha\delta\cdots=\alpha\delta\epsilon\cdots$ and $\gamma\delta\cdots=\gamma\delta\epsilon\cdots$. By $\alpha<\pi$ and the angle sum for pentagon, we get $\beta+\gamma+\delta+\epsilon>2\pi$. Then by $\alpha>\gamma$, this implies $R(\alpha\delta\epsilon)<R(\gamma\delta\epsilon)<\beta<\alpha,\gamma,\delta,2\epsilon$. Therefore $\alpha\delta\cdots,\gamma\delta\cdots$ are not vertices.

By no $\alpha\delta\cdots,\delta\thin\delta\cdots$, we know the AAD of $\thin\beta\thin\beta\thin$ is $\thin^{\delta}\beta^{\alpha}\thin^{\alpha}\beta^{\delta}\thin$. This implies a vertex $\thin^{\gamma}\alpha^{\beta}\thin^{\beta}\alpha^{\gamma}\thin\cdots=\alpha^2\beta=\thin^{\beta}\alpha^{\gamma}\thin^{\alpha}\beta^{\delta}\thin^{\gamma}\alpha^{\beta}\thin$, contradicting no $\gamma\delta\cdots$. Therefore $\beta\thin\beta\cdots$ is not a vertex. 

The AAD of $\gamma^3$ implies a vertex $\alpha\epsilon\cdots$. By no $\alpha\delta\cdots$, we get $\alpha\epsilon\cdots=\alpha\epsilon^2\cdots$. Then by $\gamma^3$ and $\gamma<\alpha,\delta,2\epsilon$, we get $\alpha\epsilon^2\cdots=\alpha\beta^k\epsilon^2(k\ge 1)$. Then by no $\beta\thin\beta\cdots,\gamma\delta\cdots$, we get $\alpha\epsilon\cdots=\alpha\beta^k\epsilon^2=\thick^{\delta}\epsilon^{\gamma}\thin^{\alpha}\beta^{\delta}\thin^{\beta}\alpha^{\gamma}\thin^{\gamma}\epsilon^{\delta}\thick=\alpha\beta\epsilon^2$. This implies a vertex $\alpha\gamma\cdots$. By $\alpha\epsilon\cdots=\alpha\beta\epsilon^2$ and no $\alpha\delta\cdots$, we know $\alpha\gamma\cdots$ is a $\hat{b}$-vertex. Then by $\gamma^3$ and $\beta<\gamma<\alpha$, we get $\alpha\gamma\cdots=\alpha\beta^k\gamma$. Then by no $\beta\thin\beta\cdots$, we get $\alpha\beta^k\gamma=\thin\alpha\thin\beta\thin\gamma\thin\beta\thin=\alpha\beta^2\gamma$. However, the angle sums of $\alpha\beta\epsilon^2,\alpha\beta^2\gamma$ imply $\beta+\gamma=2\epsilon$, contradicting the formulae for $\beta,\gamma,\epsilon$. 

\subsubsection*{Case. $\text{AVC}_3= \{\beta\delta^2,\gamma^3,\alpha\beta\gamma\}$}

The angle sums of the degree $3$ vertices, and the angle sum of $\delta\epsilon^3$ or $\epsilon^4$, and the angle sum for pentagon imply
\begin{align*}
\delta\epsilon^3 &\colon
	\alpha=(\tfrac{1}{3}+\tfrac{12}{f})\pi,\,
	\beta=(1-\tfrac{12}{f})\pi,\,
	\gamma=\tfrac{2}{3}\pi,\,
	\delta=(\tfrac{1}{2}+\tfrac{6}{f})\pi,\,
	\epsilon=(\tfrac{1}{2}-\tfrac{2}{f})\pi. \\
\epsilon^4 &\colon
	\alpha=(\tfrac{1}{3}+\tfrac{8}{f})\pi,\,
	\beta=(1-\tfrac{8}{f})\pi,\,
	\gamma=\tfrac{2}{3}\pi,\,
	\delta=(\tfrac{1}{2}+\tfrac{4}{f})\pi,\,
	\epsilon=\tfrac{1}{2}\pi.
\end{align*}
We have $\delta>\epsilon$ and $\beta<2\epsilon$. By Lemma \ref{geometry1}, this implies $\beta<\gamma$. Then by $\gamma^3,\alpha\beta\gamma$, we get $\beta<\gamma<\alpha$. 

The AAD of $\gamma^3$ implies a vertex $\alpha\epsilon\cdots$. By $\alpha\beta\gamma$, and $\gamma<2\epsilon$, and $\delta>\epsilon$, we get $R(\alpha\delta\epsilon)<R(\alpha\epsilon^2)<\beta<\alpha,\gamma,2\delta,2\epsilon$. This implies $\alpha\epsilon\cdots$ is not a vertex, a contradiction. 
\end{proof}

\begin{proposition}\label{c2d}
There is no tiling, such that $\alpha,\beta,\gamma$ have distinct values, and $\gamma\delta^2$ is the only degree $3$ $b$-vertex.
\end{proposition}

\begin{proof}
By exchanging $\alpha$ and $\gamma$, the initial discussion in the proof of Proposition \ref{a2d} can be adopted. We find that that the collection $\text{AVC}_3$ of all degree $3$ vertices is exactly one of the following:
\[
\begin{array}{lll}
\{\gamma\delta^2,\alpha^3,\alpha\beta\gamma\}, &
\{\gamma\delta^2,\alpha^3,\beta^2\gamma\}, &
\{\gamma\delta^2,\alpha^3,\beta\gamma^2\},  \\
\{\gamma\delta^2,\alpha^2\beta,\alpha\gamma^2\}, &
\{\gamma\delta^2,\alpha^2\beta,\beta^2\gamma\}, &
\{\gamma\delta^2,\alpha^2\beta,\gamma^3\}, \\
\{\gamma\delta^2,\alpha\beta^2,\alpha^2\gamma\}, &
\{\gamma\delta^2,\alpha\beta^2,\beta\gamma^2\}, &
\{\gamma\delta^2,\alpha\beta^2,\gamma^3\}, \\ 
\{\gamma\delta^2,\beta^3,\alpha\beta\gamma\}, &
\{\gamma\delta^2,\beta^3,\alpha^2\gamma\}, &
\{\gamma\delta^2,\beta^3,\alpha\gamma^2\}.
\end{array}
\]
In a special tile, the degree $3$ vertex $\delta\cdots=\gamma\delta^2$, and the vertices $\alpha\cdots$, $\beta\cdots$, $\gamma\cdots$ have degree $3$. Note that $\delta\cdots=\gamma\delta^2$ implies the degree $3$ vertex $\beta\cdots=\alpha\beta\cdots,\beta\epsilon\cdots$. Therefore $\{\gamma\delta^2,\alpha^3,\beta^2\gamma\}$, $\{\gamma\delta^2,\alpha^3,\beta\gamma^2\}$, $
\{\gamma\delta^2,\beta^3,\alpha^2\gamma\}$, $
\{\gamma\delta^2,\beta^3,\alpha\gamma^2\}$ cannot be the $\text{AVC}_3$. 

For $\text{AVC}_3=\{\gamma\delta^2,\alpha\beta^2,\beta\gamma^2\}$ or $
\{\gamma\delta^2,\alpha\beta^2,\gamma^3\}$, in a special tile, by $\delta\cdots=\gamma\delta^2$, we get $\beta\cdots=\alpha\beta^2$, with specific angle arrangement at the vertex. Then $\alpha\cdots=\alpha^2\cdots,\alpha\delta\cdots$ has high degree, a contradiction.

For $\text{AVC}_3=\{\gamma\delta^2,\beta^3,\alpha\beta\gamma\}$, in a special tile, by $\delta\cdots=\gamma\delta^2$, we get $\beta\cdots=\alpha\beta\gamma$, with specific angle arrangement at the vertex. Then $\alpha\cdots=\alpha^2\cdots,\alpha\epsilon\cdots$ has high degree, a contradiction.

We recall that $\delta\epsilon^3$ or $\epsilon^4$ is a vertex, and $\gamma^2\cdots$ is a vertex.

\subsubsection*{Case. $\text{AVC}_3= \{\gamma\delta^2,\alpha^3,\alpha\beta\gamma\}$}

The angle sums of the degree $3$ vertices, and the angle sum of $\delta\epsilon^3$ or $\epsilon^4$, and the angle sum for pentagon imply
\begin{align*}
\delta\epsilon^3 &\colon
	\alpha=\tfrac{2}{3}\pi,\,
	\beta=(\tfrac{1}{3}+\tfrac{12}{f})\pi,\,
	\gamma=(1-\tfrac{12}{f})\pi,\,
	\delta=(\tfrac{1}{2}+\tfrac{6}{f})\pi,\,
	\epsilon=(\tfrac{1}{2}-\tfrac{2}{f})\pi. \\
\epsilon^4 &\colon
	\alpha=\tfrac{2}{3}\pi,\,
	\beta=(\tfrac{1}{3}+\tfrac{8}{f})\pi,\,
	\gamma=(1-\tfrac{8}{f})\pi,\,
	\delta=(\tfrac{1}{2}+\tfrac{4}{f})\pi,\,
	\epsilon=\tfrac{1}{2}\pi.
\end{align*}
By $f\ge 24$, we get $\delta>\epsilon$ and $\beta\le 2\epsilon$. By Lemma \ref{geometry1}, this implies $\beta<\gamma$. Then by $\alpha^3,\alpha\beta\gamma$, we get $\beta<\alpha<\gamma$. 

By $\alpha\beta\gamma$ and $\alpha<\gamma$, we get $R(\gamma^2)<\beta<\alpha,\gamma,2\delta,2\epsilon$. This implies $\gamma^2\cdots$ is not a vertex, a contradiction.

\subsubsection*{Case. $\text{AVC}_3= \{\gamma\delta^2,\alpha^2\beta,\alpha\gamma^2\}$}

The angle sums of the degree $3$ vertices, and the angle sum of $\delta\epsilon^3$ or $\epsilon^4$, and the angle sum for pentagon imply
\begin{align*}
\delta\epsilon^3 &\colon
	\alpha=(\tfrac{3}{4}-\tfrac{3}{f})\pi,\,
	\beta=(\tfrac{1}{2}+\tfrac{6}{f})\pi,\,
	\gamma=(\tfrac{5}{8}+\tfrac{3}{2f})\pi,\,
	\delta=(\tfrac{11}{16}-\tfrac{3}{4f})\pi,\,
	\epsilon=(\tfrac{7}{16}+\tfrac{1}{4f})\pi. \\
\epsilon^4 &\colon
	\alpha=(\tfrac{4}{5}-\tfrac{16}{5f})\pi,\,
	\beta=(\tfrac{2}{5}+\tfrac{32}{5f})\pi,\,
	\gamma=(\tfrac{3}{5}+\tfrac{8}{5f})\pi,\,
	\delta=(\tfrac{7}{10}-\tfrac{4}{5f})\pi,\,
	\epsilon=\tfrac{1}{2}\pi.
\end{align*}
We have $\delta>\epsilon$ and $\alpha<2\beta,2\epsilon$. By Lemma \ref{geometry1}, this implies $\beta<\gamma$. Then by $\alpha^2\beta,\alpha\gamma^2$, we get $\beta<\gamma<\alpha$. 

The AAD $\thick^{\epsilon}\delta^{\beta}\thin^{\alpha}\gamma^{\epsilon}\thin^{\beta}\delta^{\epsilon}\thick$ of $\gamma\delta^2$ implies a vertex $\thin^{\alpha}\beta^{\delta}\thin^{\gamma}\epsilon^{\delta}\thick\cdots$. By the angle values, we get $2\beta+\delta+\epsilon>2\pi$ and $\beta+\gamma+2\epsilon>2\pi$. Then we get $R(\beta\delta\epsilon)<\beta<\alpha,\gamma,2\delta,2\epsilon$ and $R(\beta\epsilon^2)<\gamma<\alpha,2\beta,2\delta,2\epsilon$. This implies $\beta\epsilon\cdots=\beta^2\epsilon^2$. The angle sum of $\beta^2\epsilon^2$ further implies
\begin{align*}
\delta\epsilon^3 &\colon
	\alpha=\tfrac{18}{25}\pi,\,
	\beta=\tfrac{14}{25}\pi,\,
	\gamma=\tfrac{16}{25}\pi,\,
	\delta=\tfrac{17}{25}\pi,\,
	\epsilon=\tfrac{11}{25}\pi, \,
	f=100. \\
\epsilon^4 &\colon
	\alpha=\tfrac{3}{4}\pi,\,
	\beta=\epsilon=\tfrac{1}{2}\pi,\,
	\gamma=\tfrac{5}{8}\pi,\,
	\delta=\tfrac{11}{16}\pi,\,
	f=64.
\end{align*}

The AAD of $\delta\epsilon^3$ implies a vertex $\beta\gamma\cdots$. Since $R(\beta\gamma)=\frac{20}{25}\pi$ is not a combination of angle values, we get a contradiction. 

The angle values for $\epsilon^4$ imply $\alpha^2\cdots=\alpha\beta\cdots=\alpha^2\beta$, and $\gamma\delta\cdots=\gamma\delta^2$, and $\alpha\delta\cdots$ is not a vertex. We also know $\beta\epsilon\cdots=\beta^2\epsilon^2$. By no $\alpha\delta\cdots$, the AAD of $\thin^{\alpha}\beta^{\delta}\thin^{\gamma}\epsilon^{\delta}\thick\cdots=\beta^2\epsilon^2$ is $\thick^{\delta}\epsilon^{\gamma}\thin^{\delta}\beta^{\alpha}\thin^{\alpha}\beta^{\delta}\thin^{\gamma}\epsilon^{\delta}\thick$. This determines $T_1,T_2,T_3,T_4$ in Figure \ref{c2dA}. Then $\alpha_1\alpha_3\cdots=\alpha^2\beta$ gives $\beta_5$. Up to the horizontal flip symmetry, we may assume $T_5$ is arranged as indicated. Then $\gamma_1\delta_5\cdots=\gamma_2\delta_1\cdots=\gamma\delta^2$ determines $T_6,T_7$. Then $\beta_6\epsilon_1\epsilon_7\cdots=\beta^2\epsilon^2$ and no $\alpha\delta\cdots$ determine $T_8$. Then $\gamma_7\delta_8\cdots=\gamma\delta^2$ determines $T_9$. Then $\alpha_2\beta_7\cdots=\alpha_7\beta_9\cdots=\alpha^2\beta$ implies two $\alpha$ adjacent, a contradiction. 

\begin{figure}[htp]
\centering
\begin{tikzpicture}[>=latex,scale=1]

\foreach \a in {1,-1}
{
\begin{scope}[xscale=\a]

\draw
	(0,-0.7) -- (0,0.7) -- (0.5,1.1) -- (1,0.7) -- (1,-0.7) -- (0.5,-1.1) -- (0,-0.7)
	(0,0) -- (1,0);

\draw[line width=1.2]
	(0,0) -- (0,-0.7)
	(1,0) -- (1,0.7);
	
\node at (0.2,0.6) {\small $\alpha$};
\node at (0.2,0.2) {\small $\beta$};
\node at (0.5,0.85) {\small $\gamma$};	
\node at (0.8,0.2) {\small $\delta$}; 
\node at (0.8,0.6) {\small $\epsilon$};

\node at (0.2,-0.6) {\small $\delta$};
\node at (0.5,-0.85) {\small $\beta$};
\node at (0.2,-0.2) {\small $\epsilon$};	
\node at (0.8,-0.2) {\small $\gamma$}; 
\node at (0.8,-0.6) {\small $\alpha$};

\end{scope}
}

\draw
	(1,0.7) -- (1.8,0.3) -- (2.5,0.3) 
	(1.8,0.3) -- (1.8,-0.7)
	(1,0.7) -- (1.5,1.1) -- (2.5,1.1) -- (2.5,0.3) --(2.9,-0.2) -- (2.5,-0.7) -- (1,-0.7)
	(-0.5,1.1) -- (-0.5,1.8) -- (1.5,1.8) -- (1.5,1.1);

\draw[line width=1.2]
	(0.5,1.1) -- (0.5,1.8)
	(1.8,0.3) -- (2.5,0.3);

\node at (1.1,-0.85) {\small $\alpha$};
\node at (1.8,-0.85) {\small $\alpha$};

\node at (1.6,0.2) {\small $\gamma$};
\node at (1.6,-0.5) {\small $\alpha$};
\node at (1.2,0.4) {\small $\epsilon$};
\node at (1.2,-0.5) {\small $\beta$};
\node at (1.2,0) {\small $\delta$};

\node at (2,0.1) {\small $\delta$};
\node at (2,-0.5) {\small $\beta$};
\node at (2.45,0.1) {\small $\epsilon$};
\node at (2.45,-0.5) {\small $\alpha$};
\node at (2.65,-0.2) {\small $\gamma$};

\node at (-0.3,1.2) {\small $\alpha$};
\node at (0,0.95) {\small $\beta$};
\node at (-0.3,1.6) {\small $\gamma$};
\node at (0.3,1.2) {\small $\delta$};
\node at (0.3,1.6) {\small $\epsilon$};

\node at (1.3,0.75) {\small $\beta$};
\node at (1.8,0.5) {\small $\delta$};
\node at (1.6,0.95) {\small $\alpha$};
\node at (2.3,0.5) {\small $\epsilon$};
\node at (2.3,0.9) {\small $\gamma$};

\node at (1,0.95) {\small $\beta$};
\node at (0.7,1.2) {\small $\delta$};
\node at (1.3,1.2) {\small $\alpha$};
\node at (0.7,1.6) {\small $\epsilon$};
\node at (1.3,1.65) {\small $\gamma$};

\node[inner sep=0.5,draw,shape=circle] at (0.5,0.4) {\small $1$};
\node[inner sep=0.5,draw,shape=circle] at (0.5,-0.4) {\small $2$};
\node[inner sep=0.5,draw,shape=circle] at (-0.5,0.4) {\small $3$};
\node[inner sep=0.5,draw,shape=circle] at (-0.5,-0.4) {\small $4$};
\node[inner sep=0.5,draw,shape=circle] at (0,1.4) {\small $5$};
\node[inner sep=0.5,draw,shape=circle] at (1.5,-0.15)  {\small $7$};
\node[inner sep=0.5,draw,shape=circle] at (1,1.4) {\small 6};
\node[inner sep=0.5,draw,shape=circle] at (1.95,0.8) {\small 8};
\node[inner sep=0.5,draw,shape=circle] at (2.25,-0.2) {\small $9$};

\end{tikzpicture}
\caption{Proposition \ref{c2d}: $\text{AVC}_3= \{\gamma\delta^2,\alpha^2\beta,\alpha\gamma^2\}$.}
\label{c2dA}
\end{figure}
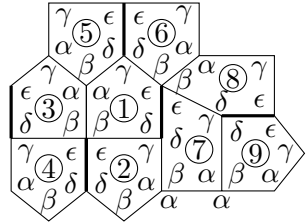

\subsubsection*{Case. $\text{AVC}_3= \{\gamma\delta^2,\alpha^2\beta,\beta^2\gamma\}$}

The angle sums of the degree $3$ vertices, and the angle sum of $\delta\epsilon^3$ or $\epsilon^4$, and the angle sum for pentagon imply
\begin{align*}
\delta\epsilon^3 &\colon
	\alpha=(\tfrac{3}{5}+\tfrac{12}{5f})\pi,\,
	\beta=\delta=(\tfrac{4}{5}-\tfrac{24}{5f})\pi,\,
	\gamma=(\tfrac{2}{5}+\tfrac{48}{5f})\pi,\,
	\epsilon=(\tfrac{2}{5}+\tfrac{8}{5f})\pi. \\
\epsilon^4 &\colon
	\alpha=(\tfrac{1}{2}+\tfrac{4}{f})\pi,\,
	\beta=\delta=(1-\tfrac{8}{f})\pi,\,
	\gamma=\tfrac{16}{f}\pi,\,
	\epsilon=\tfrac{1}{2}\pi.
\end{align*}
By $f\ge 24$, we get $\delta>\epsilon$, and $\beta<2\epsilon$, and $\beta+\epsilon>\pi$. By Lemma \ref{geometry1}, this implies $\beta<\gamma$. Then by $\alpha^2\beta,\beta^2\gamma$, we get $\beta<\alpha<\gamma$. 

The AAD $\thick^{\epsilon}\delta^{\beta}\thin^{\alpha}\gamma^{\epsilon}\thin^{\beta}\delta^{\epsilon}\thick$ of $\gamma\delta^2$ implies a vertex $\beta\epsilon\cdots$. By $\beta+\epsilon>\pi$ and $\delta>\epsilon$, we get $R(\beta\delta\epsilon)<R(\beta\epsilon^2)<\beta<\alpha,\gamma,2\delta,2\epsilon$. This implies $\beta\epsilon\cdots$ is not a vertex, a contradiction.

\subsubsection*{Case. $\text{AVC}_3= \{\gamma\delta^2,\alpha^2\beta,\gamma^3\}$}

The angle sums of the degree $3$ vertices, and the angle sum of $\delta\epsilon^3$ or $\epsilon^4$, and the angle sum for pentagon imply
\begin{align*}
\delta\epsilon^3 &\colon
	\alpha=(\tfrac{7}{9}-\tfrac{4}{f})\pi,\,
	\beta=(\tfrac{4}{9}+\tfrac{8}{f})\pi,\,
	\gamma=\delta=\tfrac{2}{3}\pi,\,
	\epsilon=\tfrac{4}{9}\pi. \\
\epsilon^4 &\colon
	\alpha=(\tfrac{5}{6}-\tfrac{4}{f})\pi,\,
	\beta=(\tfrac{1}{3}+\tfrac{8}{f})\pi,\,
	\gamma=\delta=\tfrac{2}{3}\pi,\,
	\epsilon=\tfrac{1}{2}\pi.
\end{align*}
We have $\delta>\epsilon$ and $\gamma<2\epsilon$. By Lemma \ref{geometry1}, this implies $\beta<\gamma$. Then by $\alpha^2\beta,\gamma^3$, we get $\beta<\gamma<\alpha$. 

The AAD $\thick^{\epsilon}\delta^{\beta}\thin^{\alpha}\gamma^{\epsilon}\thin^{\beta}\delta^{\epsilon}\thick$ of $\gamma^3$ implies a vertex $\alpha\epsilon\cdots$. By $\alpha^2\beta$, and $\alpha<2\epsilon$, and $\delta>\epsilon$, we get $R(\alpha\delta\epsilon)<R(\alpha\epsilon^2)<\beta<\alpha,\gamma,2\delta,2\epsilon$. This implies $\beta\epsilon\cdots$ is not a vertex, a contradiction. 

\subsubsection*{Case. $\text{AVC}_3= \{\gamma\delta^2,\alpha\beta^2,\alpha^2\gamma\}$}
 
The angle sums of the degree $3$ vertices, and the angle sum of $\delta\epsilon^3$ or $\epsilon^4$, and the angle sum for pentagon imply
\begin{align*}
\delta\epsilon^3 &\colon
	\alpha=\delta=(\tfrac{4}{5}-\tfrac{24}{5f})\pi,\,
	\beta=(\tfrac{3}{5}+\tfrac{12}{5f})\pi,\,
	\gamma=(\tfrac{2}{5}+\tfrac{48}{5f})\pi,\,
	\epsilon=(\tfrac{2}{5}+\tfrac{8}{5f})\pi. \\
\epsilon^4 &\colon 
	\alpha=\delta=(1-\tfrac{8}{f})\pi,\,
	\beta=(\tfrac{1}{2}+\tfrac{4}{f})\pi,\,
	\gamma=\tfrac{16}{f}\pi,\,
	\epsilon=\tfrac{1}{2}\pi.
\end{align*}
By $f\ge 24$, we get $\delta>\epsilon$, and $\alpha<2\epsilon$, and $\alpha+\epsilon>\pi$. By Lemma \ref{geometry1}, this implies $\beta<\gamma$. Then by $\alpha\beta^2,\alpha^2\gamma$, we get $\alpha<\beta<\gamma$. 

The AAD $\thick^{\epsilon}\delta^{\beta}\thin^{\alpha}\gamma^{\epsilon}\thin^{\beta}\delta^{\epsilon}\thick$ of $\gamma\delta^2$ implies a vertex $\beta\epsilon\cdots$. By $\alpha+\epsilon>\pi$, and $\alpha<\beta$, and $\delta>\epsilon$, we get $R(\beta\delta\epsilon)<R(\beta\epsilon^2)<\alpha<\beta,\gamma,2\delta,2\epsilon$. This implies $\beta\epsilon\cdots$ is not a vertex, a contradiction.
\end{proof}

\section{$\alpha,\beta,\gamma$ Are Not Distinct}
\label{ndistinct}

After Sections \ref{symmetric_tiling} and \ref{distinct}, it remains to classify edge-to-edge tilings by congruent non-symmetric pentagons in Figure \ref{pentagon}, such that the value of $\alpha,\beta,\gamma$ are not distinct. By Lemma \ref{geometry11}, the non-symmetry property means $\beta\ne\gamma$ and $\delta\ne\epsilon$. Then up to the exchange symmetry $(\beta,\delta)\leftrightarrow(\gamma,\epsilon)$, we may assume $\alpha=\beta\ne\gamma$. All the results in this section will assume $\alpha=\beta\ne\gamma$. 

We use $\alpha$ (and sometimes use $\dot{\alpha}$) to represent either $\alpha$ or $\beta$. This is usually the case when the discussion involves only angle values. We will use $\beta$ only in case $\alpha,\beta$ can be distinguished.

By $\alpha=\beta$, and $a\ne b$, and Lemma \ref{geometry11}, we get $\gamma\ne\delta$. 

By Lemma \ref{geometry2}, we know the following are all the possible pairs of degree $3$ $b$-vertices:
\begin{enumerate}
\item $\alpha\delta\epsilon,\gamma\epsilon^2$.
\item $\gamma\delta\epsilon,\alpha\delta^2$.
\item $\alpha\delta^2,\gamma\epsilon^2$.
\end{enumerate}

By Lemma \ref{geometry4}, the following sets of inequalities are impossible for strictly convex pentagon:
\begin{enumerate}
\item $\alpha\ge\gamma$, $\alpha\ge\epsilon$, $\delta\le\gamma$.
\item $\alpha\le\gamma$, $\alpha\le\epsilon$, $\delta\ge \gamma$.
\end{enumerate}

The following are all the degree $3$ vertices:
\begin{itemize}
\item degree $3$ $b$-vertex: $\alpha\delta^2,\alpha\epsilon^2,\gamma\delta^2,\gamma\epsilon^2,\alpha\delta\epsilon,\gamma\delta\epsilon$.
\item degree $3$ $\hat{b}$-vertex: $\alpha^3,\alpha^2\gamma,\alpha\gamma^2,\gamma^3$.
\end{itemize}
The way multiple degree $3$ $b$-vertices appear is determined by Lemma \ref{geometry2}. Moreover, by $\alpha\ne\gamma$, the four degree $3$ $\hat{b}$-vertices are mutually exclusive.

By Lemma \ref{ndegree3}, a tiling has a degree $3$ $b$-vertex (and at most two by Lemma \ref{geometry2}). By applying Lemma \ref{degree3} to $\delta,\epsilon$, which appear together twice in the pentagon, a tiling has a degree $3$ $\hat{b}$-vertex. Since the four degree $3$ $\hat{b}$-vertices are mutually exclusive, a tiling has exactly one degree $3$ $\hat{b}$-vertex. We conclude the collection of all degree $3$ vertices consists of one or two degree $3$ $b$-vertices (two degree $3$ $b$-vertices are given by Lemma \ref{geometry2}), and one degree $3$ $\hat{b}$-vertex.

\begin{proposition}\label{4a2e}
There is no tiling, such that $\alpha=\beta\ne\gamma$, and $\alpha\delta^2$ is a vertex.
\end{proposition}

\begin{proof}
By Lemma \ref{geometry2}, $\alpha\delta^2,\gamma\epsilon^2,\gamma\delta\epsilon$ are all the degree $3$ $b$-vertices. Moreover, if $\gamma\epsilon^2,\gamma\delta\epsilon$ are not vertices, then by applying Lemma \ref{ndegree3} to $\epsilon$, we know one of $\delta\epsilon^3,\epsilon^4$ is a vertex. Therefore one of $\gamma\epsilon^2,\gamma\delta\epsilon,\delta\epsilon^3,\epsilon^4$ is a vertex. By $\delta\ne\epsilon$, we may assume the four vertices are mutually exclusive.

By $\alpha\delta^2$ and $\gamma\ne\delta$ (Lemma \ref{geometry11}), we know $\alpha\gamma^2$ is not a vertex. Therefore one of $\alpha^3,\alpha^2\gamma,\gamma^3$ is a vertex.

\subsubsection*{Case. $\alpha^3$ is a vertex}

By $\alpha^3$, we know a degree $3$ vertex $\gamma\cdots$ is a $b$-vertex. By Lemma \ref{ndegree3}, either $\gamma$ or $\epsilon$ appears at a degree $3$ vertex. Then by all the degree $3$ $b$-vertices $\alpha\delta^2,\gamma\epsilon^2,\gamma\delta\epsilon$, we know one of $\gamma\epsilon^2,\gamma\delta\epsilon$ is a vertex. 

The angle sums of $\alpha\delta^2,\alpha^3,\gamma\delta\epsilon$ and the angle sum for pentagon imply $f=12$, a contradiction. 

The angle sums of $\dot{\alpha}\delta^2,\dot{\alpha}^3,\gamma\epsilon^2$ ($\dot{\alpha}$ is $\alpha$ or $\beta$) and the angle sum for pentagon imply 
\[
	\dot{\alpha}=\delta=\tfrac{2}{3}\pi,\,
	\gamma=\tfrac{8}{f}\pi,\,
	\epsilon=(1-\tfrac{4}{f})\pi. 
\]
We have $\gamma<\dot{\alpha}=\delta=\frac{2}{3}\pi<\epsilon$. This implies $\dot{\alpha}\delta^2,\gamma\epsilon^2,\gamma^k\delta^2(k\ge 2),\gamma^k\delta\epsilon(k\ge 2)$ are all the $b$-vertices. In particular, $\dot{\alpha}\epsilon\cdots,\delta\thin\epsilon\cdots,\epsilon\thin\epsilon\cdots$ are not vertices. This implies the AADs of $\thin\gamma\thin\gamma\thin,\thin\beta\thin\gamma\thin,\thin\gamma\thin\delta\thick$ are $\thin^{\epsilon}\gamma^{\alpha}\thin^{\alpha}\gamma^{\epsilon}\thin,\thin\beta\thin^{\alpha}\gamma^{\epsilon}\thin,\thin^{\epsilon}\gamma^{\alpha}\thin^{\beta}\delta^{\epsilon}\thick$. This further implies no consecutive $\gamma\gamma\gamma,\beta\gamma\gamma,\gamma\gamma\delta$. Therefore $\gamma^k\delta^2,\gamma^k\delta\epsilon$ are not vertices, and $\dot{\alpha}\delta^2,\gamma\epsilon^2$ are all the $b$-vertices. In particular, we get $\delta\cdots=\dot{\alpha}\delta^2$ and $\epsilon\cdots=\gamma\epsilon^2$.

By $\dot{\alpha}^3$ and no consecutive $\gamma\gamma\gamma$, the $\hat{b}$-vertices are $\dot{\alpha}^3,\dot{\alpha}^2\gamma^k(k\ge 2)$. By Lemma \ref{geometry6}, we get $\dot{\alpha}+2\gamma>\pi$. Therefore $k=2,3$ in $\dot{\alpha}^2\gamma^k$. Then by no consecutive $\gamma\gamma\gamma,\beta\gamma\gamma$, we know $\alpha^2\gamma^k$ is $\thin\dot{\alpha}\thin\gamma\thin\dot{\alpha}\thin\gamma\thin$, $\thin\alpha\thin\gamma\thin\gamma\thin\alpha\thin$, $\thin\alpha\thin\gamma\thin\gamma\thin\alpha\thin\gamma\thin$. In particular, we get $\beta\thin\beta\cdots=\dot{\alpha}^3$.

We have $\thin\dot{\alpha}\thin\gamma\thin\dot{\alpha}\thin\gamma\thin=\thin\dot{\alpha}\thin\gamma\thin\dot{\alpha}\thin^{\epsilon}\gamma^{\alpha}\thin$. By the AAD $\thin\beta\thin^{\alpha}\gamma^{\epsilon}\thin$ of $\thin\beta\thin\gamma\thin$, we know $\dot{\alpha}=\alpha$ in $\thin\dot{\alpha}\thin^{\epsilon}\gamma^{\alpha}\thin$. Then by no $\beta\epsilon\cdots=\dot{\alpha}\epsilon\cdots$, we get $\thin\dot{\alpha}\thin^{\epsilon}\gamma^{\alpha}\thin=\thin^{\beta}\alpha^{\gamma}\thin^{\epsilon}\gamma^{\alpha}\thin$. Then by no $\beta\epsilon\cdots$, the AAD of $\thin\gamma\thin\dot{\alpha}\thin^{\epsilon}\gamma^{\alpha}\thin$ is $\thin^{\epsilon}\gamma^{\alpha}\thin^{\beta}\alpha^{\gamma}\thin^{\epsilon}\gamma^{\alpha}\thin$. Then by the AAD $\thin^{\beta}\alpha^{\gamma}\thin^{\epsilon}\gamma^{\alpha}\thin$ of $\thin\dot{\alpha}\thin^{\epsilon}\gamma^{\alpha}\thin$, the AAD of $\thin\dot{\alpha}\thin\gamma\thin\dot{\alpha}\thin^{\epsilon}\gamma^{\alpha}\thin$ is $\thin^{\beta}\alpha^{\gamma}\thin^{\epsilon}\gamma^{\alpha}\thin^{\beta}\alpha^{\gamma}\thin^{\epsilon}\gamma^{\alpha}\thin$.

By the AAD $\thin^{\epsilon}\gamma^{\alpha}\thin^{\alpha}\gamma^{\epsilon}\thin$ of $\thin\gamma\thin\gamma\thin$, and no $\beta\epsilon\cdots$, we know the AAD of $\thin\alpha\thin\gamma\thin\gamma\thin\alpha\thin$ is $\thin^{\beta}\alpha^{\gamma}\thin^{\epsilon}\gamma^{\alpha}\thin^{\alpha}\gamma^{\epsilon}\thin^{\gamma}\alpha^{\beta}\thin$. Then the AAD of $\thin\alpha\thin\gamma\thin\gamma\thin\alpha\thin\gamma\thin$ is $\thin^{\beta}\alpha^{\gamma}\thin^{\epsilon}\gamma^{\alpha}\thin^{\alpha}\gamma^{\epsilon}\thin^{\gamma}\alpha^{\beta}\thin\gamma\thin$, which implies $\beta\epsilon\cdots$, a contradiction. Moreover, the vertex $\thin^{\beta}\alpha^{\gamma}\thin^{\epsilon}\gamma^{\alpha}\thin^{\alpha}\gamma^{\epsilon}\thin^{\gamma}\alpha^{\beta}\thin$ implies a vertex $\thin^{\delta}\beta^{\alpha}\thin^{\alpha}\beta^{\delta}\thin\cdots=\dot{\alpha}^3=\thin^{\alpha}\beta^{\delta}\thin\dot{\alpha}\thin^{\delta}\beta^{\alpha}\thin$. By $\delta\cdots=\dot{\alpha}\delta^2$, the vertex is $\thin^{\alpha}\beta^{\delta}\thin^{\dot{\alpha}}\dot{\alpha}^{\dot{\alpha}}\thin^{\delta}\beta^{\alpha}\thin$. Then $\thin^{\dot{\alpha}}\dot{\alpha}^{\dot{\alpha}}\thin$ means three $\dot{\alpha}$ in a tile, a contradiction. 

We conclude $\dot{\alpha}\delta^2,\dot{\alpha}^3,\gamma\epsilon^2,\thin^{\beta}\alpha^{\gamma}\thin^{\epsilon}\gamma^{\alpha}\thin^{\beta}\alpha^{\gamma}\thin^{\epsilon}\gamma^{\alpha}\thin$ are all the vertices. This implies $\dot{\alpha}\gamma\cdots=\thin^{\beta}\alpha^{\gamma}\thin^{\epsilon}\gamma^{\alpha}\thin^{\beta}\alpha^{\gamma}\thin^{\epsilon}\gamma^{\alpha}\thin$, and $\beta\cdots=\dot{\alpha}^3,\beta\delta^2$.

The vertex $\thin^{\beta}\alpha^{\gamma}\thin^{\epsilon}\gamma^{\alpha}\thin^{\beta}\alpha^{\gamma}\thin^{\epsilon}\gamma^{\alpha}\thin$ determines $T_1,T_2,T_3,T_4$ in  the first of Figure \ref{4a2eA}. Then $\epsilon_1\cdots=\epsilon_2\cdots=\gamma\epsilon^2$ determines $T_5,T_6$. Then $\alpha_5\gamma_6\cdots=\thin^{\beta}\alpha^{\gamma}\thin^{\epsilon}\gamma^{\alpha}\thin^{\beta}\alpha^{\gamma}\thin^{\epsilon}\gamma^{\alpha}\thin$ determines $T_7,T_8$. Then $\alpha_7\beta_5\cdots=\dot{\alpha}^3$ and $\delta_1\delta_5\cdots=\dot{\alpha}\delta^2$ give two $\dot{\alpha}$ in the same tile. This implies $\beta_1\cdots,\beta_7\cdots$ are not $\dot{\alpha}^3$. Therefore $\beta_1\cdots=\beta_7\cdots=\beta\delta^2$. Then we get two $\delta$ in the same tile, a contradiction. 

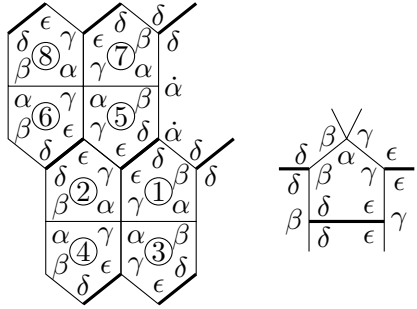
\begin{figure}[htp]
\centering
\begin{tikzpicture}[>=latex,scale=1]

\foreach \a in {1,-1}
\foreach \b/\c in {0/0, -0.5/1.8}
{
\begin{scope}[shift={(\b cm, \c cm)}, scale=\a]

\draw
	(0,-0.7) -- (0,0.7) -- (0.5,1.1) -- (1,0.7) -- (1,-0.7) -- (0.5,-1.1) -- (0,-0.7)
	(0,0) -- (1,0);
	
\draw[line width=1.2]
	(0.5,1.1) -- (0,0.7)
	(0.5,-1.1) -- (1,-0.7);
	
\node at (0.2,0.6) {\small $\epsilon$};
\node at (0.2,0.2) {\small $\gamma$};
\node at (0.5,0.85) {\small $\delta$};	
\node at (0.8,0.2) {\small $\alpha$}; 
\node at (0.8,0.6) {\small $\beta$};

\node at (0.2,-0.6) {\small $\gamma$};
\node at (0.5,-0.85) {\small $\epsilon$};
\node at (0.2,-0.2) {\small $\alpha$};	
\node at (0.8,-0.2) {\small $\beta$}; 
\node at (0.8,-0.6) {\small $\delta$};

\end{scope}
}

\draw[line width=1.2]
	(1,0.7) -- (1.5,1.1)
	(0.5,2.5) -- (1,2.9);

\node at (0.7,1.8) {\small $\dot{\alpha}$};
\node at (0.7,1.2) {\small $\dot{\alpha}$};
\node at (1.2,0.6) {\small $\delta$};
\node at (1,0.95) {\small $\delta$};
\node at (0.7,2.4) {\small $\delta$};
\node at (0.5,2.75) {\small $\delta$};

\node[inner sep=0.5,draw,shape=circle] at (0.5,0.4) {\small $1$};
\node[inner sep=0.5,draw,shape=circle] at (0.5,-0.4) {\small $3$};
\node[inner sep=0.5,draw,shape=circle] at (-0.5,0.4) {\small $2$};
\node[inner sep=0.5,draw,shape=circle] at (-0.5,-0.4) {\small $4$};

\node[inner sep=0.5,draw,shape=circle] at (0,2.2) {\small $7$};
\node[inner sep=0.5,draw,shape=circle] at (0,1.4) {\small $5$};
\node[inner sep=0.5,draw,shape=circle] at (-1,2.2) {\small $8$};
\node[inner sep=0.5,draw,shape=circle] at (-1,1.4) {\small $6$};

\begin{scope}[xshift=3cm]

\draw
	(-0.5,-0.4) -- (-0.5,0.7) -- (0,1.1) -- (0.5,0.7) -- (0.5,-0.4)
	(-0.5,0.7) -- ++(-0.4,0)
	(0,1.1) -- ++(0.2,0.4)
	(0,1.1) -- ++(-0.2,0.4);
		
\draw[line width=1.2]
	(-0.5,0) -- (0.5,0)
	(-0.5,0.7) -- ++(-0.4,0)
	(0.5,0.7) -- ++(0.4,0);
	
\node at (0,0.85) {\small $\alpha$}; 
\node at (-0.3,0.6) {\small $\beta$};
\node at (0.3,0.6) {\small $\gamma$};
\node at (-0.3,0.2) {\small $\delta$};
\node at (0.3,0.2) {\small $\epsilon$};	

\node at (-0.7,0) {\small $\beta$};
\node at (0.7,0) {\small $\gamma$};
\node at (-0.3,-0.2) {\small $\delta$};
\node at (0.3,-0.2) {\small $\epsilon$};
\node at (-0.7,0.5) {\small $\delta$};
\node at (-0.6,0.9) {\small $\delta$};
\node at (0.7,0.5) {\small $\epsilon$};
\node at (0.6,0.9) {\small $\epsilon$};

\node at (-0.25,1.1) {\small $\beta$};
\node at (0.25,1.1) {\small $\gamma$};

\end{scope}

\end{tikzpicture}
\caption{Proposition \ref{4a2e}: $\alpha^3$ or $\gamma^3$ is a vertex.}
\label{4a2eA}
\end{figure}

\subsubsection*{Case. $\gamma^3$ is a vertex}

By $\dot{\alpha}\delta^2,\gamma^3$, we know a degree $3$ vertex $\dot{\alpha}\cdots=\dot{\alpha}\delta^2$. Moreover, we know a degree $3$ vertex $\gamma\cdots=\gamma^3,\gamma\epsilon^2,\gamma\delta\epsilon$.

In the special tile in the second of Figure \ref{4a2eA}, one of the vertices $\alpha\cdots,\beta\cdots$ has degree $3$. If the vertex $\alpha\cdots$ has degree $3$, then $\alpha\cdots=\alpha\delta^2$. This implies $\gamma\cdots=\beta\gamma\cdots=\dot{\alpha}\gamma\cdots$ and $\beta\cdots=\beta^2\cdots=\dot{\alpha}^2\cdots$ have high degree, a contradiction. 

Therefore the vertex $\beta\cdots$ has degree $3$. Then $\beta\cdots=\beta\delta^2$. This implies $\alpha\cdots=\alpha\beta\cdots=\dot{\alpha}^2\cdots$ has high degree. Therefore the vertices $\delta\cdots,\epsilon\cdots$ have degree $3$. However, $\beta\cdots=\beta\delta^2$ implies the degree $3$ vertex $\delta\cdots=\beta\delta\cdots=\dot{\alpha}\delta\cdots=\dot{\alpha}\delta^2$. This implies the degree $3$ vertex $\epsilon\cdots=\epsilon^2\cdots=\gamma\epsilon^2$. Then $\gamma\delta\epsilon$ is not a vertex, and a degree $3$ vertex $\gamma\cdots=\gamma^3,\gamma\epsilon^2$. Then by the vertex $\epsilon\cdots=\gamma\epsilon^2$, we know the degree $3$ vertex $\gamma\cdots=\gamma\epsilon^2$. Then the high degree vertex $\alpha\cdots=\alpha\beta\gamma\cdots=\dot{\alpha}^2\gamma\cdots$. 

The angle sums of $\dot{\alpha}\delta^2,\gamma^3,\gamma\epsilon^2$ and the angle sum for pentagon imply
\[
\dot{\alpha}=(\tfrac{4}{9}+\tfrac{8}{3f})\pi,\,
\gamma=\epsilon=\tfrac{2}{3}\pi,\,
\delta=(\tfrac{7}{9}-\tfrac{4}{3f})\pi.
\]
We have $\dot{\alpha}<\gamma=\epsilon<\delta$. Then by $R(\dot{\alpha}^2\gamma)=(\frac{4}{9}-\frac{16}{3f})\pi<\dot{\alpha}$, we know $\dot{\alpha}^2\gamma\cdots$ is not a vertex. 

\subsubsection*{Case. $\alpha^2\gamma$ is a vertex}

The angle sums of $\alpha\delta^2,\alpha^2\gamma$, and the angle sum of one of $\gamma\epsilon^2,\gamma\delta\epsilon,\delta\epsilon^3,\epsilon^4$, and the angle sum for pentagon imply
\begin{align*}
\gamma\epsilon^2 &\colon
	\alpha=\epsilon=\tfrac{8}{f}\pi,\,
	\gamma=(2-\tfrac{16}{f})\pi,\,
	\delta=(1-\tfrac{4}{f})\pi. \\
\gamma\delta\epsilon &\colon
	\alpha=(\tfrac{1}{2}+\tfrac{2}{f})\pi,\,
	\gamma=(1-\tfrac{4}{f})\pi,\,
	\delta=(\tfrac{3}{4}-\tfrac{1}{f})\pi,\,
	\epsilon=(\tfrac{1}{4}+\tfrac{5}{f})\pi.  \\
\delta\epsilon^3 &\colon
	\alpha=(1-\tfrac{12}{f})\pi,\,
	\gamma=\tfrac{24}{f}\pi,\,
	\delta=(\tfrac{1}{2}+\tfrac{6}{f})\pi,\,
	\epsilon=(\tfrac{1}{2}-\tfrac{2}{f})\pi. \\
\epsilon^4 &\colon
	\alpha=(1-\tfrac{8}{f})\pi,\,
	\gamma=\tfrac{16}{f}\pi,\,
	\delta=(\tfrac{1}{2}+\tfrac{4}{f})\pi,\,
	\epsilon=\tfrac{1}{2}\pi. 	 
\end{align*}
We have $\delta>\epsilon$. By Lemma \ref{geometry1}, this implies $\alpha<\gamma$. By $\alpha^2\gamma$, this further implies $\gamma>\frac{2}{3}\pi$.

The AAD of $\delta\epsilon^3,\epsilon^4$  implies $\gamma^2\cdots$ is a vertex. By $\alpha^2\gamma$ and $\alpha<\gamma$, we know $R(\gamma^2)$ has no $\alpha,\gamma$. By $\gamma+\epsilon>\frac{2}{3}\pi+(\tfrac{1}{2}-\tfrac{2}{f})\pi>\pi$, and $\delta>\epsilon$, we know $R(\gamma^2)$ has no $\delta,\epsilon$. Therefore $\gamma^2\cdots$ is not a vertex, a contradiction.

For $\gamma\epsilon^2,\gamma\delta\epsilon$ is a vertex. We have the angle values in terms of $f$. Substituting the angle values into \eqref{coolsaet_eq1} and solve for $f$, we find no solution satisfying $f\ge 16$. 
\end{proof}

In Proposition \ref{4a2e}, for the case $\dot{\alpha}\delta^2,\dot{\alpha}^2\gamma,\gamma\epsilon^2$ are all the degree $3$ vertices, we may find that $\dot{\alpha}^4,\dot{\alpha}^5,\dot{\alpha}^2\delta\epsilon$ are the possible high degree vertices. The angle sum of the high degree vertex further implies
\begin{align*}
\dot{\alpha}^4 &\colon
	\dot{\alpha}=\epsilon=\tfrac{1}{2}\pi,\,
	\gamma=\pi,\,
	\delta=\tfrac{3}{4}\pi,\,
	f=16. \\
\dot{\alpha}^5,\dot{\alpha}^2\delta\epsilon &\colon
	\dot{\alpha}=\epsilon=\tfrac{2}{5}\pi,\,
	\gamma=\tfrac{6}{5}\pi,\,
	\delta=\tfrac{4}{5}\pi,\,
	f=20. 
\end{align*}
Then we may construct earth map tilings by using four or five timezones from Figure \ref{4a2eC}. In fact, we may mix the timezones in any order. 

\begin{figure}[htp]
\centering
\begin{tikzpicture}[>=latex,scale=1]

\foreach \a in {0,1}
\draw[xshift=3.5*\a cm]
	(0.5,1.5) -- (0.5,1.1) -- (1,0.7) -- (1,0) -- (1.5,-0.4) -- (1.5,-0.8)
	(-1.5,1.5) -- (-1.5,1.1) -- (-1,0.7) -- (-1,0) -- (-0.5,-0.4) -- (-0.5,-0.8) 
	(0.5,1.1) -- (0,0.7) -- (-1,0.7)
	(1,0) -- (0,0) -- (-0.5,-0.4)
	(0,0) -- (0,0.7);


\draw[line width=1.2]
	(0,0.7) -- (-1,0.7)
	(1,0) -- (0,0);

\node at (-0.5,1.4) {\small $\alpha$};
\node at (0.3,1.2) {\small $\gamma$}; 
\node at (-0.1,0.9) {\small $\epsilon$};
\node at (-0.9,0.9) {\small $\delta$};
\node at (-1.3,1.2) {\small $\beta$};

\node at (0.2,0.6) {\small $\gamma$};
\node at (0.8,0.6) {\small $\beta$};
\node at (0.2,0.2) {\small $\epsilon$}; 
\node at (0.8,0.2) {\small $\delta$};	
\node at (0.5,0.85) {\small $\alpha$};	
	
\node at (-0.2,0.5) {\small $\epsilon$}; 
\node at (-0.8,0.5) {\small $\delta$};
\node at (-0.8,0.1) {\small $\beta$};	
\node at (-0.5,-0.15) {\small $\alpha$};	
\node at (-0.2,0.1) {\small $\gamma$};

\node at (0.5,-0.7) {\small $\alpha$};
\node at (1.3,-0.5) {\small $\beta$}; 
\node at (0.9,-0.2) {\small $\delta$};
\node at (0.1,-0.2) {\small $\epsilon$};
\node at (-0.3,-0.5) {\small $\gamma$};


\begin{scope}[xshift=3.5 cm]

\draw[line width=1.2]
	(0,0.7) -- (0.5,1.1)
	(-0.5,-0.4) -- (0,0);

\node at (-0.5,1.4) {\small $\beta$};
\node at (0.3,1.2) {\small $\delta$}; 
\node at (-0.1,0.9) {\small $\epsilon$};
\node at (-0.9,0.9) {\small $\gamma$};
\node at (-1.3,1.2) {\small $\alpha$};

\node at (0.2,0.6) {\small $\epsilon$};
\node at (0.8,0.6) {\small $\beta$};
\node at (0.2,0.2) {\small $\gamma$}; 
\node at (0.8,0.2) {\small $\alpha$};	
\node at (0.5,0.85) {\small $\delta$};	
	
\node at (-0.2,0.5) {\small $\gamma$}; 
\node at (-0.8,0.5) {\small $\alpha$};
\node at (-0.8,0.1) {\small $\beta$};	
\node at (-0.5,-0.15) {\small $\delta$};	
\node at (-0.2,0.1) {\small $\epsilon$};

\node at (0.5,-0.7) {\small $\beta$};
\node at (1.3,-0.5) {\small $\alpha$}; 
\node at (0.9,-0.2) {\small $\gamma$};
\node at (0.1,-0.2) {\small $\epsilon$};
\node at (-0.3,-0.5) {\small $\delta$}; 

\end{scope}

\end{tikzpicture}
\caption{Proposition \ref{4a2e}: Geometrically non-existent tilings.}
\label{4a2eC}
\end{figure}

Similarly, for the case $\dot{\alpha}\delta^2,\dot{\alpha}^2\gamma,\gamma\delta\epsilon$ are all the degree $3$ vertices, we may find that $\dot{\alpha}^2\epsilon^2$ is the only high degree vertex. The angle sum of $\dot{\alpha}^2\epsilon^2$ further implies
\[
\dot{\alpha}=\tfrac{4}{7}\pi,\,
\gamma=\tfrac{6}{7}\pi,\,
\delta=\tfrac{5}{7}\pi,\,
\epsilon=\tfrac{3}{7}\pi,\,
f=28. 
\]
Then we may construct two tilings. Each is obtained by glueing two identical half tilings (with $14$ tiles) in Figure \ref{4a2eB}. The two copies are glue together, such that the shaded edge on the upper right are matched in twisted way. In other words, the two vertices at the ends of the shaded edge is $\gamma\delta\epsilon$ in the first picture, and is $\beta\delta^2$ in the second picture. 

\begin{figure}[htp]
\centering
\begin{tikzpicture}[>=latex,scale=1]

\foreach \b in {0,1}
\draw[xshift=6*\b cm, gray!50, line width=3]
	(1.3,1.7) -- (2.1,1.7);

\foreach \a in {1,-1}
\foreach \b in {0,1}
{
\begin{scope}[xshift=6*\b cm, scale=\a]

\draw
	(-2.1,-0.7) -- (-2.1,1.7) -- (-1.3,1.7) -- (-0.65,2.1) -- (0,1.7) -- (0,1.1) -- (0.5,0.7) -- (2.1,0.7) -- (2.5,1.25) -- (2.1,1.7) -- (0,1.7)
	(0.5,0.7) -- (0.5,-0.7) -- (0,-1.1) 
	(-2.1,0.7) -- (-0.5,0.7) -- (-1.3,1.7)
	(1.3,1.7) -- (1.3,-0.7)
	(-0.5,0) -- (0.5,0);

\end{scope}
}	


\foreach \a in {1,-1}
{
\begin{scope}[scale=\a]

\draw[line width=1.2]
	(0,-1.1) -- (0.5,-0.7)
	(1.3,1.7) -- (2.1,1.7)
	(-1.3,1.7) -- (-2.1,1.7)
	(0.5,0.7) -- (1.3,0.7)
	(2.1,0) -- (2.1,-0.7);

\node at (1.1,0.9) {\small $\epsilon$};
\node at (0.6,0.9) {\small $\delta$};
\node at (1.1,1.5) {\small $\gamma$};
\node at (0.2,1.2) {\small $\beta$};
\node at (0.2,1.5) {\small $\alpha$};

\node at (1,-1.65) {\small $\gamma$};
\node at (0.5,-1) {\small $\epsilon$};
\node at (0.65,-1.9) {\small $\alpha$};
\node at (0.2,-1.2) {\small $\delta$};
\node at (0.2,-1.6) {\small $\beta$};

\node at (0.3,0.6) {\small $\beta$};
\node at (0,0.85) {\small $\delta$};
\node at (0.3,0.2) {\small $\alpha$};	
\node at (-0.3,0.2) {\small $\gamma$}; 
\node at (-0.3,0.6) {\small $\epsilon$};

\node at (1.1,0.5) {\small $\epsilon$};
\node at (1.1,-0.5) {\small $\gamma$};
\node at (0.7,0.5) {\small $\delta$};
\node at (0.7,-0.5) {\small $\alpha$};
\node at (0.7,0) {\small $\beta$};

\node at (1.5,0.5) {\small $\alpha$};
\node at (1.5,-0.5) {\small $\beta$};
\node at (1.9,0.5) {\small $\gamma$};
\node at (1.9,-0.5) {\small $\delta$};
\node at (1.9,0) {\small $\epsilon$};

\node at (2,0.9) {\small $\alpha$};
\node at (2.25,1.2) {\small $\gamma$};
\node at (1.5,0.9) {\small $\beta$};
\node at (2,1.5) {\small $\epsilon$};
\node at (1.5,1.5) {\small $\delta$};

\node at (1.3,-0.9) {\small $\alpha$};
\node at (1.9,-0.9) {\small $\gamma$};
\node at (0.9,-0.9) {\small $\beta$};
\node at (1.9,-1.5) {\small $\epsilon$};
\node at (1.4,-1.5) {\small $\delta$};

\end{scope}
}

\draw
	(2.1,1.7) -- (2.5,2.1);
	
\node at (1.3,1.9) {\small $\epsilon$};
\node at (2,1.9) {\small $\delta$};
\node at (2.4,1.75) {\small $\gamma$};


\begin{scope}[xshift=6 cm]

\foreach \a in {1,-1}
{
\begin{scope}[scale=\a]

\draw[line width=1.2]
	(0.5,0) -- (0.5,-0.7)
	(2.1,0) -- (2.1,0.7)
	(-1.3,1.7) -- (-0.65,2.1)
	(1.3,1.7) -- (1.3,0.7)
	(-2.1,1.7) -- (-2.1,0.7);
	
\node at (1.1,0.9) {\small $\epsilon$};
\node at (0.6,0.9) {\small $\gamma$};
\node at (1.1,1.5) {\small $\delta$};
\node at (0.2,1.2) {\small $\alpha$};
\node at (0.2,1.5) {\small $\beta$};

\node at (1,-1.65) {\small $\delta$};
\node at (0.5,-1) {\small $\beta$};
\node at (0.65,-1.85) {\small $\epsilon$};
\node at (0.2,-1.2) {\small $\alpha$};
\node at (0.2,-1.55) {\small $\gamma$};

\node at (0.3,0.6) {\small $\alpha$};
\node at (0,0.85) {\small $\gamma$};
\node at (0.3,0.2) {\small $\beta$};	
\node at (-0.3,0.2) {\small $\delta$}; 
\node at (-0.3,0.6) {\small $\epsilon$};

\node at (1.1,0.5) {\small $\alpha$};
\node at (1.1,-0.5) {\small $\gamma$};
\node at (0.7,0.5) {\small $\beta$};
\node at (0.7,-0.5) {\small $\epsilon$};
\node at (0.7,0) {\small $\delta$};

\node at (1.5,0.5) {\small $\beta$};
\node at (1.5,-0.5) {\small $\alpha$};
\node at (1.9,0.5) {\small $\delta$};
\node at (1.9,-0.5) {\small $\gamma$};
\node at (1.9,0) {\small $\epsilon$};

\node at (2,0.9) {\small $\gamma$};
\node at (2.25,1.2) {\small $\alpha$};
\node at (1.5,0.9) {\small $\epsilon$};
\node at (2,1.5) {\small $\beta$};
\node at (1.5,1.5) {\small $\delta$};

\node at (1.3,-0.9) {\small $\beta$};
\node at (1.9,-0.9) {\small $\delta$};
\node at (0.9,-0.9) {\small $\alpha$};
\node at (1.9,-1.5) {\small $\epsilon$};
\node at (1.4,-1.5) {\small $\gamma$};

\end{scope}
}

\draw[line width=1.2]
	(2.1,1.7) -- (2.5,2.1);
	
\node at (1.3,1.9) {\small $\beta$};
\node at (2,1.9) {\small $\delta$};
\node at (2.35,1.7) {\small $\delta$};

\end{scope}

\end{tikzpicture}
\caption{Proposition \ref{4a2e}: Geometrically non-existent tilings.}
\label{4a2eB}
\end{figure}

The tilings given by Figures \ref{4a2eC} and \ref{4a2eB} are only combinatorially possible. Since their angle values fail Lemma \ref{geometry7}, they are geometrically impossible.

\begin{proposition}\label{4a2d}
There is no tiling, such that $\alpha=\beta\ne\gamma$, and $\alpha\epsilon^2$ is a vertex.
\end{proposition}

\begin{proof}
By $\alpha\epsilon^2$ and Lemma \ref{geometry2}, we know $\alpha\epsilon^2$ is the only degree $3$ $b$-vertex, and $\delta\cdots$ has high degree. By Lemma \ref{ndegree3}, this implies one of $\delta^4,\delta^3\epsilon$ is a vertex. Moreover, in a special tile, we know $\epsilon\cdots=\alpha\epsilon^2$, and $\alpha\cdots,\beta\cdots,\gamma\cdots$ have degree $3$, and $\delta\cdots$ has degree $4$ or $5$. Then by Lemma \ref{special_tile}, we get $f\ge 24$. Moreover, the degree $3$ vertex $\gamma\cdots$ is a $\hat{b}$-vertex. Therefore one of $\alpha^2\gamma,\alpha\gamma^2,\gamma^3$ is a vertex, and all degree $3$ vertices are $\alpha\epsilon^2$ and one of $\alpha^2\gamma,\alpha\gamma^2,\gamma^3$.

If $\alpha\epsilon^2,\gamma^3$ are all the degree $3$ vertices, then in a special tile, one of $\alpha\cdots,\beta\cdots$ has high degree, a contradiction. 

The angle sum of $\alpha\epsilon^2$, and the angle sum of one of $\alpha^2\gamma,\alpha\gamma^2$, and the angle sum of one of $\delta^4,\delta^3\epsilon$, and the angle sum for pentagon imply
\begin{align*}
\alpha^2\gamma,\delta^4 &\colon
	\alpha=(1-\tfrac{8}{f})\pi,\,
	\gamma=\tfrac{16}{f}\pi,\,
	\delta=\tfrac{1}{2}\pi,\,
	\epsilon=(\tfrac{1}{2}+\tfrac{4}{f})\pi.	\\
\alpha^2\gamma,\delta^3\epsilon &\colon
	\alpha=(1-\tfrac{12}{f})\pi,\,
	\gamma=\tfrac{24}{f}\pi,\,
	\delta=(\tfrac{1}{2}-\tfrac{2}{f})\pi,\,
	\epsilon=(\tfrac{1}{2}+\tfrac{6}{f})\pi.  \\
\alpha\gamma^2,\delta^4 &\colon
	\alpha=(\tfrac{1}{2}+\tfrac{4}{f})\pi,\,
	\gamma=\epsilon=(\tfrac{3}{4}-\tfrac{2}{f})\pi,\,
	\delta=\tfrac{1}{2}\pi.	\\
\alpha\gamma^2,\delta^3\epsilon &\colon
	\alpha=(\tfrac{4}{7}+\tfrac{24}{7f})\pi,\,
	\gamma=\epsilon=(\tfrac{5}{7}-\tfrac{12}{7f})\pi,\,
	\delta=(\tfrac{3}{7}+\tfrac{4}{7f})\pi.
\end{align*}
We have $\delta<\epsilon$ and $\alpha<2\delta$. By Lemma \ref{geometry1}, this implies $\alpha>\gamma$. 

By $\alpha^2\gamma$ or $\alpha\gamma^2$, and $\gamma<\alpha<2\delta$, we get $\alpha+\gamma+2\delta>2\pi$. Then by $\alpha>\gamma$ and $\delta<\epsilon$, this implies a $b$-vertex $\alpha\cdots=\alpha\epsilon^2$. 

If $\alpha\epsilon^2,\alpha\gamma^2$ are all the degree $3$ vertices, then by $\alpha>\gamma$, we know a $\hat{b}$-vertex $\alpha\cdots=\alpha\gamma^2$. Then we get $\alpha\cdots=\alpha\gamma^2,\alpha\epsilon^2$. This implies
\[
2f=\#\alpha=\#\alpha\gamma^2+\#\alpha\epsilon^2
\le \tfrac{1}{2}\#\gamma+\tfrac{1}{2}\#\epsilon
=\tfrac{1}{2}f+\tfrac{1}{2}f=f,
\]
a contradiction. 

Therefore $\alpha\epsilon^2,\alpha^2\gamma$ are all the degree $3$ vertices. By $\pi>\alpha>\gamma$, the pentagon is strictly convex. By $2(\alpha-\epsilon)=\alpha-\gamma>0$, we get $\alpha>\epsilon$. Then by $\alpha>\gamma$ and Lemma \ref{geometry4}, we get $\gamma>\delta$. 

For $\alpha^2\gamma,\delta^4$, we have $\alpha,\gamma,\epsilon>\delta=\tfrac{1}{2}\pi$. This implies $\delta^4$ is the only high degree vertex, and $\alpha\epsilon^2,\alpha^2\gamma,\delta^4$ are all the vertices. Then we get
\[
f=\#\gamma
=\#\alpha^2\gamma
<\tfrac{1}{2}(\#\alpha\epsilon^2+2\#\alpha^2\gamma)
=\tfrac{1}{2}\#\alpha=f,
\]
a contradiction.

Therefore $\alpha^2\gamma,\delta^3\epsilon$ are all the degree $3$ vertices. By $\alpha^2\gamma$ and $2\gamma>2\delta>\alpha>\gamma$, we know a $\hat{b}$-vertex $\alpha\cdots=\alpha^2\gamma$. Combined with $b$-vertex $\alpha\cdots=\alpha\epsilon^2$, we get $\alpha\cdots=\alpha\epsilon^2,\alpha^2\gamma$. Then by the only degree $3$ $b$-vertex $\alpha\epsilon^2$, we know $\gamma\delta\cdots,\delta\epsilon\cdots$ are high degree vertices with no $\alpha$. Then by $\delta^3\epsilon$ and $\gamma>\delta$, we get $\delta\epsilon\cdots=\delta^3\epsilon$. This implies $\gamma\delta\cdots$ has no $\epsilon$. By $f\ge 24$, we get $5\delta=(\frac{5}{2}-\frac{10}{f})\pi>2\pi$. Then by $\gamma>\delta$, we get $\gamma\delta\cdots=\gamma^2\delta^2$.

The AAD $\thick^{\epsilon}\delta^{\beta}\thin^{\beta}\delta^{\epsilon}\thick^{\epsilon}\delta^{\beta}\thin^{\gamma}\epsilon^{\delta}\thick$ of $\delta^3\epsilon$ determines $T_1,T_2,T_3,T_4$ in Figure \ref{4a2dA}. Then $\delta_4\epsilon_3\cdots=\delta^3\epsilon$ determines $T_5,T_6$. Then $\beta_1\beta_3\cdots=\dot{\alpha}^2\cdots=\dot{\alpha}^2\gamma$ and $\beta_5\gamma_3\cdots=\dot{\alpha}\gamma\cdots=\dot{\alpha}^2\gamma$ ($\dot{\alpha}$ is $\alpha$ or $\beta$) give $\gamma_7,\dot{\alpha}_8$. If $T_7$ is not arranged as indicated, then $\alpha_3\cdots=\dot{\alpha}\epsilon\cdots=\dot{\alpha}\epsilon^2$. This implies $\dot{\alpha}_8,\epsilon_8$ adjacent, a contradiction. Therefore $T_7$ is arranged as indicated. Then $\alpha_3\alpha_7\cdots=\dot{\alpha}^2\cdots=\dot{\alpha}^2\gamma$ and $\dot{\alpha}_8$ determine $T_8$. Then $\alpha_5\beta_8\cdots=\dot{\alpha}^2\cdots=\dot{\alpha}^2\gamma$ and no $\dot{\alpha}\delta\cdots$ determine $T_9$. Then $\delta_8\epsilon_9\cdots=\delta^3\epsilon$ is comparable to $\delta^3\epsilon$ we start with, and induces six tiles similar to $T_1,T_2,T_3,T_4,T_5,T_6$. The process continues, and we get more and more six tile combinations. In particular, we get more copies of $T_7$. 

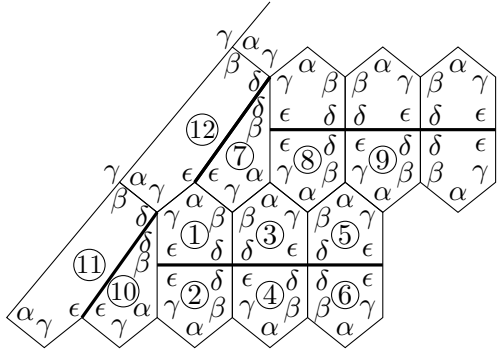
\begin{figure}[htp]
\centering
\begin{tikzpicture}[>=latex,scale=1]

\foreach \c in {0,1}
{
\begin{scope}[shift={(1.5*\c cm, 1.8*\c cm)}]

\foreach \a in {-1,0,1}
\draw[xshift=\a cm]
	(-0.5,-0.7) -- (-0.5,0.7) -- (0,1.1) -- (0.5,0.7) -- (0.5,-0.7) -- (0,-1.1) -- (-0.5,-0.7);

\foreach \b in {1,-1}
{
\begin{scope}[scale=\b]

\node at (0.3,0.6) {\small $\gamma$};
\node at (0,0.85) {\small $\alpha$};
\node at (0.3,0.2) {\small $\epsilon$};	
\node at (-0.3,0.2) {\small $\delta$}; 
\node at (-0.3,0.6) {\small $\beta$};

\node at (1.3,0.6) {\small $\gamma$};
\node at (1,0.85) {\small $\alpha$};
\node at (1.3,0.2) {\small $\epsilon$};	
\node at (0.7,0.2) {\small $\delta$}; 
\node at (0.7,0.6) {\small $\beta$};

\node at (-1.3,0.6) {\small $\gamma$};
\node at (-1,0.85) {\small $\alpha$};
\node at (-1.3,0.2) {\small $\epsilon$};	
\node at (-0.7,0.2) {\small $\delta$}; 
\node at (-0.7,0.6) {\small $\beta$};
	
\end{scope}
}

\draw
	(-1.5,-0.7) -- (-2,-1.1) -- (-2.5,-0.7) -- (-3,-1.1) -- (-3.5,-0.7) -- (-2,1.1) -- (-1.5,0.7);

\draw[line width=1.2]
	(-1.5,0) -- (1.5,0)
	(-2.5,-0.7) -- (-1.5,0.7);

\node at (-1.7,-0.6) {\small $\alpha$};
\node at (-2,-0.85) {\small $\gamma$};
\node at (-1.7,0) {\small $\beta$};
\node at (-2.25,-0.6) {\small $\epsilon$};
\node at (-1.65,0.3) {\small $\delta$};

\node at (-3.25,-0.65) {\small $\alpha$};
\node at (-3,-0.85) {\small $\gamma$};
\node at (-2,0.85) {\small $\beta$};
\node at (-2.6,-0.6) {\small $\epsilon$};
\node at (-1.7,0.65) {\small $\delta$};

\node at (-2.1,1.2) {\small $\gamma$};

\end{scope}
}

\draw
	(-0.5,2.9) -- ++(0.5,0.6);

\draw[line width=1.2]
	(-1,1.1) -- (0,2.5);

\node at (-0.25,2.95) {\small $\alpha$};
\node at (0,2.75) {\small $\gamma$};

\node[inner sep=0.5,draw,shape=circle] at (1,0.4) {\small $5$};
\node[inner sep=0.5,draw,shape=circle] at (1,-0.4) {\small $6$};
\node[inner sep=0.5,draw,shape=circle] at (0,0.4) {\small $3$};
\node[inner sep=0.5,draw,shape=circle] at (0,-0.4) {\small $4$};
\node[inner sep=0.5,draw,shape=circle] at (-1,0.4) {\small $1$};
\node[inner sep=0.5,draw,shape=circle] at (-1,-0.4) {\small $2$};
\node[inner sep=0.5,draw,shape=circle] at (-0.4,1.45) {\small $7$};
\node[inner sep=0.5,draw,shape=circle] at (0.5,1.4) {\small $8$};
\node[inner sep=0.5,draw,shape=circle] at (1.5,1.4) {\small $9$};
\node[inner sep=0,draw,shape=circle] at (-1.95,-0.35) {\footnotesize $10$};
\node[inner sep=0,draw,shape=circle] at (-2.4,0) {\footnotesize $11$};
\node[inner sep=0,draw,shape=circle] at (-0.9,1.8) {\footnotesize $12$};

\end{tikzpicture}
\caption{Proposition \ref{4a2d}: $\alpha^2\gamma,\delta^3\epsilon$ are vertices.}
\label{4a2dA}
\end{figure}

We may assume $T_1,T_2,\cdots$ is the six tile combination $T_8,T_9,\cdots$. Then $T_7$ becomes $T_{10}$, and $T_8,T_9,\cdots$ becomes the third six tile combination, which also determines $T_7$. We have more six tile combination and copies of $T_7$ further up.

The vertex $\gamma_1\delta_{10}\cdots=\gamma^2\delta^2$ determines $T_{11}$. By the same reason, we determine $T_{12}$ and more similar tiles further up. Then we get $\alpha_{12}\beta_{11}\cdots=\dot{\alpha}^2\cdots=\dot{\alpha}^2\gamma$, and similar $\dot{\alpha}^2\gamma$ further up. Then we get adjacent $\gamma$, a contradiction. 
\end{proof}

\begin{proposition}\label{4b2e}
There is no tiling, such that $\alpha=\beta\ne\gamma$, and $\gamma\delta^2$ is a vertex.
\end{proposition}

\begin{proof}
By $\gamma\delta^2$ and Lemma \ref{geometry2}, we know $\gamma\delta^2$ is the only degree $3$ $b$-vertex, and $\epsilon\cdots$ has high degree. By Lemma \ref{ndegree3}, this implies one of $\delta\epsilon^3,\epsilon^4$ is a vertex. 

We know degree $3$ vertices are $\gamma\delta^2$ and one of $\alpha^3,\alpha^2\gamma,\alpha\gamma^2,\gamma^3$. By applying Lemma \ref{degree3} to $\gamma$, we know that $\gamma\delta^2,\alpha^3$ are all the degree $3$ vertices. The angle sums of $\gamma\delta^2,\alpha^3$, and the angle sum of one of $\delta\epsilon^3,\epsilon^4$, and the angle sum for pentagon imply
\begin{align*}
\delta\epsilon^3 &\colon
	\alpha=\tfrac{2}{3}\pi,\,
	\gamma=(\tfrac{1}{2}+\tfrac{6}{f})\pi,\,
	\delta=(\tfrac{3}{4}-\tfrac{3}{f})\pi,\,
	\epsilon=(\tfrac{5}{12}+\tfrac{1}{f})\pi.  \\
\epsilon^4 &\colon
	\alpha=\tfrac{2}{3}\pi,\,
	\gamma=(\tfrac{1}{3}+\tfrac{8}{f})\pi,\,
	\delta=(\tfrac{5}{6}-\tfrac{4}{f})\pi,\,
	\epsilon=\tfrac{1}{2}\pi.	 
\end{align*}
We have $\delta>\epsilon$ and $\gamma<2\epsilon$. By Lemma \ref{geometry1}, this implies $\alpha<\gamma$. 

The AAD of $\delta\epsilon^3,\epsilon^4$  implies $\gamma^2\cdots$ is a vertex. By $\alpha^3$ and $\alpha<\gamma$, we get $R(\gamma^2)<\alpha,\gamma,2\delta,2\epsilon$. This implies $\gamma^2\cdots$ is not a vertex, a contradiction.
\end{proof}

\begin{proposition}\label{4b2d}
There is no tiling, such that $\alpha=\beta\ne\gamma$, and $\gamma\epsilon^2$ is a vertex.
\end{proposition}

\begin{proof}
By Lemma \ref{geometry2}, $\gamma\epsilon^2,\alpha\delta\epsilon,\alpha\delta^2$ are all the degree $3$ $b$-vertices. By Proposition \ref{4a2e}, we know $\alpha\delta^2$ is not a vertex. 

If $\alpha\delta\epsilon$ is not a vertex, then $\gamma\epsilon^2$ and one of $\alpha^3,\alpha^2\gamma,\alpha\gamma^2,\gamma^3$ are all the degree $3$ vertices. If $\gamma\epsilon^2,\alpha^3$ are all the degree $3$ vertices, then in any tile, one of $\alpha\cdots,\gamma\cdots$ has high degree, and $\delta\cdots$ has high degree. By Lemma \ref{special_tile}, this is a contradiction. If $\gamma\epsilon^2$ and one of $\alpha^2\gamma,\alpha\gamma^2,\gamma^3$ are all the degree $3$ vertices, then by applying Lemma \ref{degree3} to $\gamma$, we get a contradiction. 

Therefore $\alpha\delta\epsilon$ is a vertex. The angle sums of $\gamma\epsilon^2,\alpha\delta\epsilon,\alpha^2\gamma$ imply $\gamma=\delta$, contradicting Lemma \ref{geometry11}. Therefore $\gamma\epsilon^2,\alpha\delta\epsilon$ and one of $\alpha^3,\alpha\gamma^2,\gamma^3$ are all the degree $3$ vertices, and all three actually appear. The angle sums of all degree $3$ vertices and the angle sum for pentagon imply 
\begin{align*}
\alpha^3 &\colon
	\alpha=\tfrac{2}{3}\pi,\,
	\gamma=(\tfrac{1}{3}+\tfrac{4}{f})\pi,\,
	\delta=(\tfrac{1}{2}+\tfrac{2}{f})\pi,\,
	\epsilon=(\tfrac{5}{6}-\tfrac{2}{f})\pi. \\
\alpha\gamma^2 &\colon
	\alpha=\tfrac{8}{f}\pi,\,
	\gamma=(1-\tfrac{4}{f})\pi,\,
	\delta=(\tfrac{3}{2}-\tfrac{10}{f})\pi,\,
	\epsilon=(\tfrac{1}{2}+\tfrac{2}{f})\pi. \\
\gamma^3 &\colon
	\alpha=(\tfrac{1}{3}+\tfrac{4}{f})\pi,\,
	\gamma=\epsilon=\tfrac{2}{3}\pi,\,
	\delta=(1-\tfrac{4}{f})\pi. 
\end{align*}

For $\alpha\gamma^2,\gamma^3$, we have $\alpha<\gamma$, and $\alpha+2\delta>2\pi$, and $\delta>\epsilon>\frac{1}{2}\pi$. This implies $R(\delta^2)<\alpha,\gamma,2\delta,2\epsilon$. Therefore $\delta^2\cdots$ is not a vertex, contradicting $\gamma\epsilon^2$ and the balance lemma. 

For $\alpha^3$, by $\gamma\epsilon^2$ and the balance lemma, we know $\delta^2\cdots$ is a vertex. We have $\alpha>\gamma$. Then by $R(\delta^2)=(1-\frac{4}{f})\pi<\alpha+\gamma,3\gamma,2\delta,2\epsilon$, and $\alpha>\gamma$, we get $\delta^2\cdots=\gamma^2\delta^2$. The angle sum of $\gamma^2\delta^2$ further implies
\[
	\alpha=\tfrac{2}{3}\pi,\,
	\gamma=\tfrac{4}{9}\pi,\,
	\delta=\tfrac{5}{9}\pi,\,
	\epsilon=\tfrac{7}{9}\pi,\,
	f=36.
\]
By the angle values, we know $\gamma\epsilon^2,\dot{\alpha}^3,\dot{\alpha}\delta\epsilon,\dot{\alpha}\gamma^3,\gamma^2\delta^2$ ($\dot{\alpha}$ is $\alpha$ or $\beta$) are all the vertices.

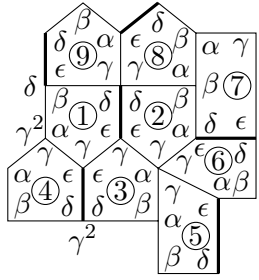
\begin{figure}[htp]
\centering
\begin{tikzpicture}[>=latex,scale=1]

\draw
	(0,-0.7) -- (0,0.7) -- (0.5,1.1) -- (1,0.7) -- (1,-0.7) -- (0.5,-1.1) -- (0,-0.7) -- (-0.5,-1.1) -- (-1,-0.7) -- (-1,0.7) -- (-0.5,1.1) -- (0,0.7)
	(-0.5,-1.1) -- (-0.5,-1.8)
	(-1,-0.7) -- (-1.5,-1.1) -- (-1.5,-1.8) -- (0.5,-1.8) 
	(1,0.7) -- (1.8,0.7) -- (1.8,-1.5) -- (1.3,-1.5) -- (0.5,-1.1) -- (0.5,-2.5) -- (1.3,-2.5) -- (1.3,-1.5)
	(1,-0.7) -- (1.8,-0.7)
	(-1,0) -- (1,0);

\draw[line width=1.2]
	(0,-0.7) -- (0,0)
	(-1,0.7) -- (-1,0)
	(0,0.7) -- (0.5,1.1)
	(-0.5,-1.1) -- (-0.5,-1.8)
	(1.3,-2.5) -- (1.3,-1.5)
	(1,-0.7) -- (1.8,-0.7);

\node at (0.2,0.6) {\small $\epsilon$};
\node at (0.2,0.2) {\small $\gamma$};
\node at (0.5,0.85) {\small $\delta$};	
\node at (0.8,0.2) {\small $\alpha$}; 
\node at (0.8,0.6) {\small $\beta$};

\node at (0.2,-0.6) {\small $\epsilon$};
\node at (0.5,-0.85) {\small $\gamma$};
\node at (0.2,-0.2) {\small $\delta$};	
\node at (0.8,-0.2) {\small $\beta$}; 
\node at (0.8,-0.6) {\small $\alpha$};

\node at (-0.2,0.6) {\small $\alpha$};
\node at (-0.5,0.85) {\small $\beta$};
\node at (-0.2,0.2) {\small $\gamma$};	
\node at (-0.8,0.2) {\small $\epsilon$}; 
\node at (-0.8,0.6) {\small $\delta$};

\node at (-0.2,-0.6) {\small $\epsilon$};
\node at (-0.5,-0.85) {\small $\gamma$};
\node at (-0.2,-0.2) {\small $\delta$};	
\node at (-0.8,-0.2) {\small $\beta$}; 
\node at (-0.8,-0.6) {\small $\alpha$};

\node at (-0.3,-1.2) {\small $\epsilon$};
\node at (0,-0.95) {\small $\gamma$};
\node at (-0.3,-1.6) {\small $\delta$};
\node at (0.3,-1.2) {\small $\alpha$};
\node at (0.3,-1.6) {\small $\beta$};

\node at (-1.3,-1.2) {\small $\alpha$};
\node at (-1,-0.95) {\small $\gamma$};
\node at (-1.3,-1.6) {\small $\beta$};
\node at (-0.7,-1.2) {\small $\epsilon$};
\node at (-0.7,-1.6) {\small $\delta$};

\node at (1.6,0.5) {\small $\gamma$};
\node at (1.6,-0.5) {\small $\epsilon$};
\node at (1.2,0.5) {\small $\alpha$};
\node at (1.2,-0.5) {\small $\delta$};
\node at (1.2,0) {\small $\beta$};

\node at (1.1,-2.3) {\small $\delta$};
\node at (1.1,-1.65) {\small $\epsilon$};
\node at (0.7,-2.3) {\small $\beta$};
\node at (0.7,-1.45) {\small $\gamma$};
\node at (0.7,-1.8) {\small $\alpha$};

\node at (1.05,-0.85) {\small $\epsilon$};
\node at (1.6,-0.9) {\small $\delta$};
\node at (0.85,-1.05) {\small $\gamma$};
\node at (1.6,-1.3) {\small $\beta$};
\node at (1.35,-1.3) {\small $\alpha$};

\node at (-0.5,-2.05) {\small $\gamma^2$};
\node at (-1.2,0) {\small $\delta$}; 
\node at (-1.2,-0.6) {\small $\gamma^2$};	

\node[inner sep=0.5,draw,shape=circle] at (-0.5,0.4) {\small 9};
\node[inner sep=0.5,draw,shape=circle] at (0.5,0.4) {\small 8};
\node[inner sep=0.5,draw,shape=circle] at (0.5,-0.4) {\small 2};
\node[inner sep=0.5,draw,shape=circle] at (-0.5,-0.4) {\small 1};
\node[inner sep=0.5,draw,shape=circle] at (0,-1.4) {\small 3};
\node[inner sep=0.5,draw,shape=circle] at (-1,-1.4) {\small 4};
\node[inner sep=0.5,draw,shape=circle] at (1.55,0) {\small 7};
\node[inner sep=0.5,draw,shape=circle] at (1.3,-1) {\small 6};
\node[inner sep=0.5,draw,shape=circle] at (1,-2) {\small 5};

\end{tikzpicture}
\caption{Proposition \ref{4b2d}: $\gamma\epsilon^2,\alpha\delta\epsilon,\alpha^3$ are vertices.}
\label{4b2dA}
\end{figure}

The AAD $\thick^{\delta}\epsilon^{\gamma}\thin^{\alpha}\gamma^{\epsilon}\thin^{\gamma}\epsilon^{\delta}\thick$ of $\gamma\epsilon^2$ determines $T_1,T_2,T_3$ in Figure \ref{4b2dA}. Then $\gamma_1\epsilon_3\cdots=\gamma\epsilon^2$ determines $T_4$. Then $\alpha_3\gamma_2\cdots=\dot{\alpha}\gamma^3$ gives $\gamma_5,\gamma_6$. If $T_5$ is not arranged as indicated, then $\beta_3\cdots=\dot{\alpha}\epsilon\cdots=\dot{\alpha}\delta\epsilon$ and $\delta_3\delta_4\cdots=\gamma^2\delta^2$ imply $\gamma,\delta$ adjacent, a contradiction. Therefore $T_5$ is arranged as indicated. Then by $\gamma_6$ and no $\epsilon\thin\epsilon\cdots$, we determine $T_6$. Then $\alpha_2\epsilon_6\cdots=\dot{\alpha}\delta\epsilon$ determines $T_7$. Then $\beta_2\beta_7\cdots=\dot{\alpha}^3$ and $\delta_1\delta_2\cdots=\gamma^2\delta^2$ determine $T_8$, and gives $\gamma_9$. Then by $\gamma_9$ and no $\epsilon\thin\epsilon\cdots$, we determine $T_9$. Then $\beta_1\epsilon_9\cdots=\dot{\alpha}\delta\epsilon$ and $\alpha_1\gamma_4\cdots=\dot{\alpha}\gamma^3$ imply $\gamma,\delta$ adjacent, a contradiction. 
\end{proof}

\begin{proposition}\label{4bde}
There is no tiling, such that $\alpha=\beta\ne\gamma$, and $\gamma\delta\epsilon$ is a vertex.
\end{proposition}

\begin{proof}
By $\alpha=\beta\ne\gamma$ and Propositions \ref{4a2d}, \ref{4a2e}, \ref{4b2d}, \ref{4b2e}, we know $\gamma\delta\epsilon$ is the only degree $3$ $b$-vertex. Then degree $3$ vertices are $\gamma\delta\epsilon$ and one of $\alpha^3,\alpha^2\gamma,\alpha\gamma^2,\gamma^3$. 

By applying Lemma \ref{degree3} to $\gamma$, we know that $\gamma\delta\epsilon,\alpha^3$ are all the degree $3$ vertices. On the other hand, the angle sums of $\gamma\delta\epsilon,\alpha^3$ and the angle sum for pentagon imply $f=12$, a contradiction. 
\end{proof}

\begin{proposition}\label{4ade}
Tilings with $\alpha=\beta\ne\gamma$, and such that $\alpha\delta\epsilon$ is a vertex, are the following:
\begin{enumerate}
\item Pentagonal subdivisions $P_{\pentagon}P_8,P_{\pentagon}P_{20}$ of the octahedron and the icosahedron.
\item The second earth map tiling {\rm $E_{\pentagon}2$} in Figure \ref{emt}, and its modifications {\rm $RE_{\pentagon}2,F_1E_{\pentagon}2,F_2E_{\pentagon}2,F'_2E_{\pentagon}2,F''_2E_{\pentagon}2$}.
\item The sporadic tiling {\rm $S_{16\pentagon}$}.
\end{enumerate}
\end{proposition}

The pentagonal subdivision $P_{\pentagon}P_4$ of the tetrahedron is not included because all propositions assume $f\ge 16$.

\begin{proof}
By $\alpha=\beta\ne\gamma$ and Propositions \ref{4a2d}, \ref{4a2e}, \ref{4b2d}, \ref{4b2e}, we know $\alpha\delta\epsilon$ is the only degree $3$ $b$-vertex. Then degree $3$ vertices are $\alpha\delta\epsilon$ and one of $\alpha^3,\alpha^2\gamma,\alpha\gamma^2,\gamma^3$. 

\subsubsection*{Case. $\alpha\delta\epsilon,\alpha^3$ are vertices}

The assumption implies $\gamma$ does not appear at degree $3$ vertices. Then by Lemma \ref{ndegree3}, we know one of $\alpha\gamma^3,\gamma^4,\gamma^5$ is a vertex. The angle sums of $\alpha\gamma\epsilon,\alpha^3$, and the angle sum of one of $\alpha\gamma^3,\gamma^4,\gamma^5$, and the angle sum for pentagon imply
\begin{align*}
\alpha\gamma^3 &\colon
	\alpha=\tfrac{2}{3}\pi,\,
	\gamma=\tfrac{4}{9}\pi,\,
	\delta+\epsilon=\tfrac{4}{3}\pi,\,
	f=36. \\
\gamma^4 &\colon
	\alpha=\tfrac{2}{3}\pi,\,
	\gamma=\tfrac{1}{2}\pi,\,
	\delta+\epsilon=\tfrac{4}{3}\pi,\,
	f=24. \\
\gamma^5 &\colon
	\alpha=\tfrac{2}{3}\pi,\,
	\gamma=\tfrac{2}{5}\pi,\,
	\delta+\epsilon=\tfrac{4}{3}\pi,\,
	f=60. 
\end{align*}

We have $\gamma<\alpha<2\gamma$. By Lemma \ref{geometry1}, this implies $\delta<\epsilon$. By $\delta+\epsilon>\pi$ and $\delta<\epsilon$, we know $R(\epsilon^2)$ has no $\delta,\epsilon$. By $\alpha\delta\epsilon$ and $\delta<\epsilon$, we get $R(\epsilon^2)<\alpha<2\gamma$. Then by the only degree $3$ $b$-vertex $\alpha\delta\epsilon$, we know $\epsilon^2\cdots$ is not a vertex. By the balance lemma, this implies a $b$-vertex is $\delta\epsilon\cdots$, with no $\delta,\epsilon$ in the remainder. Then by $\alpha\delta\epsilon$ and $\alpha<2\gamma$, we know $\alpha\delta\epsilon$ is the only $b$-vertex. By the values of $\alpha,\gamma$, we also get $\hat{b}$-vertices, and conclude ($\dot{\alpha}$ is $\alpha$ or $\beta$)
\begin{align*}
f=36 &\colon
	\text{AVC}=\{\dot{\alpha}\delta\epsilon,\dot{\alpha}^3,\dot{\alpha}\gamma^3\}. \\
f=24 &\colon
	\text{AVC}=\{\dot{\alpha}\delta\epsilon,\dot{\alpha}^3,\gamma^4\}. \\
f=60 &\colon
	\text{AVC}=\{\dot{\alpha}\delta\epsilon,\dot{\alpha}^3,\gamma^5\}.
\end{align*}

For $f=36$, we have $\dot{\alpha}\gamma^3=\alpha\gamma^3,\beta\gamma^3$. By no $\gamma\epsilon\cdots,\epsilon^2\cdots$, we know the AAD of $\alpha\gamma^3$ is $\thin^{\beta}\alpha^{\gamma}\thin^{\alpha}\gamma^{\epsilon}\thin^{\alpha}\gamma^{\epsilon}\thin^{\alpha}\gamma^{\epsilon}\thin$. This determines $T_1,T_2,T_3,T_4$ in the first of Figure \ref{4adeA}. Then $\alpha_2\gamma_1\cdots=\alpha\gamma^3$, and the AAD of $\alpha\gamma^3$ determines $T_5,T_6$. Then $\beta_1\epsilon_3\cdots=\alpha_5\epsilon_1\cdots=\alpha_6\epsilon_5\cdots=\beta_2\epsilon_6\cdots=\dot{\alpha}\delta\epsilon$ determines $T_7,T_8,T_9,T_{10}$. Then $\beta_5\beta_8\cdots=\beta_6\beta_9\cdots=\dot{\alpha}^3$ and $\delta_5\epsilon_9\cdots=\delta_6\epsilon_{10}\cdots=\dot{\alpha}\delta\epsilon$ give two $\dot{\alpha}$ in $T_{11},T_{12}$. Then by no $\gamma\delta\cdots$, we determine $T_{11},T_{12}$. Then $\alpha_7\gamma_8\cdots=\alpha\gamma^3$, and the AAD of $\alpha\gamma^3$ determines $T_{13},T_{14}$. Then $\alpha_9\delta_{12}\cdots=\alpha_{14}\epsilon_{13}\cdots=\dot{\alpha}\delta\epsilon$ determines $T_{15},T_{16}$. Then $\gamma_9\gamma_{11}\gamma_{15}\cdots=\dot{\alpha}\gamma^3$, and $\delta_{14}\epsilon_{11}\cdots=\dot{\alpha}\delta\epsilon$, and $\beta_{14}\beta_{16}\cdots=\dot{\alpha}^3$ imply three $\dot{\alpha}$ in a tile, a contradiction. 

By no $\delta\thin\epsilon\cdots,\epsilon^2\cdots$, we know the AAD of $\beta\gamma^3$ is $\thin^{\alpha}\beta^{\delta}\thin^{\alpha}\gamma^{\epsilon}\thin^{\alpha}\gamma^{\epsilon}\thin^{\alpha}\gamma^{\epsilon}\thin$. This determines $T_1,T_2,T_3,T_4$ in the second of Figure \ref{4adeA}. Then $\alpha_2\delta_1=\alpha_4\epsilon_2\cdots=\alpha_3\epsilon_4\cdots=\alpha_1\epsilon_3\cdots=\dot{\alpha}\delta\epsilon$ determines $T_5,T_6,T_7,T_8$. Then $\beta_2\gamma_5\cdots=\beta\gamma^3$, and the AAD of $\beta\gamma^3$ determines $T_9$. Moreover, $\beta_4\beta_6\cdots=\beta_3\beta_7\cdots=\dot{\alpha}^3$ and $\delta_2\epsilon_6\cdots=\delta_4\epsilon_7\cdots=\delta_3\epsilon_8\cdots=\dot{\alpha}\delta\epsilon$ give two $\dot{\alpha}$ in $T_{10},T_{11}$. Then by no $\gamma\delta\cdots$, we determine $T_{10},T_{11}$. The rest of the argument is the same as the first picture, and we get the same contradiction. 

This completes the proof that there is no tiling for $f=36$.

\begin{figure}[htp]
\centering
\begin{tikzpicture}[>=latex,yscale=-1]

\foreach \a in {0,1}
{
\begin{scope}[xshift=5.5*\a cm]

\draw
	(0,-1.1) -- (-0.5,-0.7) -- (-0.5,0.7) -- (0,1.1) -- (0.5,0.7) -- (0.5,-0.7) -- (0,-1.1)
	(-0.5,-0.7) -- (-1,-1.1) -- (-1.5,-0.7) -- (-1.5,0.7) -- (-1,1.1) -- (-0.5,0.7)
	(0.5,0) -- (-1.5,0)
	(0.5,-0.7) -- (1,-1.1) -- (1,-3.2) -- (-3.1,-3.2) -- (-3.1,-2.5) -- (0,-2.5)
	(-1.5,-0.7) -- (-3.1,-0.7) -- (-3.1,1.8) -- (-2.3,1.8)
	(1,-1.8) -- (-2.3,-1.8)
	(0,-1.1) -- (0,-3.2)
	(-2.3,-2.5) -- (-2.3,1.8) -- (0,1.8) -- (0,1.1)
	(-1,-1.1) -- (-1,-1.8)
	(-1,1.1) -- (-1,1.8)
	(-1.5,0.7) -- (-2.3,0.7);

\draw[line width=1.2]
	(-1,-1.8) -- (-2.3,-1.8)
	(0,-2.5) -- (0,-3.2)
	(-2.3,0.7) -- (-2.3,1.8)
	(-1.5,0) -- (-1.5,-0.7)
	(0,-1.1) -- (-0.5,-0.7)
	(-1,1.1) -- (-0.5,0.7);

\node at (-0.3,-0.6) {\small $\epsilon$};
\node at (0,-0.85) {\small $\delta$};
\node at (-0.3,-0.2) {\small $\gamma$};	
\node at (0.3,-0.2) {\small $\alpha$}; 
\node at (0.3,-0.6) {\small $\beta$};

\node at (-1.3,-0.6) {\small $\delta$};
\node at (-1,-0.85) {\small $\beta$};
\node at (-1.3,-0.2) {\small $\epsilon$};	
\node at (-0.7,-0.2) {\small $\gamma$}; 
\node at (-0.7,-0.6) {\small $\alpha$};

\node at (-1.3,0.6) {\small $\beta$};
\node at (-1,0.85) {\small $\delta$};
\node at (-1.3,0.2) {\small $\alpha$};	
\node at (-0.7,0.2) {\small $\gamma$}; 
\node at (-0.7,0.6) {\small $\epsilon$};

\node at (-0.5,-0.95) {\small $\delta$};
\node at (-0.8,-1.2) {\small $\beta$};
\node at (-0.2,-1.2) {\small $\epsilon$};
\node at (-0.8,-1.6) {\small $\alpha$};
\node at (-0.2,-1.6) {\small $\gamma$};

\node at (-0.5,0.95) {\small $\delta$};
\node at (-0.8,1.2) {\small $\epsilon$};
\node at (-0.2,1.2) {\small $\beta$};
\node at (-0.8,1.6) {\small $\gamma$};
\node at (-0.2,1.6) {\small $\alpha$};

\node at (-1.6,-0.9) {\small $\alpha$};
\node at (-1.2,-1.2) {\small $\beta$};
\node at (-2.1,-0.9) {\small $\gamma$};
\node at (-2.1,-1.6) {\small $\epsilon$};
\node at (-1.2,-1.6) {\small $\delta$};

\node at (-1.6,0.9) {\small $\beta$};
\node at (-1.2,1.2) {\small $\alpha$};
\node at (-2.1,0.9) {\small $\delta$};
\node at (-2.1,1.6) {\small $\epsilon$};
\node at (-1.2,1.6) {\small $\gamma$};

\node at (-2.1,-0.5) {\small $\gamma$};
\node at (-2.1,0.5) {\small $\alpha$};
\node at (-1.7,-0.5) {\small $\epsilon$};
\node at (-1.7,0.5) {\small $\beta$};
\node at (-1.7,0) {\small $\delta$};

\node at (-2.9,-0.5) {\small $\alpha$};
\node at (-2.5,-0.5) {\small $\gamma$};
\node at (-2.9,1.6) {\small $\beta$};
\node at (-2.5,0.7) {\small $\epsilon$};
\node at (-2.5,1.6) {\small $\delta$};

\node at (-2.5,-1.8) {\small $\dot{\alpha}$};
\node at (-2.5,-0.9) {\small $\dot{\alpha}$};
\node at (-2.5,-2.3) {\small $\dot{\alpha}$};

\node at (-2.9,-2.7) {\small $\alpha$};
\node at (-2.9,-3) {\small $\gamma$};
\node at (-2.3,-2.7) {\small $\beta$};
\node at (-0.2,-3) {\small $\epsilon$};
\node at (-0.2,-2.7) {\small $\delta$};

\node at (-0.2,-2.3) {\small $\alpha$};
\node at (-0.2,-2) {\small $\gamma$};
\node at (-2.1,-2.3) {\small $\beta$};
\node at (-1,-2) {\small $\epsilon$};
\node at (-2.1,-2) {\small $\delta$};

\node at (0.8,-2) {\small $\alpha$};
\node at (0.2,-2) {\small $\gamma$};
\node at (0.8,-3) {\small $\beta$};
\node at (0.2,-2.5) {\small $\epsilon$};
\node at (0.2,-3) {\small $\delta$};

\end{scope}
}


\draw
	(1,-1.1) -- (1.5,-0.7) -- (1.5,0.7) -- (1,1.1) -- (0.5,0.7)
	(1.5,0) -- (0.5,0);

\draw[line width=1.2]
	(0,1.1) -- (0.5,0.7)
	(1.5,0) -- (1.5,0.7)
	(1,-1.1) -- (0.5,-0.7)
	;

\node at (0.3,0.6) {\small $\epsilon$};
\node at (0,0.85) {\small $\delta$};
\node at (0.3,0.2) {\small $\gamma$};	
\node at (-0.3,0.2) {\small $\alpha$}; 
\node at (-0.3,0.6) {\small $\beta$};

\node at (1.3,0.6) {\small $\delta$};
\node at (1,0.85) {\small $\beta$};
\node at (1.3,0.2) {\small $\epsilon$};	
\node at (0.7,0.2) {\small $\gamma$}; 
\node at (0.7,0.6) {\small $\alpha$};

\node at (1.3,-0.6) {\small $\beta$};
\node at (1,-0.85) {\small $\delta$};
\node at (1.3,-0.2) {\small $\alpha$};	
\node at (0.7,-0.2) {\small $\gamma$}; 
\node at (0.7,-0.6) {\small $\epsilon$};

\node at (0.5,-0.95) {\small $\delta$};
\node at (0.2,-1.2) {\small $\beta$};
\node at (0.8,-1.2) {\small $\epsilon$};
\node at (0.2,-1.6) {\small $\alpha$};
\node at (0.8,-1.6) {\small $\gamma$};

\node[inner sep=0.5,draw,shape=circle] at (0,-0.4) {\small 1};
\node[inner sep=0.5,draw,shape=circle] at (0,0.4) {\small 2};
\node[inner sep=0.5,draw,shape=circle] at (1,-0.4) {\small 3};
\node[inner sep=0.5,draw,shape=circle] at (1,0.4) {\small 4};
\node[inner sep=0.5,draw,shape=circle] at (-1,-0.4) {\small 5};
\node[inner sep=0.5,draw,shape=circle] at (-1,0.4) {\small 6};
\node[inner sep=0.5,draw,shape=circle] at (0.5,-1.4) {\small 7};
\node[inner sep=0.5,draw,shape=circle] at (-0.5,-1.4) {\small 8};
\node[inner sep=0.5,draw,shape=circle] at (-2.05,0) {\small 9};
\node[inner sep=0,draw,shape=circle] at (-0.5,1.4) {\footnotesize 10};
\node[inner sep=0,draw,shape=circle] at (-1.7,-1.3) {\footnotesize 11};
\node[inner sep=0,draw,shape=circle] at (-1.7,1.3) {\footnotesize 12};
\node[inner sep=0,draw,shape=circle] at (-2.7,0) {\footnotesize 15};
\node[inner sep=0,draw,shape=circle] at (0.5,-2.5) {\footnotesize 13};
\node[inner sep=0,draw,shape=circle] at (-1.5,-2.15) {\footnotesize 14};
\node[inner sep=0,draw,shape=circle] at (-1.5,-2.85) {\footnotesize 16};


\begin{scope}[xshift=5.5cm]

\draw
	(0.5,0.7) -- (1.3,0.7) -- (1.3,-0.7) -- (0.5,-0.7);

\draw[line width=1.2]
	(0.5,0) -- (0.5,0.7)
	(1,-1.1) -- (1,-1.8);
	
\node at (0.3,0.6) {\small $\epsilon$};
\node at (0,0.85) {\small $\gamma$};
\node at (0.3,0.2) {\small $\delta$};	
\node at (-0.3,0.2) {\small $\beta$}; 
\node at (-0.3,0.6) {\small $\alpha$};

\node at (0.5,-0.95) {\small $\gamma$};
\node at (0.2,-1.2) {\small $\alpha$};
\node at (0.8,-1.2) {\small $\epsilon$};
\node at (0.2,-1.6) {\small $\beta$};
\node at (0.8,-1.6) {\small $\delta$};

\node at (1.1,-0.5) {\small $\alpha$};
\node at (1.1,0.5) {\small $\beta$};
\node at (0.7,-0.5) {\small $\gamma$};
\node at (0.7,0.5) {\small $\delta$};
\node at (0.7,0) {\small $\epsilon$};

\node at (1,-0.9) {\small $\gamma$};

\node[inner sep=0.5,draw,shape=circle] at (0,-0.4) {\small 2};
\node[inner sep=0.5,draw,shape=circle] at (0,0.4) {\small 1};
\node[inner sep=0.5,draw,shape=circle] at (-1,-0.4) {\small 4};
\node[inner sep=0.5,draw,shape=circle] at (-1,0.4) {\small 3};
\node[inner sep=0.5,draw,shape=circle] at (1.05,0) {\small 5};
\node[inner sep=0.5,draw,shape=circle] at (0.5,-1.4) {\small 9};
\node[inner sep=0.5,draw,shape=circle] at (-0.5,-1.4) {\small 6};
\node[inner sep=0.5,draw,shape=circle] at (-2.05,0) {\small 7};
\node[inner sep=0.5,draw,shape=circle] at (-0.5,1.4) {\small 8};
\node[inner sep=0,draw,shape=circle] at (-1.7,-1.3) {\footnotesize 10};
\node[inner sep=0,draw,shape=circle] at (-1.7,1.3) {\footnotesize 11};

\end{scope}


\begin{scope}[xshift=8cm]

\draw
	(0,1.1) -- (0.5,0.7) -- (0.5,-0.7) -- (0,-1.1) -- (-0.5,-0.7) -- (-0.5,0.7) -- (0,1.1)
	(0.5,0.7) -- (1,1.1) -- (1.5,0.7) -- (1.5,-0.7) -- (1,-1.1) -- (0.5,-0.7)
	(-0.5,0) -- (1.5,0)
	(0,-1.1) -- (0,-1.8) -- (2.3,-1.8) -- (2.3,0.7) -- (1.5,0.7)
	(1,-1.1) -- (1,-2.2)
	(2.3,-0.7) -- (1.5,-0.7);

\draw[line width=1.2]
	(2.3,-0.7) -- (2.3,-1.8)
	(1.5,0) -- (1.5,0.7)
	(-0.5,0) -- (-0.5,-0.7)
	(0,1.1) -- (0.5,0.7)
	(1,-1.1) -- (0.5,-0.7);
	
\node at (0.3,0.6) {\small $\epsilon$};
\node at (0,0.85) {\small $\delta$};
\node at (0.3,0.2) {\small $\gamma$};	
\node at (-0.3,0.2) {\small $\alpha$}; 
\node at (-0.3,0.6) {\small $\beta$};

\node at (1.3,0.6) {\small $\delta$};
\node at (1,0.85) {\small $\beta$};
\node at (1.3,0.2) {\small $\epsilon$};	
\node at (0.7,0.2) {\small $\gamma$}; 
\node at (0.7,0.6) {\small $\alpha$};

\node at (-0.3,-0.6) {\small $\delta$};
\node at (0,-0.85) {\small $\beta$};
\node at (-0.3,-0.2) {\small $\epsilon$};	
\node at (0.3,-0.2) {\small $\gamma$}; 
\node at (0.3,-0.6) {\small $\alpha$};

\node at (1.3,-0.6) {\small $\beta$};
\node at (1,-0.85) {\small $\delta$};
\node at (1.3,-0.2) {\small $\alpha$};	
\node at (0.7,-0.2) {\small $\gamma$}; 
\node at (0.7,-0.6) {\small $\epsilon$};

\node at (1.6,-0.9) {\small $\beta$};
\node at (1.2,-1.2) {\small $\alpha$};
\node at (2.1,-0.9) {\small $\delta$};
\node at (2.1,-1.6) {\small $\epsilon$};
\node at (1.2,-1.6) {\small $\gamma$};

\node at (2.1,0.5) {\small $\gamma$};
\node at (2.1,-0.5) {\small $\alpha$};
\node at (1.7,0.5) {\small $\epsilon$};
\node at (1.7,-0.5) {\small $\beta$};
\node at (1.7,0) {\small $\delta$};

\node at (0.5,-0.95) {\small $\delta$};
\node at (0.8,-1.2) {\small $\epsilon$};
\node at (0.2,-1.2) {\small $\beta$};
\node at (0.8,-1.6) {\small $\gamma$};
\node at (0.2,-1.6) {\small $\alpha$};

\node at (1.2,-2) {\small $\gamma$};
\node at (0.8,-2) {\small $\gamma$};

\node[inner sep=0.5,draw,shape=circle] at (0,0.4) {\small 4};
\node[inner sep=0.5,draw,shape=circle] at (1,0.4) {\small 2};
\node[inner sep=0.5,draw,shape=circle] at (1,-0.4) {\small 1};
\node[inner sep=0.5,draw,shape=circle] at (0,-0.4) {\small 3};
\node[inner sep=0.5,draw,shape=circle] at (1.7,-1.3) {\small 7};
\node[inner sep=0.5,draw,shape=circle] at (0.5,-1.4) {\small 6};
\node[inner sep=0.5,draw,shape=circle] at (2.05,0) {\small 5};

\end{scope}

\end{tikzpicture}
\caption{Proposition \ref{4ade}: $\alpha\delta\epsilon,\alpha^3$ are vertices.}
\label{4adeA}
\end{figure}
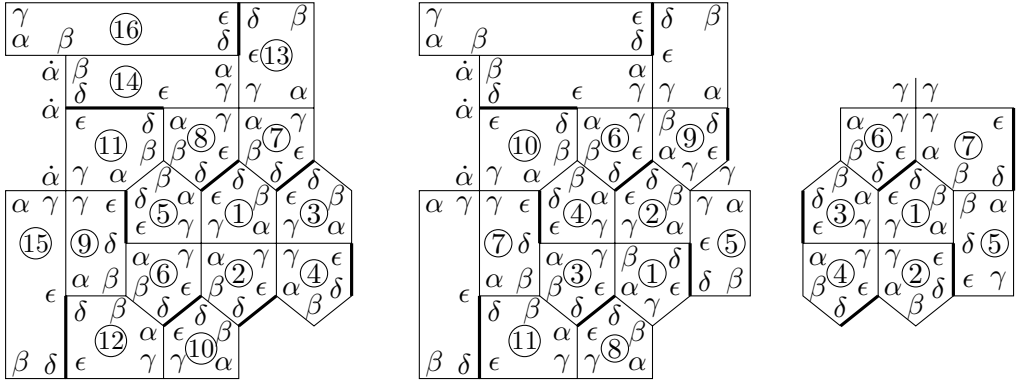

For $f=24$, by no $\epsilon^2\cdots$, we know the AAD of $\gamma^4$ is $\thin^{\alpha}\gamma^{\epsilon}\thin^{\alpha}\gamma^{\epsilon}\thin^{\alpha}\gamma^{\epsilon}\thin^{\alpha}\gamma^{\epsilon}\thin$. This determines $T_1,T_2,T_3,T_4$ in the third of Figure \ref{4adeA}. Then $\alpha_1\epsilon_2\cdots=\alpha_3\epsilon_1\cdots=\dot{\alpha}\delta\epsilon$ determines $T_5,T_6$. Then $\beta_1\beta_5\cdots=\dot{\alpha}^3$ and $\delta_1\epsilon_6\cdots=\dot{\alpha}\delta\epsilon$ give two $\dot{\alpha}$ in $T_7$. Then by no $\gamma\delta\cdots$, we determine $T_7$. Then $\gamma_6\gamma_7\cdots=\gamma^4$, and we may repeat the argument for this $\gamma^4$ in place of the initial $\gamma^4$. More repetitions imply the pentagonal subdivision $P_{\pentagon}P_8$ of the octahedron in Figure \ref{subdivision_tiling}. For $f=60$, similar argument with $\gamma^5$ in place of $\gamma^4$ implies the pentagonal subdivision $P_{\pentagon}P_{20}$ of the icosahedron in Figure \ref{subdivision_tiling}.

\medskip

\noindent{\em Geometry of Pentagon}

\medskip

For $f=24$, we have
\[
\alpha=\tfrac{2}{3}\pi,\,
\gamma=\tfrac{1}{2}\pi,\,
\delta+\epsilon=\tfrac{4}{3}\pi.
\]
By $\delta<\epsilon$, we get $\tfrac{2}{3}\pi<\epsilon<\tfrac{4}{3}\pi$. By \eqref{coolsaet_eq1}, we get 
\[
\cos^2\epsilon+2\cos\epsilon\sin\epsilon+(1-\tfrac{2}{3\sqrt{3}})\sin^2\epsilon=0.
\]
The solution satisfying $\tfrac{2}{3}\pi<\epsilon<\tfrac{4}{3}\pi$ is given by
\[
\tan\epsilon
=-\frac{3^{\frac{3}{4}}}{\sqrt{2}+3^{\frac{3}{4}}},\quad
\epsilon=0.8240\pi.
\]
Then we further find $a,b$ are given by
\begin{align*}
\cos a & ={\textstyle \tfrac{1}{3}\sqrt{2\sqrt{3}+3} }, &
a &=  0.1781\pi, \\
\cos b & ={\textstyle \tfrac{1}{3\sqrt{3}}(\sqrt{2\sqrt{3}+3}+2) }, &
b & = 0.1613\pi. 
\end{align*}
The base angles $\theta,\rho$ of isosceles triangles $\triangle ABD,\triangle ACE$ satisfy $\cos a=\cot\frac{1}{2}\gamma\cot\theta$ and $\cos a=\cot\frac{1}{2}\alpha\cot\rho$. By $\cos a,\cot\frac{1}{2}\alpha, \cot\frac{1}{2}\gamma>0$, we get $\theta,\rho<\frac{1}{2}\pi<\alpha$. Then by Lemma \ref{geometry8}, the pentagon is simple. 

For $f=60$, we get
\[
\alpha=\tfrac{2}{3}\pi,\,
\gamma=\tfrac{2}{5}\pi,\,
\delta+\epsilon=\tfrac{4}{3}\pi.
\]
Again we have $\tfrac{2}{3}\pi<\epsilon<\tfrac{4}{3}\pi$. Then \eqref{coolsaet_eq1} is a quadratic equation for $\tan\epsilon$, and the solution satisfying $\tfrac{2}{3}\pi<\epsilon<\tfrac{4}{3}\pi$ is given by 
\[
\tan\epsilon=\frac{2\sqrt{3}(6-\sqrt{15-9\sqrt{5}+\sqrt{870+114\sqrt{5}}})}{7+\sqrt{5}-3\sqrt{30+6\sqrt{5}}},\quad
\epsilon=0.9310\pi.
\]
Then we further find $a$ is given by
\[
\cos a ={\textstyle \tfrac{1}{\sqrt[4]{3^3}\sqrt[8]{5}}\sqrt{\sqrt{3}\sqrt[4]{5}+2\sqrt{2+\sqrt{5}}} }, \quad
a= 0.1206\pi.
\]
We also get $b=0.1129\pi$, and the exact expression for $\cos b$ is too complicated to be recorded here. By the same argument for $f=24$, the pentagon is simple.

\subsubsection*{Case. $\alpha\delta\epsilon,\gamma^3$ are vertices}

The angle sums of $\alpha\delta\epsilon,\gamma^3$ and the angle sum for pentagon imply
\[
\alpha=(\tfrac{1}{3}+\tfrac{4}{f})\pi,\,
\gamma=\tfrac{2}{3}\pi,\,
\delta+\epsilon=(\tfrac{5}{3}-\tfrac{4}{f})\pi.
\]
We have $\alpha<\gamma$, and $\delta+\epsilon>\pi$. By Lemma \ref{geometry1}, this implies $\delta>\epsilon$. By $\delta+\epsilon>\pi$ and $\delta>\epsilon$, we know $R(\delta^2)$ has no $\delta,\epsilon$. By $\alpha\delta\epsilon$ and $\delta>\epsilon$, we get $R(\delta^2)<\alpha<\gamma$. Therefore $\delta^2\cdots$ is not a vertex. Then by the balance lemma, we know a $b$-vertex is $\delta\epsilon\cdots$, with only $\alpha,\gamma$ in the remainder. Then by $\alpha\delta\epsilon$ and $\alpha<\gamma$, we know $\alpha\delta\epsilon$ is the only $b$-vertex. By the values of $\alpha,\gamma$, we also get $\hat{b}$-vertices, and conclude ($\dot{\alpha}$ is $\alpha$ or $\beta$)
\begin{align*}
f=36 &\colon
	\text{AVC}=\{\dot{\alpha}\delta\epsilon,\gamma^3,\dot{\alpha}^3\gamma\}. \\
f=24 &\colon
	\text{AVC}=\{\dot{\alpha}\delta\epsilon,\gamma^3,\dot{\alpha}^4\}. \\
f=60 &\colon
	\text{AVC}=\{\dot{\alpha}\delta\epsilon,\gamma^3,\dot{\alpha}^5\}.
\end{align*}

By no $\epsilon^2\cdots$, we know the AAD of $\gamma^2\cdots=\gamma^3$ is $\thin^{\alpha}\gamma^{\epsilon}\thin^{\alpha}\gamma^{\epsilon}\thin^{\alpha}\gamma^{\epsilon}\thin$. This implies the AAD of $\thin\gamma\thin\gamma\thin$ is $\thin^{\alpha}\gamma^{\epsilon}\thin^{\alpha}\gamma^{\epsilon}\thin$, and further implies no $\thin^{\beta}\alpha^{\gamma}\thin^{\gamma}\alpha^{\beta}\thin$. 

By no $\gamma\delta\cdots$, we do not have $\thin^{\beta}\alpha^{\gamma}\thin^{\delta}\beta^{\alpha}\thin$. By no $\delta^2\cdots$, we do not have $\thin^{\alpha}\beta^{\delta}\thin^{\delta}\beta^{\alpha}\thin$. Combined with no $\thin^{\beta}\alpha^{\gamma}\thin^{\gamma}\alpha^{\beta}\thin$, we do not have $\thin^{\dot{\alpha}}\dot{\alpha}\thin\dot{\alpha}^{\dot{\alpha}}\thin$. By no $\gamma\epsilon\cdots,\delta\thin\epsilon\cdots$, we also do not have $\thin^{\alpha}\gamma^{\epsilon}\thin\dot{\alpha}^{\dot{\alpha}}\thin$.

For $f=36$, by no $\thin^{\dot{\alpha}}\dot{\alpha}\thin\dot{\alpha}^{\dot{\alpha}}\thin$, we get $\thin\dot{\alpha}^{\dot{\alpha}}\thin^{\alpha}\gamma^{\epsilon}\thin\cdots=\dot{\alpha}^3\gamma=\thin\dot{\alpha}^{\dot{\alpha}}\thin\dot{\alpha}^{\dot{\alpha}}\thin\dot{\alpha}^{\dot{\alpha}}\thin^{\alpha}\gamma^{\epsilon}\thin=\thin^{\alpha}\gamma^{\epsilon}\thin\dot{\alpha}^{\dot{\alpha}}\thin\cdots$, contradicting no $\thin^{\alpha}\gamma^{\epsilon}\thin\dot{\alpha}^{\dot{\alpha}}\thin$. Therefore we do not have $\thin\dot{\alpha}^{\dot{\alpha}}\thin^{\alpha}\gamma^{\epsilon}\thin$. This further implies no $\thin^{\beta}\alpha^{\gamma}\thin^{\dot{\alpha}}\dot{\alpha}\thin$. Therefore $\thin^{\dot{\alpha}}\dot{\alpha}\thin^{\dot{\alpha}}\dot{\alpha}\thin=\thin^{\alpha}\beta^{\delta}\thin^{\dot{\alpha}}\dot{\alpha}\thin=\thin^{\alpha}\beta^{\delta}\thin^{\alpha}\beta^{\delta}\thin,\thin^{\alpha}\beta^{\delta}\thin^{\beta}\alpha^{\gamma}\thin$. The AAD $\thin^{\alpha}\beta^{\delta}\thin^{\beta}\alpha^{\gamma}\thin$ implies a vertex $\thin^{\delta}\beta^{\alpha}\thin^{\beta}\delta^{\epsilon}\thick\cdots=\dot{\alpha}\delta\epsilon=\thick^{\delta}\epsilon^{\gamma}\thin^{\delta}\beta^{\alpha}\thin^{\beta}\delta^{\epsilon}\thick$, contradicting no $\gamma\delta$. Then we conclude $\thin^{\dot{\alpha}}\dot{\alpha}\thin^{\dot{\alpha}}\dot{\alpha}\thin=\thin^{\alpha}\beta^{\delta}\thin^{\alpha}\beta^{\delta}\thin$.

By no $\thin^{\alpha}\gamma^{\epsilon}\thin\dot{\alpha}^{\dot{\alpha}}\thin,\thin\dot{\alpha}^{\dot{\alpha}}\thin^{\alpha}\gamma^{\epsilon}\thin$, we know the AAD of $\dot{\alpha}^3\gamma$ is $\thin^{\dot{\alpha}}\dot{\alpha}\thin^{\alpha}\gamma^{\epsilon}\thin^{\dot{\alpha}}\dot{\alpha}\thin\dot{\alpha}\thin$. Then by no $\thin^{\dot{\alpha}}\dot{\alpha}\thin\dot{\alpha}^{\dot{\alpha}}\thin$, the AAD is $\thin^{\dot{\alpha}}\dot{\alpha}\thin^{\dot{\alpha}}\dot{\alpha}\thin^{\dot{\alpha}}\dot{\alpha}\thin^{\alpha}\gamma^{\epsilon}\thin$. By $\thin^{\dot{\alpha}}\dot{\alpha}\thin^{\dot{\alpha}}\dot{\alpha}\thin=\thin^{\alpha}\beta^{\delta}\thin^{\alpha}\beta^{\delta}\thin$, we conclude $\dot{\alpha}^3\gamma=\thin^{\alpha}\beta^{\delta}\thin^{\alpha}\beta^{\delta}\thin^{\alpha}\beta^{\delta}\thin^{\alpha}\gamma^{\epsilon}\thin=\beta^3\gamma$. This determines $T_1,T_2,T_3,T_4$ in the first of Figure \ref{4adeB}. Then $\alpha_2\epsilon_1\cdots=\alpha_4\delta_2\cdots=\alpha_3\delta_4\cdots=\dot{\alpha}\delta\epsilon$ determines $T_5,T_6,T_7$. Then $\beta_5\gamma_2\cdots=\dot{\alpha}^3\gamma$, and the AAD of $\dot{\alpha}^3\gamma$ determines $T_8$. Then $\gamma_4\gamma_6\cdots=\gamma_3\gamma_7\cdots=\gamma^3$, and the AAD of $\gamma^3$ determines $T_9,T_{10}$. Then $\beta_6\gamma_8\cdots=\dot{\alpha}^3\gamma$, and the AAD of $\dot{\alpha}^3\gamma$ determines $T_{11},T_{12}$. Then $\alpha_7\epsilon_{10}\cdots=\alpha_{12}\delta_{11}\cdots=\dot{\alpha}\delta\epsilon$ determines $T_{13},T_{14}$. Then $\beta_7\beta_9\beta_{13}\cdots=\beta^3\gamma$, and $\delta_9\epsilon_{12}\cdots=\dot{\alpha}\delta\epsilon$, and $\gamma_{12}\gamma_{14}\cdots=\gamma^3$ imply two $\gamma$ in a tile, a contradiction. 

\begin{figure}[htp]
\centering
\begin{tikzpicture}[>=latex,scale=1]

\draw
	(0,-1.1) -- (-0.5,-0.7) -- (-0.5,0.7) -- (0,1.1) -- (0.5,0.7) -- (0.5,-0.7) -- (0,-1.1)
	(-0.5,-0.7) -- (-1,-1.1) -- (-1.5,-0.7) -- (-1.5,0.7) -- (-1,1.1) -- (-0.5,0.7)
	(0.5,0) -- (-1.5,0)
	(1,-1.8) -- (-2.3,-1.8)
	(0.5,-0.7) -- (1,-1.1) -- (1,-3.2) -- (-3.1,-3.2) -- (-3.1,-2.5) -- (0,-2.5)
	(-2.3,-2.5) -- (-2.3,1.8)
	(-2.3,0.7) -- (-1.5,0.7)
	(-1,-1.1) -- (-1,-1.8)
	(0,-1.1) -- (0,-3.2)
	(-1.5,-0.7) -- (-3.1,-0.7) -- (-3.1,1.8) -- (-1,1.8) -- (-1,1.1);

\draw[line width=1.2]
	(-1,-1.8) -- (-2.3,-1.8)
	(-2.3,0.7) -- (-2.3,1.8)
	(-1.5,0) -- (-1.5,-0.7)
	(0,-1.1) -- (-0.5,-0.7)
	(-1,1.1) -- (-0.5,0.7)
	(0,-2.5) -- (0,-3.2);

\draw
	(0.5,0.7) -- (1.3,0.7) -- (1.3,-0.7) -- (0.5,-0.7);

\draw[line width=1.2]
	(0.5,0) -- (0.5,0.7)
	(1,-1.1) -- (1,-1.8);
	
\node at (0.3,0.6) {\small $\delta$};
\node at (0,0.85) {\small $\beta$};
\node at (0.3,0.2) {\small $\epsilon$};	
\node at (-0.3,0.2) {\small $\gamma$}; 
\node at (-0.3,0.6) {\small $\alpha$};

\node at (0.3,-0.6) {\small $\gamma$};
\node at (0,-0.85) {\small $\epsilon$};
\node at (0.3,-0.2) {\small $\alpha$};	
\node at (-0.3,-0.2) {\small $\beta$}; 
\node at (-0.3,-0.6) {\small $\delta$};

\node at (-0.7,-0.6) {\small $\alpha$};
\node at (-1,-0.85) {\small $\gamma$};
\node at (-0.7,-0.2) {\small $\beta$};	
\node at (-1.3,-0.2) {\small $\delta$}; 
\node at (-1.3,-0.6) {\small $\epsilon$};

\node at (-0.7,0.6) {\small $\delta$};
\node at (-1,0.85) {\small $\epsilon$};
\node at (-0.7,0.2) {\small $\beta$};	
\node at (-1.3,0.2) {\small $\alpha$}; 
\node at (-1.3,0.6) {\small $\gamma$};

\node at (0.5,-0.95) {\small $\beta$};
\node at (0.2,-1.2) {\small $\alpha$};
\node at (0.8,-1.2) {\small $\delta$};
\node at (0.2,-1.6) {\small $\gamma$};
\node at (0.8,-1.6) {\small $\epsilon$};

\node at (1.1,-0.5) {\small $\alpha$};
\node at (1.1,0.5) {\small $\gamma$};
\node at (0.7,-0.5) {\small $\beta$};
\node at (0.7,0.5) {\small $\epsilon$};
\node at (0.7,0) {\small $\delta$};

\node at (0.95,-0.85) {\small $\beta$};

\node at (-0.5,-0.95) {\small $\epsilon$};
\node at (-0.2,-1.2) {\small $\delta$};
\node at (-0.8,-1.2) {\small $\gamma$};
\node at (-0.2,-1.6) {\small $\beta$};
\node at (-0.8,-1.6) {\small $\alpha$};

\node at (-2.1,-0.5) {\small $\beta$};
\node at (-2.1,0.5) {\small $\alpha$};
\node at (-1.7,-0.5) {\small $\delta$};
\node at (-1.7,0.5) {\small $\gamma$};
\node at (-1.7,0) {\small $\epsilon$};

\node at (-0.2,-2.3) {\small $\alpha$};
\node at (-0.2,-2) {\small $\beta$};
\node at (-2.1,-2.3) {\small $\gamma$};
\node at (-1,-2) {\small $\delta$};
\node at (-2.1,-2) {\small $\epsilon$};

\node at (-1.6,0.9) {\small $\gamma$};
\node at (-1.2,1.2) {\small $\alpha$};
\node at (-2.1,0.9) {\small $\epsilon$};
\node at (-2.1,1.6) {\small $\delta$};
\node at (-1.2,1.6) {\small $\beta$};

\node at (-1.6,-0.9) {\small $\alpha$};
\node at (-1.2,-1.2) {\small $\gamma$};
\node at (-2.1,-0.9) {\small $\beta$};
\node at (-2.1,-1.6) {\small $\delta$};
\node at (-1.2,-1.6) {\small $\epsilon$};

\node at (-2.9,-0.5) {\small $\alpha$};
\node at (-2.9,1.6) {\small $\gamma$};
\node at (-2.5,-0.5) {\small $\beta$};
\node at (-2.5,1.6) {\small $\epsilon$};
\node at (-2.5,0.7) {\small $\delta$};

\node at (-2.5,-1.8) {\small $\dot{\alpha}$};
\node at (-2.5,-0.9) {\small $\gamma$};
\node at (-2.5,-2.3) {\small $\gamma$};

\node at (-2.9,-2.7) {\small $\alpha$};
\node at (-2.9,-3) {\small $\beta$};
\node at (-2.3,-2.7) {\small $\gamma$};
\node at (-0.2,-3) {\small $\delta$};
\node at (-0.2,-2.7) {\small $\epsilon$};

\node at (0.8,-2) {\small $\alpha$};
\node at (0.2,-2) {\small $\beta$};
\node at (0.8,-3) {\small $\gamma$};
\node at (0.2,-2.5) {\small $\delta$};
\node at (0.2,-3) {\small $\epsilon$};

\node[inner sep=0.5,draw,shape=circle] at (0,0.4) {\small 1};
\node[inner sep=0.5,draw,shape=circle] at (0,-0.4) {\small 2};
\node[inner sep=0.5,draw,shape=circle] at (-1,0.4) {\small 3};
\node[inner sep=0.5,draw,shape=circle] at (-1,-0.4) {\small 4};
\node[inner sep=0.5,draw,shape=circle] at (1.05,0) {\small 5};
\node[inner sep=0.5,draw,shape=circle] at (0.5,-1.4) {\small 8};
\node[inner sep=0.5,draw,shape=circle] at (-0.5,-1.4) {\small 6};
\node[inner sep=0.5,draw,shape=circle] at (-2.05,0) {\small 7};
\node[inner sep=0.5,draw,shape=circle] at (-1.7,-1.3) {\small 9};
\node[inner sep=0,draw,shape=circle] at (-1.7,1.3) {\footnotesize 10};
\node[inner sep=0,draw,shape=circle] at (-2.7,0) {\footnotesize 13};
\node[inner sep=0,draw,shape=circle] at (0.5,-2.5) {\footnotesize 11};
\node[inner sep=0,draw,shape=circle] at (-1.5,-2.15) {\footnotesize 12};
\node[inner sep=0,draw,shape=circle] at (-1.5,-2.85) {\footnotesize 14};


\begin{scope}[xshift=2.5cm]

\draw
	(0,1.1) -- (0.5,0.7) -- (0.5,-0.7) -- (0,-1.1) -- (-0.5,-0.7) -- (-0.5,0.7) -- (0,1.1)
	(0.5,0.7) -- (1,1.1) -- (1.5,0.7) -- (1.5,-0.7) -- (1,-1.1) -- (0.5,-0.7)
	(-0.5,0) -- (1.5,0)
	(0,-1.1) -- (0,-1.8) -- (2.3,-1.8) -- (2.3,0.7) -- (1.5,0.7)
	(1,-1.1) -- (1,-2.2)
	(2.3,-0.7) -- (1.5,-0.7);

\draw[line width=1.2]
	(2.3,-0.7) -- (2.3,-1.8)
	(1.5,0) -- (1.5,0.7)
	(-0.5,0) -- (-0.5,-0.7)
	(0,1.1) -- (0.5,0.7)
	(1,-1.1) -- (0.5,-0.7);
	
\node at (0.3,0.6) {\small $\delta$};
\node at (0,0.85) {\small $\epsilon$};
\node at (0.3,0.2) {\small $\beta$};	
\node at (-0.3,0.2) {\small $\alpha$}; 
\node at (-0.3,0.6) {\small $\gamma$};

\node at (1.3,0.6) {\small $\epsilon$};
\node at (1,0.85) {\small $\gamma$};
\node at (1.3,0.2) {\small $\delta$};	
\node at (0.7,0.2) {\small $\beta$}; 
\node at (0.7,0.6) {\small $\alpha$};

\node at (-0.3,-0.6) {\small $\epsilon$};
\node at (0,-0.85) {\small $\gamma$};
\node at (-0.3,-0.2) {\small $\delta$};	
\node at (0.3,-0.2) {\small $\beta$}; 
\node at (0.3,-0.6) {\small $\alpha$};

\node at (1.3,-0.6) {\small $\gamma$};
\node at (1,-0.85) {\small $\epsilon$};
\node at (1.3,-0.2) {\small $\alpha$};	
\node at (0.7,-0.2) {\small $\beta$}; 
\node at (0.7,-0.6) {\small $\delta$};

\node at (1.6,-0.9) {\small $\gamma$};
\node at (1.2,-1.2) {\small $\alpha$};
\node at (2.1,-0.9) {\small $\epsilon$};
\node at (2.1,-1.6) {\small $\delta$};
\node at (1.2,-1.6) {\small $\beta$};

\node at (2.1,0.5) {\small $\beta$};
\node at (2.1,-0.5) {\small $\alpha$};
\node at (1.7,0.5) {\small $\delta$};
\node at (1.7,-0.5) {\small $\gamma$};
\node at (1.7,0) {\small $\epsilon$};

\node at (0.5,-0.95) {\small $\epsilon$};
\node at (0.8,-1.2) {\small $\delta$};
\node at (0.2,-1.2) {\small $\gamma$};
\node at (0.8,-1.6) {\small $\beta$};
\node at (0.2,-1.6) {\small $\alpha$};

\node at (1.2,-2) {\small $\beta$};
\node at (0.8,-2) {\small $\beta$};

\node[inner sep=0.5,draw,shape=circle] at (0,0.4) {\small 3};
\node[inner sep=0.5,draw,shape=circle] at (1,0.4) {\small 1};
\node[inner sep=0.5,draw,shape=circle] at (1,-0.4) {\small 2};
\node[inner sep=0.5,draw,shape=circle] at (0,-0.4) {\small 4};
\node[inner sep=0.5,draw,shape=circle] at (1.7,-1.3) {\small 7};
\node[inner sep=0.5,draw,shape=circle] at (0.5,-1.4) {\small 6};
\node[inner sep=0.5,draw,shape=circle] at (2.05,0) {\small 5};

\end{scope}
	
\end{tikzpicture}
\caption{Proposition \ref{4ade}: $\alpha\delta\epsilon,\gamma^3$ are vertices.}
\label{4adeB}
\end{figure}
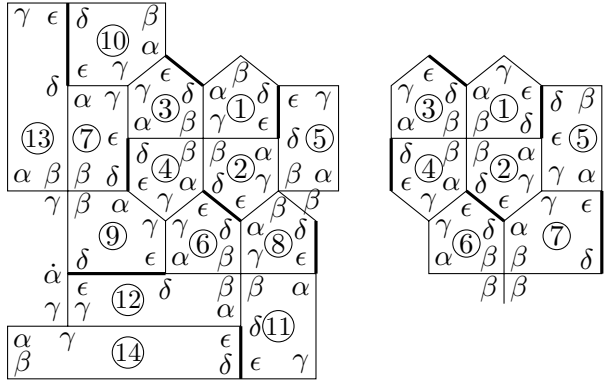

For $f=24$, by no $\thin^{\dot{\alpha}}\dot{\alpha}\thin\dot{\alpha}^{\dot{\alpha}}\thin$, we know the AAD of $\dot{\alpha}^4$ is $\thin^{\dot{\alpha}}\dot{\alpha}\thin^{\dot{\alpha}}\dot{\alpha}\thin^{\dot{\alpha}}\dot{\alpha}\thin^{\dot{\alpha}}\dot{\alpha}\thin$. Then by no $\dot{\alpha}\gamma\cdots$, we get $\dot{\alpha}^4=\thin^{\alpha}\beta^{\delta}\thin^{\alpha}\beta^{\delta}\thin^{\alpha}\beta^{\delta}\thin^{\alpha}\beta^{\delta}\thin=\beta^4$. This determines $T_1,T_2,T_3,T_4$ in the second of Figure \ref{4adeB}. Then $\alpha_2\delta_1\cdots=\alpha_4\delta_2\cdots=\dot{\alpha}\delta\epsilon$ determines $T_5,T_6$. Then $\gamma_2\gamma_5\cdots=\gamma^3$ and $\epsilon_2\delta_5\cdots=\dot{\alpha}\delta\epsilon$ determine $T_7$. Then $\beta_6\beta_7\cdots=\dot{\alpha}^4$, and we may repeat the argument for this $\dot{\alpha}^4$ in place of the initial $\dot{\alpha}^4$. More repetitions imply the pentagonal subdivision of the octahedron. For $f=60$, similar argument with $\dot{\alpha}^5$ in place of $\dot{\alpha}^4$ implies the pentagonal subdivision of the icosahedron. After the exchange $(\beta,\delta)\leftrightarrow(\gamma,\delta)$, these are the pentagonal subdivisions $P_{\pentagon}P_8,P_{\pentagon}P_{20}$ in Figure \ref{subdivision_tiling}.

\medskip

\noindent{\em Geometry of Pentagon}

\medskip

For $f=24$, we get
\[
\alpha=\tfrac{1}{2}\pi,\,
\gamma=\tfrac{2}{3}\pi,\,
\delta+\epsilon=\tfrac{3}{2}\pi.
\]
By $\delta>\epsilon$, we get $\epsilon<\tfrac{3}{4}\pi$. By \eqref{coolsaet_eq2} and \eqref{coolsaet_eq3}, we get 
\[
\cos a=\cot\tfrac{1}{2}\alpha\cot\delta
=\cot\tfrac{1}{2}\gamma\cot\epsilon.
\]
Substituting the angle values, we find the solution satiswfying $\epsilon<\tfrac{3}{4}\pi$ is given by
\[
\cos\epsilon=\tfrac{1}{\sqrt{2}}{\textstyle \sqrt{3-\sqrt{3}}},\quad
\epsilon= 0.2068\pi.
\]
Then we further find $a,b$ are given by
\begin{align*}
\cos a & =\tfrac{1}{\sqrt[4]{3}}, &
a &= 0.2252\pi, \\
\cos b & =\tfrac{1}{2\sqrt{3}}(\sqrt{3}+2\sqrt[4]{3}-1), &
b & = 0.0766\pi. 
\end{align*}
The base angles of $\triangle ABD$ and $\triangle ACE$ are $\epsilon,\delta-\pi<\alpha$. Then by Lemma \ref{geometry8}, we know the pentagon is simple.

For $f=60$, we get
\[
\alpha=\tfrac{2}{5}\pi,\,
\gamma=\tfrac{2}{3}\pi,\,
\delta+\epsilon=\tfrac{8}{5}\pi.
\]
By $\delta>\epsilon$, we get $\epsilon<\tfrac{4}{5}\pi$. By \eqref{coolsaet_eq1}, we get a quadratic equation for $\tan\epsilon$ 
\begin{align*}
&(3+4\sqrt{5}-{\textstyle \sqrt{15(5-2\sqrt{5})} })\tan^2\epsilon
+2{\textstyle \sqrt{10(25-11\sqrt{5})} }\tan\epsilon \\
&+(5(5-2\sqrt{5})-{\textstyle \sqrt{15(5-2\sqrt{5})} }
=0.
\end{align*}
The solution of $\tan\epsilon$ satisfying $\epsilon<\tfrac{4}{5}\pi$ is the smaller of the two roots. This determines $\epsilon$, which has the approximate value $\epsilon=0.0126\pi$. Then we find 
\[
\cos a
=\tfrac{1}{\sqrt{15}}{\textstyle \sqrt{-3\sqrt{5}+\sqrt{6(25+11\sqrt{5})}  }},\quad
a= 0.1835\pi,
\]
and $b= 0.1349\pi$ (the exact expression for $\cos b$ is too complicated). Then we calculate the base angles of isosceles triangles $\triangle ABD,\triangle ACE$ and get $0.1919\pi, 0.3258\pi<\alpha$. Then by Lemma \ref{geometry8}, the pentagon is simple. 

\subsubsection*{Case. $\alpha\delta\epsilon,\alpha\gamma^2$ are vertices}

The angle sums of $\alpha\delta\epsilon,\alpha\gamma^2$ and the angle sum for pentagon imply
\[
\alpha=\tfrac{8}{f}\pi,\,
\gamma=(1-\tfrac{4}{f})\pi,\,
\delta+\epsilon=(2-\tfrac{8}{f})\pi.
\]
We have $\alpha<\gamma$, and $\delta+\epsilon>\pi$. By Lemma \ref{geometry1} and $\delta+\epsilon=2\gamma$, this implies $\delta>\gamma>\epsilon$. By the same argument for the case $\alpha\delta\epsilon,\gamma^3$ are vertices, we know $\alpha\delta\epsilon$ is the only $b$-vertex. This implies that all companion pairs are twisted, and a tile determines its companion.

Applying Lemma \ref{geometry5} to $\alpha,\beta$, we get $3\alpha>\pi$. This means $f<24$, or $f=16,18,20,22$. Then by the values of $\alpha,\gamma$, we also get $\hat{b}$-vertices, and conclude ($\dot{\alpha}$ is $\alpha$ or $\beta$)
\begin{align*}
f=16 &\colon
	\text{AVC}=\{\dot{\alpha}\delta\epsilon,\dot{\alpha}\gamma^2,\dot{\alpha}^4\}. \\
f=20 &\colon
	\text{AVC}=\{\dot{\alpha}\delta\epsilon,\dot{\alpha}\gamma^2,\dot{\alpha}^3\gamma,\dot{\alpha}^5\}. 
\end{align*}
The cases $f=18,20$ are dismissed because all vertices have degree $3$.

For $f=20$, we get
\[
\dot{\alpha}=\tfrac{2}{5}\pi,\,
\gamma=\tfrac{4}{5}\pi,\,
\delta+\epsilon=\tfrac{8}{5}\pi.
\]
Then by \eqref{coolsaet_eq1}, we get 
\[
(\sqrt{5}(\sqrt{5}-1)\sin\epsilon+{\textstyle \sqrt{50-22\sqrt{5}} }\cos\epsilon)\sin\epsilon=0.
\]
There is no solution satisfying $\epsilon<\gamma=\tfrac{4}{5}\pi$.

Therefore $f=16$. By no $\gamma\epsilon\cdots,\delta\thin\epsilon\cdots$, we do not have $\thin^{\dot{\alpha}}\dot{\alpha}\thin^{\epsilon}\gamma^{\alpha}\thin$. This implies $\dot{\alpha}\gamma^2=\thin\gamma\thin^{\dot{\alpha}}\dot{\alpha}\thin\gamma\thin=\thin\gamma\thin^{\dot{\alpha}}\dot{\alpha}\thin^{\alpha}\gamma^{\epsilon}\thin=\thin^{\dot{\alpha}}\dot{\alpha}\thin^{\alpha}\gamma^{\epsilon}\thin\gamma\thin$. Then by no $\epsilon^2\cdots$, we get $\dot{\alpha}\gamma^2=\thin^{\dot{\alpha}}\dot{\alpha}\thin^{\alpha}\gamma^{\epsilon}\thin^{\alpha}\gamma^{\epsilon}\thin$. By $\dot{\alpha}\gamma\cdots=\gamma^2\cdots=\dot{\alpha}\gamma^2$, this implies the AAD of $\thin\gamma\thin\gamma\thin$ is $\thin^{\alpha}\gamma^{\epsilon}\thin^{\alpha}\gamma^{\epsilon}\thin$, and we do not have $\thin\dot{\alpha}^{\dot{\alpha}}\thin^{\alpha}\gamma^{\epsilon}\thin$. This implies no $\thin^{\beta}\alpha^{\gamma}\thin^{\gamma}\alpha^{\beta}\thin,\thin^{\beta}\alpha^{\gamma}\thin^{\dot{\alpha}}\dot{\alpha}\thin$. 

By no $\gamma\delta\cdots$, we do not have $\thin^{\beta}\alpha^{\gamma}\thin^{\delta}\beta^{\alpha}\thin$. By no $\delta^2\cdots$, we do not have $\thin^{\alpha}\beta^{\delta}\thin^{\delta}\beta^{\alpha}\thin$. Combined with no $\thin^{\beta}\alpha^{\gamma}\thin^{\gamma}\alpha^{\beta}\thin$, we do not have $\thin^{\dot{\alpha}}\dot{\alpha}\thin\dot{\alpha}^{\dot{\alpha}}\thin$. Therefore $\dot{\alpha}^4=\thin^{\dot{\alpha}}\dot{\alpha}\thin^{\dot{\alpha}}\dot{\alpha}\thin^{\dot{\alpha}}\dot{\alpha}\thin^{\dot{\alpha}}\dot{\alpha}\thin$. Then by no $\thin^{\beta}\alpha^{\gamma}\thin^{\dot{\alpha}}\dot{\alpha}\thin$, we get $\thin^{\dot{\alpha}}\dot{\alpha}\thin^{\dot{\alpha}}\dot{\alpha}\thin=\thin^{\alpha}\beta^{\delta}\thin^{\dot{\alpha}}\dot{\alpha}\thin$. This implies $\dot{\alpha}^4=\thin^{\alpha}\beta^{\delta}\thin^{\alpha}\beta^{\delta}\thin^{\alpha}\beta^{\delta}\thin^{\alpha}\beta^{\delta}\thin=\beta^4$.

In Figure \ref{4adeC}, the AAD $\thin^{\alpha}\beta^{\delta}\thin^{\alpha}\beta^{\delta}\thin$ determines $T_1,T_2$. The two tiles determine their companions $T_3,T_4$. Then $\delta_3\epsilon_1\cdots=\dot{\alpha}\delta\epsilon$ and $\gamma_1\gamma_4=\dot{\alpha}\gamma^2$ determine two $\dot{\alpha}$ in $T_5$. Then $T_5$ has two angle arrangements as indicated in the two pictures. In each case, $T_5$ further determines its companion $T_6$.

We apply the argument to all four $\thin^{\alpha}\beta^{\delta}\thin^{\alpha}\beta^{\delta}\thin$ in $\beta^4$, and get two earth map tilings with four timezones. The first timezone is given by $T_1,T_4,T_5,T_6$ in the first of Figure \ref{4adeC}, and the second timezone is given by $T_1,T_3,T_5,T_6$ in the second of Figure \ref{4adeC}. Both earth map tilings can be divided into the upper half consisting of tiles similar to $T_1,T_2,T_3,T_4$, and the lower half consisting of tiles similar to $T_5,T_6$. The two tilings have the same upper half, and the lower half are related by a flip. The first is an $f=16$ case of the earth map tiling $E_{\pentagon}2$ in Figure \ref{emt}. The second is the sporadic tiling $S_{16\pentagon}$ in Figure \ref{emt2mod7}.

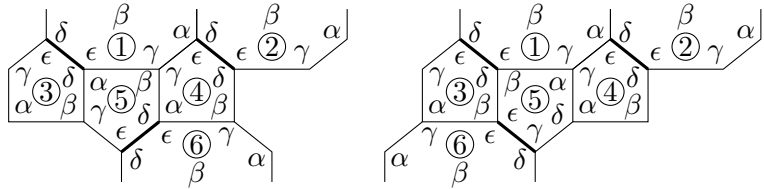
\begin{figure}[htp]
\centering
\begin{tikzpicture}[>=latex,scale=1]


\foreach \a in {0,1}
\foreach \b in {0,1}
{
\begin{scope}[xshift=5.5*\b cm, xshift=2*\a cm]

\draw
	(0,1.5) -- (0,1.1) -- (0.5,0.7) -- (1.5,0.7) -- (2,1.1) -- (2,1.5)
	(0,1.1) -- (-0.5,0.7) -- (-0.5,0) -- (0.5,0) -- (0.5,0.7);

\draw[line width=1.2]
	(0,1.1) -- (0.5,0.7);

\node at (1,1.4) {\small $\beta$};
\node at (1.8,1.2) {\small $\alpha$}; 
\node at (1.4,0.9) {\small $\gamma$};
\node at (0.6,0.9) {\small $\epsilon$};
\node at (0.2,1.2) {\small $\delta$};

\node at (-0.3,0.6) {\small $\gamma$};
\node at (0.3,0.6) {\small $\delta$};
\node at (-0.3,0.2) {\small $\alpha$}; 
\node at (0.3,0.2) {\small $\beta$};	
\node at (0,0.85) {\small $\epsilon$};	
	
\end{scope}
}


\draw
	(0.5,0) -- (1,-0.4) -- (1.5,0)
	(1,-0.8) -- (1,-0.4) 
	(2.5,0) -- (3,-0.4) -- (3,-0.8);
	
\draw[line width=1.2]
	(1,-0.4) -- (1.5,0);
	
\node at (1.3,0.5) {\small $\beta$}; 
\node at (0.7,0.5) {\small $\alpha$};
\node at (0.7,0.1) {\small $\gamma$};	
\node at (1,-0.15) {\small $\epsilon$};	
\node at (1.3,0.1) {\small $\delta$};

\node at (2,-0.7) {\small $\beta$};
\node at (2.8,-0.5) {\small $\alpha$}; 
\node at (2.4,-0.2) {\small $\gamma$};
\node at (1.6,-0.2) {\small $\epsilon$};
\node at (1.2,-0.5) {\small $\delta$}; 

\node[inner sep=0.5,draw,shape=circle] at (1,1) {\small 1};
\node[inner sep=0.5,draw,shape=circle] at (3,1) {\small 2};
\node[inner sep=0.5,draw,shape=circle] at (0,0.4) {\small 3};
\node[inner sep=0.5,draw,shape=circle] at (2,0.4) {\small 4};
\node[inner sep=0.5,draw,shape=circle] at (1,0.3) {\small 5};
\node[inner sep=0.5,draw,shape=circle] at (2,-0.3) {\small 6};


\begin{scope}[xshift=5.5 cm]

\draw
	(0.5,0) -- (1,-0.4) -- (1.5,0)
	(1,-0.8) -- (1,-0.4) 
	(-0.5,0) -- (-1,-0.4) -- (-1,-0.8);
	
\draw[line width=1.2]
	(1,-0.4) -- (0.5,0);
	
\node at (1.3,0.5) {\small $\alpha$}; 
\node at (0.7,0.5) {\small $\beta$};
\node at (0.7,0.1) {\small $\epsilon$};	
\node at (1,-0.15) {\small $\gamma$};	
\node at (1.3,0.1) {\small $\delta$};

\node at (0,-0.7) {\small $\beta$};
\node at (0.8,-0.5) {\small $\delta$}; 
\node at (0.4,-0.2) {\small $\epsilon$};
\node at (-0.4,-0.2) {\small $\gamma$};
\node at (-0.8,-0.5) {\small $\alpha$}; 

\node[inner sep=0.5,draw,shape=circle] at (1,1) {\small 1};
\node[inner sep=0.5,draw,shape=circle] at (3,1) {\small 2};
\node[inner sep=0.5,draw,shape=circle] at (0,0.4) {\small 3};
\node[inner sep=0.5,draw,shape=circle] at (2,0.4) {\small 4};
\node[inner sep=0.5,draw,shape=circle] at (1,0.3) {\small 5};
\node[inner sep=0.5,draw,shape=circle] at (0,-0.3) {\small 6};

\end{scope}
	
\end{tikzpicture}
\caption{Proposition \ref{4ade}: $\alpha\delta\epsilon,\alpha\gamma^2$ are vertices.}
\label{4adeC}
\end{figure}

Combinatorially, the argument can also be applied to $\beta^5$, for the case $f=20$, and produces earth map tilings with five timezones. However, the pentagon does not exist for geometrical reason.

\medskip

\noindent{\em Geometry of Pentagon}

\medskip

For $f=16$, we get
\[
\alpha=\tfrac{1}{2}\pi,\,
\gamma=\tfrac{3}{4}\pi,\,
\delta+\epsilon=\tfrac{3}{2}\pi.
\]
Then by \eqref{coolsaet_eq2} and \eqref{coolsaet_eq3}, we get 
\[
\cos a=\cot\tfrac{1}{2}\alpha\cot\delta
=\cot\tfrac{1}{2}\gamma\cot\epsilon.
\]
Substituting the angle values, we find the solution satisfying $\epsilon<\gamma=\tfrac{3}{4}\pi$ is given by
\[
\cos\epsilon=\frac{1}{\sqrt[4]{2}},\quad
\epsilon=0.1820\pi.
\]
Then we further find $a,b$ are given by
\begin{align*}
\cos a & ={\textstyle \sqrt{\sqrt{2}-1} }, &
a &= 0.2774\pi, \\
\cos b & ={\textstyle \sqrt{2}\sqrt{\sqrt{2}-1} }, &
b & = 0.1359\pi. 
\end{align*}
The base angles of $\triangle ABD$ and $\triangle ACE$ are $\epsilon,\delta-\pi<\alpha$. Then by Lemma \ref{geometry8}, we know the pentagon is simple.

\subsubsection*{Case. $\alpha\delta\epsilon,\alpha^2\gamma$ are vertices.}

The angle sums of $\alpha\delta\epsilon,\alpha^2\gamma$ and the angle sum for pentagon imply
\begin{equation}\label{4ade_eq1}
\alpha=(1-\tfrac{4}{f})\pi,\,
\gamma=\tfrac{8}{f}\pi,\,
\delta+\epsilon=(1+\tfrac{4}{f})\pi.
\end{equation}
We have $\alpha>\gamma$. By Lemma \ref{geometry1}, this implies $\delta<\epsilon$. 

We have $\alpha+\gamma=\delta+\epsilon$. If $\alpha>\epsilon$, then this implies $\gamma<\delta$. By $\alpha<\pi$, we know the pentagon is strictly convex. Then by Lemma \ref{geometry4}, we get a contradiction. Therefore we have $\alpha<\epsilon$. 

By $\delta+\epsilon>\pi$ and $\delta<\epsilon$, we know $R(\epsilon^2)$ has no $\delta,\epsilon$. By $\alpha^2\gamma$ and $\alpha<\epsilon$, we get $R(\epsilon^2)<\gamma<\alpha$. Therefore $\epsilon^2\cdots$ is not a vertex. Then by $\alpha\delta\epsilon$ and the balance lemma, we know $\alpha\delta\epsilon,\gamma^k\delta\epsilon$ are all the $b$-vertices. By $\alpha^2\gamma$ and $\alpha>\gamma$, we know $\alpha^2\gamma,\alpha\gamma^k,\gamma^k$ are all the $\hat{b}$-vertices. Then we get all the vertices ($\dot{\alpha}$ is $\alpha$ or $\beta$)
\[
\text{AVC}
=\{\dot{\alpha}\delta\epsilon,\dot{\alpha}^2\gamma,\dot{\alpha}\gamma^k,\gamma^k,\gamma^k\delta\epsilon\}.
\]
This implies all companion pairs are twisted, and a tile determines its companion.

The AAD $\thin^{\alpha}\gamma^{\epsilon}\thin^{\alpha}\gamma^{\epsilon}\thin$ determines $T_1,T_2$ in Figure \ref{4adeD1}. Then $T_1,T_2$ determine their companions $T_3,T_4$. Then $\beta_1\beta_4\cdots=\dot{\alpha}^2\gamma$ and $\delta_1\epsilon_3\cdots=\dot{\alpha}\delta\epsilon,\gamma^k\delta\epsilon$ determine $T_5$. Then $T_5$ determines its companion $T_6$.

By no $\epsilon^2\cdots$, the AAD of $\gamma^k$ is $\thin^{\alpha}\gamma^{\epsilon}\thin^{\alpha}\gamma^{\epsilon}\thin\cdots\thin^{\alpha}\gamma^{\epsilon}\thin$. Applying the argument above to all $\thin^{\alpha}\gamma^{\epsilon}\thin^{\alpha}\gamma^{\epsilon}\thin$ in the vertex, we get the earth map tiling with $T_1,T_3,T_5,T_6$ as the timezone. We may also choose $T_1,T_4,T_5,T_6$ as the timezone. The second choice shows the tiling is the second earth map tiling $E_{\pentagon}2$ in Figure \ref{emt}.

\begin{figure}[htp]
\centering
\begin{tikzpicture}[>=latex,scale=1]

\draw[gray!50, line width=3]
	(0,1.5) -- (0,1.1) -- (0.5,0.7) -- (1.5,0.7) -- (1.5,0) -- (1,-0.4) -- (1,-0.8);
	
\draw[gray!50, xshift=8 cm, line width=3]
	(0,1.5) -- (0,1.1) -- (0.5,0.7) -- (0.5,0) -- (-0.5,0) -- (-1,-0.4) -- (-1,-0.8);
		
\foreach \a in {0,...,3}
{
\begin{scope}[xshift=2*\a cm]

\draw
	(0,1.5) -- (0,1.1) -- (0.5,0.7) -- (1.5,0.7) -- (2,1.1) -- (2,1.5)
	(1,-0.8) -- (1,-0.4) -- (0.5,0) -- (-0.5,0) -- (-1,-0.4) -- (-1,-0.8)
	(2,1.1) -- (2.5,0.7) -- (2.5,0) -- (1.5,0) -- (1,-0.4)
	(1.5,0) -- (1.5,0.7)
	(0.5,0) -- (0.5,0.7);
		
\draw[line width=1.2]
	(2,1.1) -- (1.5,0.7)
	(1,-0.4) -- (0.5,0);

\node at (1,1.4) {\small $\gamma$};
\node at (1.8,1.2) {\small $\epsilon$}; 
\node at (1.4,0.9) {\small $\delta$};
\node at (0.6,0.9) {\small $\beta$};
\node at (0.2,1.2) {\small $\alpha$};
	
\node at (1.7,0.6) {\small $\epsilon$};
\node at (2.3,0.6) {\small $\beta$};
\node at (1.7,0.2) {\small $\gamma$}; 
\node at (2.3,0.2) {\small $\alpha$};	
\node at (2,0.85) {\small $\delta$};	

\node at (1.3,0.5) {\small $\alpha$}; 
\node at (0.7,0.5) {\small $\gamma$};
\node at (0.7,0.1) {\small $\epsilon$};	
\node at (1,-0.15) {\small $\delta$};	
\node at (1.3,0.1) {\small $\beta$};

\node at (0,-0.7) {\small $\gamma$};
\node at (0.8,-0.5) {\small $\epsilon$}; 
\node at (0.4,-0.2) {\small $\delta$};
\node at (-0.4,-0.2) {\small $\beta$};
\node at (-0.8,-0.5) {\small $\alpha$}; 

\end{scope}
}

\node at (0,1.7) {\footnotesize $l_2$}; 
\node at (8,1.7) {\footnotesize $l_1$};

\node[inner sep=0.5,draw,shape=circle] at (7,1) {\small 1};
\node[inner sep=0.5,draw,shape=circle] at (5,1) {\small 2};
\node[inner sep=0.5,draw,shape=circle] at (8,0.4) {\small 3};
\node[inner sep=0.5,draw,shape=circle] at (6,0.4) {\small 4};
\node[inner sep=0.5,draw,shape=circle] at (7,0.3) {\small 5};
\node[inner sep=0.5,draw,shape=circle] at (6,-0.3) {\small 6};
\node[inner sep=0.5,draw,shape=circle] at (1,1) {\small 7};
\node[inner sep=0.5,draw,shape=circle] at (2,0.4) {\small 8};
\node[inner sep=0.5,draw,shape=circle] at (1,0.3) {\small 9};
\node[inner sep=0,draw,shape=circle] at (0,-0.3) {\footnotesize 10};
	
\end{tikzpicture}
\caption{Proposition \ref{4ade}: ${\mc P}^k_{k-1}$ and ${\mc P}^k_k$.}
\label{4adeD1}
\end{figure}
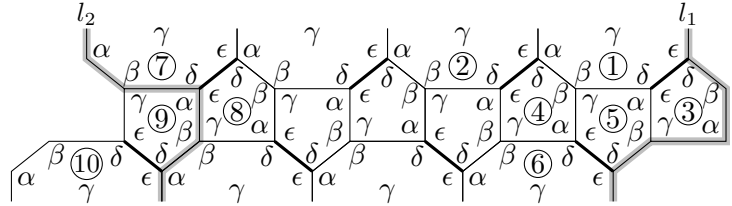

In general, if we have the the AAD $(\thin^{\alpha}\gamma^{\epsilon}\thin)^k=\thin^{\alpha}\gamma^{\epsilon}\thin\cdots\thin^{\alpha}\gamma^{\epsilon}\thin$ as part of a vertex, then we get a partial earth map tiling ${\mc P}^k_{k-1}$, with $T_1,T_3,T_5,T_6$ as the timezone, and with $(\thin^{\alpha}\gamma^{\epsilon}\thin)^{k-1}$ at the other end. The boundary of ${\mc P}^k_{k-1}$ is given by the shaded edges $l_1,l_2$ in Figure \ref{4adeD1} (the picture has $k=4$). We note that ${\mc P}^k_{k-1}$ actually consists of $k-1$ full timezones, and one half timezone consisting of $T_7,T_8$ (similar to $T_2,T_4$). 

The full earth map tiling is obtained by glueing two copies of ${\mc P}^{q+1}_q$ together, as in Figure \ref{emt2mod2}. The inner hemisphere in Figure \ref{emt2mod2} is the vertical flip of the partial earth map tiling in Figure \ref{4adeD1}. Each ${\mc P}^{q+1}_q$ tiles a hemisphere, and is actually a half earth map tiling. By $\alpha=\beta$, the angle values along the boundary between the two hemispheres is given by the first of Figure \ref{emt2mod3}.

In the subsequent tiling construction, we also denote by ${\mc P}^k_k$ the partial earth map tiling by adding two tiles $T_9,T_{10}$ to ${\mc P}^k_{k-1}$. This makes the half timezone $T_7,T_8$ in ${\mc P}^k_{k-1}$ into one full timezone. Therefore ${\mc P}^k_k$ consists of $k$ full timezones, with $(\thin^{\alpha}\gamma^{\epsilon}\thin)^k$ at both ends. 

We already proved that, if $\gamma^k$ is a vertex, then the tiling is the earth map tiling. In fact, we have $k=\frac{f}{4}$ in $\gamma^k$. Next we assume $\gamma^k$ is not a vertex. Then we get all the vertices
\[
\text{AVC}
=\{\dot{\alpha}\delta\epsilon,\dot{\alpha}^2\gamma,\dot{\alpha}\gamma^{q+1},\gamma^q\delta\epsilon\},
\quad q=\tfrac{f-4}{8}.
\]

We have consecutive $\gamma^k=\gamma\cdots\gamma$ at $\alpha\gamma^{q+1},\gamma^q\delta\epsilon$. By no $\epsilon^2\cdots$, the AAD of the consecutive $\gamma^k$ has two possibilities 
\begin{align*}
(\thin^{\alpha}\gamma^{\epsilon}\thin)^k
&=\thin^{\alpha}\gamma^{\epsilon}\thin\cdots\thin^{\alpha}\gamma^{\epsilon}\thin
=\thin^{\alpha}\gamma\cdots\gamma^{\epsilon}\thin, \\
(\thin^{\epsilon}\gamma^{\alpha}\thin)^s(\thin^{\alpha}\gamma^{\epsilon}\thin)^t
&=\thin^{\epsilon}\gamma^{\alpha}\thin\cdots\thin^{\epsilon}\gamma^{\alpha}\thin^{\alpha}\gamma^{\epsilon}\thin\cdots\thin^{\alpha}\gamma^{\epsilon}\thin
=\thin^{\epsilon}\gamma\cdots\gamma^{\epsilon}\thin.
\end{align*}
Here we introduce the abbreviated notations on the right side. Then by no $\delta\thin\epsilon\cdots$,  the AAD of $\dot{\alpha}\gamma^{q+1}$ has the following possibilities 
\begin{equation}\label{4ade_eqa}
\thin^{\beta}\alpha^{\gamma}\thin^{\alpha}\gamma\cdots\gamma^{\epsilon}\thin,\;
\thin^{\beta}\alpha^{\gamma}\thin^{\epsilon}\gamma\cdots\gamma^{\alpha}\thin,\;
\thin^{\beta}\alpha^{\gamma}\thin^{\epsilon}\gamma\cdots\gamma^{\epsilon}\thin,\;
\thin^{\alpha}\beta^{\delta}\thin^{\alpha}\gamma\cdots\gamma^{\epsilon}\thin.
\end{equation}
and the AAD of $\gamma^q\delta\epsilon$ has the following possibilities
\begin{equation}\label{4ade_eqde}
\thick^{\epsilon}\delta^{\beta}\thin^{\alpha}\gamma\cdots\gamma^{\epsilon}\thin^{\gamma}\epsilon^{\delta}\thick,\;
\thick^{\epsilon}\delta^{\beta}\thin^{\epsilon}\gamma\cdots\gamma^{\alpha}\thin^{\gamma}\epsilon^{\delta}\thick,\;
\thick^{\epsilon}\delta^{\beta}\thin^{\epsilon}\gamma\cdots\gamma^{\epsilon}\thin^{\gamma}\epsilon^{\delta}\thick.
\end{equation}

We know the AAD $\thin^{\alpha}\gamma\cdots\gamma^{\epsilon}\thin$ induces the partial earth map tiling ${\mc P}^k_{k-1}$ in the first of Figure \ref{4adeD1}. The AAD $\thin^{\epsilon}\gamma\cdots\gamma^{\epsilon}\thin$ is $\thin^{\epsilon}\gamma\cdots\gamma^{\alpha}\thin^{\alpha}\gamma\cdots\gamma^{\epsilon}\thin$, which is two AADs $\thin^{\alpha}\gamma\cdots\gamma^{\epsilon}\thin$ glued at $\thin^{\epsilon}\gamma^{\alpha}\thin^{\alpha}\gamma^{\epsilon}\thin$. The AAD $\thin^{\epsilon}\gamma^{\alpha}\thin^{\alpha}\gamma^{\epsilon}\thin$ determines $T_1,T_2$ in Figure \ref{4adeD2}. Then $T_1,T_2$ determine their companions $T_3,T_4$, and $\alpha_1\alpha_2\cdots=\dot{\alpha}^2\gamma$ gives $\gamma_5$. Up to the horizontal flip symmetry, we may assume $T_5$ is arranged as indicated. Then $T_5$ determines its companion $T_6$. Moreover, $\alpha_5\beta_2\cdots=\dot{\alpha}^2\gamma$ and $\delta_2\epsilon_4\cdots=\dot{\alpha}\delta\epsilon,\gamma^q\delta\epsilon$ determine $T_7$. Then $T_7$ determines its companion $T_8$.

The tiles $T_5,T_6,T_7,T_8$ form the four tile patch ${\mc A}$ in Figure \ref{emt2mod8}. If $T_5$ is arranged in another way, which is the horizontal flip of $T_5$ as indicated in the picture, then we get the horizontal flip of the patch ${\mc A}$. The flip is also described in Figure \ref{emt2mod8}.

The two AADs $\thin^{\alpha}\gamma\cdots\gamma^{\epsilon}\thin$ in $\thin^{\epsilon}\gamma\cdots\gamma^{\epsilon}\thin$ induce two partial earth map tilings ${\mc P}^s_{s-1}$ and ${\mc P}^t_{t-1}$. Their meeting place $\thin^{\epsilon}\gamma^{\alpha}\thin^{\alpha}\gamma^{\epsilon}\thin$ is filled by a copy of ${\mc A}$, and we get a {\em partial tiling} ${\mc P}^s_{s-1}\cup {\mc A}\cup {\mc P}^t_{t-1}$ in the first of Figure \ref{4adeD2} (the picture has $s=2$ and $t=2$), with $\gamma^k,\gamma^{k-1}$ ($k=s+t$) at the two ends, and the boundary given by the shaded edges $l_1,l_2$. 

\begin{figure}[htp]
\centering
\begin{tikzpicture}[>=latex,scale=1]


\foreach \a in {1,-1}
{
\draw[xshift=4*\a cm, xscale=\a, gray!50, line width=3]
	(0,1.5) -- (0,1.1) -- (0.5,0.7) -- (0.5,0) -- (-0.5,0) -- (-1,-0.4) -- (-1,-0.8);
	
\draw[gray!50, xscale=\a, line width=3]	
	(0,1.1) -- (0.5,0.7) -- (1.5,0.7) -- (1.5,0) -- (1,-0.4) -- (1,-0.8); 

\foreach \b in {1,2}
\draw[shift={(-3cm, -2.6*\b cm)}, xscale=\a, gray!50, line width=3]
	(1,1.5) -- (1,1.1) -- (1.5,0.7) -- (1.5,0) -- (0.5,0) -- (0,-0.4);

}


\foreach \b in {0,1,2}
{
\begin{scope}[yshift=-2.6*\b cm]

\foreach \a in {-1,...,1}
\draw[xshift=2*\a cm]
	(-1,-0.8) -- (-1,-0.4) -- (-0.5,0) -- (0.5,0) -- (1,-0.4) -- (1,-0.8);

\foreach \a in {-2,...,1}
\draw[xshift=2*\a cm]
	(0,1.5) -- (0,1.1) -- (0.5,0.7) -- (1.5,0.7) -- (2,1.1) -- (2,1.5)
	(0.5,0) -- (1,-0.4) -- (1.5,0);

\foreach \a in {-2,...,2}
\draw[xshift=2*\a cm]
	(0.5,0) -- (0.5,0.7) -- (0,1.1) -- (-0.5,0.7) -- (-0.5,0) -- (0.5,0);

\end{scope}
}


\foreach \a/\b in {0/0,1/0,-1/1,0/1,1/1,-2/2,-1/2,0/2,1/2}
{
\begin{scope}[yshift=-2.6*\b cm, xshift=2*\a cm]

\draw[line width=1.2]
	(2,1.1) -- (1.5,0.7);
	
\node at (1,1.4) {\small $\gamma$};
\node at (1.8,1.2) {\small $\epsilon$}; 
\node at (1.4,0.9) {\small $\delta$};
\node at (0.6,0.9) {\small $\beta$};
\node at (0.2,1.2) {\small $\alpha$};

\node at (1.7,0.6) {\small $\epsilon$};
\node at (2.3,0.6) {\small $\beta$};
\node at (1.7,0.2) {\small $\gamma$}; 
\node at (2.3,0.2) {\small $\alpha$};	
\node at (2,0.85) {\small $\delta$};

\end{scope}
}


\foreach \a/\b in {0/0,-1/0,-1/1}
{
\begin{scope}[yshift=-2.6*\b cm, xshift=2*\a cm]

\draw[line width=1.2]
	(-2,1.1) -- (-1.5,0.7);
	
\node at (-1,1.4) {\small $\gamma$};
\node at (-1.8,1.2) {\small $\epsilon$}; 
\node at (-1.4,0.9) {\small $\delta$};
\node at (-0.6,0.9) {\small $\beta$};
\node at (-0.2,1.2) {\small $\alpha$};

\node at (-1.7,0.6) {\small $\epsilon$};
\node at (-2.3,0.6) {\small $\beta$};
\node at (-1.7,0.2) {\small $\gamma$}; 
\node at (-2.3,0.2) {\small $\alpha$};	
\node at (-2,0.85) {\small $\delta$};

\end{scope}
}

	
\foreach \a/\b in {0/0,1/0,-1/1,0/1,1/1,-1/2,0/2,1/2}
{
\begin{scope}[yshift=-2.6*\b cm, xshift=2*\a cm]

\draw[line width=1.2]
	(1,-0.4) -- (0.5,0);
	
\node at (1.3,0.5) {\small $\alpha$}; 
\node at (0.7,0.5) {\small $\gamma$};
\node at (0.7,0.1) {\small $\epsilon$};	
\node at (1,-0.15) {\small $\delta$};	
\node at (1.3,0.1) {\small $\beta$};

\node at (0,-0.7) {\small $\gamma$};
\node at (0.8,-0.5) {\small $\epsilon$}; 
\node at (0.4,-0.2) {\small $\delta$};
\node at (-0.4,-0.2) {\small $\beta$};
\node at (-0.8,-0.5) {\small $\alpha$}; 

\end{scope}
}


\draw[line width=1.2]
	(-2.5,0) -- (-3,-0.4);
	
\node at (-2,-0.7) {\small $\gamma$};
\node at (-1.2,-0.5) {\small $\alpha$}; 
\node at (-1.6,-0.2) {\small $\beta$};
\node at (-2.4,-0.2) {\small $\delta$};
\node at (-2.8,-0.5) {\small $\epsilon$}; 

\node at (-3.3,0.5) {\small $\alpha$}; 
\node at (-2.7,0.5) {\small $\gamma$};
\node at (-2.7,0.1) {\small $\epsilon$};	
\node at (-3,-0.15) {\small $\delta$};	
\node at (-3.3,0.1) {\small $\beta$};


\foreach \a/\b in {0/0,-1/1}
{
\begin{scope}[yshift=-2.6*\b cm, xshift=2*\a cm]

\draw[line width=1.2]
	(-0.5,0) -- (-0.5,0.7);

\node at (-0.3,0.6) {\small $\epsilon$};
\node at (0.3,0.6) {\small $\alpha$};
\node at (-0.3,0.2) {\small $\delta$}; 
\node at (0.3,0.2) {\small $\beta$};	
\node at (0,0.85) {\small $\gamma$};		

\node at (-0.7,0.5) {\small $\delta$}; 
\node at (-1.3,0.5) {\small $\beta$};
\node at (-1.3,0.1) {\small $\alpha$};	
\node at (-1,-0.15) {\small $\gamma$};	
\node at (-0.7,0.1) {\small $\epsilon$};

\node[inner sep=0.5,draw,shape=circle] at (-1,1) {\small 1};
\node[inner sep=0.5,draw,shape=circle] at (1,1) {\small 2};
\node[inner sep=0.5,draw,shape=circle] at (-2,0.4) {\small 3};
\node[inner sep=0.5,draw,shape=circle] at (2,0.4) {\small 4};
\node[inner sep=0.5,draw,shape=circle] at (0,0.4) {\small 5};
\node[inner sep=0.5,draw,shape=circle] at (-1,0.3) {\small 6};
\node[inner sep=0.5,draw,shape=circle] at (1,0.3) {\small 7};
\node[inner sep=0.5,draw,shape=circle] at (0,-0.3) {\small 8};

\end{scope}
}


\begin{scope}[shift={(-4cm,-5.2 cm)}]

\draw[line width=1.2]
	(0.5,0) -- (0.5,0.7);

\node at (-0.3,0.6) {\small $\alpha$};
\node at (0.3,0.6) {\small $\epsilon$};
\node at (-0.3,0.2) {\small $\beta$}; 
\node at (0.3,0.2) {\small $\delta$};	
\node at (0,0.85) {\small $\gamma$};		

\node at (0.7,0.5) {\small $\delta$}; 
\node at (1.3,0.5) {\small $\beta$};
\node at (1.3,0.1) {\small $\alpha$};	
\node at (1,-0.15) {\small $\gamma$};	
\node at (0.7,0.1) {\small $\epsilon$};

\node[inner sep=0.5,draw,shape=circle] at (0,0.4) {\small 3};
\node[inner sep=0.5,draw,shape=circle] at (1,0.3) {\small 6};

\end{scope}

\node at (-4,1.7) {\footnotesize $l_2$}; 
\node at (4,1.7) {\footnotesize $l_1$};

\end{tikzpicture}
\caption{Proposition \ref{4ade}: $\thin^{\epsilon}\gamma\cdots\gamma^{\epsilon}\thin$ gives ${\mc P}^s_{s-1}\cup {\mc A}\cup {\mc P}^t_{t-1}$, and is modified to ${\mc A}\cup {\mc P}^k_k$ with $\thin^{\alpha}\gamma\cdots\gamma^{\epsilon}\thin$.}
\label{4adeD2}
\end{figure}
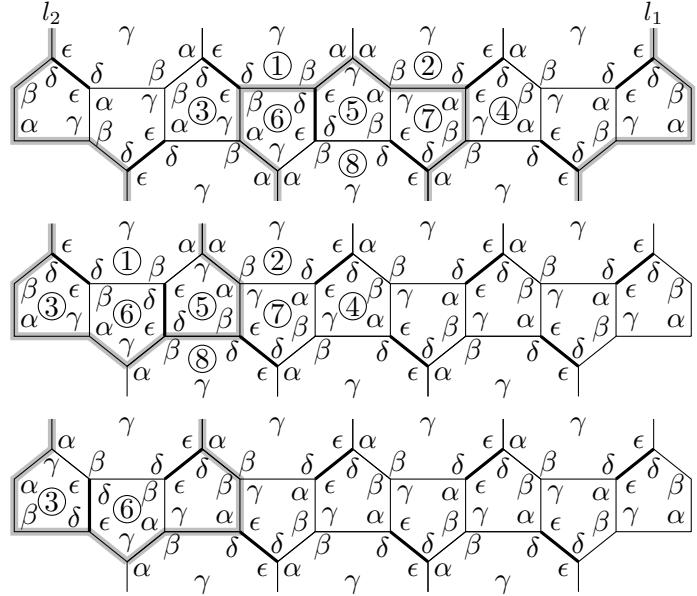

The partial tiling ${\mc P}^s_{s-1}\cup {\mc A}\cup {\mc P}^t_{t-1}$ can be modified without changing the angle values along the boundaries $l_1,l_2$. We may flip ${\mc A}$, or change $s,t$ as along as $s+t=k$. The change of $s,t$ shifts the location of ${\mc A}$ inside the partial tiling. The second of Figure \ref{4adeD2} shows that the patch ${\mc A}$ in the first is shifted to one side of the partial tiling. 

We note that $T_1,T_3,T_5,T_6$ in the first of Figure \ref{4adeD2} also form an ``upper'' patch ${\mc A}$. The shift of the original ``lower'' ${\mc A}$ induces the corresponding shift of the upper ${\mc A}$. In the second of Figure \ref{4adeD2}, the upper ${\mc A}$ (the boundary is outlined in the second picture) is also shifted to one side of the partial earth map tiling. This is ${\mc A}\cup {\mc P}^k_k$, with ${\mc A}$ being the upper one. 

If we further flip the patch ${\mc A}$ in ${\mc A}\cup {\mc P}^k_k$, then we get the third of Figure \ref{4adeD2}. This changes the original $\thin^{\epsilon}\gamma\cdots\gamma^{\epsilon}\thin$ to $\thin^{\alpha}\gamma\cdots\gamma^{\epsilon}\thin$. Of course we may also shift the upper ${\mc A}$ in the second of Figure \ref{4adeD2} to the other side of the partial tiling. Then we get ${\mc P}^k_k\cup {\mc A}$. This changes the original $\thin^{\epsilon}\gamma\cdots\gamma^{\epsilon}\thin$ to $\thin^{\epsilon}\gamma\cdots\gamma^{\alpha}\thin$.

The three partial tilings in Figure \ref{4adeD2} have the same boundary and angle values along the boundary. This means that, if a vertex contains $\thin^{\epsilon}\gamma\cdots\gamma^{\epsilon}\thin$, then we may modify the tiling by the flip and shift of the upper ${\mc A}$ in Figure \ref{4adeD2}, such that $\thin^{\epsilon}\gamma\cdots\gamma^{\epsilon}\thin$ is changed to $\thin^{\alpha}\gamma\cdots\gamma^{\epsilon}\thin$ or $\thin^{\epsilon}\gamma\cdots\gamma^{\alpha}\thin$. Therefore the third AAD in \eqref{4ade_eqa} gives a relation between the first and second AADs via the flip and shift of ${\mc A}$
\begin{equation}\label{4ade_eqa_shift}
\thin^{\beta}\alpha^{\gamma}\thin^{\alpha}\gamma\cdots\gamma^{\epsilon}\thin
\longleftarrow
\thin^{\beta}\alpha^{\gamma}\thin^{\epsilon}\gamma\cdots\gamma^{\epsilon}\thin
\longrightarrow
\thin^{\beta}\alpha^{\gamma}\thin^{\epsilon}\gamma\cdots\gamma^{\alpha}\thin.
\end{equation}
The third in \eqref{4ade_eqde} also gives the similar relation
\begin{equation}\label{4ade_eqde_shift}
\thick^{\epsilon}\delta^{\beta}\thin^{\alpha}\gamma\cdots\gamma^{\epsilon}\thin^{\gamma}\epsilon^{\delta}\thick
\longleftarrow
\thick^{\epsilon}\delta^{\beta}\thin^{\epsilon}\gamma\cdots\gamma^{\epsilon}\thin^{\gamma}\epsilon^{\delta}\thick
\longrightarrow
\thick^{\epsilon}\delta^{\beta}\thin^{\epsilon}\gamma\cdots\gamma^{\alpha}\thin^{\gamma}\epsilon^{\delta}\thick.
\end{equation}

Conversely, a vertex containing $\thin^{\alpha}\gamma\cdots\gamma^{\epsilon}\thin$ fits into the relation above if it induces a partial tiling with ${\mc A}$ on a side, like the third of Figure \ref{4adeD2}. We note that the third of Figure \ref{4adeD2} can also be regarded as obtained by adding $T_3,T_6$ to ${\mc P}^{k+1}_k$.

\subsubsection*{Subcase. $\dot{\alpha}\gamma^{q+1}$ is a vertex.}

Suppose the first AAD $\thin^{\beta}\alpha^{\gamma}\thin^{\alpha}\gamma\cdots\gamma^{\epsilon}\thin$ of \eqref{4ade_eqa} is a vertex. The $\thin^{\alpha}\gamma\cdots\gamma^{\epsilon}\thin$ part of the vertex gives ${\mc P}^{q+1}_q$ (indicated by $\gamma^{q+1}$ and $\gamma^q$ in the first picture) on the left of $l_1$ and the right of $l_2$. The $\thin^{\beta}\alpha^{\gamma}\thin$ part of the vertex determines $T_1$. Then $T_1$ determines its companion $T_2$. Then $\gamma_2\cdots=\alpha\gamma\cdots=\dot{\alpha}^2\gamma,\dot{\alpha}\gamma^{q+1}$. 

The first of Figure \ref{4adeF} is the case $\gamma_2\cdots=\dot{\alpha}^2\gamma$. The vertex $\thin^{\delta}\beta^{\alpha}\thin^{\epsilon}\gamma^{\alpha}\thin\cdots$ on the left of $l_1$ is $\dot{\alpha}\gamma\cdots=\dot{\alpha}^2\gamma,\dot{\alpha}\gamma^{q+1}$. If the vertex is $\dot{\alpha}\gamma^{q+1}$, then by $\gamma_2\cdots=\dot{\alpha}^2\gamma$, the AAD of the vertex is $\thin^{\delta}\beta^{\alpha}\thin^{\epsilon}\gamma^{\alpha}\thin^{\alpha}\gamma^{\epsilon}\thin\cdots$, contradicting the possible AADs \eqref{4ade_eqa} of $\dot{\alpha}\gamma^{q+1}$. Therefore the vertex is $\dot{\alpha}^2\gamma$, and we get two adjacent $\dot{\alpha}$ indicated by the arrows. The vertex $\delta\epsilon\cdots$ on the lower left of $l_1$ is $\dot{\alpha}\delta\epsilon,\gamma^q\delta\epsilon$. By no three adjacent $\dot{\alpha}$, we know the vertex is $\gamma^q\delta\epsilon$. This determines which of the two adjacent $\dot{\alpha}$ indicated by the arrows is $\alpha$ or $\beta$. Then by the possible AADs \eqref{4ade_eqde} of $\gamma^q\delta\epsilon$, the vertex is $\thick^{\epsilon}\delta^{\beta}\thin^{\alpha}\gamma^{\epsilon}\thin\cdots=\thick^{\epsilon}\delta^{\beta}\thin^{\alpha}\gamma\cdots\gamma^{\epsilon}\thin^{\gamma}\epsilon^{\delta}\thick$. This gives $\epsilon$ at the bottom of the picture. Then the bottom vertex is $\gamma^q\epsilon\cdots=\gamma^q\delta\epsilon$. This gives the bottom $\delta$. Then the vertex $\alpha\cdots$ on the lower right of $l_2$ is $\alpha\beta\cdots=\dot{\alpha}^2\gamma$. This gives $\gamma_3$, and implies $\beta\gamma\cdots$ on the right of $l_2$ is not $\dot{\alpha}\gamma^{q+1}$, and must be $\dot{\alpha}^2\gamma$. By $\gamma_3$, this determines $T_3$. Then the vertex $\beta\cdots$ on the right of $l_2$ is $\beta\delta\cdots=\dot{\alpha}\delta\epsilon$. This implies $\gamma_1\cdots=\dot{\alpha}\gamma^2\cdots=\dot{\alpha}\gamma^{q+1}$. Moreover, by the possible AADs \eqref{4ade_eqa} of $\dot{\alpha}\gamma^{q+1}$, the AAD of the vertex is $\thin^{\epsilon}\gamma_1^{\alpha}\thin^{\gamma}\alpha^{\beta}\thin^{\epsilon}\gamma^{\alpha}\thin\cdots=\thin^{\beta}\alpha^{\gamma}\thin^{\alpha}\gamma\cdots\gamma^{\epsilon}\thin$. The $\thin^{\alpha}\gamma\cdots\gamma^{\epsilon}\thin$ part of the vertex induces a partial earth map tiling ${\mc P}^{q+1}_q$, which consists of $P,T_1,T_2,T_3$, and fills the middle part of the tiling. The tiling is the third flip modification $F_3E_{\pentagon}2$ in Figure \ref{emt2mod3}. 

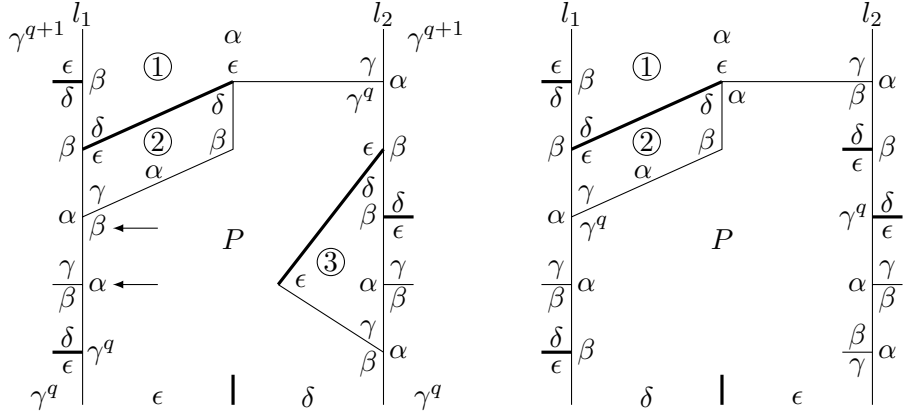
\begin{figure}[htp]
\centering
\begin{tikzpicture}[>=latex,scale=1]

\foreach \a in {0,1}
{
\begin{scope}[xshift=6.5*\a cm]

\draw
	(2,1.8) -- (0,1.8)
	(-2,2.5) -- (-2,-2.5)
	(2,2.5) -- (2,-2.5)
	(-2,0) -- (0,0.9) -- (0,1.8)
	(-2,-0.9) -- ++(-0.4,0)
	(2,-0.9) -- ++(0.4,0);

\draw[line width=1.2]
	(2,0) -- ++(0.4,0) 
	(-2,1.8) -- ++(-0.4,0)
	(-2,-1.8) -- ++(-0.4,0)
	(0,1.8) -- (-2,0.9);

\node at (-2,2.7) {\small $l_1$};
\node at (2,2.7) {\small $l_2$};

\node at (2.2,1.8) {\small $\alpha$};
\node at (2.2,0.9) {\small $\beta$};
\node at (2.2,0.2) {\small $\delta$};
\node at (2.2,-0.2) {\small $\epsilon$};
\node at (2.2,-0.7) {\small $\gamma$};
\node at (2.2,-1.1) {\small $\beta$};
\node at (2.2,-1.8) {\small $\alpha$};

\node at (-2.2,2) {\small $\epsilon$};
\node at (-2.2,1.6) {\small $\delta$};
\node at (-2.2,0.9) {\small $\beta$};
\node at (-2.2,0) {\small $\alpha$};
\node at (-2.2,-0.7) {\small $\gamma$};
\node at (-2.2,-1.1) {\small $\beta$};
\node at (-2.2,-1.6) {\small $\delta$};
\node at (-2.2,-2) {\small $\epsilon$};

\node at (0,2.4) {\small $\alpha$};
\node at (1.8,2) {\small $\gamma$};
\node at (-1.8,1.8) {\small $\beta$};
\node at (0,2) {\small $\epsilon$};
\node at (-1.8,1.2) {\small $\delta$};

\node at (-0.2,1) {\small $\beta$};
\node at (-1.05,0.6) {\small $\alpha$};
\node at (-1.8,0.3) {\small $\gamma$};
\node at (-1.8,0.8) {\small $\epsilon$};
\node at (-0.2,1.5) {\small $\delta$};

\node at (0,-0.3) {\small $P$};

\node[inner sep=0.5,draw,shape=circle] at (-1,2) {\small 1};
\node[inner sep=0.5,draw,shape=circle] at (-1,1) {\small 2};

\end{scope}
}


\draw
	(2,-1.8) -- (0.6,-0.9);

\draw[line width=1.2]
	(2,0.9) -- (0.6,-0.9)
	(0,-2.5) -- ++(0,0.4);
	
\draw [->]
	(-1,-0.15) -- ++(-0.6,0);
\draw [->]
	(-1,-0.9) -- ++(-0.6,0);

\node at (1.7,1.55) {\small $\gamma^q$};
\node at (1.8,0.9) {\small $\epsilon$};
\node at (1.8,-1.95) {\small $\beta$};
\node at (1,-2.4) {\small $\delta$};
\node at (-1,-2.4) {\small $\epsilon$};
\node at (-1.75,-1.8) {\small $\gamma^q$};
\node at (-1.8,-0.9) {\small $\alpha$};
\node at (-1.8,-0.15) {\small $\beta$};

\node at (1.8,-0.9) {\small $\alpha$};
\node at (1.8,-1.5) {\small $\gamma$};
\node at (1.8,0) {\small $\beta$};
\node at (0.9,-0.85) {\small $\epsilon$};
\node at (1.8,0.4) {\small $\delta$};

\node at (2.7,2.4) {\small $\gamma^{q+1}$};
\node at (2.6,-2.4) {\small $\gamma^q$};
\node at (-2.6,2.4) {\small $\gamma^{q+1}$};
\node at (-2.5,-2.4) {\small $\gamma^q$};

\node[inner sep=0.5,draw,shape=circle] at (1.3,-0.6) {\small 3};


\begin{scope}[xshift=6.5 cm]

\draw
	(2,-1.8) -- ++(-0.4,0);

\draw[line width=1.2]
	(2,0.9) -- ++(-0.4,0)
	(0,-2.5) -- ++(0,0.4);

\node at (-1.7,-0.15) {\small $\gamma^q$};
\node at (-1.8,-0.9) {\small $\alpha$};
\node at (-1.8,-1.8) {\small $\beta$};
\node at (1.8,-2) {\small $\gamma$};
\node at (1.8,-1.6) {\small $\beta$};
\node at (1.8,-0.9) {\small $\alpha$};
\node at (1.75,0) {\small $\gamma^q$};
\node at (1.8,0.7) {\small $\epsilon$};
\node at (1.8,1.1) {\small $\delta$};
\node at (1.8,1.6) {\small $\beta$};
\node at (0.2,1.6) {\small $\alpha$};
\node at (-1,-2.4) {\small $\delta$};
\node at (1,-2.4) {\small $\epsilon$};

\end{scope}

\end{tikzpicture}
\caption{Proposition \ref{4ade}: Tiling for $\thin^{\beta}\alpha^{\gamma}\thin^{\alpha}\gamma\cdots\gamma^{\epsilon}\thin$.}
\label{4adeF}
\end{figure}

The second of Figure \ref{4adeF} is the case $\gamma_2\cdots=\dot{\alpha}\gamma^{q+1}$. This implies the vertex $\beta\gamma\cdots$ on the left of $l_1$ is not $\dot{\alpha}\gamma^{q+1}$, and must be $\dot{\alpha}^2\gamma$. Then by the possible AADs \eqref{4ade_eqa} of $\dot{\alpha}\gamma^{q+1}$, the AAD of $\gamma_2\cdots$ is $\thin^{\alpha}\gamma_2^{\epsilon}\thin^{\beta}\alpha^{\gamma}\thin^{\alpha}\gamma^{\epsilon}\thin\cdots=\thin^{\beta}\alpha^{\gamma}\thin^{\alpha}\gamma\cdots\gamma^{\epsilon}\thin$. Then the $\thin^{\alpha}\gamma\cdots\gamma^{\epsilon}\thin$ part of the vertex induces a partial earth map tiling ${\mc P}^{q+1}_q$, which consists of $P,T_1,T_2$, and fills the middle part of the tiling. The tiling is the first flip modification $F_1E_{\pentagon}2$ in Figure \ref{emt2mod3}. 

We conclude the tiling with $\thin^{\beta}\alpha^{\gamma}\thin^{\alpha}\gamma\cdots\gamma^{\epsilon}\thin$ is either $F_1E_{\pentagon}2$ or $F_3E_{\pentagon}2$. Moreover, if the tiling contains the patch ${\mc A}$ in Figure \ref{emt2mod8}, then it is $F_3E_{\pentagon}2$.

Suppose the second AAD $\thin^{\beta}\alpha^{\gamma}\thin^{\epsilon}\gamma\cdots\gamma^{\alpha}\thin$ of \eqref{4ade_eqa} is a vertex. The AAD gives ${\mc P}^{q+1}_q$ on the left of $l_1$ and the right of $l_2$, and determines $T_1$ in Figure \ref{4adeE}. Then $\gamma_1\cdots=\gamma\delta\epsilon\cdots=\gamma^q\delta\epsilon$. By the possible AADs \eqref{4ade_eqde} of $\gamma^q\delta\epsilon$, we get $\gamma_1\cdots=\thin^{\epsilon}\gamma_1^{\alpha}\thin^{\gamma}\epsilon^{\delta}\thick\cdots=\thick^{\epsilon}\delta^{\beta}\thin^{\epsilon}\gamma\cdots\gamma^{\alpha}\thin^{\gamma}\epsilon^{\delta}\thick$. The $\thin^{\epsilon}\gamma\cdots\gamma^{\alpha}\thin$ part of the vertex induces a partial earth map tiling ${\mc P}^q_{q-1}$, which consists of $P$ and $T_1$. Then $\beta_1\cdots=\dot{\alpha}^2\cdots=\dot{\alpha}^2\gamma$ and $\delta_1\cdots=\delta\epsilon\cdots=\dot{\alpha}\delta\epsilon,\gamma^q\delta\epsilon$ determine $T_2$. Then $T_2$ determines its companion $T_3$. Then $\gamma_3\cdots=\dot{\alpha}\gamma^p\cdots=\dot{\alpha}\gamma^{p+1}$. By the possible AADs \eqref{4ade_eqa} of $\dot{\alpha}\gamma^{q+1}$, we get $\gamma_3\cdots=\thin^{\gamma}\alpha^{\beta}\thin^{\alpha}\gamma_3^{\epsilon}\thin\gamma\thin\cdots=\thin^{\beta}\alpha^{\gamma}\thin^{\epsilon}\gamma\cdots\gamma^{\alpha}\thin$. This determines $T_4$. Then $T_4$ determines its companion $T_5$. The tiling is the second flip modification $F_2E_{\pentagon}2$ in Figure \ref{emt2mod3}.

\begin{figure}[htp]
\centering
\begin{tikzpicture}[>=latex,scale=1]

\draw
	(-2,1.8) -- (2,1.8)
	(-2,2.5) -- (-2,-2.5)
	(2,2.5) -- (2,-2.5)
	(-2,-0.9) -- ++(-0.4,0)
	(2,-0.9) -- ++(0.4,0)
	(2,-1.8) -- (0.2,-1.8) -- (-0.2,0) -- (-2,0)
	(0.6,1.8) -- (0.6,0) -- (2,-1.8)
	(-0.2,0) -- ++(0,0.4)
	(0.6,0.9) -- (0.2,0.9);

\draw[line width=1.2]
	(2,0) -- ++(0.4,0) 
	(-2,1.8) -- ++(-0.4,0)
	(-2,-1.8) -- ++(-0.4,0)
	(-0.6,1.8) -- (0.6,1.8)
	(-2,0.9) -- ++(0.4,0)
	(0.6,0) -- (2,0.9)
	(-1,0) -- (-1,-2.5)
	(0.2,-1.8) -- ++(0.3,0.4);

\node at (-2,2.7) {\small $l_1$};
\node at (2,2.7) {\small $l_2$};

\node at (2.2,1.8) {\small $\alpha$};
\node at (2.2,0.9) {\small $\beta$};
\node at (2.2,0.2) {\small $\delta$};
\node at (2.2,-0.2) {\small $\epsilon$};
\node at (2.2,-0.7) {\small $\gamma$};
\node at (2.2,-1.1) {\small $\beta$};
\node at (2.2,-1.8) {\small $\alpha$};

\node at (-2.2,2) {\small $\epsilon$};
\node at (-2.2,1.6) {\small $\delta$};
\node at (-2.2,0.9) {\small $\beta$};
\node at (-2.2,0) {\small $\alpha$};
\node at (-2.2,-0.7) {\small $\gamma$};
\node at (-2.2,-1.1) {\small $\beta$};
\node at (-2.2,-1.6) {\small $\delta$};
\node at (-2.2,-2) {\small $\epsilon$};

\node at (0,2.5) {\small $\alpha$};
\node at (1.8,2) {\small $\beta$};
\node at (-1.8,2) {\small $\gamma$};
\node at (-0.6,2) {\small $\epsilon$};
\node at (0.6,2) {\small $\delta$};

\node at (0.1,-1.9) {\small $\alpha$};
\node at (1.8,-2) {\small $\gamma$};
\node at (0.5,-2.4) {\small $\epsilon$};
\node at (-0.8,-0.2) {\small $\delta$};
\node at (-0.35,-0.2) {\small $\beta$};

\node at (-1.55,1.55) {\small $\gamma^{q-1}$};
\node at (-1.8,1.1) {\small $\epsilon$};
\node at (-1.8,0.7) {\small $\delta$};
\node at (-1.8,0.2) {\small $\beta$};
\node at (-1,0.2) {\small $\alpha$};
\node at (-0.4,0.2) {\small $\gamma$};
\node at (-0,0) {\small $\alpha$};
\node at (0.3,-1.4) {\small $\delta$};
\node at (0.55,-1.6) {\small $\epsilon$};
\node at (1.3,-1.55) {\small $\gamma^{q-1}$};
\node at (0.4,0) {\small $\alpha$};
\node at (0.4,1.6) {\small $\epsilon$};
\node at (0.4,1.1) {\small $\gamma$};
\node at (0.4,0.7) {\small $\beta$};

\node at (-0.8,0.9) {\small $P$};

\node at (0.8,1.6) {\small $\alpha$};
\node at (0.8,0.9) {\small $\beta$};
\node at (1.8,1.6) {\small $\gamma$};
\node at (1.8,1) {\small $\epsilon$};
\node at (0.8,0.35) {\small $\delta$};

\node at (1.8,0) {\small $\beta$};
\node at (1.8,-0.9) {\small $\alpha$};
\node at (1.8,-1.3) {\small $\gamma$};
\node at (0.9,-0.05) {\small $\epsilon$};
\node at (1.8,0.55) {\small $\delta$};

\node at (-1.8,-0.9) {\small $\alpha$};
\node at (-1.8,-1.8) {\small $\beta$};
\node at (-1.8,-0.2) {\small $\gamma$};
\node at (-1.2,-0.2) {\small $\epsilon$};
\node at (-1.5,-2.4) {\small $\delta$};

\node[inner sep=0.5,draw,shape=circle] at (0,2.1) {\small 1};
\node[inner sep=0.5,draw,shape=circle] at (1.3,1.2) {\small 2};
\node[inner sep=0.5,draw,shape=circle] at (1.3,-0.3) {\small 3};
\node[inner sep=0.5,draw,shape=circle] at (-0.5,-1.35) {\small 4};
\node[inner sep=0.5,draw,shape=circle] at (-1.5,-1.35) {\small 5};

\end{tikzpicture}
\caption{Proposition \ref{4ade}: Tiling for $\thin^{\beta}\alpha^{\gamma}\thin^{\epsilon}\gamma\cdots\gamma^{\alpha}\thin$.}
\label{4adeE}
\end{figure}
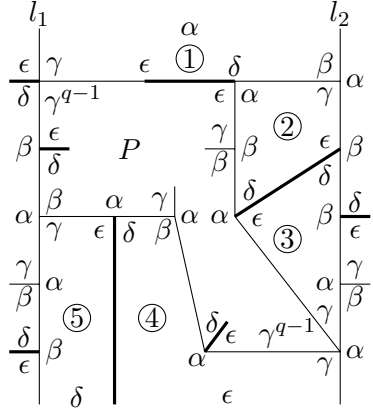

We conclude the tiling with $\thin^{\beta}\alpha^{\gamma}\thin^{\epsilon}\gamma\cdots\gamma^{\alpha}\thin$ is $F_2E_{\pentagon}2$. In fact, if we use \eqref{4ade_eqa_shift} to shift and flip ${\mc A}$ to get $\thin^{\beta}\alpha^{\gamma}\thin^{\alpha}\gamma\cdots\gamma^{\epsilon}\thin$, then we get a tiling with ${\mc A}$. By the earlier argument, we know the tiling after the change is $F_3E_{\pentagon}2$. Therefore $F_2E_{\pentagon}2,F_3E_{\pentagon}2$ are related by the shift and flip of ${\mc A}$.

Suppose the third AAD $\thin^{\beta}\alpha^{\gamma}\thin^{\epsilon}\gamma\cdots\gamma^{\epsilon}\thin$ of \eqref{4ade_eqa} is a vertex. By \eqref{4ade_eqde_shift}, we may shift and flip ${\mc A}$ to change the AAD to $\thin^{\beta}\alpha^{\gamma}\thin^{\epsilon}\gamma\cdots\gamma^{\alpha}\thin$, which gives the tiling $F_2E_{\pentagon}2$. Therefore the tiling with $\thin^{\beta}\alpha^{\gamma}\thin^{\epsilon}\gamma\cdots\gamma^{\epsilon}\thin$ is $F'_2E_{\pentagon}2$, obtained from $F_2E_{\pentagon}2$ by the shift and flip of ${\mc A}$.

Suppose the fourth AAD $\thin^{\alpha}\beta^{\delta}\thin^{\alpha}\gamma\cdots\gamma^{\epsilon}\thin$ of \eqref{4ade_eqa} is a vertex. The AAD gives ${\mc P}^{q+1}_q$ on the left of $l_1$ and the right of $l_2$, and determines $T_1$ in Figure \ref{4adeH}. Then $T_1$ determines its companion $T_2$. Then $\gamma_2\cdots=\beta\gamma\cdots=\dot{\alpha}^2\gamma,\dot{\alpha}\gamma^{q+1}$.

The first of Figure \ref{4adeH} is the case $\gamma_2\cdots=\dot{\alpha}\gamma^{q+1}$. By the possible AADs \eqref{4ade_eqa} of $\dot{\alpha}\gamma^{q+1}$, the AAD of the vertex is $\thin^{\alpha}\beta^{\delta}\thin^{\alpha}\gamma\cdots\gamma^{\epsilon}\thin$. The $\thin^{\alpha}\gamma\cdots\gamma^{\epsilon}\thin$ part of the vertex gives a partial earth map tiling ${\mc P}^{q+1}_q$, which consists of $P,T_1,T_2$, and fills the middle part of the tiling. The tiling is the $120^{\circ}$ degree rotation modification $RE_{\pentagon}2$ in Figure \ref{emt2mod2}.

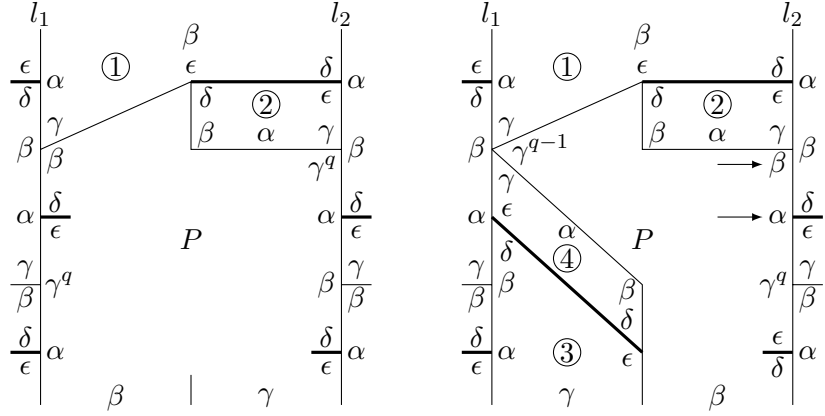
\begin{figure}[htp]
\centering
\begin{tikzpicture}[>=latex,scale=1]

\foreach \a in {0,1}
{
\begin{scope}[xshift=6*\a cm]

\draw
	(-2,0.9) -- (0,1.8) -- (0,0.9) -- (2,0.9)
	(-2,2.5) -- (-2,-2.5)
	(2,2.5) -- (2,-2.5)
	(-2,-0.9) -- ++(-0.4,0)
	(2,-0.9) -- ++(0.4,0);

\draw[line width=1.2]
	(2,1.8) -- (0,1.8)
	(2,0) -- ++(0.4,0) 
	(-2,1.8) -- ++(-0.4,0)
	(-2,-1.8) -- ++(-0.4,0);

\node at (-2,2.7) {\small $l_1$};
\node at (2,2.7) {\small $l_2$};

\node at (2.2,1.8) {\small $\alpha$};
\node at (2.2,0.9) {\small $\beta$};
\node at (2.2,0.2) {\small $\delta$};
\node at (2.2,-0.2) {\small $\epsilon$};
\node at (2.2,-0.7) {\small $\gamma$};
\node at (2.2,-1.1) {\small $\beta$};
\node at (2.2,-1.8) {\small $\alpha$};

\node at (-2.2,2) {\small $\epsilon$};
\node at (-2.2,1.6) {\small $\delta$};
\node at (-2.2,0.9) {\small $\beta$};
\node at (-2.2,0) {\small $\alpha$};
\node at (-2.2,-0.7) {\small $\gamma$};
\node at (-2.2,-1.1) {\small $\beta$};
\node at (-2.2,-1.6) {\small $\delta$};
\node at (-2.2,-2) {\small $\epsilon$};

\node at (0,2.4) {\small $\beta$};
\node at (1.8,2) {\small $\delta$};
\node at (-1.8,1.8) {\small $\alpha$};
\node at (0,2) {\small $\epsilon$};
\node at (-1.8,1.2) {\small $\gamma$};

\node at (1,1.1) {\small $\alpha$};
\node at (1.8,1.1) {\small $\gamma$};
\node at (0.2,1.1) {\small $\beta$};
\node at (1.8,1.6) {\small $\epsilon$};
\node at (0.2,1.6) {\small $\delta$};

\node at (0,-0.3) {\small $P$};

\node[inner sep=0.5,draw,shape=circle] at (-1,2) {\small 1};
\node[inner sep=0.5,draw,shape=circle] at (1,1.5) {\small 2};

\end{scope}
}


\draw
	(0,-2.5) -- ++(0,0.4);

\draw[line width=1.2]
	(2,-1.8) -- ++(-0.4,0)
	(-2,0) -- ++(0.4,0);

\node at (1.75,0.65) {\small $\gamma^q$};
\node at (1.8,0) {\small $\alpha$};
\node at (1.8,-0.9) {\small $\beta$};
\node at (1.8,-1.6) {\small $\delta$};
\node at (1.8,-2) {\small $\epsilon$};

\node at (-1.8,0.75) {\small $\beta$};
\node at (-1.8,0.2) {\small $\delta$};
\node at (-1.8,-0.2) {\small $\epsilon$};
\node at (-1.75,-0.9) {\small $\gamma^q$};
\node at (-1.8,-1.8) {\small $\alpha$};

\node at (-1,-2.4) {\small $\beta$};
\node at (1,-2.4) {\small $\gamma$};


\begin{scope}[xshift=6cm]

\draw
	(0,-2.5) -- (0,-0.9) -- (-2,0.9);

\draw[line width=1.2]
	(2,-1.8) -- ++(-0.4,0)
	(-2,0) -- (0,-1.8);

\draw[->]
	(1,0) -- ++(0.6,0);
\draw[->]
	(1,0.7) -- ++(0.6,0);
	
\node at (1.8,0.7) {\small $\beta$};
\node at (1.8,0) {\small $\alpha$};
\node at (1.75,-0.9) {\small $\gamma^q$};
\node at (1.8,-1.6) {\small $\epsilon$};
\node at (1.8,-2) {\small $\delta$};

\node at (-1.35,0.9) {\small $\gamma^{q-1}$};

\node at (-1,-0.2) {\small $\alpha$};
\node at (-0.2,-1) {\small $\beta$};
\node at (-1.8,0.45) {\small $\gamma$};
\node at (-0.2,-1.35) {\small $\delta$};
\node at (-1.8,0.1) {\small $\epsilon$};

\node at (-1.8,-0.45) {\small $\delta$};
\node at (-1.8,-0.9) {\small $\beta$};
\node at (-1.8,-1.8) {\small $\alpha$};
\node at (-0.2,-1.9) {\small $\epsilon$};
\node at (-1,-2.4) {\small $\gamma$};

\node at (1,-2.4) {\small $\beta$};

\node[inner sep=0.5,draw,shape=circle] at (-1,-1.8) {\small 3};
\node[inner sep=0.5,draw,shape=circle] at (-1,-0.55) {\small 4};

\end{scope}

\end{tikzpicture}
\caption{Proposition \ref{4ade}: Tiling for $\thin^{\alpha}\beta^{\delta}\thin^{\alpha}\gamma\cdots\gamma^{\epsilon}\thin$.}
\label{4adeH}
\end{figure}

The second of Figure \ref{4adeH} is the case $\gamma_2\cdots=\dot{\alpha}^2\gamma$. The vertex $\delta\epsilon\cdots$ on the right of $l_2$ is $\dot{\alpha}\delta\epsilon,\gamma^q\delta\epsilon$. If the vertex is $\gamma^q\delta\epsilon$, then by the possible AADs \eqref{4ade_eqde} of $\gamma^q\delta\epsilon$, its AAD is $\thin^{\epsilon}\gamma^{\alpha}\thin^{\beta}\delta^{\epsilon}\thick\cdots=\thick^{\epsilon}\delta^{\beta}\thin^{\alpha}\gamma\cdots\gamma^{\epsilon}\thin^{\gamma}\epsilon^{\delta}\thick$. This implies the vertex $\beta\gamma\cdots$ on the right of $l_2$ is $\beta\gamma\epsilon\cdots=\dot{\alpha}\gamma\epsilon\cdots$, a contradiction. Therefore the vertex $\delta\epsilon\cdots$ on the right of $l_2$ is $\dot{\alpha}\delta\epsilon$, and we get two adjacent $\dot{\alpha}$ indicated by the arrows. The vertex $\beta\gamma\cdots$ on the right of $l_2$ is $\dot{\alpha}^2\gamma,\dot{\alpha}\gamma^{q+1}$. By no three adjacent $\dot{\alpha}$, we know the vertex is $\dot{\alpha}\gamma^{q+1}$. This determines which of the two adjacent $\dot{\alpha}$ indicated by the arrows is $\alpha$ or $\beta$. Then by the possible AADs \eqref{4ade_eqa} of $\dot{\alpha}\gamma^{q+1}$, the vertex is $\thin^{\alpha}\beta^{\delta}\thin^{\alpha}\gamma^{\epsilon}\thin\cdots=\thin^{\alpha}\beta^{\delta}\thin^{\alpha}\gamma\cdots\gamma^{\epsilon}\thin$. This implies the vertex $\alpha\cdots$ on the lower right of $l_2$ is $\alpha\epsilon\cdots=\dot{\alpha}\delta\epsilon$. This gives $\beta$ at the bottom of the picture. Then the bottom vertex is $\beta\gamma^q\cdots=\dot{\alpha}\gamma^{q+1}$. Moreover, the vertex $\delta\epsilon\cdots$ on the lower left of $l_1$ is $\dot{\alpha}\delta\epsilon,\gamma^q\delta\epsilon$. The two vertices determine $T_3$. Then $T_3$ determines its companion $T_4$. Then $\gamma_1\gamma_4\cdots=\beta\gamma^2\cdots=\dot{\alpha}\gamma^{q+1}$. By the possible AADs \eqref{4ade_eqa} of $\dot{\alpha}\gamma^{q+1}$, we get $\gamma_1\gamma_4\cdots=\thin^{\alpha}\gamma_4^{\epsilon}\thin^{\alpha}\beta^{\delta}\thin^{\alpha}\gamma_1^{\epsilon}\thin\cdots=\thin^{\alpha}\beta^{\delta}\thin^{\alpha}\gamma\cdots\gamma^{\epsilon}\thin$. The $\thin^{\alpha}\gamma\cdots\gamma^{\epsilon}\thin$ part of the vertex gives a partial earth map tiling ${\mc P}^{q+1}_q$, which consists of $P,T_1,T_2,T_3,T_4$, and fills the middle part of the tiling. The tiling is the $240^{\circ}$ degree rotation modification $RE_{\pentagon}2$ in Figure \ref{emt2mod2}.

\subsubsection*{Subcase. $\dot{\alpha}\gamma^{q+1}$ is not a vertex, and $\gamma^q\delta\epsilon$ is a vertex.}

We get updated list of all the vertices
\[
\text{AVC}
=\{\dot{\alpha}\delta\epsilon,\dot{\alpha}^2\gamma,\gamma^q\delta\epsilon\},
\quad q=\tfrac{f-4}{8}.
\]
In particular, we have $\dot{\alpha}\gamma\cdots=\dot{\alpha}^2\gamma$. Moreover, since $f>12$ implies $q>2$, we have $\gamma^q\cdots=\gamma^q\delta\epsilon$, $\gamma\delta\epsilon\cdots=\gamma^q\delta\epsilon$, and $\gamma^{q-1}\delta\epsilon\cdots=\gamma^q\delta\epsilon$. 

Suppose the second AAD $\thick^{\epsilon}\delta^{\beta}\thin^{\epsilon}\gamma\cdots\gamma^{\alpha}\thin^{\gamma}\epsilon^{\delta}\thick$ of \eqref{4ade_eqde} is a vertex. The AAD gives ${\mc P}^q_{q-1}$ on the left of $l_1$ and the right of $l_2$, and determines $T_1,T_2$ in Figure \ref{4adeL}. Then $\alpha_1\cdots=\alpha\beta\cdots=\dot{\alpha}^2\gamma$ gives $\gamma_3$. There are two ways of arranging $T_3$, indicated by the two pictures. 

In the first of Figure \ref{4adeL}, $T_3$ determines its companion $T_4$. Then $\gamma_4\cdots=\gamma^q\cdots=\gamma^q\delta\epsilon=\thick^{\epsilon}\delta^{\beta}\thin^{\epsilon}\gamma\cdots\gamma^{\alpha}\thin^{\gamma}\epsilon^{\delta}\thick$ or $\thick^{\epsilon}\delta^{\beta}\thin^{\alpha}\gamma\cdots\gamma^{\epsilon}\thin^{\gamma}\epsilon^{\delta}\thick$. 

If $\gamma_4\cdots=\thick^{\epsilon}\delta^{\beta}\thin^{\epsilon}\gamma\cdots\gamma^{\alpha}\thin^{\gamma}\epsilon^{\delta}\thick$, then we get $T_5,T_6$. Then $\alpha_3\gamma_1\cdots=\alpha_5\beta_3\cdots=\dot{\alpha}^2\gamma$ determines $T_7$. Then $T_7$ determines its companion $T_8$. 

If $\gamma_4\cdots=\thick^{\epsilon}\delta^{\beta}\thin^{\alpha}\gamma\cdots\gamma^{\epsilon}\thin^{\gamma}\epsilon^{\delta}\thick$, then $\delta_5,\epsilon_6$ are switched, and $\delta_3\epsilon_4\cdots=\gamma\delta\epsilon\cdots=\gamma^q\delta\epsilon=\thick^{\epsilon}\delta^{\beta}\thin^{\alpha}\gamma\cdots\gamma^{\epsilon}\thin^{\gamma}\epsilon^{\delta}\thick$ or $\thick^{\epsilon}\delta^{\beta}\thin^{\epsilon}\gamma\cdots\gamma^{\epsilon}\thin^{\gamma}\epsilon^{\delta}\thick$. If $\delta_3\epsilon_4\cdots=\thick^{\epsilon}\delta^{\beta}\thin^{\epsilon}\gamma\cdots\gamma^{\epsilon}\thin^{\gamma}\epsilon^{\delta}\thick$, then by \eqref{4ade_eqde_shift}, we may shift and flip ${\mc A}$ to change the AAD to $\thick^{\epsilon}\delta^{\beta}\thin^{\alpha}\gamma\cdots\gamma^{\epsilon}\thin^{\gamma}\epsilon^{\delta}\thick$. Therefore we only need to consider the case $\delta_3\epsilon_4\cdots=\thick^{\epsilon}\delta^{\beta}\thin^{\alpha}\gamma\cdots\gamma^{\epsilon}\thin^{\gamma}\epsilon^{\delta}\thick$. Then we still get $\alpha_5$ as in the picture. Then we may use the same reason to get the same $T_7,T_8$. 

Now $T_1,T_2,T_7,T_8$ form the patch ${\mc A}$ in Figure \ref{emt2mod8}. If we flip the patch, then the vertex $\delta_1\epsilon_2\cdots=\gamma^q\delta\epsilon$ changes to $\alpha\gamma^{q+1}$, and the AAD of $\alpha\gamma^{q+1}$ is the first $\thin^{\beta}\alpha^{\gamma}\thin^{\alpha}\gamma\cdots\gamma^{\epsilon}\thin$ of \eqref{4ade_eqa}. Since the tiling contains ${\mc A}$, it is $F_3E_{\pentagon}2$, or $F'_2E_{\pentagon}2$. Considering the possible shift and flip of ${\mc A}$ needed to change the AAD of $\delta_3\epsilon_4\cdots$, we conclude the tiling in the first of Figure \ref{4adeL} is $F'_2E_{\pentagon}2$ or $F''_2E_{\pentagon}2$.

\begin{figure}[htp]
\centering
\begin{tikzpicture}[>=latex,scale=1]

\foreach \a in {0,1}
{
\begin{scope}[xshift=6*\a cm]

\draw
	(2,1.8) -- (0,1.8) -- (-1,0.9) -- (-2,0.9)
	(-2,2.5) -- (-2,-2.5)
	(2,2.5) -- (2,-2.5)
	(-1,0.9) -- (-1,-2.5)
	(-1,-1.8) -- (-2,0)
	(2,-1.8) -- (0,-1.8)
	(2,-0.9) -- ++(0.4,0)
	(-2,-0.9) -- ++(-0.4,0);

\draw[line width=1.2]
	(0,1.8) -- (0,2.5)
	(2,0) -- ++(0.4,0)
	(-2,1.8) -- ++(-0.4,0)
	(-2,-1.8) -- ++(-0.4,0);

\node at (-2,2.7) {\small $l_1$};
\node at (2,2.7) {\small $l_2$};

\node at (2.2,1.8) {\small $\alpha$};
\node at (2.2,0.9) {\small $\beta$};
\node at (2.2,0.2) {\small $\delta$};
\node at (2.2,-0.2) {\small $\epsilon$};
\node at (2.2,-0.7) {\small $\gamma$};
\node at (2.2,-1.1) {\small $\beta$};
\node at (2.2,-1.8) {\small $\alpha$};

\node at (-2.2,2) {\small $\epsilon$};
\node at (-2.2,1.6) {\small $\delta$};
\node at (-2.2,0.9) {\small $\beta$};
\node at (-2.2,0) {\small $\alpha$};
\node at (-2.2,-0.7) {\small $\gamma$};
\node at (-2.2,-1.1) {\small $\beta$};
\node at (-2.2,-1.6) {\small $\delta$};
\node at (-2.2,-2) {\small $\epsilon$};

\node at (1.3,2) {\small $\alpha$};
\node at (0.7,2) {\small $\beta$};
\node at (1.8,2) {\small $\gamma$};
\node at (0.2,2) {\small $\delta$};
\node at (1,2.4) {\small $\epsilon$};

\node at (-1.8,1.1) {\small $\alpha$};
\node at (-1.1,1.1) {\small $\gamma$};
\node at (-1.8,1.8) {\small $\beta$};
\node at (-0.2,1.9) {\small $\epsilon$};
\node at (-1,2.4) {\small $\delta$};

\node[inner sep=0.5,draw,shape=circle] at (-1,1.8) {\small 1};
\node[inner sep=0.5,draw,shape=circle] at (1.5,2.3) {\small 2};
\node[inner sep=0.5,draw,shape=circle] at (-1.5,0.3) {\small 3};
\node[inner sep=0.5,draw,shape=circle] at (-1.6,-1.4) {\small 4};

\end{scope}
}

\node[inner sep=0.5,draw,shape=circle] at (-0.5,-1.4) {\small 5};
\node[inner sep=0.5,draw,shape=circle] at (1.5,-2.3) {\small 6};
\node[inner sep=0.5,draw,shape=circle] at (-0.5,0.32) {\small 7};
\node[inner sep=0.5,draw,shape=circle] at (0.7,0.3) {\small 8};


\draw
	(0,-1.8) -- (0,-0.7)
	(-1,-0.7) -- (1.3,-0.7) -- (1.3,1.8);

\draw[line width=1.2]
	(0,-1.8) -- (0,-2.5)
	(0.7,1.8) -- (0,-0.7)
	(-1,-1.8) -- (-2,0);

\node at (-0.8,0.8) {\small $\alpha$};
\node at (-0.8,-0.5) {\small $\gamma$};
\node at (-0,1.6) {\small $\beta$};
\node at (-0.1,-0.5) {\small $\epsilon$};
\node at (0.45,1.6) {\small $\delta$};

\node at (1.1,0.3) {\small $\alpha$};
\node at (1.1,1.6) {\small $\gamma$};
\node at (1.1,-0.5) {\small $\beta$};
\node at (0.8,1.6) {\small $\epsilon$};
\node at (0.3,-0.5) {\small $\delta$};

\node at (-1.2,0.7) {\small $\alpha$};
\node at (-1.8,0.7) {\small $\gamma$};
\node at (-1.2,-0.7) {\small $\beta$};
\node at (-1.8,0) {\small $\epsilon$};
\node at (-1.15,-1.2) {\small $\delta$};

\node at (-1.8,-1.8) {\small $\alpha$};
\node at (-1.5,-2.4) {\small $\gamma$};
\node at (-1.8,-0.9) {\small $\beta$};
\node at (-1.2,-1.8) {\small $\epsilon$};
\node at (-1.85,-0.6) {\small $\delta$};

\node at (-0.8,-0.9) {\small $\alpha$};
\node at (-0.2,-0.9) {\small $\gamma$};
\node at (-0.8,-1.8) {\small $\beta$};
\node at (-0.2,-1.8) {\small $\epsilon$};
\node at (-0.5,-2.4) {\small $\delta$};

\node at (1.3,-2) {\small $\alpha$};
\node at (0.7,-2) {\small $\beta$};
\node at (1.8,-2) {\small $\gamma$};
\node at (0.2,-2) {\small $\delta$};
\node at (1,-2.4) {\small $\epsilon$};


\begin{scope}[xshift=6cm]

\draw
	(-1,0.9) -- (0,0.9) -- (0,-2.5)
	(-1,-0.7) -- (0,-0.7);

\draw[line width=1.2]
	(-1,0.9) -- (-1,-0.7)
	(-1,-1.8) -- (-1,-2.5);

\node at (-1.2,0.7) {\small $\epsilon$};
\node at (-1.8,0.7) {\small $\gamma$};
\node at (-1.2,-0.7) {\small $\delta$};
\node at (-1.8,0) {\small $\alpha$};
\node at (-1.15,-1.2) {\small $\beta$};

\node at (-1.8,-1.8) {\small $\beta$};
\node at (-1.5,-2.4) {\small $\delta$};
\node at (-1.8,-0.9) {\small $\alpha$};
\node at (-1.2,-1.8) {\small $\epsilon$};
\node at (-1.85,-0.6) {\small $\gamma$};

\node at (-0.2,0) {\small $\alpha$};
\node at (-0.2,-0.5) {\small $\gamma$};
\node at (-0.2,0.7) {\small $\beta$};
\node at (-0.8,-0.5) {\small $\epsilon$};
\node at (-0.8,0.7) {\small $\delta$};

\node at (-0.8,-0.9) {\small $\beta$};
\node at (-0.2,-0.9) {\small $\alpha$};
\node at (-0.8,-1.8) {\small $\delta$};
\node at (-0.2,-1.8) {\small $\gamma$};
\node at (-0.5,-2.4) {\small $\epsilon$};

\node at (1,-2.4) {\small $\gamma$};

\node[inner sep=0.5,draw,shape=circle] at (-0.5,0.3) {\small 5};
\node[inner sep=0.5,draw,shape=circle] at (-0.5,-1.4) {\small 6};

\end{scope}

\end{tikzpicture}
\caption{Proposition \ref{4ade}: Tiling for $\thick^{\epsilon}\delta^{\beta}\thin^{\epsilon}\gamma\cdots\gamma^{\alpha}\thin^{\gamma}\epsilon^{\delta}\thick$.}
\label{4adeL}
\end{figure}
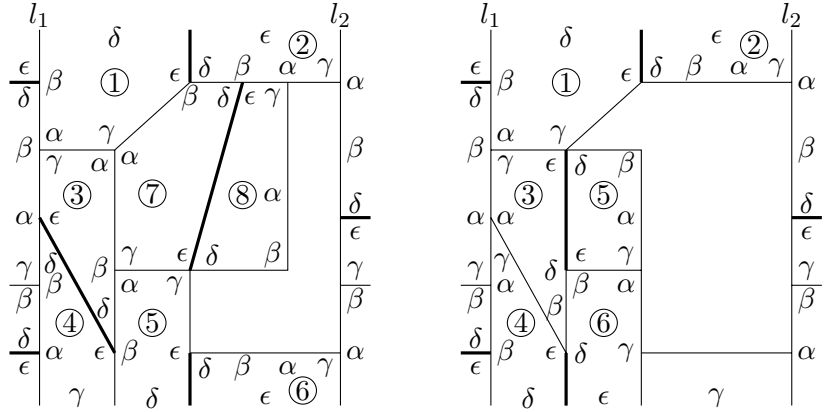

In the second of Figure \ref{4adeL}, we have $\alpha_3\cdots=\alpha^2\cdots=\dot{\alpha}^2\gamma$, and we know the vertex $\beta\gamma\cdots$ on the left of $l_1$ is $\beta\gamma\cdots=\dot{\alpha}^2\gamma$. This determines $T_4$. Then $T_3,T_4$ determine their companions $T_5,T_6$. Now $T_3,T_4,T_5,T_6$ form the patch ${\mc A}$ in Figure \ref{emt2mod8}. If we flip the patch, then the vertex $\delta_4\epsilon_6\cdots=\gamma^q\delta\epsilon$ changes to $\alpha\gamma^{q+1}$, and the AAD of $\alpha\gamma^{q+1}$ is the first $\thin^{\beta}\alpha^{\gamma}\thin^{\alpha}\gamma\cdots\gamma^{\epsilon}\thin$ or the third $\thin^{\beta}\alpha^{\gamma}\thin^{\epsilon}\gamma\cdots\gamma^{\epsilon}\thin$ of \eqref{4ade_eqa}. Since the tiling contains ${\mc A}$, the tiling in the second of Figure \ref{4adeL} is $F'_2E_{\pentagon}2$. 

Suppose the third AAD $\thick^{\epsilon}\delta^{\beta}\thin^{\epsilon}\gamma\cdots\gamma^{\epsilon}\thin^{\gamma}\epsilon^{\delta}\thick$ of \eqref{4ade_eqde} is a vertex. By \eqref{4ade_eqde_shift}, we may shift and flip ${\mc A}$ to change the AAD to the second $\thick^{\epsilon}\delta^{\beta}\thin^{\epsilon}\gamma\cdots\gamma^{\alpha}\thin^{\gamma}\epsilon^{\delta}\thick$ of \eqref{4ade_eqde}. Since the later gives the tiling $F'_2E_{\pentagon}2$ or $F''_2E_{\pentagon}2$, we know the AAD $\thick^{\epsilon}\delta^{\beta}\thin^{\epsilon}\gamma\cdots\gamma^{\epsilon}\thin^{\gamma}\epsilon^{\delta}\thick$ gives the tiling $F'_2E_{\pentagon}2$ or $F''_2E_{\pentagon}2$.

Suppose the first AAD $\thick^{\epsilon}\delta^{\beta}\thin^{\alpha}\gamma\cdots\gamma^{\epsilon}\thin^{\gamma}\epsilon^{\delta}\thick$ in \eqref{4ade_eqde} is a vertex. The $\thin^{\alpha}\gamma\cdots\gamma^{\epsilon}\thin$ part of the vertex gives ${\mc P}^q_{q-1}$ on the left of $l_1$ and the right of $l_2$ in Figure \ref{4adeK}, and the $\thick^{\epsilon}\delta^{\beta}\thin$ and $\thin^{\gamma}\epsilon^{\delta}\thick$ parts determine $T_1,T_2$. We have $\gamma_1\cdots=\thin^{\alpha}\gamma_1^{\epsilon}\thin^{\gamma}\epsilon^{\delta}\thick^{\epsilon}\delta^{\beta}\thin\cdots=\gamma^q\delta\epsilon=\thick^{\epsilon}\delta^{\beta}\thin^{\alpha}\gamma\cdots\gamma^{\epsilon}\thin^{\gamma}\epsilon^{\delta}\thick$ or $\thick^{\epsilon}\delta^{\beta}\thin^{\epsilon}\gamma\cdots\gamma^{\epsilon}\thin^{\gamma}\epsilon^{\delta}\thick$. If the vertex is $\thick^{\epsilon}\delta^{\beta}\thin^{\epsilon}\gamma\cdots\gamma^{\epsilon}\thin^{\gamma}\epsilon^{\delta}\thick$, then the tiling is $F'_2E_{\pentagon}2$ or $F''_2E_{\pentagon}2$. If the vertex is $\thick^{\epsilon}\delta^{\beta}\thin^{\alpha}\gamma\cdots\gamma^{\epsilon}\thin^{\gamma}\epsilon^{\delta}\thick$, then the $\thin^{\alpha}\gamma\cdots\gamma^{\epsilon}\thin$ part of the vertex induces a partial earth map tiling ${\mc P}^q_{q-1}$, which consists of $P,T_1,T_2$. Moreover, $\beta_2\cdots=\alpha\beta\cdots=\dot{\alpha}^2\gamma$ gives $\gamma_3$. Then we have two possible arrangements of $T_3$, given by the two pictures.

In the first of Figure \ref{4adeK}, $T_3$ determines its companion $T_4$. Then $\alpha_3\cdots=\alpha\beta\cdots=\dot{\alpha}^2\gamma$ and no $\delta\epsilon^2\cdots$ determine $T_5$. Then $T_5$ determines its companion $T_6$. Then $\alpha_6\gamma_4\cdots=\alpha\gamma^q\cdots=\dot{\alpha}\gamma^{q+1}=\thin^{\beta}\alpha^{\gamma}\thin^{\alpha}\gamma\cdots\gamma^{\epsilon}\thin$ or $\thin^{\beta}\alpha^{\gamma}\thin^{\epsilon}\gamma\cdots\gamma^{\epsilon}\thin$. Since the tiling contains the patch ${\mc A}$ consisting of $T_3,T_4,T_5,T_6$, the tiling is $F'_2E_{\pentagon}2$ or $F''_2E_{\pentagon}2$.

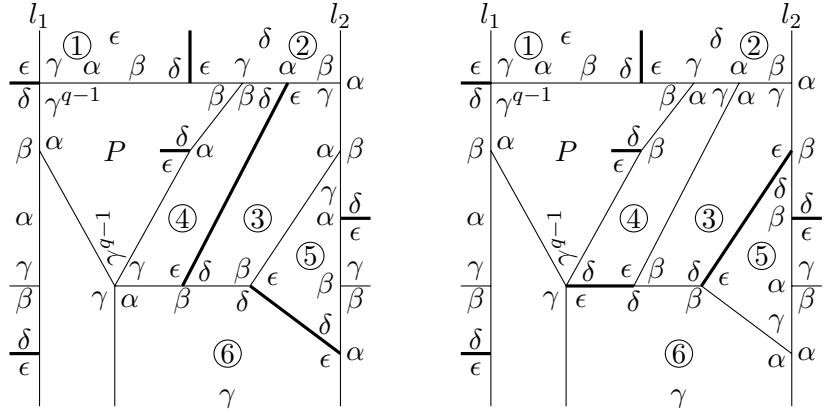
\begin{figure}[htp]
\centering
\begin{tikzpicture}[>=latex,scale=1]

\foreach \a in {0,1}
{
\begin{scope}[xshift=6*\a cm]

\draw
	(-2,1.8) -- (2,1.8)
	(-2,2.5) -- (-2,-2.5)
	(2,2.5) -- (2,-2.5)
	(0.7,1.8) -- (0,0.9) -- (-1,-0.9) -- (-1,-2.5)
	(2,0.9) -- (0.8,-0.9) -- (-1,-0.9) -- (-2,0.9)
	(0.8,-0.9) -- (2,-1.8)
	(1.3,1.8) -- (-0.1,-0.9)
	(-2,-0.9) -- ++(-0.4,0)
	(2,-0.9) -- ++(0.4,0);

\draw[line width=1.2]
	(0,1.8) -- (0,2.5)
	(2,0) -- ++(0.4,0)
	(-2,1.8) -- ++(-0.4,0)
	(-2,-1.8) -- ++(-0.4,0)
	(0,0.9) -- ++(-0.4,0);

\node at (-2,2.7) {\small $l_1$};
\node at (2,2.7) {\small $l_2$};

\node at (1.3,2) {\small $\alpha$};
\node at (0.7,2) {\small $\gamma$};
\node at (1.8,2) {\small $\beta$};
\node at (0.2,2) {\small $\epsilon$};
\node at (1,2.4) {\small $\delta$};

\node at (-1.3,2) {\small $\alpha$};
\node at (-0.7,2) {\small $\beta$};
\node at (-1.8,2) {\small $\gamma$};
\node at (-0.2,2) {\small $\delta$};
\node at (-1,2.4) {\small $\epsilon$};

\node at (2.2,1.8) {\small $\alpha$};
\node at (2.2,0.9) {\small $\beta$};
\node at (2.2,0.2) {\small $\delta$};
\node at (2.2,-0.2) {\small $\epsilon$};
\node at (2.2,-0.7) {\small $\gamma$};
\node at (2.2,-1.1) {\small $\beta$};
\node at (2.2,-1.8) {\small $\alpha$};

\node at (-2.2,2) {\small $\epsilon$};
\node at (-2.2,1.6) {\small $\delta$};
\node at (-2.2,0.9) {\small $\beta$};
\node at (-2.2,0) {\small $\alpha$};
\node at (-2.2,-0.7) {\small $\gamma$};
\node at (-2.2,-1.1) {\small $\beta$};
\node at (-2.2,-1.6) {\small $\delta$};
\node at (-2.2,-2) {\small $\epsilon$};

\node at (-1.55,1.55) {\small $\gamma^{q-1}$};
\node at (-1.8,1) {\small $\alpha$};
\node at (0.35,1.6) {\small $\beta$};
\node at (-0.1,1.1) {\small $\delta$};
\node at (-0.3,0.7) {\small $\epsilon$};
\node[rotate=90] at (-1.05,-0.25) {\small $\gamma^{q-1}$};
\node at (-1.2,-1.1) {\small $\gamma$};

\node at (-1,0.9) {\small $P$};

\node[inner sep=0.5,draw,shape=circle] at (-1.5,2.3) {\small 1};
\node[inner sep=0.5,draw,shape=circle] at (1.5,2.3) {\small 2};
\node[inner sep=0.5,draw,shape=circle] at (0.9,0) {\small 3};
\node[inner sep=0.5,draw,shape=circle] at (-0.1,0) {\small 4};
\node[inner sep=0.5,draw,shape=circle] at (1.6,-0.5) {\small 5};
\node[inner sep=0.5,draw,shape=circle] at (0.5,-1.8) {\small 6};

\end{scope}
}


\draw[line width=1.2]
	(1.3,1.8) -- (-0.1,-0.9)
	(0.8,-0.9) -- (2,-1.8);

\node at (1.8,0.9) {\small $\alpha$};
\node at (1.8,1.6) {\small $\gamma$};
\node at (0.7,-0.7) {\small $\beta$};
\node at (1.4,1.6) {\small $\epsilon$};
\node at (0.2,-0.7) {\small $\delta$};		

\node at (0.2,0.9) {\small $\alpha$};
\node at (-0.7,-0.7) {\small $\gamma$};
\node at (0.75,1.6) {\small $\beta$};	
\node at (-0.2,-0.7) {\small $\epsilon$};
\node at (1,1.6) {\small $\delta$};	

\node at (1.8,0) {\small $\alpha$};
\node at (1.85,0.3) {\small $\gamma$};
\node at (1.8,-0.9) {\small $\beta$};
\node at (1.1,-0.8) {\small $\epsilon$};
\node at (1.8,-1.4) {\small $\delta$};

\node at (-0.1,-1.1) {\small $\beta$};
\node at (0.5,-2.4) {\small $\gamma$};
\node at (-0.8,-1.1) {\small $\alpha$};
\node at (1.8,-1.9) {\small $\epsilon$};
\node at (0.7,-1.1) {\small $\delta$};


\begin{scope}[xshift=6cm]

\draw[line width=1.2]
	(-0.1,-0.9) -- (-1,-0.9)
	(0.8,-0.9) -- (2,0.9);

\node at (1.8,0.9) {\small $\epsilon$};
\node at (1.8,1.6) {\small $\gamma$};
\node at (0.7,-0.7) {\small $\delta$};
\node at (1.4,1.6) {\small $\alpha$};
\node at (0.2,-0.7) {\small $\beta$};		

\node at (0.2,0.9) {\small $\beta$};
\node at (-0.7,-0.7) {\small $\delta$};
\node at (0.75,1.6) {\small $\alpha$};	
\node at (-0.2,-0.7) {\small $\epsilon$};
\node at (1.05,1.6) {\small $\gamma$};	

\node at (1.8,0) {\small $\beta$};
\node at (1.85,0.4) {\small $\delta$};
\node at (1.8,-0.9) {\small $\alpha$};
\node at (1.1,-0.8) {\small $\epsilon$};
\node at (1.8,-1.4) {\small $\gamma$};

\node at (-0.8,-1.1) {\small $\epsilon$};
\node at (0.5,-2.4) {\small $\gamma$};
\node at (-0.1,-1.1) {\small $\delta$};
\node at (1.8,-1.9) {\small $\alpha$};
\node at (0.7,-1.1) {\small $\beta$};

\end{scope}

\end{tikzpicture}
\caption{Proposition \ref{4ade}: Tiling for $\thick^{\epsilon}\delta^{\beta}\thin^{\alpha}\gamma\cdots\gamma^{\epsilon}\thin^{\gamma}\epsilon^{\delta}\thick$.}
\label{4adeK}
\end{figure}

In the second of Figure \ref{4adeK}, $\alpha_2\alpha_3\cdots=\alpha^2\gamma$ and $\gamma_2\cdots=\alpha\gamma\cdots=\dot{\alpha}^2\gamma$ determine $T_4$. Then $T_3,T_4$ determine their companions $T_5,T_6$. Now $T_3,T_4,T_5,T_6$ form the patch ${\mc A}$. If we flip the patch, then we get the first of Figure \ref{4adeK}. Since the first tiling is $F'_2E_{\pentagon}2$ or $F''_2E_{\pentagon}2$, the second tiling is also $F'_2E_{\pentagon}2$ or $F''_2E_{\pentagon}2$.

\medskip

\noindent{\em Geometry of Pentagon}

\medskip

We justify the existence of the pentagon with the angle values given by \eqref{4ade_eq1}. By $\delta<\epsilon$, we also know $(\frac{1}{2}+\tfrac{2}{f})\pi<\epsilon<(1+\tfrac{4}{f})\pi$. 

The equality \eqref{coolsaet_eq1} becomes $A\cos^2\epsilon+B\cos\epsilon\sin\epsilon+C\sin^2\epsilon=0$, where $\theta=\tfrac{4}{f}\pi$ and
\begin{align*}
A &= 2+2\cos\theta-\cos3\theta-2\cos4\theta-\cos5\theta, \\
B &= 4\sin\theta+4\sin2\theta+4\sin3\theta+2\sin4\theta, \\
C &= -6\cos\theta-4\cos2\theta-\cos3\theta-4\cos4\theta-\cos5\theta.
\end{align*}
Figure \ref{4ade_geometry1} gives the graphs of $A$ and $C$ as functions of $\frac{1}{f}$. We see that $C\ne 0$ for $f\ge 16$ and multiple of $4$. This implies $\cos\epsilon\ne 0$, and we get a quadratic equation for $\tan\epsilon$. For each $f$, there are two solutions for $\tan\epsilon$. For each solution, we may further use $(\frac{1}{2}+\tfrac{2}{f})\pi<\epsilon<(1+\tfrac{4}{f})\pi$ to determine $\epsilon$. 

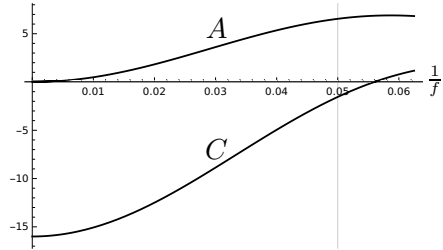
\begin{figure}[htp]
\centering
\begin{tikzpicture}[>=latex,scale=1]

\node[scale=0.6] at (0,0) {\includegraphics{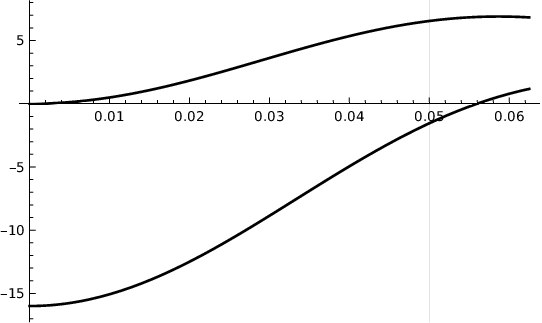}};

\node at (0,1.3) {\small $A$};
\node at (0,-0.3) {\small $C$};
\node at (2.9,0.6) {\scriptsize $\tfrac{1}{f}$};

\end{tikzpicture}
\caption{Proposition \ref{4ade}: The signs of $A$ and $C$.}
\label{4ade_geometry1}
\end{figure}

The graphs of $A$ and $C$ show that $\frac{C}{A}>0$ for $f=16$ and $\frac{C}{A}<0$ for $f\ge 20$. For $f=16$, only one solution satisfies $\frac{5}{8}\pi=(\frac{1}{2}+\tfrac{2}{16})\pi<\epsilon<(1+\tfrac{4}{16})\pi=\frac{5}{4}\pi$, and is given by
\[
\tan\epsilon=\tfrac{1}{\sqrt{2}}{\textstyle \sqrt{11+8\sqrt{2}} }-2-\tfrac{3}{\sqrt{2}},\quad
\epsilon=0.7889\pi.
\]
Then we get
\begin{align*}
\cos a &={\textstyle \sqrt{4\sqrt{2}-5} }, &
a &= 0.1992\pi, \\
\cos b &=3-2\sqrt{2}+{\textstyle \sqrt{44\sqrt{2}-62} }, &
b &= 0.2763\pi.
\end{align*}

For $f\ge 20$, by $\frac{C}{A}<0$, there is one positive solution $\tan\epsilon_+$ and one negative solution $\tan\epsilon_-$. Their signs mean $\pi<\epsilon_+<(1+\frac{4}{f})\pi$ and $(\frac{1}{2}+\frac{2}{f})\pi<\epsilon_-<\pi$. Therefore
\[
\tan\epsilon_+<\tan(1+\tfrac{4}{f})\pi, \quad
\tan\epsilon_->\tan(\tfrac{1}{2}+\tfrac{2}{f})\pi.
\]
Figure \ref{4ade_geometry2} shows the graphs of $\tan\epsilon_+,\tan\epsilon_-,\tan(1+\tfrac{4}{f})\pi$ as functions of $\frac{1}{f}$. We find that $\tan\epsilon_+>\tan(1+\tfrac{4}{f})\pi$, a contradiction. We also find
\[
\tan\epsilon_-
>-1
=\tan(\tfrac{1}{2}+\tfrac{2}{8})\pi
>\tan(\tfrac{1}{2}+\tfrac{2}{f})\pi.
\]
This means $(\frac{1}{2}+\frac{2}{f})\pi<\epsilon_-<\pi$, and only $\epsilon_-$ is acceptable for tiling.

\begin{figure}[htp]
\centering
\begin{tikzpicture}[>=latex,scale=1]

\node[scale=0.6] at (-3.5,0) {\includegraphics{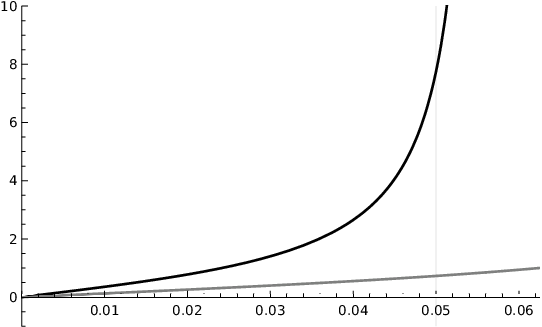}};

\node[scale=0.6] at (3.5,0) {\includegraphics{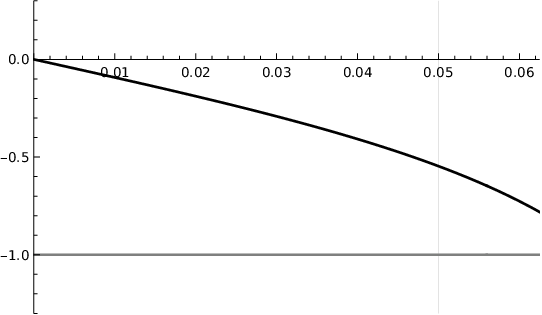}};

\node at (-0.6,-1.35) {\scriptsize $\tfrac{1}{f}$};
\node at (6.4,1) {\scriptsize $\tfrac{1}{f}$};
\node at (-3.6,-0.7) {\scriptsize $\tan\epsilon_+$};
\node at (-2,-0.9) {\scriptsize $\tan(1+\tfrac{4}{f})\pi$};
\node at (5,0.25) {\scriptsize $\tan\epsilon_-$};

\end{tikzpicture}
\caption{Proposition \ref{4ade}: $\tan\epsilon_+$ and $\tan\epsilon_-$ as functions of $\frac{1}{f}$.}
\label{4ade_geometry2}
\end{figure}
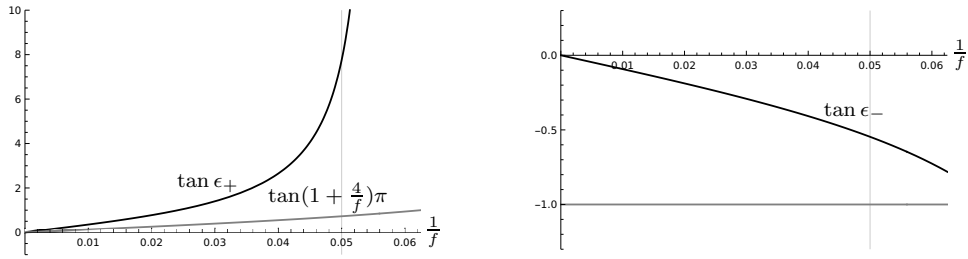

Then we use \eqref{coolsaet_eq2} and \eqref{coolsaet_eq3} to express $\cos a$ as a function of $\frac{1}{f}$. The graph of this function is given by Figure \ref{4ade_geometry3}. We find that the value of the function is always in $(0,1)$. Therefore we can find $0<a<\frac{1}{2}\pi$ satisfying \eqref{coolsaet_eq2} and \eqref{coolsaet_eq3}. By the converse in Lemma \ref{geometry7}, we know the pentagon with the given $\epsilon_-$ (and the corresponding values for the other angles) exist.

\begin{figure}[htp]
\centering
\begin{tikzpicture}[>=latex,scale=1]

\node[scale=0.6] at (0,0) {\includegraphics{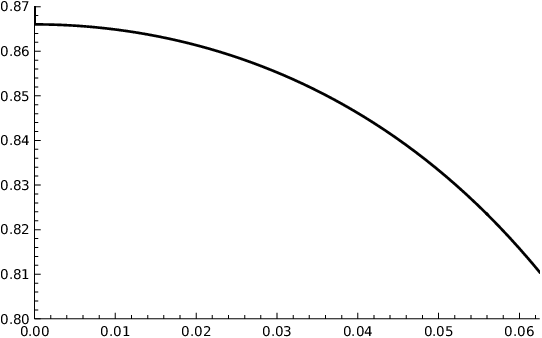}};

\node at (2.9,-1.5) {\scriptsize $\tfrac{1}{f}$};
\node at (0.8,0.95) {\scriptsize $\cos a$};

\end{tikzpicture}
\caption{Proposition \ref{4ade}: $\cos a$ as a function of $\frac{1}{f}$.}
\label{4ade_geometry3}
\end{figure}

It remains to argue that the pentagon is simple. We apply Lemma \ref{geometry8}. The base angles $\theta,\rho$ of the isosceles triangles $\triangle ABD,\triangle ACE$ satisfy $\cos a=\cot\frac{1}{2}\gamma\cot\theta$ and $\cos a=\cot\frac{1}{2}\alpha\cot\rho$. By $\cos a,\cot\frac{1}{2}\alpha, \cot\frac{1}{2}\gamma>0$, we get $\theta,\rho<\frac{1}{2}\pi<\alpha$. Then by Lemma \ref{geometry8}, the pentagon is simple. 
\end{proof}

\end{document}